\documentclass[12pt]{article}

\usepackage{graphicx} 
\usepackage{amsmath,amssymb,amsthm}
\usepackage{color}
\usepackage[nohead,margin=1.0in]{geometry}
\usepackage{booktabs}
\usepackage{threeparttable}
\usepackage{algorithm,algorithmic}
\usepackage{verbatim}
\usepackage{url}
\usepackage{epstopdf}
\usepackage{mathtools}
\usepackage{float}
\allowdisplaybreaks
\usepackage[colorlinks,linktocpage,linkcolor=blue]{hyperref}

\graphicspath{{./pics/}}
\usepackage{subfigure}
\usepackage{array}
\theoremstyle{plain}
\newtheorem{theorem}{Theorem}[section]
\newtheorem{remark}{Remark}[section]

\newtheorem{lemma}{Lemma}[section]
\newtheorem{assumption}{Assumption}[section]

\newtheorem{example}{Example}[section]

\numberwithin{equation}{section}

\renewcommand{\d}{\mathrm{d}}
\def\II{(\Omega)}
\def\bsgamma{\boldsymbol{\gamma}}

\title{Imaging Anisotropic Conductivity from Internal Measurements with Mixed Least-Squares Deep Neural Networks\thanks{The work of B. Jin is supported by Hong Kong RGC General Research Fund (Project
14306423), and a Direct Grant for Research (2024/25) and a start-up fund, both from The Chinese University of Hong Kong. The work of Z. Zhou is supported by Hong Kong
Research Grants Council (15303021) and an internal grant of Hong Kong Polytechnic University (Project ID: P0038888, Work
Programme: 1-ZVX3).}}

\author{Siyu Cen\thanks{Department of Applied Mathematics, The Hong Kong Polytechnic University, Kowloon, Hong Kong, P.R. China. (\texttt{siyu2021.cen@connect.polyu.hk; zhizhou@polyu.edu.hk})}
\and Bangti Jin\thanks{Department of Mathematics, The Chinese University of Hong Kong, Shatin, New Territories, Hong Kong, P.R. China (\texttt{bangti.jin@gmail.com, b.jin@cuhk.edu.hk}).}
\and Xiyao Li\thanks{Department of Computer Science, University College London, Gower Street, London WC1E 6BT, UK (\texttt{xiyao.li.20@ucl.ac.uk})}
\and Zhi Zhou\footnotemark[2]}

\begin{document}

\maketitle

\begin{abstract}
In this work we develop a novel algorithm, termed as mixed least-squares deep neural network (MLS-DNN), to recover an anisotropic conductivity tensor from the internal measurements of the solutions. It is based on applying the least-squares formulation to the mixed form of the elliptic problem, and approximating the internal flux and conductivity tensor simultaneously using deep neural networks. We provide error bounds on the approximations obtained via both population and empirical losses. The analysis relies on the canonical source condition, approximation theory of deep neural networks and statistical learning theory. We also present multiple numerical experiments to illustrate the performance of the method, and conduct a comparative study with the standard Galerkin finite element method and physics informed neural network. The results indicate that the method can accurately recover the anisotropic conductivity in both two- and three-dimensional cases, up to 10\% noise in the data.   
\end{abstract}

\section{Introduction}\label{sec:intro}
The conductivity value of an object varies significantly with the composition and type of materials, and its accurate imaging can provide valuable structural information about the object. This fact underpins several important medical imaging modalities, e.g., electrical impedance tomography (EIT), current density impedance imaging, and acousto-electrical tomography. It is widely accepted that the conductivity values of most biological tissues, e.g., cell membranes, muscle fibers and nerve fiber structures, must be investigated using an anisotropic tensor model. See the works \cite{nicholson1965specific,bronzino2006electrical,ammari2016spectroscopic} for comprehensive overviews on mathematical models and relevant theoretical results for related inverse problems in bio-engineering. In this work, we aim at identifying an anisotropic conductivity tensor $A(x)$ in elliptic problems from internal data using deep neural networks (DNNs). Let $\Omega\subset \mathbb{R}^d$ be a simply connected open bounded domain with a smooth boundary $\partial\Omega$. Consider the following boundary value problem for the function $u$ with a Neumann boundary condition
\begin{equation}\label{eqn:PDE-elliptic}
 \left\{\begin{aligned}
   -\nabla\cdot(A\nabla u)  &= f,&&  \mbox{in }\Omega ,\\
    \vec{n} \cdot(A\nabla u) &=  g,&&\mbox{on }\partial\Omega,
   \end{aligned}\right.
\end{equation}
where $\vec{n} $ denotes the unit outward normal direction to the boundary $\partial\Omega$. The functions {$f\in L^2(\Omega)$ and $g\in H^{-\frac12}(\partial\Omega)$} in \eqref{eqn:PDE-elliptic} are the given source and flux, respectively, and are assumed to satisfy the standard compatibility condition $\int_\Omega f\ {\rm d} x + \int_{\partial\Omega}g\ {\rm d} S = 0$ in order to ensure solvability. 
In addition to the Neumann setting, we also consider the Dirichlet setting of the elliptic problem, i.e., $u=g\ \mbox{on}\ \partial\Omega$. Let $\mathcal{S}_d$ denote the set of all symmetric $d \times d$ matrices equipped with the inner product $A \cdot B = \text{trace}(AB)$, and define the set \begin{equation*}
    \mathcal{M}:=\{A\in \mathcal{S}_d:\  c_0|\xi|^2\le A\xi \cdot \xi \le c_1|\xi|^2\text{ for all }\xi \in \mathbb{R}^{d}  \},
\end{equation*}
with $0<c_0,c_1<\infty$. The anisotropic conductivity  $A$ is assumed to belong to the following admissible set
\begin{equation*}
    \mathcal{K}=\{A:\ A(x)\in\mathcal{M},\text{ a.e. }x\in \Omega  \}.
\end{equation*}
Then the conductivity tensor $A\in L^\infty(\Omega)^{d\times d}$ is uniformly elliptic. This and standard elliptic regularity theory ensure that problem \eqref{eqn:PDE-elliptic} has a unique solution $u\in H^1(\Omega)$ with the normalization condition $\int_\Omega u\ {\rm d}x=0$ \cite[Lemma 2.1]{liu2022imaging}.
We use the notation $u(A)$ to indicate the dependence of the solution $u$ to problem \eqref{eqn:PDE-elliptic} on $A\in\mathcal{K}$.

The concerned inverse problem is to recover the anisotropic conductivity $A$ from suitable
internal observation of the solution $u$. This class of inverse problems has been extensively studied in both mathematics and bio-engineering communities. Given one current density magnitude, the reconstruction problem in a known conformal class was studied in \cite{hoell2014current}. In a series of interesting works, Bal et al investigated the uniqueness of the reconstruction of an anisotropic resistivity from current densities (i.e., $A\nabla u$) \cite{bal2014imaging,bal2014inverse} and power densities (i.e., $A\nabla u_i\cdot \nabla u_j$ for solutions $u_i,u_j$ of problem \eqref{eqn:PDE-elliptic}) \cite{monard2012inverse,monard2013inverse,monard2018imaging}. For example, Bal et al \cite{bal2014inverse} established a minimum number of current densities that is sufficient to guarantee a unique local reconstruction of the conductivity tensor $A$. Ammari et al. \cite{ammari2017determining} studied the recovery of an anisotropic conductivity proportional to a known conductivity tensor with data acquired from diffusion tensor imaging. Theoretically, Hsiao and Sprekels \cite{hsiao1988stability} derived the stability estimate for the reconstruction of matrix-valued coefficients of the form $A = \nabla p \otimes \nabla u$; Rannacher and Vexler \cite{rannacher2005priori} analyzed the reconstruction error of matrix-valued coefficients (with a finite number of unknown parameters) in case of finitely many point-wise observations. Numerically, the reconstruction can be obtained using asymptotic regularization \cite{hoffmann1985identification}, regularized output least-squares approach \cite{Deckelnick:2012,liu2022imaging} and mixed-type formulation approach \cite{KohnLowe:1988}. Decklnick \cite{Deckelnick:2012} investigated the numerical approximation with the Galerkin finite element method (FEM), and established the convergence of discrete approximations using the concept of discrete $H$-convergence. Liu et al \cite{liu2022imaging} revisited the convergence issue of the Galerkin FEM approximation for nonzero Dirichlet boundary conditions. All the above works except \cite{KohnLowe:1988} focus on the Dirichlet boundary conditions.

In this work, we develop a new numerical approach to reconstruct anisotropic conductivity using DNNs. It is based on optimizing a loss defined through a mixed-type least-squares reformulation of the governing equation \eqref{eqn:PDE-elliptic}, with an $L^2(\Omega)^{d,d}$ penalty on the conductivity $A$. The mixed least-squares formulation was first proposed by Kohn and Lowe in \cite{KohnLowe:1988} for identifying isotropic conductivities, where both the conductivity $q$ and the current flux $\sigma  := q\nabla u(q)$ were approximated using the Galerkin FEM. Kohn and Lowe \cite{KohnLowe:1988} extended the mixed formulation to the anisotropic case by introducing variational forms for the least-squares objectives, resulting in a convex functional in the anisotropic conductivity $A$ and the flux $\sigma:= A\nabla u(A) $. However, no  numerical experiments were presented for the reconstruction of $A$. Inspired by the work \cite{KohnLowe:1988}, we approximate both current flux $\sigma$ and conductivity $A$ separately using two DNNs, adopt a least-squares objective for all the equality constraints, and minimize the overall population or empirical loss with respect to the DNN parameters. The use of DNNs in place of the Galerkin FEM effectively exploits the inductive bias and expressivity of DNNs for function approximations, which can be highly beneficial for numerical recovery. By leveraging the existing approximation theory of DNNs \cite{guhring2021approximation} and tools from statistical learning theory \cite{AnthonyBartlett:1999,BartlettMendelson:2002}, we derive
rigorous error bounds on the approximations in terms of the accuracy of
the noisy observational data $z^\delta$, DNN architecture (depth, width and parameter bound), and
the numbers of sampling points in the domain $\Omega$ and on the boundary $\partial\Omega$, etc. This is carried out for both population and empirical losses (resulting from the Monte Carlo quadrature of the integrals). These error bounds provide solid theoretical underpinnings for the approach. In practice, the proposed approach is easy to implement, robust against data noise, and capable of handling inverse problems in three-dimensional spaces. To the best of our knowledge, previous works have not numerically validated the reconstruction of anisotropic conductivities in 3D. These distinct features make the proposed method highly attractive, and the numerical results clearly show its significant potential. The systematic development of the DNN formulations, rigorous error analysis and extensive numerical validation represent the main contributions of this work.

In recent years, the use of DNNs for solving direct and inverse problems for PDEs has received a lot of attention; see \cite{EHanJentzen:2022,TanyuMaass:2022} for overviews. Existing neural inverse schemes using DNNs can roughly be divided into two groups: supervised approaches (see, e.g., \cite{HamiltonHauptmann:2018,KhooYing:2019,WeiChen:2019,GuoJiang:2021,CenJinShin:2023,TanyuNing:2023,ZhouHan:2023jcp}) and unsupervised approaches (see, e.g., \cite{BarSochen:2021,PakravanMistani:2021,XuDarve:2022,JagtapMaoKarniadakis:2022,JinLiLu:2022,PokkunuruKhan:2023,ZhangLiLiu:2023,Jin:2024conductivity}). Supervised methods exploit the availability of (abundant) paired training data to extract problem-specific features, and are concerned with learning forward operators or their (regularized) inverses. See \cite{ChenYang:2023} for one recent mathematical theory of the supervised approach. In contrast, unsupervised approaches exploit extraordinary expressivity of DNNs as universal function approximators for the unknown coefficient and the state, i.e., using DNNs as ansatz functions in the approximation scheme, which enjoy excellent approximation properties for high-dimensional functions (in favorable situations), and the associated inductive biases \cite{Rahaman:2019}. The work \cite{BarSochen:2021} investigated image reconstruction in the classical EIT problem, using the strong formulation (also with the $L^\infty$ norm consistency); see \cite{PokkunuruKhan:2023} for further improvement using an energy-based model. The work \cite{JinLiLu:2022} applied the deep Ritz method to a least-gradient reformulation for the current density impedance imaging, and derived a generalization error for the loss function. The approach performs reasonably well for both full and partial interior current density data, and shows remarkable robustness against data noise. Pakravan et al \cite{PakravanMistani:2021} developed a hybrid approach, blending high expressivity of DNNs with the accuracy and reliability of traditional numerical methods for PDEs, and demonstrated the approach for recovering the isotropic conductivity in one- and two-dimensional elliptic PDEs; see also \cite{CenJinQuanZhou:2023} for a hybrid DNN-FEM approach for recovering the conductivity coefficient in elliptic problems, where the conductivity and state are approximated using DNNs and Galerkin FEM, respectively.
Very recently, Jin et al. \cite{Jin:2024conductivity} developed a novel approach for recovering isotropic conductivity using a least-squares mixed-type reformulation of the governing equations, with DNNs employed as ansatz functions to approximate both the conductivity and flux simultaneously. Our approach follows the unsupervised paradigm, but recovers anisotropic conductivity, thereby extending the work \cite{Jin:2024conductivity} from the case of isotropic conductivity to the anisotropic case. Compared with the work \cite{Jin:2024conductivity}, the anisotropic case is more involved, and requires different techniques in the development of the algorithm as well as the error analysis, due to a lack of stability results comparable with the isotropic case. More precisely,
the key challenges and obstacles in extending the study \cite{Jin:2024conductivity} to the anisotropic
case include: (i) The unavailability of conditional stability estimates for the error analysis, and instead we have to employ the canonical projection type source condition for deriving the error estimate; (ii) We have
to employ a matrix-valued projection operator to preserve the box constraint on the admissible conductivity tensor (and the coercivity of the elliptic operator); (iii) The
formulation involves more unknowns and thus the implementation is much more challenging.

The rest of the paper is organized as follows. In Section \ref{sec:alg}, we develop the reconstruction algorithm for imaging anisotropic conductivity
based on a mixed formulation of the governing equation. In Section \ref{sec:conv}, we conduct an error analysis of the approach for both population and empirical losses. In Section \ref{sec:numer}, we present several numerical experiments to show the effectiveness of the approach.
We conclude the introduction with useful notation. Throughout, we denote the standard matrix norm induced by Euclidean norm for vectors by $\|A\|_{\ell^2}$, and define the $L^{\infty}(\Omega)^{d,d}$ and $L^2(\Omega)^{d,d}$ norms, respectively, by
\begin{equation*}
    \|A\|_{L^{\infty}(\Omega)^{d,d}}:=\sup_{x\in\Omega} \|A(x)\|_{\ell^2}\quad\mbox{and}\quad
    \|A\|_{L^{2}(\Omega)^{d,d}}:=\left(\int_{\Omega} \|A(x)\|^2_{\ell^2}\ \d x\right)^{\frac{1}{2}}.
\end{equation*}
The notation $W^{k,p}(\Omega)$ and $W_0^{k,p}(\Omega)$ denote the standard Sobolev spaces of order $k $ for any integer $k\ge 0$ and real $p\ge 1$, equipped with the norm $\|\cdot\|_{W^{k,p}(\Omega)}$. We also write $H^k(\Omega) $ and $H_0^k(\Omega) $ with the norm $\|\cdot\|_{H^{k}(\Omega)}$ if $p=2$ and write $L^p(\Omega)$ with the norm $ \|\cdot\|_{L^p(\Omega)}$ if $k=0$. The spaces
on the boundary $\partial\Omega$ are defined similarly. The notation $(\cdot,\cdot)$ and $(\cdot,\cdot)_{L^2(\partial\Omega)}$ denote the $L^2(\Omega)$ / $L^2(\Omega)^d$ / $L^2(\Omega)^{d,d}$ and $L^2(\partial\Omega)$ inner products, respectively. For a Banach space $B$, the notation $B^{\otimes d}$ represents the $d$-fold product space. The notation $c$ denotes a generic constant not necessarily the same at each occurrence but it is always independent of the approximation accuracy $\epsilon$ of the given DNN, the noise level $\delta$, the penalty parameters ($\gamma_\sigma$, $\gamma_b$ and $\gamma_q$) and the statistical error ($e_d$, $e_\sigma$, $e_b$, $e_{b'}$ and $e_A$).

\section{Reconstruction algorithm: mixed least-squares DNN}\label{sec:alg}
In this section, we develop the reconstruction algorithm, i.e., mixed least-squares deep neural network (MLS-DNN), using the mixed least-squares formulation and deep neural networks.

\subsection{Deep neural networks}

First we describe notation and properties of the concerned class of deep neural networks (DNNs). Let $\{d_\ell\}_{\ell=0}^L \subset\mathbb{N}$ be fixed at $d_0=d$, and a parameterization $\Theta=
\{({W}^{(\ell)},b^{(\ell)})_{\ell=1}^L\}$ consisting of weight matrices and bias vectors, with
${W}^{(\ell)}=[W_{ij}^{(\ell)}]\in \mathbb{R}^{d_\ell\times d_{\ell-1}}$ and $b^{(\ell)}=
[b^{(\ell)}_i]\in\mathbb{R}^{d_{\ell}}$ being the weight matrix and bias vector at the $\ell$-th layer.
Then a DNN $v_\theta:= v^{(L)}:\Omega\subset\mathbb{R}^d\to\mathbb{R}^{d_L}$
realized by a parameter $\theta\in\Theta$ is defined recursively by
\begin{equation}\label{eqn:NN-realization}
\mbox{DNN realization:}\quad
\left\{\begin{aligned}
v^{(0)}&=x,\quad x\in\Omega\subset\mathbb{R}^d,\\		
v^{(\ell)}&=\rho({W}^{(\ell)}v^{(\ell-1)}+b^{(\ell)}),\quad \ell=1,2,\dots,L-1,\\
		v^{(L)}&={W}^{(L)}v^{(L-1)}+b^{(L)},
	\end{aligned}\right.
\end{equation}
where the nonlinear activation function $\rho:\mathbb{R}\to\mathbb{R}$ is applied
componentwise to a vector, and below we take $\rho\equiv\tanh$: $x\to
\frac{e^x-e^{-x}}{e^x+e^{-x}}$. The DNN $v_\theta$ has a depth $L$ and width ${W}:= \max_{\ell=0,\dots,L}(d_{\ell})$.
Given the parametrization $\Theta$ (i.e., architecture), we denote the associated DNN function class $\mathcal{N}_\Theta$ by $\mathcal{N}_\Theta:=
\{v_\theta:\ \theta\in \Theta\}$.  The
following approximation result holds \cite[Proposition 4.8]{guhring2021approximation}.
\begin{lemma}\label{lem:tanh-approx}
Let $s\in\mathbb{N}\cup\{0\}$ and $p\in[1,\infty]$ be fixed, and $v\in W^{k,p}
(\Omega)$ with $\mathbb{N}\ni k\geq s+1$. Then for any $\epsilon>0$, there exists at least one
 $\theta\in\Theta$ with depth $O(\log(d+k))$ and number of nonzero weights $O
\big(\epsilon^{-\frac{d}{k-s-\mu (s=2)}}\big)$ and weight parameters bounded by $O(\epsilon^{-2-\frac{2(d/p+d+k+\mu(s=2))+d/p+d}{k-s-\mu(s=2)}})$ in the maximum norm, where $\mu>0$ is arbitrarily small, such that the DNN realization $v_\theta\in \mathcal{N}_\Theta$ satisfies
$\|v-v_\theta\|_{W^{s,p}(\Omega)} \leq \epsilon.$
\end{lemma}
We denote the set of DNNs of depth $L$, the number $N$ of nonzero weights, and maximum weight bound $R$ by
\begin{equation*}
    \mathcal{N}(L,N,R) =: \{v_\theta \mbox{ is a DNN of depth }L: \|\theta\|_{\ell^0}\leq N, \|\theta\|_{\ell^\infty}\leq R\},
\end{equation*}
where $\|\cdot\|_{\ell^0}$ and $\|\cdot\|_{\ell^\infty}$ denote, respectively, the number of nonzero entries in and the maximum norm of a vector. Furthermore, for any $\epsilon>0$ and $p\geq 1$, we denote the DNN parameter set by $\mathfrak{P}_{p,\epsilon}$ for
\begin{equation*}
\mathcal{N}\Big(c\log(d+2), c\epsilon^{-\frac{d}{1-\mu}}, c \epsilon^{-2-\frac{4p+3d+3pd}{p(1-\mu)}}\Big).
\end{equation*}

Below, for a vector-valued function, we use DNNs to approximate its components. This can be achieved by
parallelization, which assembles multiple DNNs into one larger DNN.
Moreover,
the new DNN does not change the depth, and its width equals to the sum of that of the subnetworks.

\subsection{MLS-DNN in the Neumann case}\label{ssec:neum}
Now we discuss the inverse conductivity problem for the Neumann problem \eqref{eqn:PDE-elliptic}.
One single internal observation $u$ is insufficient to uniquely determine an anisotropic conductivity $A$, and thus in practice, we take $N$ measurements,
corresponding to $N$ tuples $(f_i,g_i)_{i=1}^N$ of the  source $f$ and boundary flux $g$. Upon letting $\sigma_i  = A\nabla u_i$, we
rewrite problem \eqref{eqn:PDE-elliptic} into $N$ first-order systems
\begin{equation}\label{eqn:PDE-elliptic-system}
 \left\{\begin{aligned}
   \sigma_i &=A\nabla u_i,&&  \mbox{in }\Omega ,\\
   -\nabla\cdot \sigma_i &=f_i, && \mbox{in }\Omega ,\\
     \vec{n}\cdot \sigma_i &=g_i,  &&\mbox{on }\partial\Omega,
   \end{aligned}\right. \qquad i=1,\dots,N.
\end{equation}
To recover the anisotropic conductivity $A$, we employ noisy internal measurements $z_i^\delta$ of the solutions  $u_i^\dag = u_i(A^\dag)$ of problem $\eqref{eqn:PDE-elliptic}$ with the exact conductivity matrix $A^{\dag}$ and excitations $(f_i,g_i)$:
\begin{equation}\label{eqn:noise_level}
    \delta_i := \lVert u_i^\dagger-z_i^{\delta}\rVert_{H^{1}(\Omega)},\quad i=1,\ldots, N,\quad \mbox{and}\quad \delta =\max_{i=1,\ldots,N}\delta_i.
\end{equation}
The noise level $\delta$ imposes a mild regularity condition on the noisy data $z_i^\delta$, which may be obtained by presmoothing the raw data beforehand, e.g., via filtering or denoising. Moreover, let $P_\mathcal{K}$ be the orthogonal projection onto the admissible set $\mathcal{K}$, which is characterized by the following variational inequality
\begin{equation}\label{eqn:prjchar}
    (A-P_{\mathcal{K}}(A)):(B-P_{\mathcal{K}}(A))\le 0,\quad \forall B\in\mathcal{K},
\end{equation}
where $:$ denotes the Frobenius inner product b e tween two matrices.
The orthogonal projection $P_{\mathcal{K}}(A)$ can be implemented pointwise by first computing the eigenvalue decomposition of $A$, i.e., $A=U\Lambda U^\top$, with $\Lambda=\mathrm{diag}(\lambda_i)\in\mathbb{R}^{d\times d}$ being diagonal and $U\in\mathbb{R}^{d\times d}$ being column orthonormal, and then the projection $P_\mathcal{K}(A)$ is given by $P_\mathcal{K}(A)=U\mathrm{diag}(P_{[c_0,c_1]}(\lambda_i)) U^\top$, with $P_{[c_0,c_1]}(\lambda)=\max(c_0,\min(\lambda,c_1))$.

We employ a residual minimization scheme based on \eqref{eqn:PDE-elliptic-system} and approximate both $A$ and $\sigma_i$ using DNNs. For the anisotropic conductivity $A$, there are $d(d+1)/2$ independent entries, and we employ $d(d+1)/2 $ identical DNN function classes (of depth $L_A$ and width $W_A$) with the parameterization $\mathfrak{P}_{\infty,\epsilon_A}$, and stack them into one DNN. Similarly, we stack $Nd$ identical DNN function classes (of depth $L_{\sigma}$ and width $W_{\sigma}$) with the parametrization $\mathfrak{P}_{2,\epsilon_{\sigma}}$ into one DNN to approximate $\{\sigma_i\}_{i=1}^{N}$. Throughout, $\theta$ and $\kappa$ denote the parameters of DNN approximations to $A$ and $\{\sigma_i\}_{i=1}^N$, respectively, and
$(A_\theta,\{\sigma_{i,\kappa}\}_{i=1}^N)$ are the DNN realization of the parameter configuration $(\theta,\kappa)$. Using the least-squares formulation on the equalities in \eqref{eqn:PDE-elliptic-system}, we obtain the population loss
\begin{equation}\label{eqn:loss_NN}
\begin{aligned}
    J_{\gamma}(\theta,\kappa)=\sum_{i=1}^{N}&\left(\|\sigma_{i,\kappa}-P_{\mathcal{K}}(A_{\theta})\nabla z_i^{\delta}\|_{L^2(\Omega)^d}^2+\gamma_{\sigma}\|\nabla\cdot\sigma_{i,\kappa}+f_i\|_{L^2(\Omega)}^2\right.\\
    &\left.+\gamma_b\| \vec{n}\cdot\sigma_{i,\kappa}-g_i \|_{L^2(\partial\Omega)}^2+\gamma_A\|P_{\mathcal{K}}(A_{\theta})\|_{L^2(\Omega)^{d,d}}^2\right).
    \end{aligned}
\end{equation}
The projection $P_{\mathcal{K}}$ explicitly enforces the constraint of the admissible set $\mathcal{K}$.
Then the reconstruction problem at the population level reads
\begin{equation}\label{eqn:obj-Neum}
	\min_{(\theta,\kappa)\in (\mathfrak{P}_{\infty,\epsilon_A}^{\otimes d(d+1)/2},\mathfrak{P}_{2,\epsilon_\sigma}^{\otimes Nd})} J_{\bsgamma}(\theta,\kappa),
\end{equation}
where the superscript $\otimes Nd$ denotes the $Nd$-fold direct product via the parallelization technique, $\gamma_\sigma,\gamma_b,\gamma_A>0$
are penalty parameters, and $\bsgamma=(\gamma_\sigma,\gamma_b,\gamma_q)\in\mathbb{R}_+^3$.
The $L^2(\Omega)$ penalty on $A_\theta$ is to stabilize the minimization process, cf. \cite{Deckelnick:2012,liu2022imaging} (in the sense of H-convergence).
The well-posedness of problem
\eqref{eqn:obj-Neum} follows by a standard argument in calculus of variation. Indeed, the
compactness of the parameterizations $\mathfrak{P}_{\infty,\epsilon_A}^{\otimes d(d+1)/2}$ and $\mathfrak{P}_{2,\epsilon_\sigma}^{\otimes Nd}$
holds due to the uniform boundedness on the parameter vectors and finite-dimensionality of the spaces. The smoothness
of $\tanh$ implies the continuity of the loss $J_{\bsgamma}$
in the DNN parameters $(\theta,\kappa)$. These two properties ensure the existence of a global minimizer $(\theta^*,\kappa^*)$.

The population loss $J_{\gamma}$ involves high-dimensional integrals that need to be appropriately approximated, e.g., using the Monte Carlo method. Let $\mathcal{U}(\Omega)$ and $\mathcal{U}(\partial\Omega)$ be uniform distributions over the domain $\Omega$ and boundary $\partial\Omega$, respectively. Using the expectation $\mathbb{E}_{\mathcal{U}(\cdot)}[\cdot]$ with respect to $\mathcal{U}(\cdot)$ and likewise $\mathbb{E}_{\mathcal{U}(\partial\Omega)}[\cdot]$, we can rewrite the loss $J_{\gamma}$ as
\begin{align*}
    J_{\gamma}(\theta,\kappa)=&\sum_{i=1}^{N}|\Omega|\mathbb{E}_{X\sim\mathcal{U}(\Omega)}\left[\|\sigma_{i,\kappa}(X)-P_{\mathcal{K}}(A_{\theta}(X))\nabla z_i^{\delta}(X)\|_{\ell^2}^2\right]\\
&+\gamma_{\sigma}\sum_{i=1}^{N}|\Omega|\mathbb{E}_{X\sim\mathcal{U}(\Omega)}\left[(\nabla\cdot \sigma_{i,\kappa}(X)+f_i(X))^2\right]\\
    &+\gamma_{b}\sum_{i=1}^{N}|\partial\Omega|\mathbb{E}_{Y\sim\mathcal{U}(\partial\Omega)}\left[(\vec{n}\cdot \sigma_{i,\kappa}(Y)-g_i(Y))^2\right]\\
  &+N\gamma_A|\Omega|\mathbb{E}_{X\sim\mathcal{U}(\Omega)}\left[ \|(P_{\mathcal{K}}A_{\theta})(X)\|_{\ell^2}^2 \right]\\
    =:&\mathcal{E}_{d}(\sigma_{\kappa},A_{\theta})+\gamma_{\sigma}\mathcal{E}_{\sigma}(\sigma_{\kappa})+\gamma_b\mathcal{E}_{b}(\sigma_{\kappa})+\gamma_A\mathcal{E}_{A}(A_{\theta}),
\end{align*}
where the notation $|\Omega|$ and $|\partial\Omega|$ denote the Lebesgue measures of $\Omega$ and $\partial\Omega$. Let
$X=\{X_j\}_{j=1}^{n_r}$ and $Y=\{Y_{j}\}_{j=1}^{n_b}$ be independent and identically
distributed (i.i.d.) samples drawn from $\mathcal{U}(\Omega)$ and $\mathcal{U}(\partial\Omega)$, respectively, where $n_r$ and $n_b$ are the numbers of sampling points in $\Omega$ and on $\partial\Omega$, respectively.
The empirical loss $\widehat{J}_{\bsgamma}(\theta,\kappa)$ reads
\begin{align}\label{eqn:loss_NN_empirical}
    \widehat{J}_{\gamma}(\theta,\kappa)=\widehat{\mathcal{E}}_{d}(\sigma_{\kappa},A_{\theta})+\gamma_{\sigma}\widehat{\mathcal{E}}_{\sigma}(\sigma_{\kappa})+\gamma_b\widehat{\mathcal{E}}_{b}(\sigma_{\kappa})+\gamma_A\widehat{\mathcal{E}}_{A}(A_{\theta}),
\end{align}
where the Monte Carlo approximations $\widehat{\mathcal{E}}_i$, $i\in\{d,\sigma,b,A\}$, are given by
\begin{align*}
    \widehat{\mathcal{E}}_{d}(\sigma_{\kappa},A_{\theta}):=& n_r^{-1}|\Omega|\sum_{i=1}^{N}\sum_{j=1}^{n_r}\|\sigma_{i,\kappa}(X_j)-P_{\mathcal{K}}(A_{\theta}(X_j))\nabla z_i^{\delta}(X_j)\|_{\ell^2}^2,\\
    \widehat{\mathcal{E}}_{\sigma}(\sigma_{\kappa}):=& n_r^{-1}|\Omega|\sum_{i=1}^{N}\sum_{j=1}^{n_r}(\nabla\cdot\sigma_{i,\kappa}(X_j)+f_i(X_j))^2,\\
    \widehat{\mathcal{E}}_{b}(\sigma_{\kappa}):=& n_b^{-1}|\partial\Omega|\sum_{i=1}^{N}\sum_{j=1}^{n_b}(\vec{n}\cdot\sigma_{i,\kappa}(Y_j)-g_i(Y_j))^2,\\
    \widehat{\mathcal{E}}_{A}(A_{\theta}):=& N n_r^{-1}|\Omega|\sum_{j=1}^{n_r}\|(P_{\mathcal{K}}A_{\theta})(X_j)\|_{\ell^2}^2.
\end{align*}
The practical reconstruction scheme reads
\begin{equation}\label{eqn:obj-Neum-emp}
	\min_{(\theta,\kappa)\in (\mathfrak{P}_{\infty,\epsilon_A}^{\otimes d(d+1)/2},\mathfrak{P}_{2,\epsilon_\sigma}^{\otimes Nd})}\widehat{J}_{\bsgamma}(\theta,\kappa).
\end{equation}
In practice, the gradients of the loss $\widehat{J}_{\bsgamma}(\theta,\kappa)$ with respect to the input $x$ (i.e., the spatial derivatives) and the DNN parameters $(\theta,\kappa)$ can be evaluated efficiently using automatic differentiation, and the optimization is commonly done via gradient type algorithms, e.g. Adam \cite{KingmaBa:2015} and L-BFGS \cite{ByrdLu:1995}, which are implemented in many standard software platforms, e.g., PyTorch and TensorFlow.

\subsection{MLS-DNN in the Dirichlet case}\label{ssec:diri}
Now we extend the reconstruction algorithm to the Dirichlet case
\begin{equation}\label{eqn:PDE-diri}
 \left\{\begin{aligned}
   -\nabla\cdot(A\nabla u)  &= f,&&  \mbox{in }\Omega ,\\
    u &=  g,&&\mbox{on }\partial\Omega.
   \end{aligned}\right.
\end{equation}
Like before, we recast problem \eqref{eqn:PDE-diri} into the following first-order system for $i=1,\ldots,N$,
\begin{equation}\label{eqn:PDE-elliptic-system_Diri}
 \left\{\begin{aligned}
   \sigma_i &=A\nabla u_i,&&  \mbox{in }\Omega ,\\
   -\nabla\cdot \sigma_i &=f_i, && \mbox{in }\Omega ,\\
    u_i &=g_i,  &&\mbox{on }\partial\Omega.
   \end{aligned}\right.
\end{equation}
To recover the anisotropic conductivity $A$, we discretize $A$ and $\sigma_i$ by two DNN function classes $\mathfrak{P}_{\infty,\epsilon_A}^{\otimes d(d+1)/2}$ and $\mathfrak{P}_{2,\epsilon_\sigma}^{\otimes Nd} $, and employ a residual fitting scheme based on the first-order system \eqref{eqn:PDE-elliptic-system_Diri} with $u_i$ replaced by noisy measurement $z_i^\delta$:
\begin{equation}\label{eqn:noise_level_Diri}
    \delta_i=\lVert u_i^\dagger(A^\dag )-z_i^{\delta}\rVert_{H^{\frac{3}{2}}(\Omega)},\quad i=1,\ldots, N,\quad \mbox{and} \quad \delta = \max_{i=1,\ldots,N}\delta_i.
\end{equation}
Following Section \ref{ssec:neum}, we employ the following population loss
\begin{equation}\label{eqn:loss_NN_Diri}
    \begin{aligned}
            J_{\gamma}(\theta,\kappa)=\sum_{i=1}^{N}&\left(\|\sigma_{i,\kappa}-P_{\mathcal{K}}(A_{\theta})\nabla z_i^{\delta}\|_{L^2(\Omega)^d}^2+\gamma_{\sigma}\|\nabla\cdot\sigma_{i,\kappa}+f_i\|_{L^2(\Omega)}^2\right.\\
            &\left.+\gamma_b\| \sigma_{i,\kappa}-A^{\dag}\nabla z_i^{\delta} \|_{L^2(\partial\Omega)^d}^2+\gamma_A\|P_{\mathcal{K}}(A_{\theta})\|_{L^2(\Omega)^{d,d}}^2\right).
    \end{aligned}
\end{equation}
Below we denote by $ \mathcal{E}_{b'}(\sigma_{\kappa})= \sum_{i=1}^{N}|\partial\Omega|\mathbb{E}_{Y\sim\mathcal{U}(\partial\Omega)}\left[\| \sigma_{i,\kappa}(Y)-A^\dag(Y)\nabla z_i^\delta(Y)\|_{\ell^2}^2\right]$.
Note that the formulation \eqref{eqn:loss_NN_Diri} requires a knowledge of  $A^\dag|_{\partial\Omega}$. This assumption is frequently made when reconstructing scalar conductivity \cite{Alessandrini:1986,Richter:1981,BarSochen:2021}. Similar to the Neumann case, the population loss \eqref{eqn:loss_NN_Diri} admits a minimizer $(\theta^*,\kappa^*)$ and we denote by $(A_\theta^*,\sigma_{i,\kappa}^*)$ its DNN realization. In practice, the integrals in the loss \eqref{eqn:loss_NN_Diri} are approximated by Monte Carlo methods. Let $X=\{X_j\}_{j=1}^{n_r}\sim \mathcal{U}(\Omega)$ and $Y=\{Y_{j}\}_{j=1}^{n_b}\sim\mathcal{U}(\partial\Omega)$ be i.i.d. samples in the domain $\Omega$ and on the boundary $\partial\Omega$, respectively. Then we obtain the following empirical loss
\begin{align}\label{eqn:loss_NN_empirical_diri}
    \widehat{J}_{\gamma}(\theta,\kappa)=\widehat{\mathcal{E}}_{d}(\sigma_{\kappa},A_{\theta})+\gamma_{\sigma}\widehat{\mathcal{E}}_{\sigma}(\sigma_{\kappa})+\gamma_b\widehat{\mathcal{E}}_{b'}(\sigma_{\kappa})+\gamma_A\widehat{\mathcal{E}}_{A}(A_{\theta}),
\end{align}
where  $\widehat{\mathcal{E}}_{d}(\sigma_{\kappa},A_{\theta}), \widehat{\mathcal{E}}_{\sigma}(\sigma_{\kappa}), \widehat{\mathcal{E}}_{A}(A_{\theta}) $ are defined in Section \ref{ssec:neum}, and
\begin{align*}
     \widehat{\mathcal{E}}_{b'}(\sigma_{\kappa})&= n_b^{-1}|\partial\Omega|\sum_{i=1}^{N}\sum_{j=1}^{n_b}\| \sigma_{i,\kappa}(Y_j)-A^\dag(Y_j)\nabla z_i^\delta(Y_j) \|_{\ell^2}^2.
\end{align*}

\section{Error analysis}\label{sec:conv}

In this section, we provide an error analysis of the reconstruction algorithm for both population and empirical losses. The
analysis below focuses on one single data set $(f,g)$, and
the extension to multiple datasets is direct.

\subsection{Error estimate for the population loss}
First we derive an $L^2(\Omega )$ error bound on the DNN approximation $P_\mathcal{K}(A_\theta^\ast)$  via the population loss $J_\gamma (\theta ,\kappa )$. The following assumption ensures that $A^\dagger$ and $\sigma^\dagger$ have sufficient regularity to be effectively approximated by DNNs. \begin{assumption}\label{ass:regularity}
$A^{\dag}\in W^{2,\infty}(\Omega)^{d\times d}\cap \mathcal{K}$,  $f\in H^1(\Omega)$. Moreover, $g\in H^{\frac{3}{2}}(\partial\Omega)$ in the Neumann case and $g\in H^{\frac{5}{2}}(\partial\Omega)$ in the Dirichlet case. 
\end{assumption}

Let $u^{\dag}:=u(A^{\dag})$ and $\sigma_i^{\dag}:=A^{\dag}\nabla u(A^{\dag})$ below.  Assumption \ref{ass:regularity} and the standard elliptic regularity theory \cite{Grisvard:2011} imply $u^{\dag}\in H^3(\Omega)$ and $\sigma^{\dag}\in H^2(\Omega)^d$.
The next lemma on the minimal value of the population loss $J_{\gamma}(\theta^\ast,\kappa^\ast)$ is crucial to deriving the $L^2(\Omega )$ error bound in Theorem \ref{thm:error_A}.
\begin{lemma}\label{lem:loss_bound}
Let Assumption \ref{ass:regularity} hold. Fix small $\epsilon_A$, $\epsilon_\sigma>0$ and let $(\theta^*,\kappa^*)\in \mathfrak{P}_{\infty,\epsilon_A}^{\otimes d(d+1)/2}\times\mathfrak{P}_{2,\epsilon_\sigma}^{\otimes d}$ be a minimizer of the population loss \eqref{eqn:loss_NN} or \eqref{eqn:loss_NN_Diri}. Then the following estimate holds
    \begin{equation*}
        J_\gamma(\theta^*,\kappa^*)\le
     \left\{\begin{aligned}
                c\left(\epsilon_A^2+(1+\gamma_{\sigma}+\gamma_{b})\epsilon_{\sigma}^2+\delta^2+\gamma_{A}\right), & \quad \mbox{Neumann case},\\
                c\left(\epsilon_A^2+(1+\gamma_{\sigma}+\gamma_{b})\epsilon_{\sigma}^2+(1+\gamma_b)\delta^2+\gamma_{A}\right), &\quad\mbox{Dirichlet case}.\\
      \end{aligned}\right.
    \end{equation*}
\end{lemma}
\begin{proof}
We give the proof only in the Neumann case, since the proof in the Dirichlet case is identical, except using condition \eqref{eqn:noise_level_Diri} to bound the boundary loss $\mathcal{E}_{b'}(\sigma_\kappa)$. First, Assumption \ref{ass:regularity} and the standard elliptic regularity theory \cite{Grisvard:2011} ensures $A^\dag \in W^{2,\infty}(\Omega)^{d\times d}$ and $\sigma^\dag\in H^2(\Omega)^d$. Then Lemma  \ref{lem:tanh-approx} implies that there exists one tuple $(\theta_{\epsilon}, \kappa_{\epsilon})\in\mathfrak{P}_{\infty,\epsilon_A}^{\otimes d(d+1)/2}\times\mathfrak{P}_{2,\epsilon_\sigma}^{\otimes d}$ such that its DNN realization $(A_{\theta_{\epsilon}},\sigma_{\kappa_{\epsilon}})$ satisfies
\begin{equation}\label{eqn:approx_A_sigma}
        \| A^{\dag}-A_{\theta_{\epsilon}}\|_{L^{\infty}(\Omega)^{d,d}}\le \epsilon_A\quad \text{and}\quad \| \sigma^{\dag}-\sigma_{\kappa_{\epsilon}}\|_{H^{1}(\Omega)^d}\le \epsilon_{\sigma}.
\end{equation}
Moreover, by the characterization \eqref{eqn:approx_A_sigma}, we have
\begin{align}\label{eqn:contraction}
        \sup_{x\in\Omega}  \| A^{\dag}(x) - P_{\mathcal{K}} (A_{\theta_{\epsilon}}(x))\|_{\ell^2} \le & \sup_{x\in\Omega}  \| A^{\dag}(x) -  A_{\theta_{\epsilon}}(x) \|_{\ell^2}\le \epsilon_A.
\end{align}
The minimizing property of $(\theta^*,\kappa^*)$ and the triangle inequality imply
\begin{align*}
        &J_\gamma(\theta^*,\kappa^*)\le   J_\gamma(\theta_{\epsilon},\kappa_{\epsilon})\\
        =& \|\sigma_{\kappa_{\epsilon}}-P_{\mathcal{K}}(A_{\theta_{\epsilon}})\nabla z^{\delta}\|_{L^2(\Omega)^d}^2+\gamma_{\sigma}\|\nabla\cdot\sigma_{\kappa_{\epsilon}}+f\|_{L^2(\Omega)}^2\\
        &+\gamma_b\|\vec{n} \cdot\sigma_{\kappa_{\epsilon}}-g \|_{L^2(\partial\Omega)}^2+\gamma_A\|P_{\mathcal{K}}(A_{\theta_{\epsilon}})\|_{L^2(\Omega)^{d,d}}^2\\
        \le & c\big( \|\sigma_{\kappa_{\epsilon}}-\sigma^{\dag} \|_{L^2(\Omega)^d}^2 +\|\sigma^{\dag}-P_{\mathcal{K}}(A_{\theta_{\epsilon}})\nabla u^{\dag}\|_{L^2(\Omega)^d}^2 +\|P_{\mathcal{K}}(A_{\theta_{\epsilon}})(\nabla u^{\dag}-\nabla z^{\delta})\|_{L^2(\Omega)^d}^2\big)\\
    &+\gamma_{\sigma} \|\nabla\cdot\sigma_{\kappa_{\epsilon}}-\nabla\cdot\sigma^{\dag} \|_{L^2(\Omega)}^2+\gamma_{b}\|\vec{n}\cdot \sigma_{\kappa_{\epsilon}}-\vec{n}\cdot \sigma^{\dag}  \|_{L^2(\partial\Omega)}^2+\gamma_{A}\|P_{\mathcal{K}}(A_{\theta_{\epsilon}})\|_{L^2(\Omega)^{d,d}}^2.
\end{align*}
The inequality $\|Au\|_{L^2(\Omega)^{d,d}}        \leq \|A \|_{L^{\infty}(\Omega)^{d,d} }\|u\|_{L^2(\Omega)^d}$,
the estimate \eqref{eqn:approx_A_sigma} and the definition of the noise level $\delta$ yield
\begin{align}
        \|\sigma^{\dag}-P_{\mathcal{K}}(A_{\theta_{\epsilon}})\nabla u^{\dag}\|_{L^2(\Omega)^d}&\le \|A^{\dag}-P_{\mathcal{K}}(A_{\theta_{\epsilon}})\|_{L^{\infty}(\Omega)^{d,d} }\|\nabla u^{\dag}\|_{L^2(\Omega)^d}\le c \epsilon_{A},\label{eqn:stab1}\\
        \|P_{\mathcal{K}}(A_{\theta_{\epsilon}})(\nabla u^{\dag}-\nabla z^{\delta})\|_{L^2(\Omega)^d}&\le \|P_{\mathcal{K}}(A_{\theta_{\epsilon}})\|_{L^{\infty}(\Omega)^{d,d}}\| \nabla(u^{\dag}-z^{\delta})\|_{L^2(\Omega)^d}\le c\delta.\label{eqn:stab2}
\end{align}
Finally, combining the preceding estimates yields the desired estimate.
\end{proof}

To derive an $L^2(\Omega)^{d,d}$ error estimate, we make the following assumption on the exact anisotropic conductivity $A^\dagger$. This condition imposes a certain structural assumption on the exact conductivity $A^\dag$, and it is commonly known as the projected source condition in the literature \cite{Deckelnick:2012}. For two vectors $p,q\in\mathbb{R}^d$, $p\otimes q \in\mathbb{R}^{d\times d}$ denotes the symmetric tensor product, i.e., $(p\otimes q)_{jk}=\frac{1}{2}(p_j q_k+p_k q_j) $, $j,k=1,\dots,d$.
\begin{assumption}\label{ass:coeff}
There exists $\psi\in H_0^1(\Omega)\cap W^{2,\infty}(\Omega)$ such that
\begin{equation}\label{eqn:ass_coeff}
    A^{\dag}(x)=P_{\mathcal{K}}\big(\nabla u^\dag(x)\otimes\nabla\psi(x)\big) \quad \text{a.e. in }\Omega.
\end{equation}
\end{assumption}

Now we can state an error estimate on the approximation $P_\mathcal{K}(A_{\theta}^*)$ using the population loss.
\begin{theorem}\label{thm:error_A}
Let Assumptions \ref{ass:regularity} and \ref{ass:coeff} hold. Fix small $\epsilon_A$, $\epsilon_\sigma>0$, let $(\theta^*,\kappa^*)\in \mathfrak{P}_{\infty,\epsilon_A}^{\otimes d(d+1)/2}\times\mathfrak{P}_{2,\epsilon_\sigma}^{\otimes d} $ be a minimizer of the population loss \eqref{eqn:loss_NN} or \eqref{eqn:loss_NN_Diri} and $(A_\theta^*,\sigma_{\kappa}^*)$ its DNN realization. Then with
$$
\eta:=\left\{\begin{aligned}
(\epsilon_A^2 +(1+\gamma_{\sigma}+\gamma_{b})\epsilon_{\sigma}^2 +\delta^2)^\frac12,&\quad \mbox{Neumann case},\\
(\epsilon_A^2 +(1+\gamma_{\sigma}+\gamma_{b})\epsilon_{\sigma}^2 +(1+\gamma_b)\delta^2)^\frac12,&\quad \mbox{Dirichlet case},
\end{aligned}\right.
$$
the following error estimate holds
\begin{equation*}
   \|A^{\dag}-P_{\mathcal{K}}(A_{\theta}^*)\|_{L^2(\Omega)^{d,d}}\le c(\gamma_A^{-1}\eta^2+\eta +(1+\gamma_\sigma^{-1}+\gamma_b^{-1})\gamma_A)^\frac{1}{2}.
\end{equation*}
\end{theorem}

\begin{proof}
First we prove the Neumann case. Obviously the following identity holds
\begin{align}\label{eqn:A-PA0}
\|A^{\dag}-P_{\mathcal{K}}(A_{\theta}^*)\|_{L^2(\Omega)^{d,d}}^2
=\|P_{\mathcal{K}}(A_{\theta}^*)\|_{L^2(\Omega)^{d,d}}^2+2\left(A^{\dag},A^{\dag}-P_{\mathcal{K}}(A_{\theta}^*)\right)-\|A^{\dag}\|_{L^2(\Omega)^{d,d}}^2.
    \end{align}
The minimizing property $J_{\gamma}(\theta^*,\kappa^*)\le J_{\gamma}(\theta_{\epsilon},\kappa_{\epsilon})$, with $(\theta_{\epsilon},\kappa_{\epsilon})$ from Lemma \ref{lem:loss_bound}, the approximation property \eqref{eqn:approx_A_sigma}, the estimate \eqref{eqn:contraction} and the argument of Lemma \ref{lem:loss_bound} imply
\begin{align}
&\gamma_A \|P_{\mathcal{K}}(A_\theta^*)\|_{L^2(\Omega)^{d,d}}^2+ \|\sigma_{\kappa}^*-P_{\mathcal{K}}(A_{\theta}^*)\nabla z^{\delta}\|_{L^2(\Omega)^d}^2\nonumber\\
&\quad +\gamma_{\sigma}\|\nabla\cdot\sigma_{\kappa}^*+f\|_{L^2(\Omega)}^2+\gamma_b\|\vec{n} \cdot\sigma_{\kappa}^*-g \|_{L^2(\partial\Omega)}^2\nonumber\\
\leq&  \|\sigma_{\kappa_{\epsilon}}-P_{\mathcal{K}}(A_{\theta_{\epsilon}})\nabla z^{\delta}\|_{L^2(\Omega)^d}^2+\gamma_{\sigma}\|\nabla\cdot\sigma_{\kappa_{\epsilon}}\nonumber\\
&\quad +f\|_{L^2(\Omega)}^2+\gamma_b\| \vec{n}\cdot(\sigma_{\kappa_{\epsilon}}-\sigma^\dag) \|_{L^2(\partial\Omega)}^2+ \gamma_A\|P_{\mathcal{K}} (A_{\theta_{\epsilon}})\|_{L^2(\Omega)^{d,d}}^2\nonumber\\
       \le & c \left(\epsilon_A^2+(1+\gamma_{\sigma}+\gamma_{b})\epsilon_{\sigma}^2+\delta^2 \right)+c\gamma_A\epsilon_A+\gamma_A\|A^\dag\|_{L^2(\Omega)^{d,d}}^2\nonumber\\
       \le&  c\eta^2+c\gamma_A \epsilon_A+\gamma_A\|A^{\dag}\|_{L^2(\Omega)^{d,d}}^2.\label{eqn:bdd-PA}
\end{align}
From the source representation of $A^\dag$ in Assumption \ref{ass:coeff}, we derive
\begin{align*}
    &\big(A^{\dag},A^{\dag}-P_{\mathcal{K}}(A_{\theta}^*)\big)=\big((A^{\dag}-P_{\mathcal{K}}(A_{\theta}^*)) \nabla u^{\dag},\nabla \psi\big)
    =\big(\sigma^{\dag},\nabla \psi\big) -\big(P_{\mathcal{K}}(A_{\theta}^*)\nabla u^{\dag},\nabla \psi\big).
\end{align*}
By integration by parts, we deduce
\begin{align*}
  \big(\sigma^{\dag},\nabla \psi\big) =& -\big(\nabla\cdot\sigma^{\dag},\psi\big)+\big(\vec{n}\cdot \sigma^{\dag},\psi\big)_{L^2(\partial\Omega)} = (f,\psi) +\big(g,\psi\big)_{L^2(\partial\Omega)}\\
  = & \big(f+\nabla\cdot\sigma_{\kappa}^*,\psi\big) +\big( \sigma_{\kappa}^*,\nabla \psi\big) -\big(\vec{n}\cdot \sigma_{\kappa}^*,\psi\big)_{L^2(\partial\Omega)}
        +\big(g,\psi\big)_{L^2(\partial\Omega)}.
\end{align*}
Consequently,
\begin{align*}
  \big(A^{\dag},A^{\dag}-P_{\mathcal{K}}(A_{\theta}^*)\big) =&\big(f+\nabla\cdot\sigma_{\kappa}^*,\psi\big)
        +\big( \sigma_{\kappa}^*-P_{\mathcal{K}}(A_{\theta}^*)  \nabla z^{\delta},\nabla \psi\big)\\
        &+\big(P_{\mathcal{K}}(A_{\theta}^*)  (\nabla z^{\delta}- \nabla u^{\dag}),\nabla \psi\big)
        +\big(g-\vec{n}\cdot \sigma_{\kappa}^*,\psi\big)_{L^2(\partial\Omega)}.
\end{align*}
Since $\|\psi\|_{H^1(\Omega)}$ is bounded, we have
    \begin{align}\label{eqn:est-A-PA}
    \left|\big(A^{\dag},A^{\dag}-P_{\mathcal{K}}(A_{\theta}^*) \big)\right| \le & c \|f+\nabla\cdot \sigma_\kappa^* \|_{L^2\II}+c\| \sigma_\kappa^*-P_\mathcal{K}(A_\theta^*)\nabla z^{\delta} \|_{L^2\II}\nonumber\\
    &+c\|  P_\mathcal{K}(A_\theta^*)(\nabla z^\delta-\nabla u^\dag) \|_{L^2\II}+c\| g-\vec{n}\cdot \sigma_\kappa^* \|_{L^2(\partial\Omega)}.
\end{align}
Combing the identity \eqref{eqn:A-PA0} with the estimates \eqref{eqn:bdd-PA} and \eqref{eqn:est-A-PA} yields
\begin{align*}
    &\gamma_A \|A^{\dag}-P_{\mathcal{K}}(A_{\theta}^*)\|_{L^2(\Omega)^{d,d}}^2+ \|\sigma_{\kappa}^*-P_{\mathcal{K}}(A_{\theta}^*)\nabla z^{\delta}\|_{L^2(\Omega)^d}^2\\
&\quad +\gamma_{\sigma}\|\nabla\cdot\sigma_{\kappa}^*+f\|_{L^2(\Omega)}^2+\gamma_b\|\vec{n} \cdot\sigma_{\kappa}^*-g \|_{L^2(\partial\Omega)}^2\\
    \le &c\eta^2+c\gamma_A \eta +c \gamma_A\left(\|\sigma_{\kappa}^*-P_{\mathcal{K}}(A_{\theta}^*)\nabla z^{\delta}\|_{L^2(\Omega)^d}\right.\\
&\quad \left.+\|\nabla\cdot\sigma_{\kappa}^*+f\|_{L^2(\Omega)}+\|\vec{n} \cdot\sigma_{\kappa}^*-g \|_{L^2(\partial\Omega)}\right)\\
    \le & c\eta^2+c\gamma_A\eta+c(1+\gamma_\sigma^{-1}+\gamma_b^{-1})\gamma_A^2 +\tfrac{1}{2}\big(\|\sigma_{\kappa}^*-P_{\mathcal{K}}(A_{\theta}^*)\nabla z^{\delta}\|_{L^2(\Omega)^d}^2\\
&+\gamma_{\sigma}\|\nabla\cdot\sigma_{\kappa}^*+f\|_{L^2(\Omega)}^2+\gamma_b\|\vec{n} \cdot\sigma_{\kappa}^*-g \|_{L^2(\partial\Omega)}^2\big).
\end{align*}
Thus we conclude 
\begin{equation*}
     \|A^{\dag}-P_{\mathcal{K}}(A_{\theta}^*)\|_{L^2(\Omega)^{d,d}}\le c(\gamma_A^{-1}\eta^2+\eta +(1+\gamma_\sigma^{-1}+\gamma_b^{-1})\gamma_A)^\frac{1}{2},
\end{equation*}
i.e., the assertion in the Neumann case. 
Next, we prove the Dirichlet case. The only difference lies in bounding $\big(A^{\dag},A^{\dag}-P_{\mathcal{K}}(A_{\theta}^*)\big)$.
Indeed, Assumption \ref{ass:coeff} and integration by parts give
\begin{align*}
   \big(A^{\dag},A^{\dag}-P_{\mathcal{K}}(A_{\theta}^*)\big)
   =&\big(f+\nabla\cdot\sigma_{\kappa}^*,\psi\big)
   +\big( \sigma_{\kappa}^*-P_{\mathcal{K}}(A_{\theta}^*) \nabla z^{\delta},\nabla \psi\big)
   +\big(P_{\mathcal{K}}(A_{\theta}^*) \nabla (z^{\delta}- u^{\dag}),\nabla \psi\big)\\
    &+\big(\vec{n}\cdot (\sigma^{\dag}-A^{\dag}\nabla z^{\delta}+A^{\dag}\nabla z^{\delta}-\sigma_{\kappa}^*),\psi\big)_{L^2(\partial\Omega)}.
\end{align*}
Then Lemma \ref{lem:loss_bound} and the trace inequality imply again an estimate analogous to \eqref{eqn:est-A-PA}, and completes the proof of the theorem. 
\end{proof}
\begin{remark}\label{rem:error_A}
Theorem \ref{thm:error_A} suggests choosing the algorithmic parameters $ \gamma_{\sigma},\gamma_b\sim \mathcal{O}(1)$ and $\gamma_A\sim \mathcal{O}(\delta)$ in the loss $J_\gamma$, and $\epsilon_A=\epsilon_\sigma=O(\delta)$. With these choices, we obtain a convergence rate $O(\delta^{\frac12})$ for the approximation $P_\mathcal{K}(A_\theta^*)$. Note that this rate is comparable with that for the continuous output least-squares formulation or its Galerkin FEM approximation \cite[Theorem 4.2]{Deckelnick:2012}.
\end{remark}

\subsection{Error estimate for the empirical loss}
Now we analyze the approximation $P_\mathcal{K}(\widehat{A}_\theta^*)$ (i.e., the DNN realization of the minimizer $(\widehat{\theta}^*,\widehat{\kappa}^*)$ to the empirical loss $\widehat{J}_{\bsgamma}(\theta,\kappa)$). Compared with the population loss $J_{\bsgamma}$, the empirical loss $\widehat{J}_{\bsgamma}$ additionally involves the quadrature error due to Monte Carlo methods. The following regularity assumption is needed for applying the standard Rademacher complexity argument \cite{AnthonyBartlett:1999,BartlettMendelson:2002}; see, e.g., \cite{JiaoLaiLoWang:2024} for relevant studies in the context of direct problems.
\begin{assumption}\label{ass:regularity_high}
$f\in L^\infty(\Omega)$, and $g\in L^{\infty}(\partial\Omega)$. Moreover, $z^\delta\in W^{1,\infty}(\Omega)$  and $z^\delta \in W^{1,\infty}(\Omega)\cap W^{1,\infty}(\partial\Omega)$ in the Neumann and Dirichlet cases, respectively.
\end{assumption}
The key of the analysis is to bound the following statistical error:
\begin{equation*}
    \sup_{ \theta \in \mathfrak{P}_{\infty,\epsilon_A}^{\otimes d(d+1)/2},  \kappa \in\mathfrak{P}_{2,\epsilon_\sigma}^{\otimes d} } |\widehat{J}_{\gamma}(\theta,\kappa)-J_{\gamma}(\theta,\kappa) |\le \Delta\mathcal{E}_d+ \gamma_\sigma\Delta\mathcal{E}_\sigma+ \gamma_b\Delta\mathcal{E}_b+ \gamma_A\Delta\mathcal{E}_A,
\end{equation*}
with the error components given, respectively, by
\begin{align*}
    &\Delta\mathcal{E}_d=\sup_{ \theta \in \mathfrak{P}_{\infty,\epsilon_A}^{\otimes d(d+1)/2},  \kappa \in\mathfrak{P}_{2,\epsilon_\sigma}^{\otimes d } } |\widehat{\mathcal{E}}_d(\sigma_{\kappa},A_{\theta})-\mathcal{E}_{d}(\sigma_{\kappa},A_{\theta}) |, &&\Delta\mathcal{E}_\sigma=\sup_{   \kappa \in\mathfrak{P}_{2,\epsilon_\sigma}^{\otimes d } } |\widehat{\mathcal{E}}_\sigma(\sigma_{\kappa})-\mathcal{E}_{\sigma}(\sigma_{\kappa}) |,\\
    &\Delta\mathcal{E}_i=\sup_{   \kappa \in\mathfrak{P}_{2,\epsilon_\sigma}^{\otimes d } } |\widehat{\mathcal{E}}_i(\sigma_{\kappa})-\mathcal{E}_{i}(\sigma_{\kappa}) |,\quad i\in\{b,b'\},&&\Delta\mathcal{E}_A=\sup_{ \theta \in \mathfrak{P}_{\infty,\epsilon_A}^{\otimes d(d+1)/2} } |\widehat{\mathcal{E}}_A( A_{\theta})-\mathcal{E}_{A}(A_{\theta}) |.
\end{align*}
Now we give the quadrature error for each term in high probability. The proof follows identically as that of \cite[Theorem 3.4]{Jin:2024conductivity}, using Rademacher complexity, Dudley's lemma and the continuity of the DNN functions with respect to their parameters, and hence the detailed proof is omitted.
\begin{lemma}\label{lem:error_stat_bound}
Let Assumption \ref{ass:regularity_high} hold, $A_{\theta}\in \mathcal{N}(L_{A},N_{A},R_{A})$, and $\sigma_{\kappa}\in \mathcal{N}(L_{\sigma},N_\sigma,R_{\sigma})$. Fix $\tau\in (0,\frac{1}{8})$ and let
     \begin{align*}
         e_d=& c n_{r}^{-\frac{1}{2}}R_\sigma^2 N_\sigma^2 (N_{\sigma}+N_A)^{\frac{1}{2}} \big(\log R_\sigma+\log N_{\sigma}+\log R_A+\log N_A  +\log n\big)^{\frac{1}{2}}\\
           &+ \tilde c n_r^{-\frac12} R_{\sigma}^2 N_{\sigma}^2|\log \tau|^{\frac12},\\
         e_\sigma=& c n_{r}^{-\frac{1}{2}}R_\sigma^{2L_{\sigma}} N_\sigma^{2L_{\sigma}-\frac{3}{2}}   \big(\log R_\sigma+\log N_{\sigma}+ \log n\big)^{\frac{1}{2}}+ \tilde cn_r^{-\frac12} R_{\sigma}^{2L_{\sigma}} N_{\sigma}^{2L_{\sigma}-2}|\log \tau|^{\frac12},\\
         e_b=&e_{b'}= c n_{b}^{-\frac{1}{2}}R_\sigma^{2 } N_\sigma^{ \frac{5}{2}} \big(\log R_\sigma+\log N_{\sigma}+\log n\big)^{\frac{1}{2}}+ \tilde c n_b^{-\frac12}R_{\sigma}^{2 } N_{\sigma}^{2 }|\log \tau|^{\frac12},\\
         e_A=& c n_{r}^{-\frac{1}{2}} N_A^{ \frac{1}{2}}   \big(\log R_A+\log N_A+\log n\big)^{\frac{1}{2}}+ \tilde cn_r^{-\frac12} |\log \tau|^{\frac12},
     \end{align*}
where the constants $c $ and $\tilde c$ may depend on $|\Omega|$, $|\partial\Omega|$, $c_1$, $d$, $\|z^{\delta}\|_{W^{1,\infty}(\Omega)}$ $\|f\|_{L^{\infty}(\Omega)}$ and $\|g\|_{L^{\infty}(\partial\Omega)}$ at most polynomially. Then, with probability at least $1-\tau$, each of the following statements holds
     \begin{equation*}
         \Delta\mathcal{E}_i\le e_i,\quad i\in\{d,\sigma,b,b',A\}.
     \end{equation*}
\end{lemma}

The next result gives the error estimate on the numerical approximation $P_\mathcal{K}(\widehat{A}_\theta^*)$. The result indicates that with sufficiently large numbers   $n_r,n_b$ of sampling points in the domain $\Omega$ and on the boundary $\partial\Omega$, the overall approximation error is $O(\delta^{1/2})$, which is comparable with that for the population loss in Theorem \ref{thm:error_A}.
\begin{theorem}\label{thm:error_A_empirical}
Let Assumptions \ref{ass:regularity}--\ref{ass:regularity_high} hold. Fix small $\epsilon_A$, $\epsilon_\sigma>0$ and let $(\widehat{\theta}^*,\widehat{\kappa}^*)\in \mathfrak{P}_{\infty,\epsilon_A}^{\otimes d(d+1)/2}\times\mathfrak{P}_{2,\epsilon_\sigma}^{\otimes d}$ be a minimizer of the empirical loss \eqref{eqn:loss_NN_empirical} or \eqref{eqn:loss_NN_empirical_diri}, and $(\widehat{A}_{\theta}^*,\widehat{\sigma}_{\kappa}^*)$ its NN realization. Fix $\tau\in (0,\frac{1}{8})$, let the bounds $e_d$, $e_\sigma$, $e_b$ and $e_A$ be defined in Lemma \ref{lem:error_stat_bound} and let
$$
\eta=\left\{\begin{aligned}
    (\epsilon_A^2+(1+\gamma_\sigma+\gamma_b)\epsilon_{\sigma}^2+\delta^2+e_d+\gamma_\sigma e_\sigma+\gamma_b e_b+\gamma_A e_A)^\frac12, & \quad \mbox{Neumann case},\\
    (\epsilon_A^2+(1+\gamma_\sigma+\gamma_b)\epsilon_{\sigma}^2+(1+\gamma_b)\delta^2+e_d+\gamma_\sigma e_\sigma+\gamma_b e_{b}+\gamma_A e_A)^\frac12,& \quad\mbox{Dirichlet case}.
    \end{aligned}\right.
$$ Then with probability at least $1-4\tau$, there holds
    \begin{equation*}
        \|A^{\dag}-P_{\mathcal{K}}(\widehat{A}_{\theta}^*)\|_{L^2(\Omega)^{d,d}}\le c\big(\gamma_A^{-1}\eta^2+\eta+(1+\gamma_\sigma^{-1}+\gamma_b^{-1})\gamma_A\big)^\frac12.
    \end{equation*}
\end{theorem}
\begin{proof}
We prove only the Neumann case. 
The proof is again based on the identity \eqref{eqn:A-PA0}, with $P_{\mathcal{K}}(\widehat{A}_{\theta}^*)$ in place of $P_{\mathcal{K}}({A}_{\theta}^*)$.
Lemma \ref{lem:error_stat_bound}  implies that with probability at least $1-4\tau$,
\begin{align*}
    J_\gamma(\widehat{\theta}^*,\widehat{\kappa}^*) 
 =&\mathcal{E}_{d}(\widehat\sigma_{\kappa}^*,\widehat{A}_{\theta}^*)+\gamma_{\sigma}\mathcal{E}_{\sigma}(\widehat{\sigma}^*_{\kappa})+\gamma_b{\mathcal{E}}_{b}(\widehat{\sigma}^*_{\kappa})+\gamma_A{\mathcal{E}}_{A}(\widehat{A}^*_{\theta}) \\
 =& \big(\widehat{\mathcal{E}}_{d}(\widehat\sigma_{\kappa}^*,\widehat{A}_{\theta}^*)+{\mathcal{E}}_{d}(\widehat\sigma_{\kappa}^*,\widehat{A}_{\theta}^*)- \widehat{\mathcal{E}}_{d}(\widehat\sigma_{\kappa}^*,\widehat{A}_{\theta}^*)\big)+\gamma_{\sigma}\big(\widehat{\mathcal{E}}_{\sigma}(\widehat{\sigma}^*_{\kappa})+{\mathcal{E}}_{\sigma}(\widehat{\sigma}^*_{\kappa})-\widehat{\mathcal{E}}_{\sigma}(\widehat{\sigma}^*_{\kappa})\big)\\
 &+\gamma_b(\widehat{\mathcal{E}}_{b}(\widehat{\sigma}^*_{\kappa})+{\mathcal{E}}_{b}(\widehat{\sigma}^*_{\kappa})-\widehat{\mathcal{E}}_{b}(\widehat{\sigma}^*_{\kappa})\big)+\gamma_A(\widehat{\mathcal{E}}_{A}(\widehat{A}^*_{\theta})+{\mathcal{E}}_{A}(\widehat{A}^*_{\theta})-\widehat{\mathcal{E}}_{A}(\widehat{A}^*_{\theta})\big)\\
        \le & \widehat{J}_{\gamma}(\widehat\theta^*,\widehat\kappa^*)+(e_d+\gamma_\sigma e_\sigma + \gamma_be_b+\gamma_Ae_A).
\end{align*}
By the minimizing property of $(\widehat\theta^*,\widehat\kappa^*)$ and Lemma \ref{lem:error_stat_bound}, with probability at least $1-4\tau$, there holds
\begin{align}\label{eqn:loss_bound_empirical}
       \widehat{J}_{\gamma}(\widehat{\theta}^*,\widehat{\kappa}^*)\le & [ \widehat{J}_{\gamma}({\theta}^*,{\kappa}^*)-{J}_{\gamma}({\theta}^*,{\kappa}^*)]+{J}_{\gamma}({\theta}^*,{\kappa}^*)\nonumber\\
       \le& {J}_{\gamma}({\theta}^*,{\kappa}^*)+\sup_{ (\theta,\kappa)\in \mathfrak{P}_{\infty,\epsilon_A}^{\otimes d(d+1)/2}\times\mathfrak{P}_{2,\epsilon_\sigma}^{\otimes d } } |J_{\gamma}(\theta,\kappa)-\widehat{J}_{\gamma}(\theta,\kappa) |\nonumber\\
       \le & {J}_{\gamma}({\theta}^*,{\kappa}^*)+(e_d+\gamma_\sigma e_\sigma+\gamma_b e_b+\gamma_A e_A).
   \end{align}
Meanwhile, the minimizing property of $({\theta}^*,{\kappa}^*)$, the construction of the tuple $(\theta_\epsilon,\kappa_\epsilon)$ in the proof of Lemma \ref{lem:loss_bound} and the argument of Lemma \ref{lem:loss_bound} imply
\begin{align*}
  {J}_{\gamma}({\theta}^*,{\kappa}^*) & \leq {J}_{\gamma}(\theta_{\epsilon},\kappa_{\epsilon})\leq c(\epsilon_A^2+(1+\gamma_\sigma+\gamma_b)\epsilon_{\sigma}^2+\delta^2) + \gamma_A\|(P_{\mathcal{K}}A_{\theta_\epsilon})\|_{L^2(\Omega)^{d,d}}^2.
\end{align*}
The last three estimates and the definition of $\eta$ yield
\begin{align*}
  J_\gamma(\widehat{\theta}^*,\widehat{\kappa}^*) \leq c\eta^2 + \gamma_A\|(P_{\mathcal{K}}A_{\theta_\epsilon})\|_{L^2(\Omega)^{d,d}}^2.
\end{align*}
Now with the splitting \eqref{eqn:A-PA0} and the estimate \eqref{eqn:contraction}, we obtain
 \begin{align*}
       &\gamma_A \|A^{\dag}-P_{\mathcal{K}}(\widehat{A}_{\theta}^*)\|_{L^2(\Omega)^{d,d}}^2 + \| \widehat{\sigma}_{\kappa}^*-P_{\mathcal{K}}(\widehat{A}_{\theta}^*)  \nabla z^{\delta} \|_{L^2(\Omega)^d}^2\\
       &\quad +\gamma_\sigma\| f+\nabla\cdot\widehat{\sigma}_{\kappa}^*\|_{L^2(\Omega)}^2+\gamma_b \|\vec{n}\cdot (\sigma^{\dag}-\widehat{\sigma}_{\kappa}^* ) \|_{L^2(\partial\Omega)}^2\\
        \le& c(\eta^2 + \gamma_A\epsilon_A) +2\gamma_A\big(A^{\dag},A^{\dag}-P_{\mathcal{K}}(\widehat{A}_{\theta}^*)\big).
\end{align*}
Repeating the argument of Theorem \ref{thm:error_A} yields
\begin{align*}
   \big(A^{\dag},A^{\dag}-P_{\mathcal{K}}(\widehat{A}_{\theta}^*)\big)=&\big(f+\nabla\cdot\widehat{\sigma}_{\kappa}^*,\psi\big) +\big( \widehat{\sigma}_{\kappa}^*-P_{\mathcal{K}}(\widehat{A}_{\theta}^*)  \nabla z^{\delta},\nabla \psi\big)\\
  &+\big(P_{\mathcal{K}}(\widehat{A}_{\theta}^*) (\nabla z^{\delta}- \nabla u^{\dag}),\nabla \psi\big)
   +\big(\vec{n}\cdot (\sigma^{\dag}-\widehat{\sigma}_{\kappa}^*),\psi\big)_{L^2(\partial\Omega)}.
    \end{align*}
Thus, upon noting the $H^1(\Omega)$ bound on $\psi$, it follows directly that 
\begin{align*}
          & \gamma_A \|A^{\dag}-P_{\mathcal{K}}(\widehat{A}_{\theta}^*)\|_{L^2(\Omega)^{d,d}}^2 + \| \widehat{\sigma}_{\kappa}^*-P_{\mathcal{K}}(\widehat{A}_{\theta}^*)  \nabla z^{\delta} \|_{L^2(\Omega)^d}^2\\
          &\quad +\gamma_\sigma\| f+\nabla\cdot\widehat{\sigma}_{\kappa}^*\|_{L^2(\Omega)}^2 +\gamma_b \|\vec{n}\cdot (\sigma^{\dag}-\widehat{\sigma}_{\kappa}^* ) \|_{L^2(\partial\Omega)}^2\\
        \le &c(\eta^2 + \gamma_A\epsilon_A) +c\gamma_A\big(\|f+\nabla\cdot\widehat{\sigma}_{\kappa}^*\|_{L^2(\Omega)} +\|\widehat{\sigma}_{\kappa}^*-P_{\mathcal{K}}(\widehat{A}_{\theta}^*)  \nabla z^{\delta}\|_{L^2(\Omega)^d}\\
  &+\|P_{\mathcal{K}}(\widehat{A}_{\theta}^*) (\nabla z^{\delta}- \nabla u^{\dag})\|_{L^2(\Omega)^d}+\|\vec{n}\cdot (\sigma^{\dag}-\widehat{\sigma}_{\kappa}^*)\|_{L^2(\partial\Omega)}\big).
\end{align*}
By Young's inequality, we deduce 
\begin{equation*}
     \|A^{\dag}-P_{\mathcal{K}}(\widehat{A}_{\theta}^*)\|_{L^2(\Omega)^{d,d}}\leq (c\gamma_A^{-1}\eta^2+\eta + (1+\gamma_\sigma^{-1}+\gamma_b^{-1})\gamma_A)^{\frac12}. 
\end{equation*}
This proves the assertion in the Neumann case. The proof in the Dirichlet case is similar and hence omitted.
\end{proof}

\section{Numerical experiments and discussions}\label{sec:numer}
Now we showcase the performance of the proposed mixed least-squares DNN approach. All computations are carried out on the software platform TensorFlow
1.15.0 using Intel Core i7-11700K Processor with 16 CPUs. We measure the accuracy of a reconstruction $\hat A$ (with respect to the exact one $A^\dag$) by the relative $L^2(\Omega)^{d,d}$ error $e(\hat A)$ defined by
$$e(\hat A)=\|A^\dag-\hat A\|_{L^2(\Omega)^{d,d}}/\|A^\dag\|_{L^2(\Omega)^{d,d}},$$ 
computed by the equivalent $L^2(\Omega)^{d,d}$ norm defined by
$\|A\|_{L^{2}(\Omega)^{d,d}}:=(\int_{\Omega} \|A(x)\|^2_{F}\ \d x)^{\frac{1}{2}}$, where $\|A\|_F$ denotes the Frobenius norm of a matrix $A$. For an elliptic problem in $\mathbb{R}^d$, we use DNNs with an output dimension $d(d+1)/2$ and $d$ to approximate the conductivity
$A$ and flux $\sigma$, respectively. Unless otherwise stated, the DNN for $\sigma$ has 4 hidden layers with 26, 26, 26, and 10
neurons on each layer (i.e., depth $L_{\sigma}=5$, width $W_{\sigma}=26$), and the DNN for $A$ has 9 hidden layers with 32 neurons on each layer (i.e., depth $L_{A}=10$, width $W_{A}=32$). For Example \ref{exam:neu2d1} and \ref{exam:diri2d1}, $L_A$ is taken to be 5. Increasing the depth / width can further enhance the expressivity of the DNN, but also increases the total number of DNN parameters to be learned. Empirically this choice of depth and width balances well the expressivity and computational efficiency of the DNN. The numbers of sampling points in the domain $\Omega$ and
on the boundary $\partial\Omega$ are denoted by $n_r$ and $n_b$, respectively. The penalty parameters $\gamma_\sigma$, $\gamma_b$ and $\gamma_A$ are for the divergence term, boundary term and $L^2(\Omega)^{d,d}$ penalty
term, respectively, in the losses \eqref{eqn:loss_NN} and \eqref{eqn:loss_NN_Diri}. The empirical loss $\widehat{J}_{\bsgamma}$
is  minimized by the popular ADAM \cite{KingmaBa:2015}. In the experiments, we adopt a stagewise decaying learning rate schedule, which is determined by the starting learning rate (lr), decay rate (dr) and the epoch index at which the decay takes place (step). The term epoch refers to the total number of epochs used for obtaining the reconstruction. The hyper-parameters of the DNN approach include penalty parameters and optimizer parameters, and are determined in a trial-and-error manner. That is, the algorithm is run with different hyper-parameter settings, and then the set of hyper-parameters attaining the smallest reconstruction error $e(\hat A)$ is reported. This strategy is often used for tuning hyper-parameters in deep learning \cite{YuZhu:2020}. A systematic study of the hyper-parameter selection is beyond the scope of this work. Table \ref{tab:algpara} summarizes the algorithmic parameters for the experiments where the numbers in the brackets indicate the parameters used for noisy data. For most examples, $A^\dag$ and $u^\dag$ are known explicitly and the values of $\nabla u^\dag$, $f$ and $g$ at the sampling points are computed directly from the analytic expressions. Then i.i.d. Gaussian noise is added pointwise to $\nabla u^\dag$ at the sampling points drawn uniformly in the domain (and also on the boundary in the Dirichlet case) to form the observation 
\begin{align*}
    \nabla z^\delta(x)=\nabla u^\dag(x)+\delta\cdot\iota(x)\cdot\max_{x\in\Omega}\|\nabla u^\dag(x)\|_{\ell^\infty},
\end{align*}
where $\iota(x)\in\mathbb{R}^d$ follows the standard
Gaussian distribution (with zero mean and the covariance being the identity matrix), and $\delta>0$ denotes the relative noise level (slightly abusing the notation). For Example \ref{exam:neu2d4}, the solution $u^\dag$ and the gradient $\nabla u^\dag$ are approximated using the standard FEM solver FreeFEM++ \cite{Hecht:2012}, and the observation $\nabla z^\delta$ is obtained by adding Gaussian noise pointwise. For the purpose of comparison, we shall also present reconstruction results by the Galerkin FEM and physics informed neural network for two 2D examples with partial internal data. The FEM reconstruction is obtained using a $50\times 50$ mesh, with the spaces $[P_1]^3$ and $[P_1]^2$ approximating the conductivity tensor $A$ and the internal flux $\sigma$, respectively. The resulting optimization problem is solved using the conjugate-gradient method \cite[Section 5]{Nocedal:2006}. The Python source code for reproducing all the experiments will be made available at github \url{https://github.com/jlqzprj/matrix-cond-mixed-DNN}.

\begin{table}
\begin{center}
\begin{threeparttable} 
\caption{The algorithmic parameters used for the examples. The top block is for Neumann problems, and the bottom block for Dirichlet problems.}\label{tab:algpara}
{
\setlength{\tabcolsep}{2pt}
\begin{tabular}{c|cccccc}
\toprule
  $\backslash$No.& \ref{exam:neu2d1} & \ref{exam:neu2d2}& \ref{exam:neu2d3}& \ref{exam:neu2d4} & \ref{exam:neu3d1} &  \ref{exam:neu3d2}\\
\midrule
     $\gamma_\sigma$& 1(10) & 1(10) & 1(10) & 1 & 1 & 1(10) \\
     $\gamma_b$& 1(10) & 1(2.5) & 1(2.5) & 1(2.5) & 1/6(5/3) & 1/6(5/3) \\
     $\gamma_A$& 1e-5 & 1e-5 & 1e-5 & 1e-5 & 1e-5 & 1e-5\\
     $n_r$& 1e4 & 1.5e4 & 1.5e4 & 1.5e4 & 2e4 & 2.4e4 \\
     $n_b$& 1e3 & 4e3 & 4e3 & 4e3 & 4.8e3 & 6e3 \\
     lr& 2e-3 & 3e-3 & 3e-3 & 2e-3 & 5e-3 & 5e-3 \\
     dr& 0.7 & 0.8 & 0.8 & 0.8 & 0.8 & 0.8 \\
     step& 3e3 & 3e3 & 3e3 & 3e3 & 3e3 & 3e3 \\
     epoch& 8e4(1e4) & 8e4(1.44e4) & 8e4(6e3) & 8e4(2e4) & 8e4(3e4) & 8e4(3e4) \\
\midrule
     &\ref{exam:diri2d1}&\ref{exam:diri2d2}& \ref{exam:diri2d3}& \ref{exam:diri2d4}&\ref{exam:diri3d1}&\ref{exam:diri3d2}\\
     $\gamma_\sigma$& 1(10) & 1(10) & 1(10) & 1 & 1(10) & 1 \\
     $\gamma_b$& 1(10) & 1(2.5) & 1(2.5) & 1(1/4) & 1/6(5/3) & 1/6(5/3) \\
     $\gamma_A$& 1e-5 & 1e-5 & 1e-5 & 1e-5 & 1e-5 & 1e-5 \\
     $n_r$& 1e4 & 1.5e4 & 1.5e4 & 1.5e4 & 2e4 & 2.4e4 \\
     $n_b$& 1e3 & 4e3 & 4e3 & 4e3 & 4.8e3 & 6e3 \\
     lr& 2e-3 & 3e-3 & 3e-3 & 3e-3 & 5e-3 & 5e-3 \\
     dr& 0.7 & 0.8 & 0.8 & 0.8 & 0.8 & 0.8 \\
     step& 3e3 & 3e3 & 3e3 & 4e3 & 3e3 & 3e3 \\
     epoch& 8e4(1.5e3) & 8e4(6.6e3) & 8e4(2.3e3) & 8e4(9e3) & 8e4(3e4) & 8e4(1e4) \\
\bottomrule
\end{tabular}}
\end{threeparttable}
\end{center}
\end{table}

\subsection{The Neumann problem} First we illustrate the approach in the Neumann case. The first example is about recovering an anisotropic conductivity matrix with polynomial entries.
\begin{example}\label{exam:neu2d1}
The domain $\Omega = (0,1)^2$, $A^\dag= \begin{pmatrix}
    2+x_1^2+x_2^2&1+(x_1-\frac12)(x_2-\frac12)\\
    1+(x_1-\frac12)(x_2-\frac12)&2+(x_1-\frac12)^2+(x_2-\frac12)^2\\
\end{pmatrix},$ $u_1^\dag=x_1+x_2+\frac{1}{3}(x_1^3+x_2^3)$, $u_2^\dag=x_1-x_2+\frac{1}{3}(x_1^3-x_2^3)$, $u_3^\dag=-u_1^\dagger$ and $u_4^\dag=-u_2^\dagger$.
\end{example}

\begin{figure}[htb!]
\centering
\setlength{\tabcolsep}{0em}
\begin{tabular}{ccccc}
\includegraphics[width=0.199\textwidth]{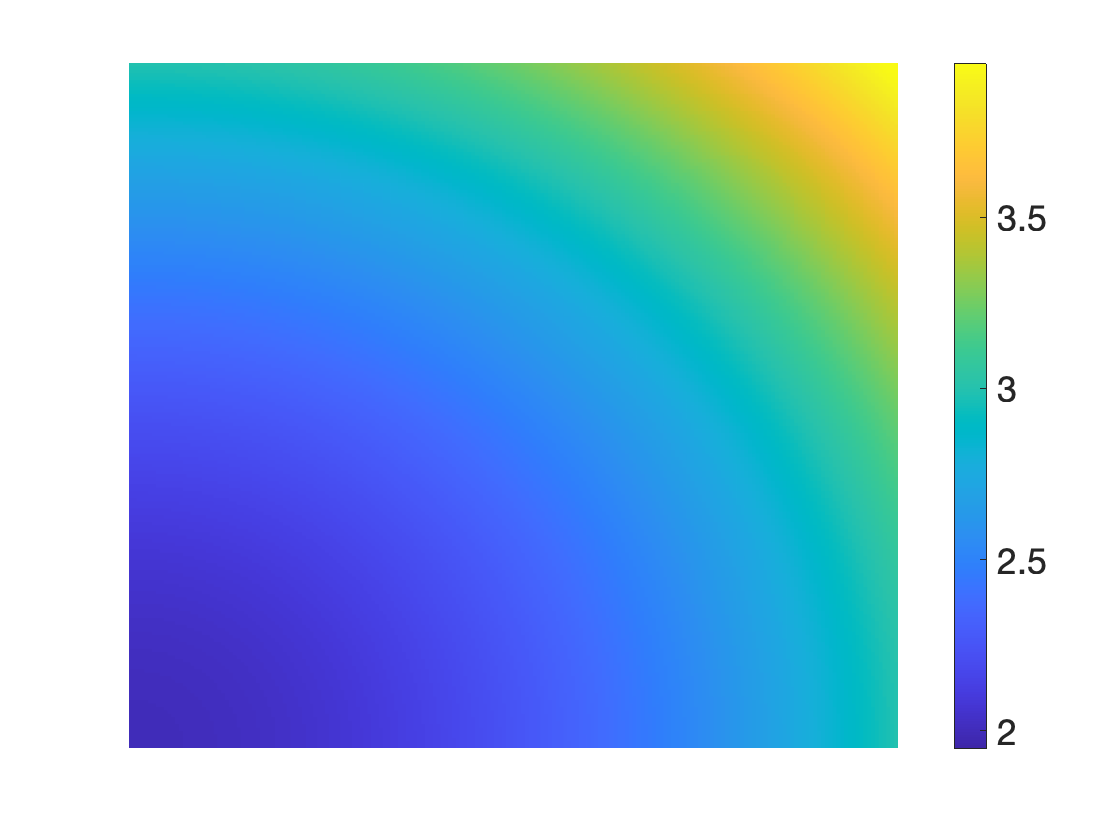} &
\includegraphics[width=0.199\textwidth]{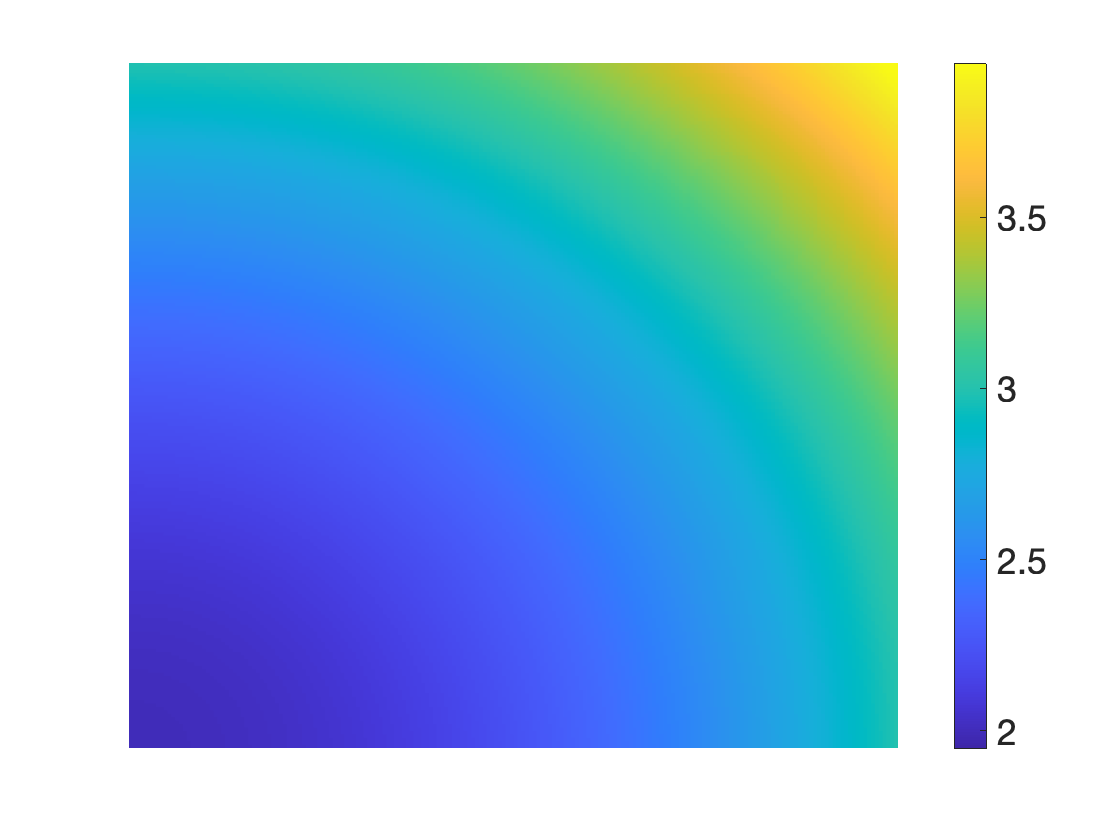} &
\includegraphics[width=0.199\textwidth]{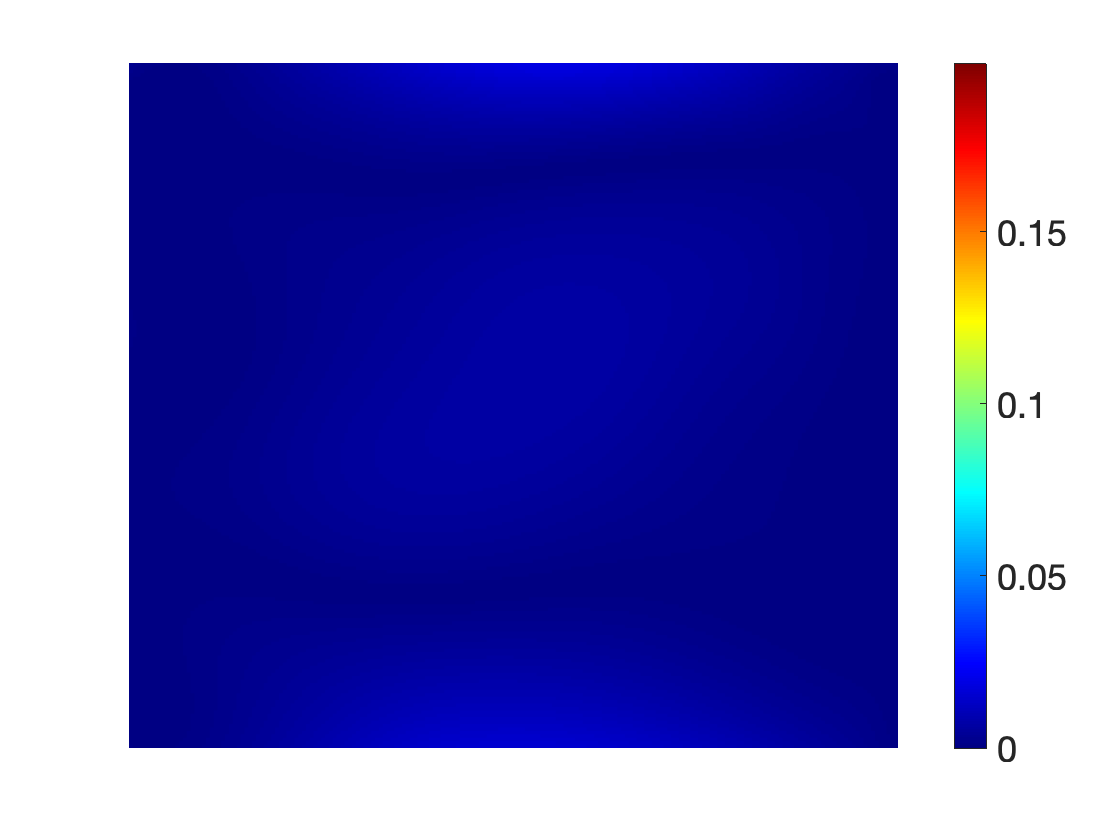} &
\includegraphics[width=0.199\textwidth]{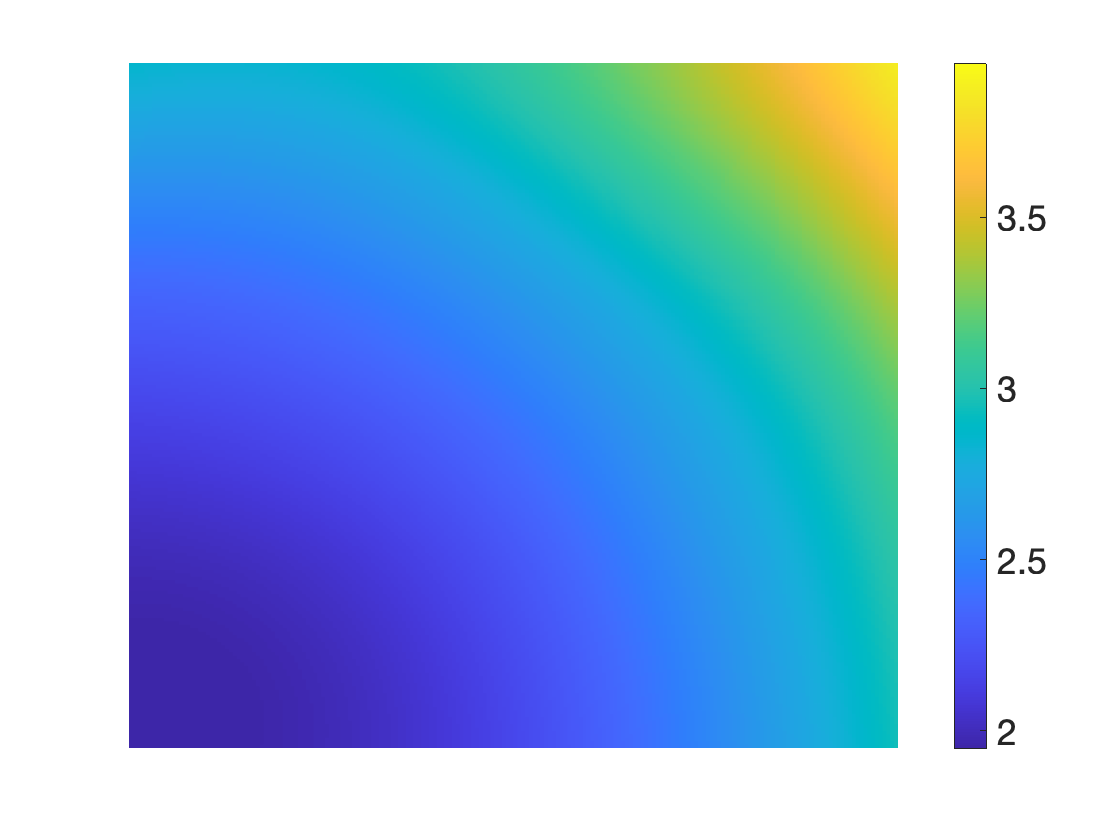} &
\includegraphics[width=0.199\textwidth]{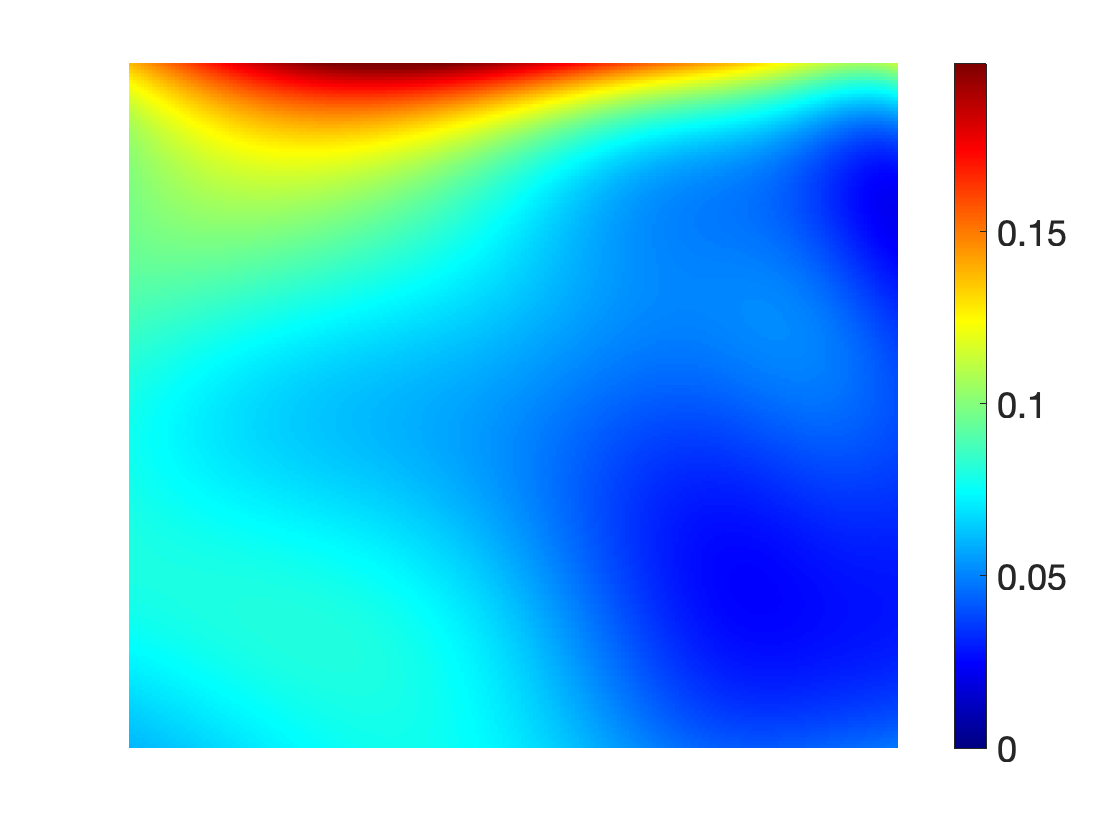} \\
\includegraphics[width=0.199\textwidth]{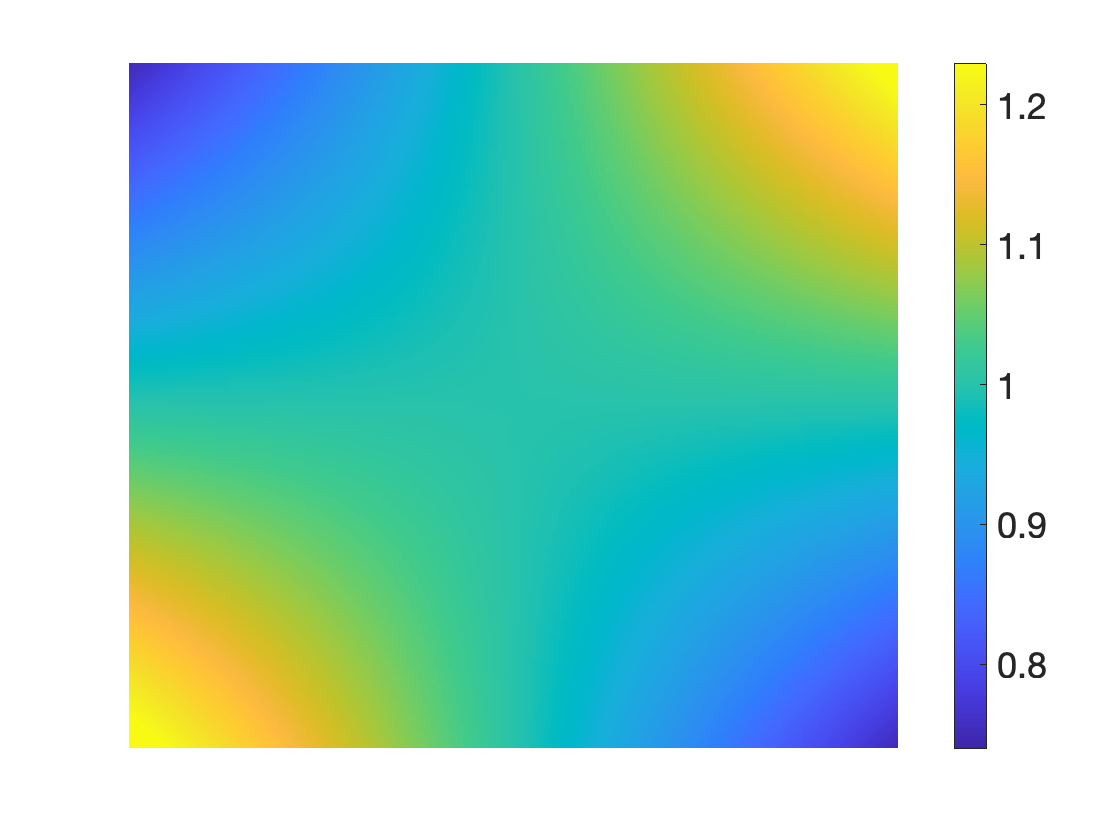} &
\includegraphics[width=0.199\textwidth]{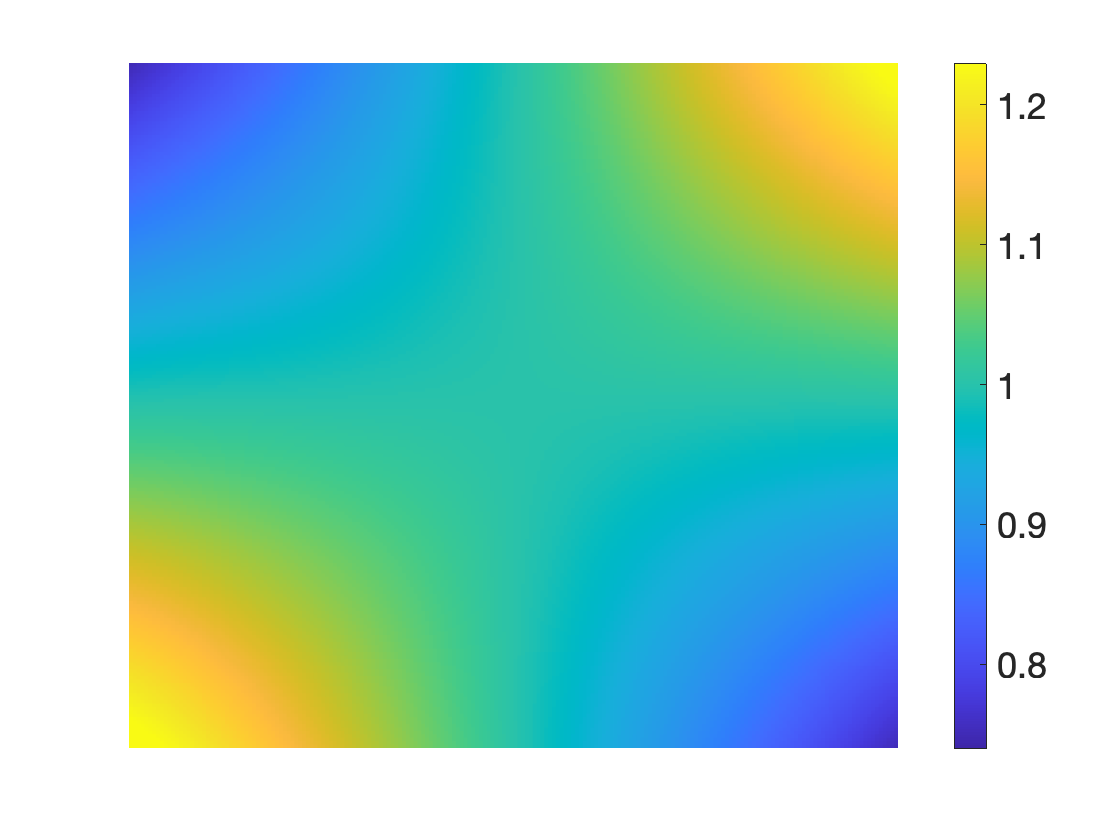} &
\includegraphics[width=0.199\textwidth]{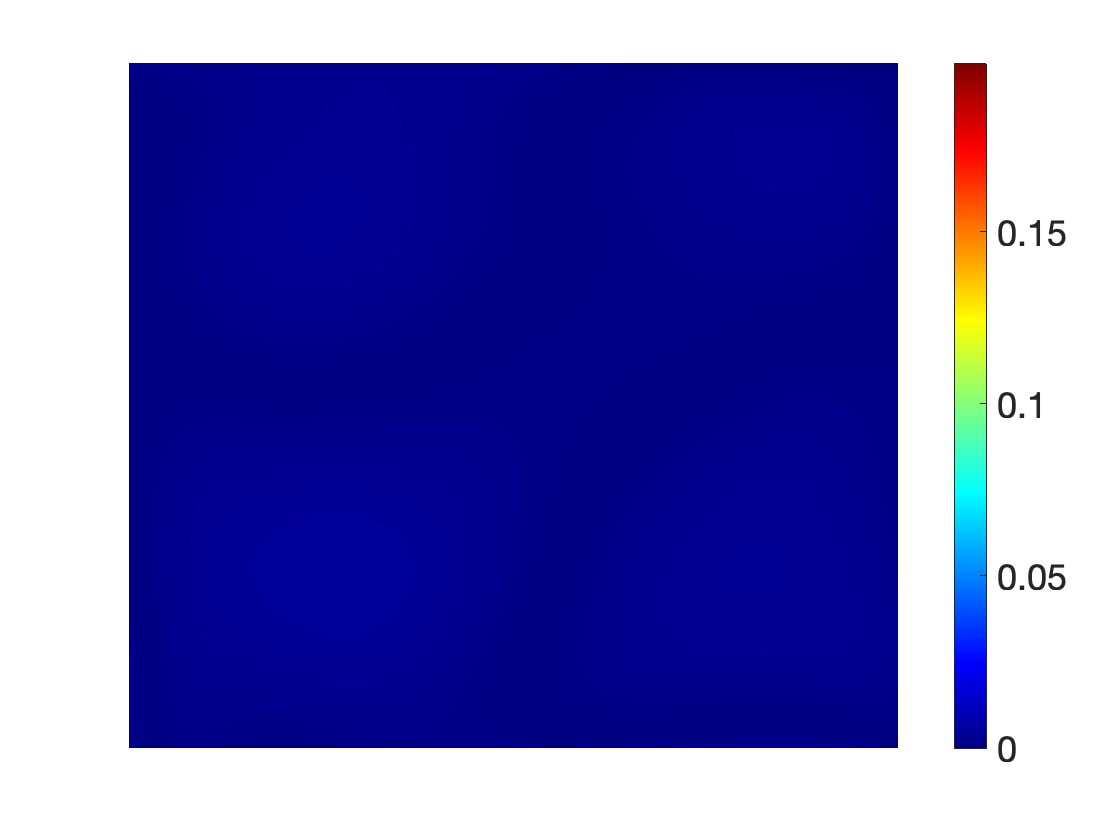} &
\includegraphics[width=0.199\textwidth]{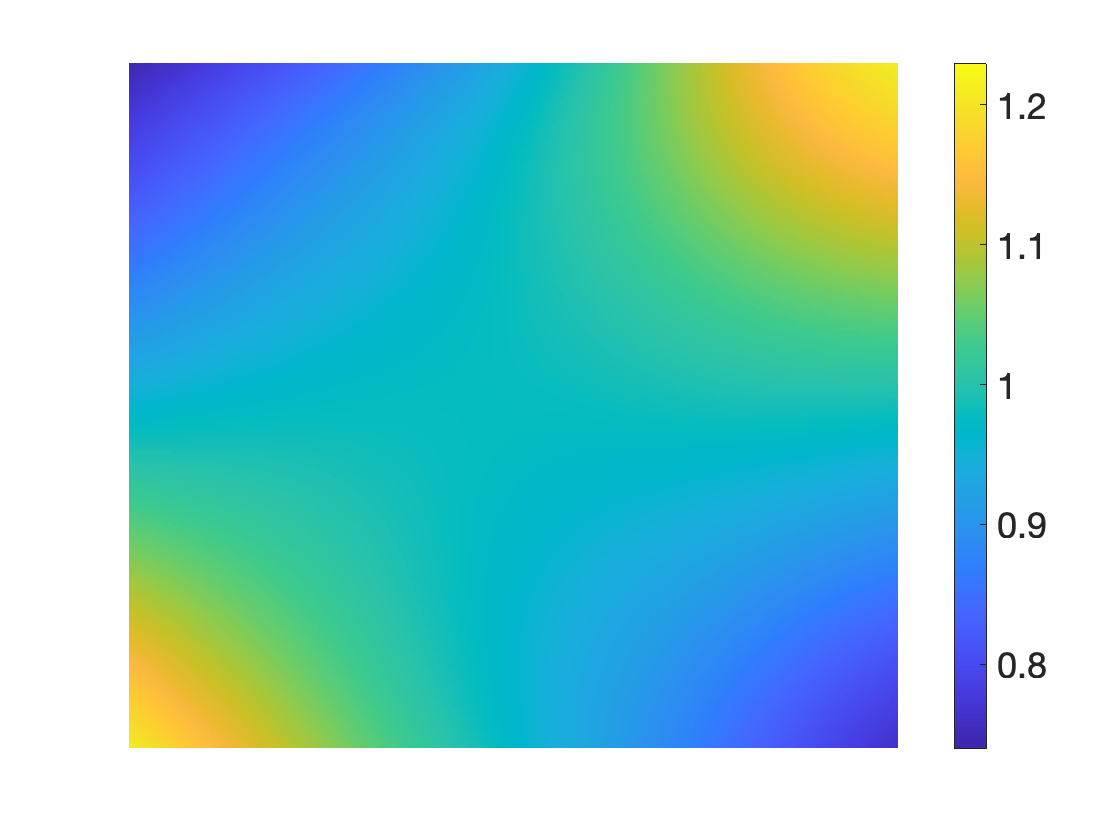} &
\includegraphics[width=0.199\textwidth]{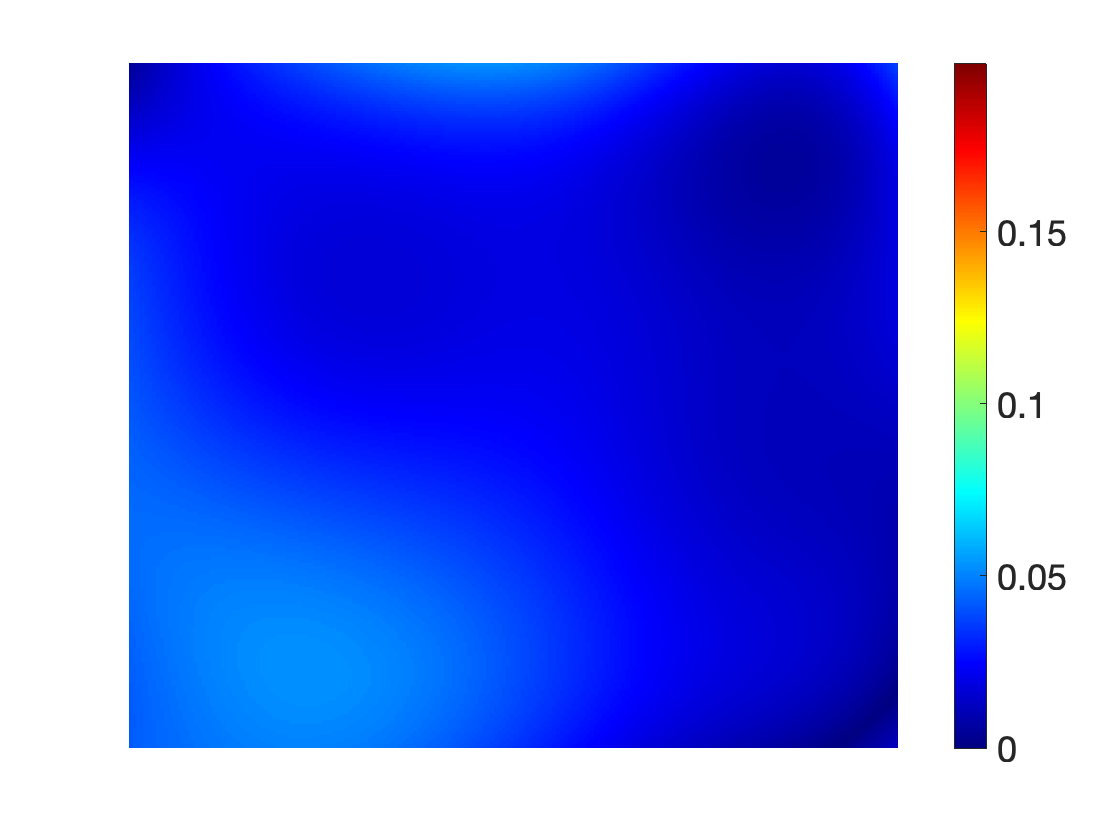} \\
\includegraphics[width=0.199\textwidth]{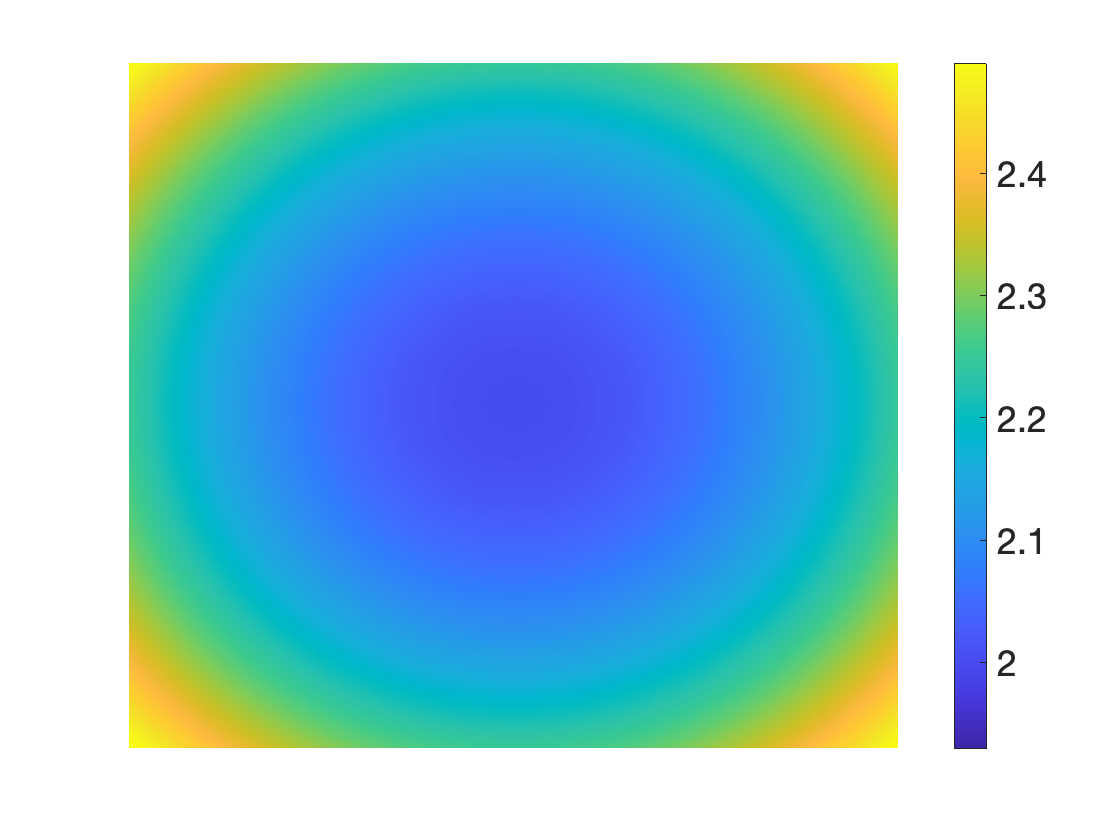} &
\includegraphics[width=0.199\textwidth]{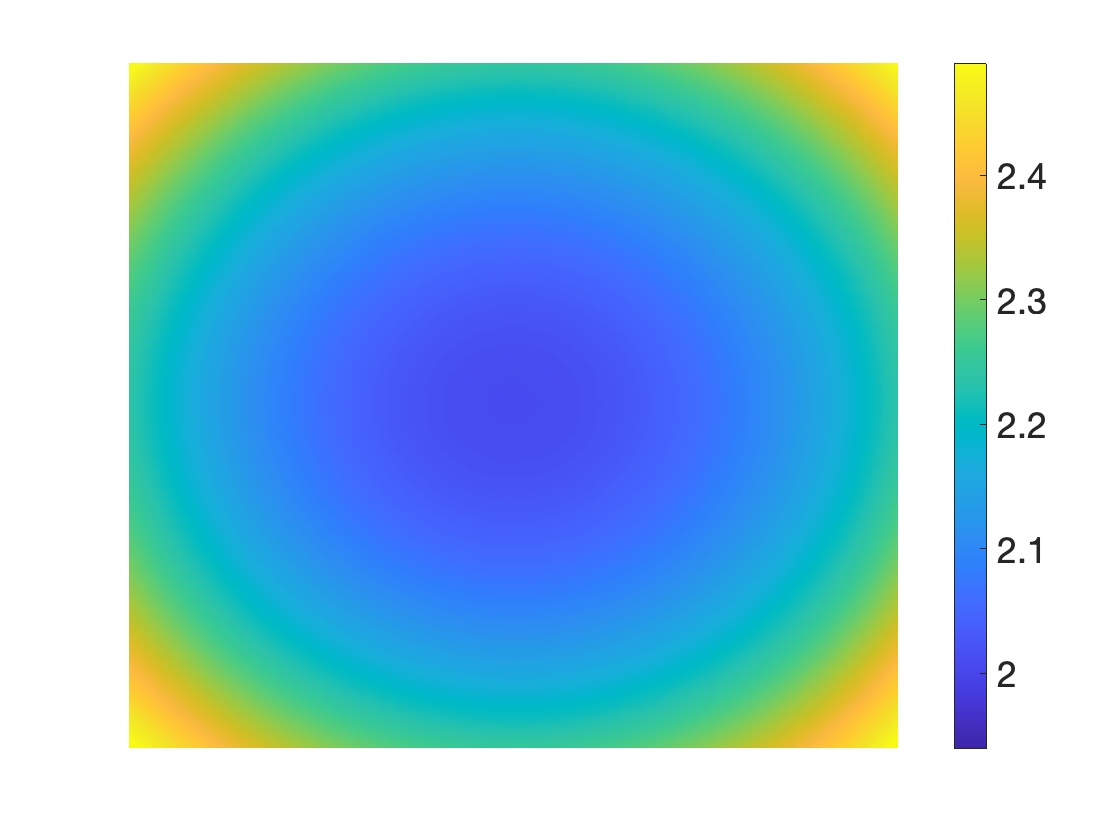} &
\includegraphics[width=0.199\textwidth]{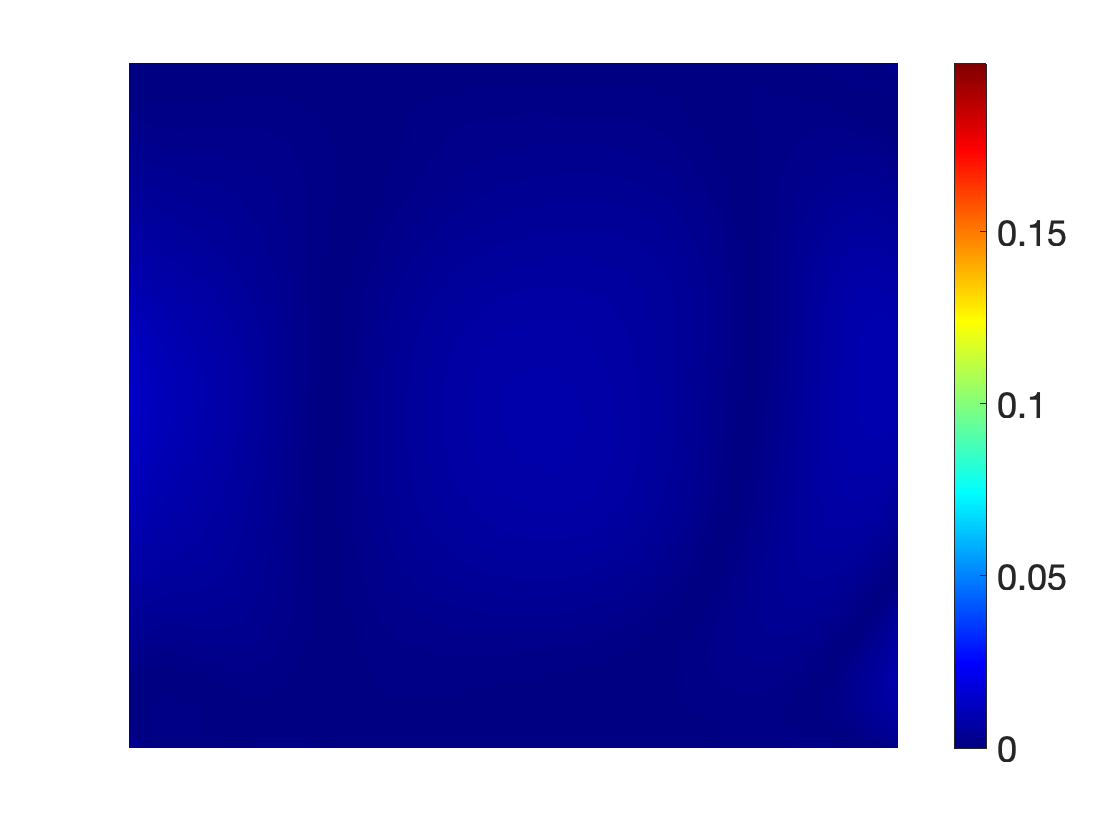} &
\includegraphics[width=0.199\textwidth]{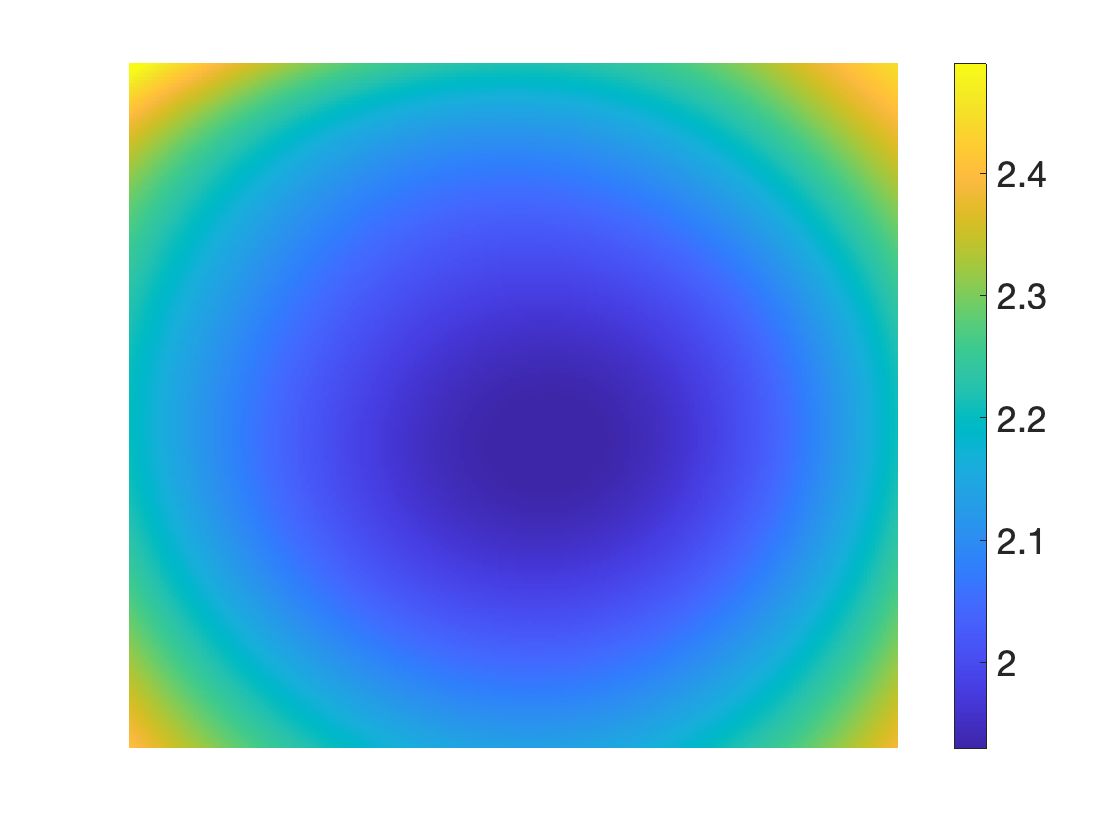} &
\includegraphics[width=0.199\textwidth]{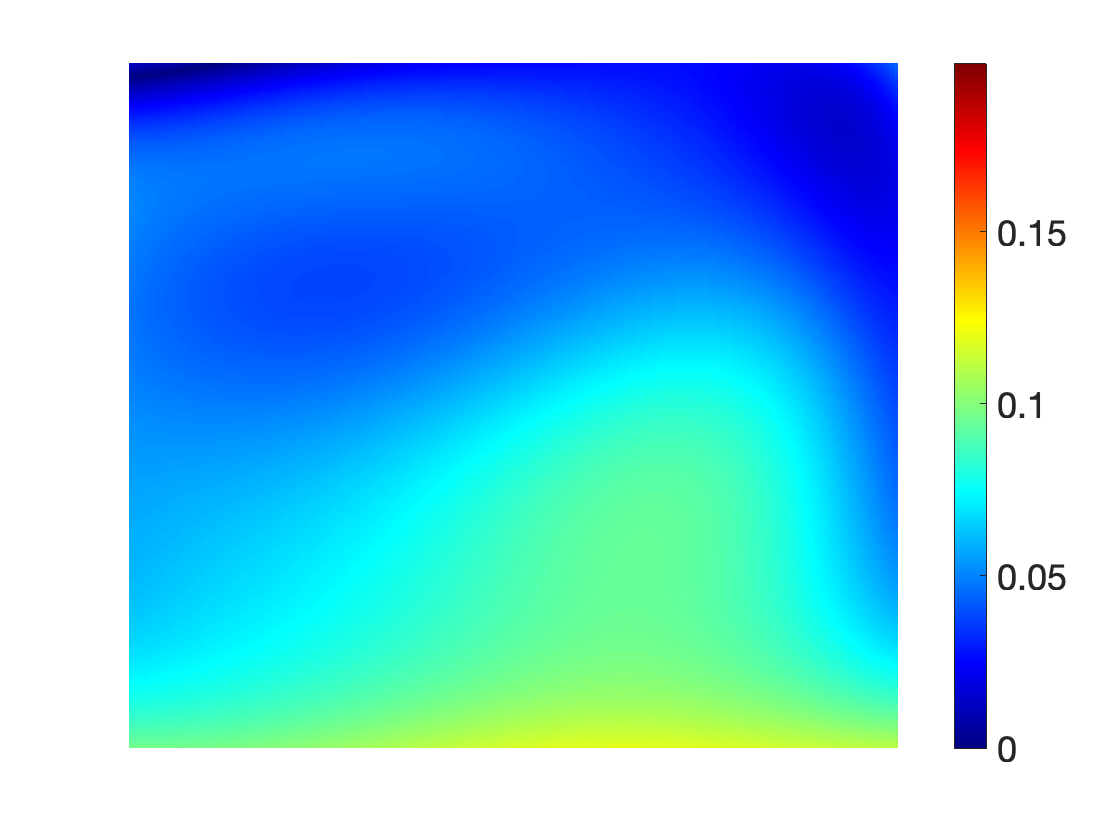} \\
(a) $A^\dag$  & (b) $\hat A$ & (c) $|\hat A-A^\dag|$ & (d) $\hat A$ & (e) $|\hat A-A^\dag|$
\end{tabular}
\caption{The reconstructions for Example \ref{exam:neu2d1} with exact data in (b) and noisy data $(\delta=10\%)$ in (d). From the top to bottom, the results are for $A_{11}$, $A_{12}$ and $A_{22}$, respectively.}
\label{fig:neu2d1}
\end{figure}
 
\begin{figure}[htb!]
\centering
\setlength{\tabcolsep}{0em}
\begin{tabular}{ccc}
\includegraphics[width=0.32\textwidth]{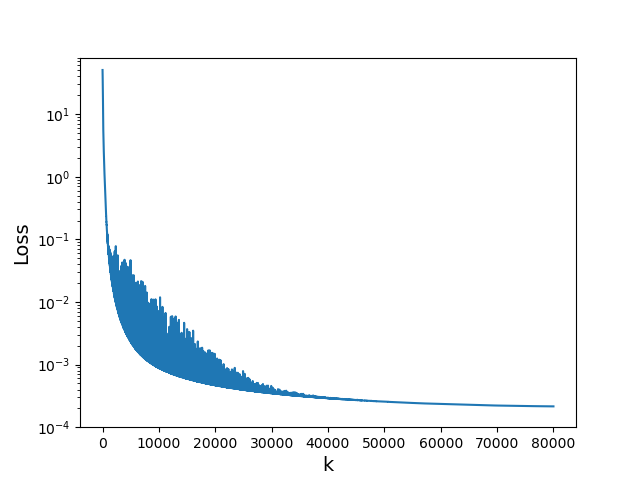} &
\includegraphics[width=0.32\textwidth]{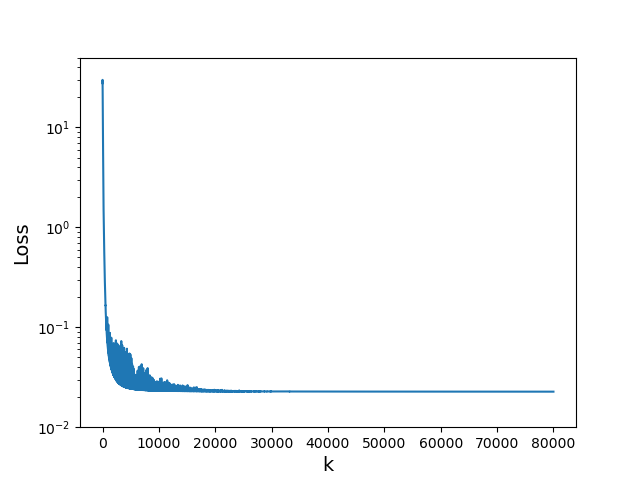} &
\includegraphics[width=0.32\textwidth]{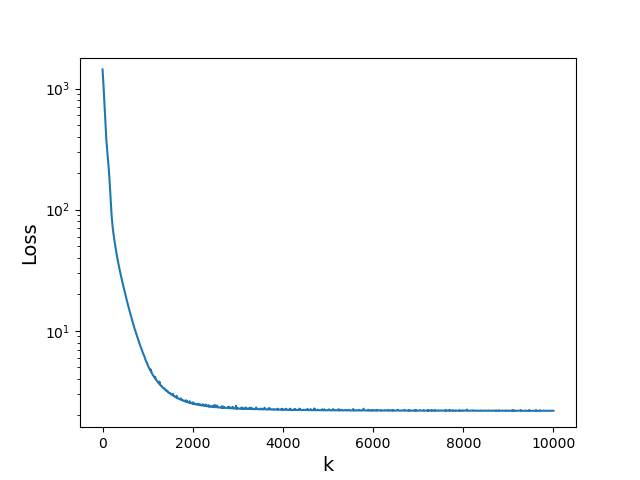}\\
\includegraphics[width=0.32\textwidth]{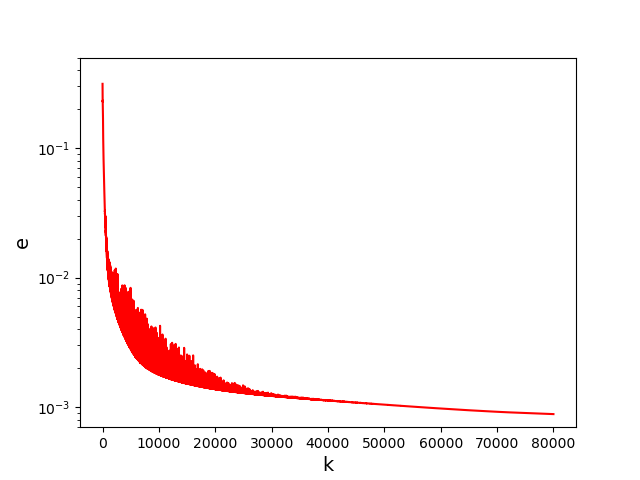} &
\includegraphics[width=0.32\textwidth]{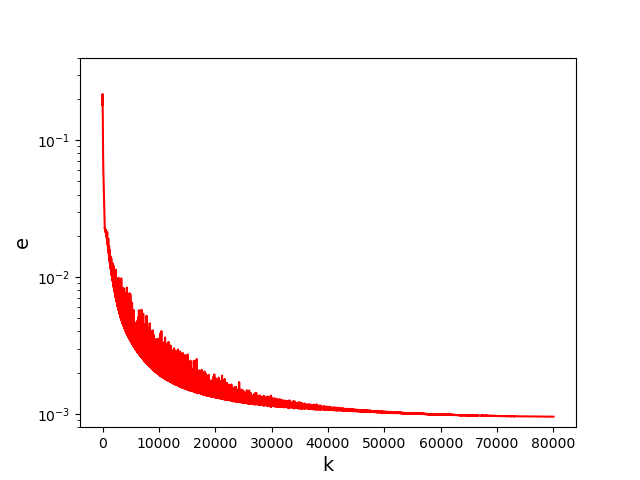} &
\includegraphics[width=0.32\textwidth]{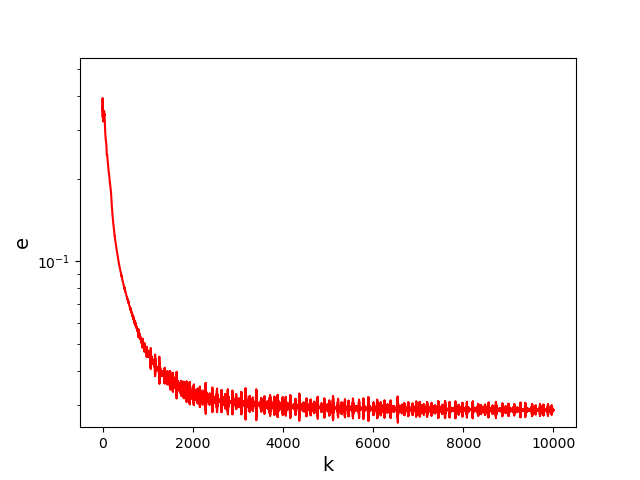}\\
(a) $\delta=0\%$  & (b) $\delta=1\%$ & (c) $\delta=10\%$
\end{tabular}
\caption{The evolution of the loss $($top$)$ and the reconstruction error $e(\hat A)$ $($bottom$)$ during the training process for Example \ref{exam:neu2d1} at different noise levels.}
\label{fig:neu1losse}
\end{figure}

Fig. \ref{fig:neu2d1} shows the DNN reconstructions and plots of the pointwise errors $|\hat A-A^\dag|$ for exact and noisy data. The features of each entry of the anisotropic conductivity tensor are well resolved, even for up to $10\%$ noise, indicating the robustness of the DNN approach with respect to noise. This robustness is further supported by the convergence behavior of the loss and error $e(\hat A)$ at different noise levels, cf. Fig. \ref{fig:neu1losse}. The DNN approach seems fairly stable with respect to the iteration index, and the final errors are close to each other for all noise levels. Note that early stopping must be adopted for reconstructions with noisy data, as neural networks are known to exhibit spectral bias \cite{Ulyanov:2018,dolean2024multilevel}, where low-frequency features are learned much faster than high-frequency features. Early stopping prevents the neural networks from overfitting to the data noise. However, a practical yet provable criterion for choosing the stopping index is still missing, so this index has been determined using a trial-and-error approach.

There are several algorithmic parameters influencing the overall quality of the DNN reconstruction $\hat A$, e.g., the numbers of sampling points ($n_r$ and $n_b$), and DNN architectural parameters {(width, depth, and activation function)}. A practical guideline for choosing these parameters is still missing, and we explore the issue empirically.
In Table \ref{neut1}, we study the impact of the regularization parameter $\gamma_A$ on the reconstruction error $e(\hat{A})$. The DNN method demonstrates remarkable robustness in the presence of data noise, maintaining fairly accurate results even for noise level $\delta$ up to 10\%. Additionally, the results show little sensitivity with respect to the variation of $\gamma_A$, indicating that this parameter plays a less critical role than suggested by the theoretical analysis. This is attributed to the inductive bias of the DNN representation, which favors smooth functions, therefore providing a built-in regularizing effect that is not accounted for in the theoretical analysis. Note also that there is an inherent accuracy limitation of the DNN approach, i.e., the reconstruction cannot be made
arbitrarily accurate for exact data $\nabla u^\dag$. This may be attributed to the optimization error: the nonconvexity of the loss landscape can cause the optimizer to fail in finding a global minimizer of the empirical loss, settling instead for an approximate local minimizer. This phenomenon has been observed across a broad range of neural solvers based on DNNs
\cite{RaissiPerdikarisKarniadakis:2019,EYu:2018,JinLiLu:2022,Jin:2024conductivity,krishnapriyan2021characterizing}. Tables \ref{neut2}-\ref{neut3}
show that the reconstruction error  does not vary much
with different DNN architectures and numbers of sampling points. This also
agrees well with the convergence behavior of ADAM in Fig.
\ref{fig:neu1losse}: both loss $\widehat J_{\bsgamma}$ and error $e(\hat A)$ stagnate
at a certain level.

In the experiments, we have adopted the off-shelf optimizer ADAM for the training of the loss, which has been widely employed in deep learning. Empirically, it takes thousands of iterations to reach convergence, and the efficiency still deserves further improvement. Even worse, the optimizer may get trapped in a local minimum, leading to pronounced optimization errors. This fact greatly complicates the numerical verification of error estimates, which
so far has not been numerically realized for neural PDE solvers and represents one outstanding challenge in scientific machine learning \cite{SiegelXu:2023}.
It is an important theoretical question to analyze the optimization error and its impact on the reconstruction. The current approaches include neural tangent kernel \cite{Jacot:2018} and mean field analysis \cite{Mei:2018}, but it remains to extend these approaches to the context of PDE inverse problems.

\begin{table}[htb!]
  \centering
  \caption{The variation of the reonstruction error $e(\hat{A})$ with respect to various algorithmic parameters. }
\begin{threeparttable}
\subfigure[$e$ v.s. $\gamma_A$ and $\delta$\label{neut1}]{\begin{tabular}{c|ccc}
\toprule
  $\gamma_A\backslash\delta$&   0\% &  1\%& 10\%\\
\midrule
     1.00e-2 &  3.07e-3 & 3.52e-3 & 3.12e-2 \\
     1.00e-3 &  1.00e-3 & 1.32e-3 & 2.98e-2\\
     1.00e-4 &  9.56e-4 & 1.19e-3 & 2.98e-2\\
     1.00e-5 &  8.84e-4 & 9.54e-4 & 2.87e-2\\
\bottomrule
\end{tabular}} \\
\subfigure[$e$ v.s. $W_A$ and  $L_A$\label{neut2}]{
\begin{tabular}{c|ccc}
\toprule
${W_A}\backslash L_A$&   5 &  10& 15\\
\midrule
     8 &1.18e-3&9.49e-4&8.61e-4 \\
     16&6.95e-4&1.38e-3&1.58e-3\\
     24&1.28e-3&1.86e-3&1.66e-3\\
     32&8.84e-4&1.88e-3&1.62e-3\\
\bottomrule
\end{tabular}}\quad
\subfigure[$e$ v.s. $n_b$ and $n_r$\label{neut3}]{
\begin{tabular}{c|ccc}
\toprule
$n_b\backslash n_r$&   5000 &  10000& 15000\\
\midrule
     500&  1.30e-3&1.40e-3&1.05e-3\\
     1000& 1.15e-3&8.84e-4&1.77e-3\\
     2000& 1.52e-3&1.81e-3&1.03e-3\\
     4000& 1.94e-3&1.39e-3&1.67e-3\\

\bottomrule
\end{tabular}}
\end{threeparttable}
\end{table}

The second example is about recovering an anisotropic conductivity matrix with mixed oscillatory and polynomial entries.
\begin{example}
The domain $\Omega = (0,1)^2$, $A^\dag= \begin{pmatrix}
    2+\frac{\sin(4\pi x_2)}{2}&1+\frac{\sin^2(2\pi x_2)}{2}\\
    1+\frac{\sin^2(2\pi x_2)}{2}&2\\
\end{pmatrix},$ $u_1^\dag=x_1+x_2+\frac{1}{3}(x_1^3+x_2^3)$, $u_2^\dag=x_1-x_2+\frac{1}{3}(x_1^3-x_2^3)$, $u_3^\dag=-u_1^\dagger$ and $u_4^\dag=-u_2^\dagger$.
\label{exam:neu2d2}
\end{example}

\begin{figure}[htb!]
\centering
\setlength{\tabcolsep}{0em}
\begin{tabular}{ccccc}
\includegraphics[width=0.199\textwidth]{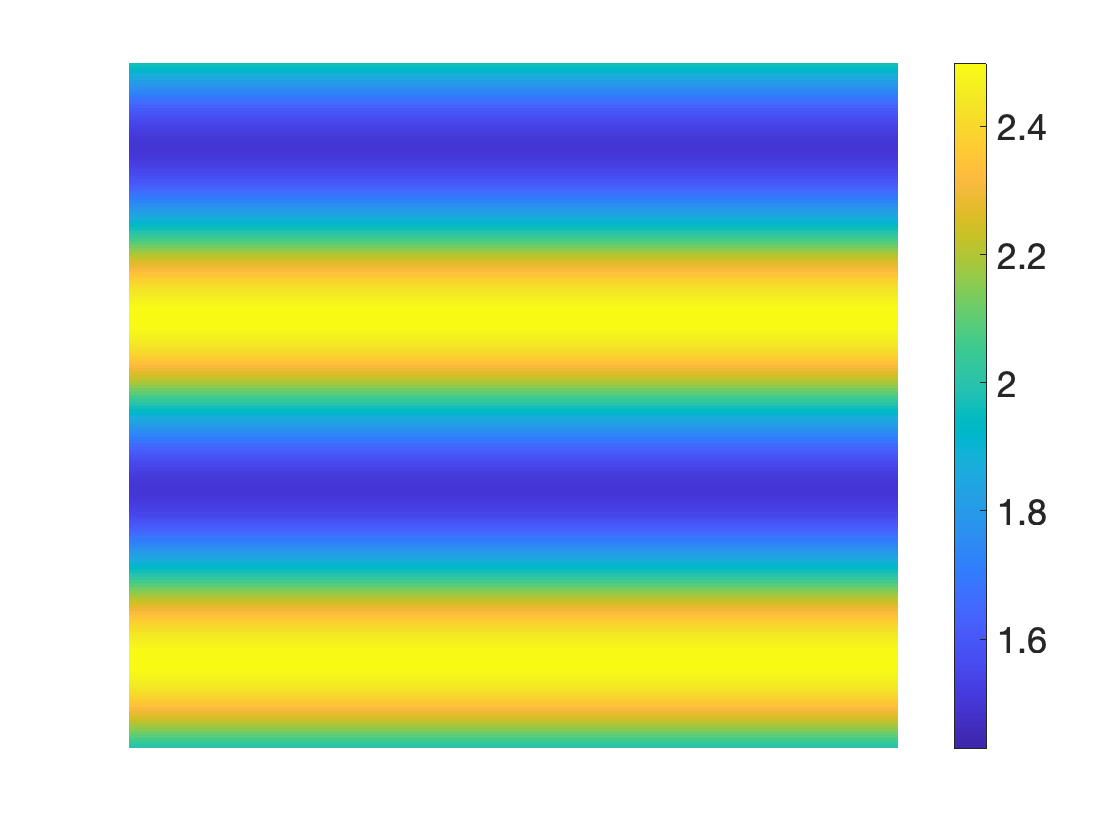} &
\includegraphics[width=0.199\textwidth]{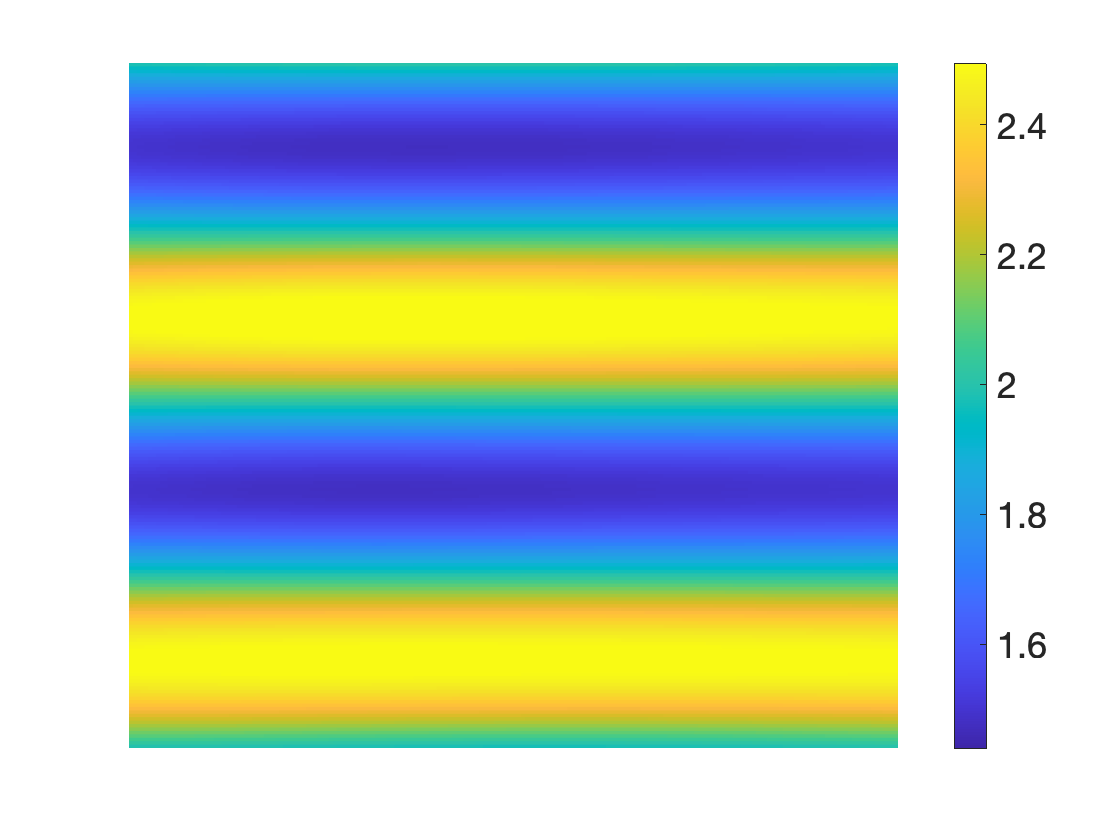} &
\includegraphics[width=0.199\textwidth]{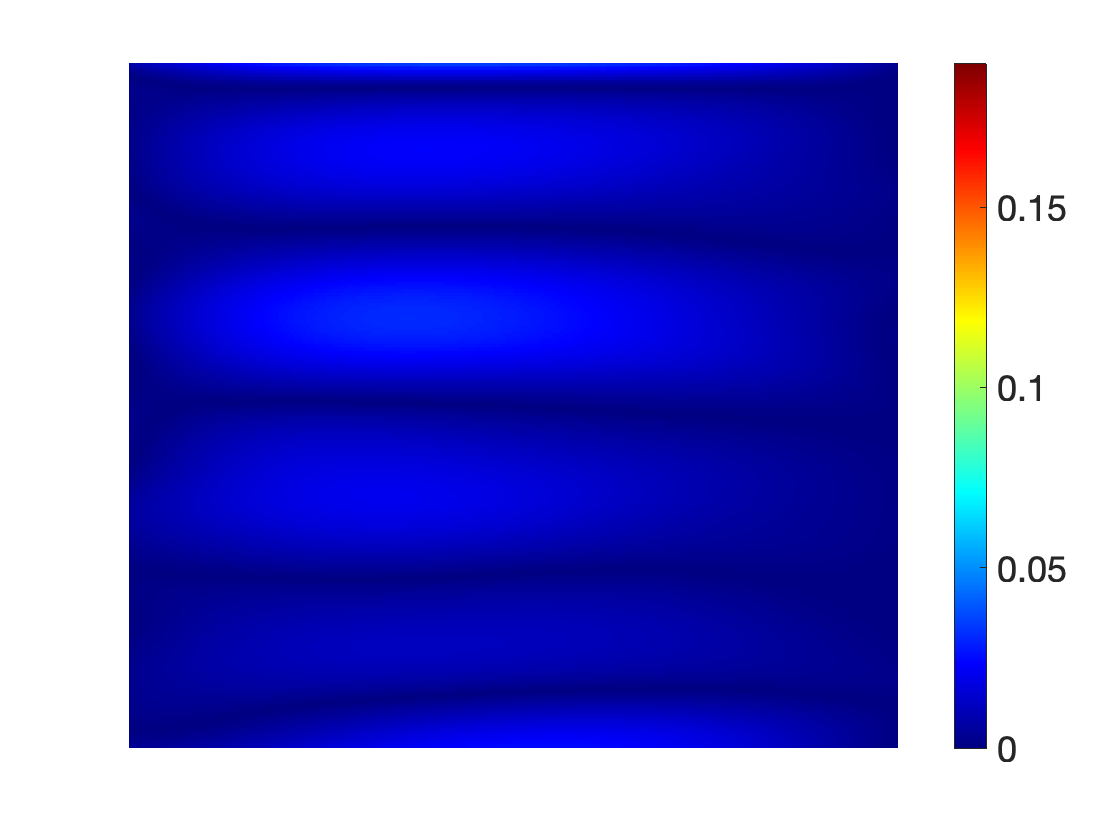} &
\includegraphics[width=0.199\textwidth]{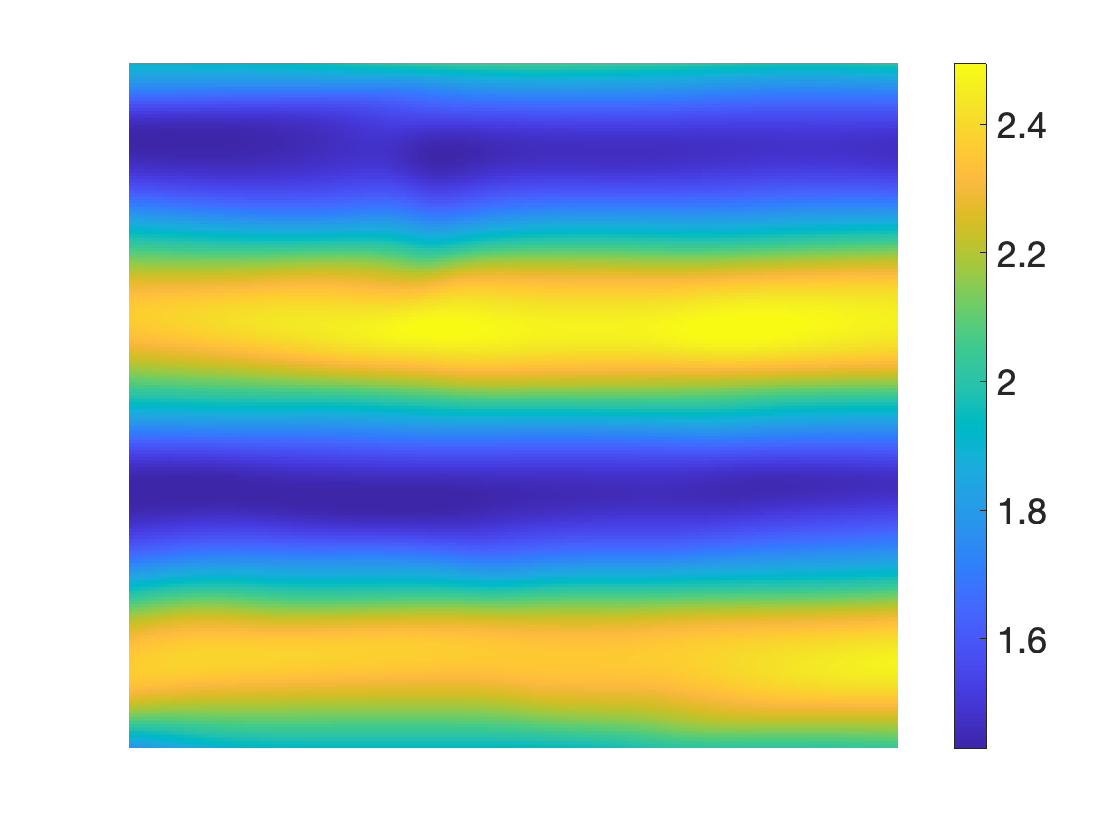} &
\includegraphics[width=0.199\textwidth]{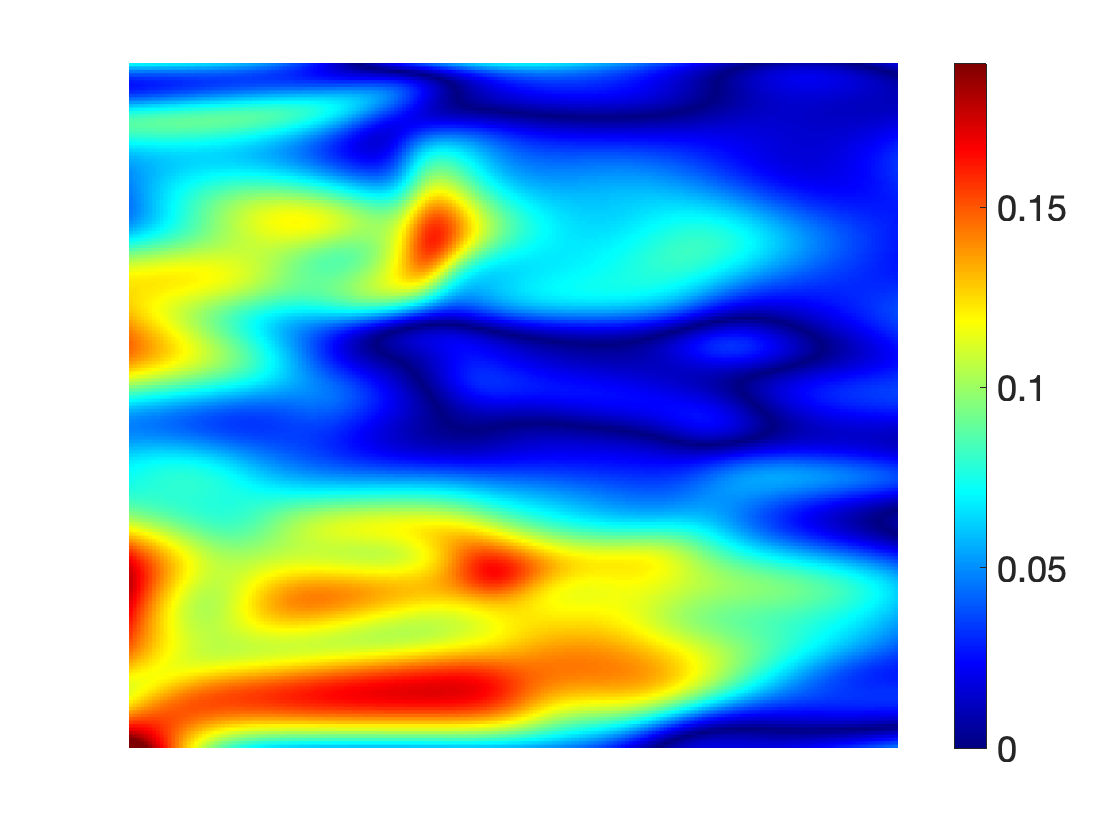} \\
\includegraphics[width=0.199\textwidth]{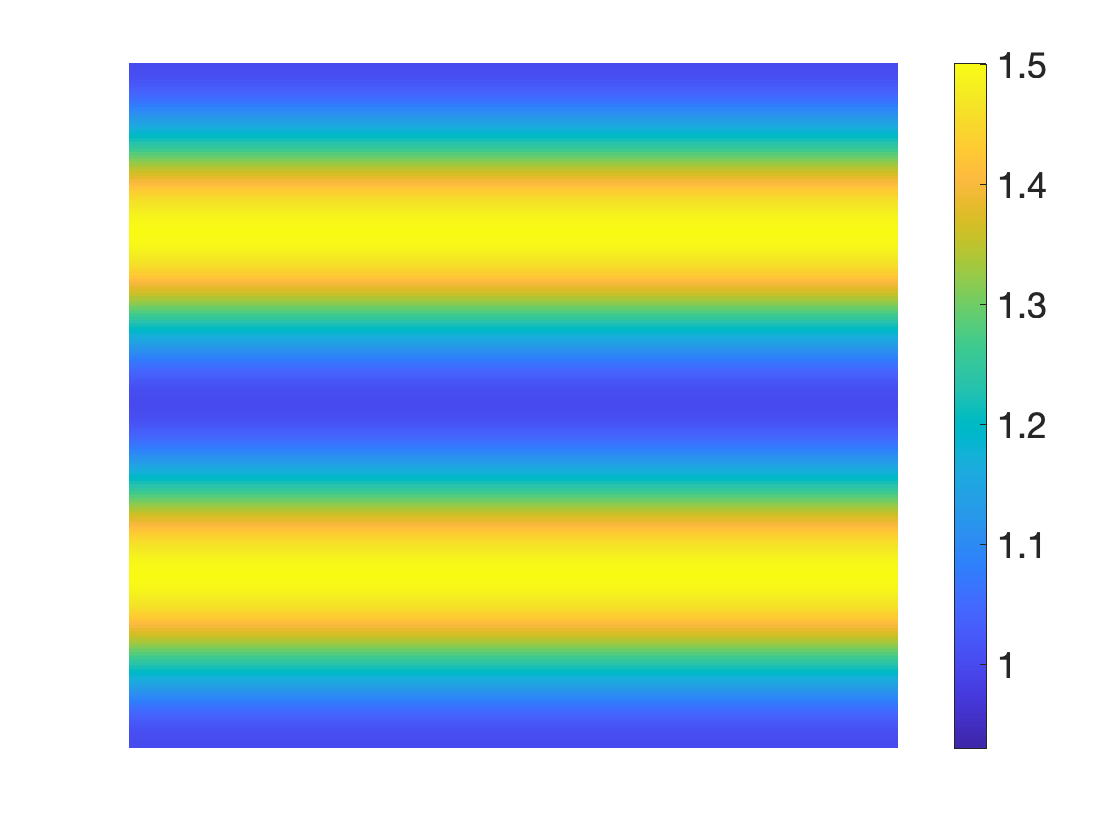} &
\includegraphics[width=0.199\textwidth]{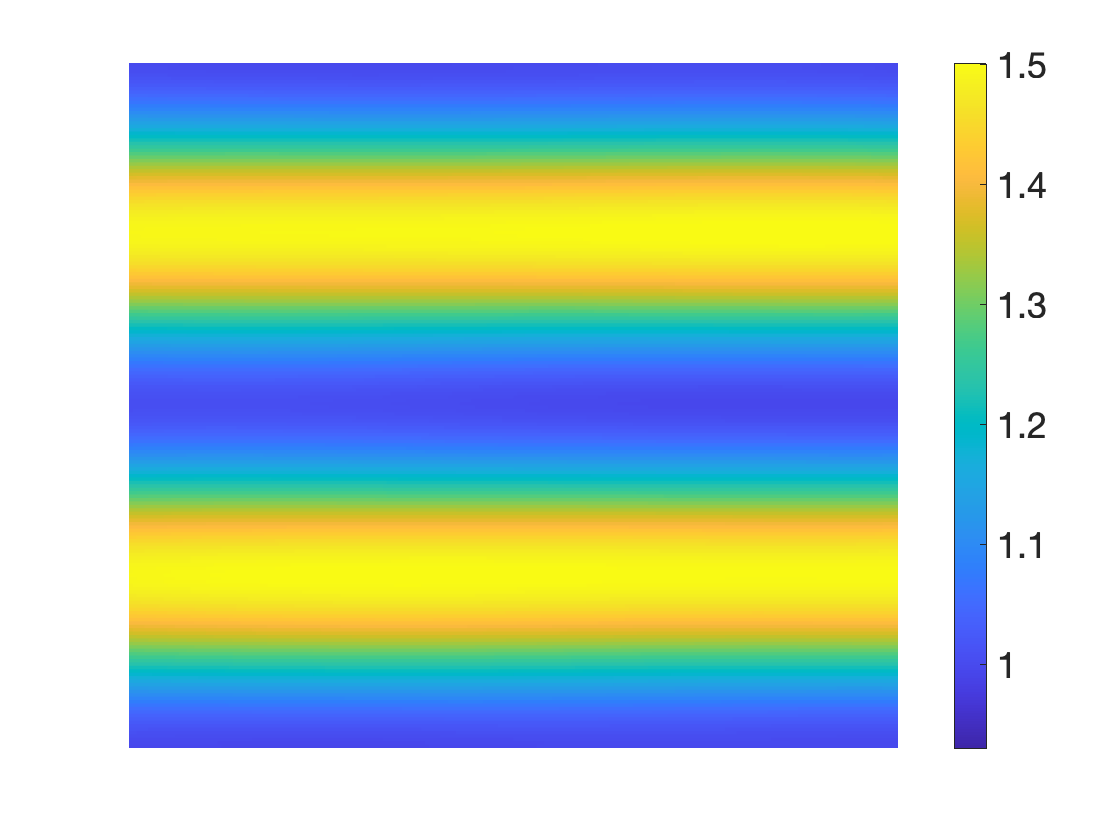} &
\includegraphics[width=0.199\textwidth]{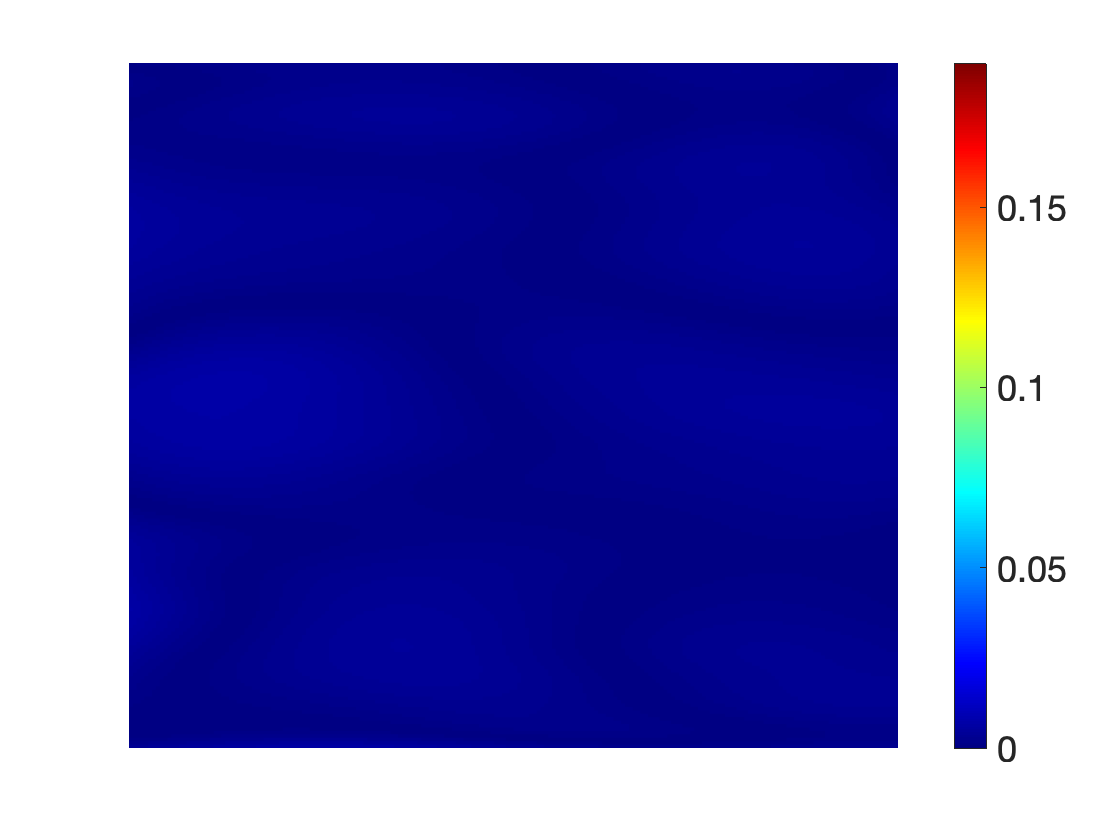} &
\includegraphics[width=0.199\textwidth]{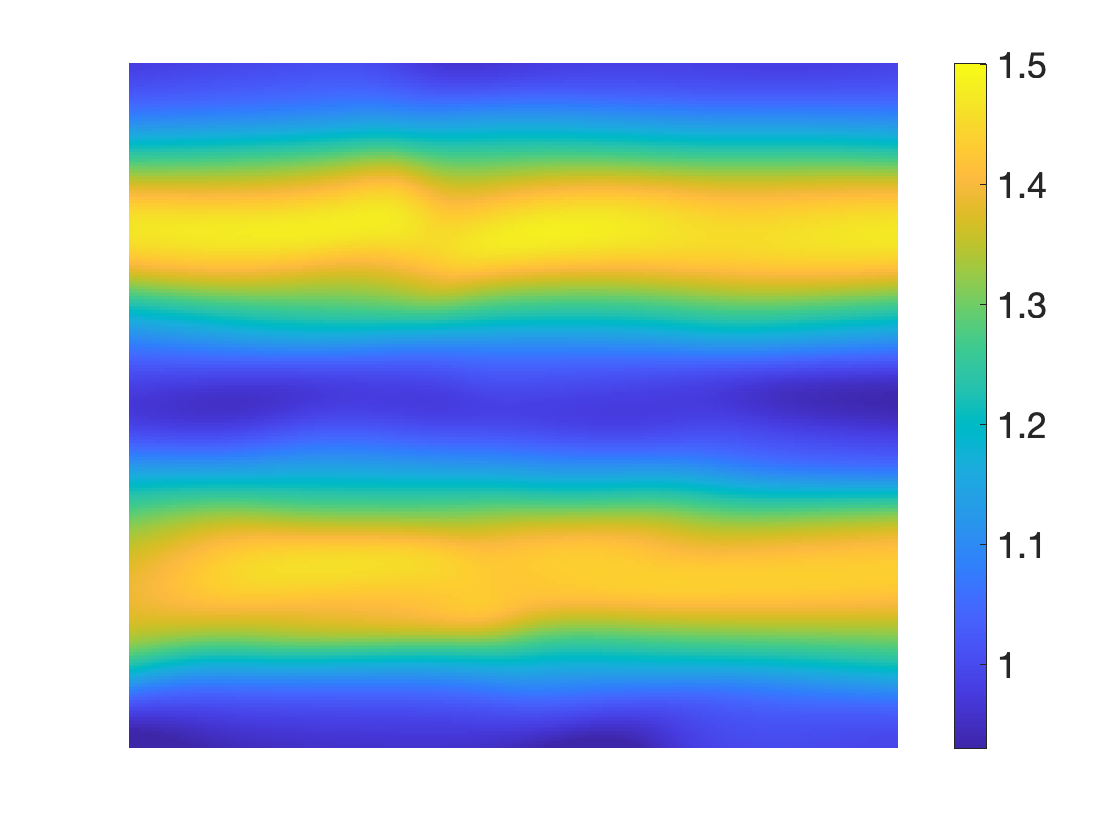} &
\includegraphics[width=0.199\textwidth]{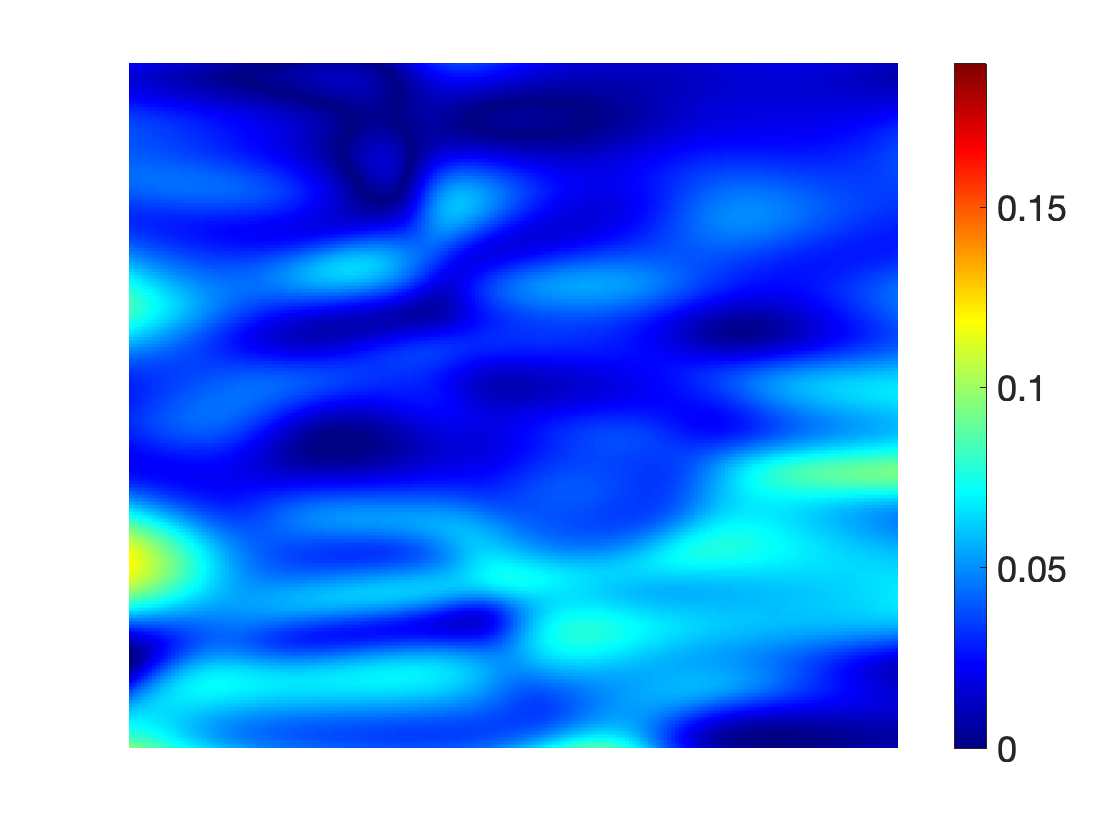} \\
\includegraphics[width=0.199\textwidth]{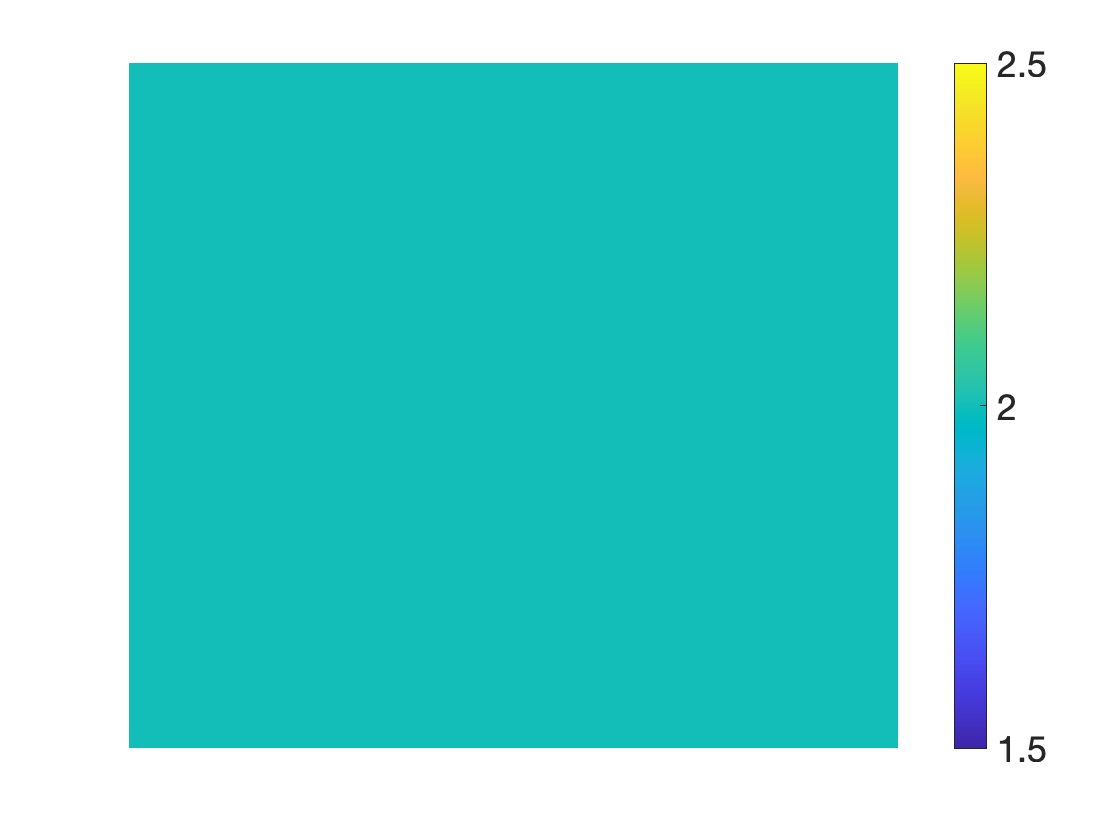} &
\includegraphics[width=0.199\textwidth]{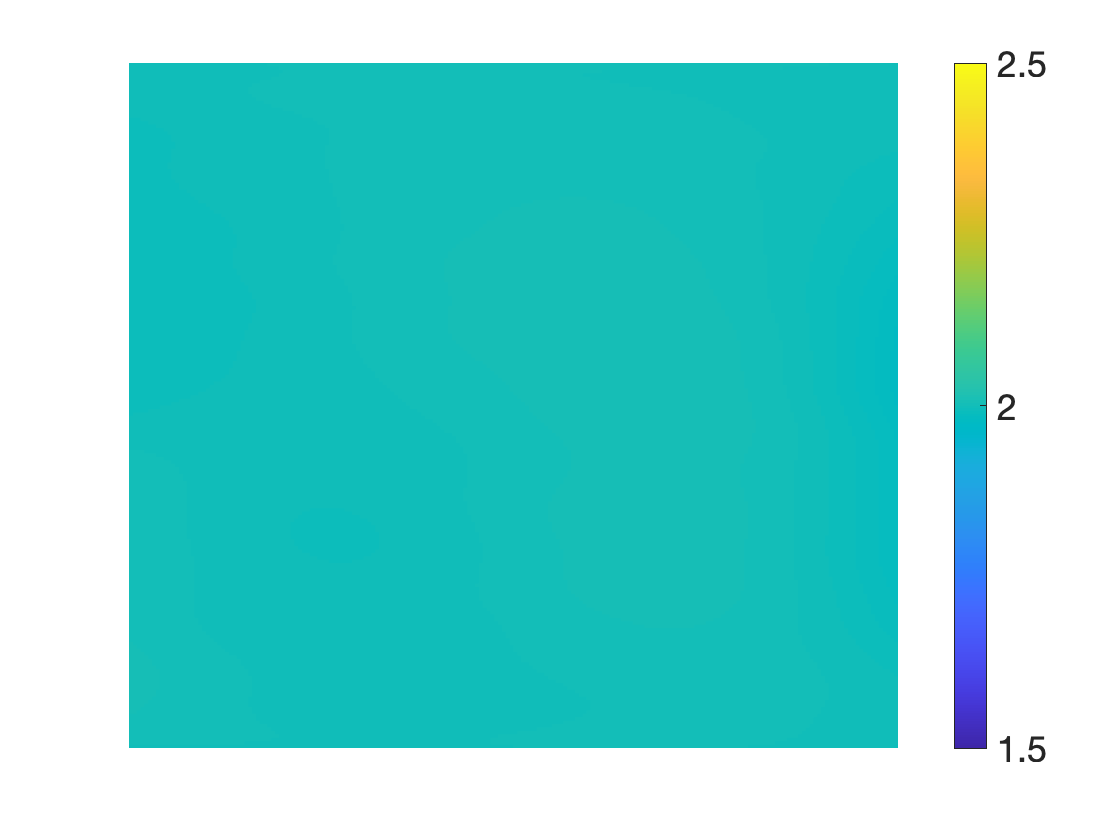} &
\includegraphics[width=0.199\textwidth]{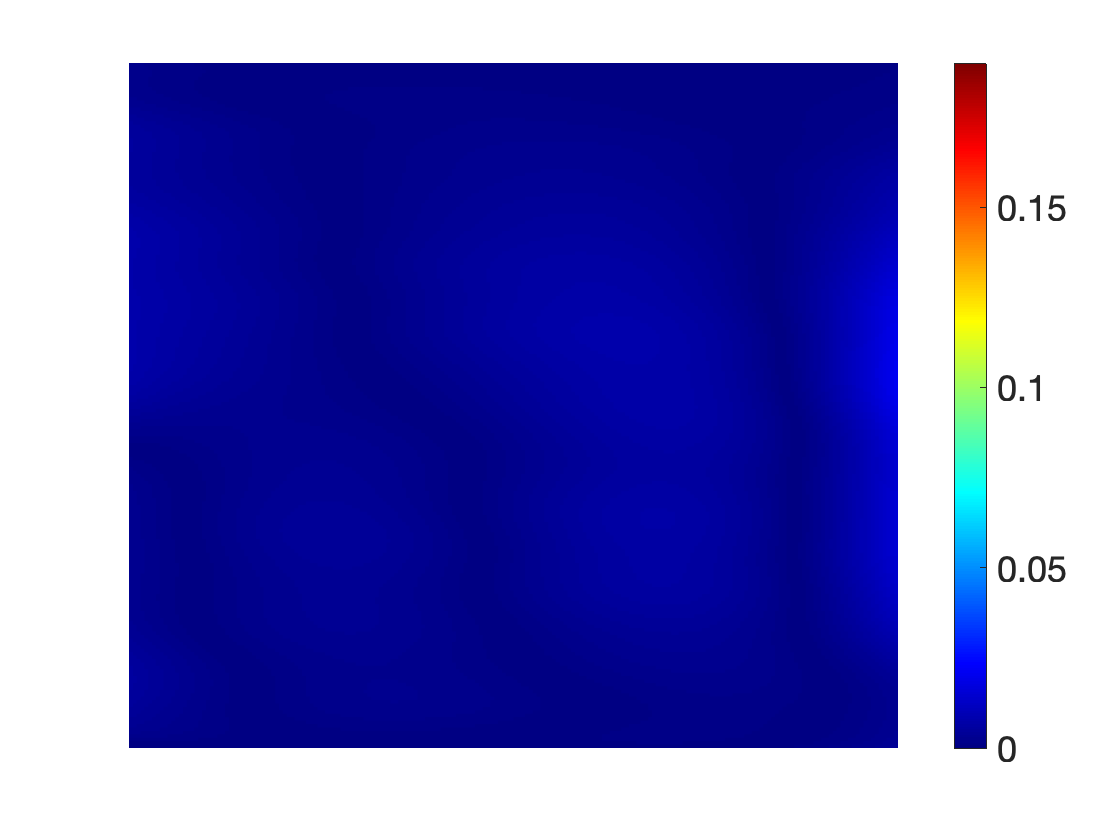} &
\includegraphics[width=0.199\textwidth]{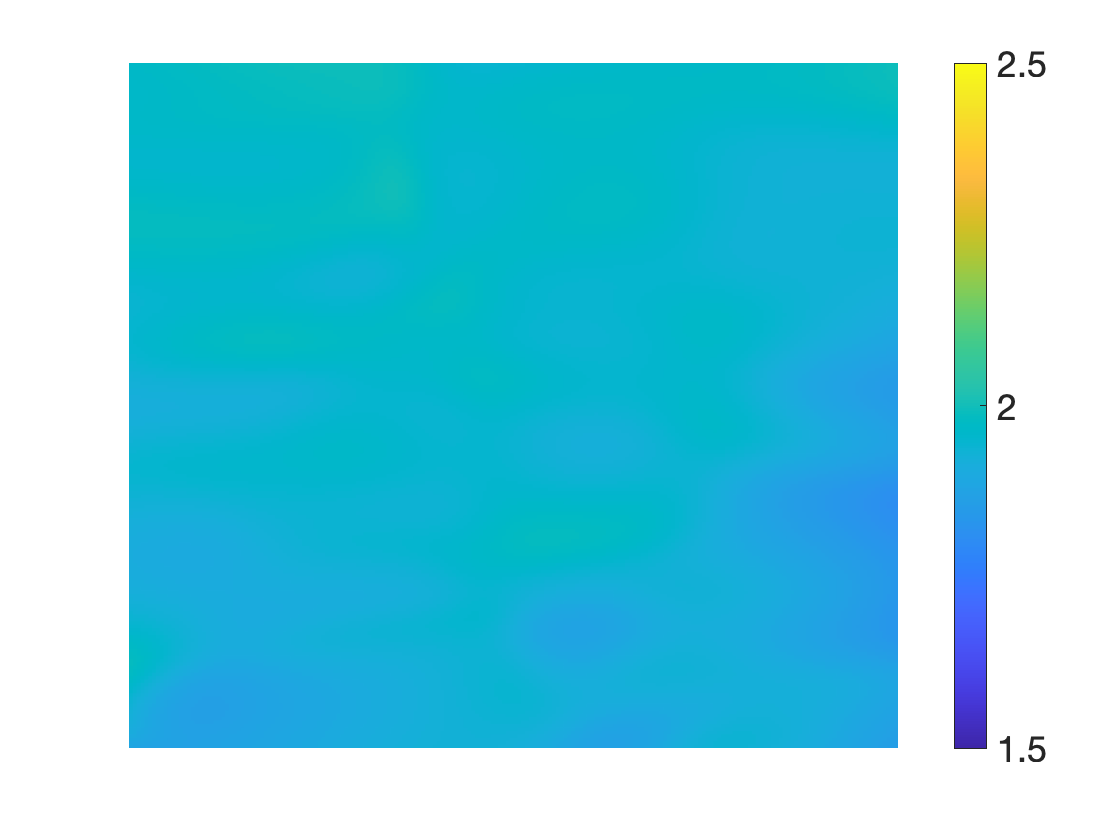} &
\includegraphics[width=0.199\textwidth]{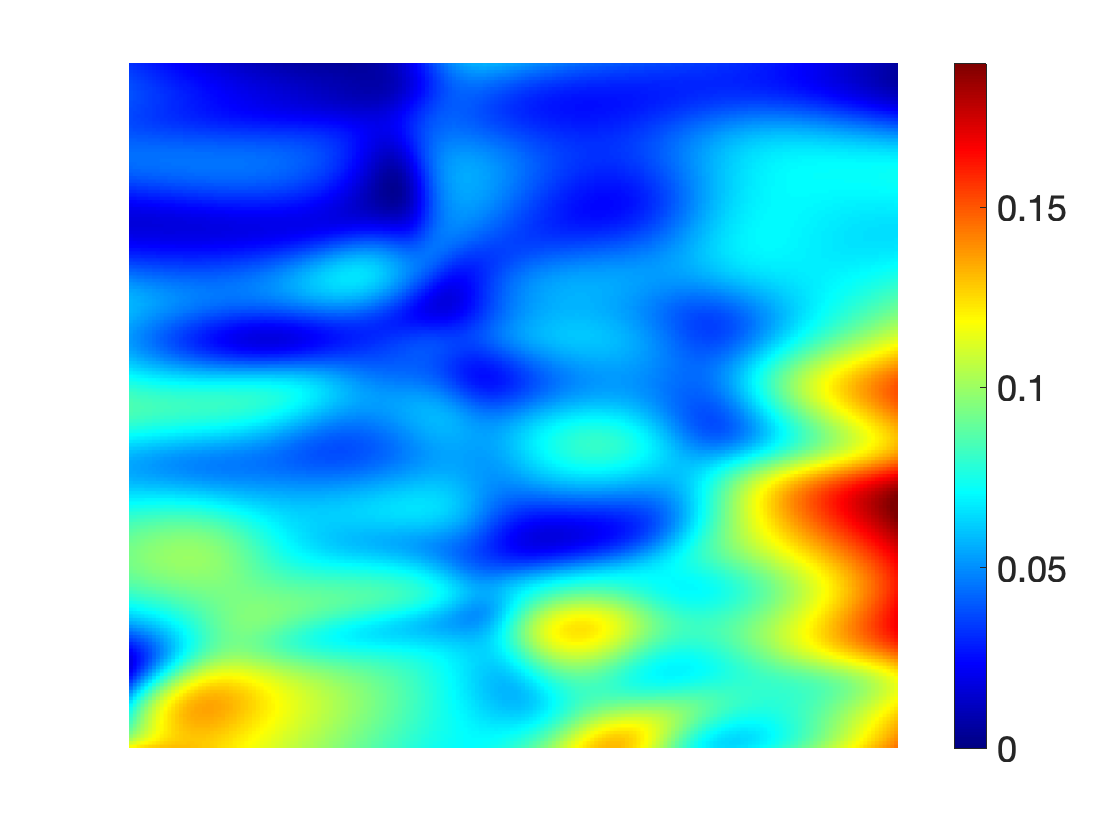} \\
(a) $A^\dag$  & (b) $\hat A$ & (c) $|\hat A-A^\dag|$ & (d) $\hat A$ & (e) $|\hat A-A^\dag|$
\end{tabular}
\caption{The reconstructions for Example \ref{exam:neu2d2} with exact data in (b) and noisy data $(\delta=10\%)$ in (d). From the top to bottom, the results are for $A_{11}$, $A_{12}$ and $A_{22}$, respectively.}
\label{fig:neu2d2}
\end{figure}

Fig. \ref{fig:neu2d2} presents the reconstructions for both exact and noisy data ($\delta = 10\%$). The DNN approximations successfully capture the oscillatory features of $A^\dagger$ in both cases, with only a mild deterioration in reconstruction quality when 10\% data noise is present. This shows its robustness against data noise.

\begin{figure}[htb!]
\centering
\setlength{\tabcolsep}{0em}
\begin{tabular}{ccccc}
\includegraphics[width=0.199\textwidth]{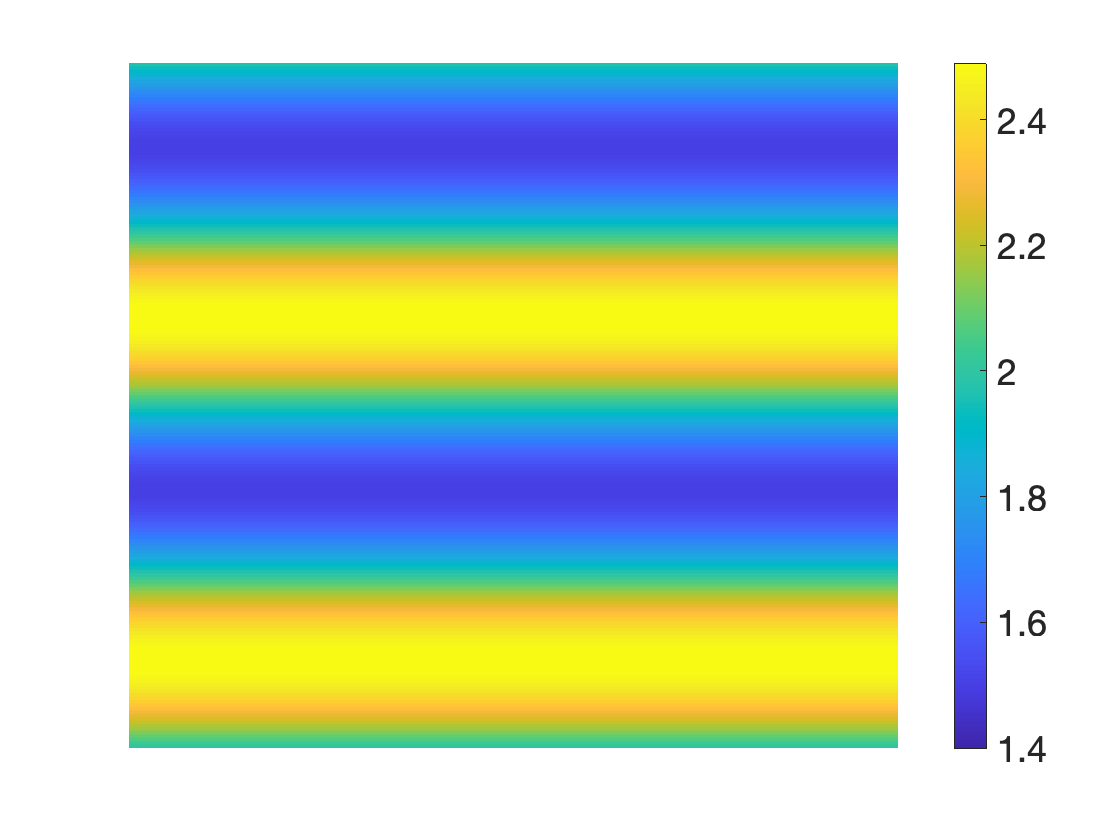} &
\includegraphics[width=0.199\textwidth]{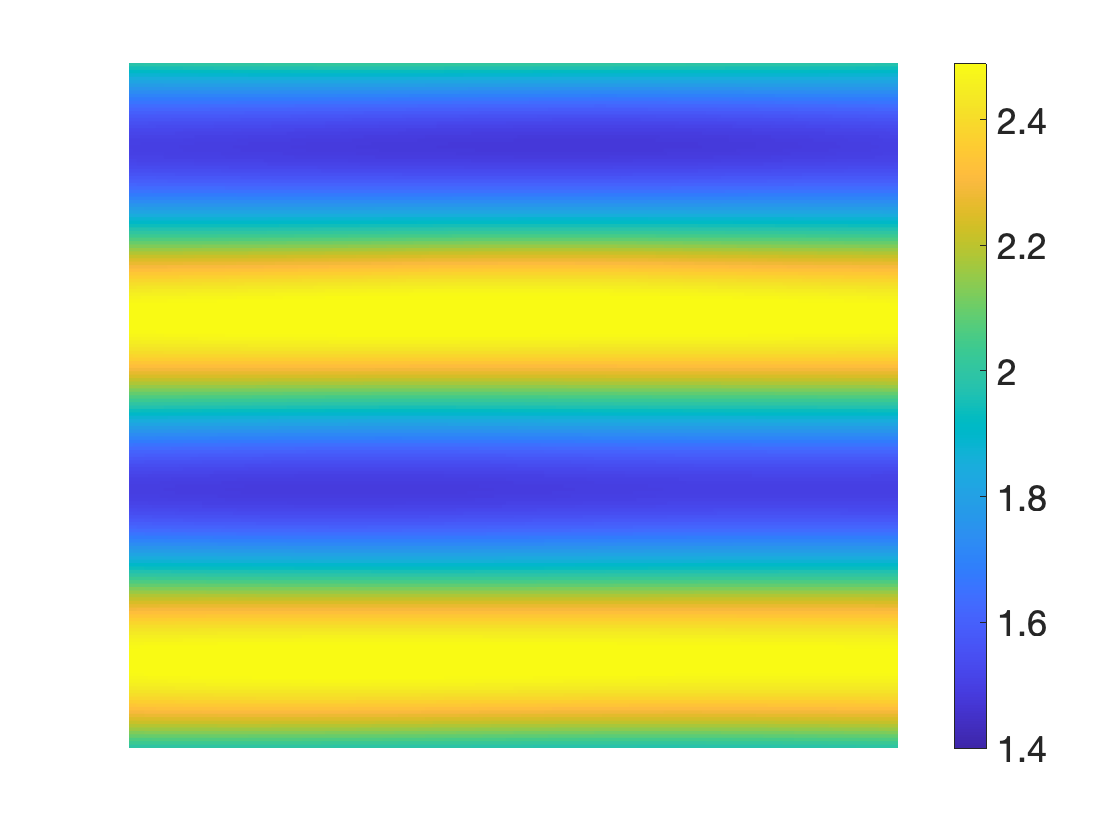} &
\includegraphics[width=0.199\textwidth]{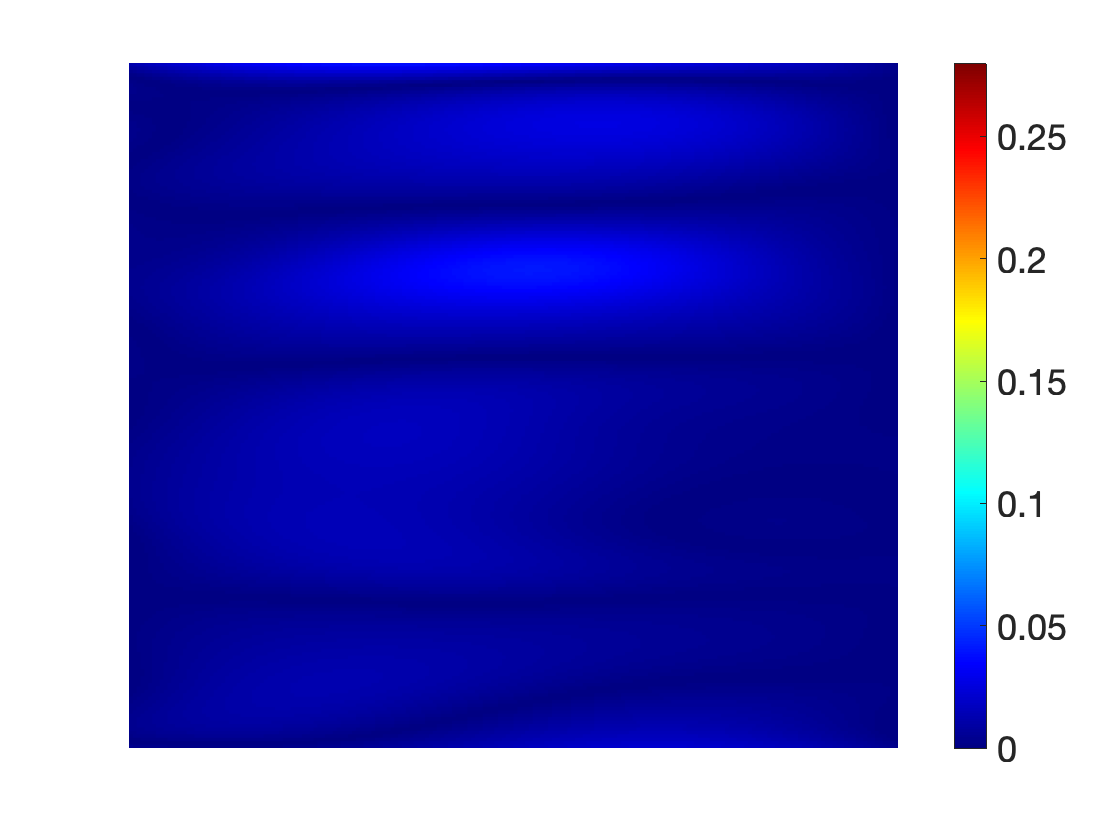} &
\includegraphics[width=0.199\textwidth]{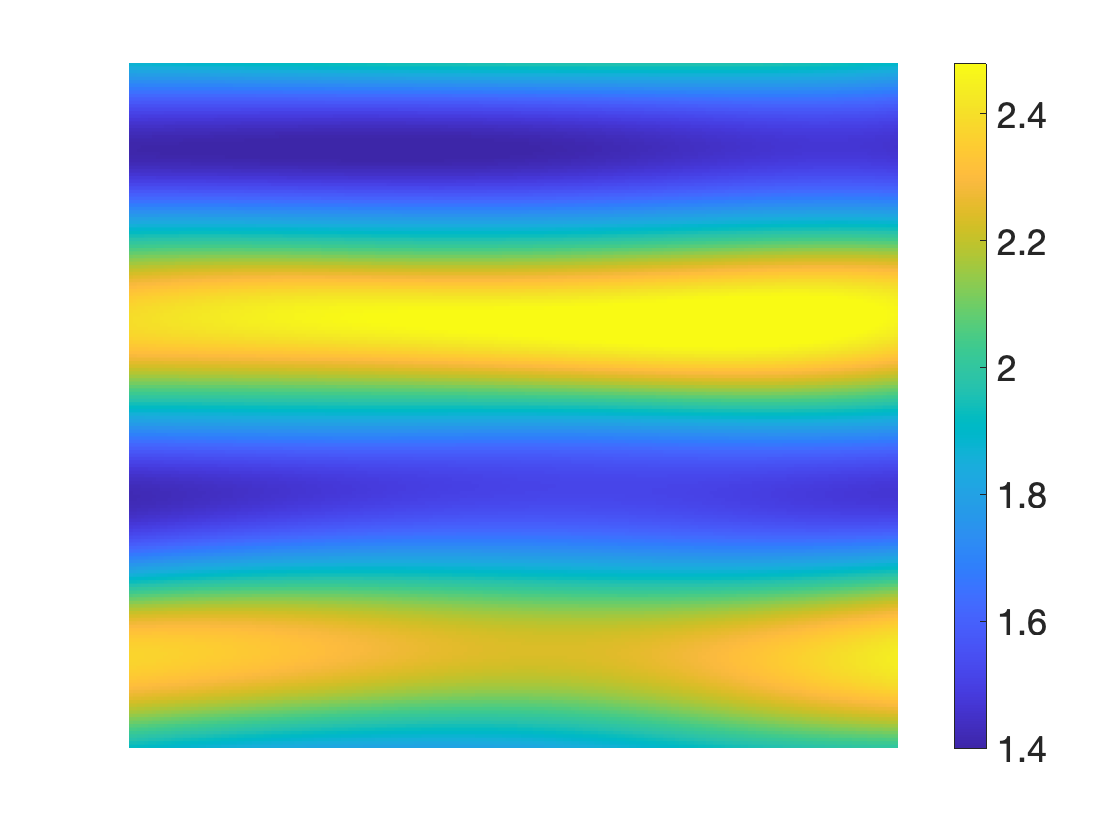} &
\includegraphics[width=0.199\textwidth]{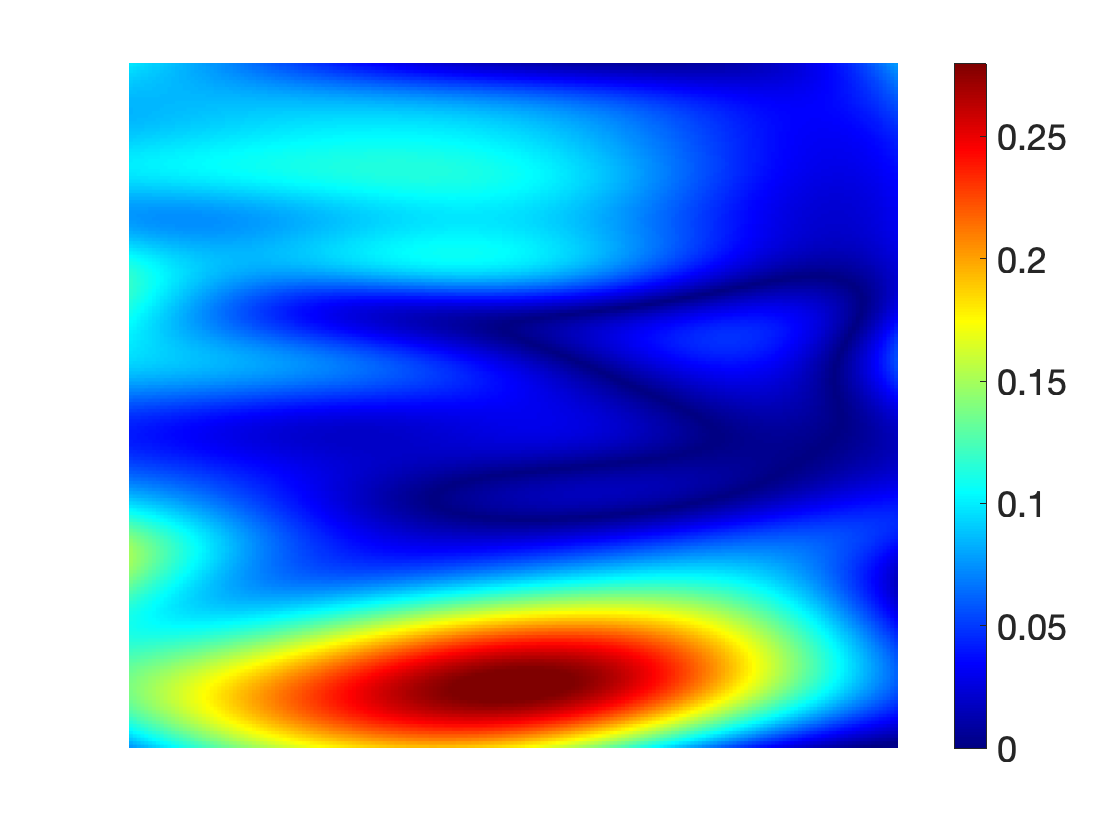} \\
\includegraphics[width=0.199\textwidth]{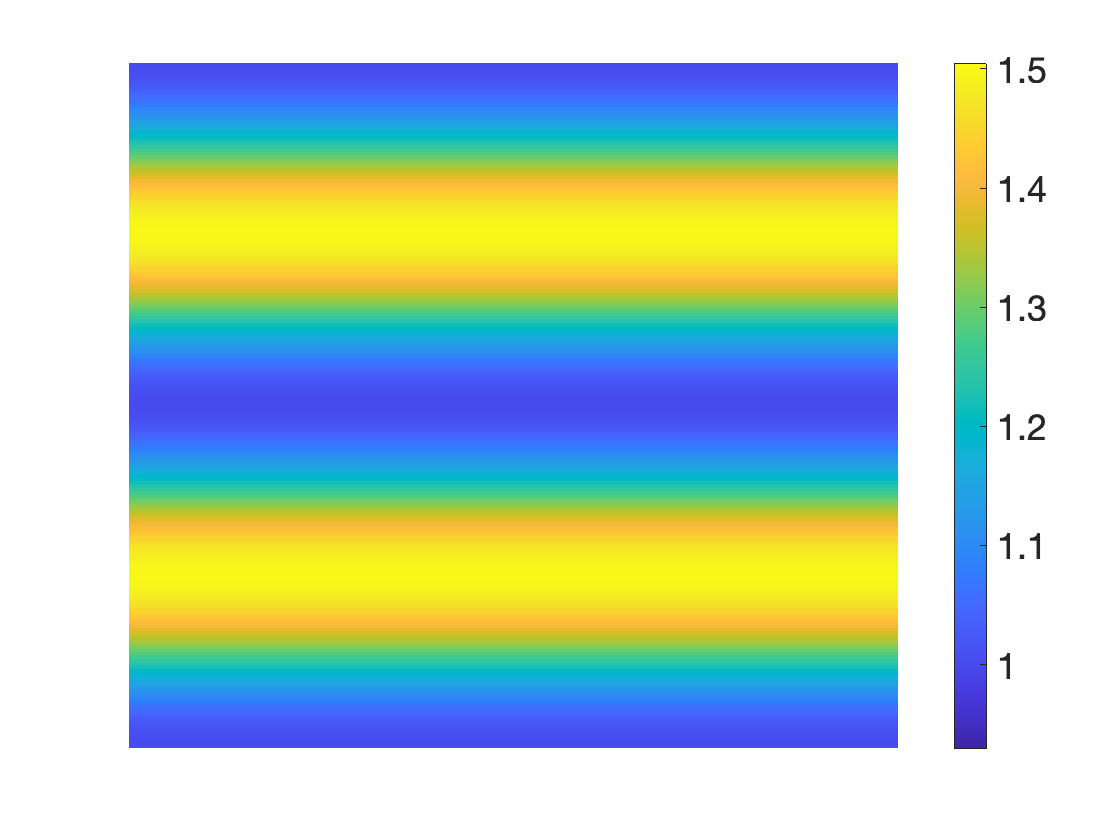} &
\includegraphics[width=0.199\textwidth]{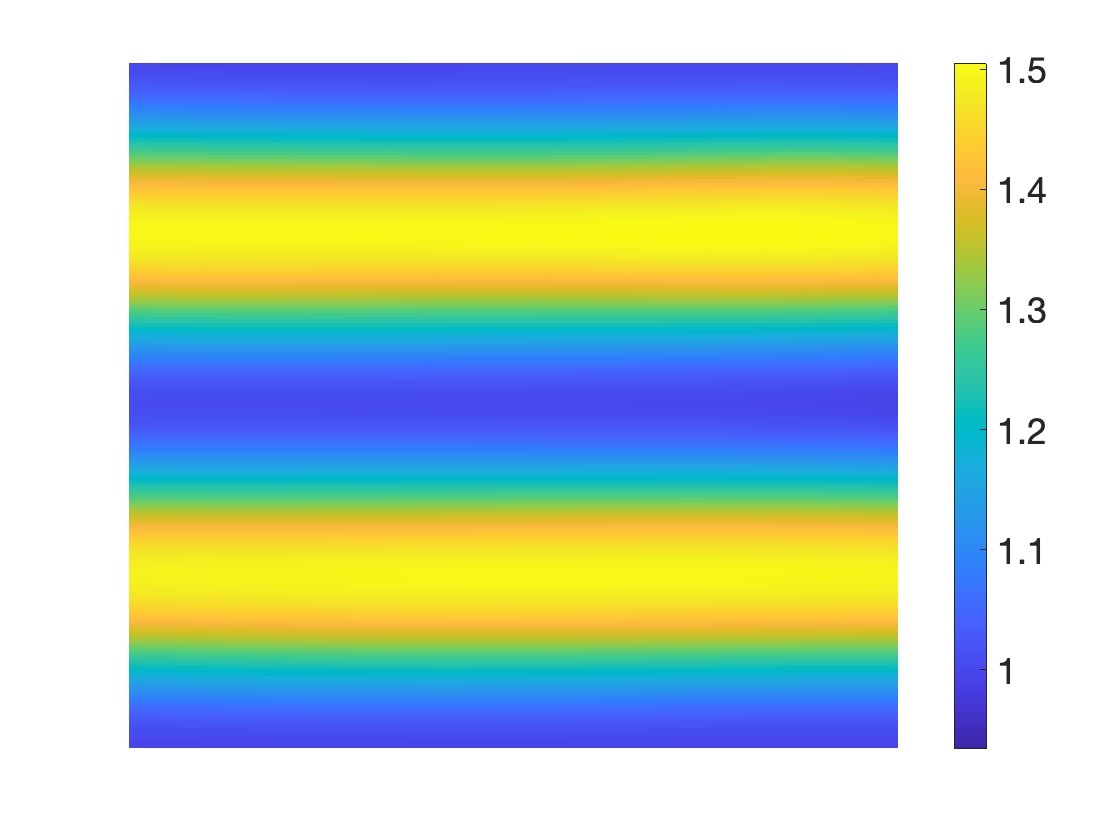} &
\includegraphics[width=0.199\textwidth]{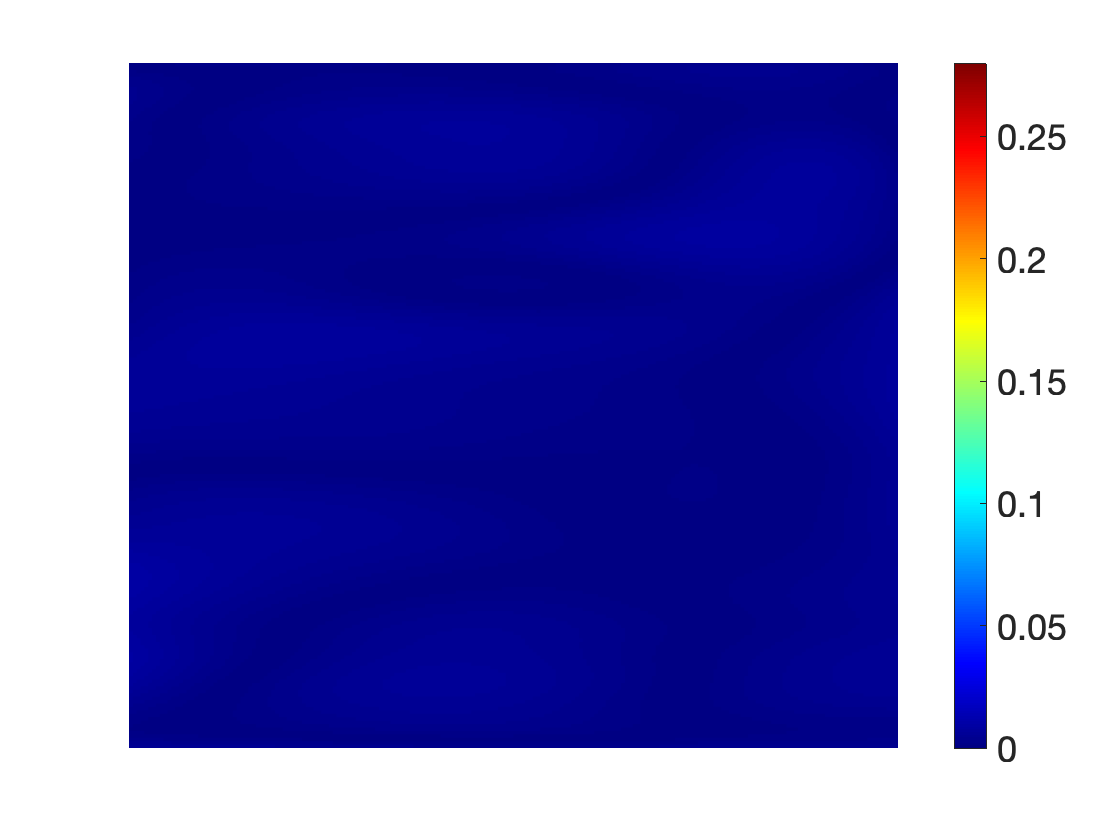} &
\includegraphics[width=0.199\textwidth]{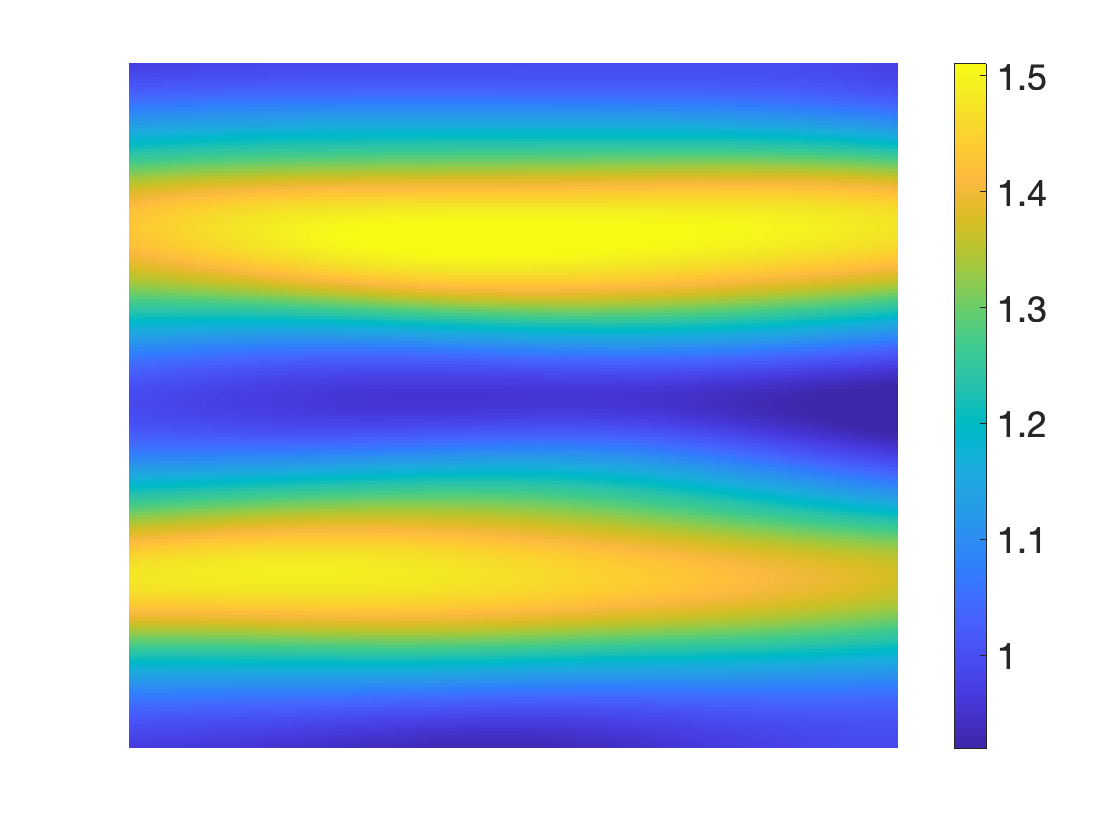} &
\includegraphics[width=0.199\textwidth]{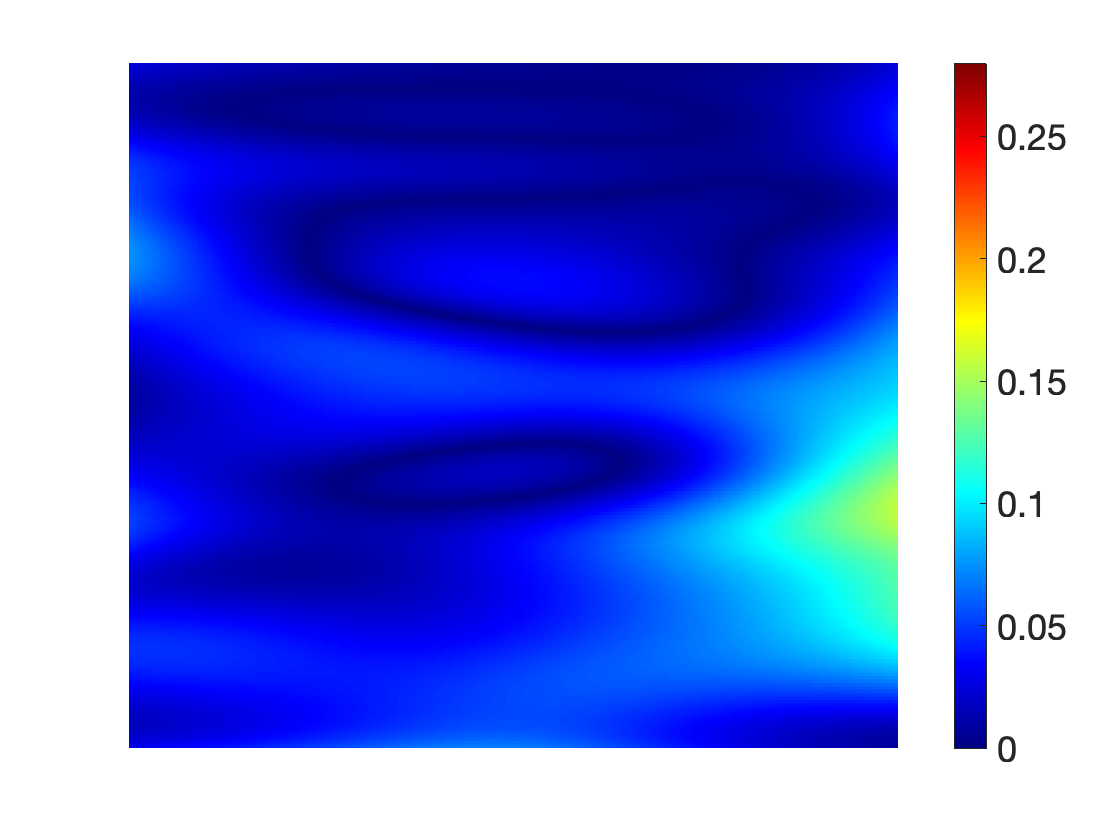} \\
\includegraphics[width=0.199\textwidth]{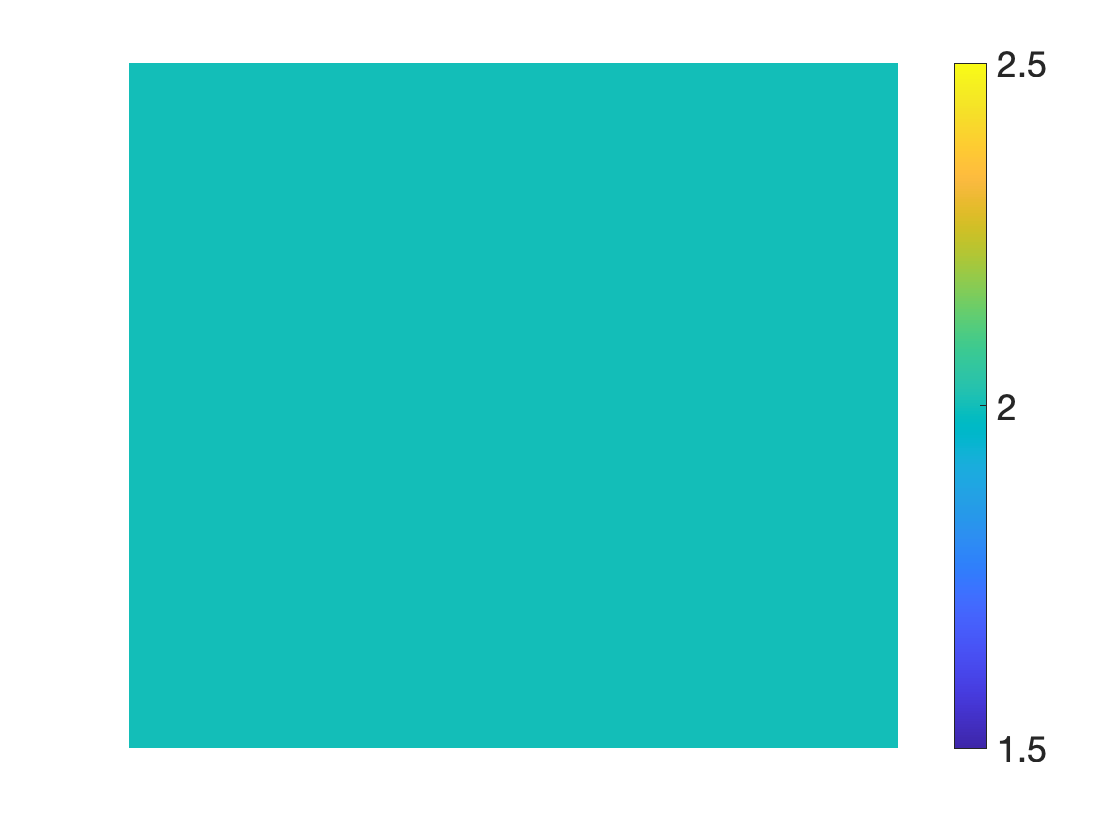} &
\includegraphics[width=0.199\textwidth]{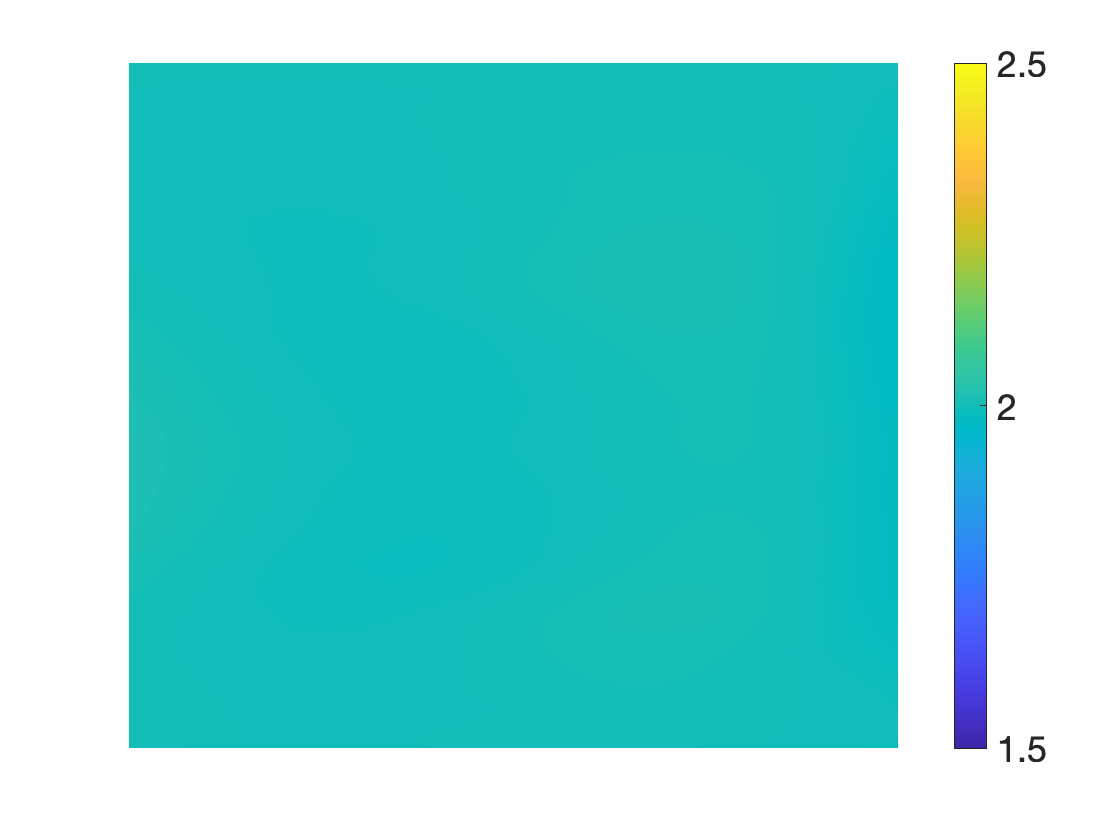} &
\includegraphics[width=0.199\textwidth]{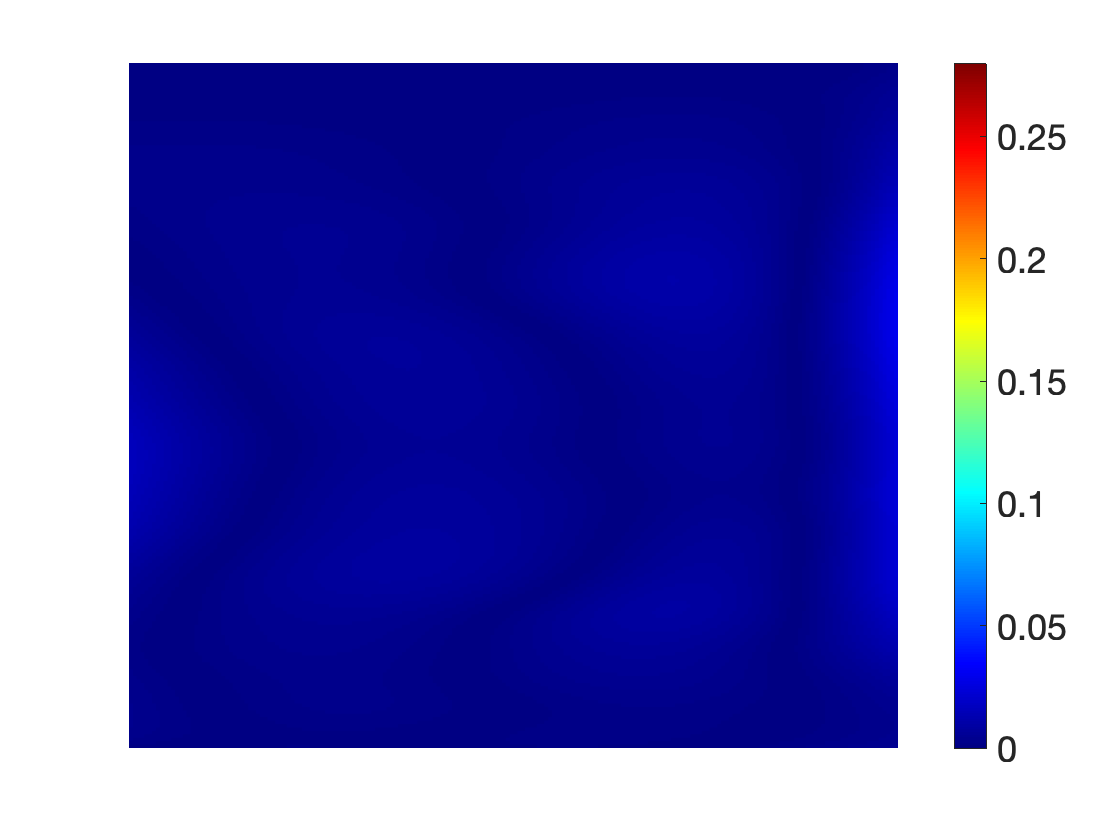} &
\includegraphics[width=0.199\textwidth]{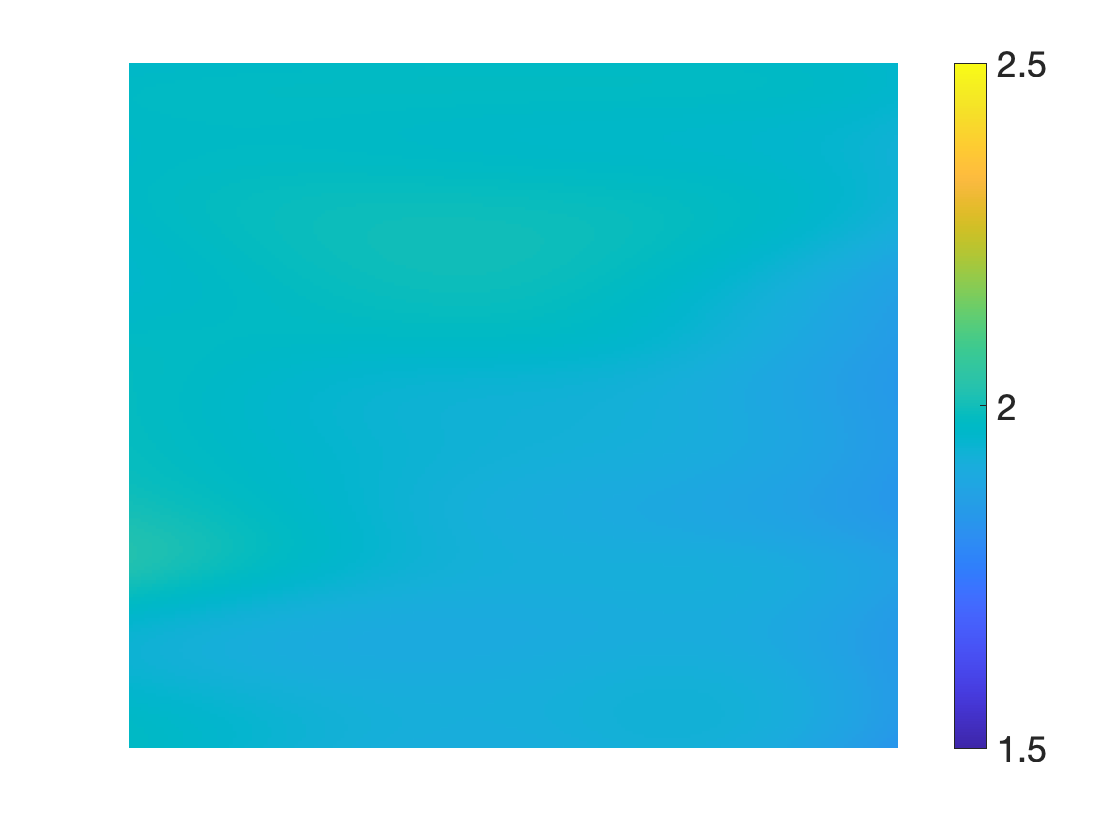} &
\includegraphics[width=0.199\textwidth]{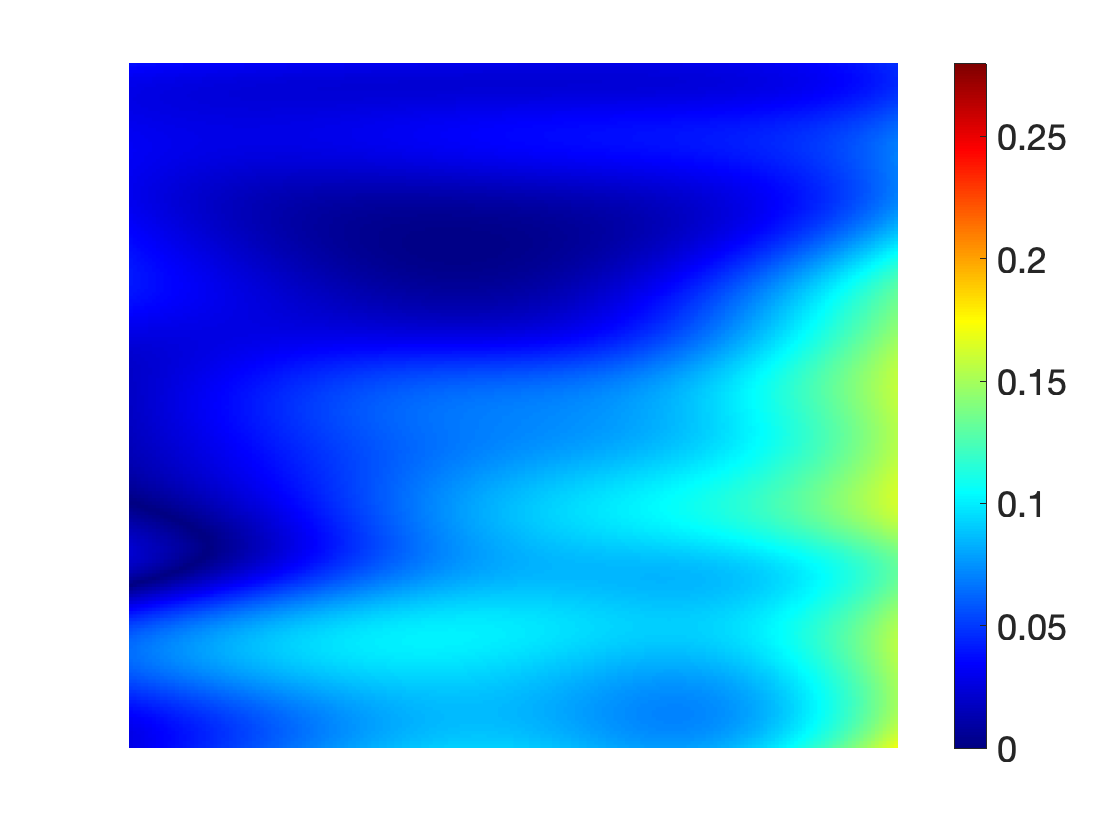} \\
(a) $A^\dag$  & (b) $\hat A$ & (c) $|\hat A-A^\dag|$ & (d) $\hat A$ & (e) $|\hat A-A^\dag|$
\end{tabular}
\caption{The reconstructions for Example \ref{exam:neu2d3} with exact data in (b) and noisy data $(\delta=10\%)$ in (d). From the top to bottom, the results are for $A_{11}$, $A_{12}$ and $A_{22}$, respectively.}
\label{fig:neu2d3}
\end{figure}

Next, we present a more ill-posed problem of using partial internal data.
\begin{example}
The domain $\Omega = (0,1)^2$, the measurement $\nabla  z_i^\delta$ on the region $\omega=\Omega\setminus(0.2,0.8)^2$, $A^\dag= \begin{pmatrix}
    2+\frac{\sin(4\pi x_2)}{2}&1+\frac{\sin^2(2\pi x_2)}{2}\\
    1+\frac{\sin^2(2\pi x_2)}{2}&2\\
\end{pmatrix},$ $u_1^\dag=x_1+x_2+\frac{1}{3}(x_1^3+x_2^3)$, $u_2^\dag=x_1-x_2+\frac{1}{3}(x_1^3-x_2^3)$, $u_3^\dag=-u_1^\dagger$ and $u_4^\dag=-u_2^\dagger$.
\label{exam:neu2d3}
\end{example}
Since the observational data $\nabla  z_i^\delta$ is only given over a subdomain $\omega$ - a narrow band near the domain boundary $\partial\Omega$ with a band width $0.2$, we modify the loss in \eqref{eqn:loss_NN} as
\begin{equation}\label{eqn:loss_neupartial}
\begin{aligned}
    J_{\gamma}(\theta,\kappa)=\sum_{i=1}^{N}&\left(\|\sigma_{i,\kappa}-P_{\mathcal{K}}(A_{\theta})\nabla z_i^{\delta}\|_{L^2(\omega)^d}^2+\gamma_{\sigma}\|\nabla\cdot\sigma_{i,\kappa}+f_i\|_{L^2(\Omega)}^2\right.\\
    &\left.+\gamma_b\| \vec{n}\cdot\sigma_{i,\kappa}-g_i \|_{L^2(\partial\Omega)}^2+\gamma_A\|P_{\mathcal{K}}(A_{\theta})\|_{L^2(\Omega)^{d,d}}^2\right).
    \end{aligned}
\end{equation}
Fig. \ref{fig:neu2d3} shows that the reconstruction results are accurate across the entire domain $\Omega$, including the central region where no observational data $\nabla z_i^\delta$ is provided. When 10\% data noise is present, there is only a mild deformation and underestimation in the oscillation peaks. This highlights the stability of the proposed MLS-DNN approach for partial internal data.

\begin{figure}[htb!]
\centering
\setlength{\tabcolsep}{0em}
\begin{tabular}{ccccc}
\includegraphics[width=0.21\textwidth]{a11neu2dpex.png} &
\includegraphics[width=0.199\textwidth]{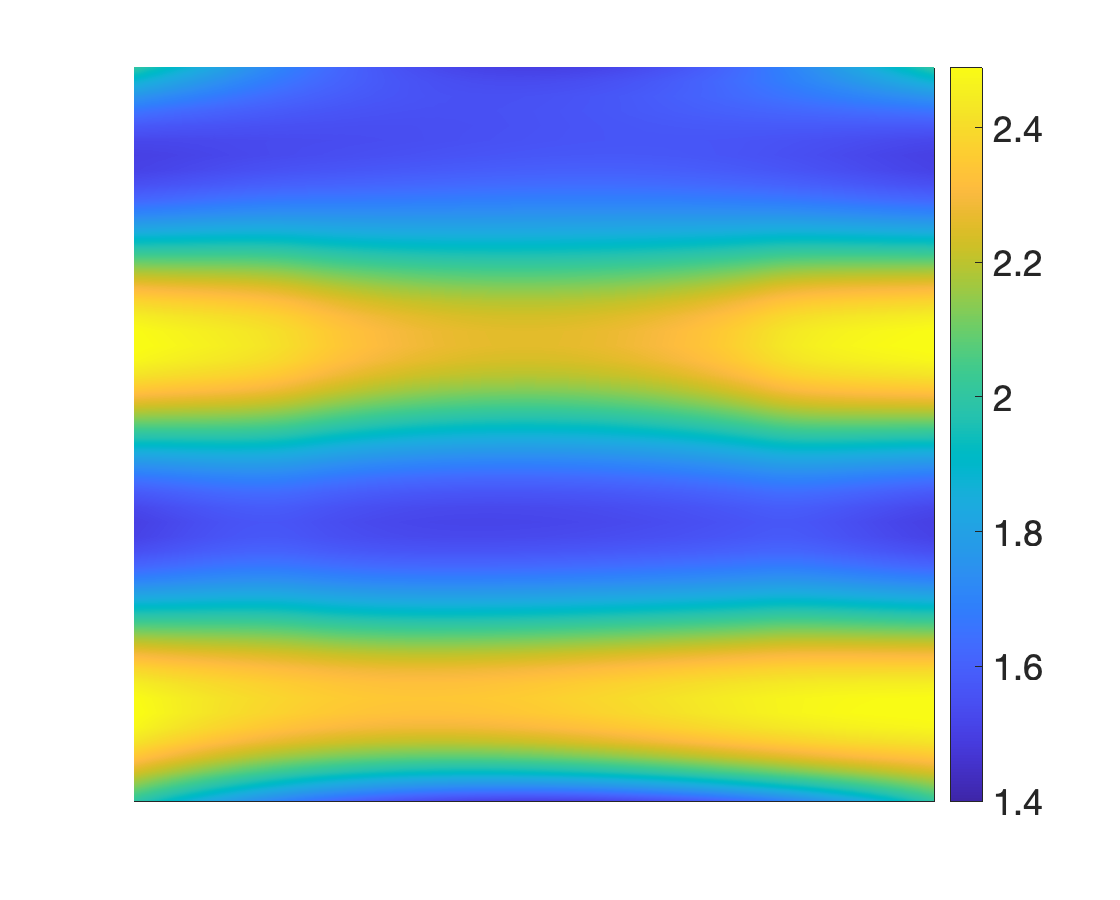} &
\includegraphics[width=0.199\textwidth]{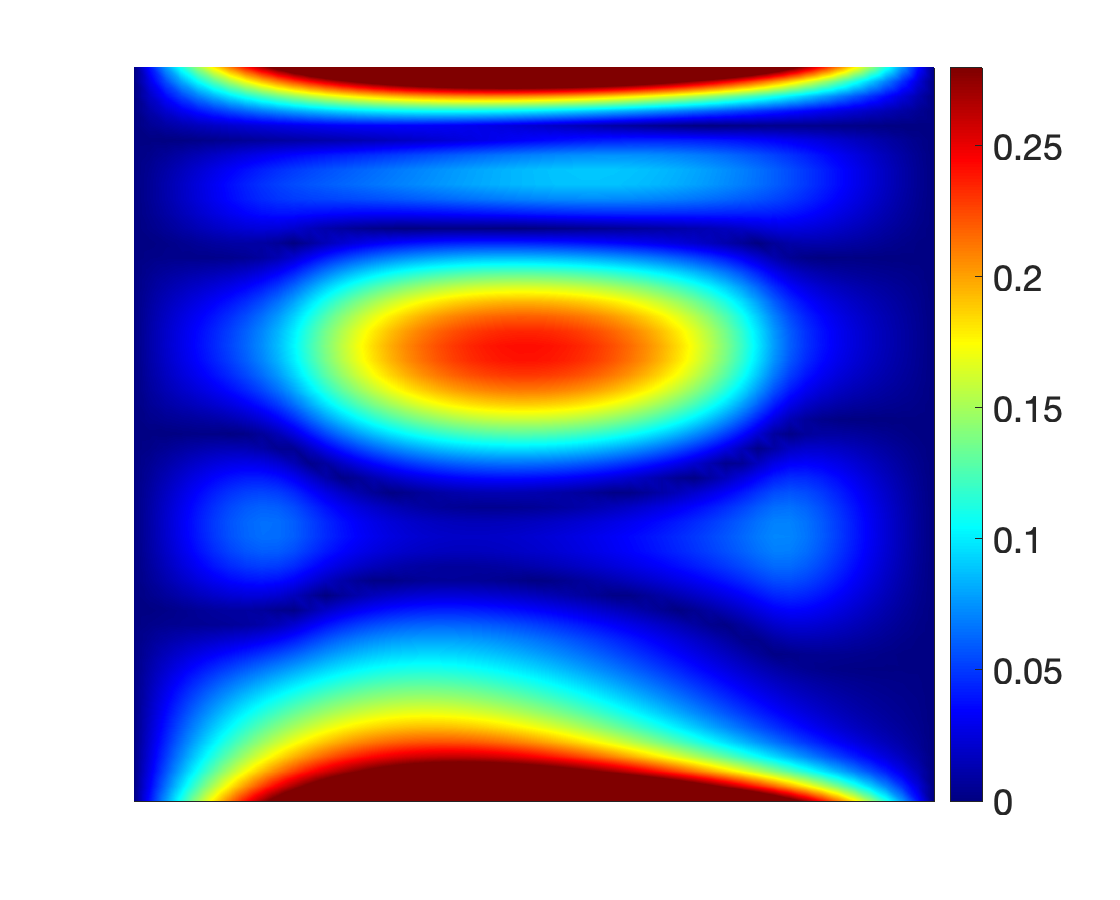} &
\includegraphics[width=0.199\textwidth]{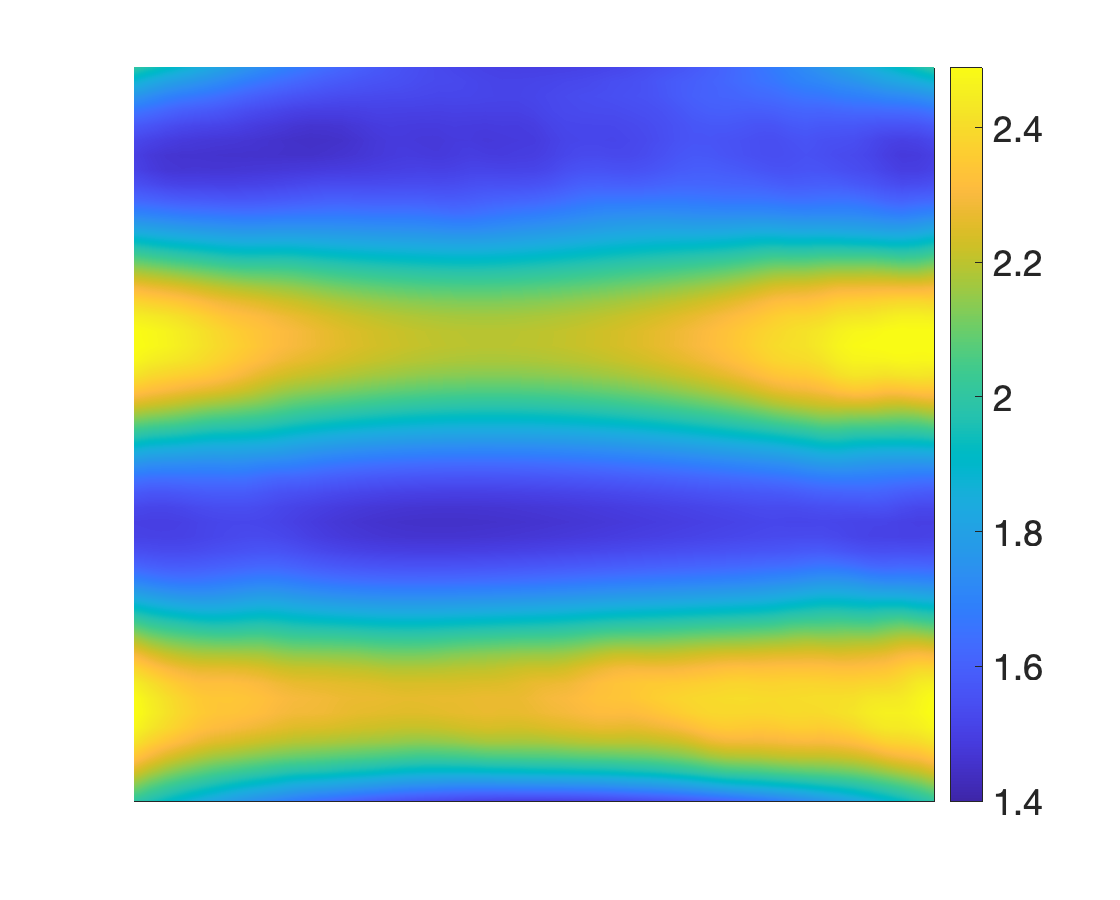} &
\includegraphics[width=0.199\textwidth]{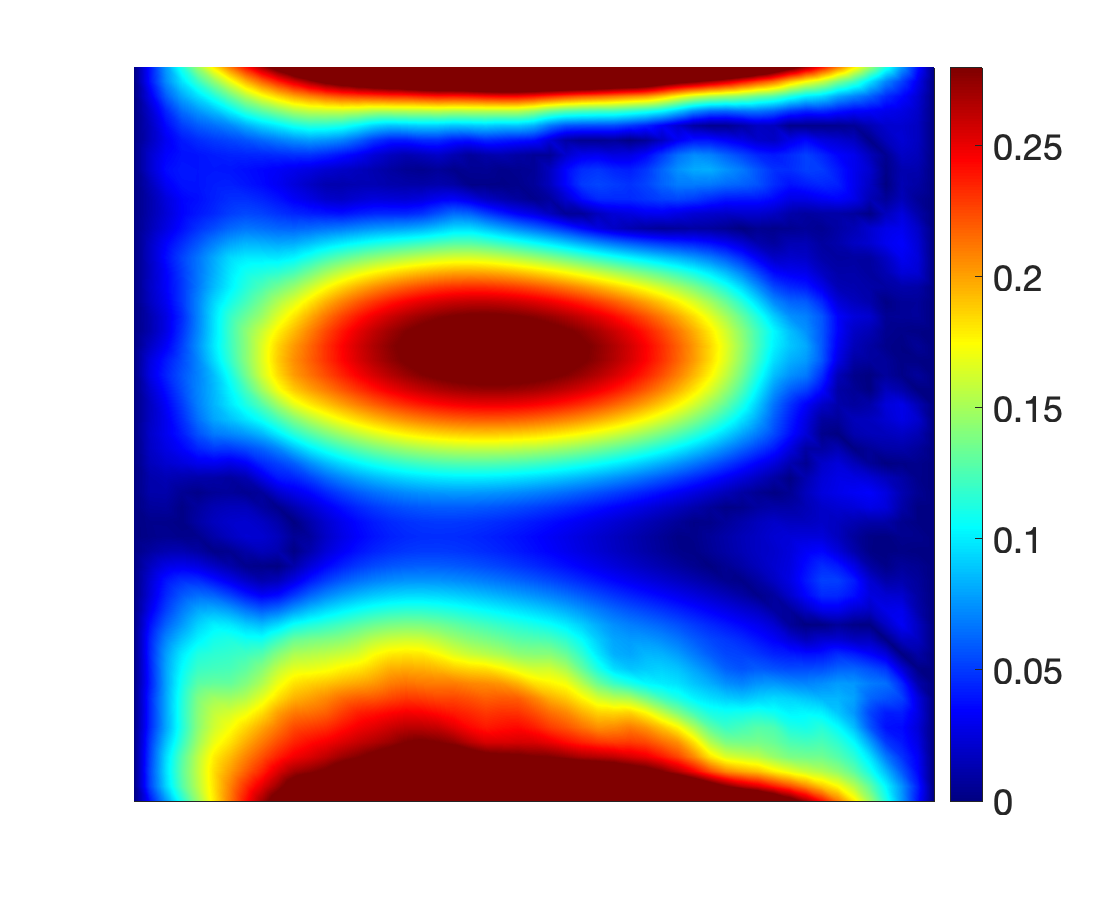} \\
\includegraphics[width=0.21\textwidth]{a12neu2dpex.png} &
\includegraphics[width=0.199\textwidth]{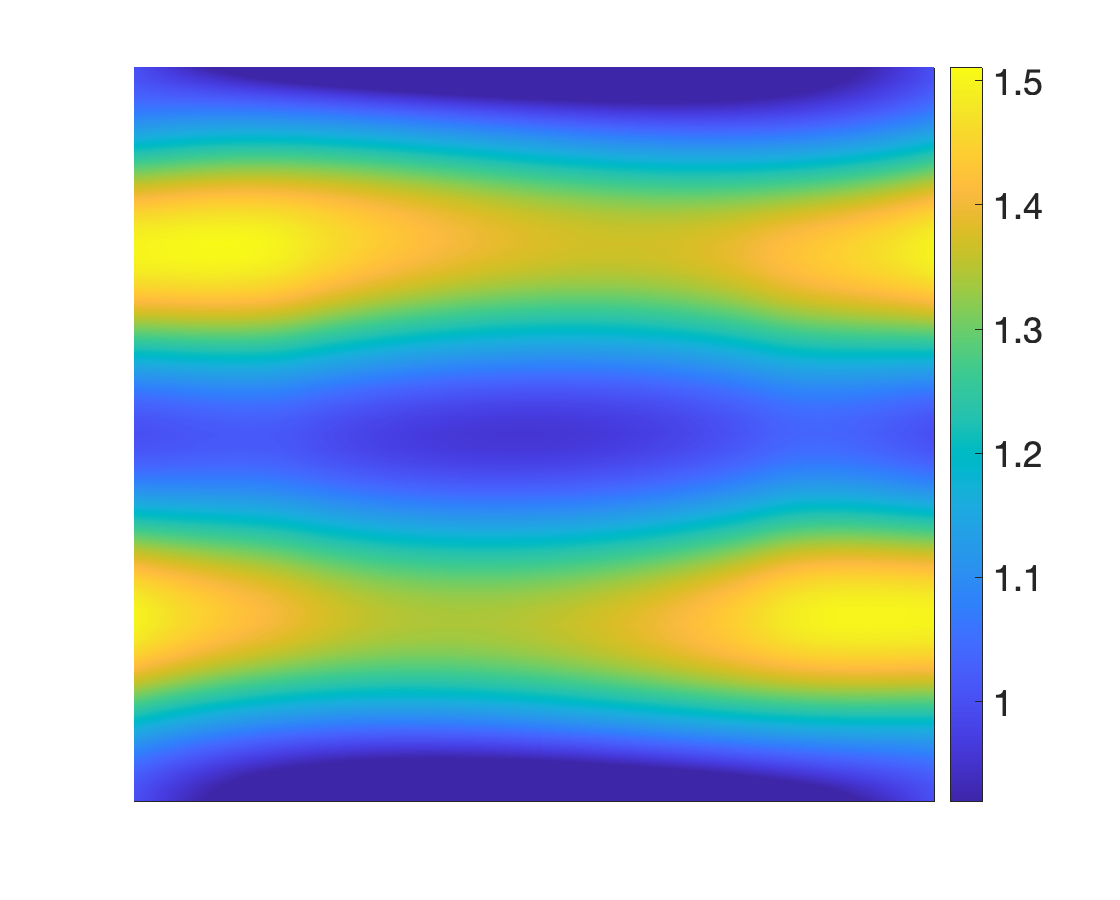} &
\includegraphics[width=0.199\textwidth]{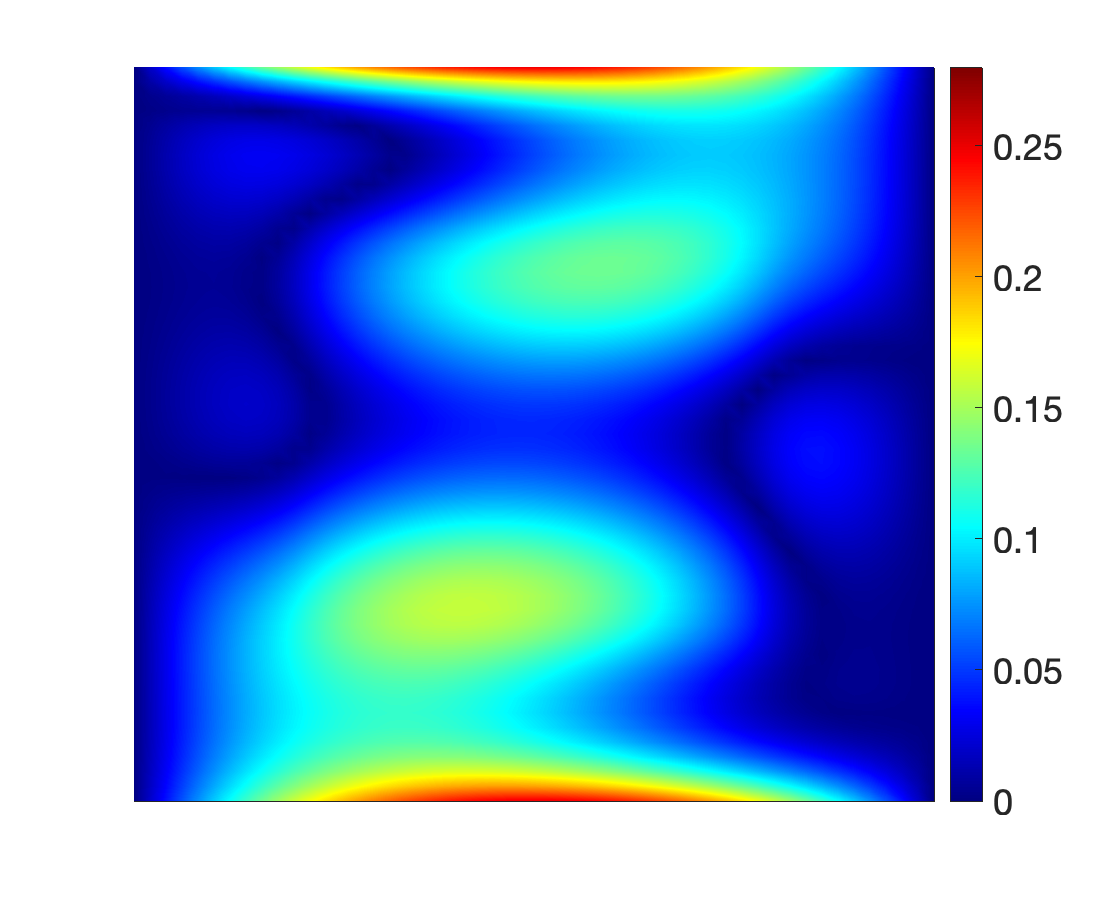} &
\includegraphics[width=0.199\textwidth]{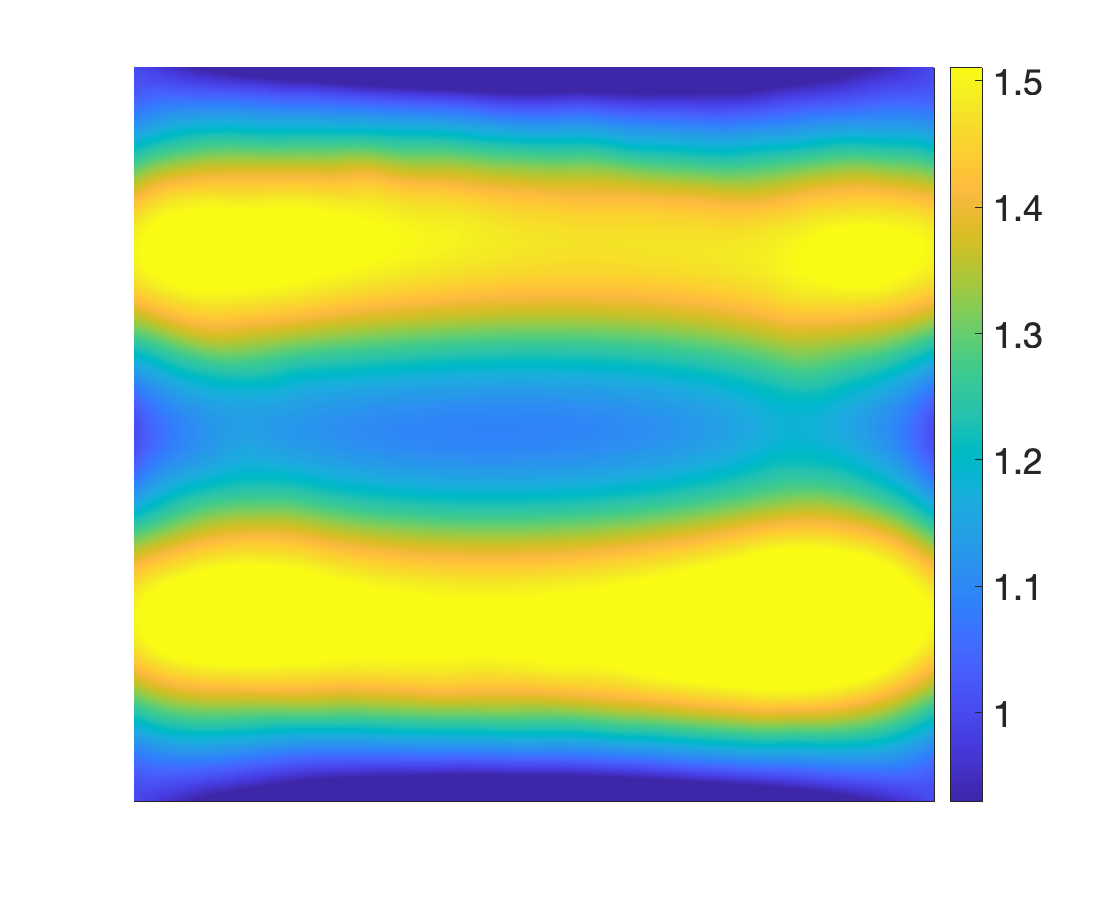} &
\includegraphics[width=0.199\textwidth]{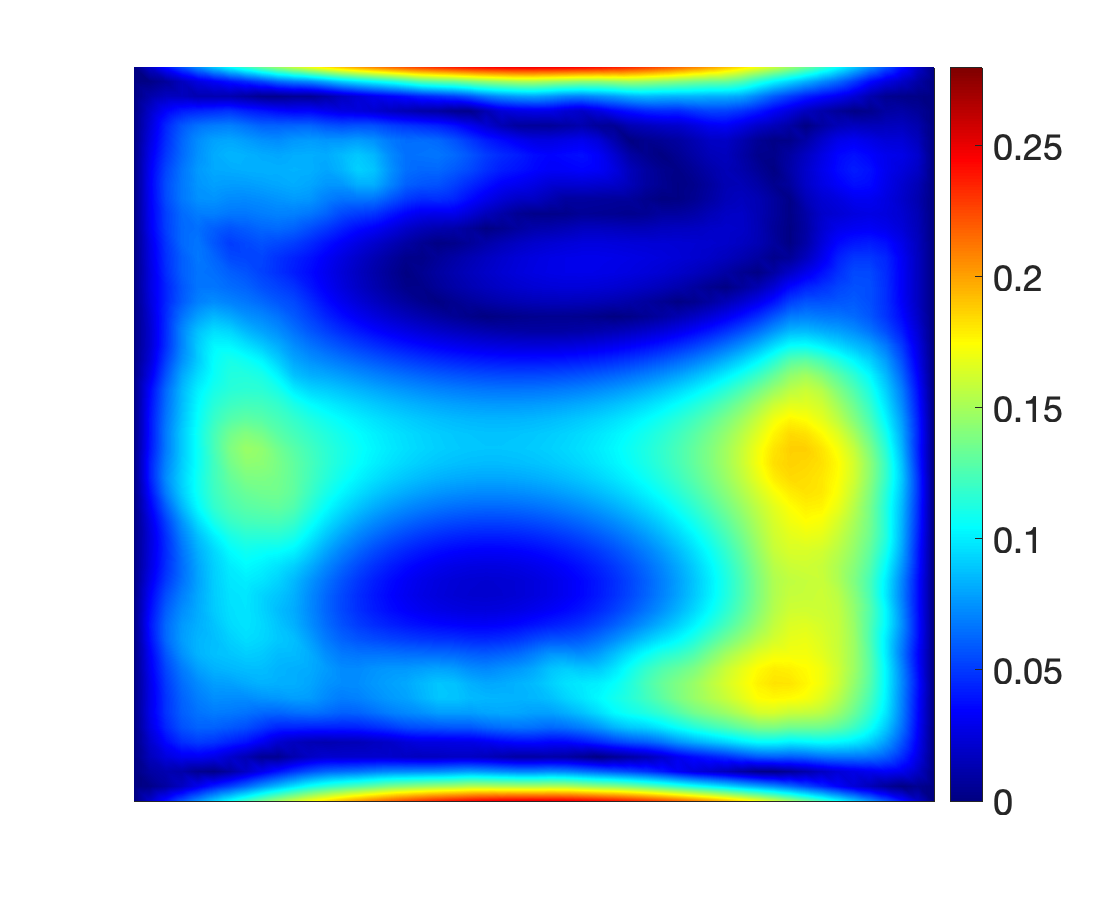} \\
\includegraphics[width=0.21\textwidth]{a22neu2dpex.png} &
\includegraphics[width=0.199\textwidth]{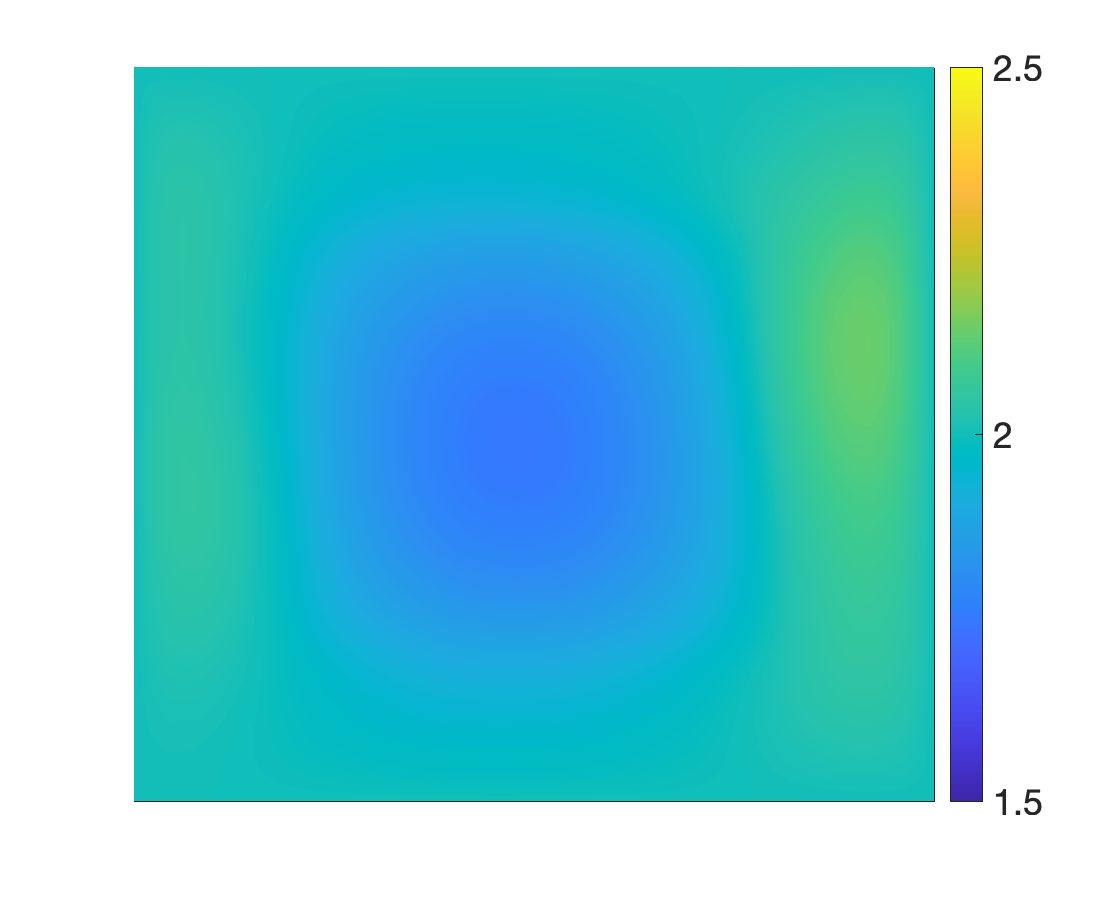} &
\includegraphics[width=0.199\textwidth]{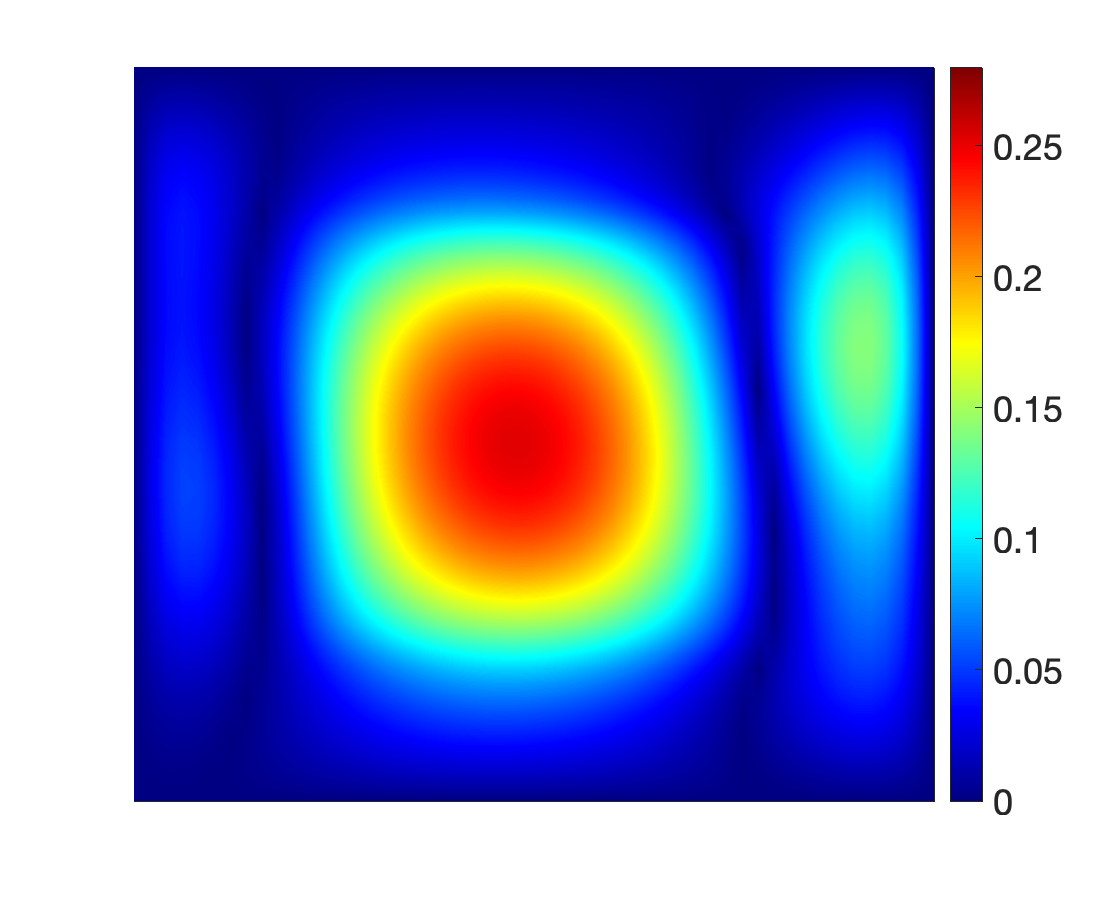} &
\includegraphics[width=0.199\textwidth]{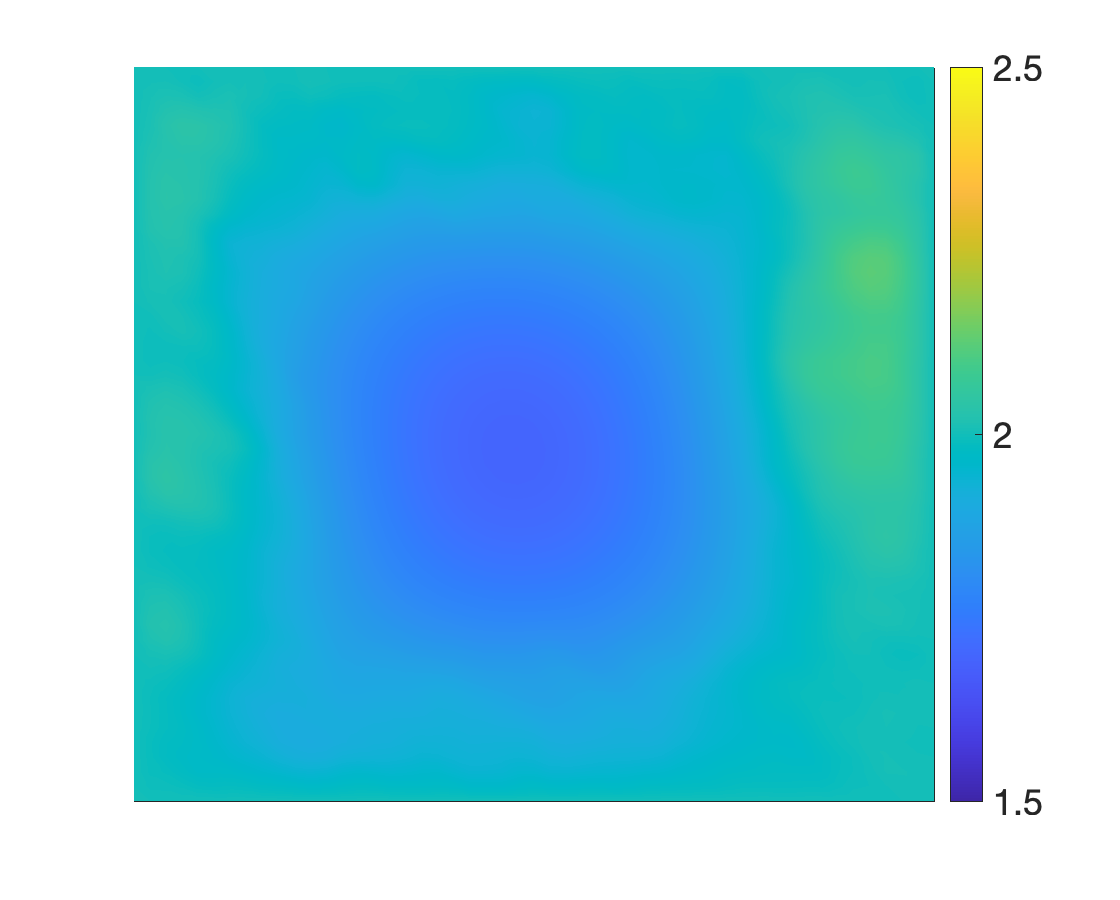} &
\includegraphics[width=0.199\textwidth]{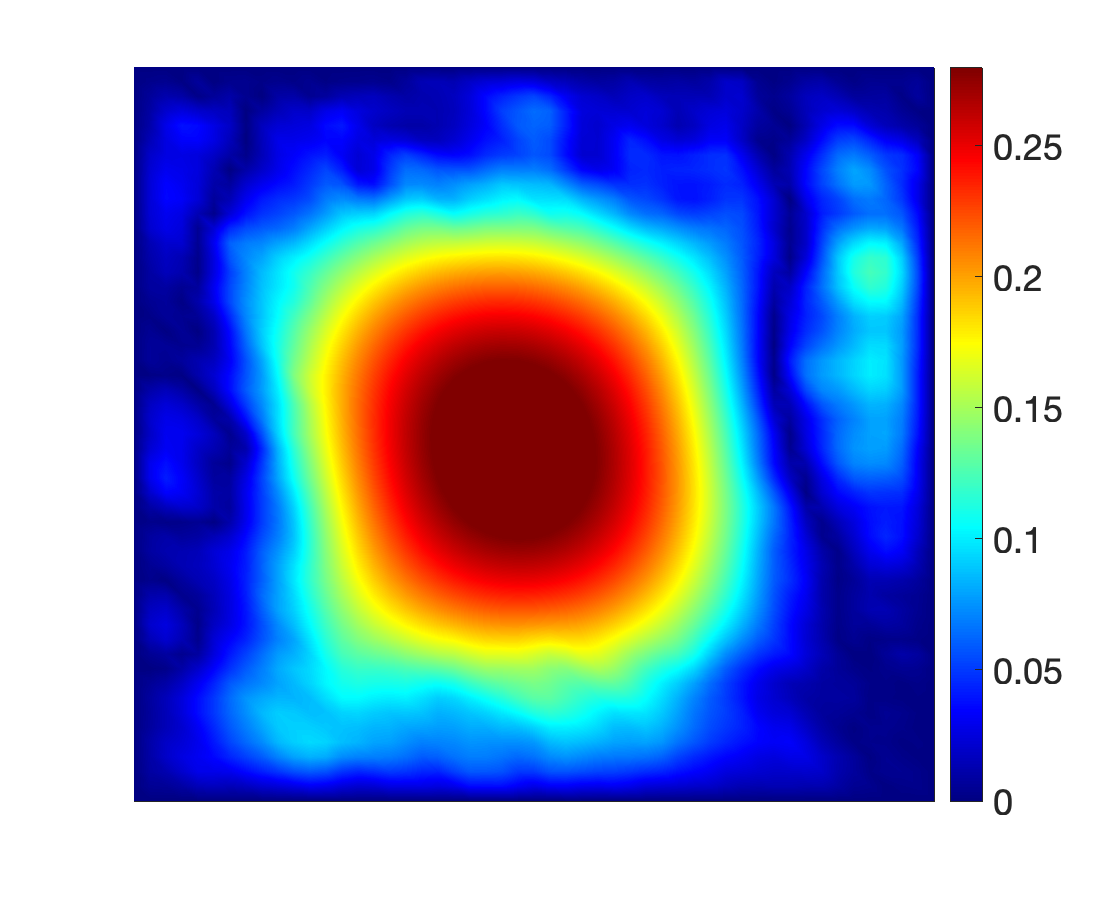} \\
(a) $A^\dag$  & (b) $\hat A$ & (c) $|\hat A-A^\dag|$ & (d) $\hat A$ & (e) $|\hat A-A^\dag|$
\end{tabular}
\caption{The reconstructions for Example \ref{exam:neu2d3} with exact data in (b) and noisy data $(\delta=10\%)$ in (d) using the FEM. From the top to bottom, the results are for $A_{11}$, $A_{12}$ and $A_{22}$, respectively.}
\label{fig:neufem}
\end{figure}

We also present the numerical results obtained using the Galerkin FEM for Example \ref{exam:neu2d3}. Specifically, we discretize the loss function \eqref{eqn:loss_neupartial} using the standard P1 FEM. Our initial guess is set as follows: $A_{11}^0= (2+\sin(4\pi x_2)/2)(1+x_1(x_1-1)) $, $A_{12}^0= (1+\sin(2\pi x_2)^2/2)(1+x_1(x_1-1)) $, $A_{22}^0= 2(1+4x_1x_2(x_1-1)(1-x_2)) $ and $\sigma_i^0=A^0\nabla u_i^\dagger$. The penalty parameters are chosen to be $\gamma_\sigma=\gamma_b=1$, $\gamma_A=$1e-6 for reconstructions using both exact and noisy data. Fig. \ref{fig:neufem} shows the reconstruction results obtained using the FEM. Due to the highly ill-posed nature of the inverse problem, the FEM reconstructions exhibit poor stability, with pronounced errors in resolving the conductivity features, even in the noise-free case.

\begin{figure}[htb!]
\centering
\setlength{\tabcolsep}{0em}
\begin{tabular}{ccccc}
\includegraphics[width=0.199\textwidth]{a11neu2dpex.png} &
\includegraphics[width=0.199\textwidth]{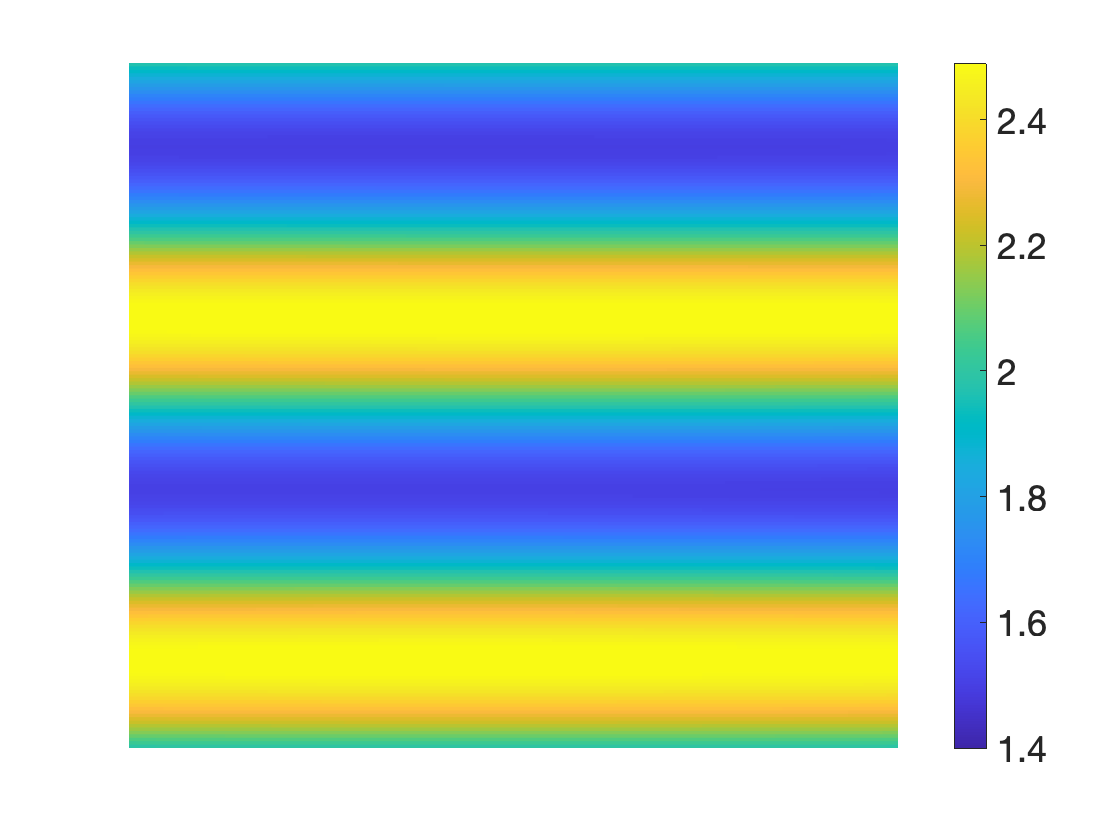} &
\includegraphics[width=0.199\textwidth]{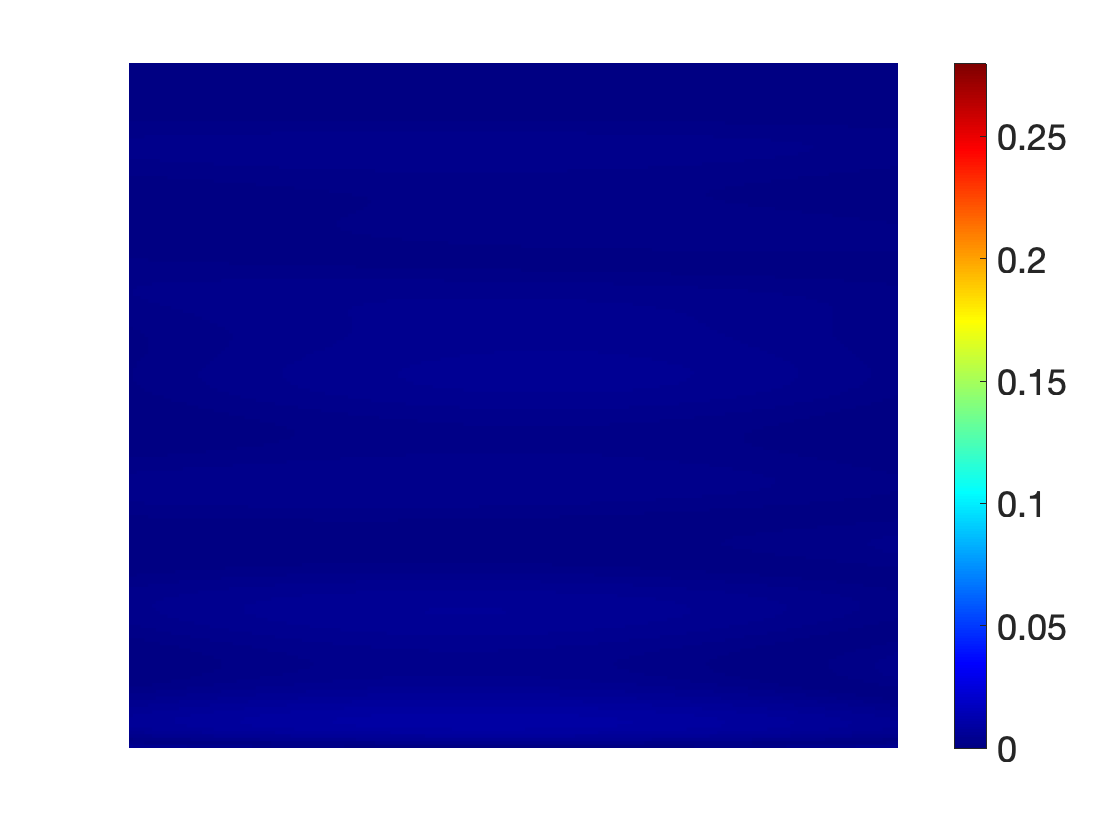} &
\includegraphics[width=0.199\textwidth]{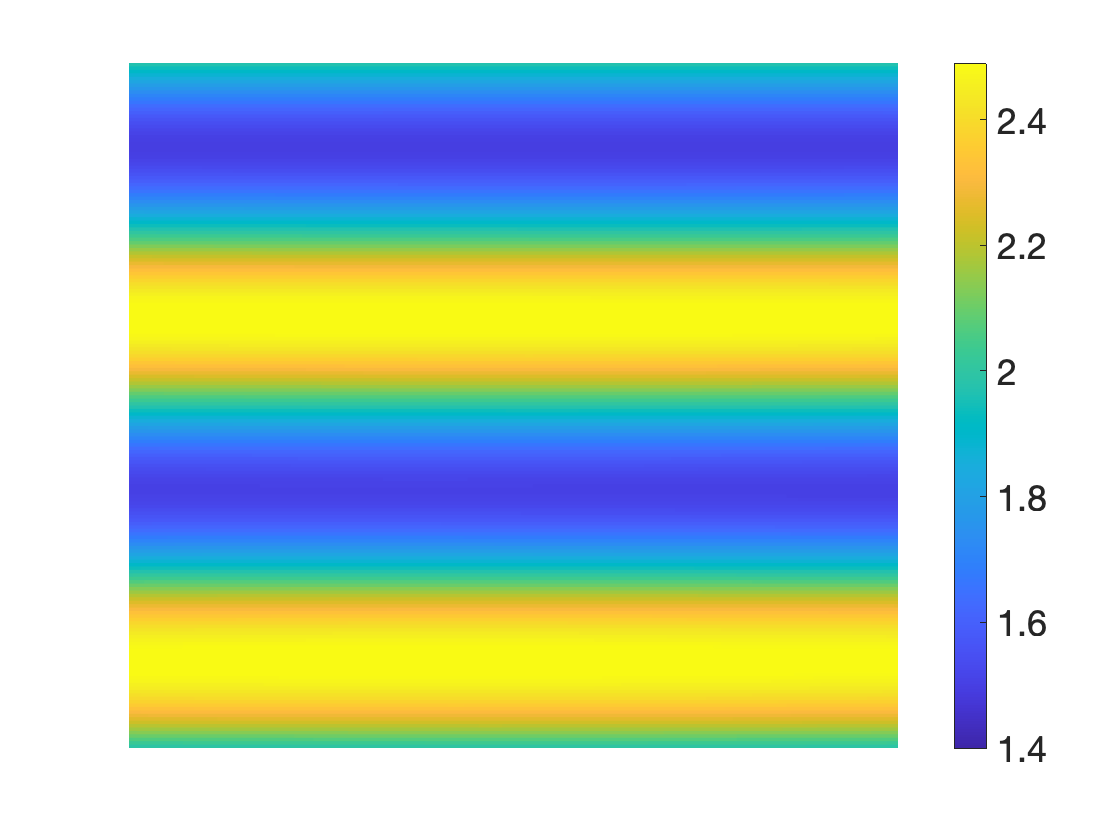} &
\includegraphics[width=0.199\textwidth]{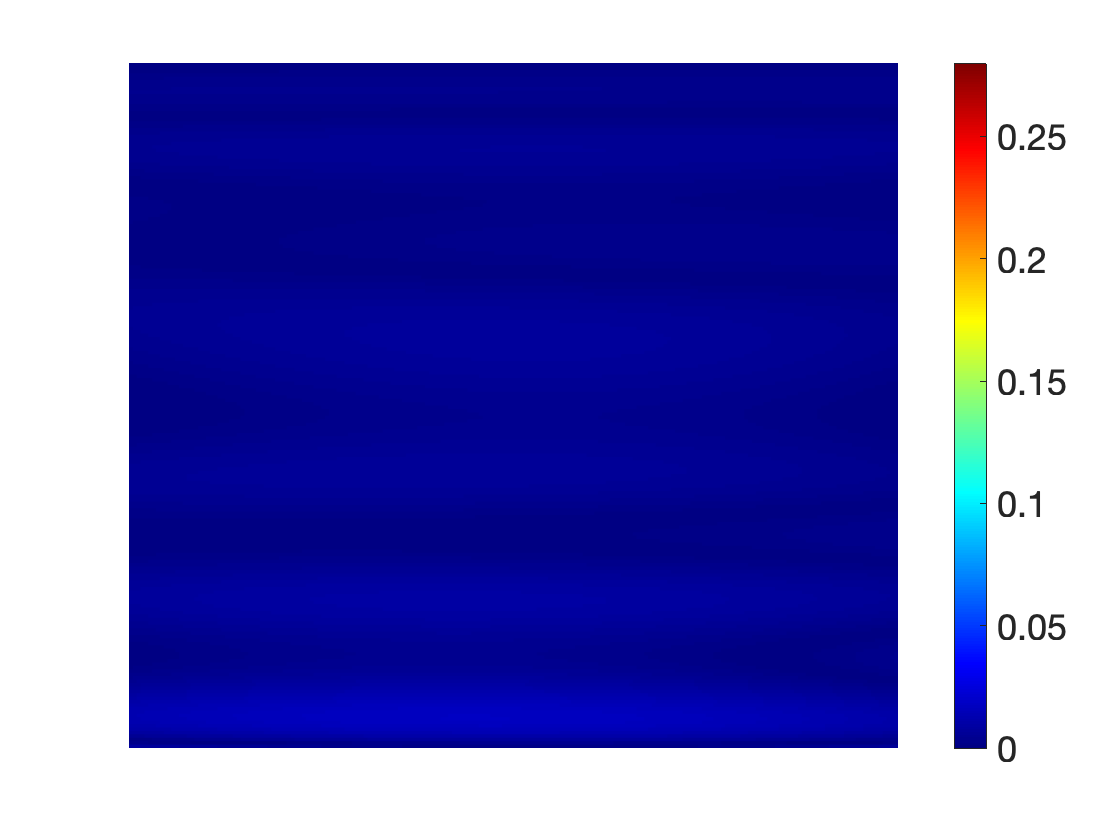} \\
\includegraphics[width=0.199\textwidth]{a12neu2dpex.png} &
\includegraphics[width=0.199\textwidth]{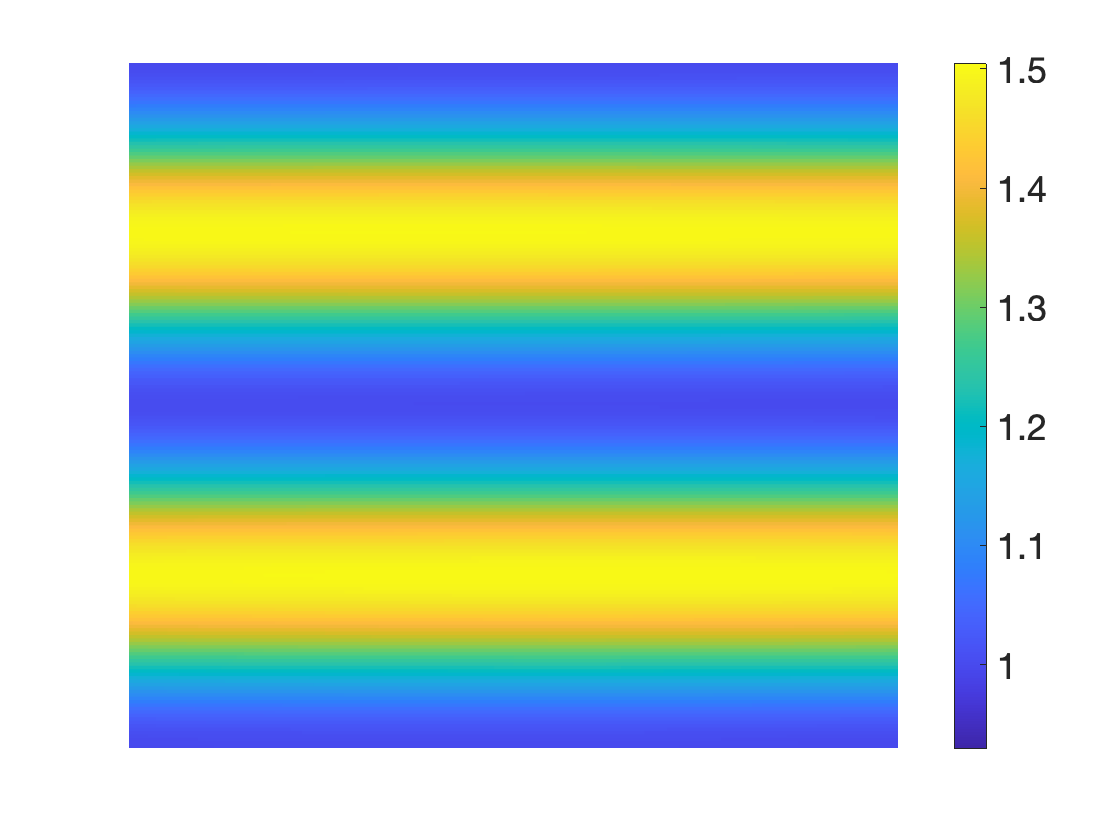} &
\includegraphics[width=0.199\textwidth]{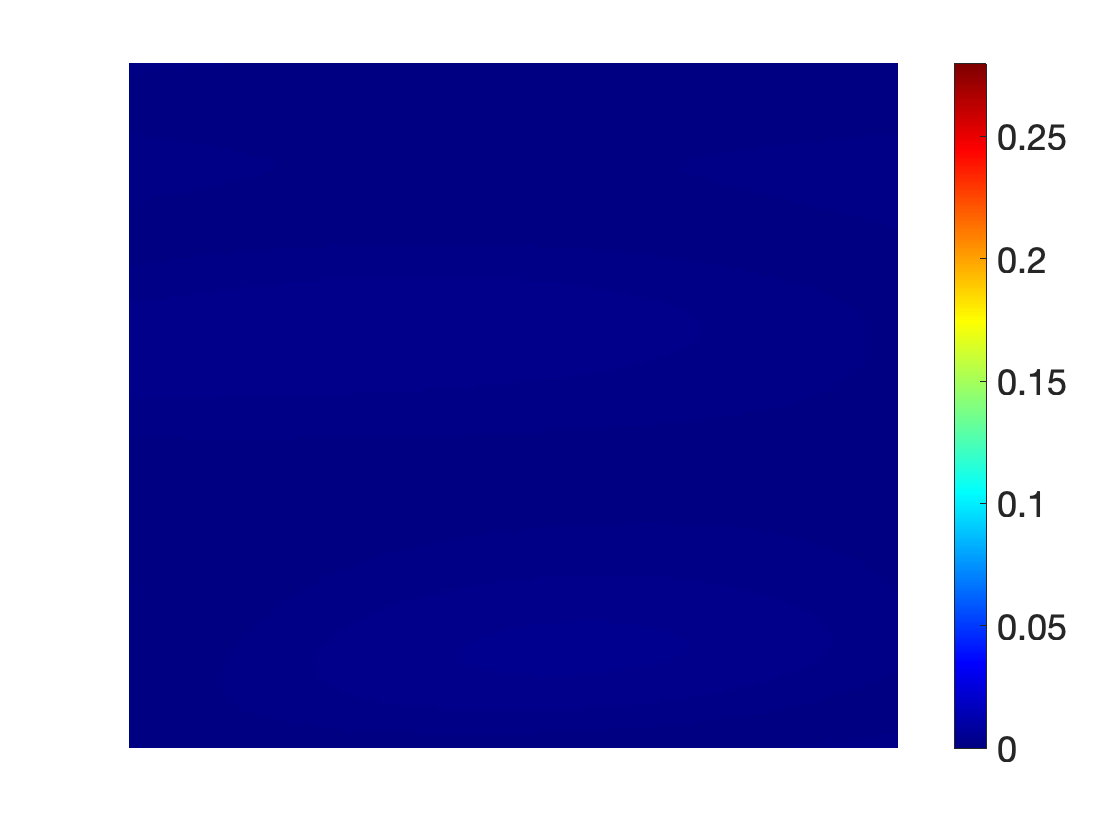} &
\includegraphics[width=0.199\textwidth]{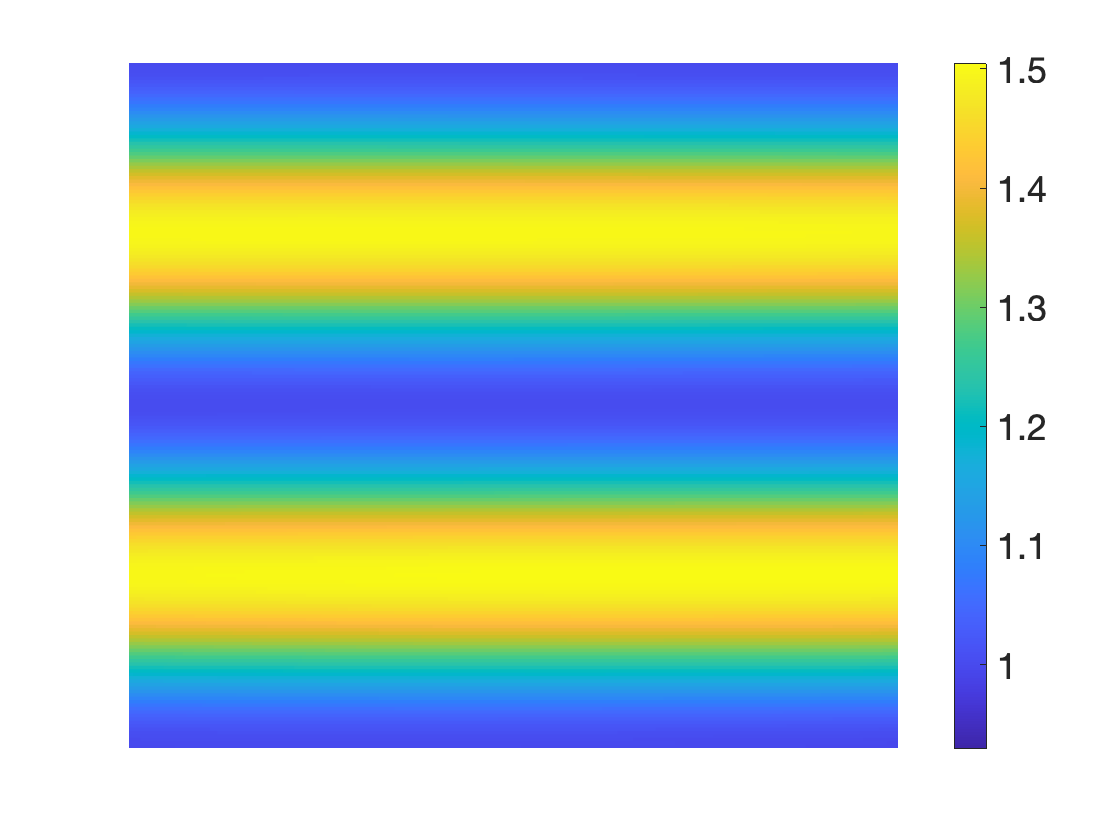} &
\includegraphics[width=0.199\textwidth]{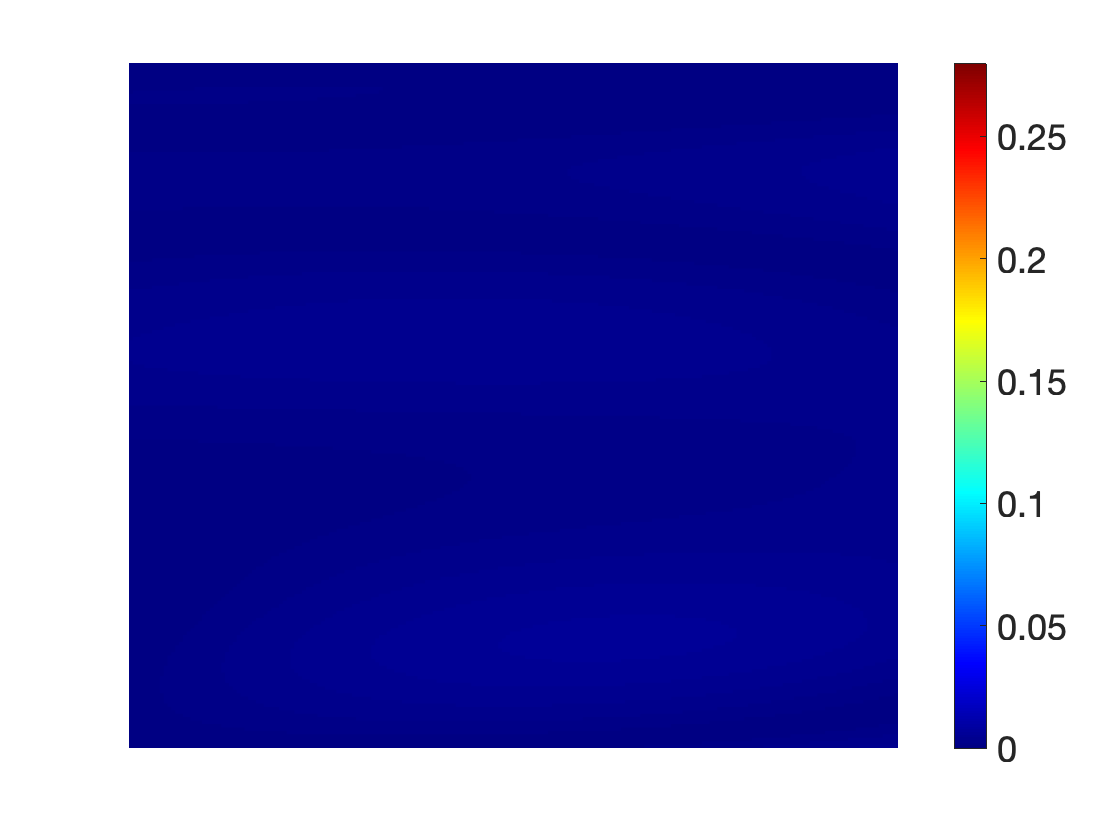} \\
\includegraphics[width=0.199\textwidth]{a22neu2dpex.png} &
\includegraphics[width=0.199\textwidth]{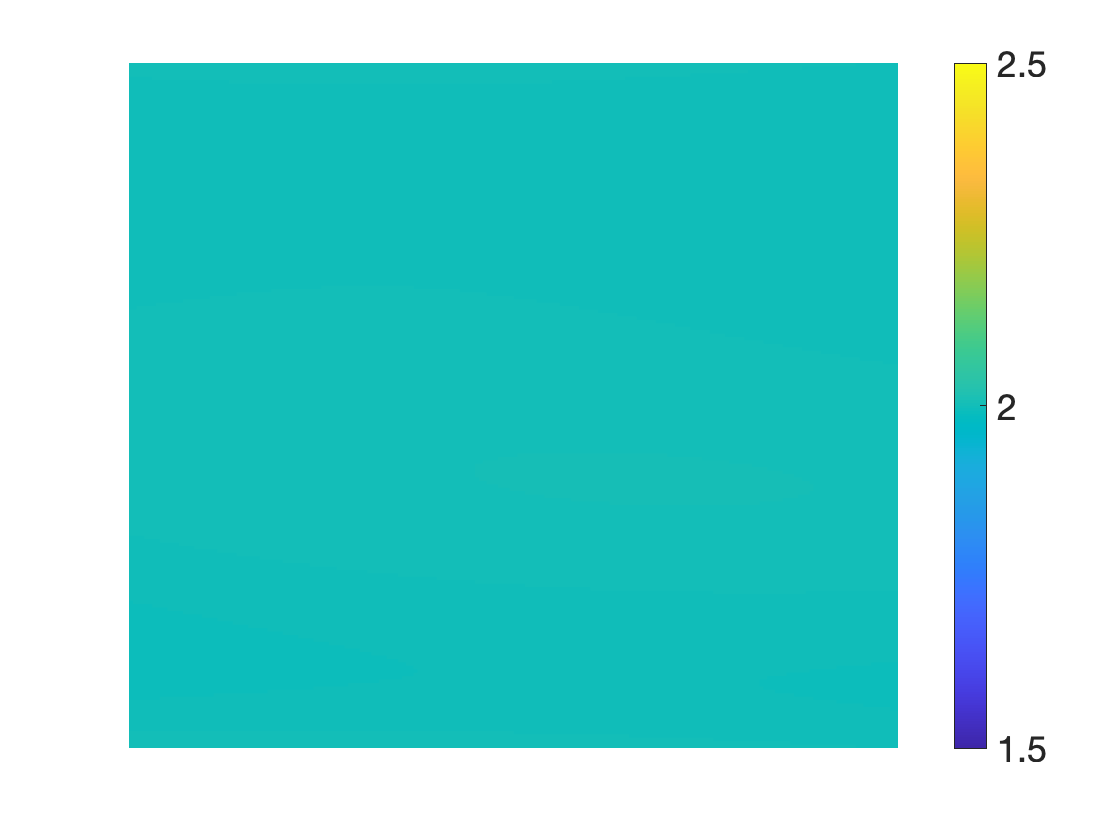} &
\includegraphics[width=0.199\textwidth]{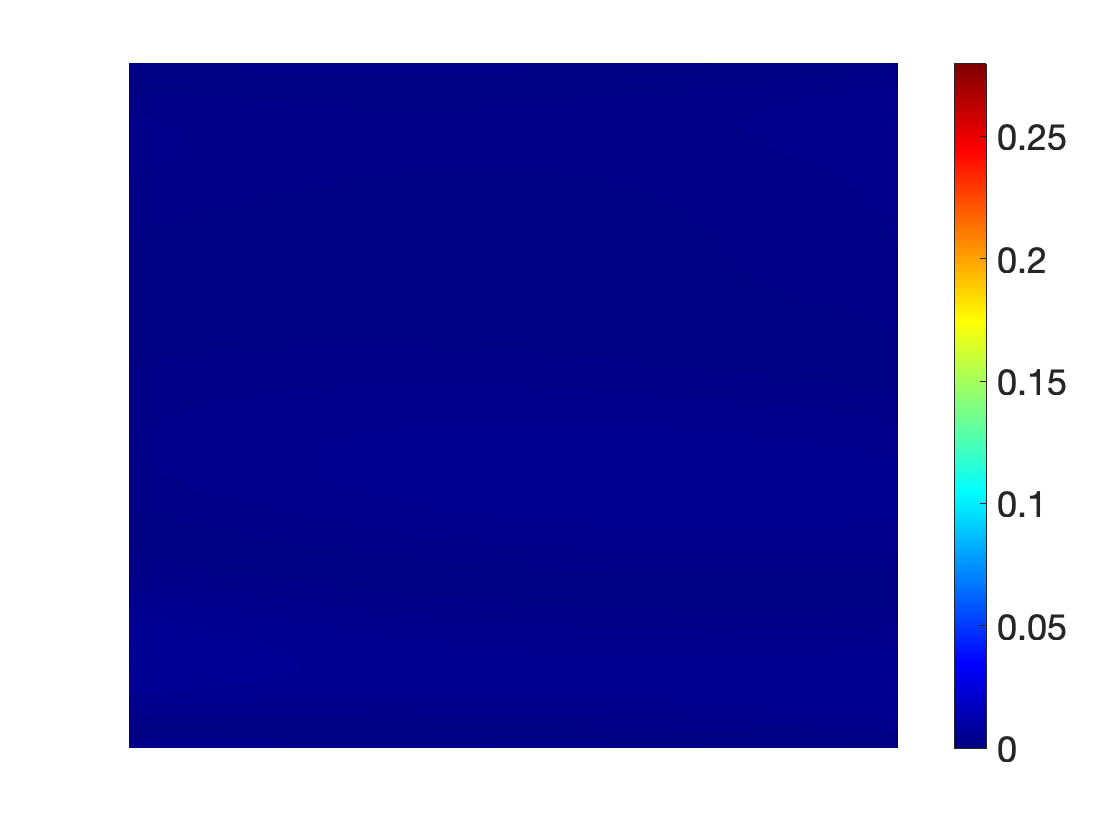} &
\includegraphics[width=0.199\textwidth]{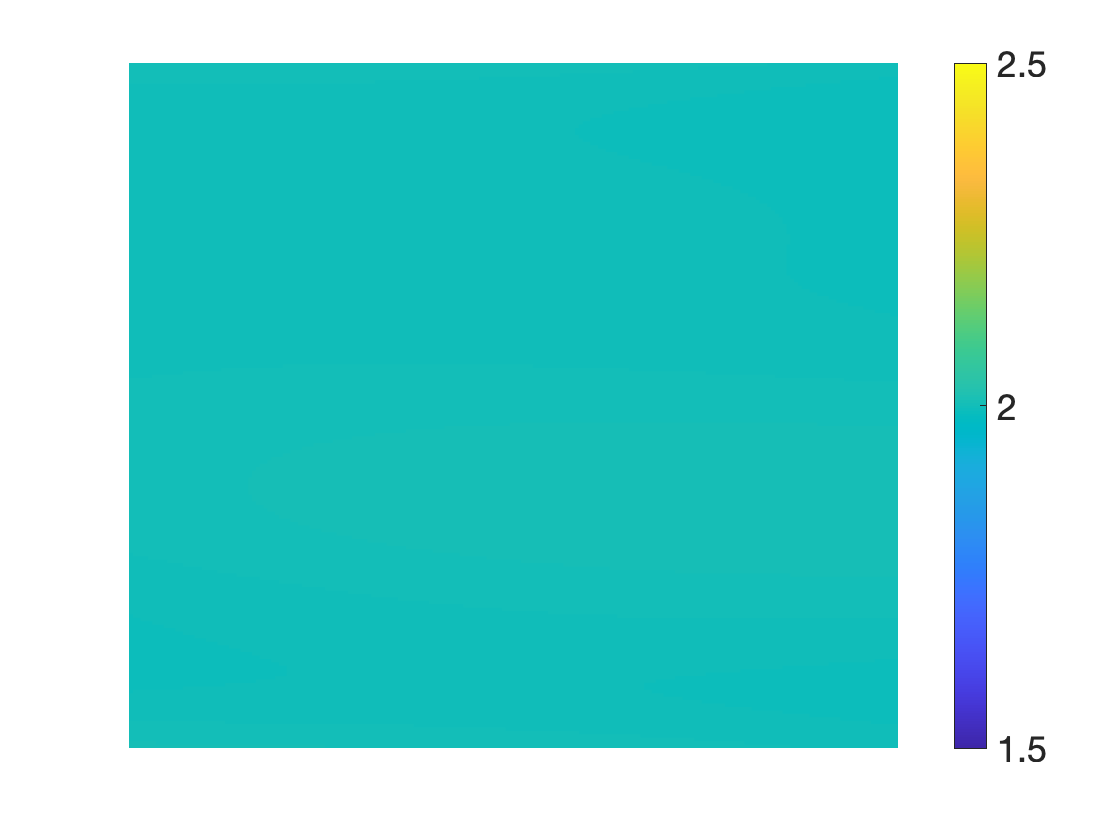} &
\includegraphics[width=0.199\textwidth]{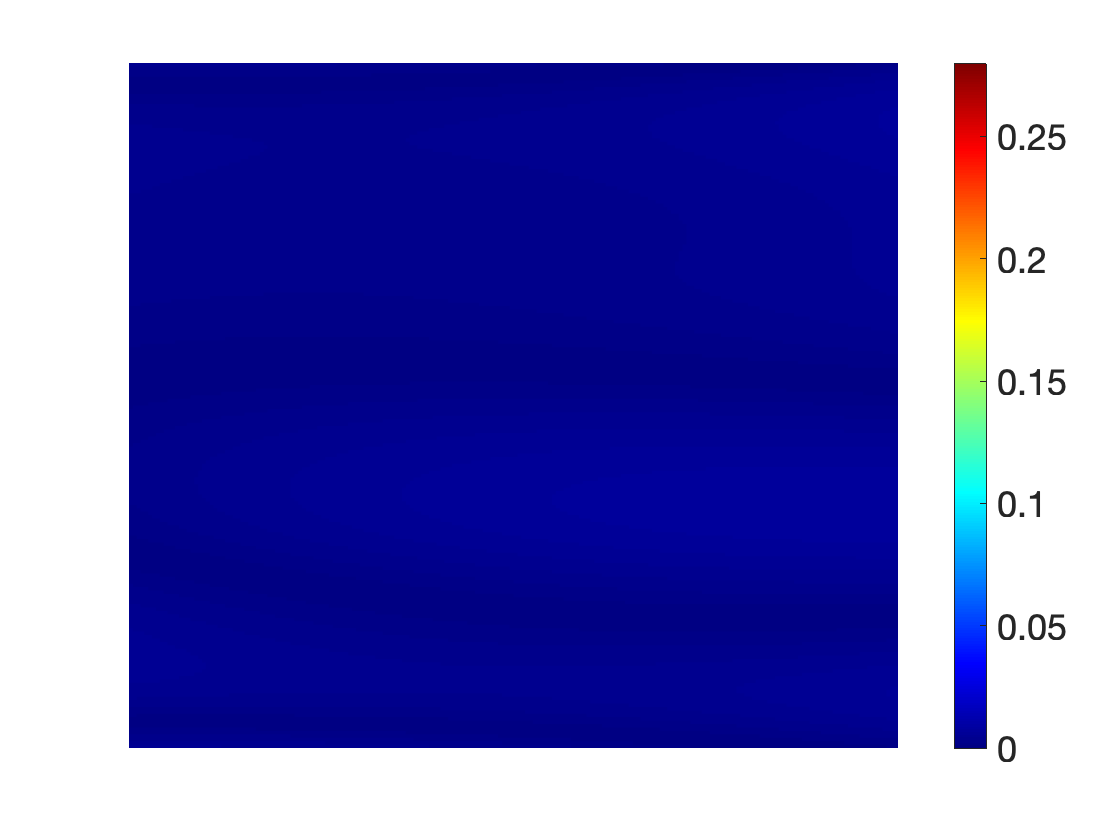} \\
(a) $A^\dag$  & (b) $\hat A$ & (c) $|\hat A-A^\dag|$ & (d) $\hat A$ & (e) $|\hat A-A^\dag|$
\end{tabular}
\caption{The reconstructions for Example \ref{exam:neu2d3} with exact data in (b) and noisy data $(\delta=10\%)$ in (d) using PINN. From the top to bottom, the results are for $A_{11}$, $A_{12}$ and $A_{22}$, respectively.}
\label{fig:neupinn}
\end{figure}

Additionally, we present also the results by the standard PINN for Example \ref{exam:neu2d3}. Specifically, we approximate both the conductivity tensor $A$ and the state $u$ separately using two DNNs, and minimize an empirical version of the following population loss
\begin{equation}\label{obj:neupinnloss}
    \begin{aligned}
            J_{\gamma}(\theta,\kappa)=&\sum_{i=1}^{N}\left(\| \nabla u_{i,\kappa} - \nabla z_i^{\delta} \|_{L^2(\omega)^d}^2+\gamma_{1}\| \nabla \cdot (A_{\theta} \nabla u_{i,\kappa}) + f_i \|_{L^2(\Omega)}^2\right.\\
            &\left.+\gamma_2 \|\vec{n}\cdot(A_{\theta} \nabla u_{i,\kappa} )-g_i\|_{L^2(\partial \Omega)}^2 +\gamma_3 \| A_{\theta}\|_{L^2(\Omega)^{d,d}}^2+\gamma_4 \|A_\theta-A^\dagger\|_{L^2(\partial \Omega)^{d,d}}^2 \right).
    \end{aligned}
\end{equation}
Note that the projection operator $P_{\mathcal{K}}$ is not included to enforce the box constraint of the admissible set $\mathcal{K}$, since the presence of $P_{\mathcal{K}}$ poses the challenge of computing the derivative with respect to the DNN parameters $(\theta,\kappa)$. Numerically, it is observed that compared to the formulation \eqref{eqn:loss_NN}, the additional knowledge $A^\dagger|_{\partial\Omega}$ in \eqref{obj:neupinnloss} is crucial: in the absence of the boundary penalization term $\|A_\theta-A^\dag\|^2_{L^2(\partial\Omega)^{d,d}}$, the PINN loss \eqref{obj:neupinnloss} fails to learn a reliable reconstruction for the conductivity tensor $A$. Fig. \ref{fig:neupinn} indicates that the reconstructions by PINN are largely comparable to that by MLS-DNN for exact data and slightly surpasses the MLS-DNNs for noisy data. This improvement is attributed to the \textit{a priori} knowledge $A^\dagger|_{\partial\Omega}$ used in the formulation. Also while the conductivity tensor $A$ is reconstructed accurately, the state $u$ is actually not very accurately resolved. Numerically, we also observe that introducing a boundary penalization term for $u$ (i.e., $\|u_{i,\kappa}-u^{\dagger}_i\|_{L^2(\partial \Omega)}^2$) into the loss \eqref{obj:neupinnloss} can successfully address this issue. However, we leave a thorough investigation of the mechanism of the PINN loss and its variants to future work, since PINN is not the primary focus of this study. Last, we briefly comment on the computing time of PINN versus MLS-DNNs. PINN takes longer per epoch to compute $J_{\gamma}(\theta,\kappa)$, primarily due to the storage and computation of second-order spatial derivatives in the term $\| \nabla \cdot (A_{\theta} \nabla u_{i,\kappa}) + f_i \|_{L^2(\Omega)}^2$.

\begin{figure}[htb!]
\centering
\setlength{\tabcolsep}{0em}
\begin{tabular}{ccccc}
\includegraphics[width=0.199\textwidth]{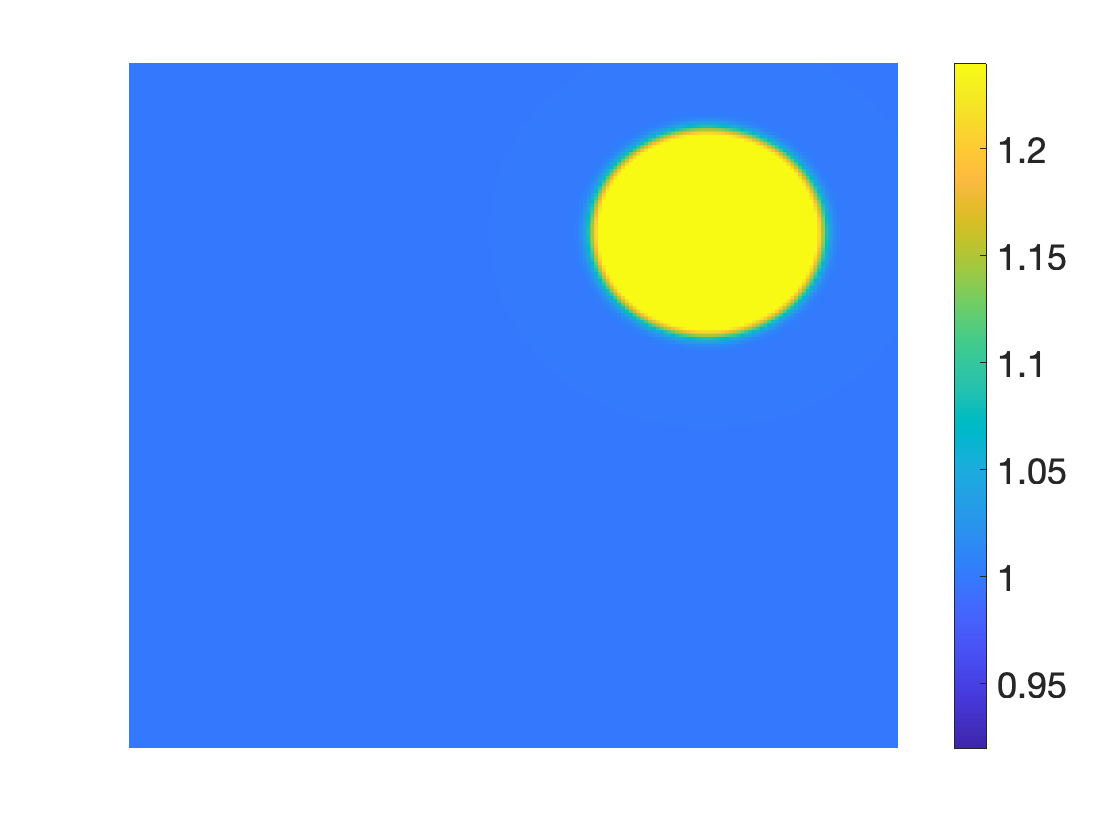} &
\includegraphics[width=0.199\textwidth]{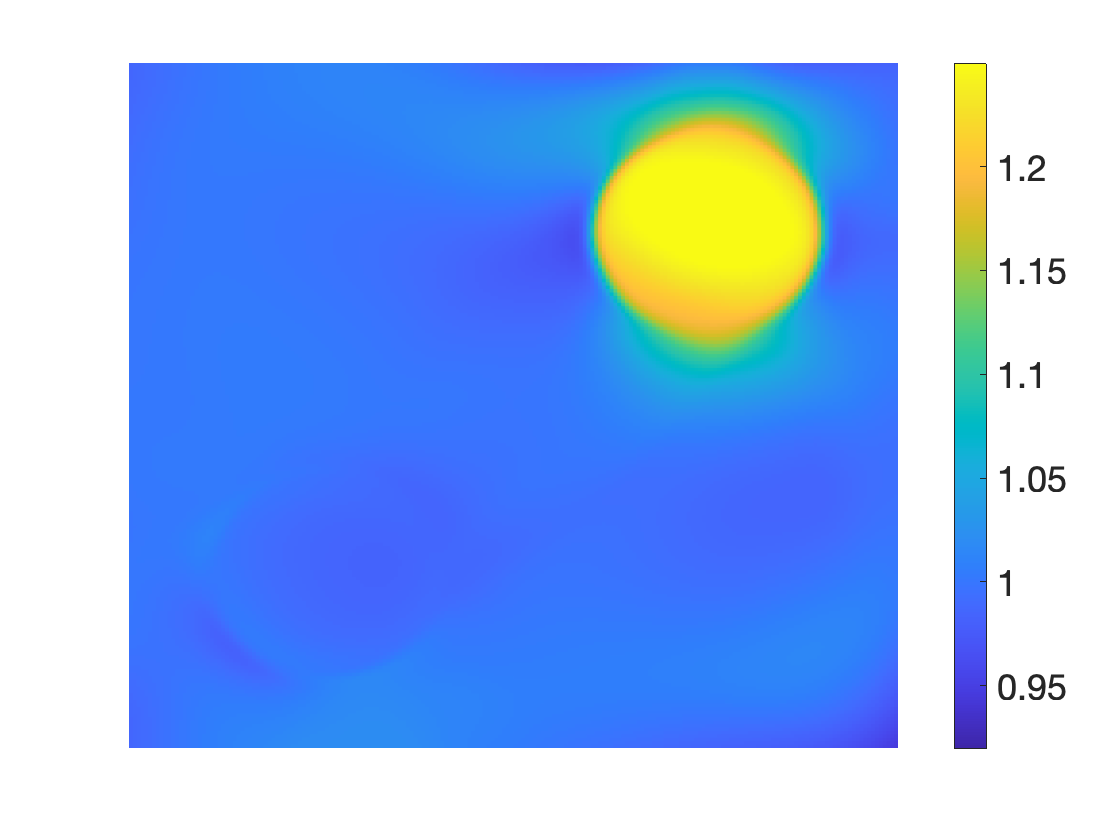} &
\includegraphics[width=0.199\textwidth]{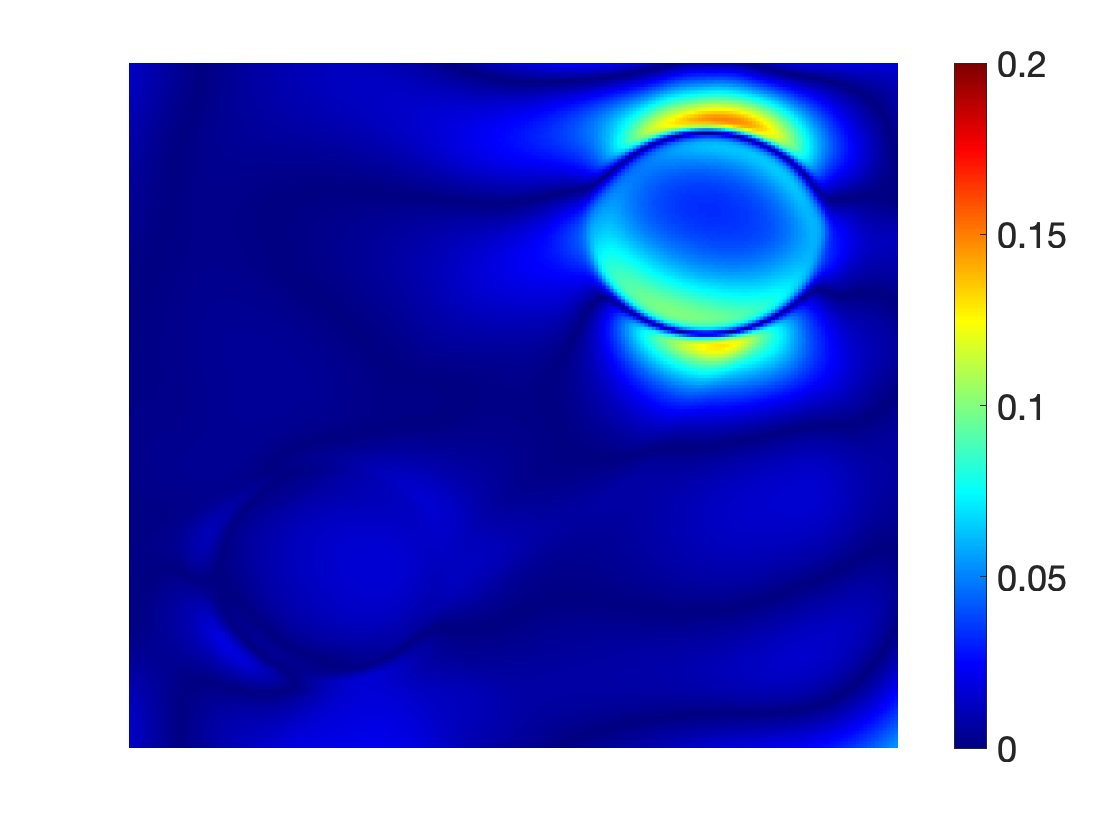} &
\includegraphics[width=0.199\textwidth]{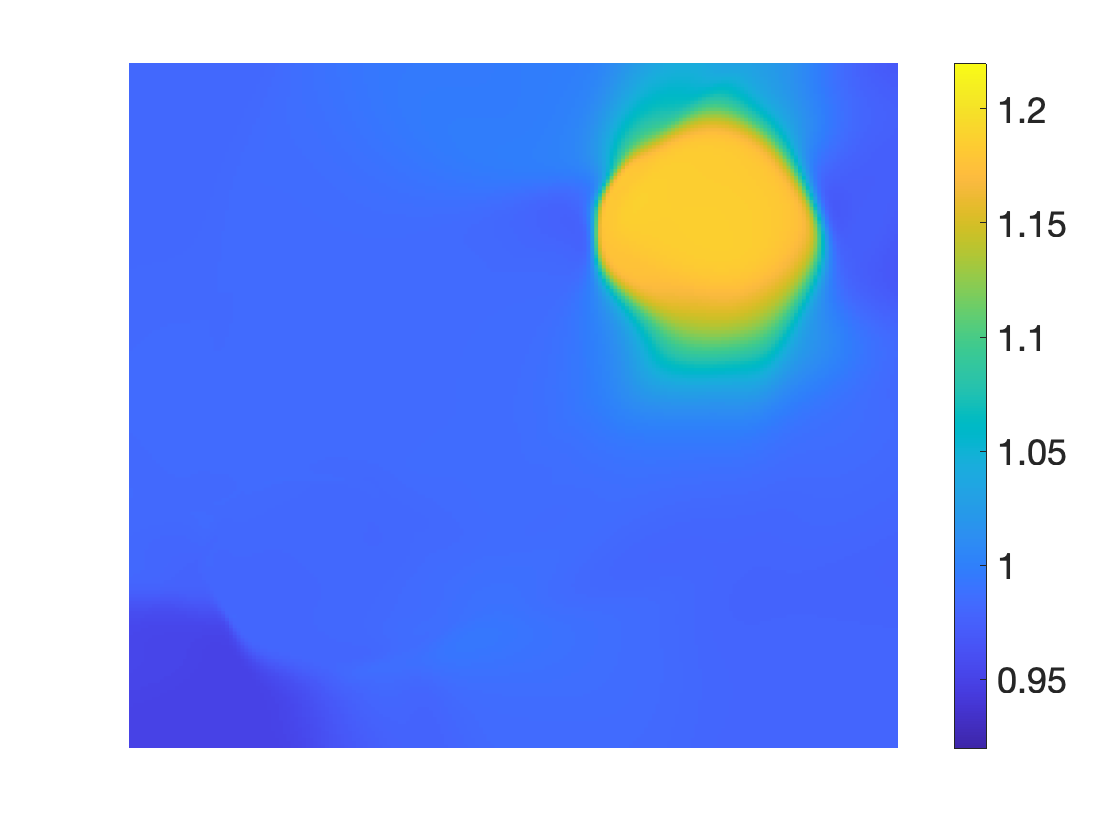} &
\includegraphics[width=0.199\textwidth]{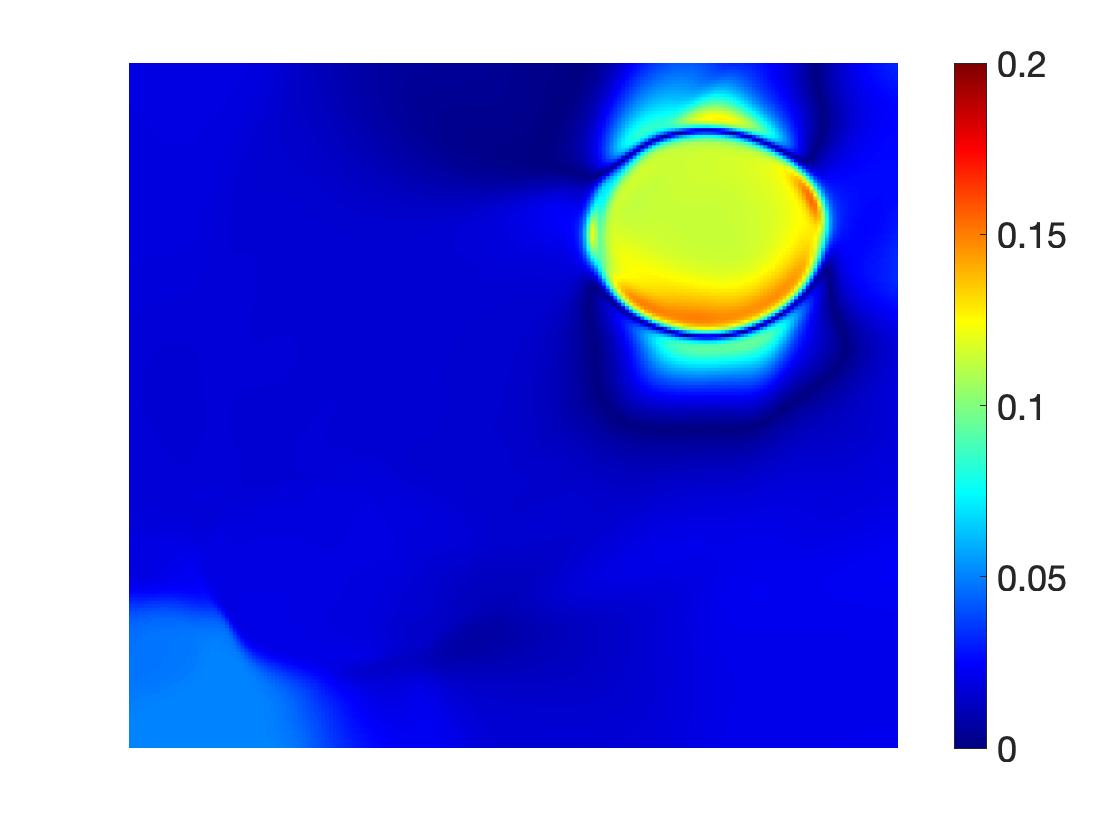} \\
\includegraphics[width=0.199\textwidth]{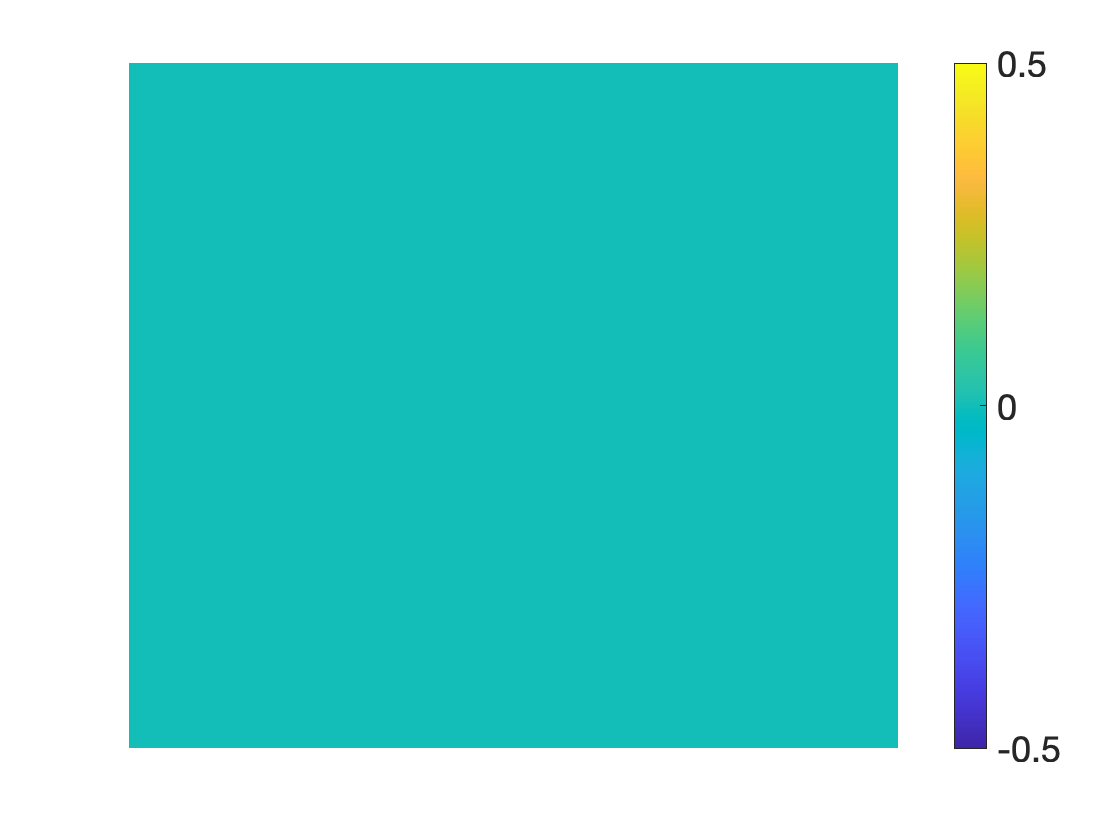} &
\includegraphics[width=0.199\textwidth]{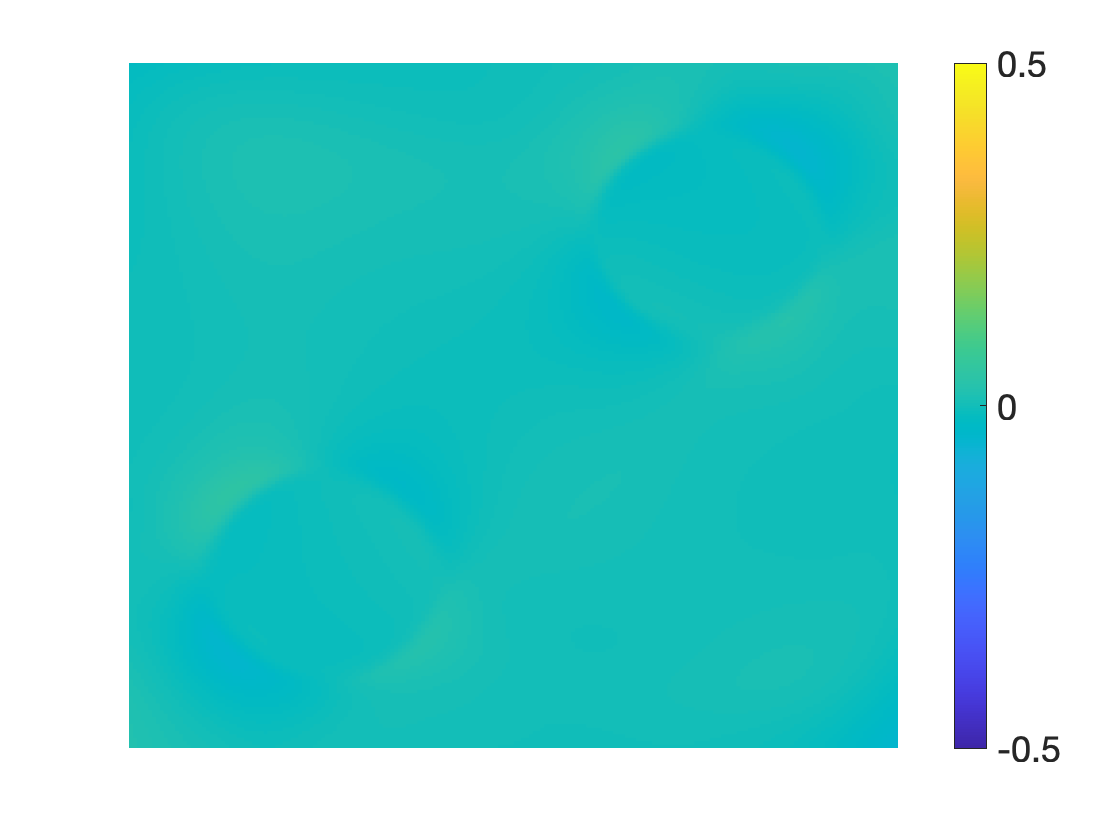} &
\includegraphics[width=0.199\textwidth]{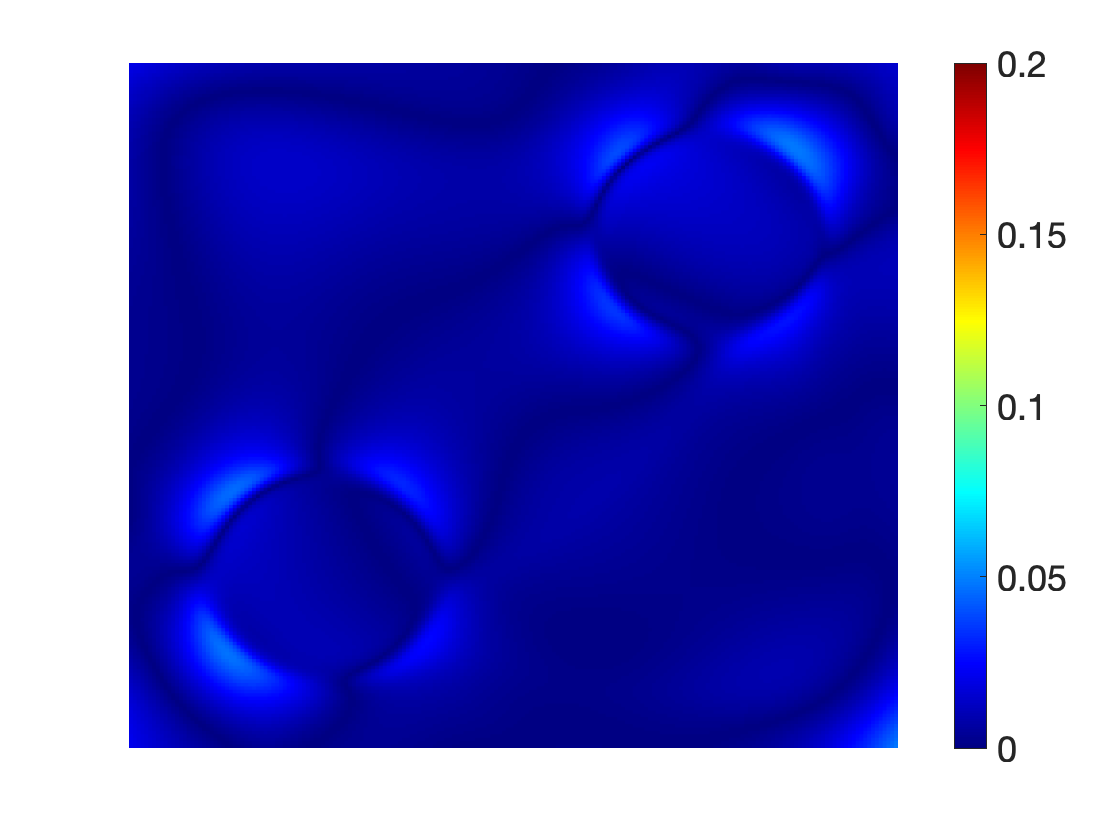} &
\includegraphics[width=0.199\textwidth]{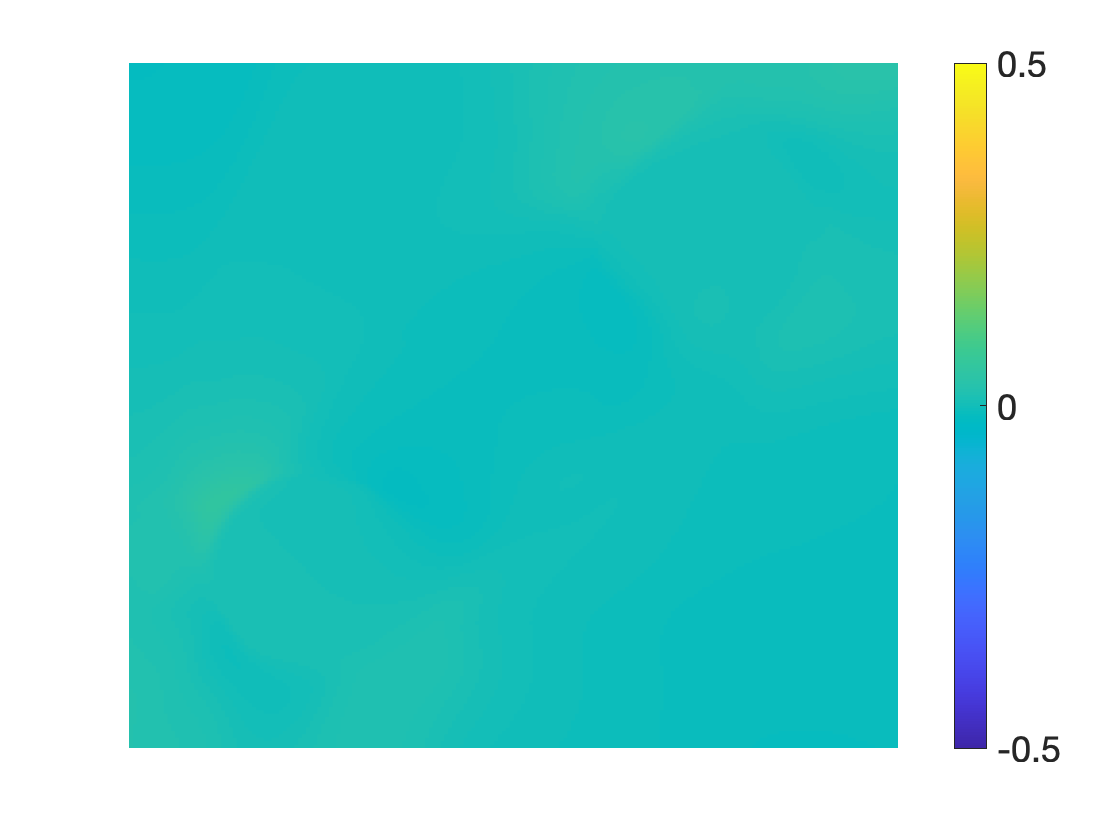} &
\includegraphics[width=0.199\textwidth]{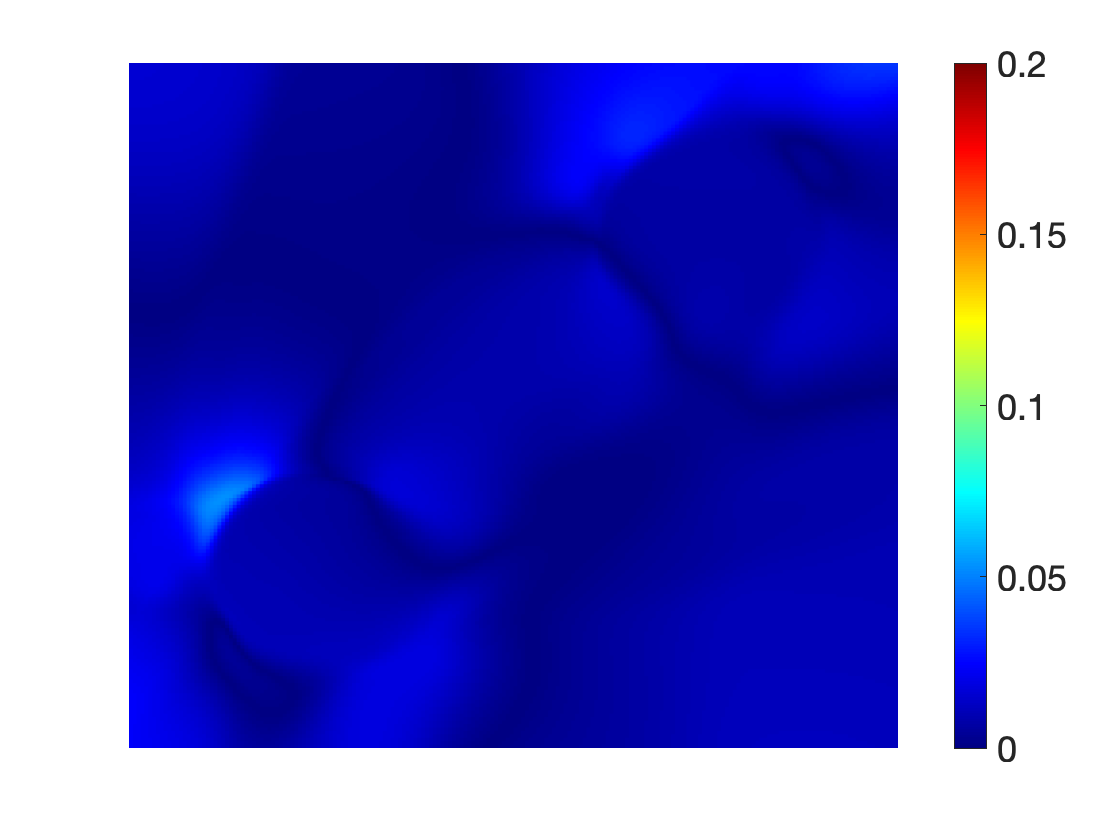} \\
\includegraphics[width=0.199\textwidth]{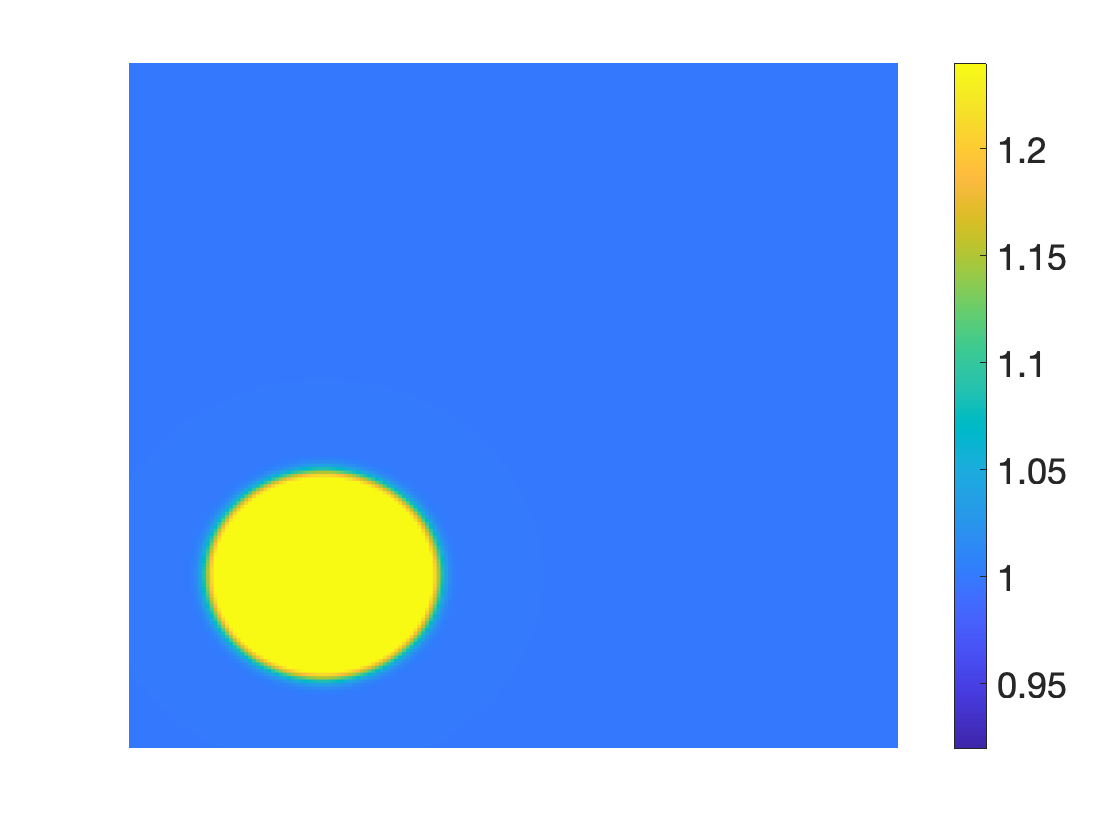} &
\includegraphics[width=0.199\textwidth]{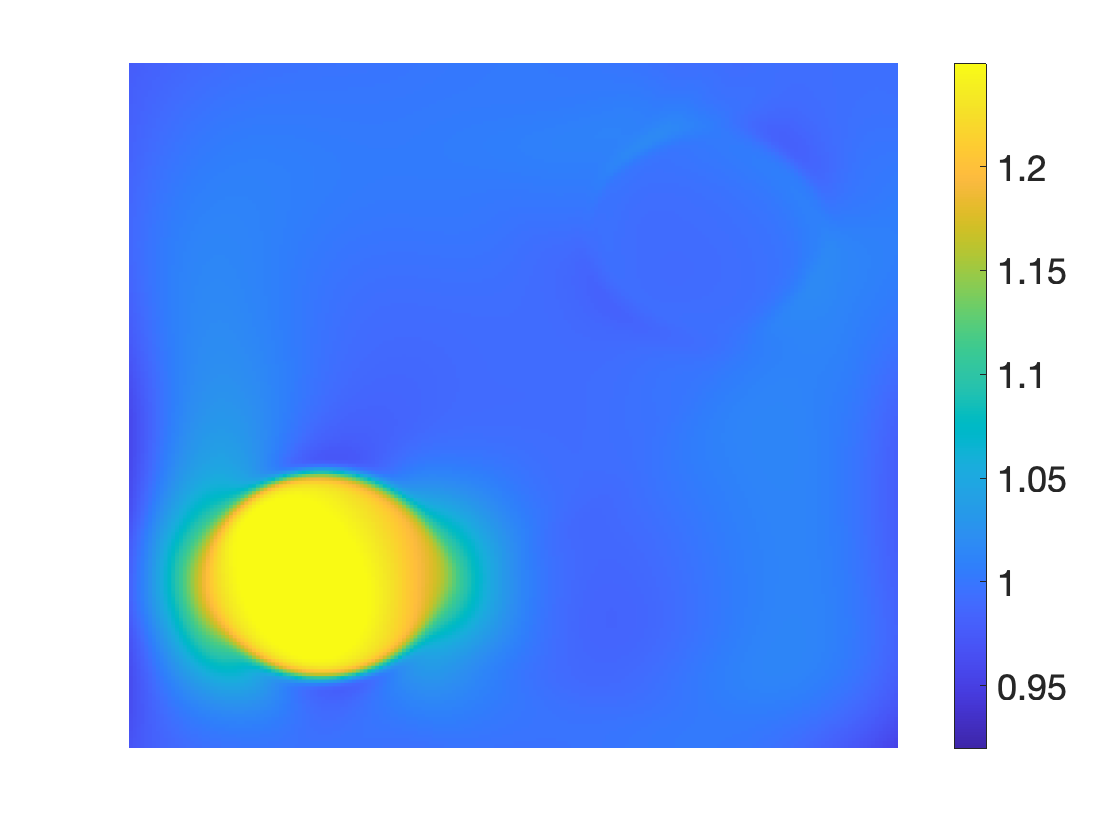} &
\includegraphics[width=0.199\textwidth]{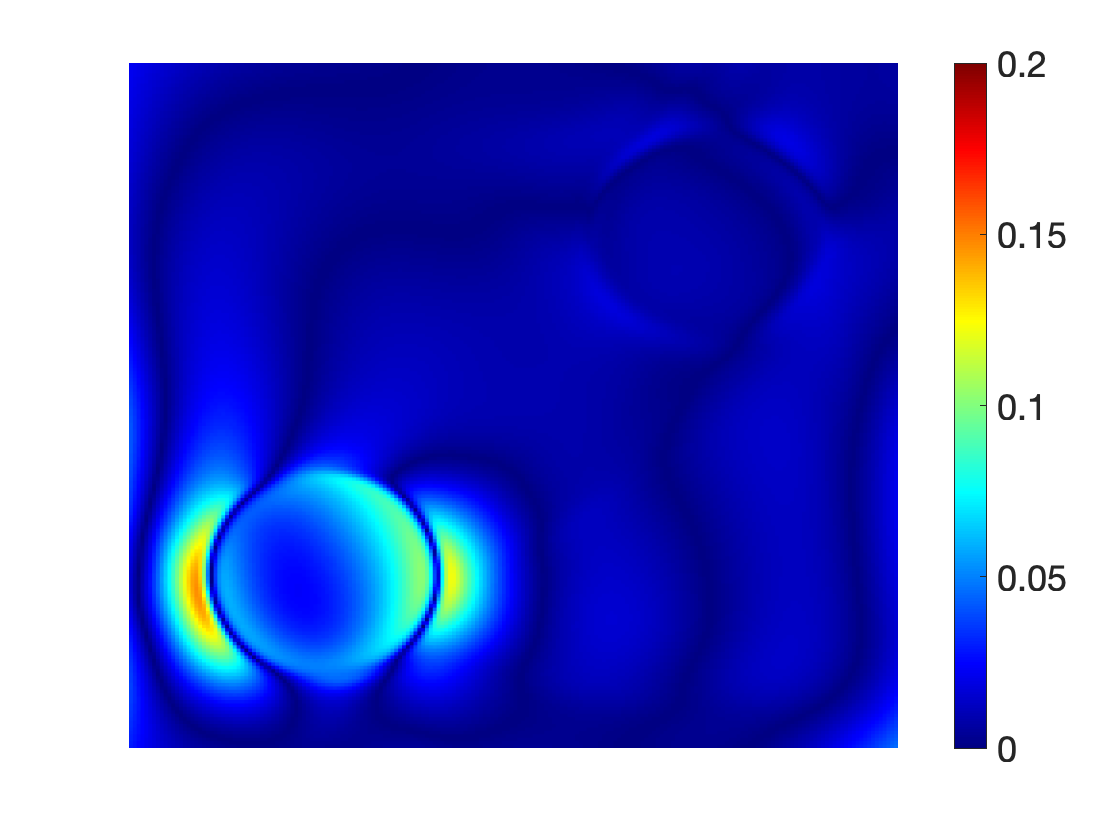} &
\includegraphics[width=0.199\textwidth]{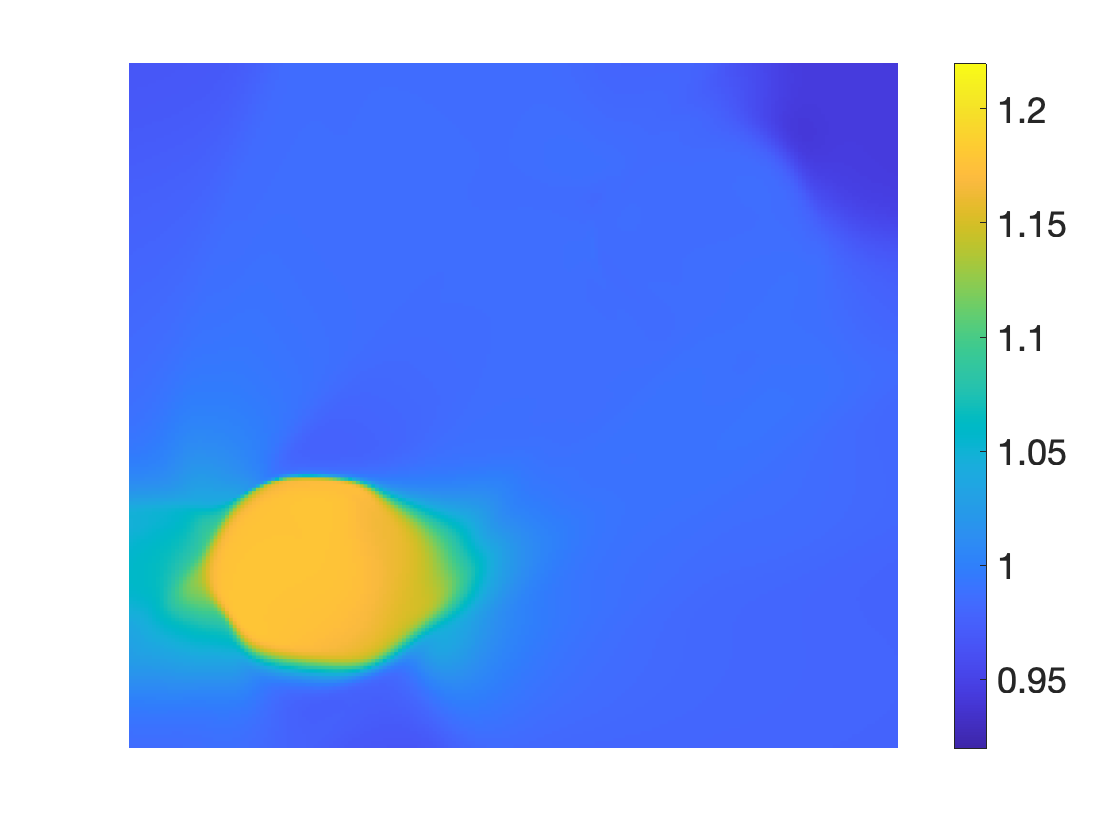} &
\includegraphics[width=0.199\textwidth]{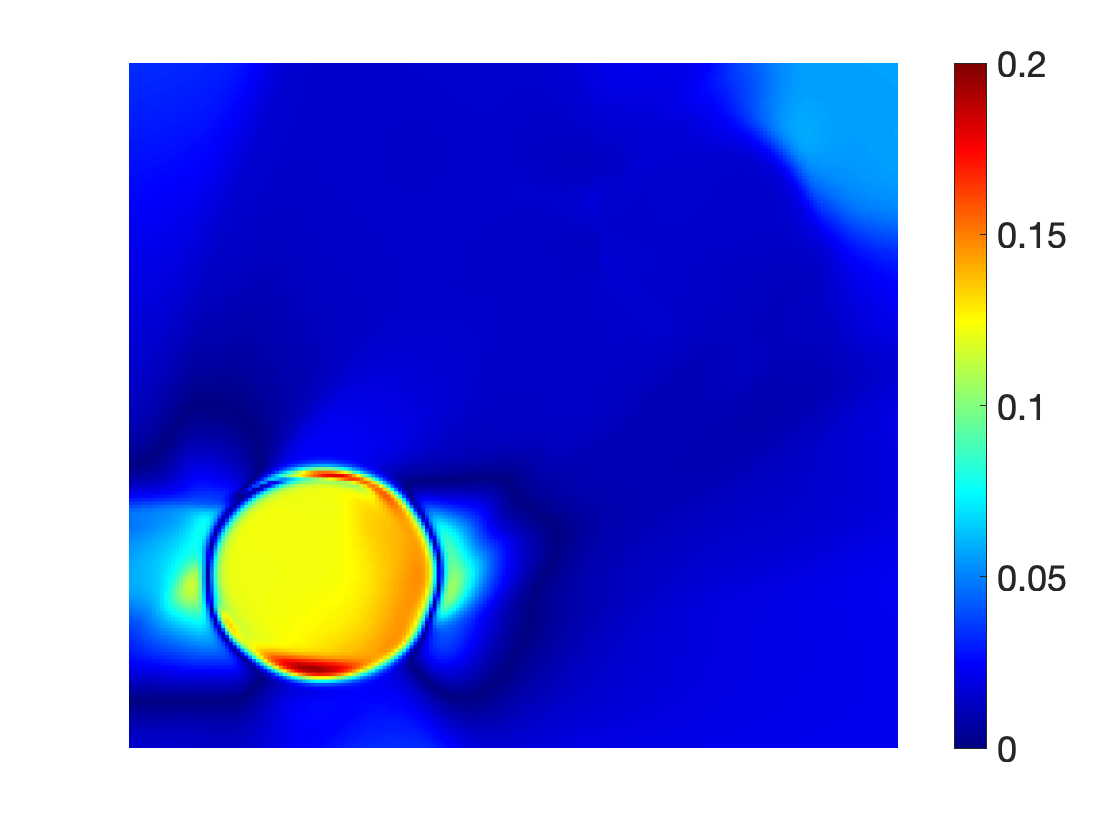} \\
(a) $A^\dag$  & (b) $\hat A$ & (c) $|\hat A-A^\dag|$ & (d) $\hat A$ & (e) $|\hat A-A^\dag|$
\end{tabular}
\caption{The reconstructions for Example \ref{exam:neu2d4} with exact data in (b) and noisy data $(\delta=5\%)$ in (d). From the top to bottom, the results are for $A_{11}$, $A_{12}$ and $A_{22}$, respectively.}
\label{fig:neu2d4}
\end{figure}

The fourth example is about recovering a diagonal conductivity matrix $A^\dag$ with nearly piecewise constant component functions.
\begin{example}
The domain $\Omega = (0,1)^2$, $A^\dag= \begin{pmatrix}
    \zeta_1(x_1,x_2) & 0\\
   0 & \zeta_2(x_1,x_2)\\
\end{pmatrix},$ where $\zeta_1(x_1,x_2)=1+0.3/(1+\exp(600((x_1-0.75)^2+(x_2-0.75)^2-0.0225)))$ and $\zeta_2(x_1,x_2)= 1+0.3/(1+\exp(600((x_1-0.25)^2+(x_2-0.25)^2-0.0225)))$. $f_{i=1,2,3,4}\equiv0$ and $g_1=x_1-0.5$, $g_2=x_2-0.5$, $g_3=\cos(\pi x_1)$ and $g_4=\cos(\pi x_2)$.
\label{exam:neu2d4}
\end{example}
Since the entries of $A^\dag$
are nearly piecewise constant, we also include the  total variation penalty \cite{RudinOsherFatemi:1992},
i.e., $\gamma_{tv}\sum_{i,j=1}^2| A_{ij}|_{\rm TV}$, to the loss $J_{\bsgamma}(\theta,\kappa)$, in order to
promote piecewise constancy in the reconstruction. Fig. \ref{fig:neu2d4} presents the reconstruction results for both exact and noisy data ($\delta = 5\%$), where $\gamma_{tv}=\ $1e-4 and 3e-3, respectively. For the exact data, the reconstruction accurately captures the location and shape of the circular bumps. When 5\% noise is present, there is a slight underestimation of the bump peaks and a mild blurring effect in the region near the interface.

The fifth example is about recovering a 3D conductivity matrix with mixed oscillatory and polynomial entries.
\begin{example}\label{exam:neu3d1}
    The domain $\Omega=(0,1)^3$, $A^\dagger = \begin{pmatrix}
    3& 1+\frac{\sin(4\pi x_3)}{2}& 1+\frac{\sin(4\pi x_2)}{2}\\
    1+\frac{\sin(4\pi x_3)}{2}& 2& 1+\frac{\sin(4\pi x_1)}{2}\\
    1+\frac{\sin(4\pi x_2)}{2}& 1+\frac{\sin(4\pi x_1)}{2}& 2+\frac{\sin(2\pi x_1)\sin(2\pi x_2)}{2}\\
    \end{pmatrix}$, $u_1^\dag=x_1+x_2+x_3+\frac{1}{3}(x_1^3+x_2^3+x_3^3)$, $u_2^\dag=x_1-x_2+x_3+\frac{1}{3}(x_1^3-x_2^3+x_3^3)$, $u_3^\dag=x_1+x_2-x_3+\frac{1}{3}(x_1^3+x_2^3-x_3^3)$, $u_4^\dag=-x_1+x_2+x_3+\frac{1}{3}(-x_1^3+x_2^3+x_3^3)$, $u_5^\dagger=-u_3^\dagger$ and $u_6^\dagger=-u_2^\dagger$.
\end{example}

\begin{figure}[htb!]
\centering
\setlength{\tabcolsep}{0em}
\begin{tabular}{ccccc}
\includegraphics[width=0.199\textwidth]{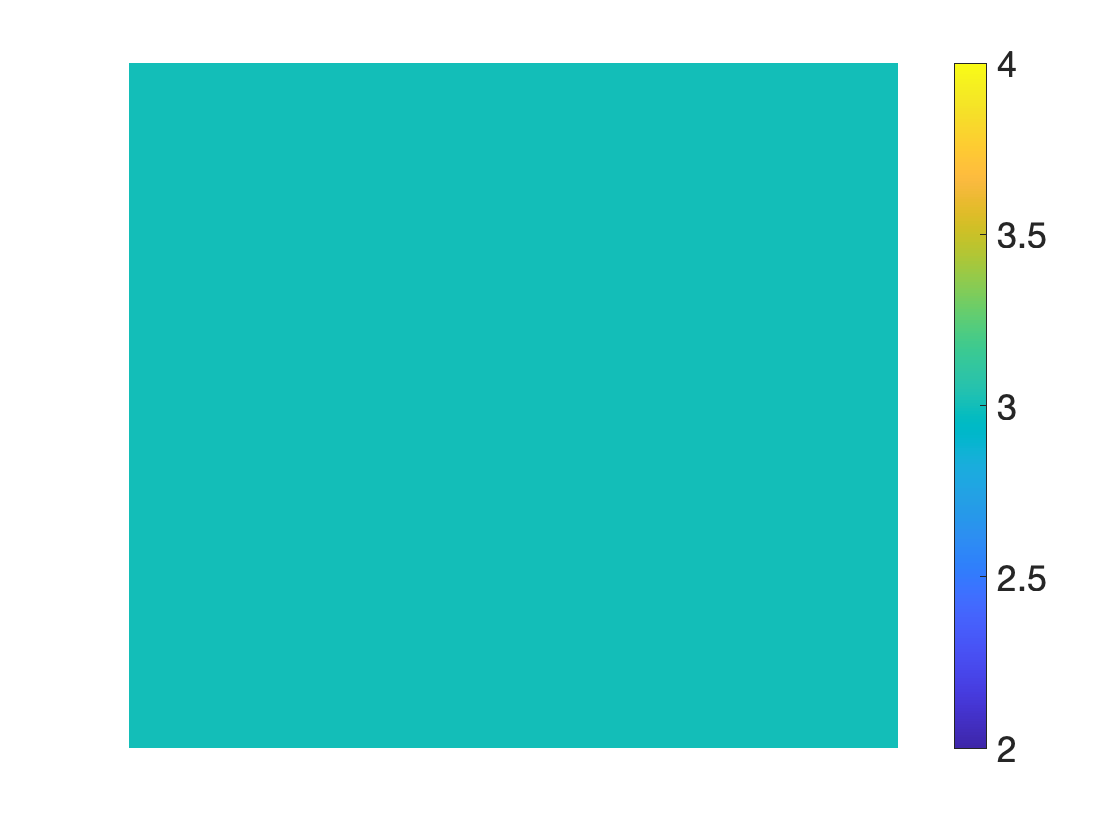} &
\includegraphics[width=0.199\textwidth]{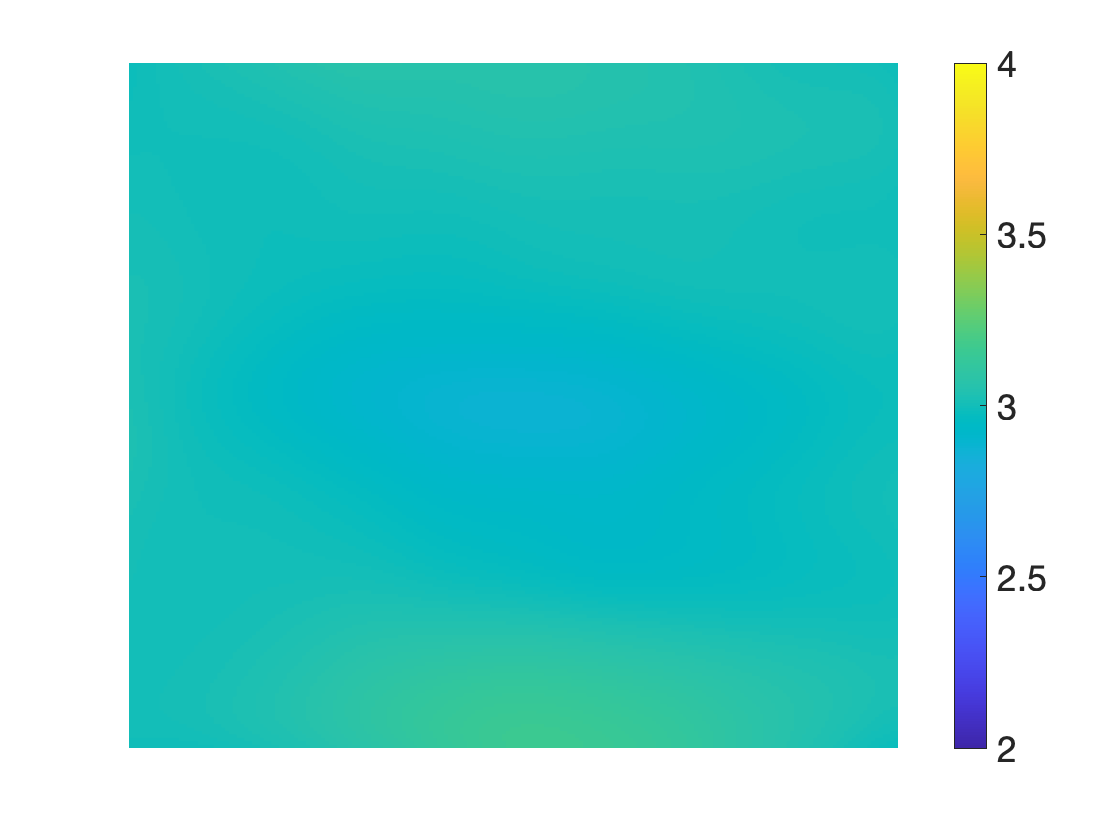} &
\includegraphics[width=0.199\textwidth]{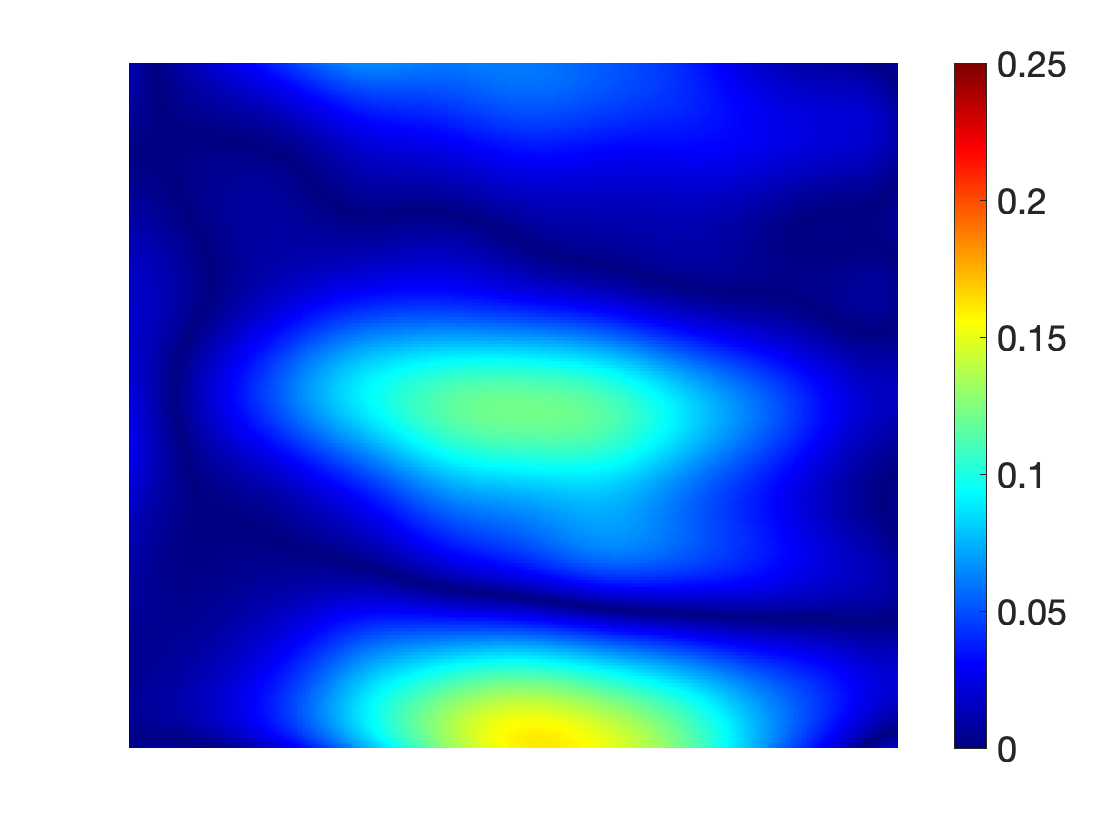} &
\includegraphics[width=0.199\textwidth]{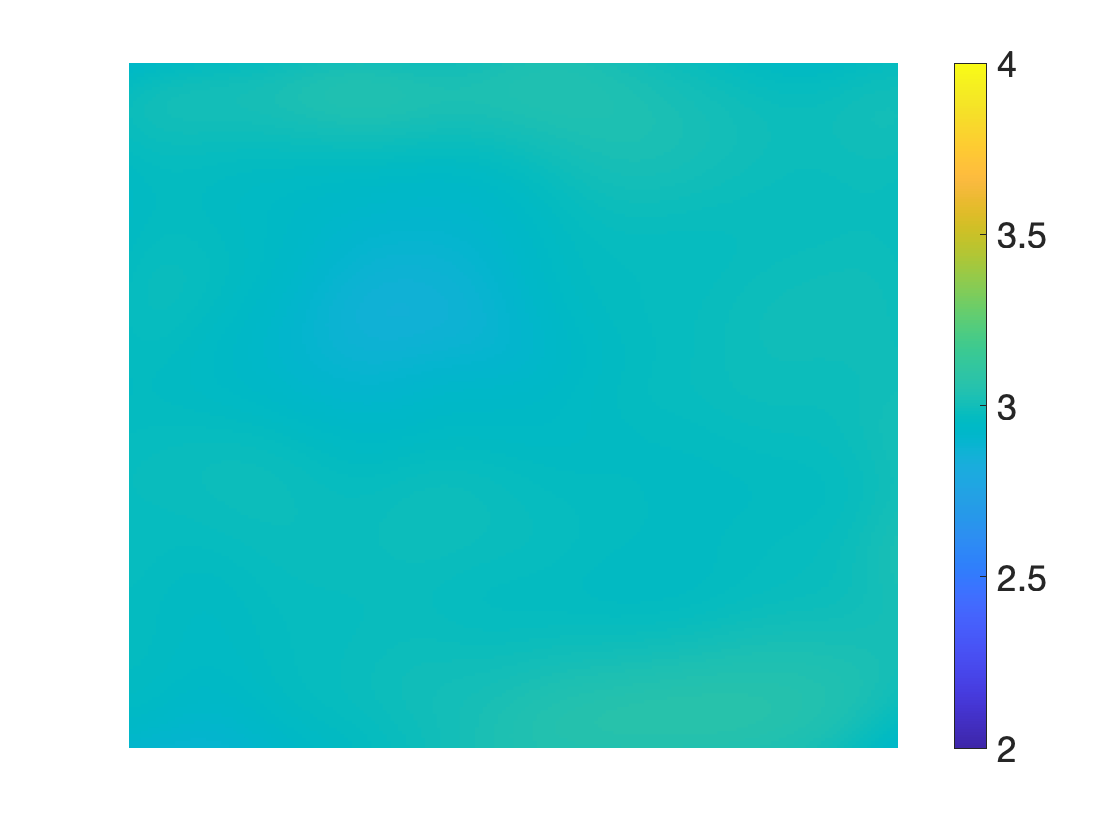} &
\includegraphics[width=0.199\textwidth]{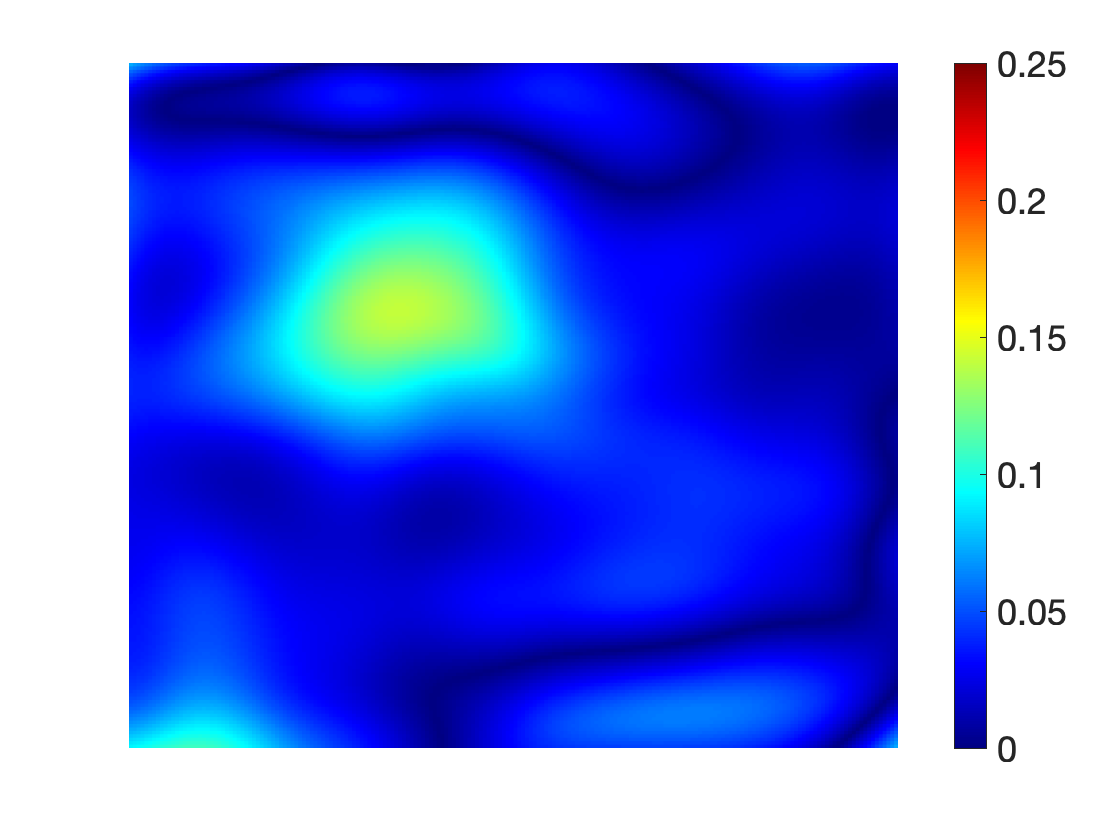}\\
\includegraphics[width=0.199\textwidth]{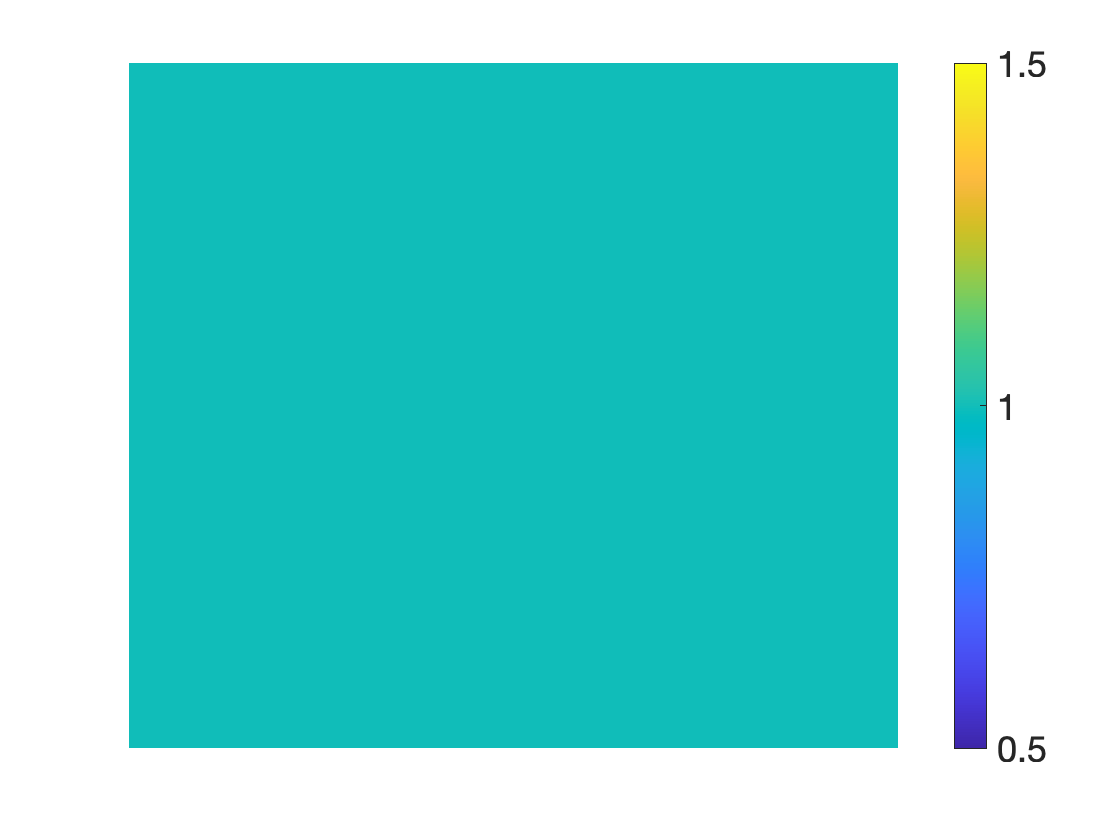} &
\includegraphics[width=0.199\textwidth]{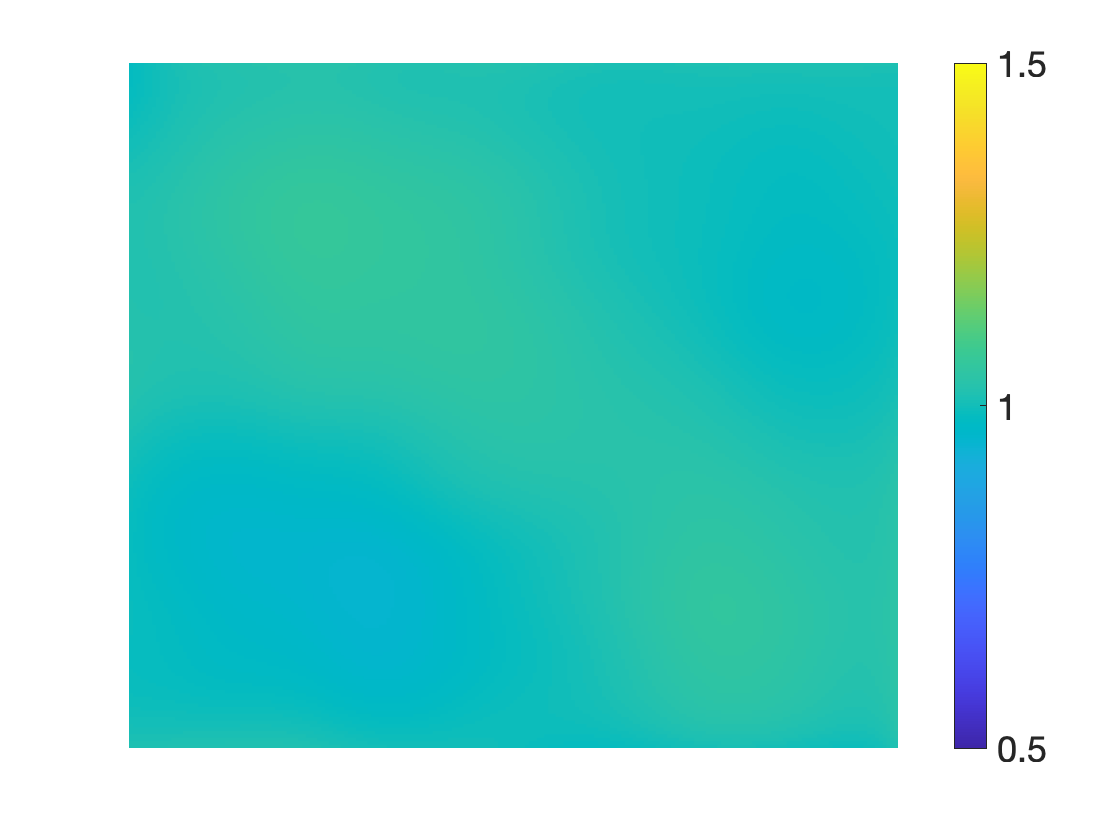} &
\includegraphics[width=0.199\textwidth]{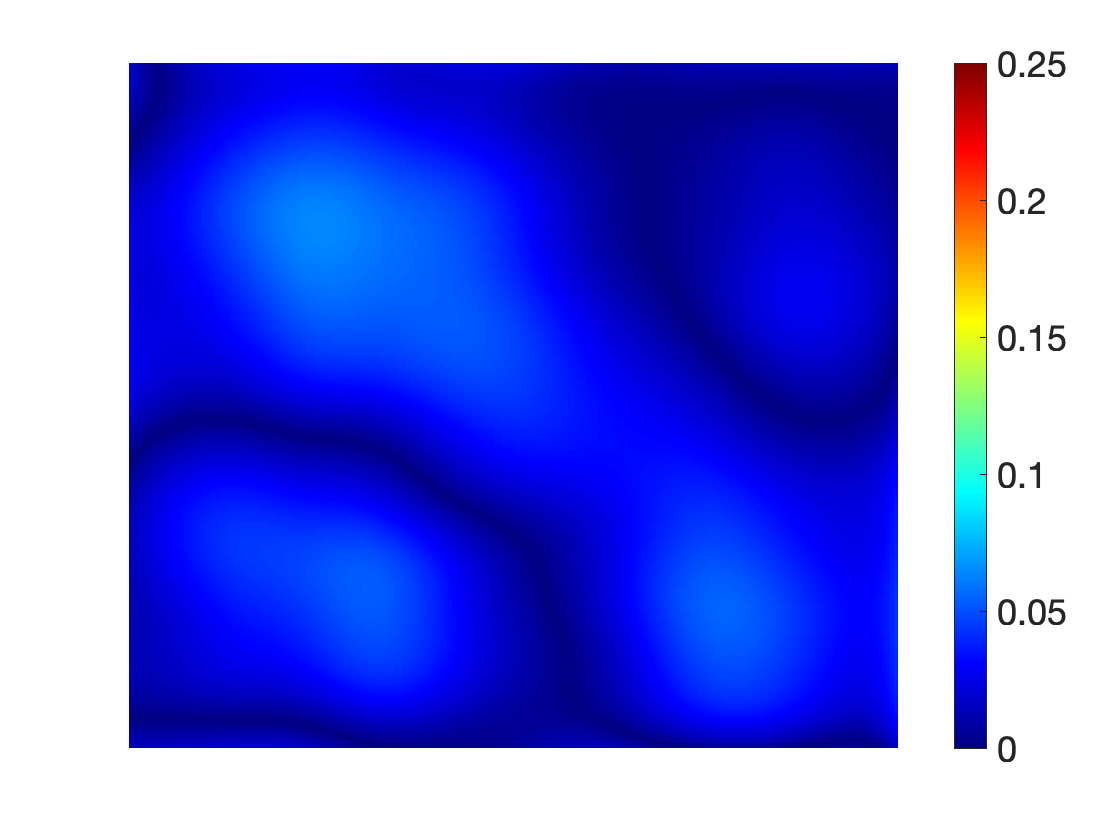}&
\includegraphics[width=0.199\textwidth]{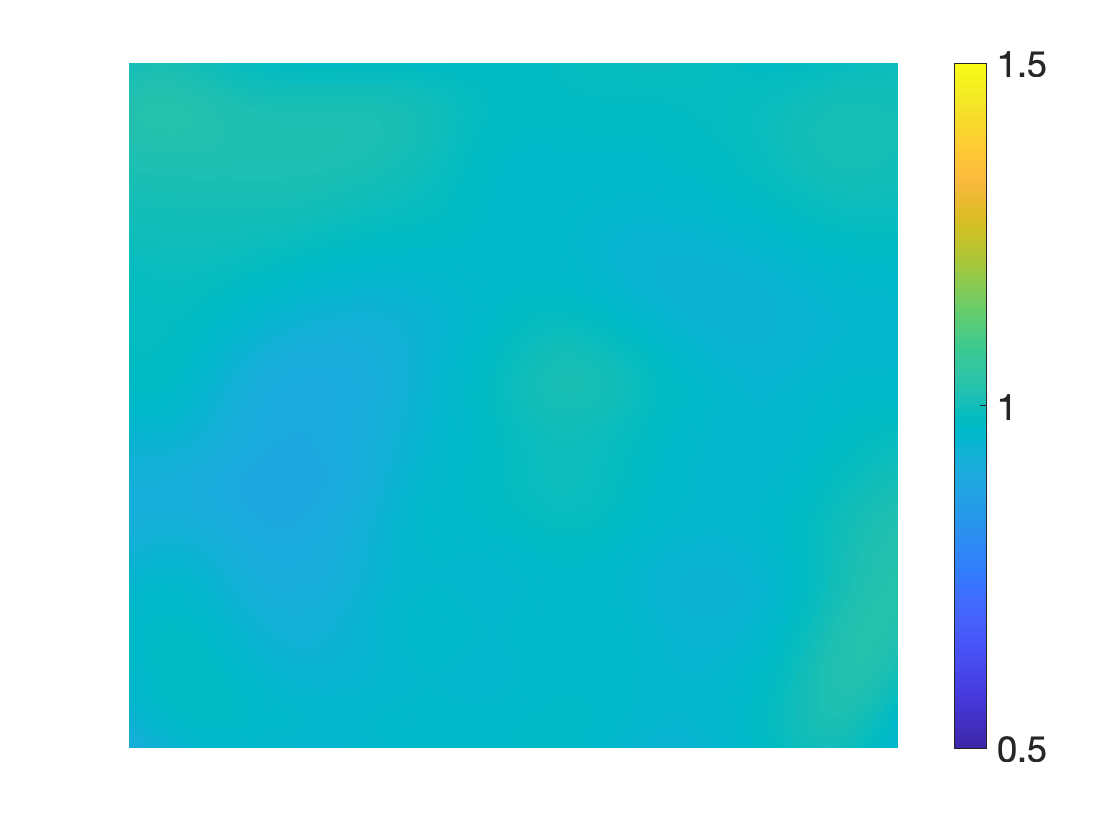} &
\includegraphics[width=0.199\textwidth]{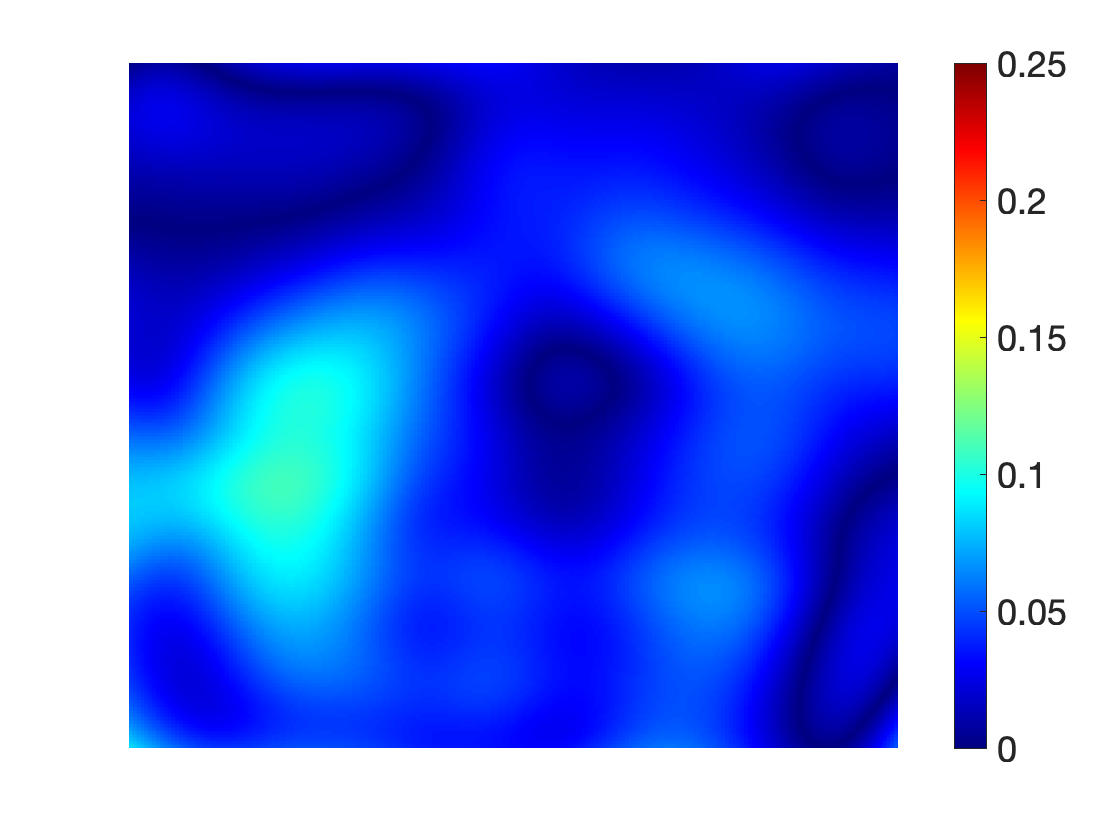}\\
\includegraphics[width=0.199\textwidth]{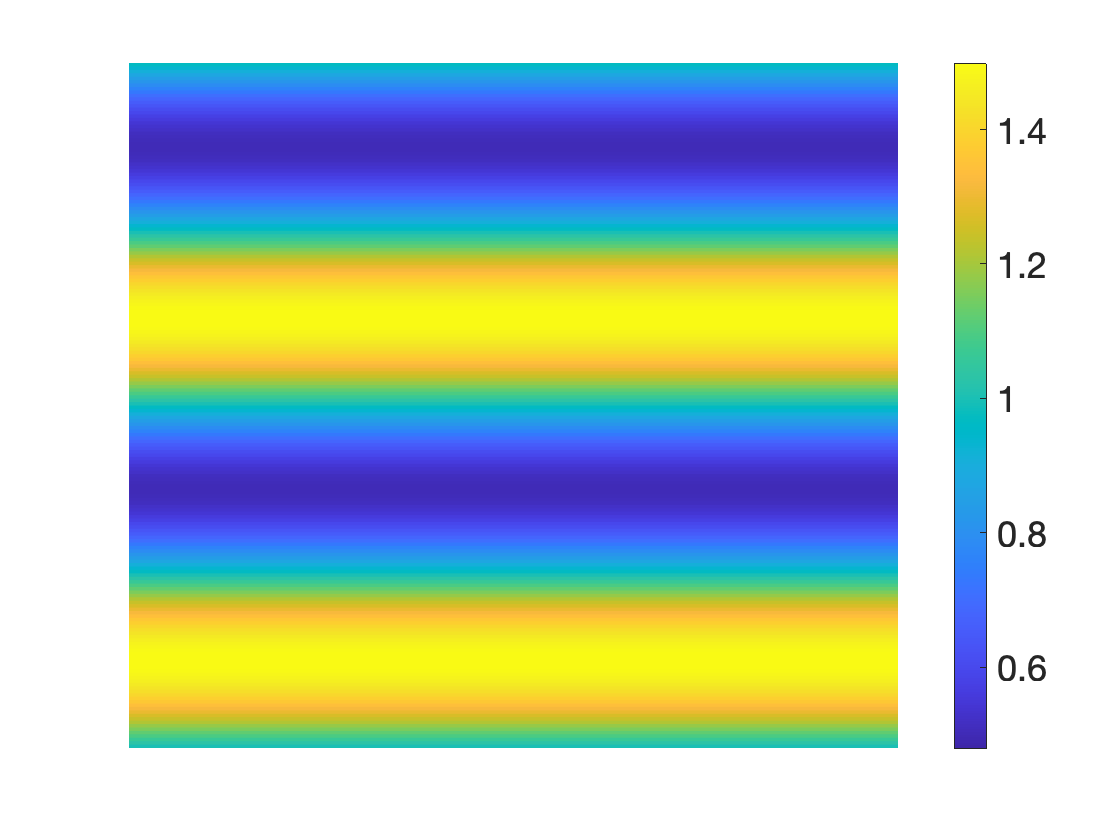} &
\includegraphics[width=0.199\textwidth]{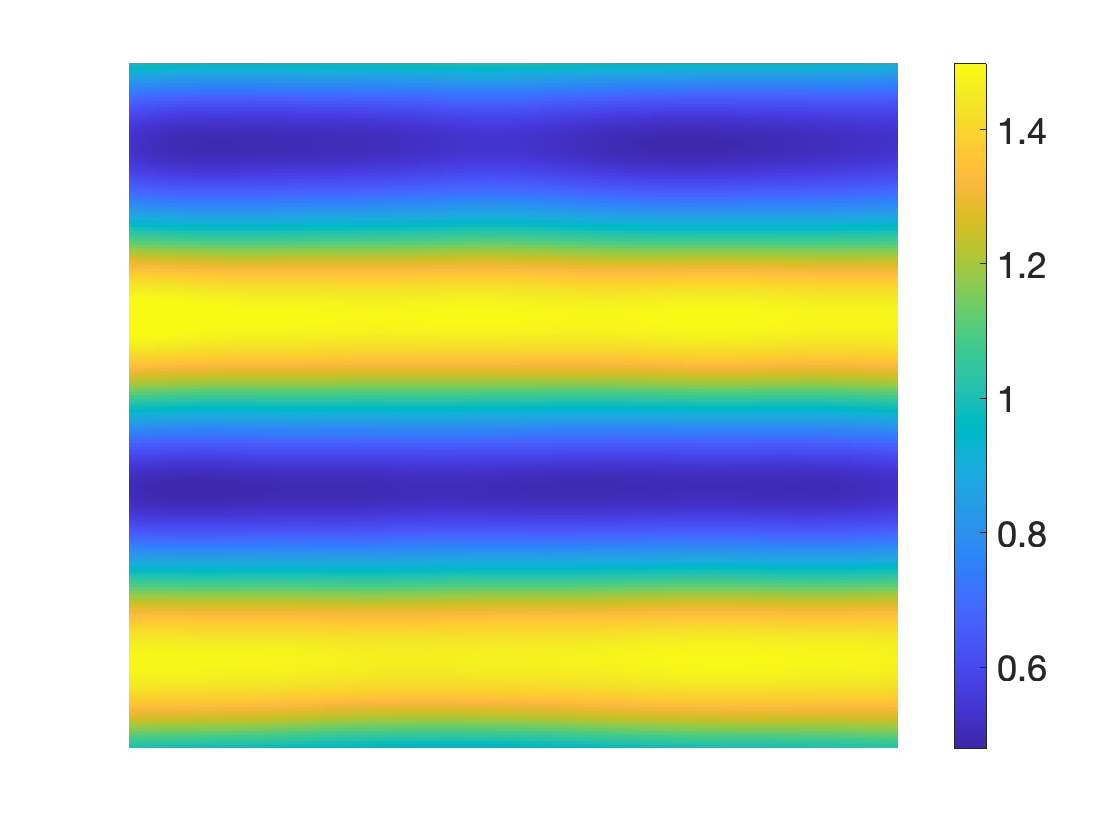} &
\includegraphics[width=0.199\textwidth]{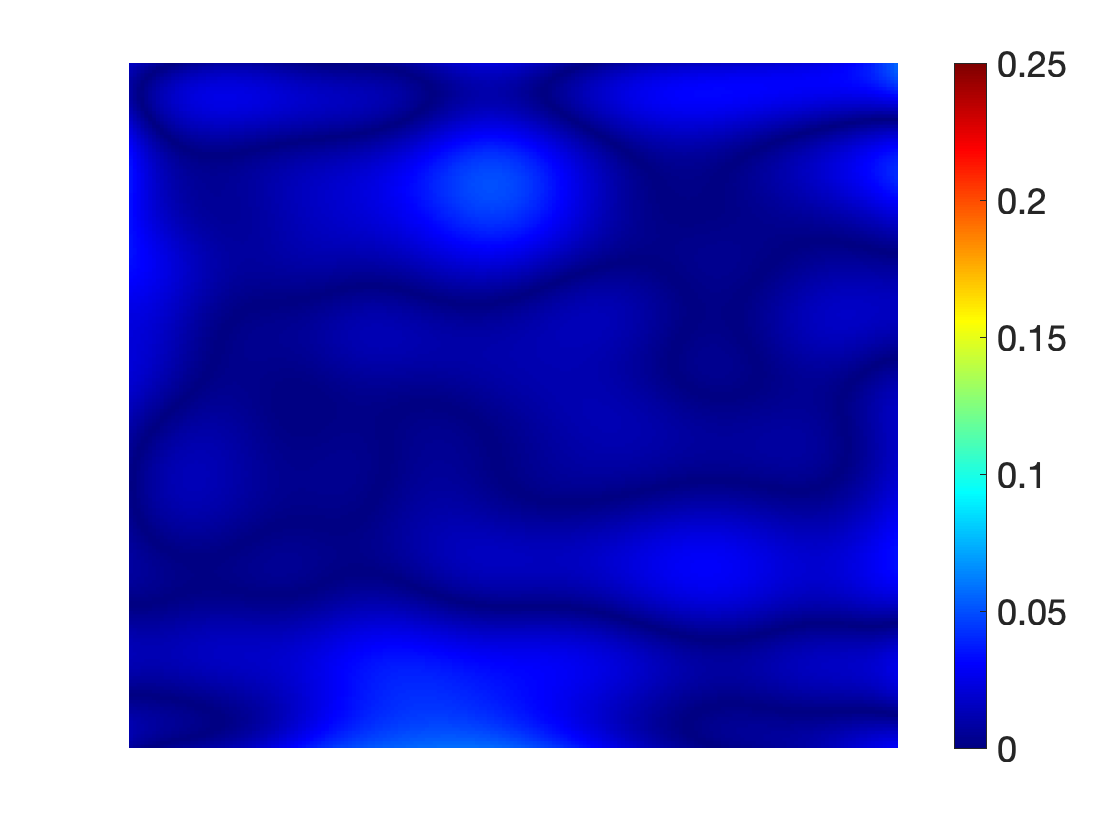} &
\includegraphics[width=0.199\textwidth]{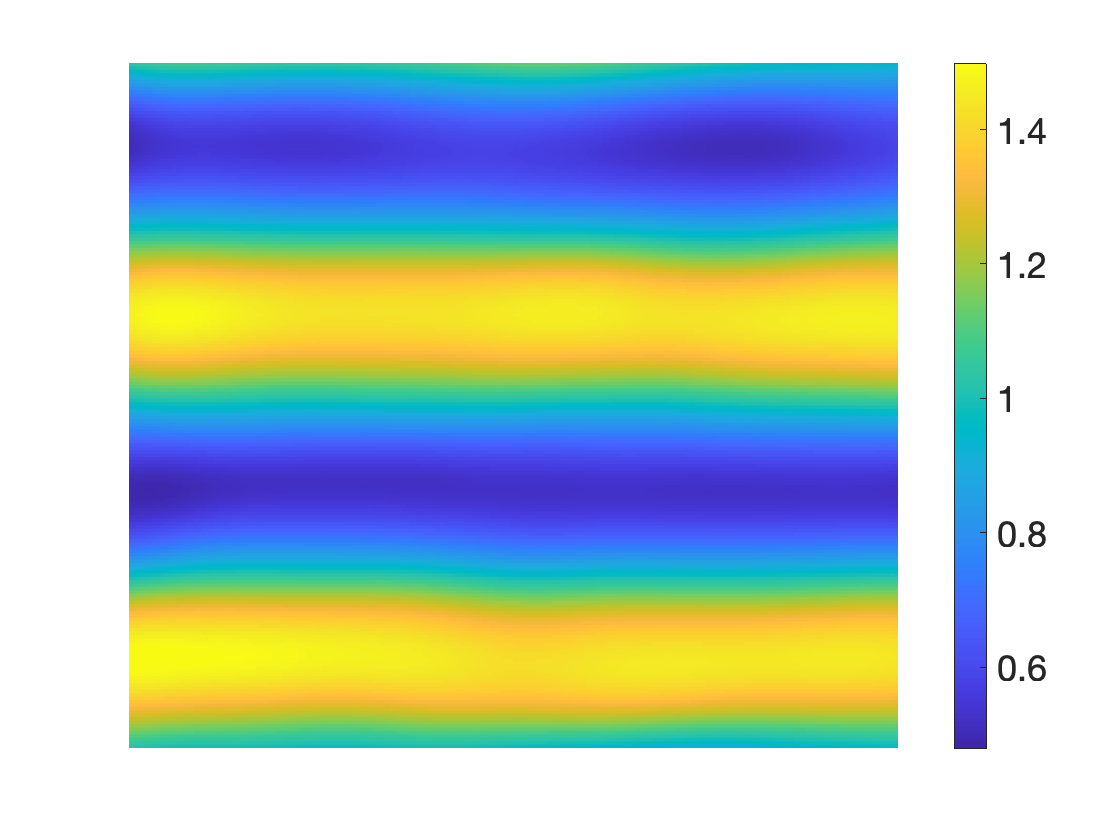} &
\includegraphics[width=0.199\textwidth]{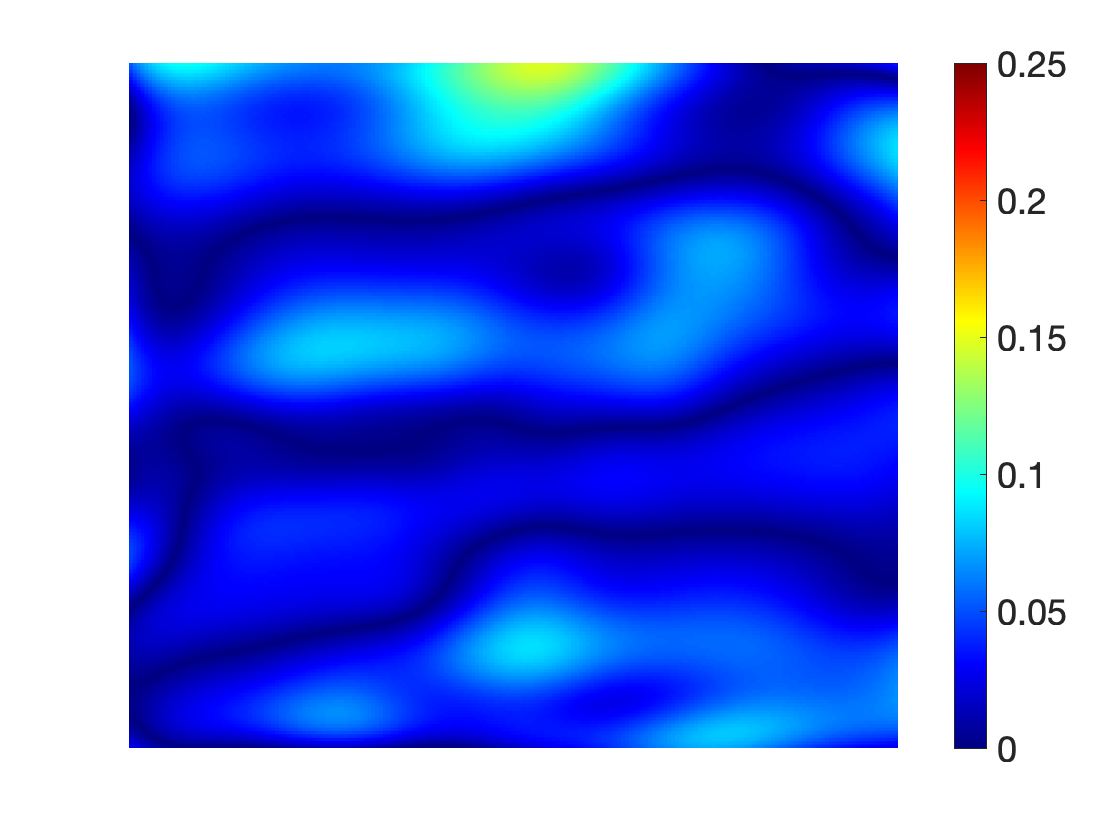}\\
\includegraphics[width=0.199\textwidth]{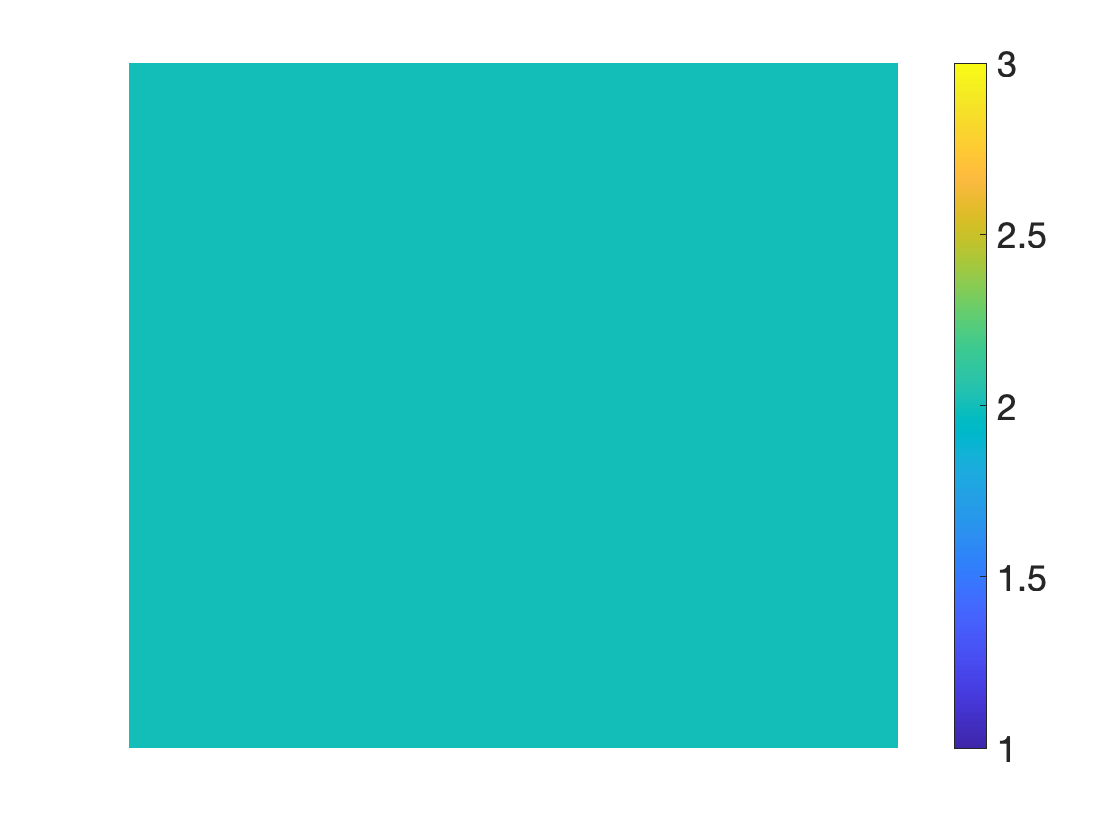} &
\includegraphics[width=0.199\textwidth]{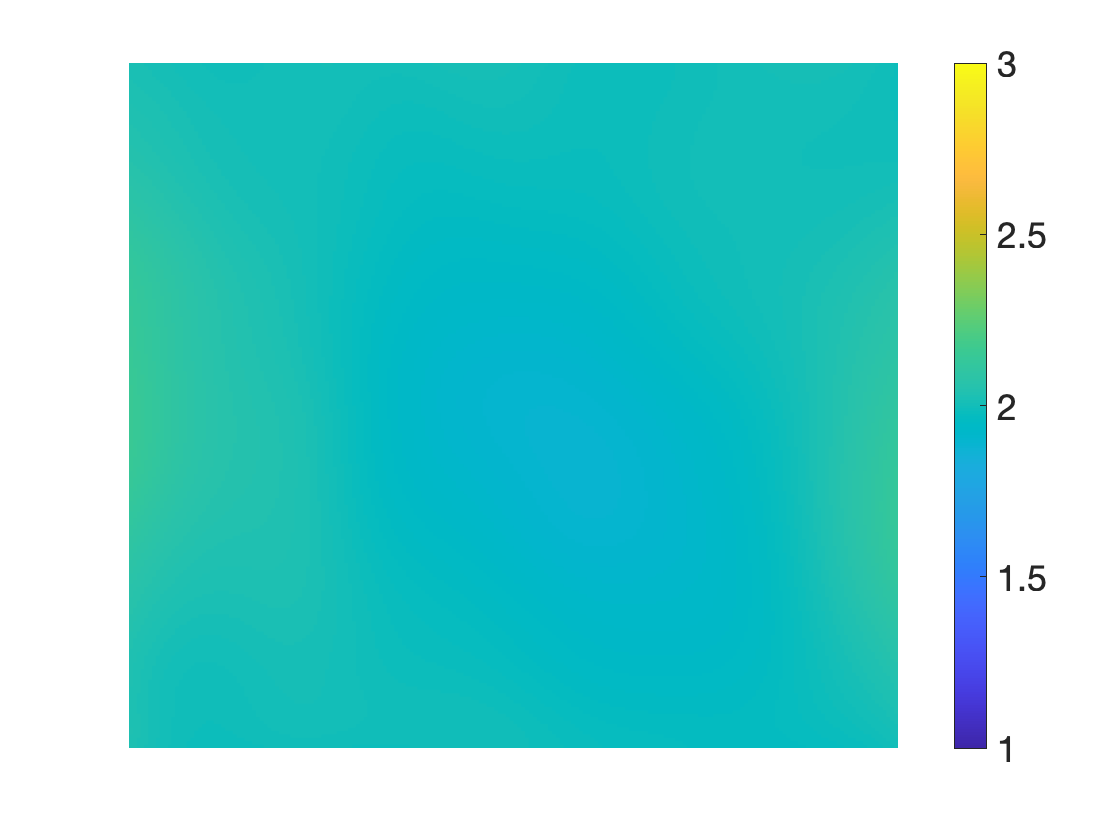} &
\includegraphics[width=0.199\textwidth]{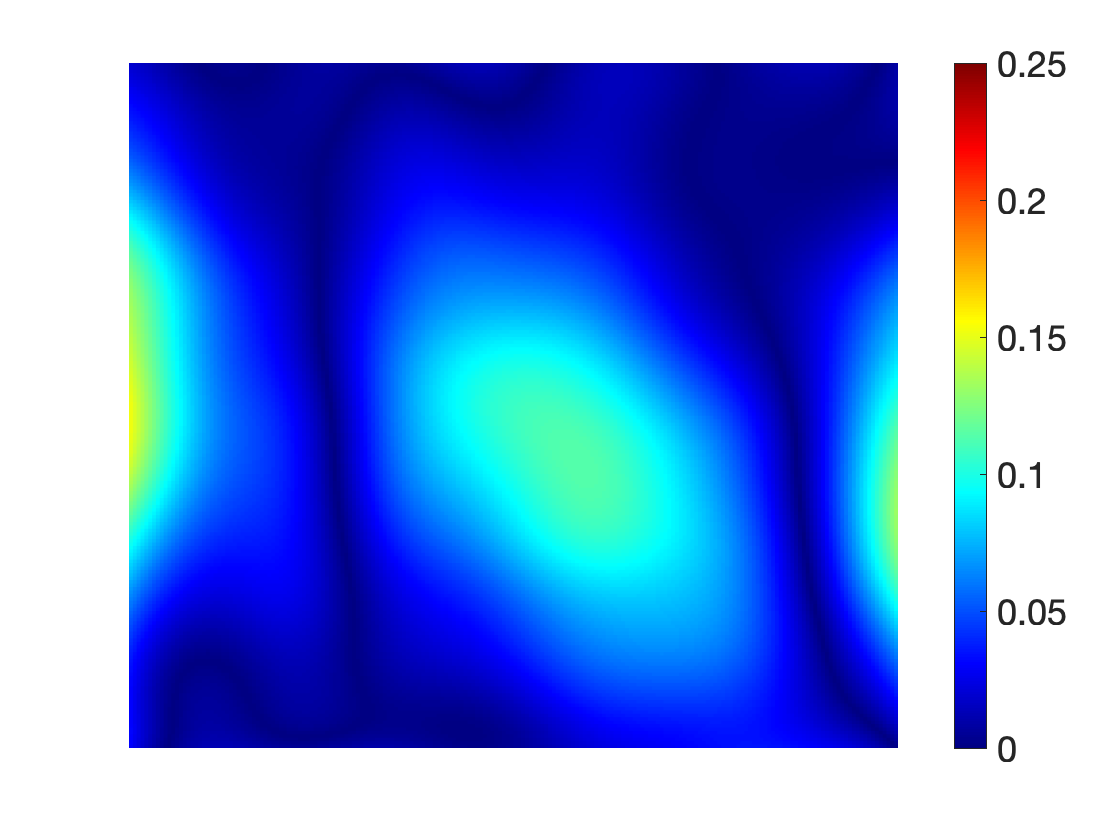} &
\includegraphics[width=0.199\textwidth]{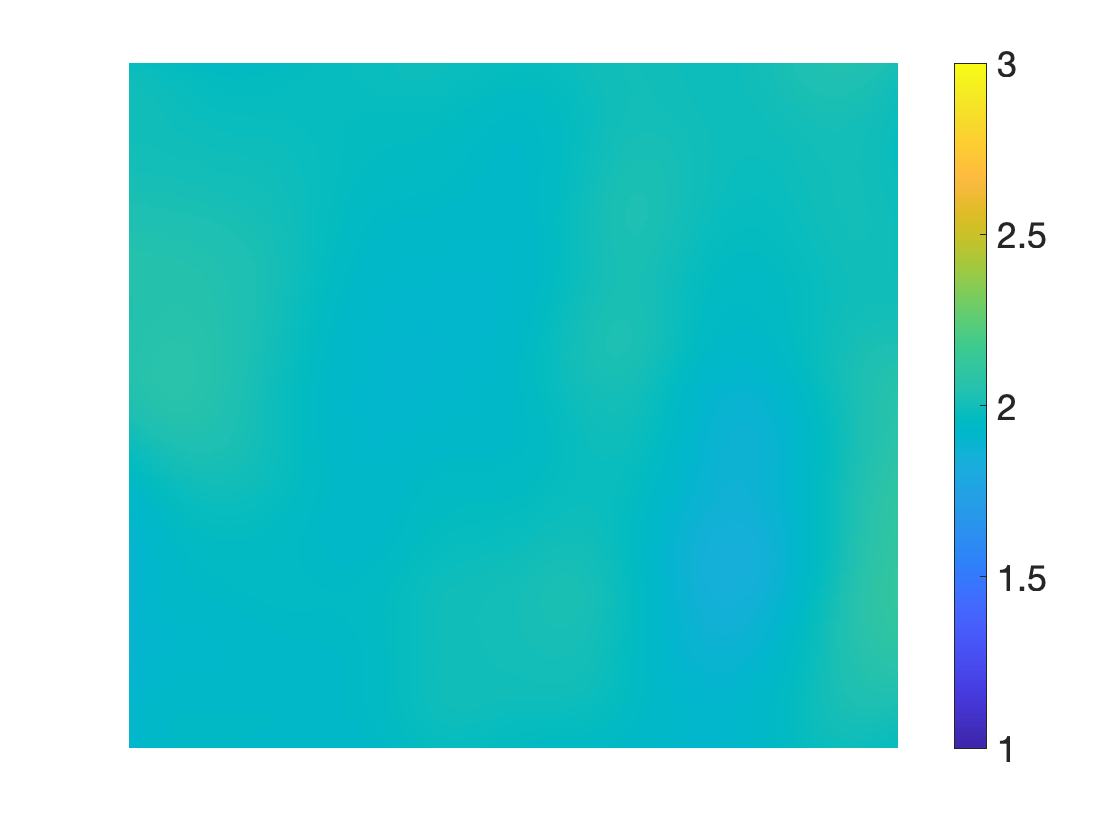} &
\includegraphics[width=0.199\textwidth]{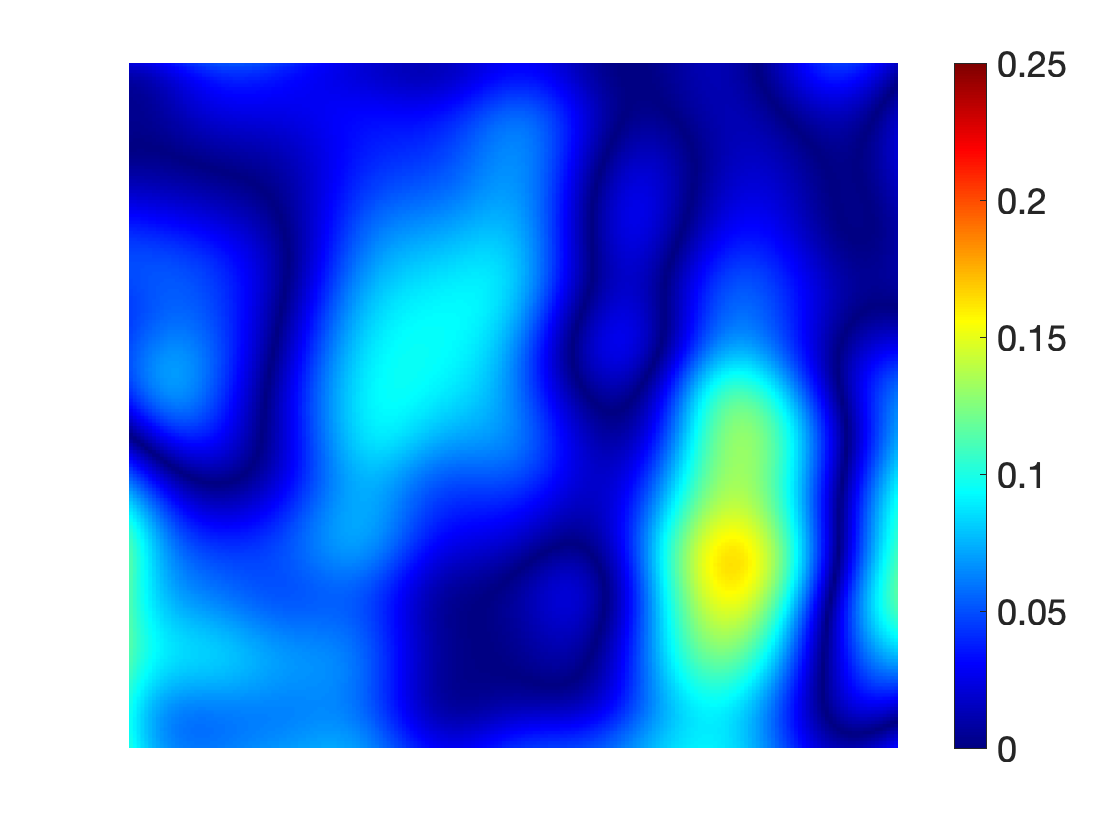}\\
\includegraphics[width=0.199\textwidth]{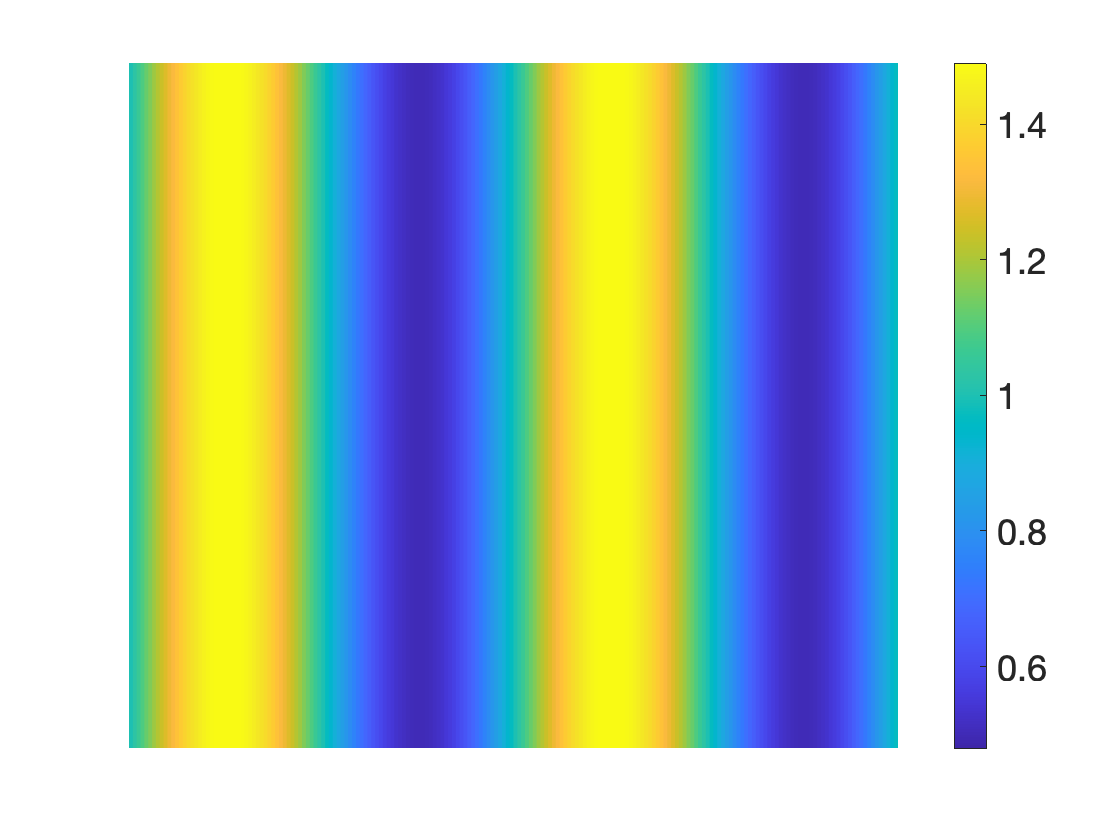} &
\includegraphics[width=0.199\textwidth]{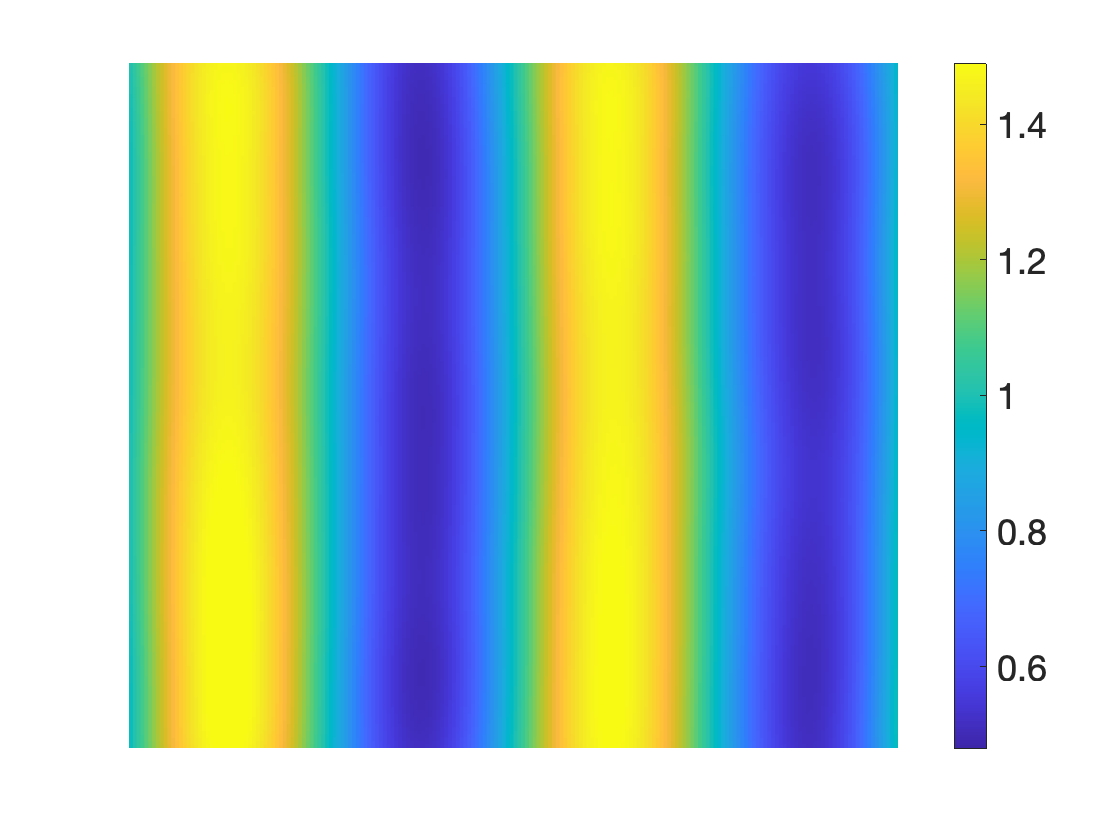} &
\includegraphics[width=0.199\textwidth]{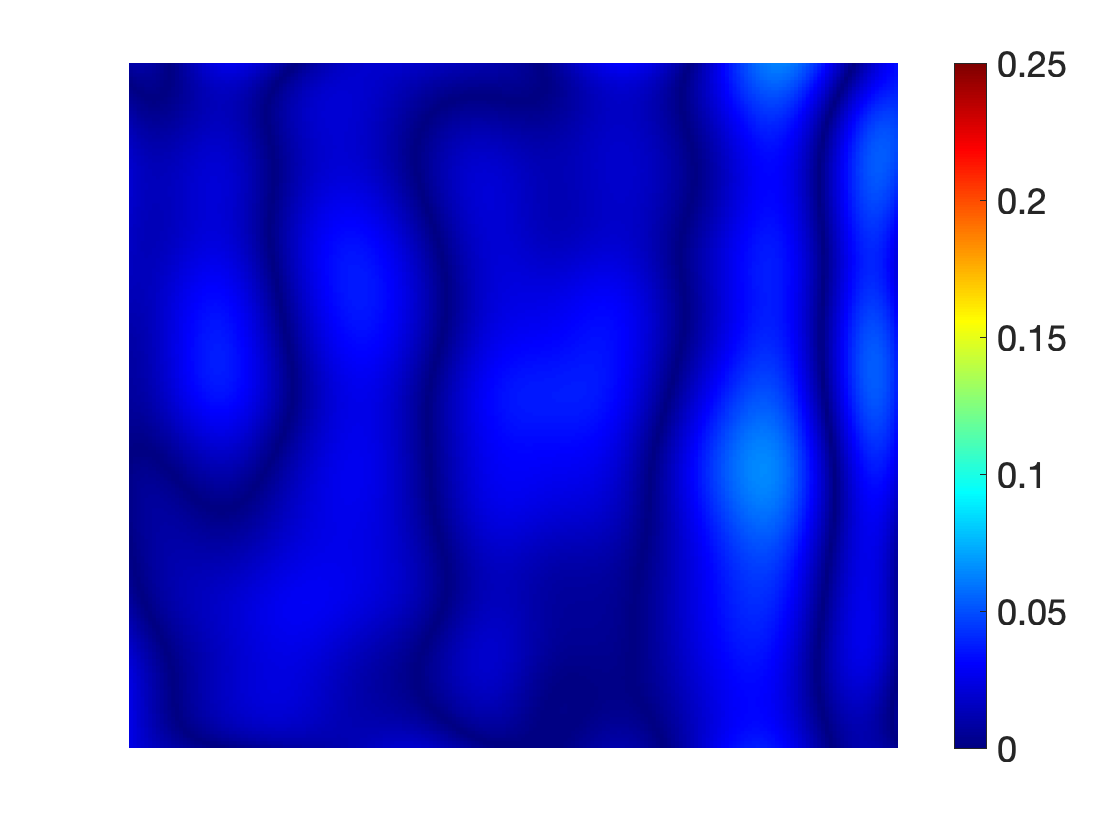} &
\includegraphics[width=0.199\textwidth]{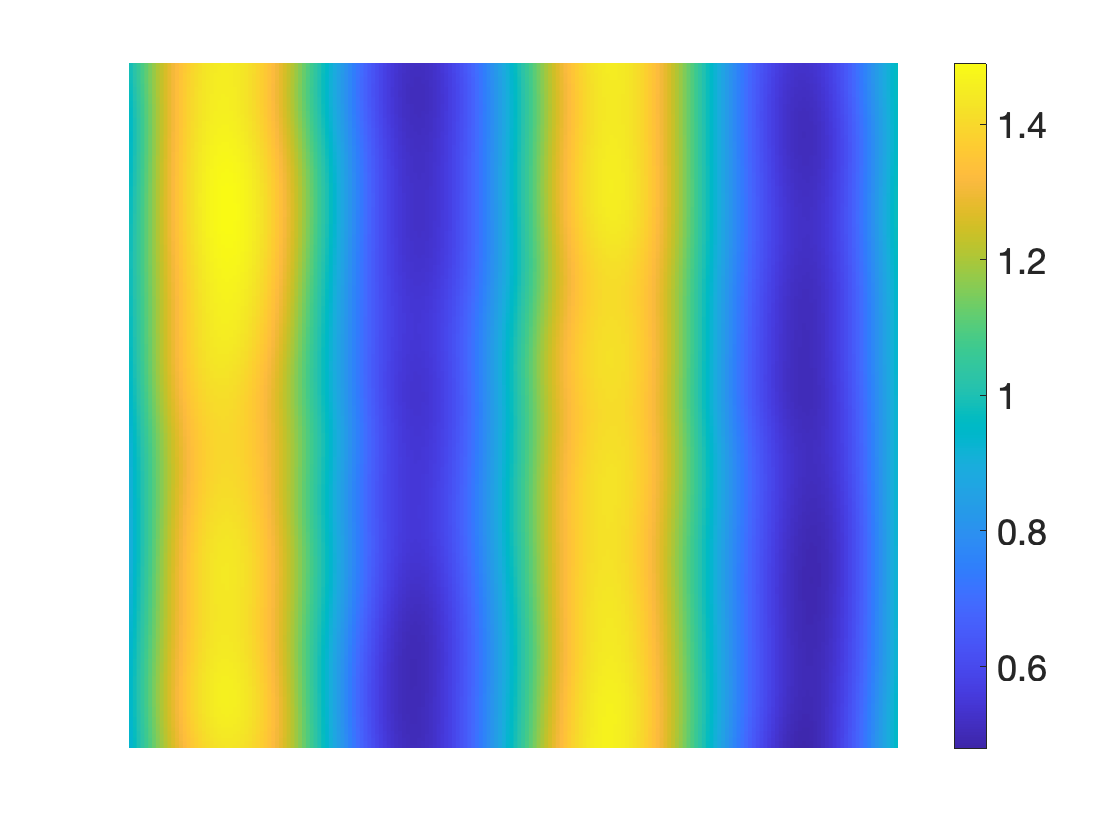} &
\includegraphics[width=0.199\textwidth]{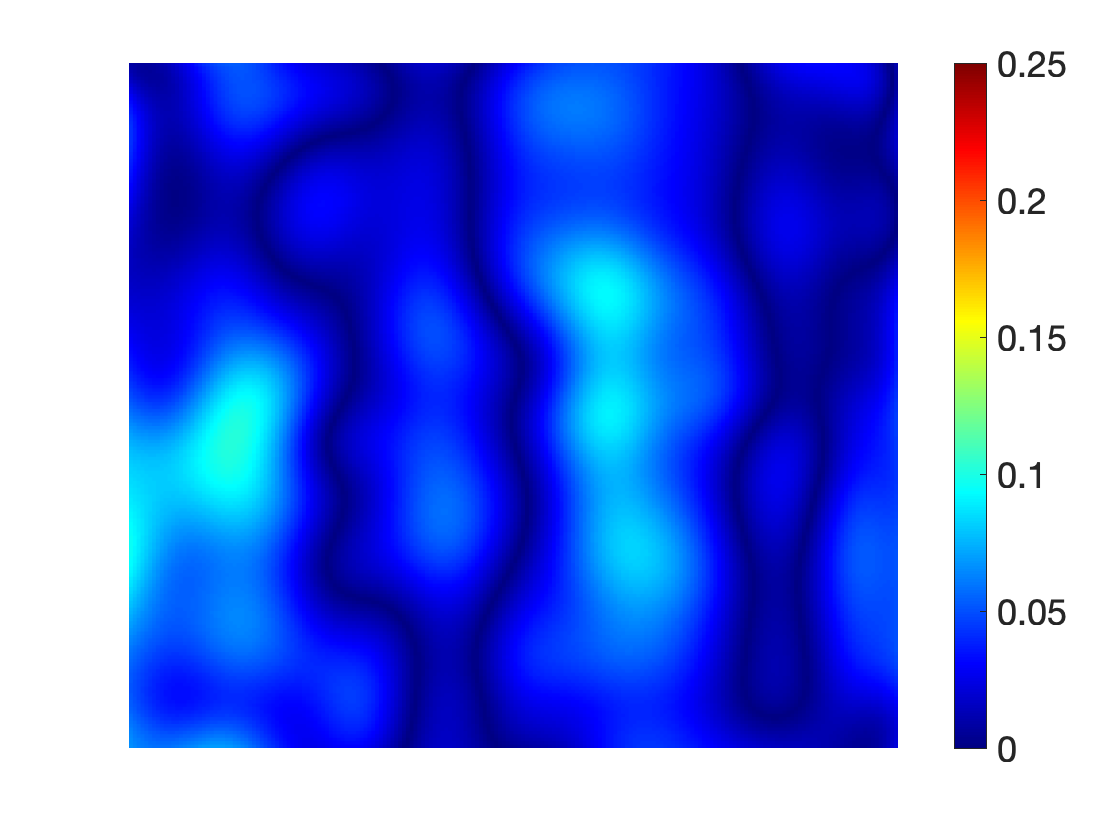}\\
\includegraphics[width=0.199\textwidth]{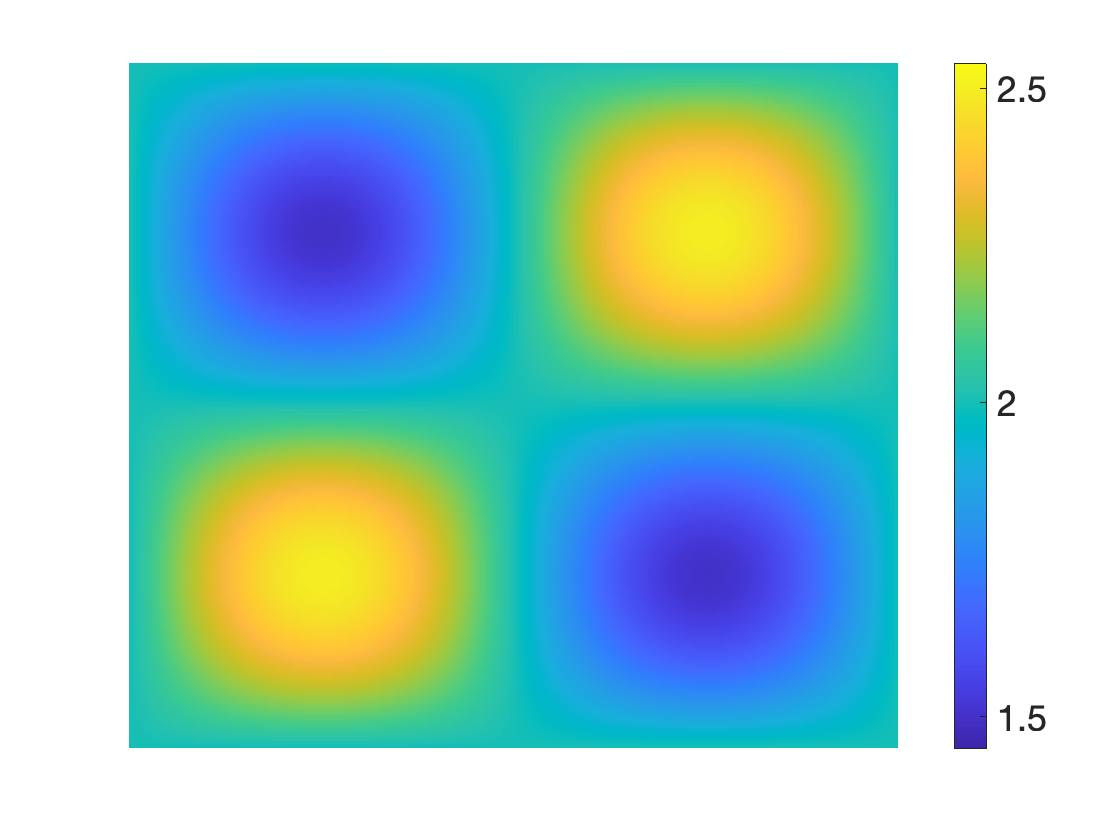} &
\includegraphics[width=0.199\textwidth]{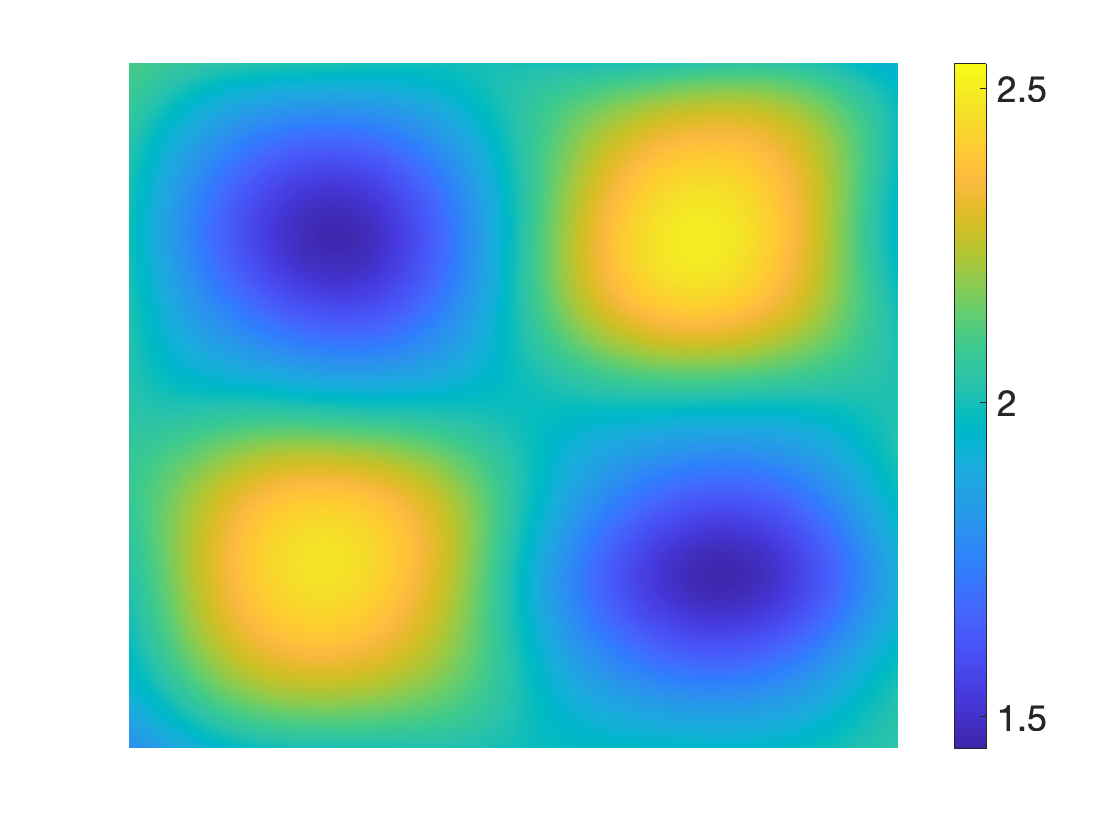} &
\includegraphics[width=0.199\textwidth]{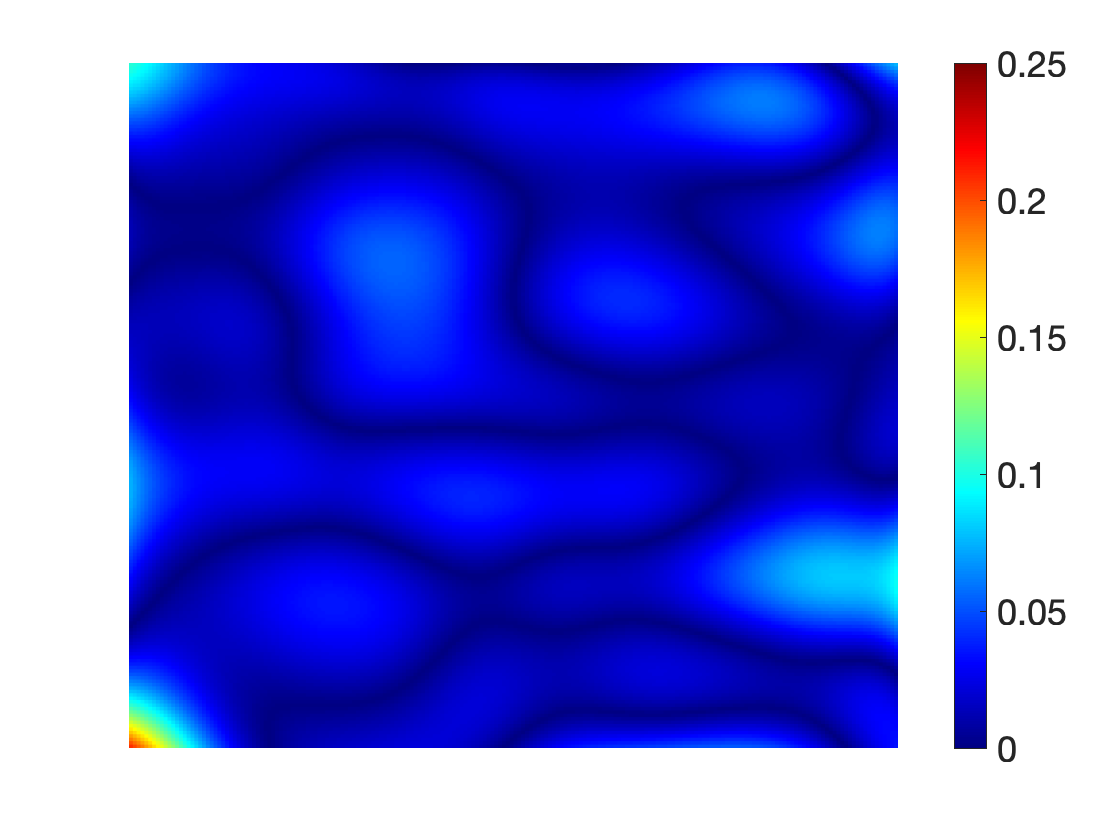} &
\includegraphics[width=0.199\textwidth]{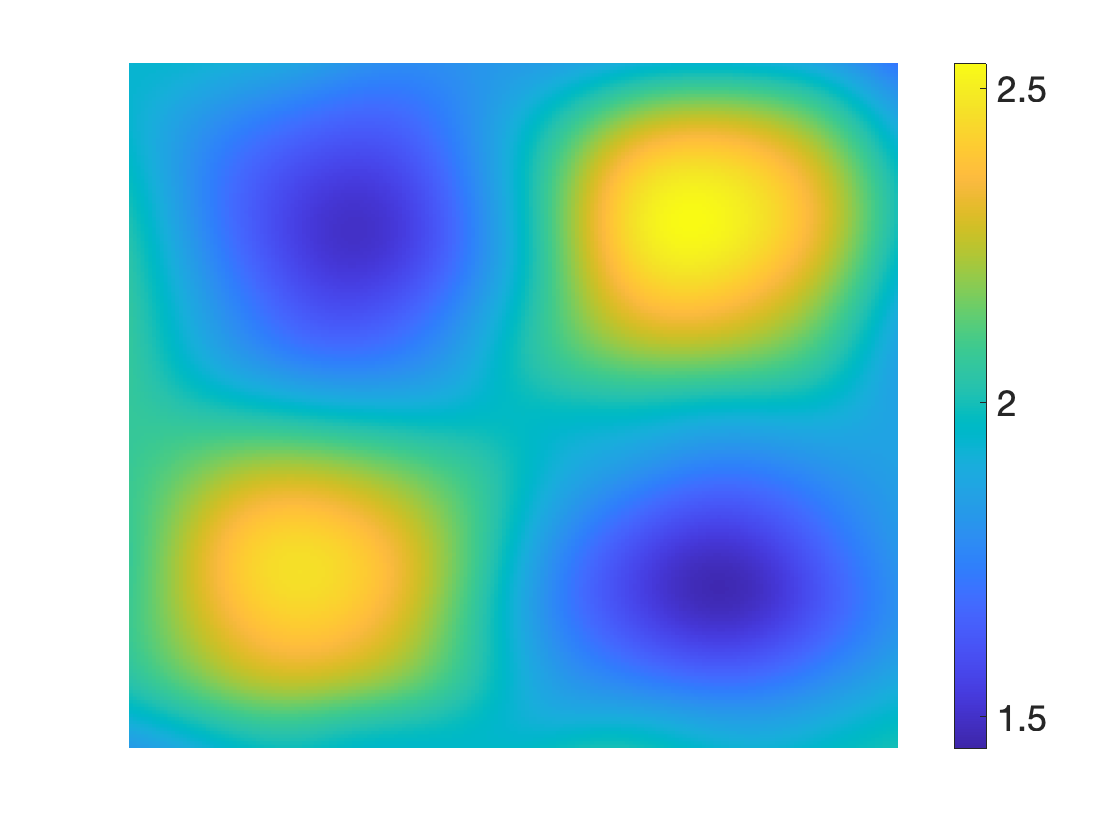} &
\includegraphics[width=0.199\textwidth]{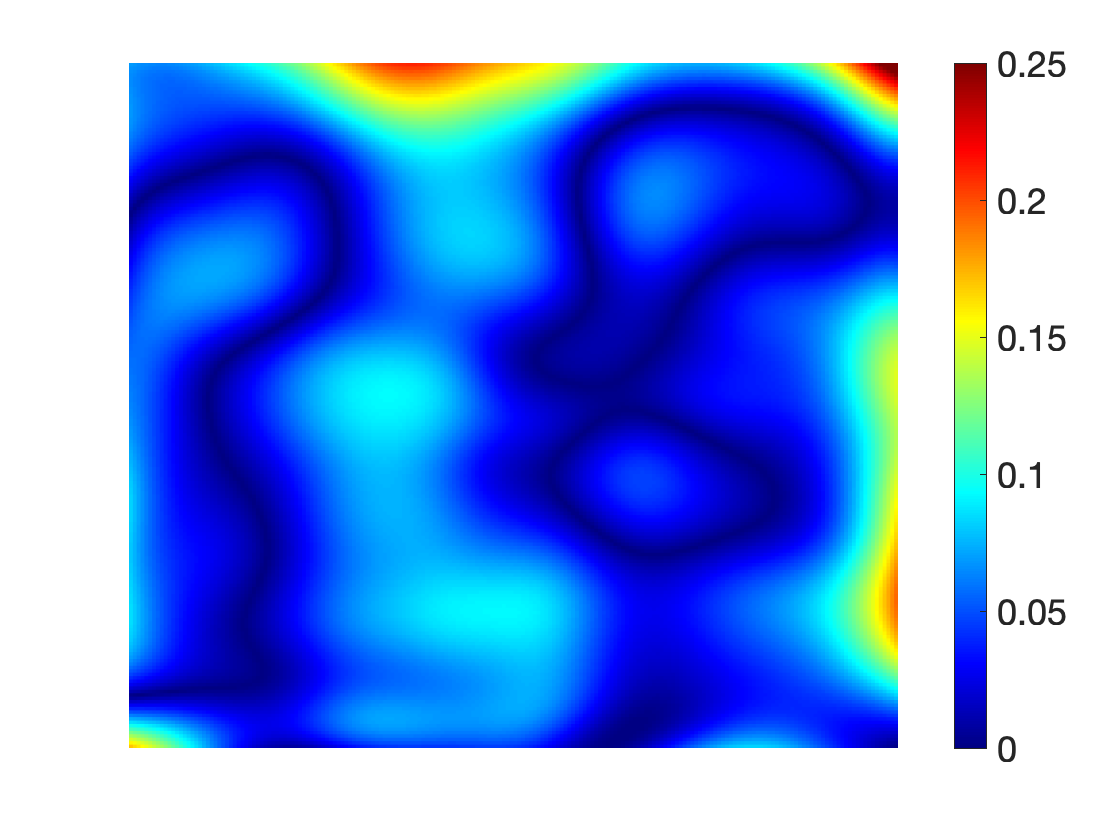}\\
(a) $A^\dag$  & (b) $\hat A$ & (c) $|\hat A-A^\dag|$ & (d) $\hat A$ & (e) $|\hat A-A^\dag|$
\end{tabular}
\caption{The reconstructions for Example \ref{exam:neu3d1} with exact data in (b) and noisy data $(\delta=5\%)$ in (d). From the top to bottom, the results are for $A_{11}$, $A_{12}$, $A_{13}$, $A_{22}$, $A_{23}$ and $A_{33}$, respectively.}
\label{fig:neu3d1}
\end{figure}

Fig. \ref{fig:neu3d1} shows the results on a 2D cross section at $x_3=0.5$ for both exact and noisy data ($\delta=5\%$). The overall features of the anisotropic conductivity are well resolved at both noise levels, showing the high robustness of the DNN approach to data noise. In addition, the grid-free nature of the DNN discretization makes the method easy to implement in the 3D case when compared to the traditional FEM.

The last example is about recovering a 3D conductivity matrix from partial internal data.
\begin{example}\label{exam:neu3d2}
    The domain $\Omega=(0,1)^3$, the measurement $\nabla z_i^\delta$ on the region $\omega=\Omega\setminus (0.2,0.8)^3$, $A^\dagger = \begin{pmatrix}
    2+x_2^2& \frac{1}{2}+x_3^2& 1+\frac{\sin(4\pi x_2)}{2}\\
    \frac12+x_3^2& 2+x_1^2& 1+\frac{\sin(4\pi x_1)}{2}\\
    1+\frac{\sin(4\pi x_2)}{2}& 1+\frac{\sin(4\pi x_1)}{2}& 2+x_1^2+x_2^2\\
    \end{pmatrix}$, $u_1^\dag=x_1+x_2+x_3+\frac{1}{3}(x_1^3+x_2^3+x_3^3)$, $u_2^\dag=x_1-x_2+x_3+\frac{1}{3}(x_1^3-x_2^3+x_3^3)$, $u_3^\dag=x_1+x_2-x_3+\frac{1}{3}(x_1^3+x_2^3-x_3^3)$, $u_4^\dag=-x_1+x_2+x_3+\frac{1}{3}(-x_1^3+x_2^3+x_3^3)$, $u_5^\dagger=-u_3^\dagger$ and $u_6^\dagger=-u_2^\dagger$.
\end{example}

\begin{figure}[htb!]
\centering
\setlength{\tabcolsep}{0em}
\begin{tabular}{ccccc}
\includegraphics[width=0.199\textwidth]{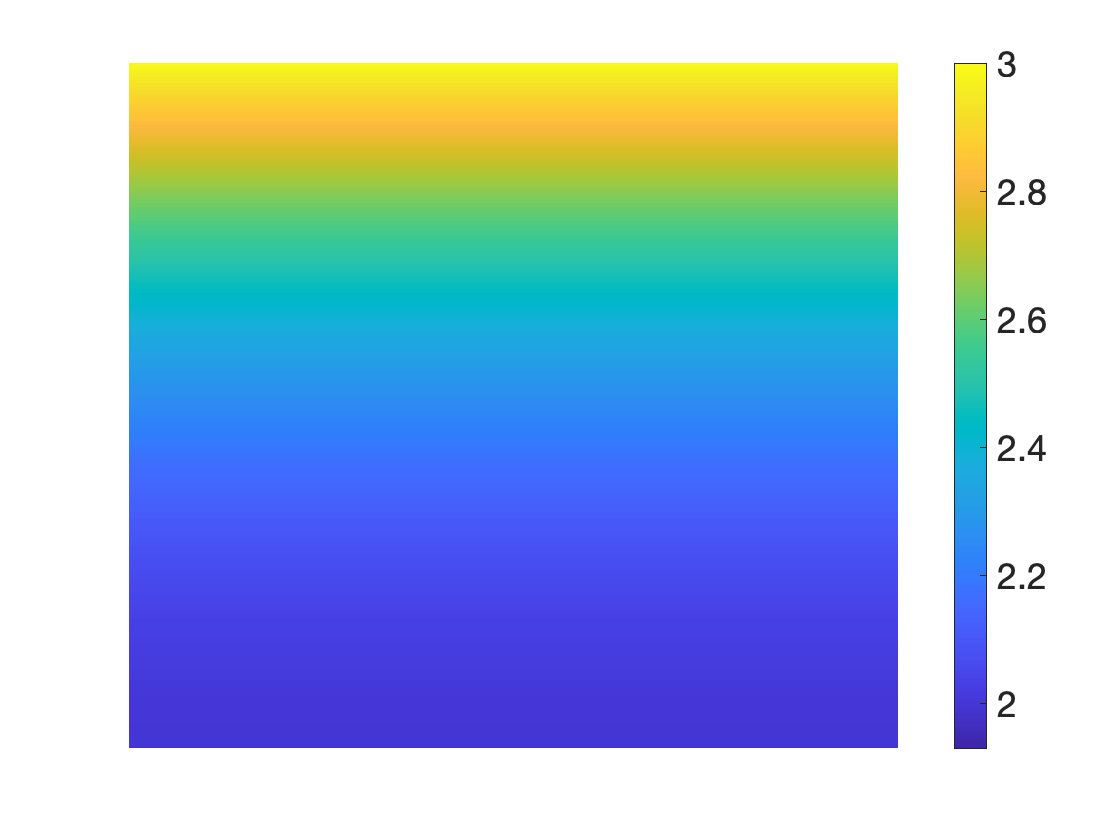} &
\includegraphics[width=0.199\textwidth]{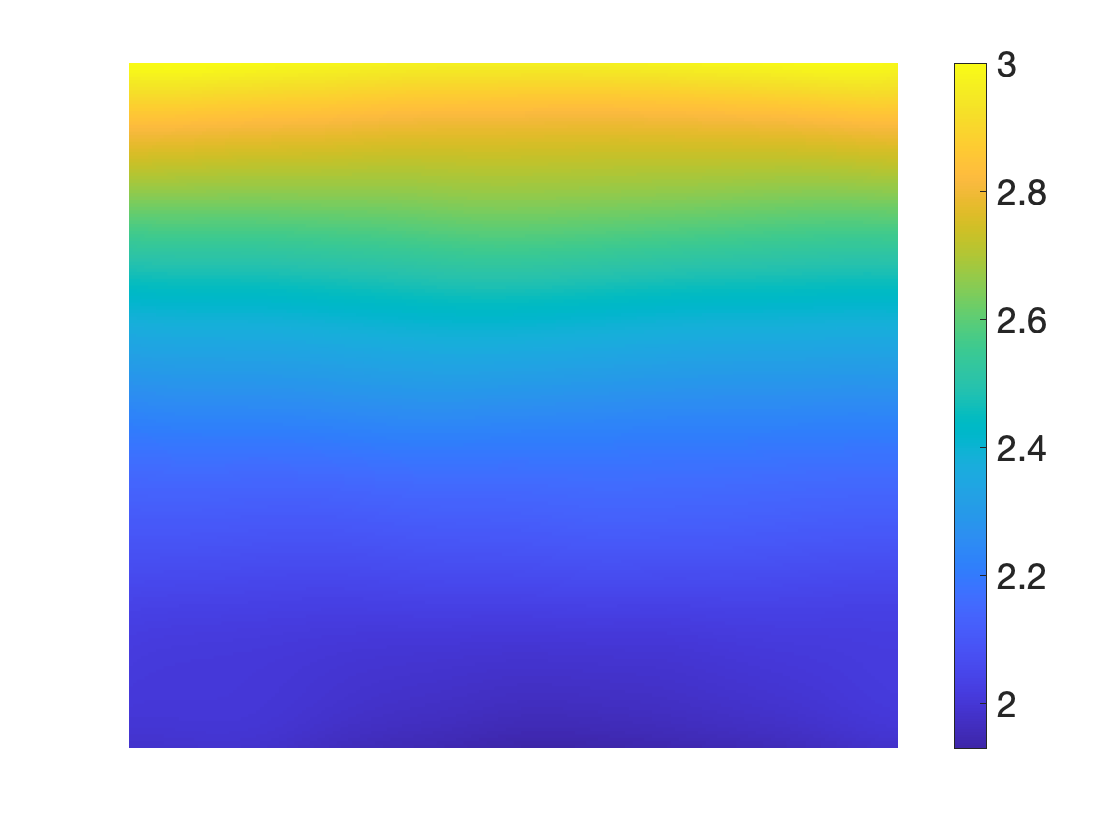} &
\includegraphics[width=0.199\textwidth]{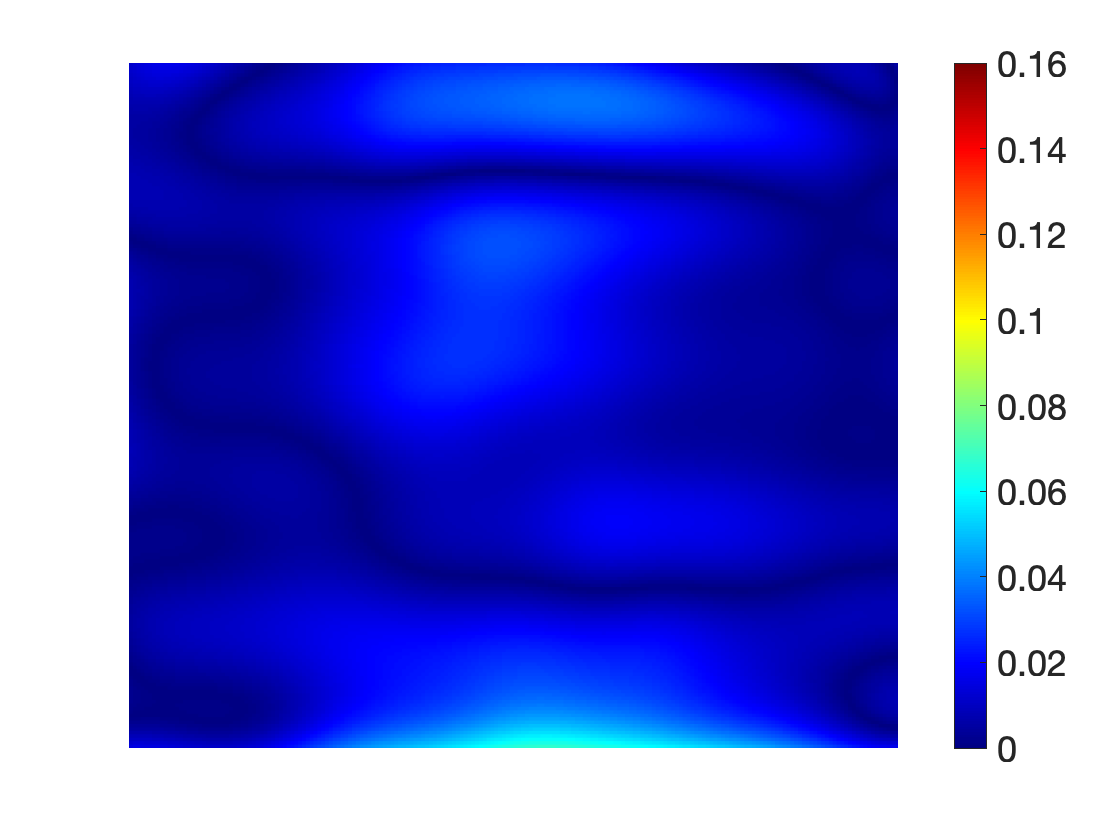} &
\includegraphics[width=0.199\textwidth]{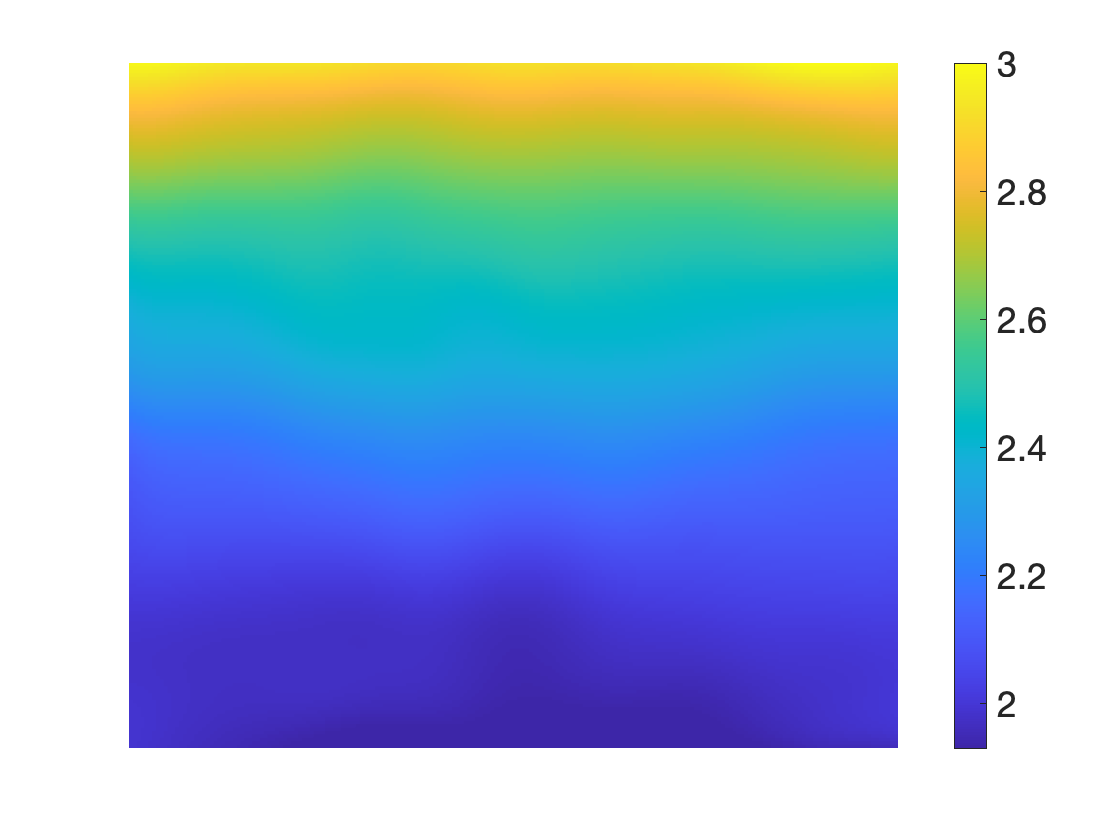} &
\includegraphics[width=0.199\textwidth]{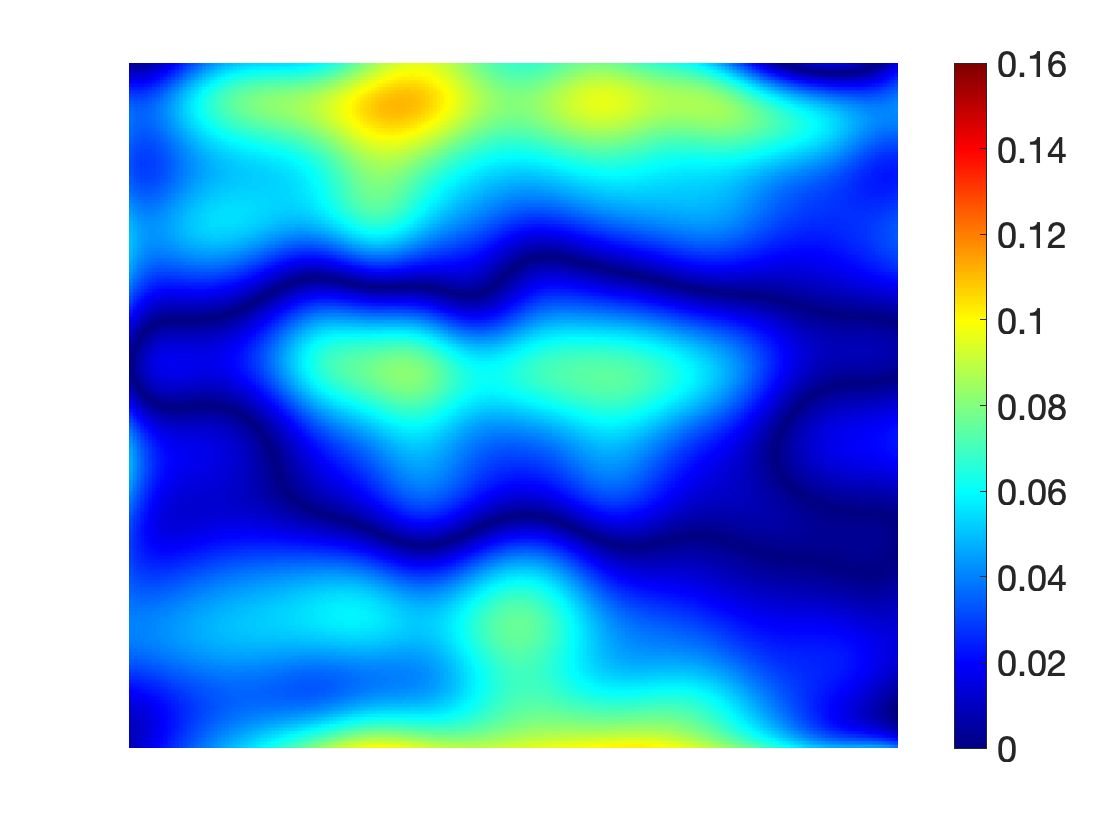}\\
\includegraphics[width=0.199\textwidth]{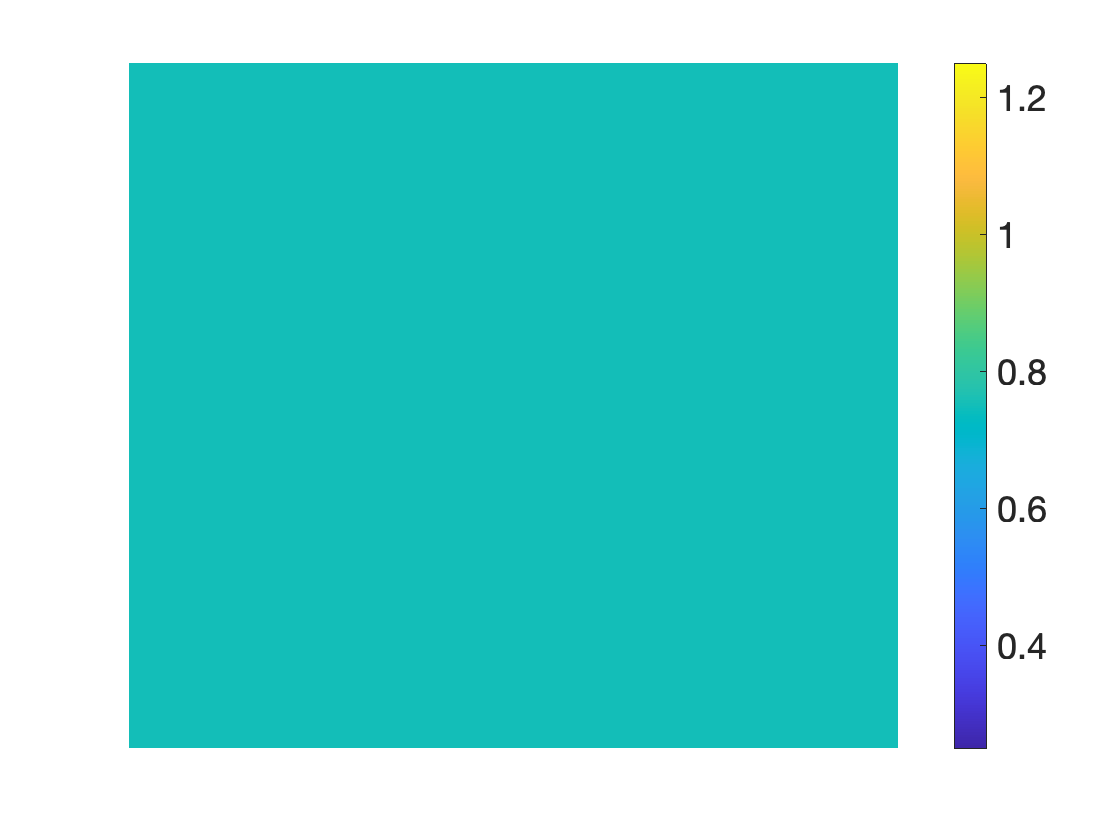} &
\includegraphics[width=0.199\textwidth]{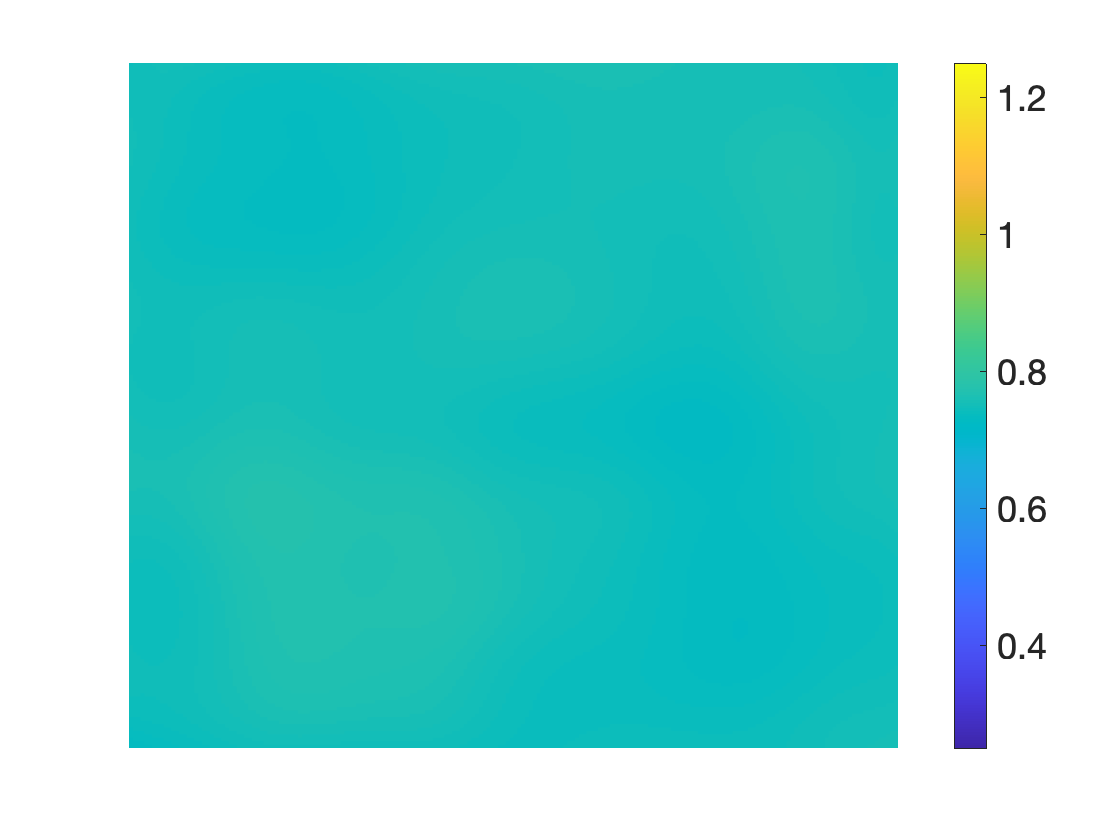} &
\includegraphics[width=0.199\textwidth]{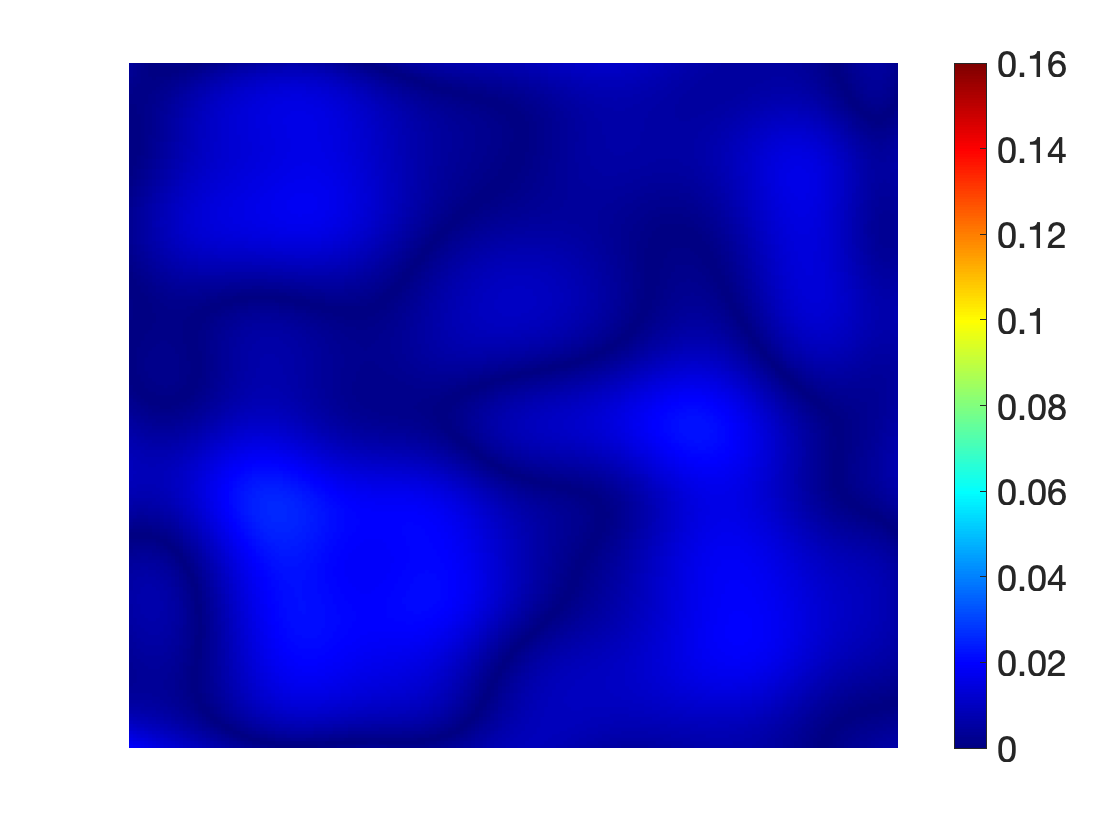}&
\includegraphics[width=0.199\textwidth]{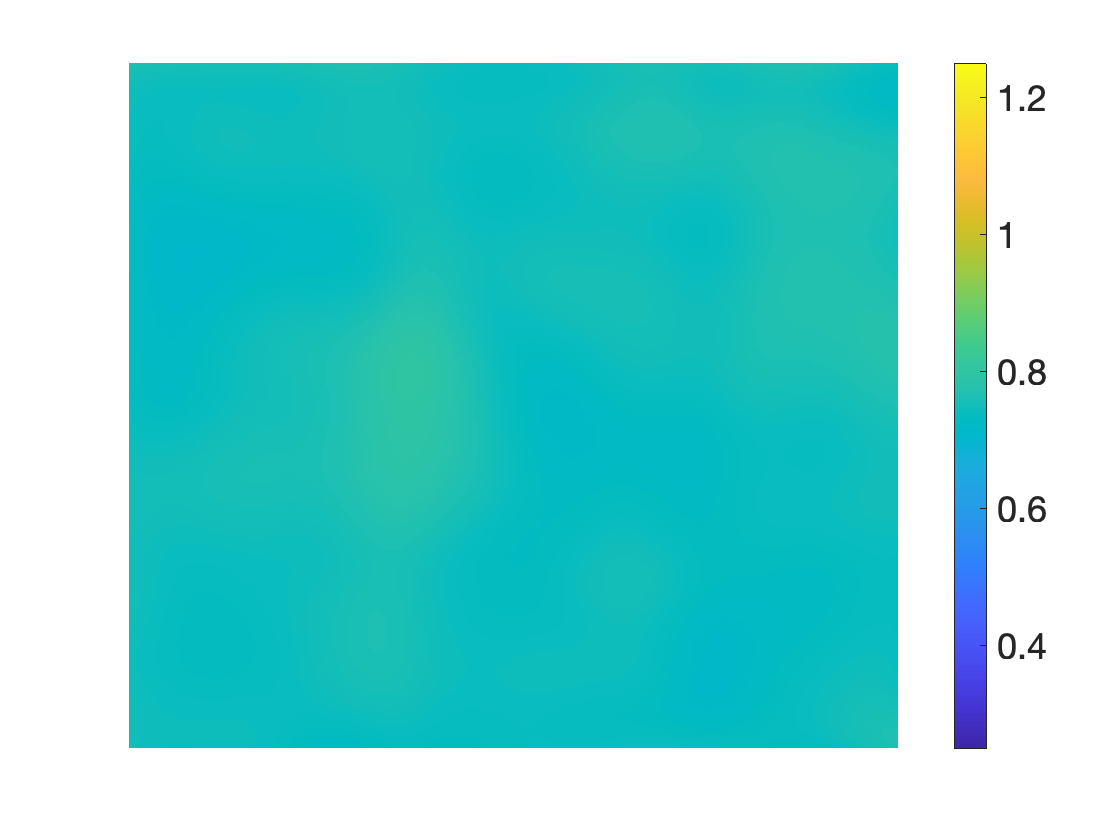} &
\includegraphics[width=0.199\textwidth]{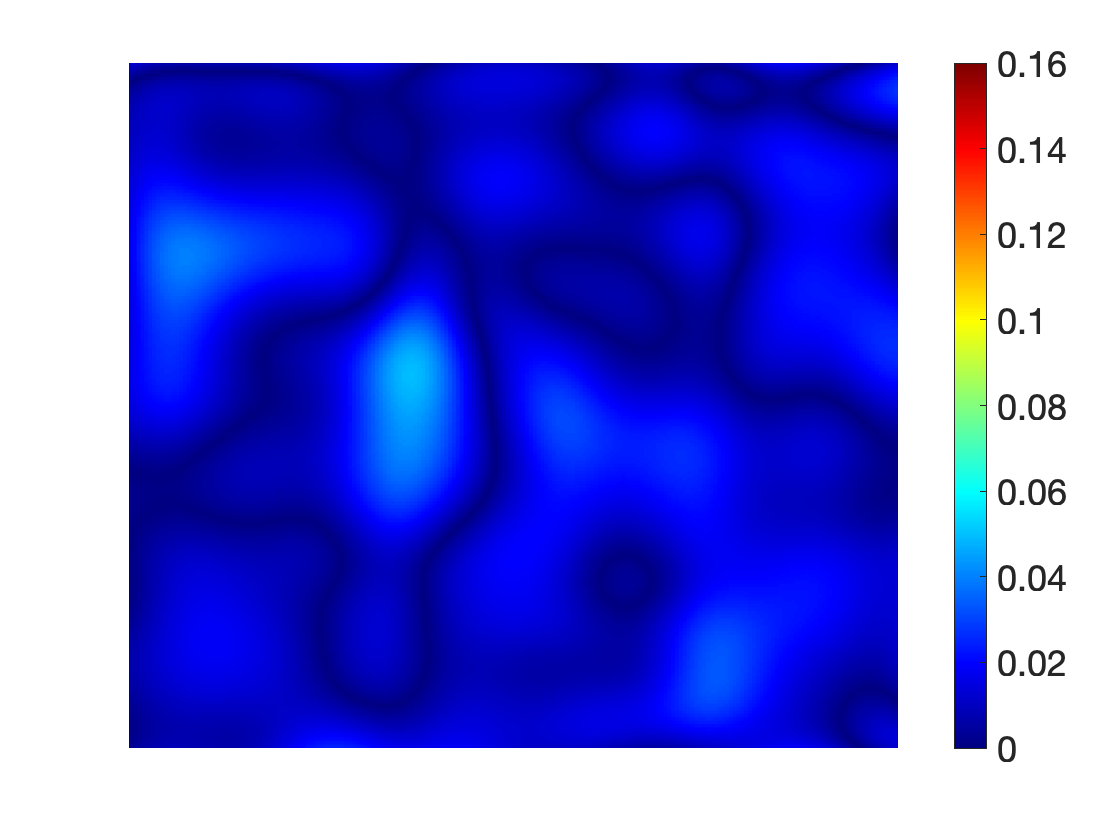}\\
\includegraphics[width=0.199\textwidth]{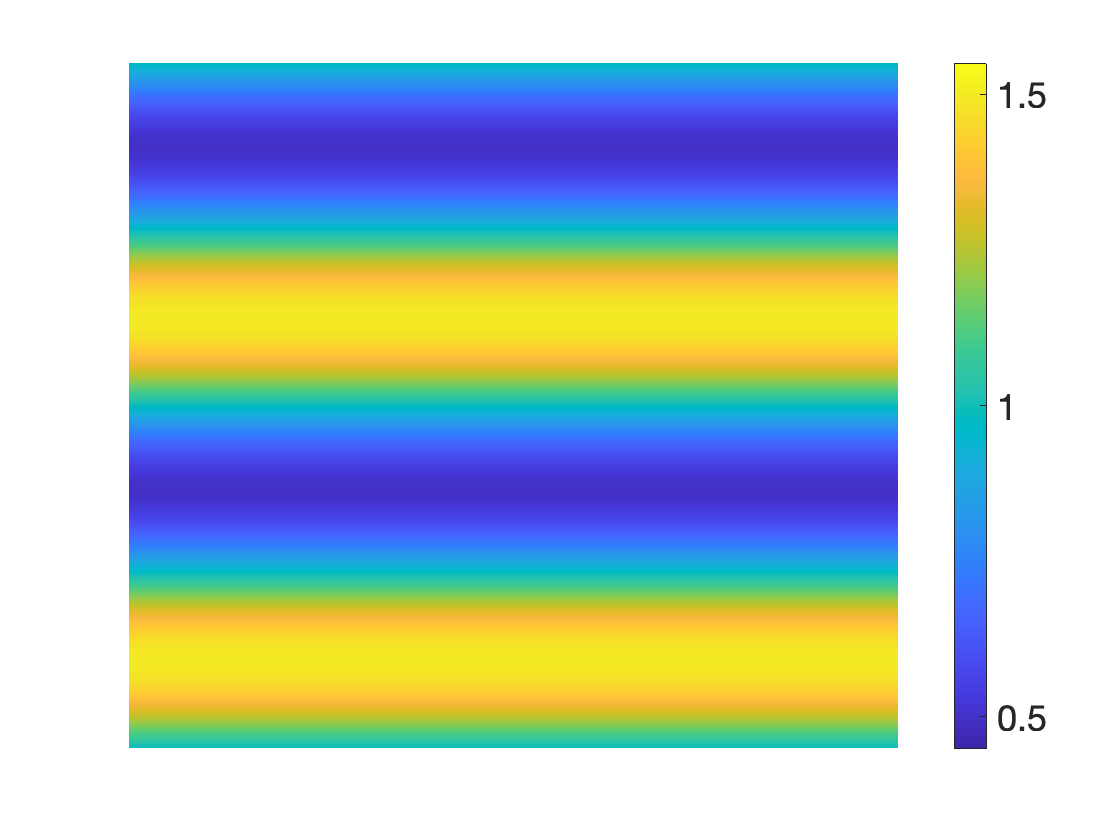} &
\includegraphics[width=0.199\textwidth]{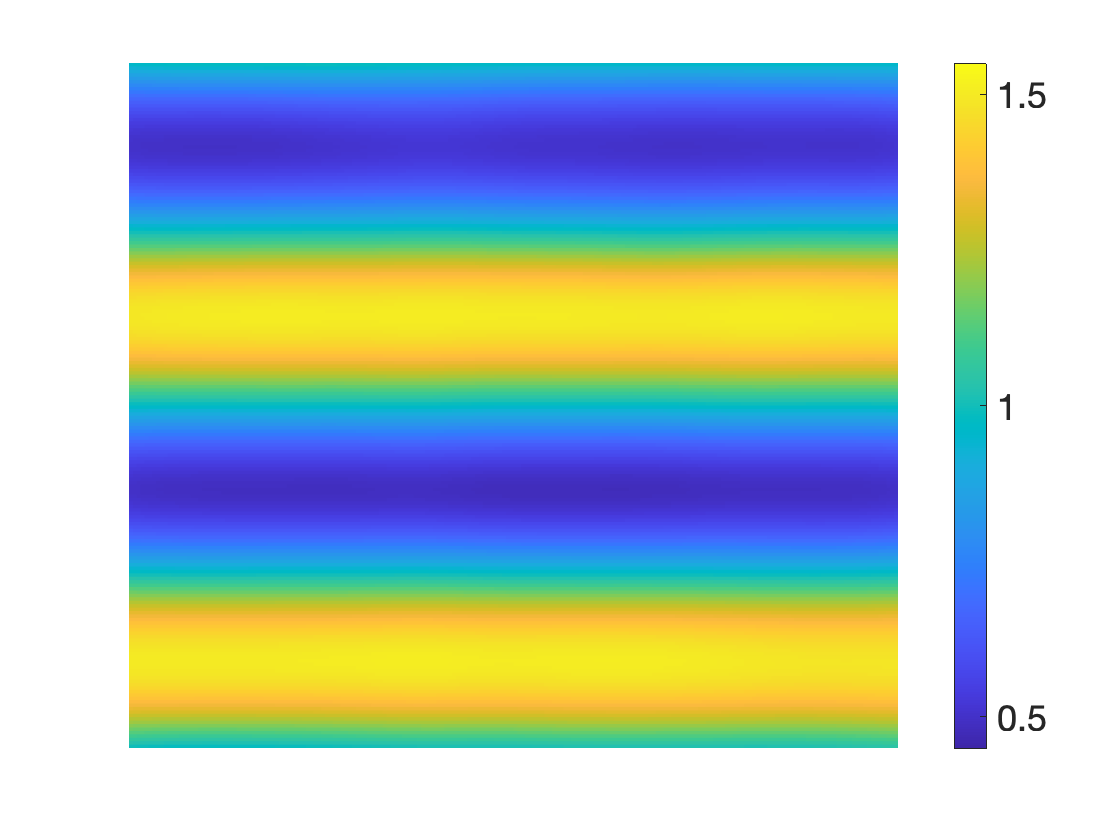} &
\includegraphics[width=0.199\textwidth]{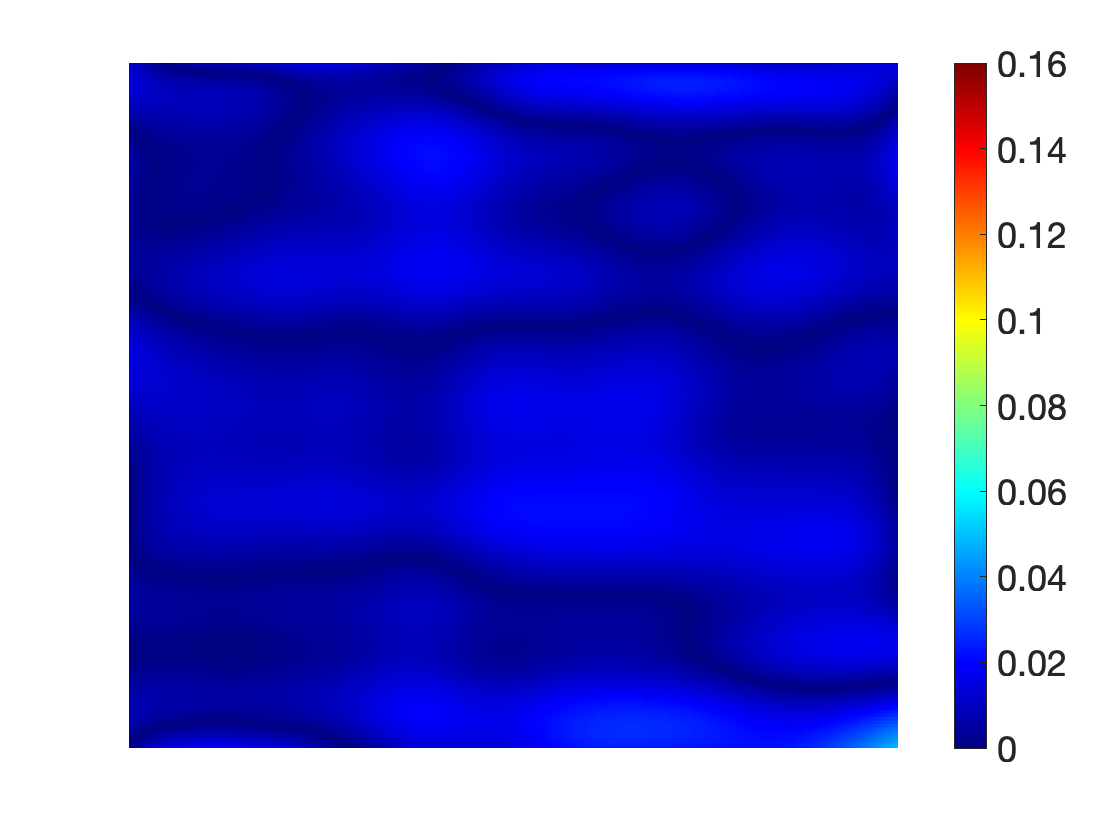} &
\includegraphics[width=0.199\textwidth]{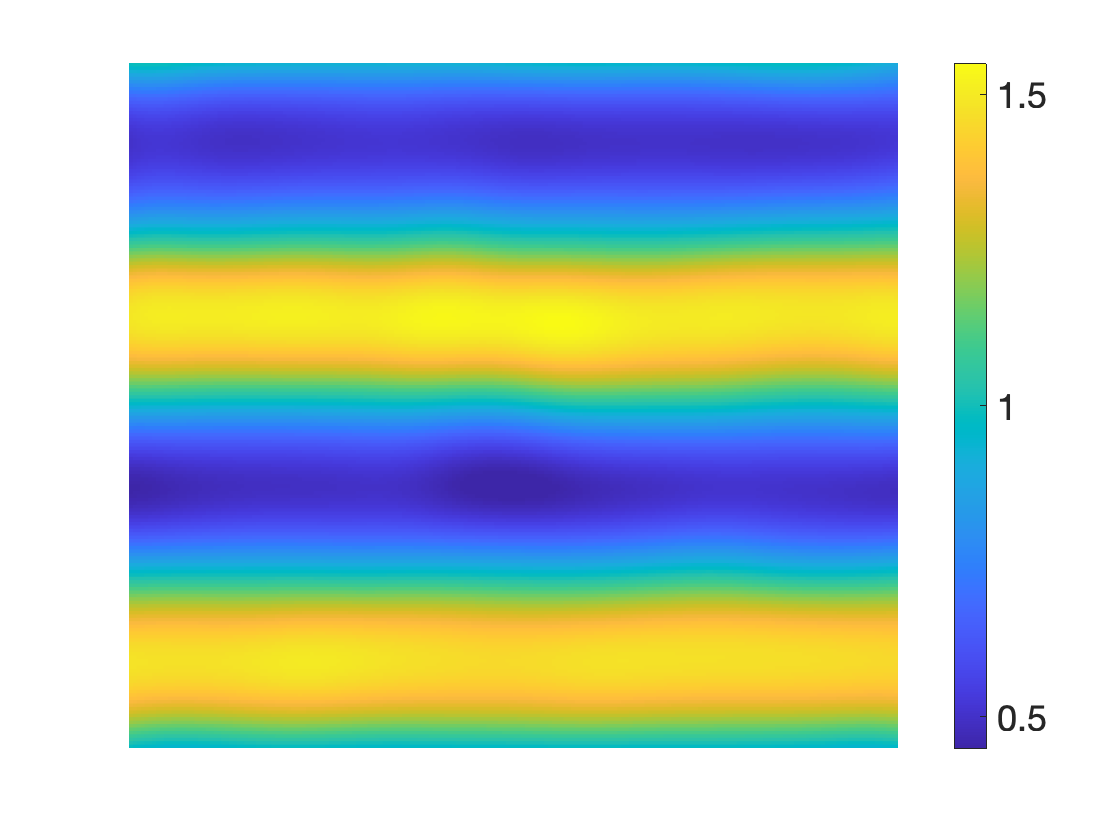} &
\includegraphics[width=0.199\textwidth]{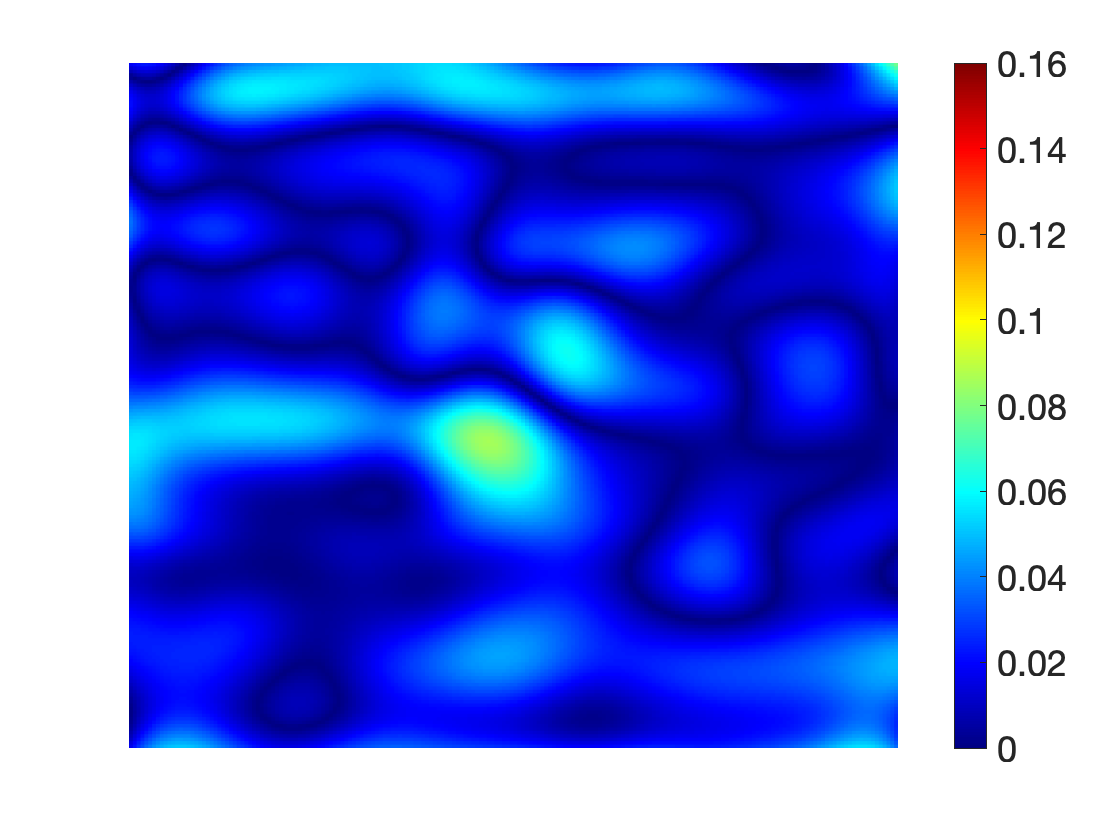}\\
\includegraphics[width=0.199\textwidth]{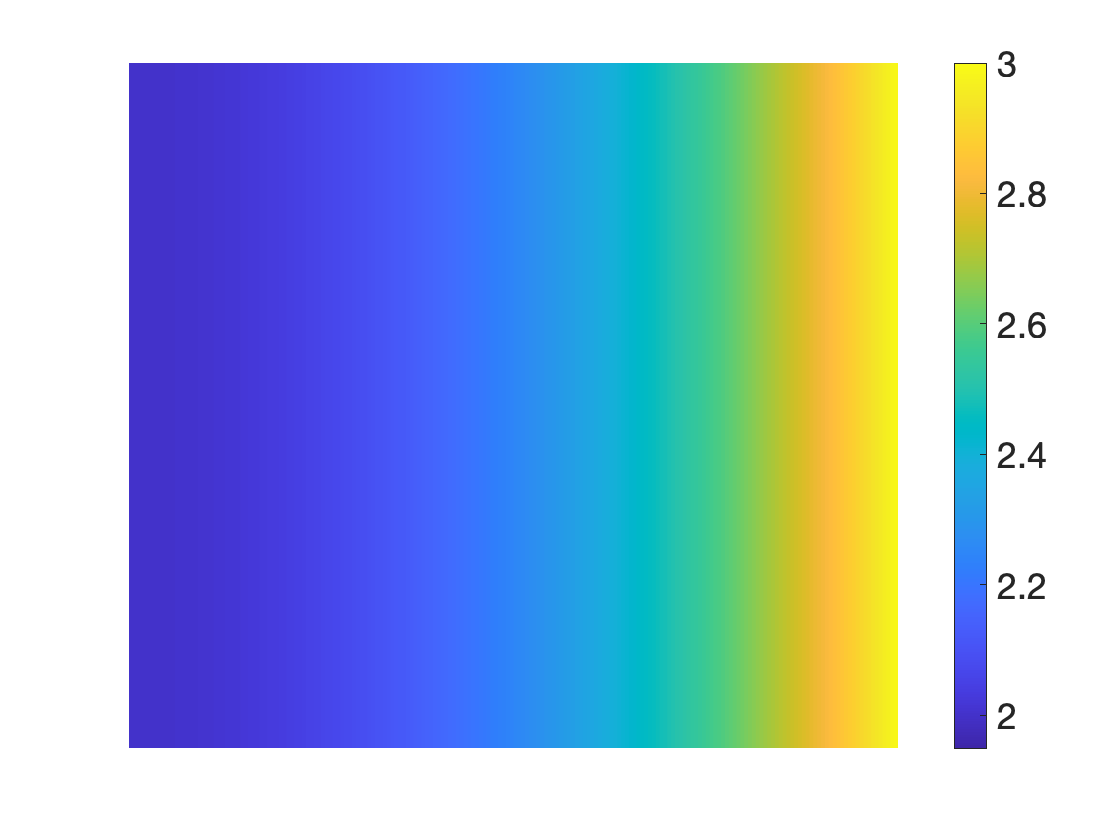} &
\includegraphics[width=0.199\textwidth]{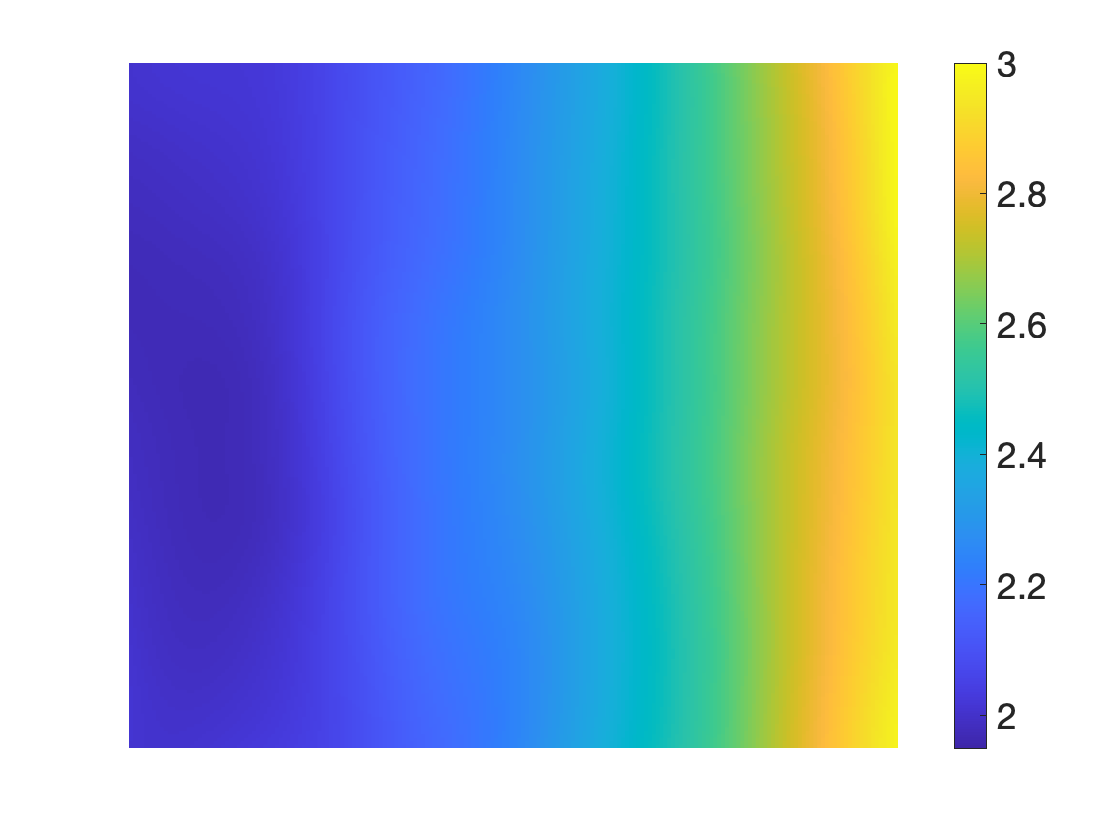} &
\includegraphics[width=0.199\textwidth]{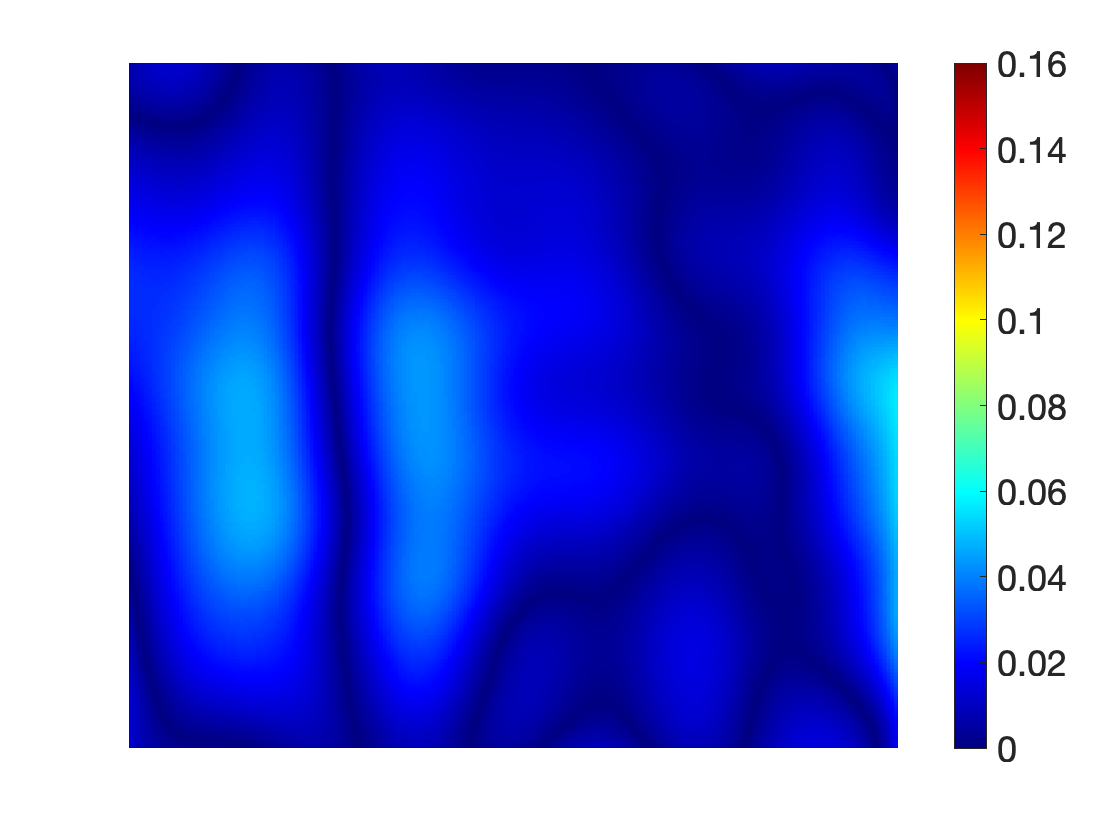} &
\includegraphics[width=0.199\textwidth]{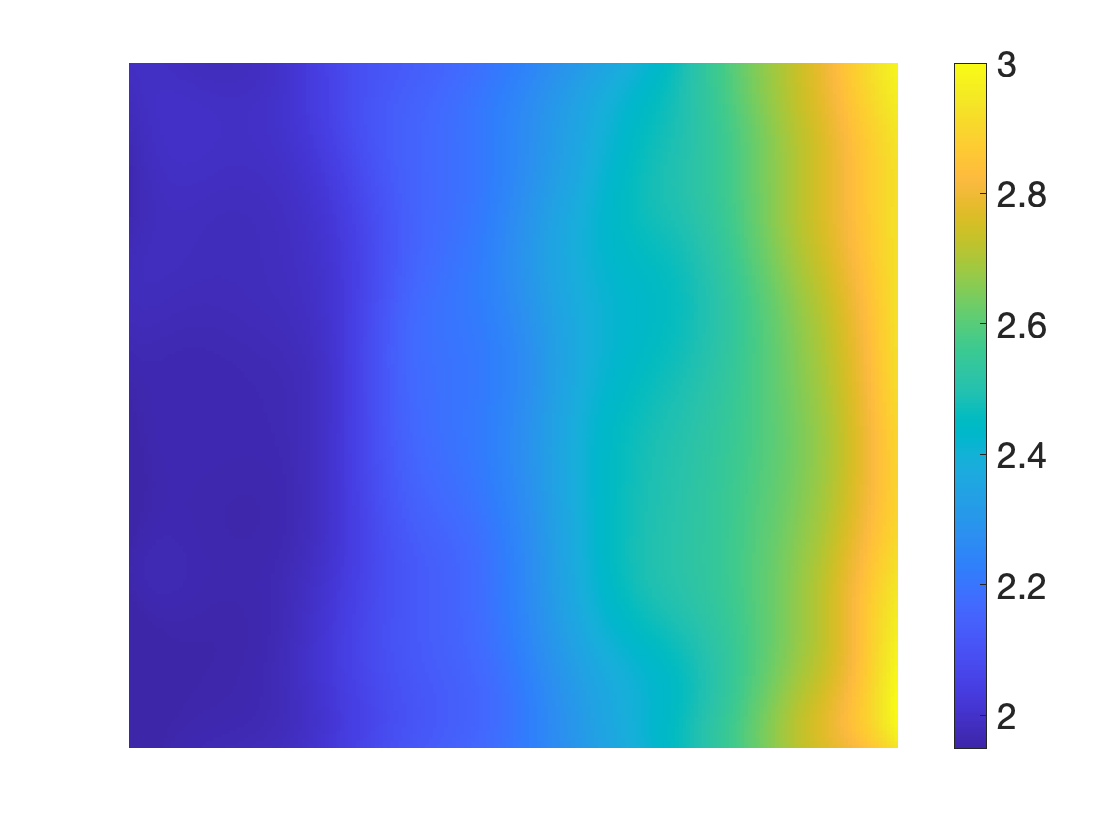} &
\includegraphics[width=0.199\textwidth]{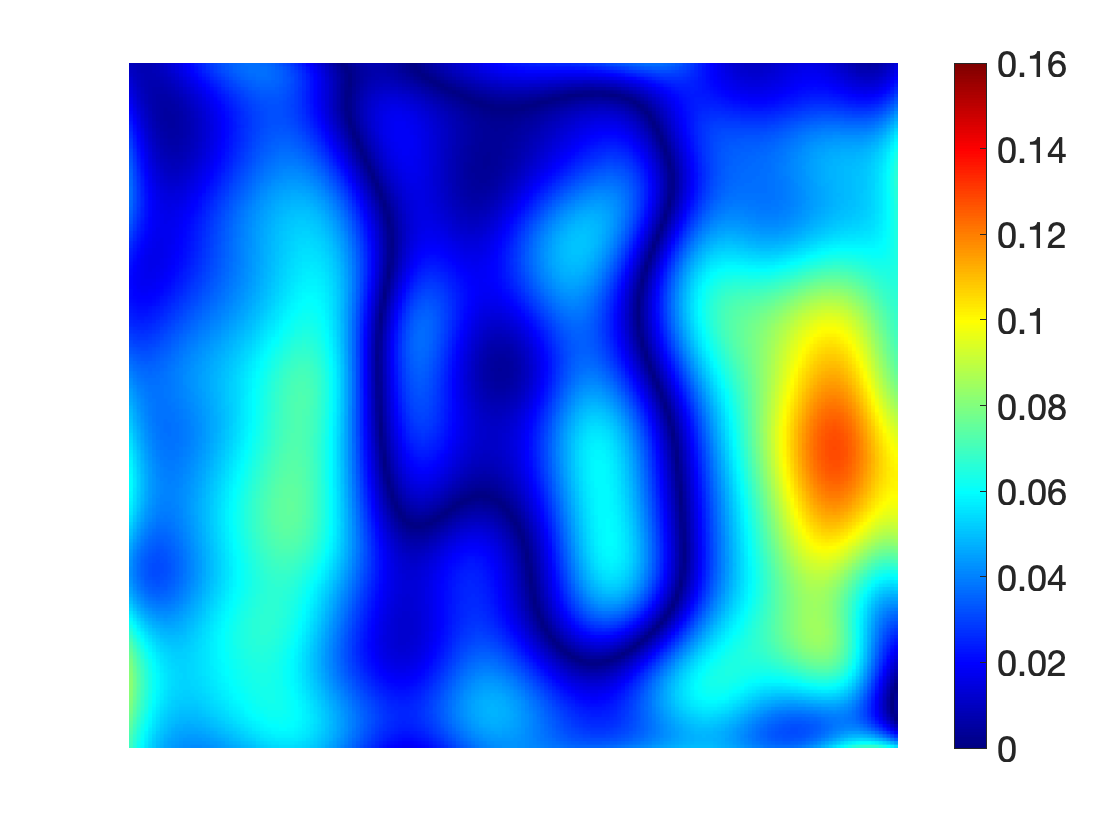}\\
\includegraphics[width=0.199\textwidth]{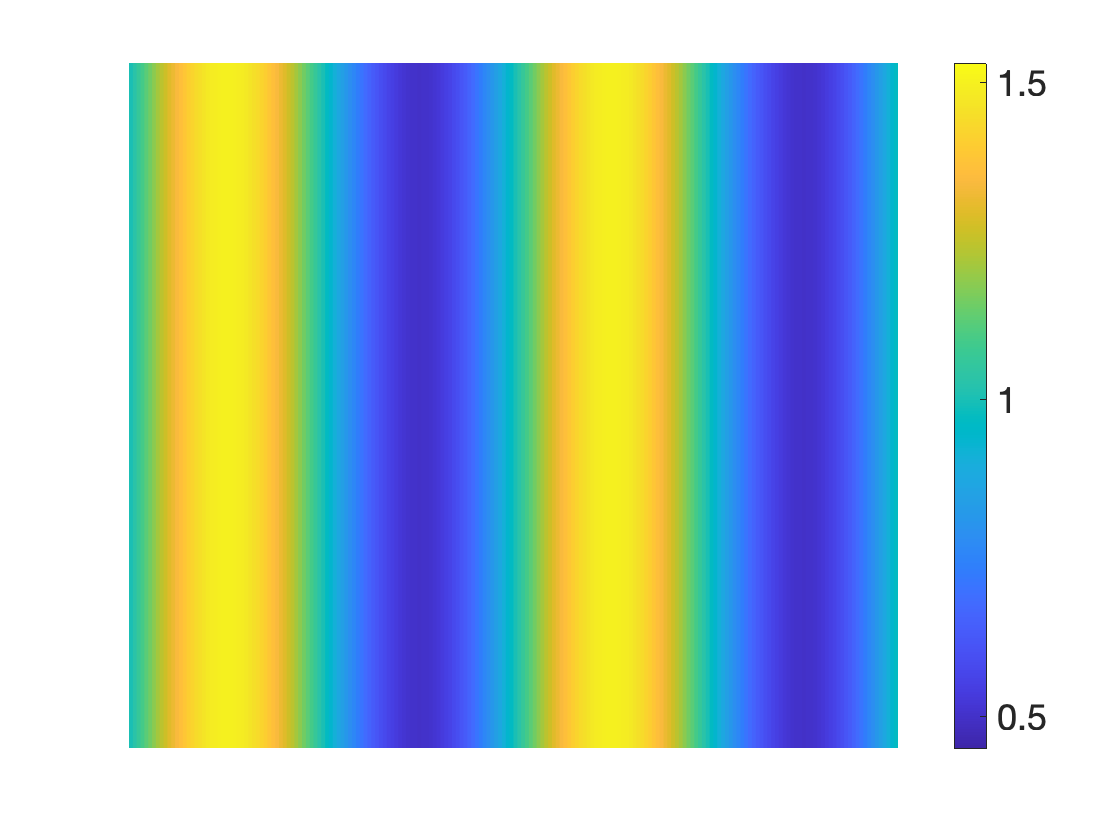} &
\includegraphics[width=0.199\textwidth]{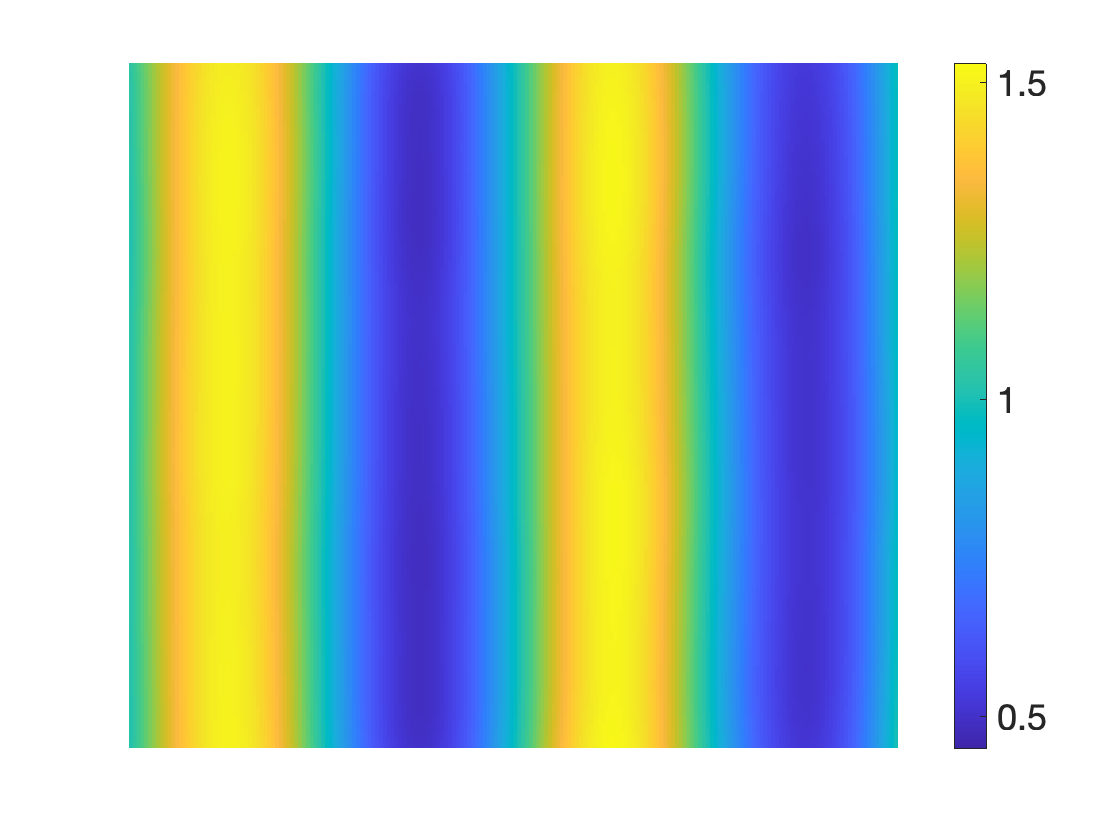} &
\includegraphics[width=0.199\textwidth]{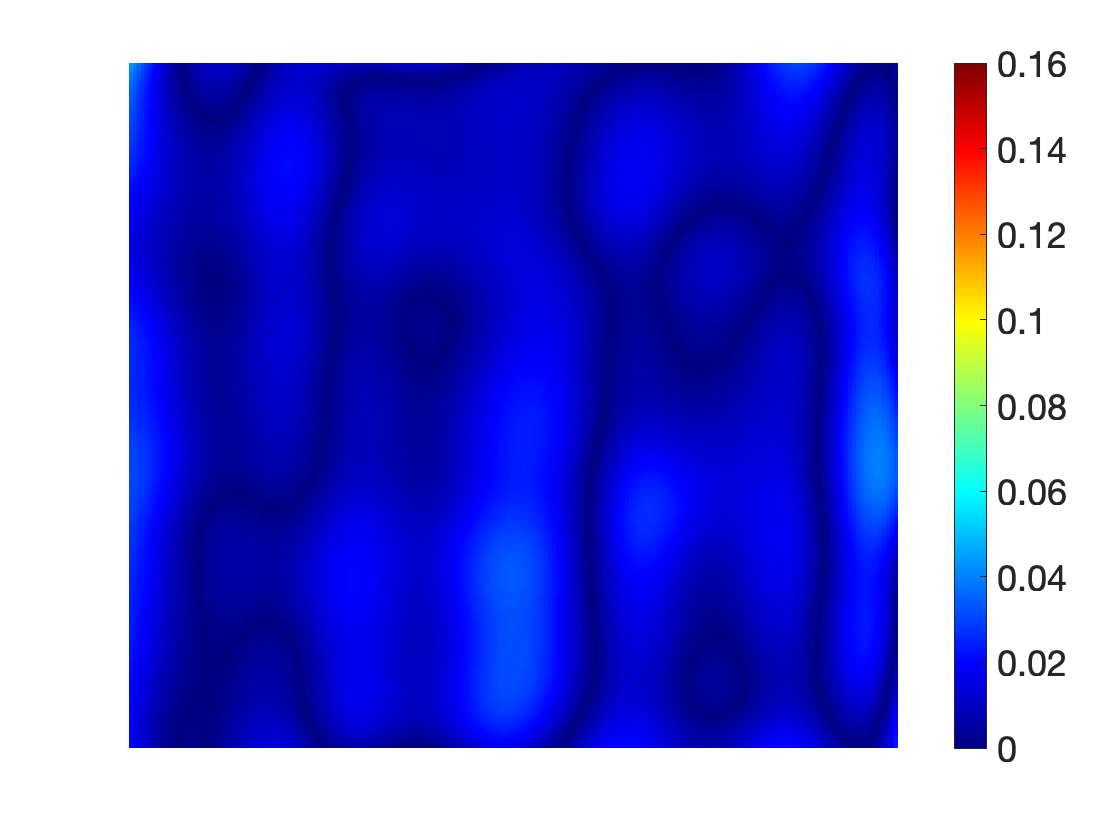} &
\includegraphics[width=0.199\textwidth]{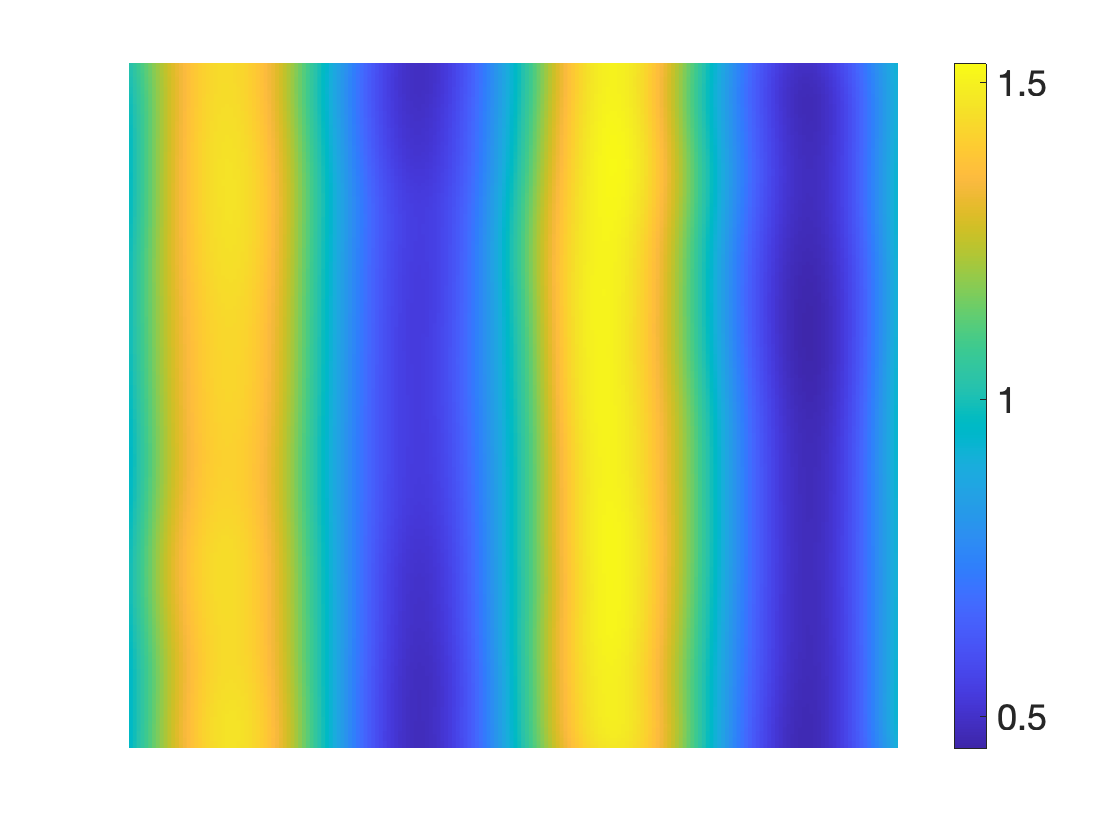} &
\includegraphics[width=0.199\textwidth]{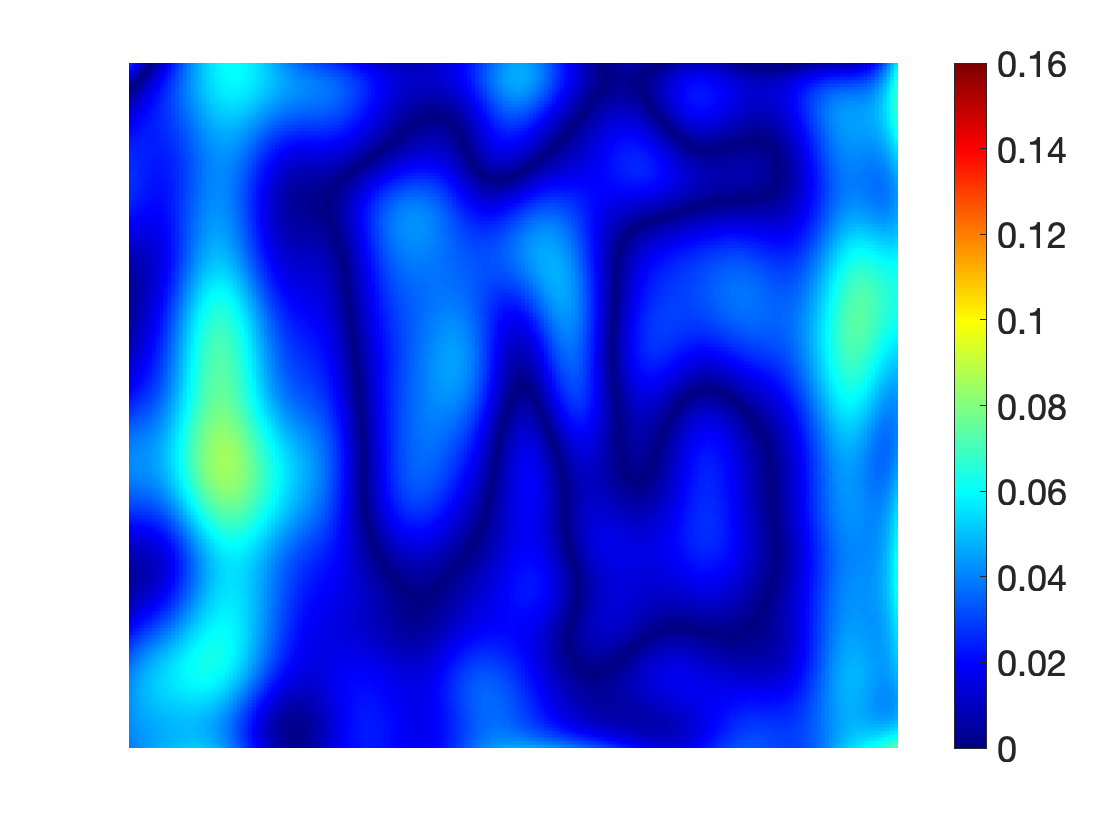}\\
\includegraphics[width=0.199\textwidth]{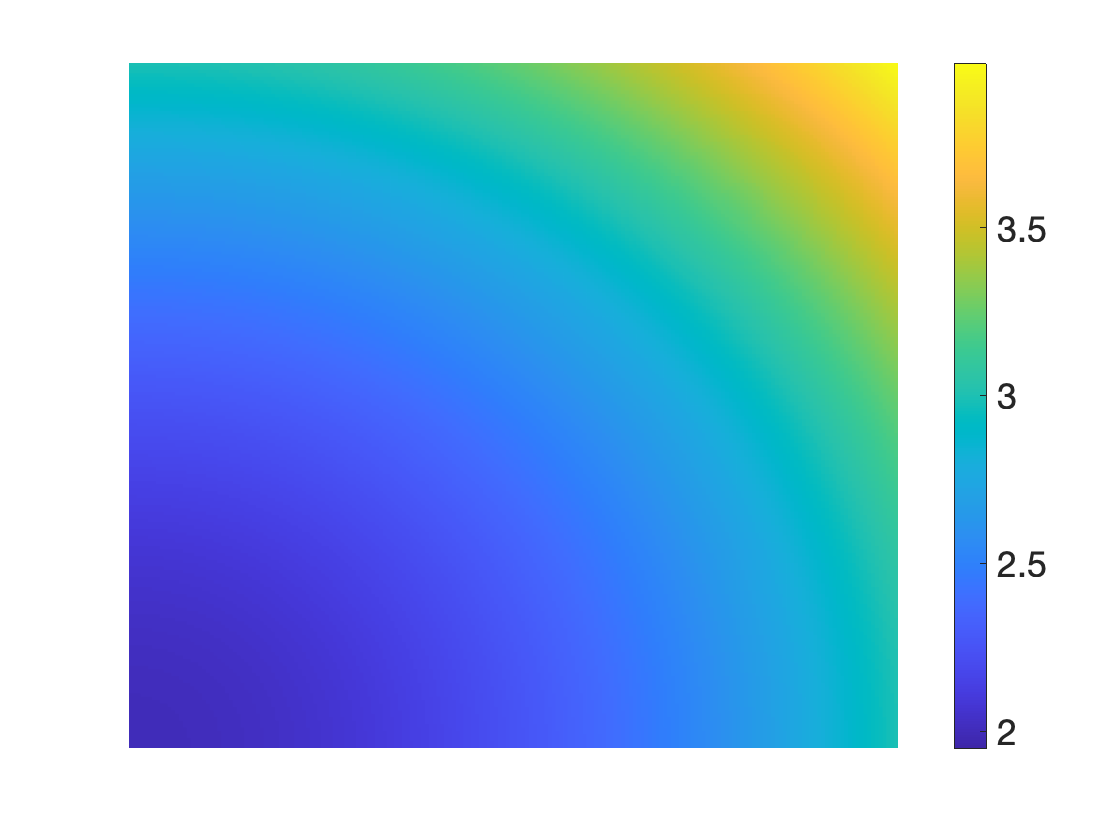} &
\includegraphics[width=0.199\textwidth]{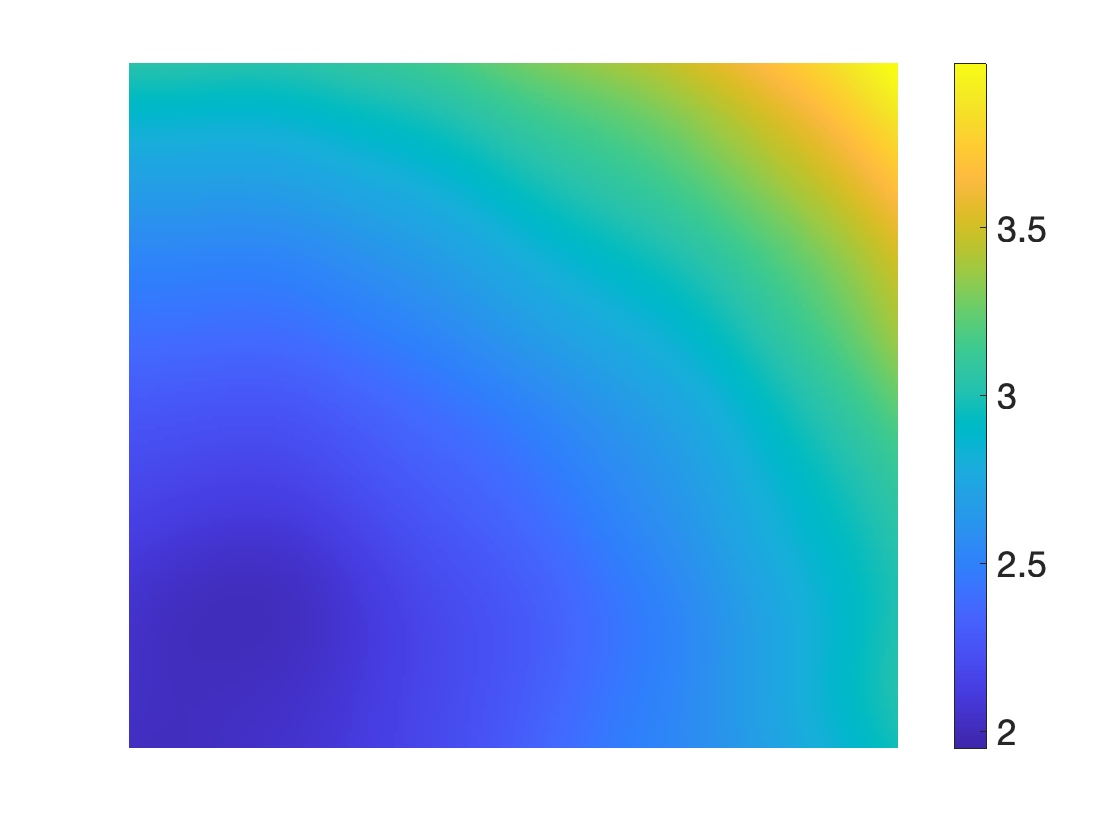} &
\includegraphics[width=0.199\textwidth]{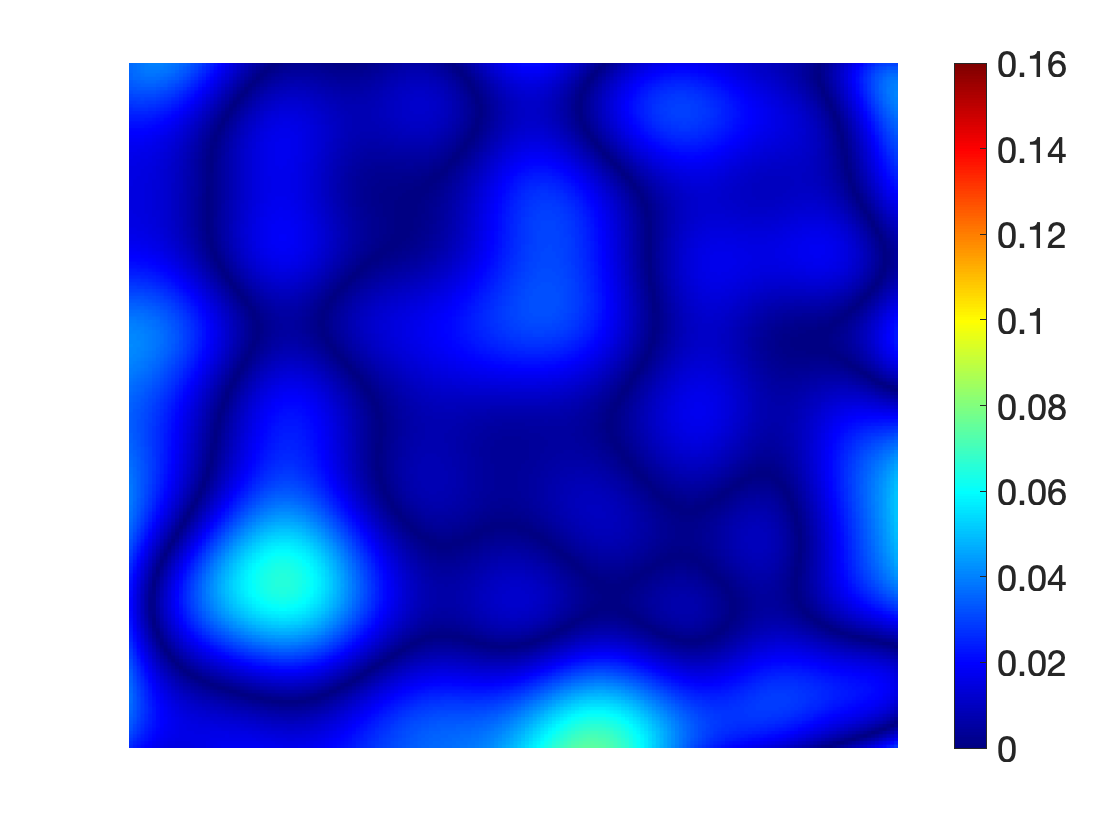} &
\includegraphics[width=0.199\textwidth]{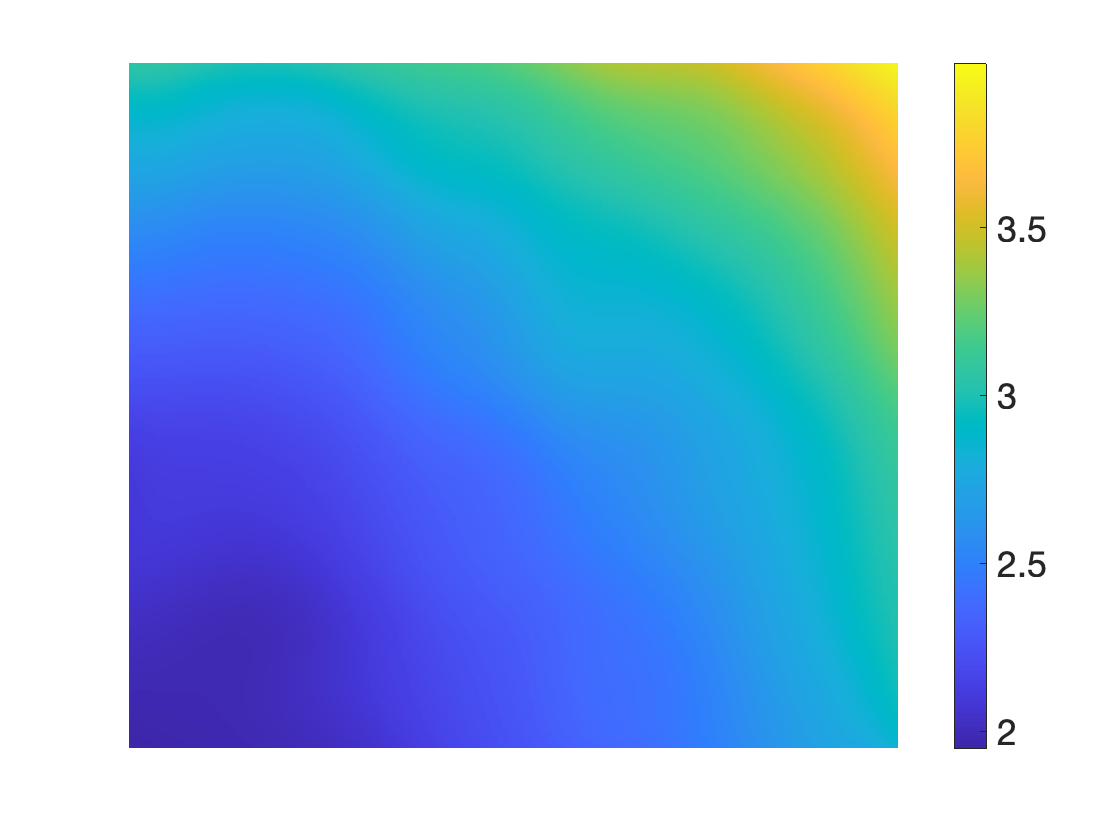} &
\includegraphics[width=0.199\textwidth]{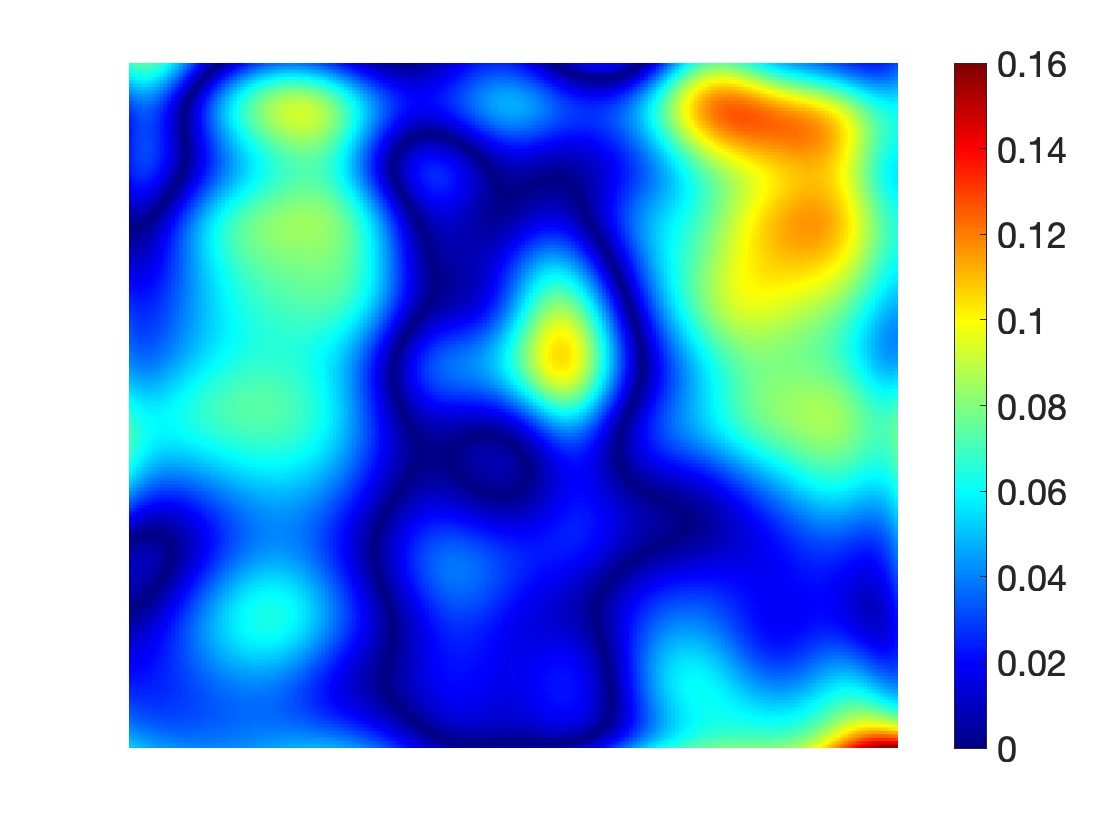}\\
(a) $A^\dag$  & (b) $\hat A$ & (c) $|\hat A-A^\dag|$ & (d) $\hat A$ & (e) $|\hat A-A^\dag|$
\end{tabular}
\caption{The reconstructions for Example \ref{exam:neu3d2} with exact data in (b) and noisy data $(\delta=5\%)$ in (d). From the top to bottom, the results are for $A_{11}$, $A_{12}$, $A_{13}$, $A_{22}$, $A_{23}$ and $A_{33}$, respectively.}
\label{fig:neu3d2}
\end{figure}

Fig. \ref{fig:neu3d2} shows the reconstruction results on a 2D cross section at $x_3=0.5$. Note that the reconstruction is accurate over the whole domain $\Omega$,
including the central region $\Omega\setminus \omega$ where there is no observational data,
and the quality does not deteriorate much for up to $5\%$ noise. This observation agrees well with that of the two-dimensional case in Example \ref{exam:neu2d3}. 

\subsection{The Dirichlet problem} Now we present a set of six numerical experiments for
the Dirichlet case. The first example is about recovering an anisotropic conductivity matrix with polynomial entries.
\begin{example}\label{exam:diri2d1}
    The domain $\Omega = (0,1)^2$, $A^\dag= \begin{pmatrix}
    2+x_1^2+x_2^2&1+(x_1-\frac12)(x_2-\frac12)\\
    1+(x_1-\frac12)(x_2-\frac12)&2+(x_1-\frac12)^2+(x_2-\frac12)^2\\
\end{pmatrix},$ $u_1^\dag=x_1+x_2+\frac{1}{3}(x_1^3+x_2^3)$, $u_2^\dag=x_1-x_2+\frac{1}{3}(x_1^3-x_2^3)$, $u_3^\dag=-u_1^\dagger$ and $u_4^\dag=-u_2^\dagger$.
\end{example}

\begin{figure}[htb!]
\centering
\setlength{\tabcolsep}{0em}
\begin{tabular}{ccccc}
\includegraphics[width=0.199\textwidth]{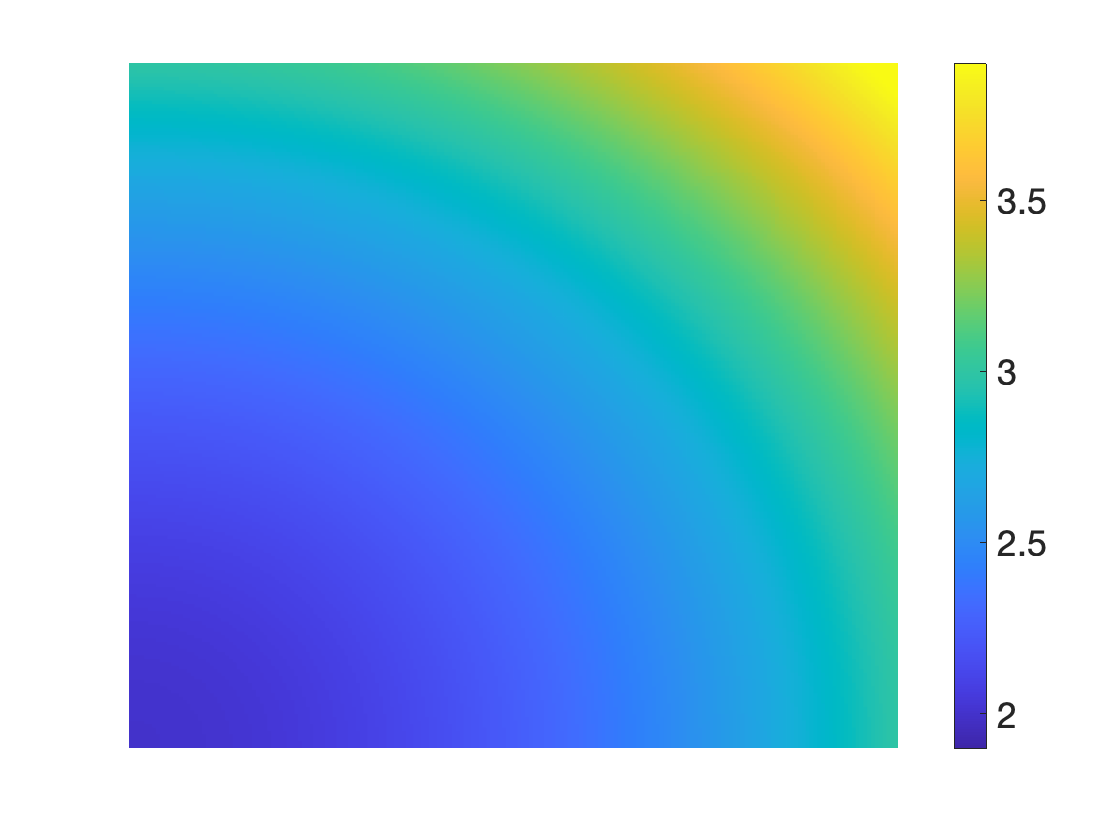} &
\includegraphics[width=0.199\textwidth]{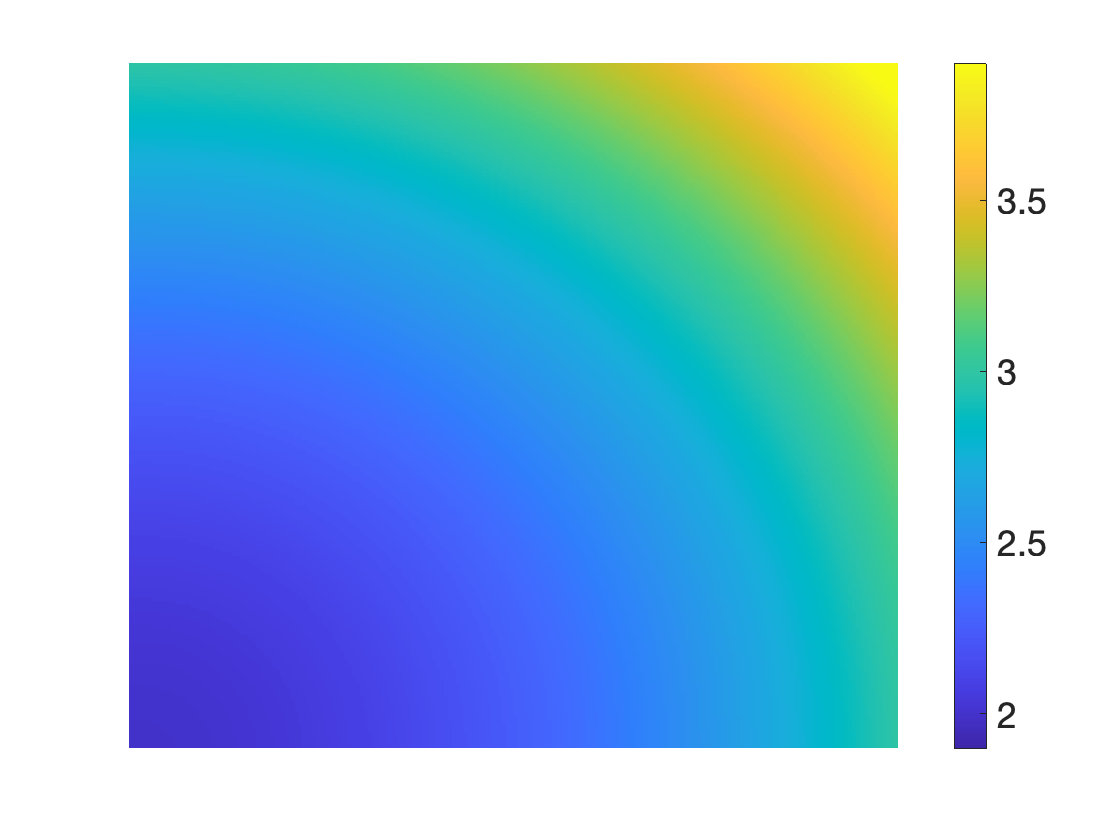} &
\includegraphics[width=0.199\textwidth]{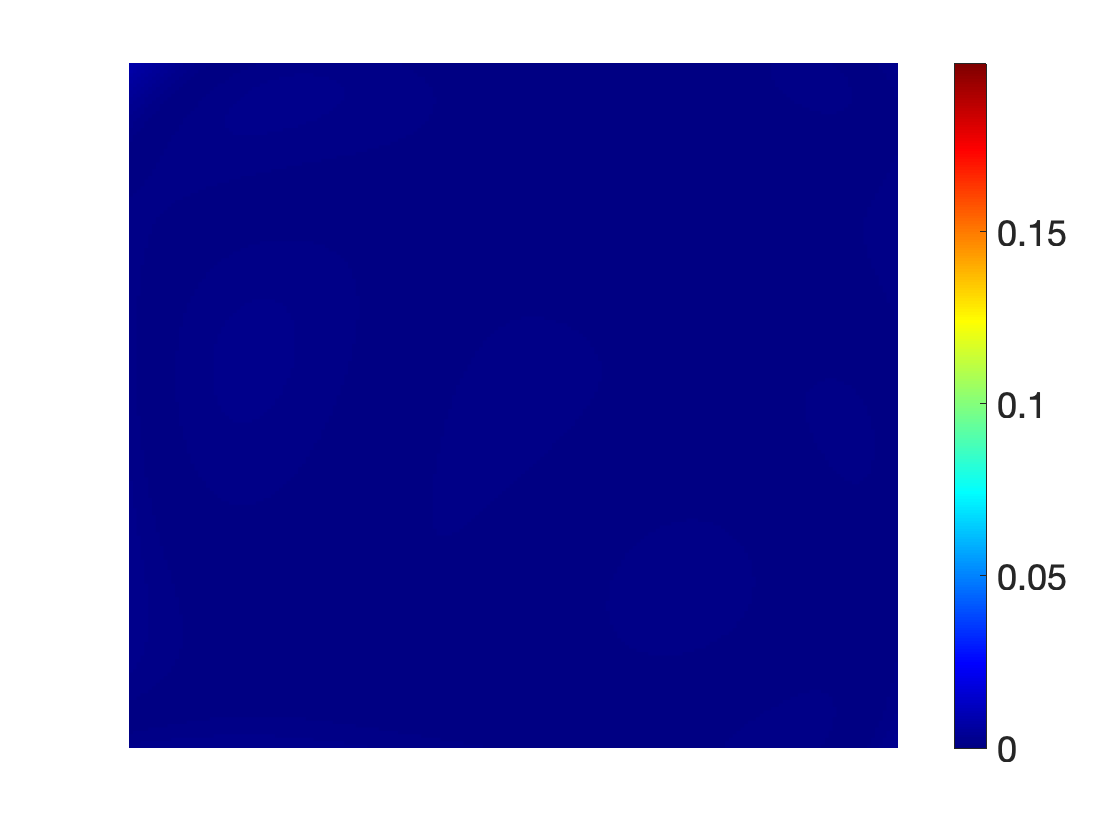} &
\includegraphics[width=0.199\textwidth]{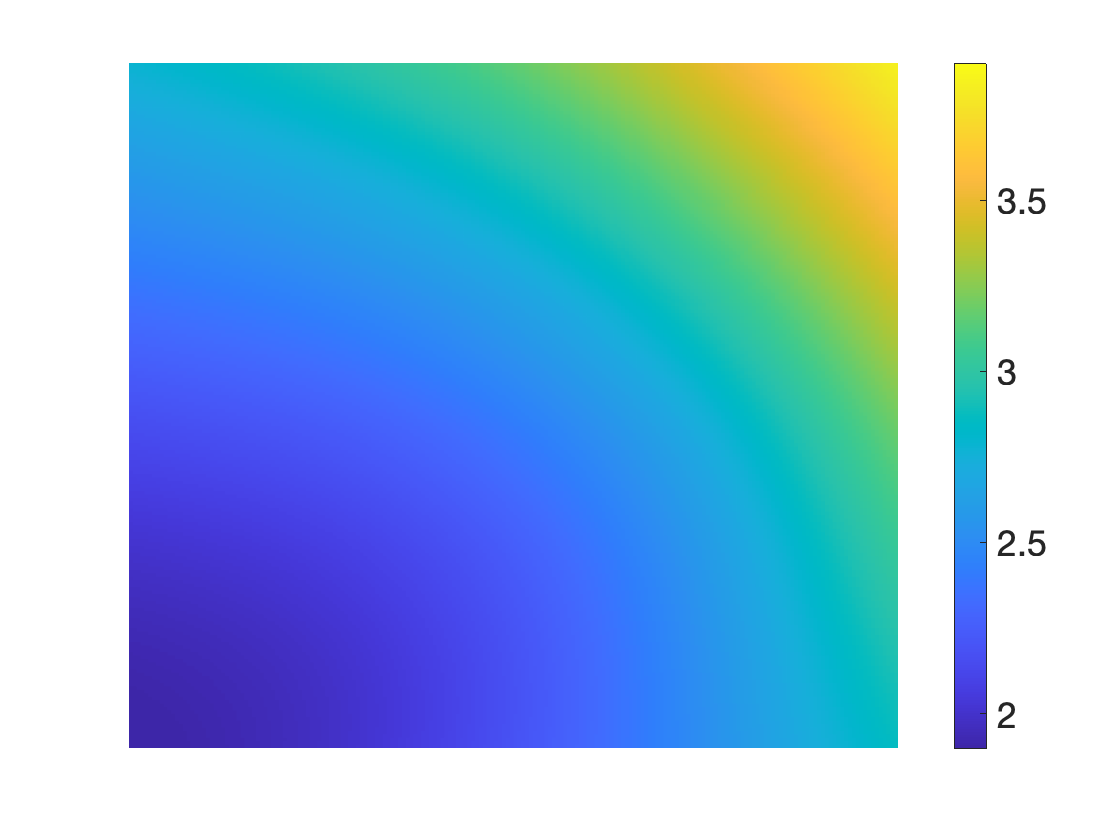} &
\includegraphics[width=0.199\textwidth]{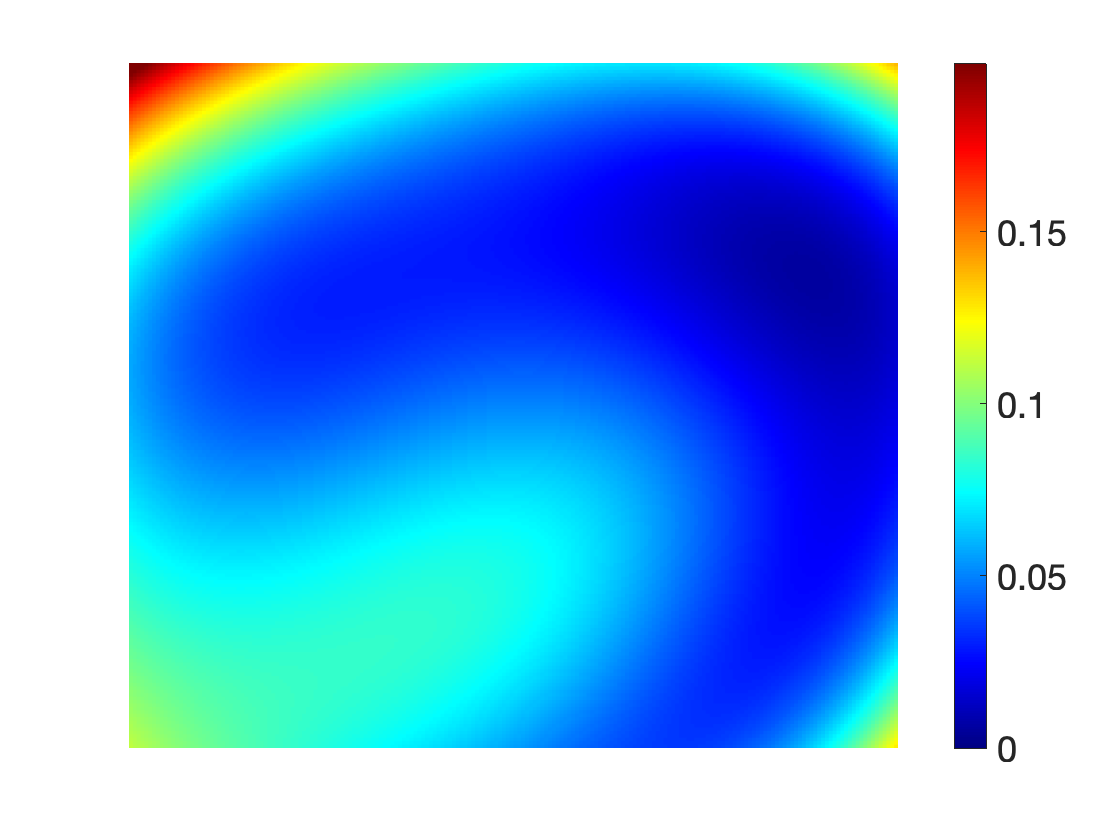} \\
\includegraphics[width=0.199\textwidth]{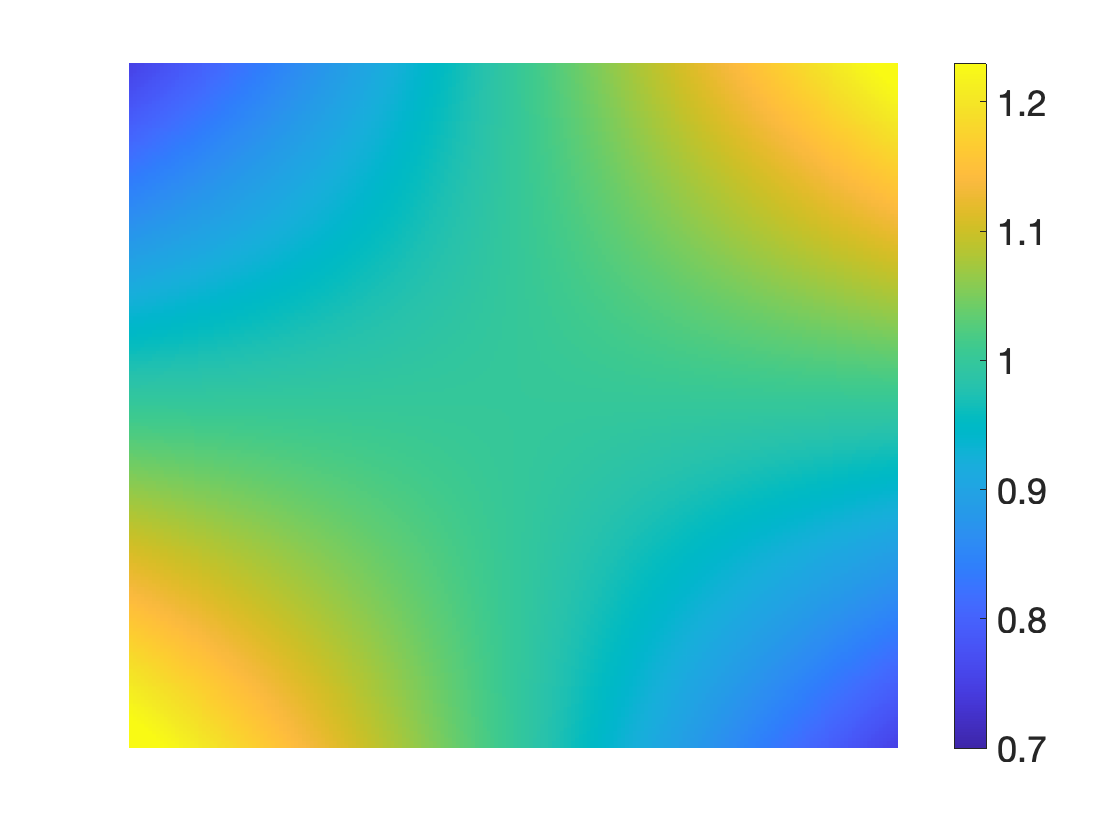} &
\includegraphics[width=0.199\textwidth]{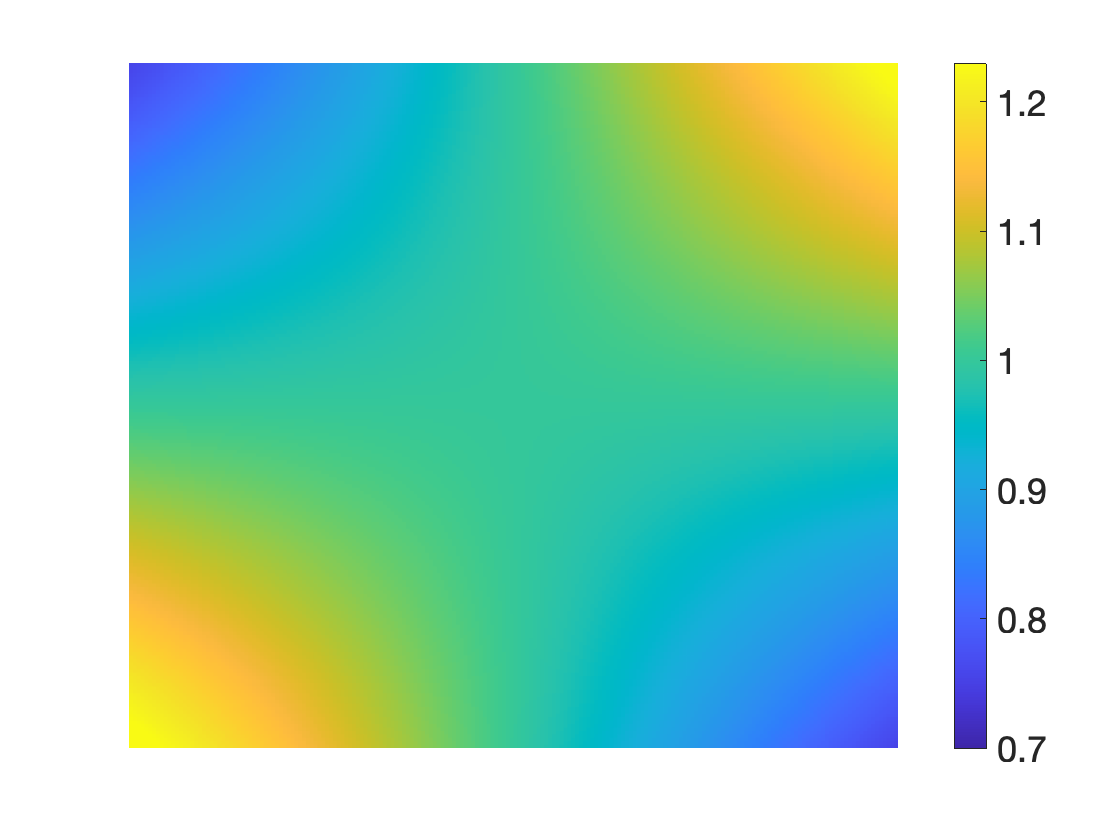} &
\includegraphics[width=0.199\textwidth]{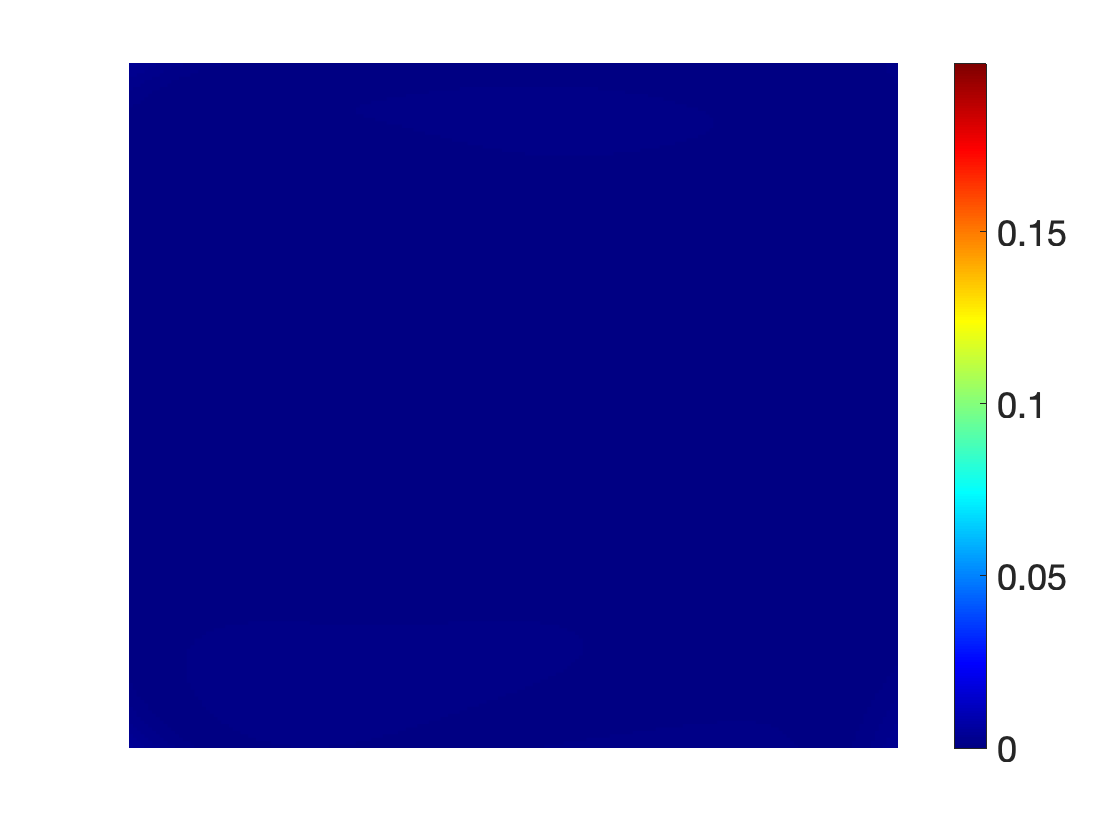} &
\includegraphics[width=0.199\textwidth]{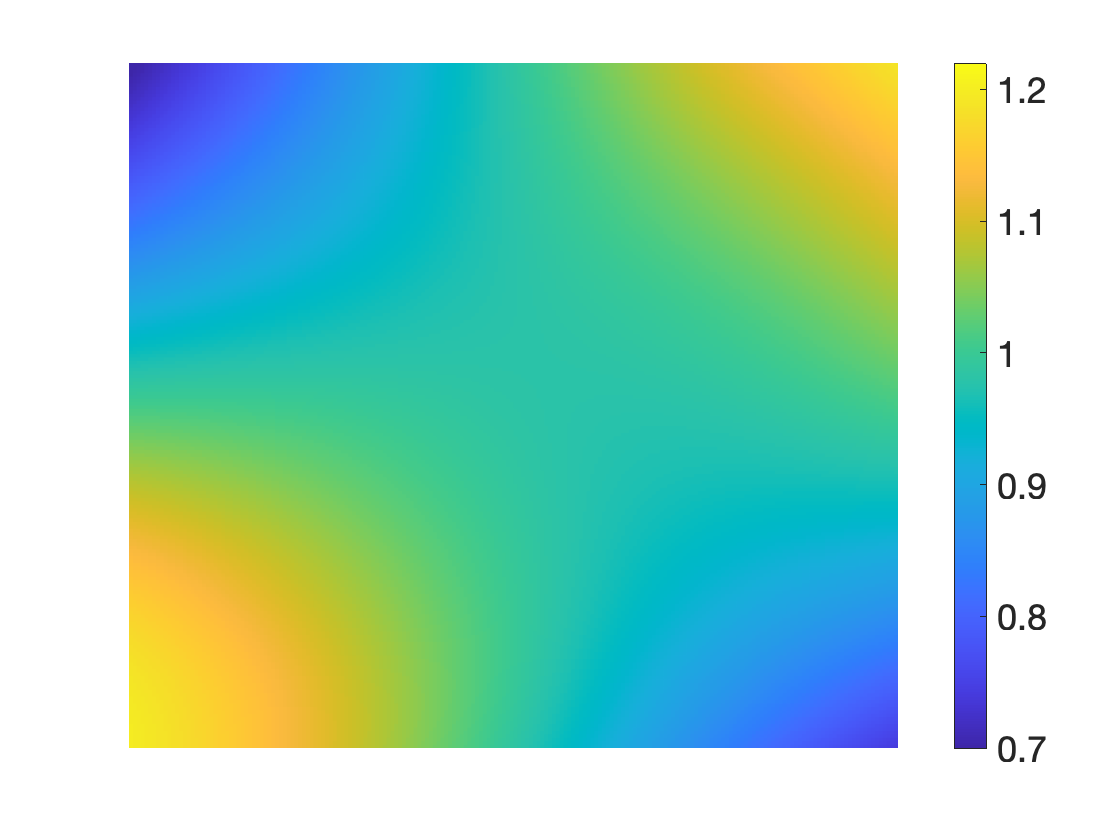} &
\includegraphics[width=0.199\textwidth]{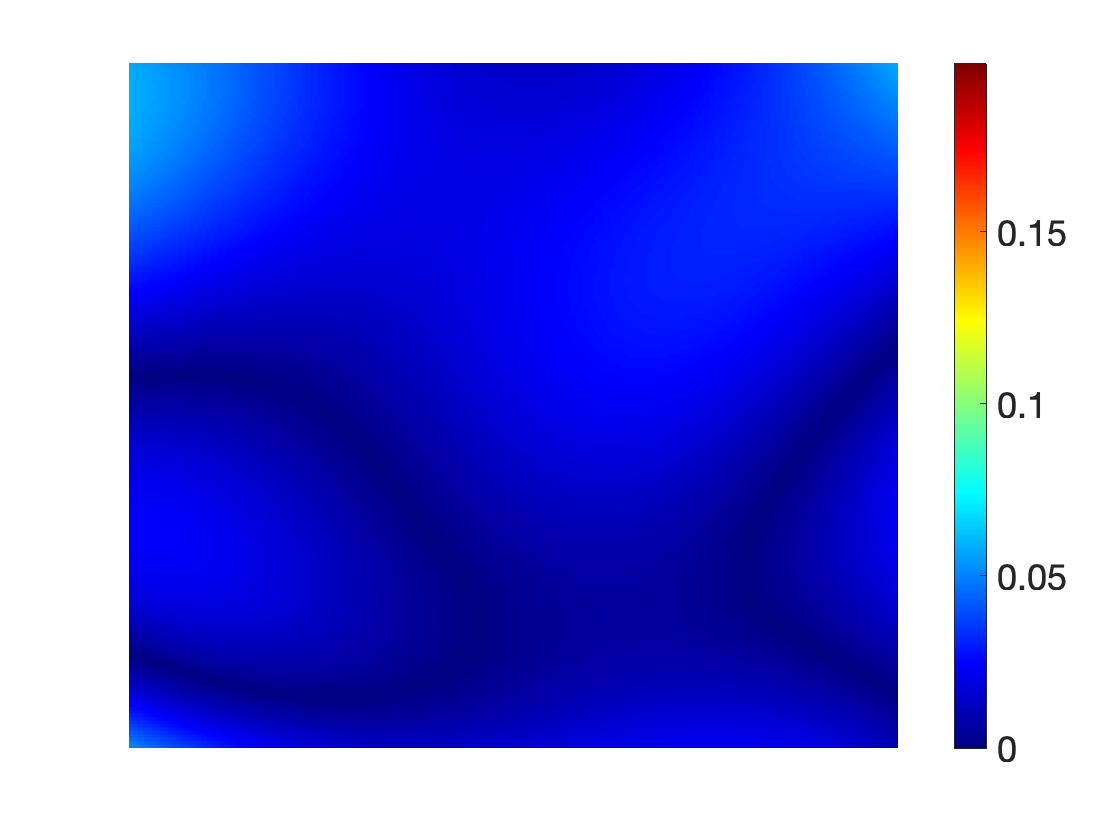} \\
\includegraphics[width=0.199\textwidth]{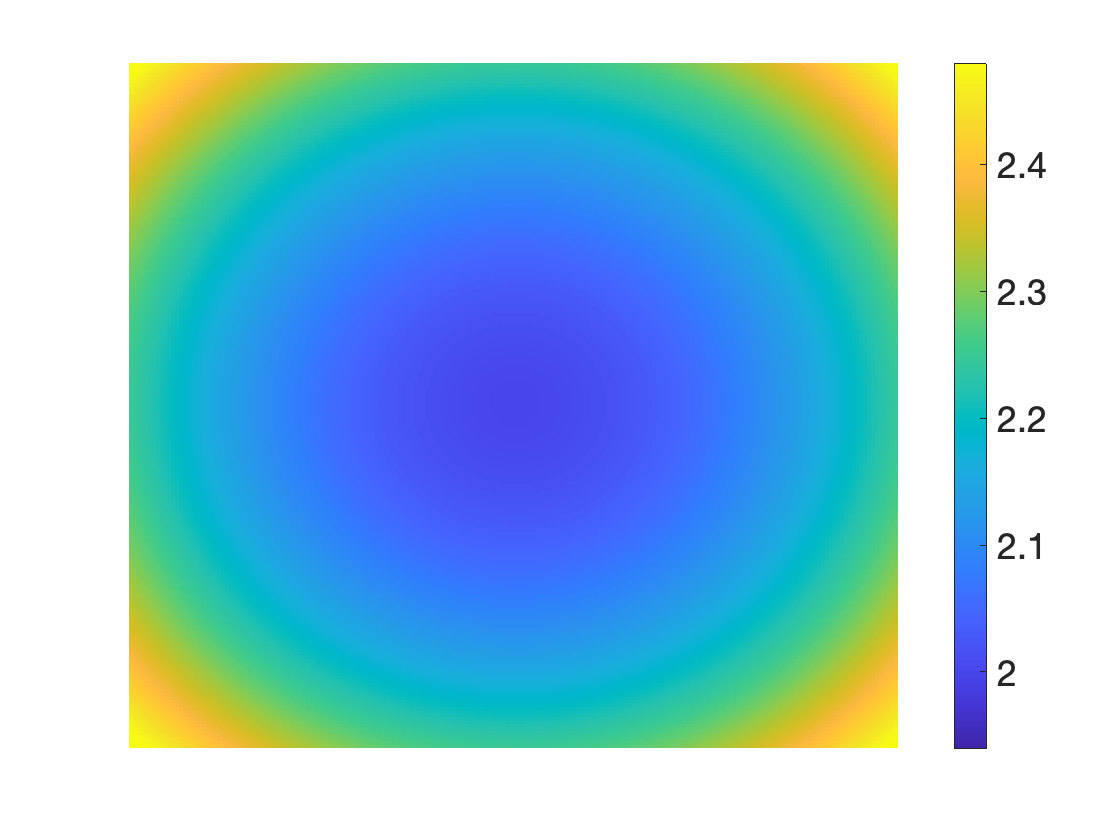} &
\includegraphics[width=0.199\textwidth]{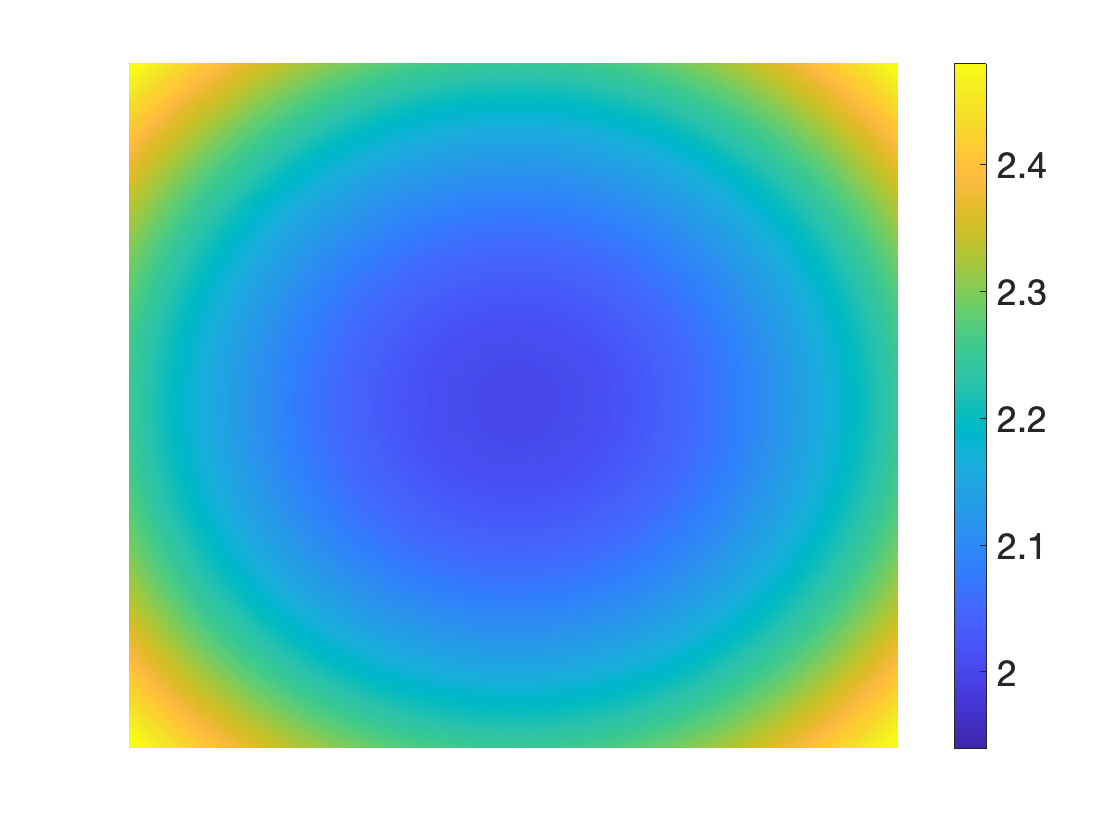} &
\includegraphics[width=0.199\textwidth]{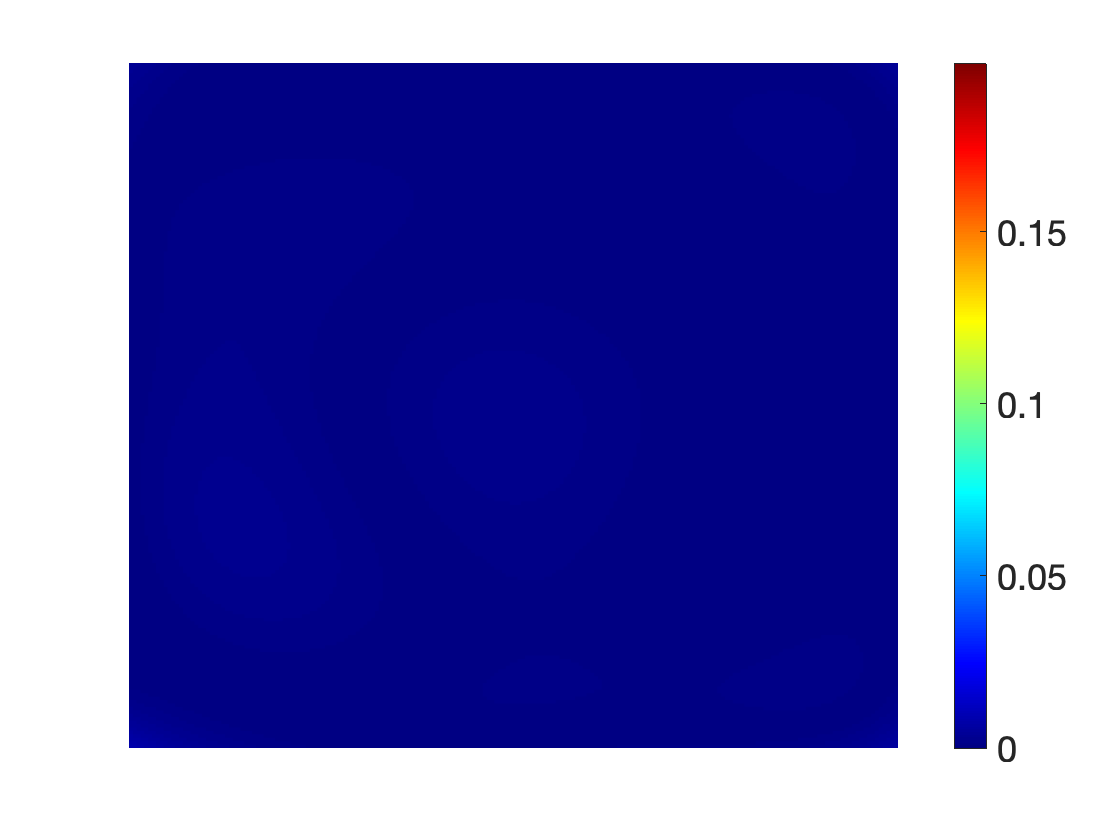} &
\includegraphics[width=0.199\textwidth]{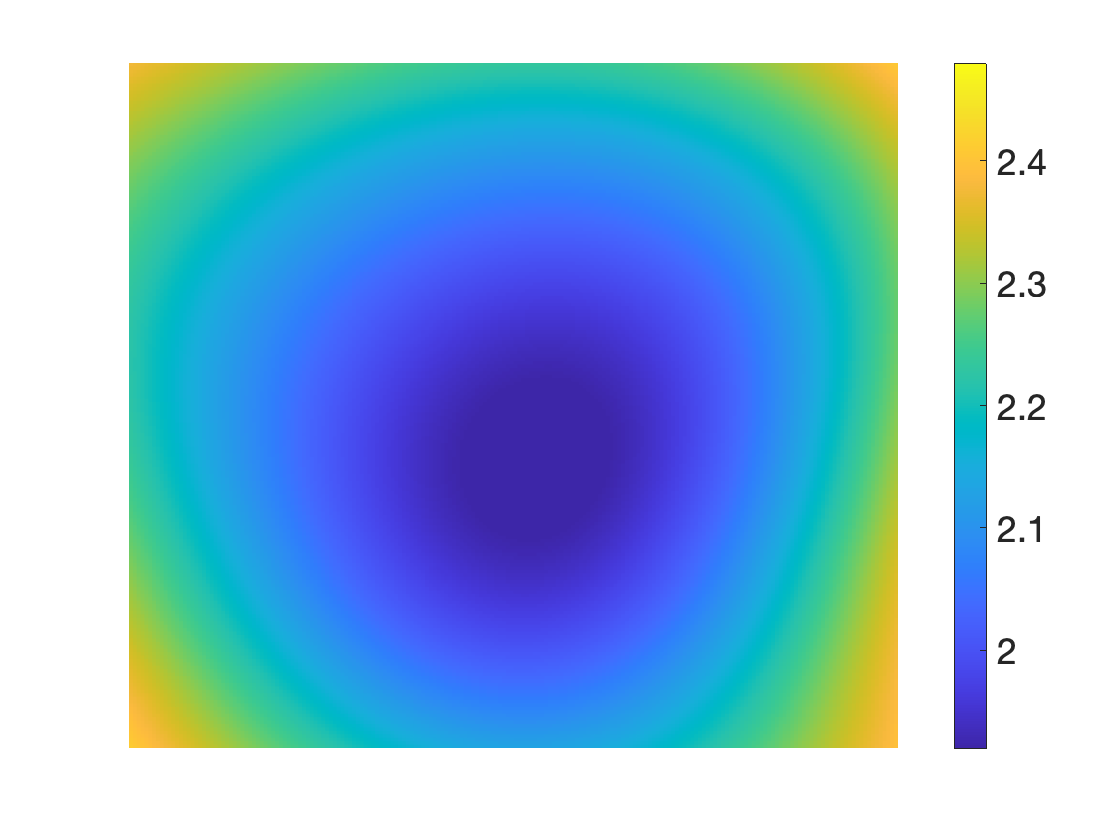} &
\includegraphics[width=0.199\textwidth]{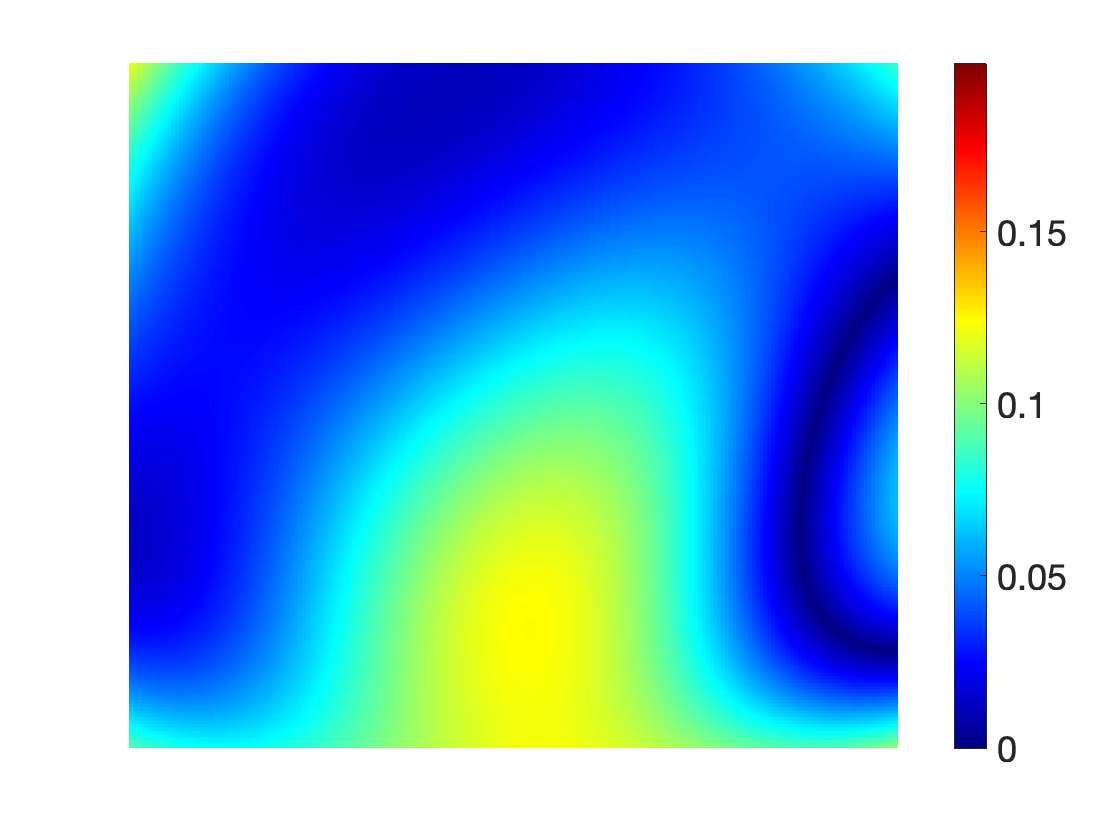} \\
(a) $A^\dag$  & (b) $\hat A$ & (c) $|\hat A-A^\dag|$ & (d) $\hat A$ & (e) $|\hat A-A^\dag|$
\end{tabular}
\caption{The reconstructions for Example \ref{exam:diri2d1} with exact data in (b) and noisy data $(\delta=10\%)$ in (d). From the top to bottom, the results are for $A_{11}$, $A_{12}$ and $A_{22}$, respectively.}
\label{fig:diri2d1}
\end{figure}

\begin{figure}[htb!]
\centering
\setlength{\tabcolsep}{0em}
\begin{tabular}{ccc}
\includegraphics[width=0.32\textwidth]{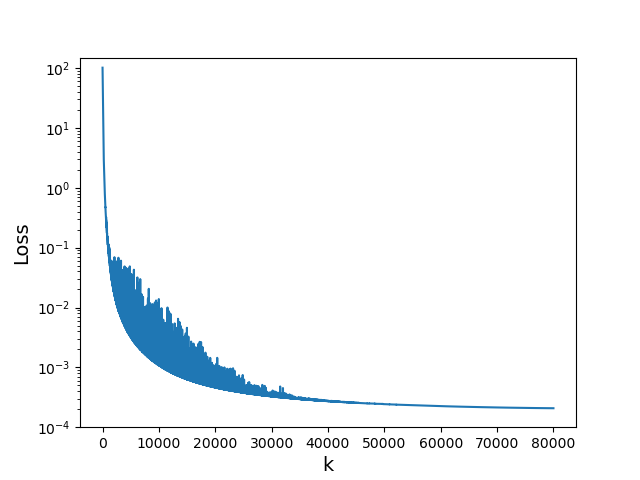} &
\includegraphics[width=0.32\textwidth]{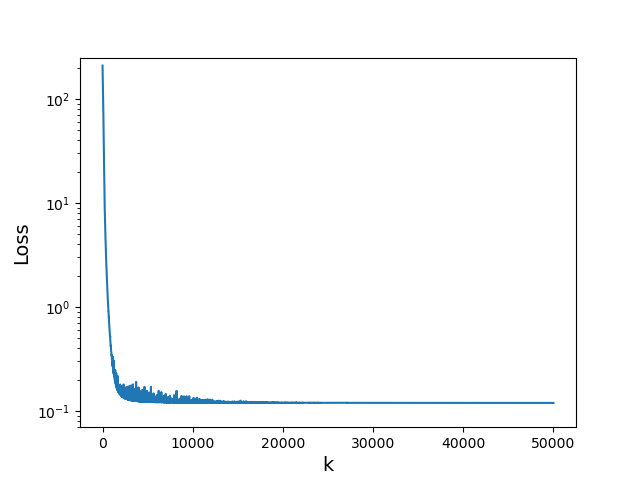} &
\includegraphics[width=0.32\textwidth]{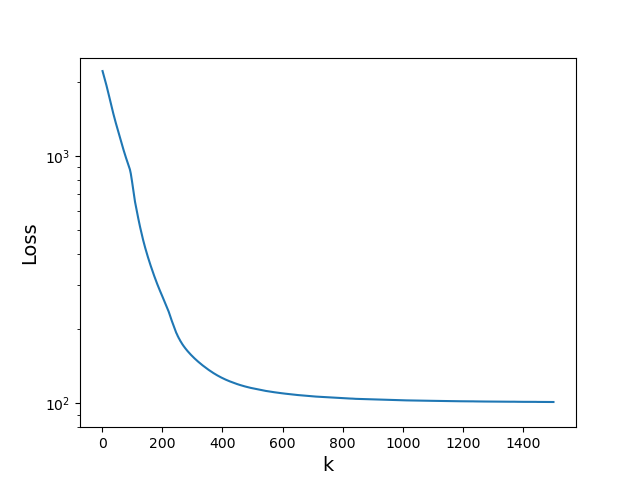}\\
\includegraphics[width=0.32\textwidth]{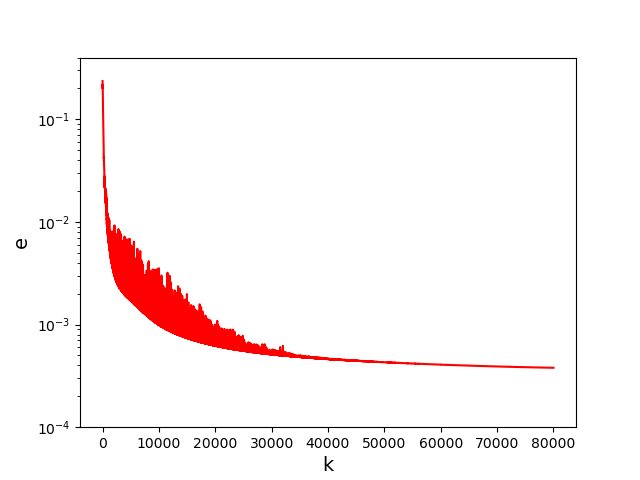} &
\includegraphics[width=0.32\textwidth]{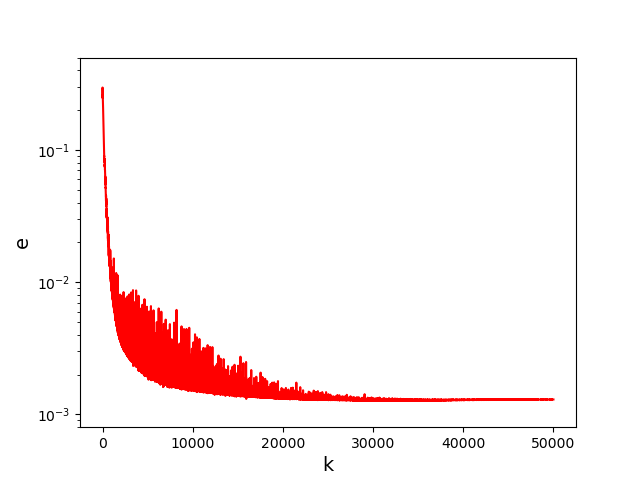} &
\includegraphics[width=0.32\textwidth]{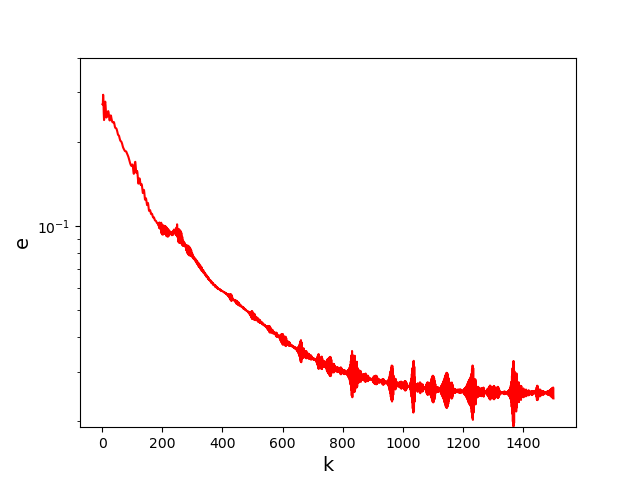}\\
(a) $\delta=0\%$  & (b) $\delta=1\%$ & (c) $\delta=10\%$
\end{tabular}
\caption{The evolution of the loss $($top$)$ and the reconstruction error $e(\hat A)$ $($bottom$)$ during the training process for Example \ref{exam:diri2d1} at different noise levels.}
\label{fig:diri1losse}
\end{figure}

Fig. \ref{fig:diri2d1} presents the results for both exact and noisy data. Like before, the DNN approach enjoys excellent robustness to data noise. With $\delta=10\%$, there is minimal degradation in reconstruction quality, except for a slight deformation in the conductivity landscape of entry $A_{22}$. Tables \ref{dirit1}-\ref{dirit3} illustrate the impact of the regularization parameter $\gamma_A$, DNN architecture, and the numbers of random sampling points on reconstruction quality. These parameters do not significantly affect the reconstruction error as long as they are within a suitable range, showing the robustness of the approach. The convergence behavior of the optimizer in Fig. \ref{fig:diri1losse} shows that both the loss $\widehat{J}_{\bsgamma}$ and error $e(\hat{A})$ stagnate at comparable levels regardless of the noise level $\delta$, further corroborating the robustness of the approach.

\begin{table}[htp!]
  \centering
  \caption{The variation of the reconstruction error $e(\hat{A})$ with respect to various algorithmic parameters. }
\begin{threeparttable}
\subfigure[$e$ v.s. $\gamma_A$ and $\delta$\label{dirit1}]{\begin{tabular}{c|ccc}
\toprule
  $\gamma_A\backslash\delta$&   0\% &  1\%& 10\%\\
\midrule
     1.00e-2 &  2.76e-3 & 3.16e-3 & 2.68e-2 \\
     1.00e-3 &  4.77e-4 & 1.35e-3 & 2.62e-2\\
     1.00e-4 &  3.80e-4 & 1.32e-3 & 2.87e-2\\
     1.00e-5 &  3.80e-4 & 1.30e-3 & 2.41e-2\\
\bottomrule
\end{tabular}} \\
\subfigure[$e$ v.s. $W_A$ and  $L_A$\label{dirit2}]{
\begin{tabular}{c|ccc}
\toprule
${W_A}\backslash L_A$&   5 &  10& 15\\
\midrule
     8 &1.12e-3&6.91e-4&5.66e-4 \\
     16&5.16e-4&4.71e-4&4.69e-4\\
     24&4.29e-4&4.21e-4&4.97e-4\\
     32&3.80e-4&3.98e-4&4.65e-4\\
\bottomrule
\end{tabular}}\quad
\subfigure[$e$ v.s. $n_b$ and $n_r$\label{dirit3}]{
\begin{tabular}{c|ccc}
\toprule
$n_b\backslash n_r$&   5000 &  10000& 15000\\
\midrule
     500&  4.81e-4&5.46e-4&4.23e-4\\
     1000& 4.41e-4&3.80e-4&3.76e-4\\
     2000& 3.87e-4&3.68e-4&4.54e-4\\
     4000& 5.17e-4&4.15e-4&4.57e-4\\

\bottomrule
\end{tabular}}
\end{threeparttable}
\end{table}

The second example is about recovering an anisotropic conductivity matrix with mixed oscillatory and polynomial entries.
\begin{example}
The domain $\Omega = (0,1)^2$, $A^\dag= \begin{pmatrix}
    2+\frac{\sin(4\pi x_2)}{2}&1+\frac{\sin^2(2\pi x_2)}{2}\\
    1+\frac{\sin^2(2\pi x_2)}{2}&2+x_1^2\\
\end{pmatrix},$ $u_1^\dag=x_1+x_2+\frac{1}{3}(x_1^3+x_2^3)$, $u_2^\dag=x_1-x_2+\frac{1}{3}(x_1^3-x_2^3)$, $u_3^\dag=-u_1^\dagger$, $u_4^\dag=-u_2^\dagger$.
\label{exam:diri2d2}
\end{example}

\begin{figure}[htb!]
\centering
\setlength{\tabcolsep}{0em}
\begin{tabular}{ccccc}
\includegraphics[width=0.199\textwidth]{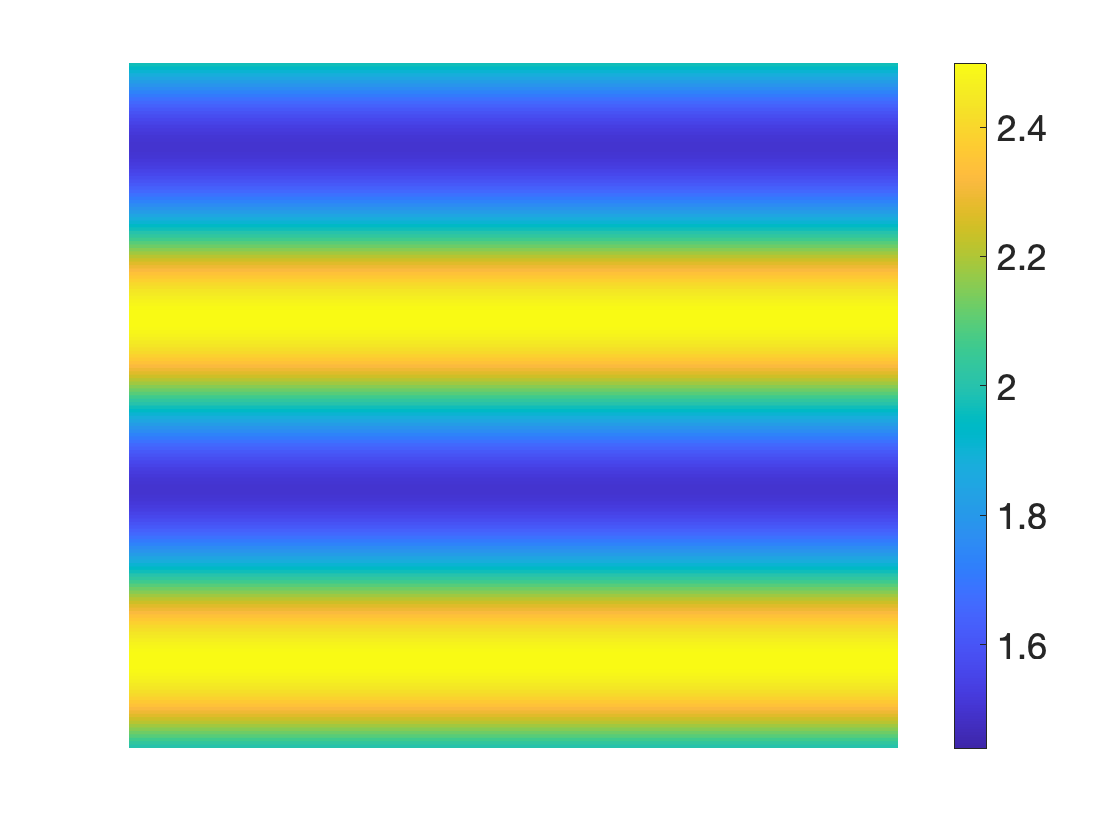} &
\includegraphics[width=0.199\textwidth]{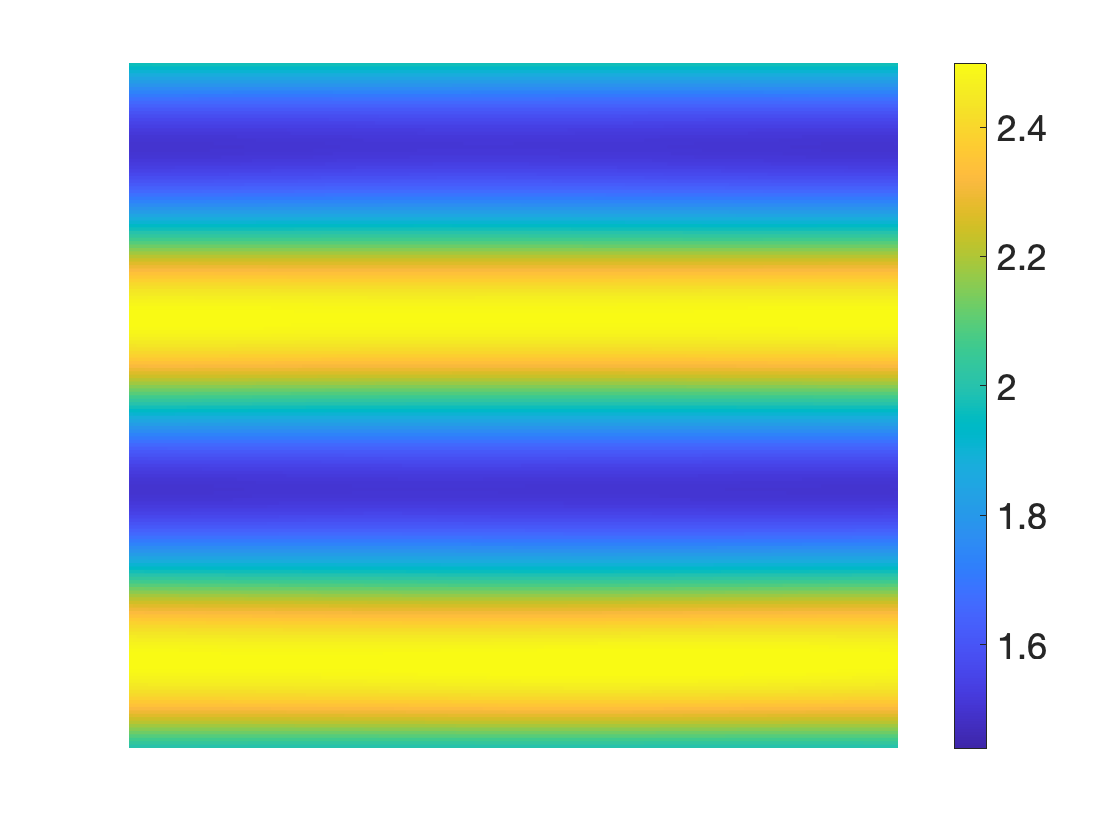} &
\includegraphics[width=0.199\textwidth]{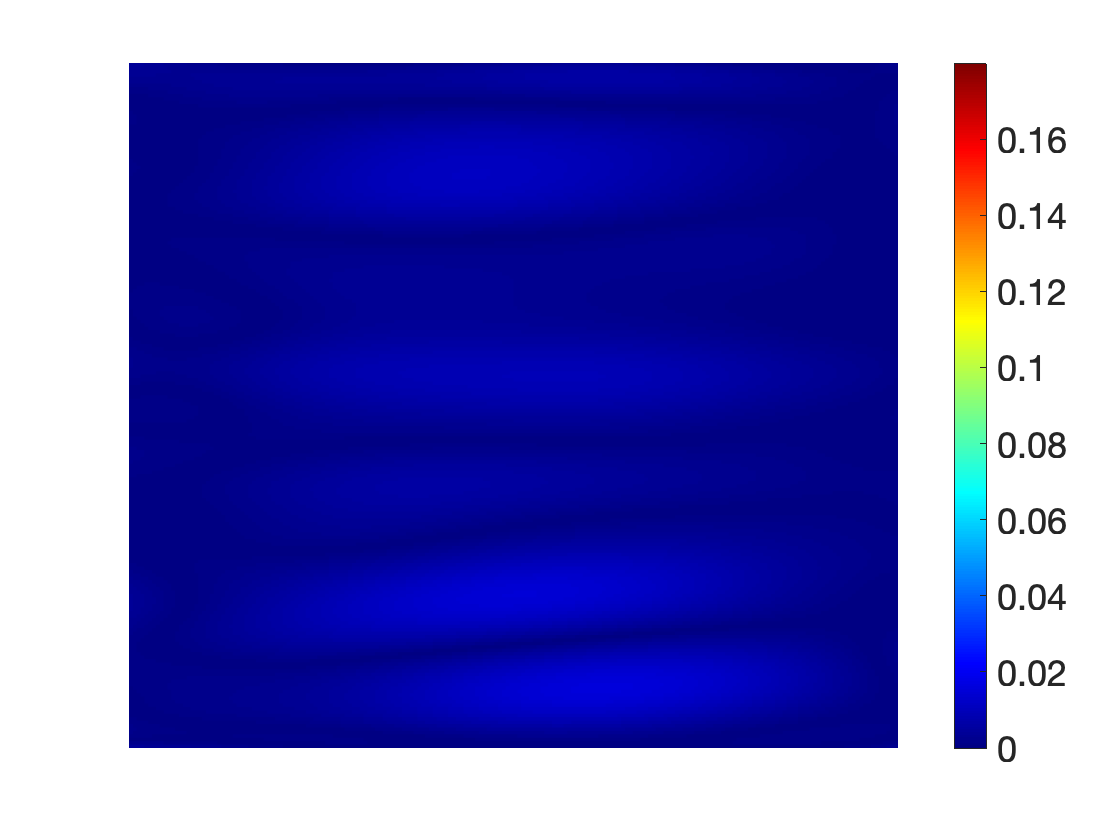} &
\includegraphics[width=0.199\textwidth]{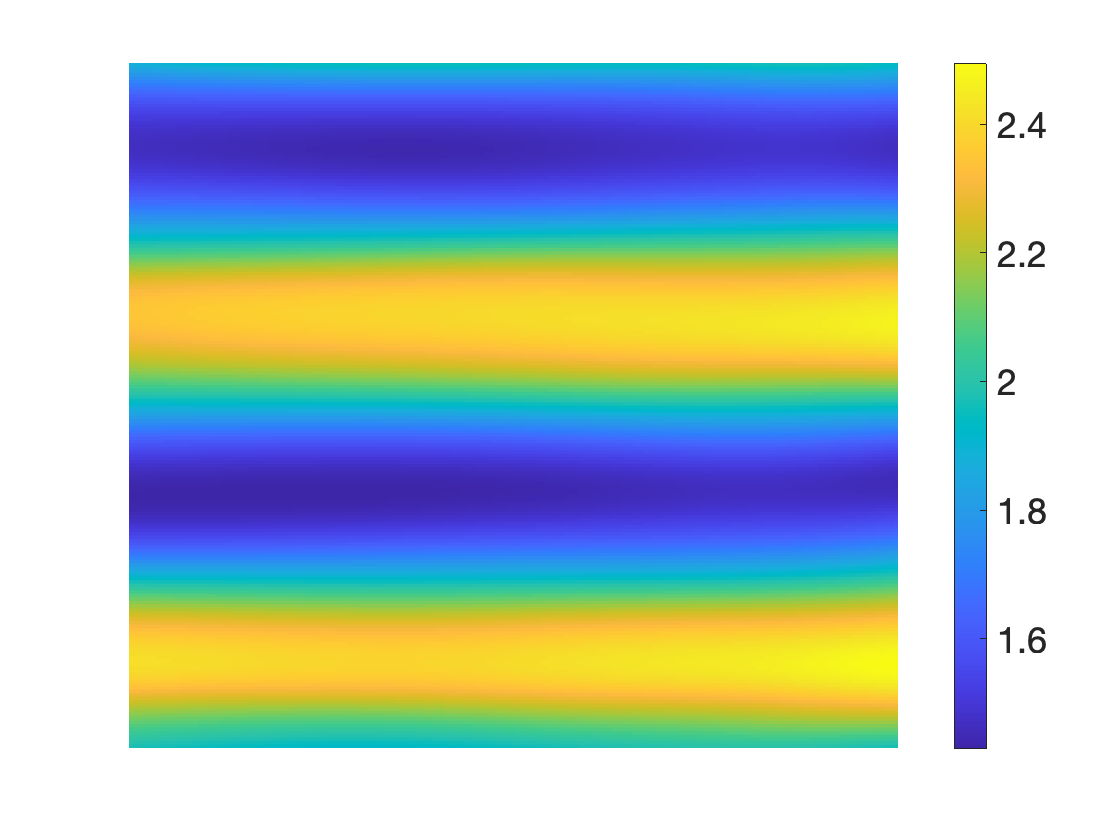} &
\includegraphics[width=0.199\textwidth]{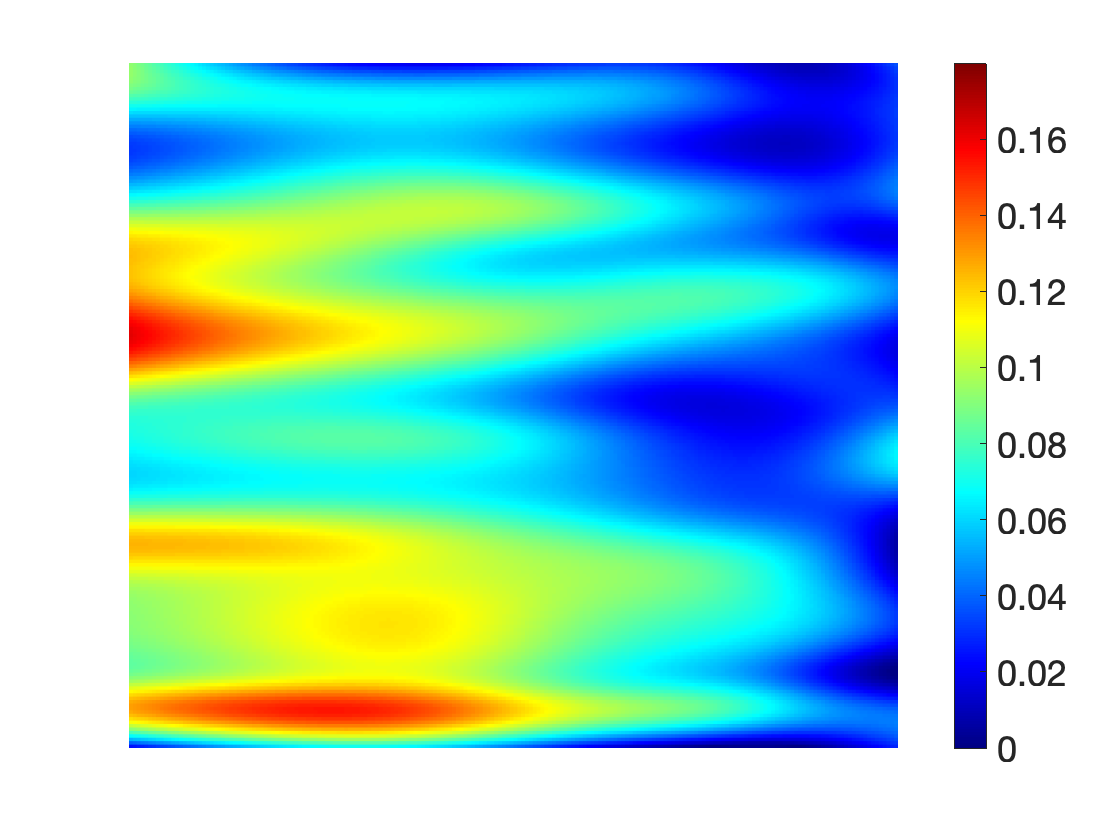} \\
\includegraphics[width=0.199\textwidth]{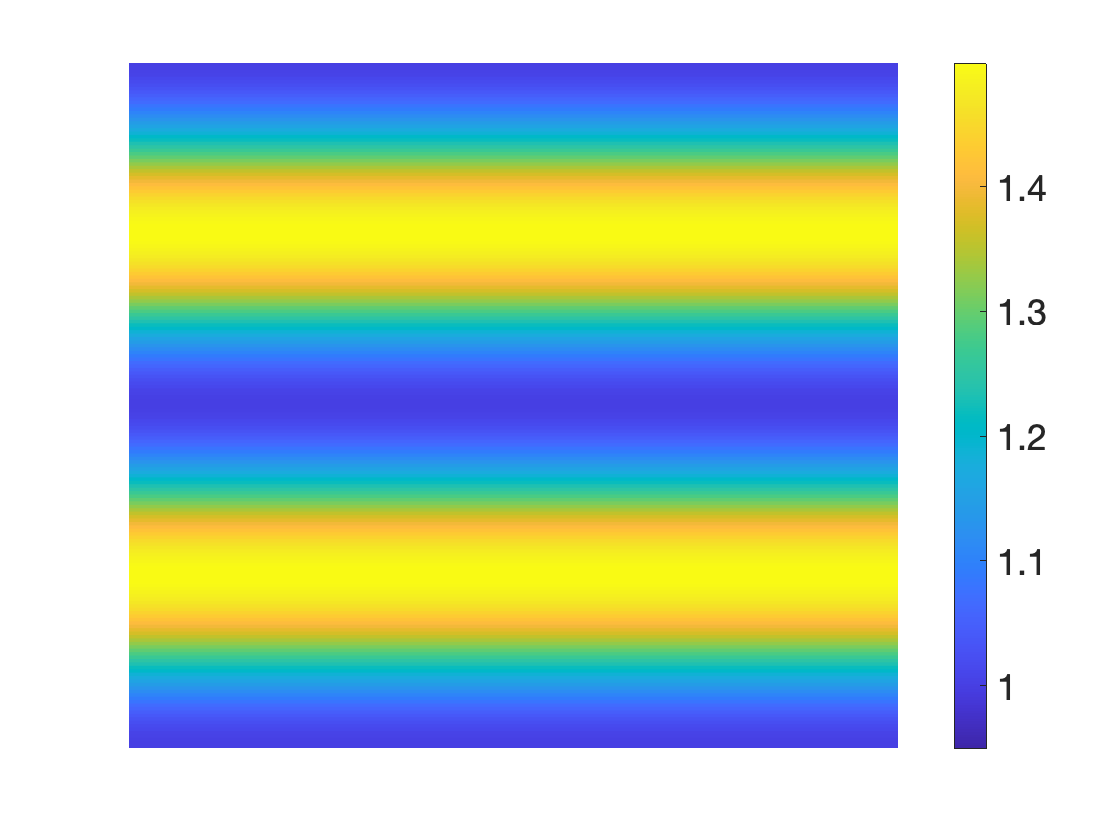} &
\includegraphics[width=0.199\textwidth]{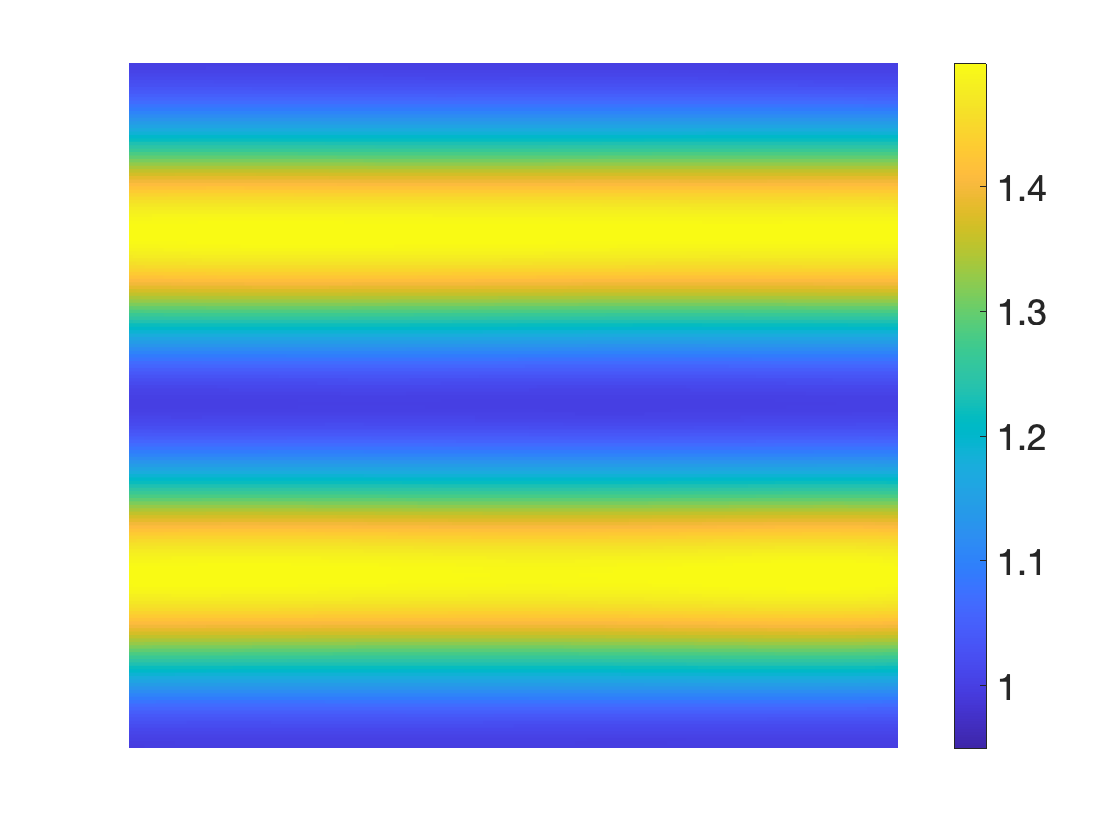} &
\includegraphics[width=0.199\textwidth]{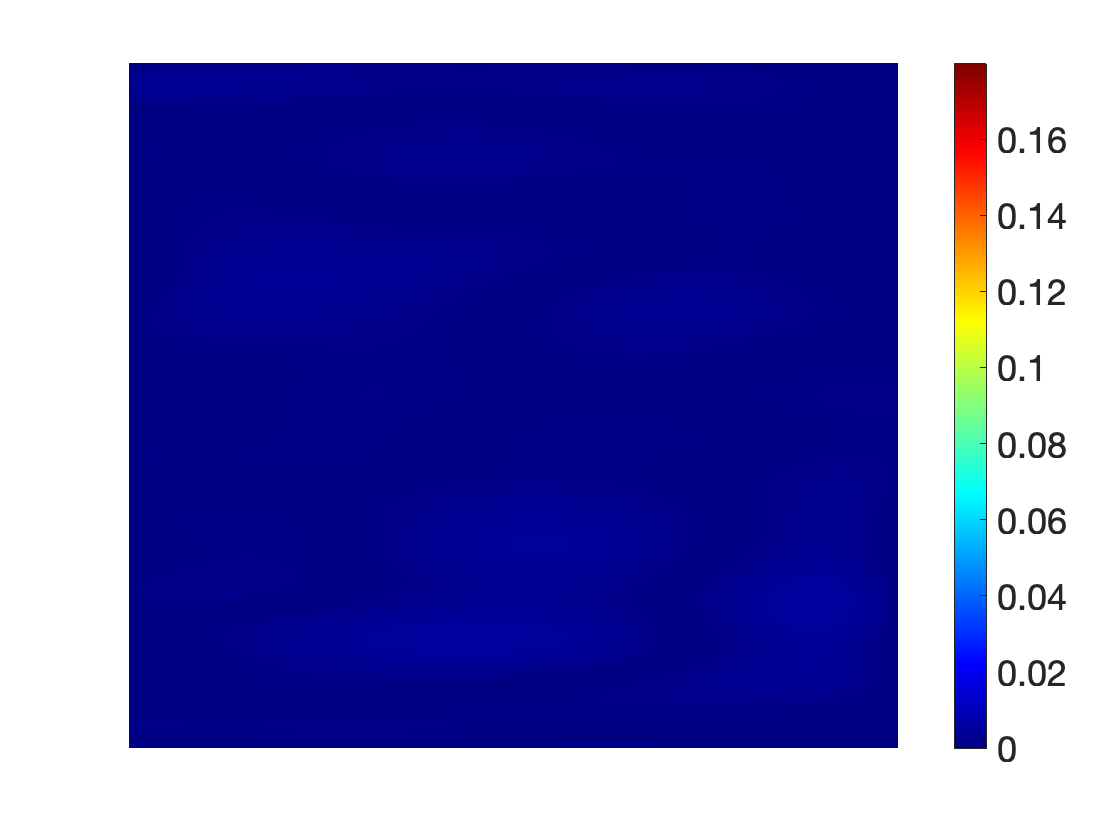} &
\includegraphics[width=0.199\textwidth]{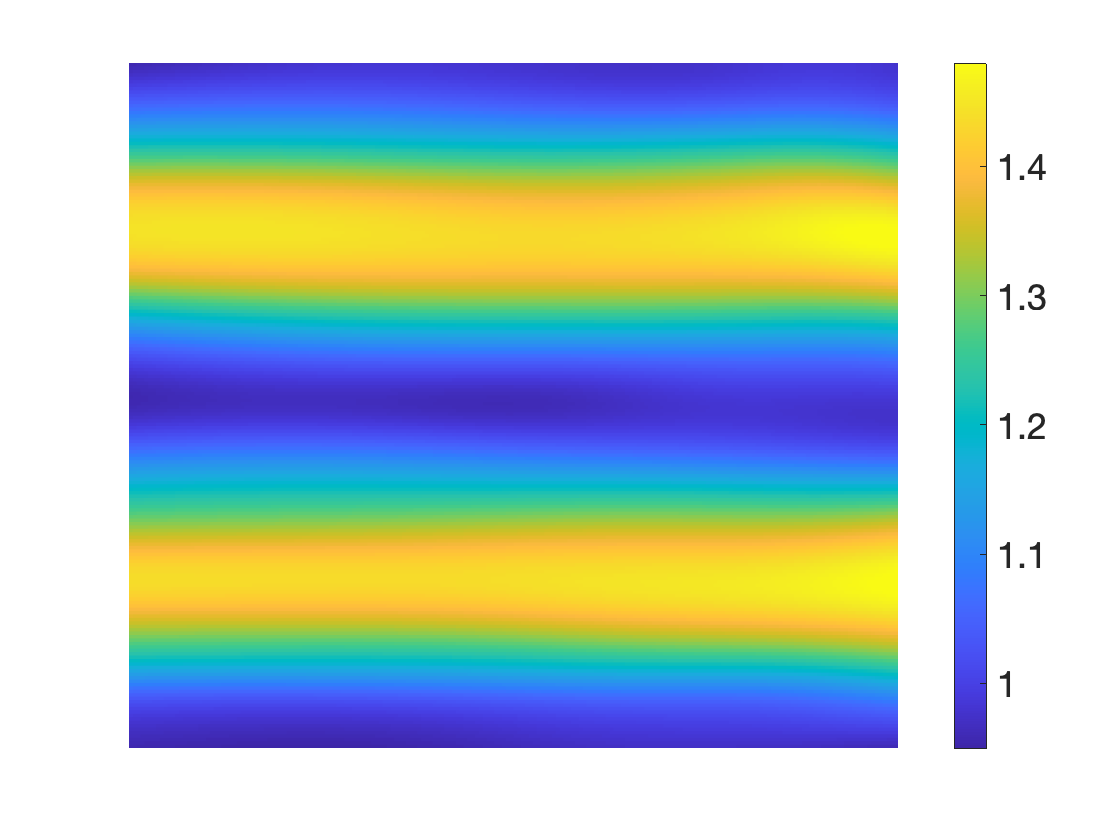} &
\includegraphics[width=0.199\textwidth]{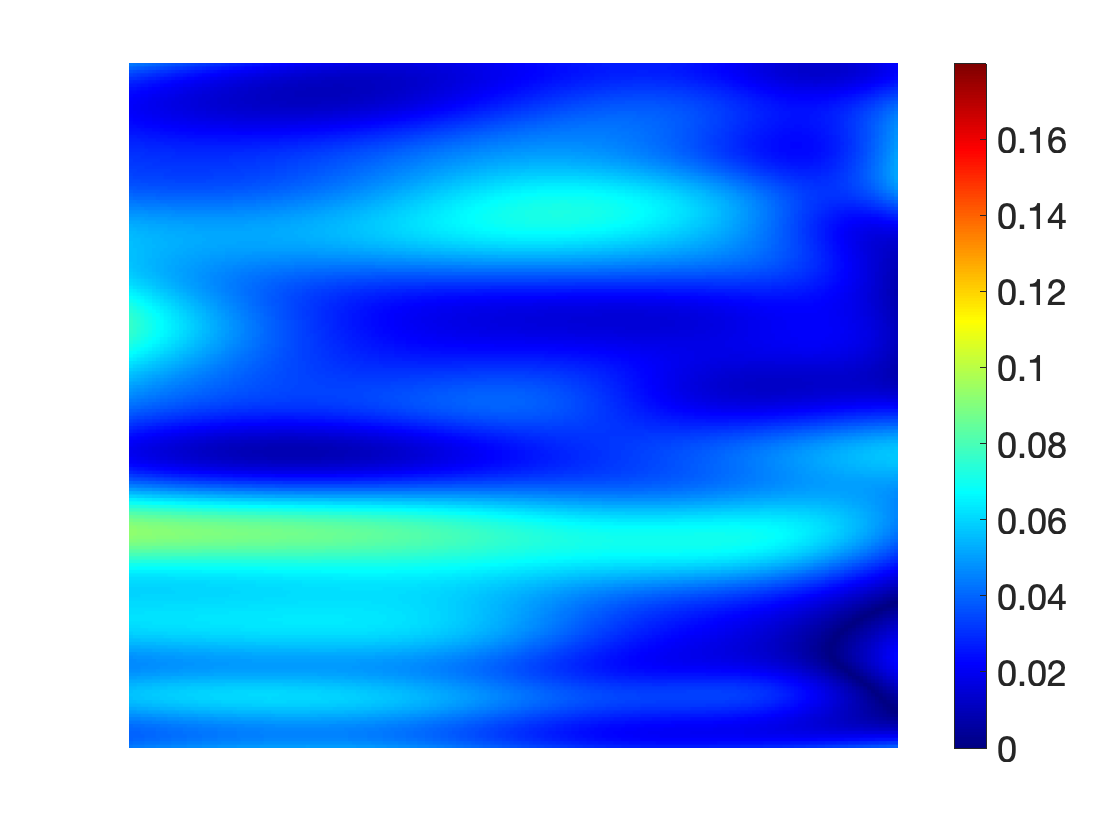} \\
\includegraphics[width=0.199\textwidth]{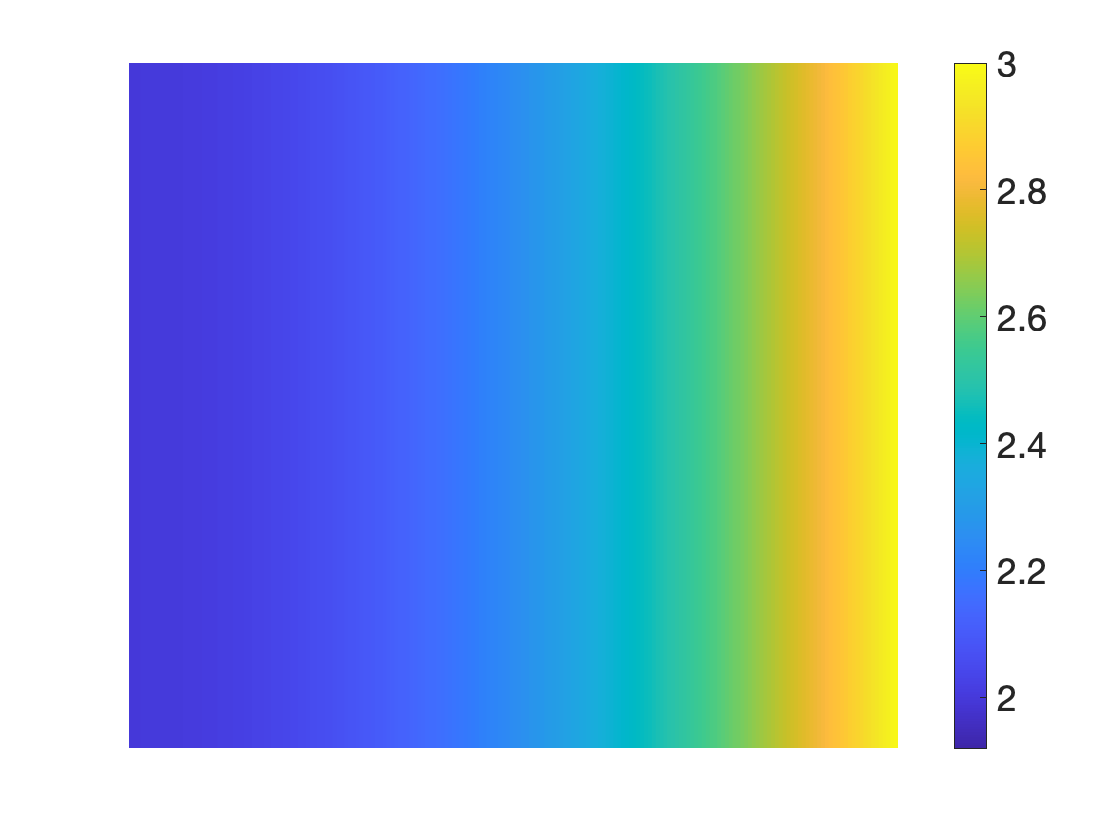} &
\includegraphics[width=0.199\textwidth]{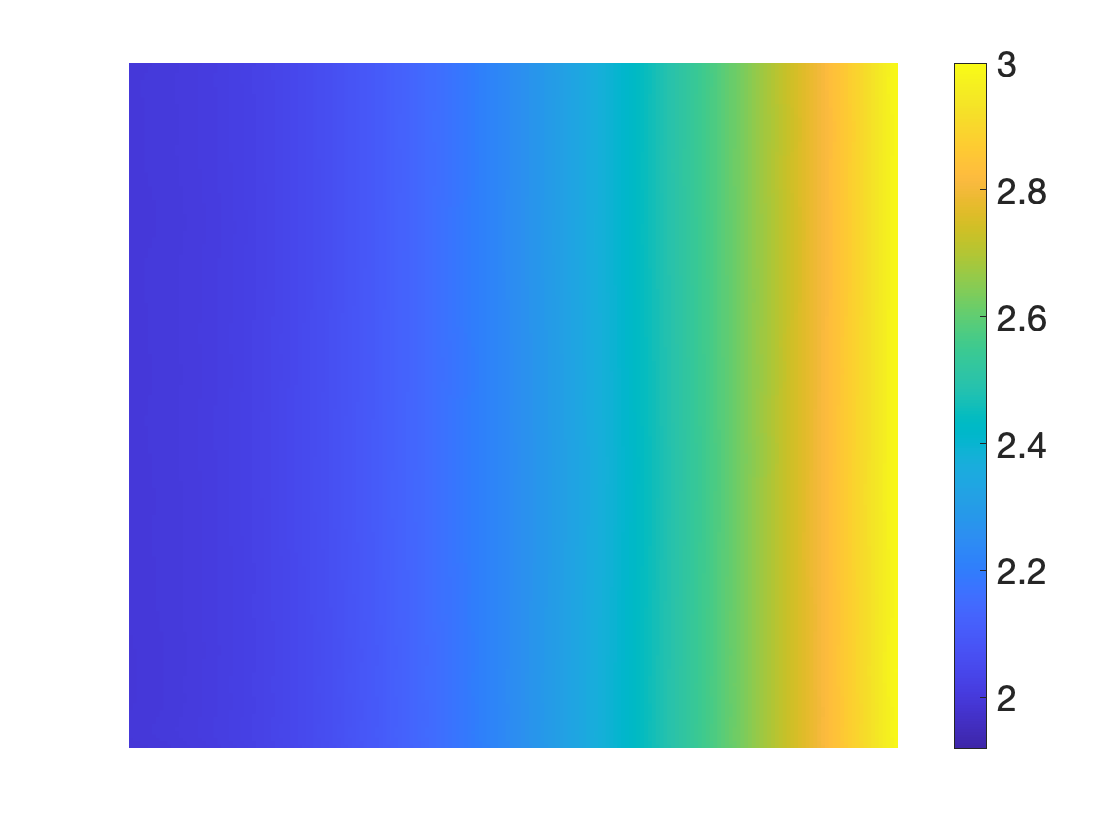} &
\includegraphics[width=0.199\textwidth]{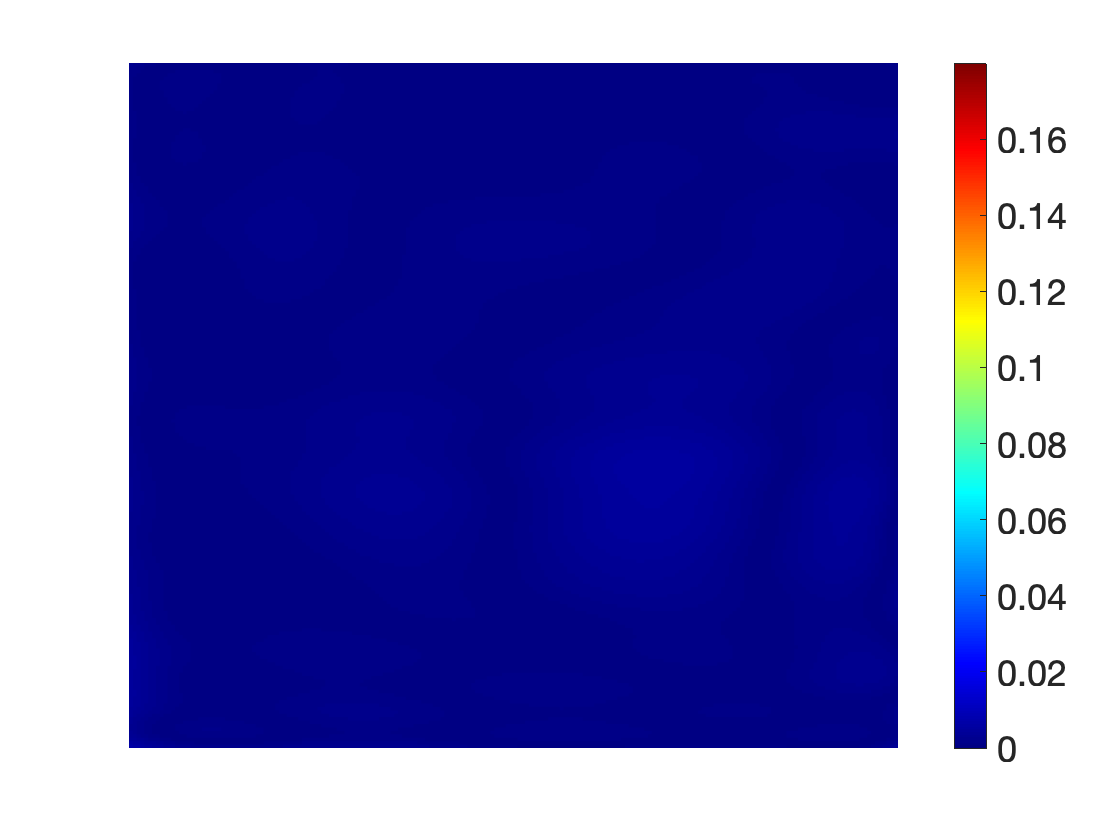} &
\includegraphics[width=0.199\textwidth]{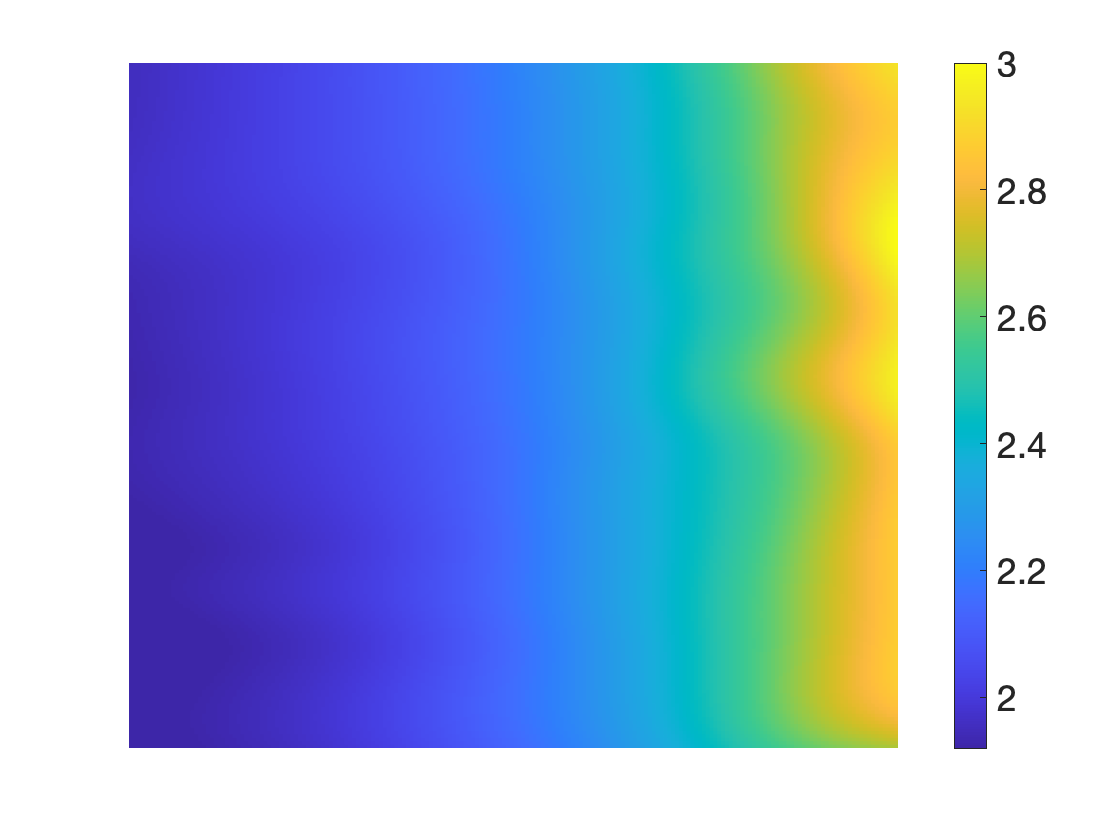} &
\includegraphics[width=0.199\textwidth]{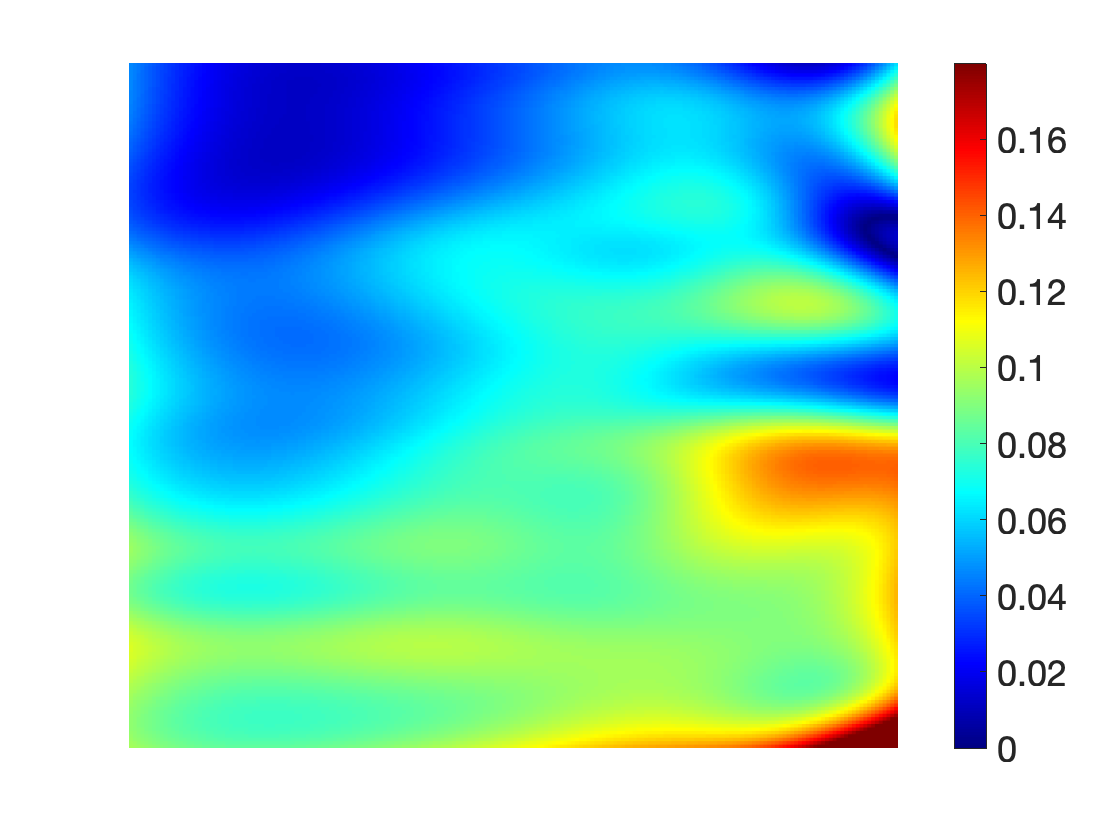} \\
(a) $A^\dag$  & (b) $\hat A$ & (c) $|\hat A-A^\dag|$ & (d) $\hat A$ & (e) $|\hat A-A^\dag|$
\end{tabular}
\caption{The reconstructions for Example \ref{exam:diri2d2} with exact data in (b) and noisy data $(\delta=10\%)$ in (d). From the top to bottom, the results are for $A_{11}$, $A_{12}$ and $A_{22}$, respectively.}
\label{fig:diri2d2}
\end{figure}

Fig. \ref{fig:diri2d2} shows the results for both exact and noisy data ($\delta = 10\%$). The DNN approximations can accurately capture the overall features of $A^\dagger$ in both cases, with almost no deterioration in the reconstruction quality even when 10\% noise is present in the data. This again shows the high robustness of the MLS-DNN approach against data noise.

\begin{figure}[htb!]
\centering
\setlength{\tabcolsep}{0em}
\begin{tabular}{ccccc}
\includegraphics[width=0.199\textwidth]{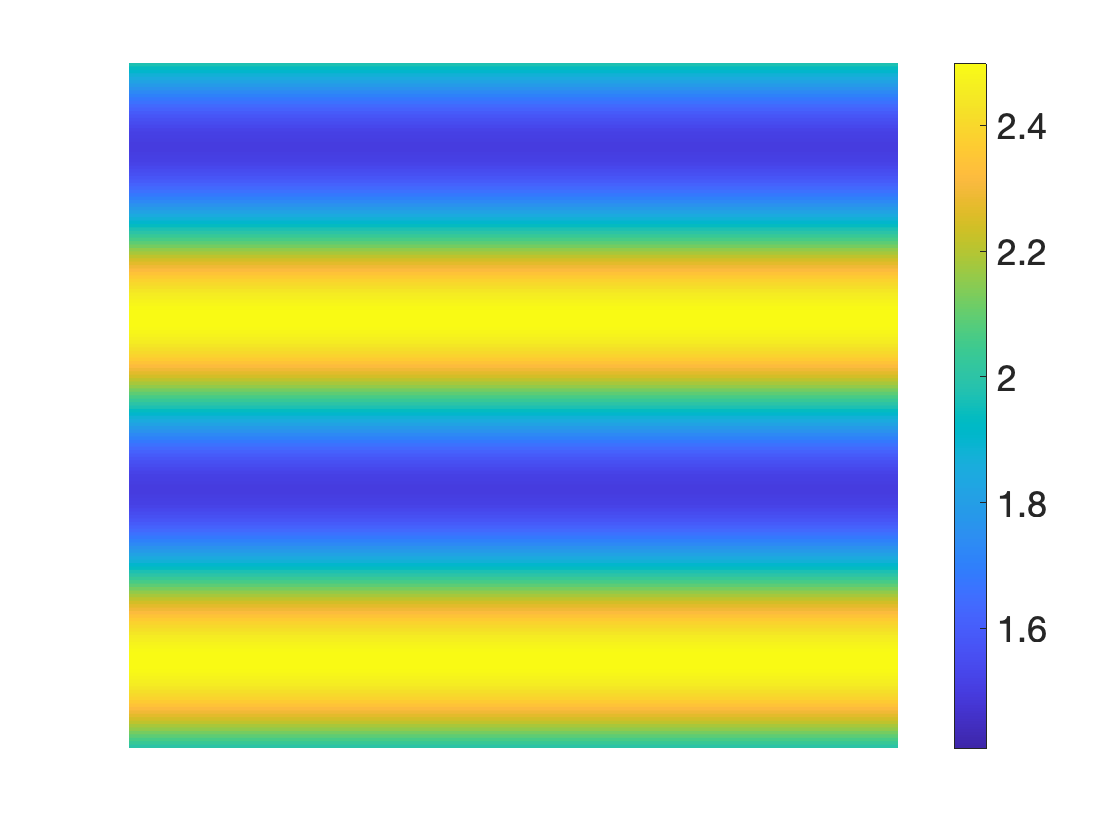} &
\includegraphics[width=0.199\textwidth]{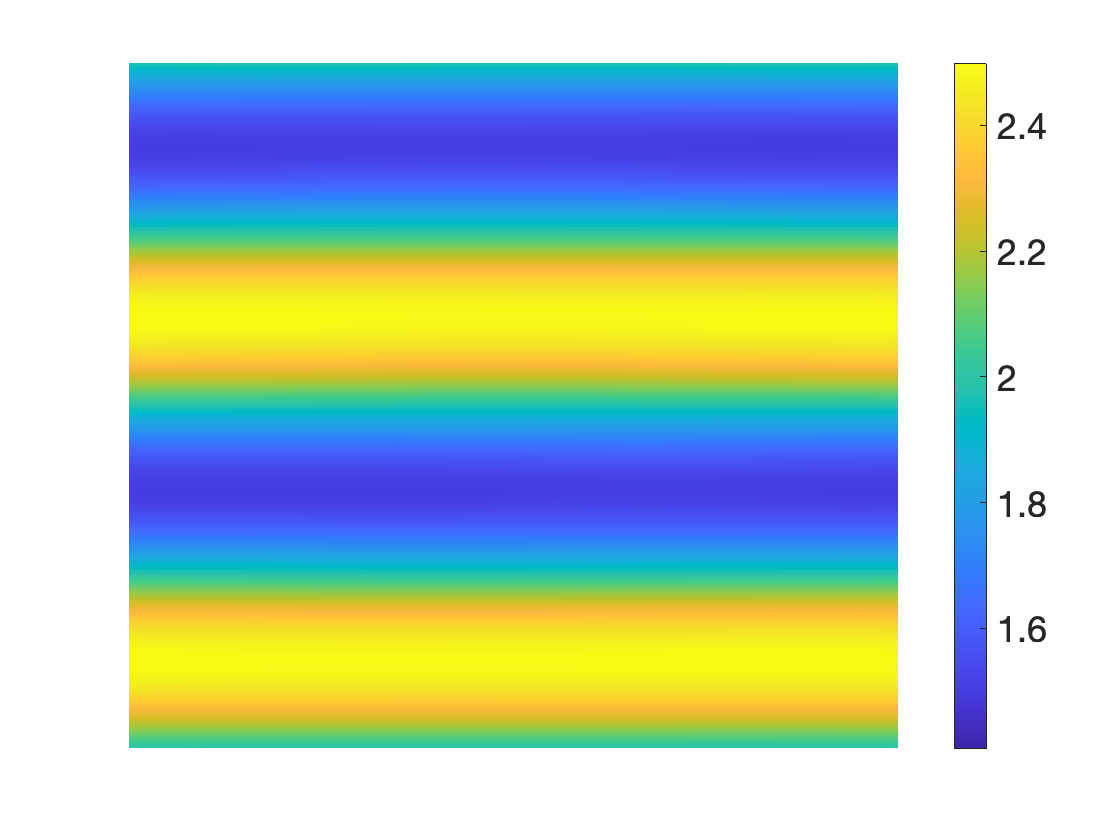} &
\includegraphics[width=0.199\textwidth]{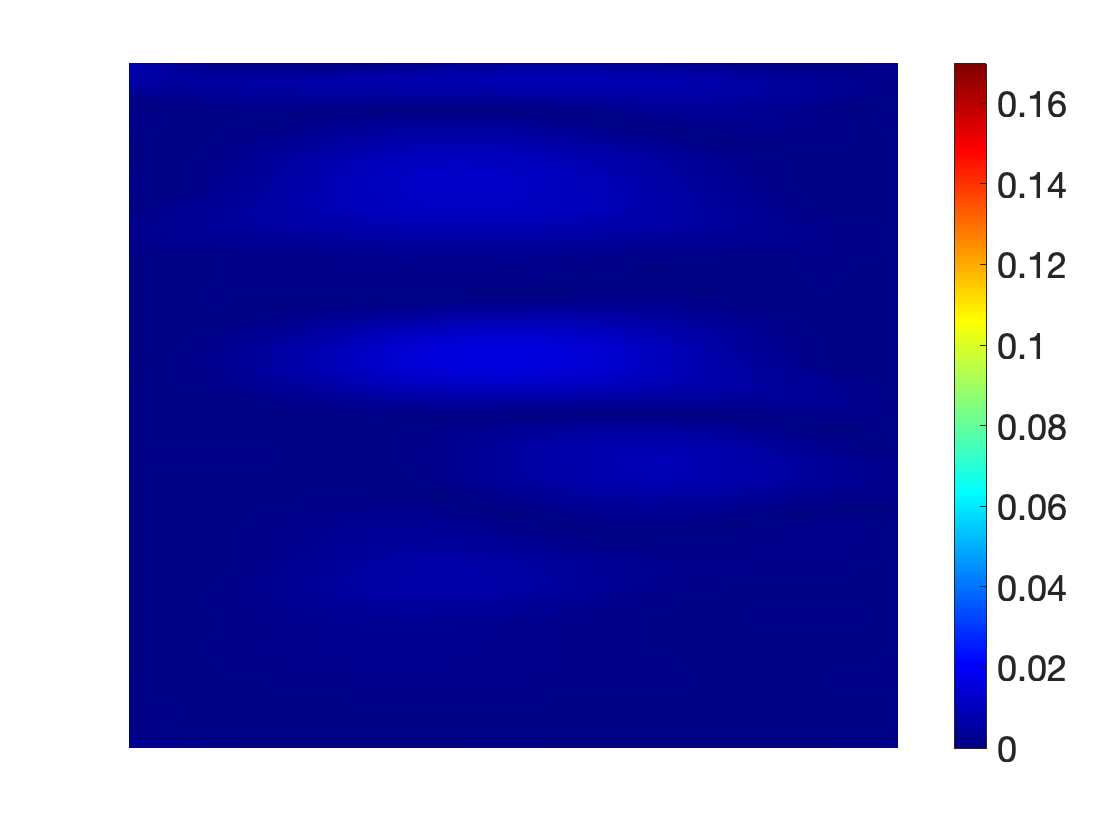} &
\includegraphics[width=0.199\textwidth]{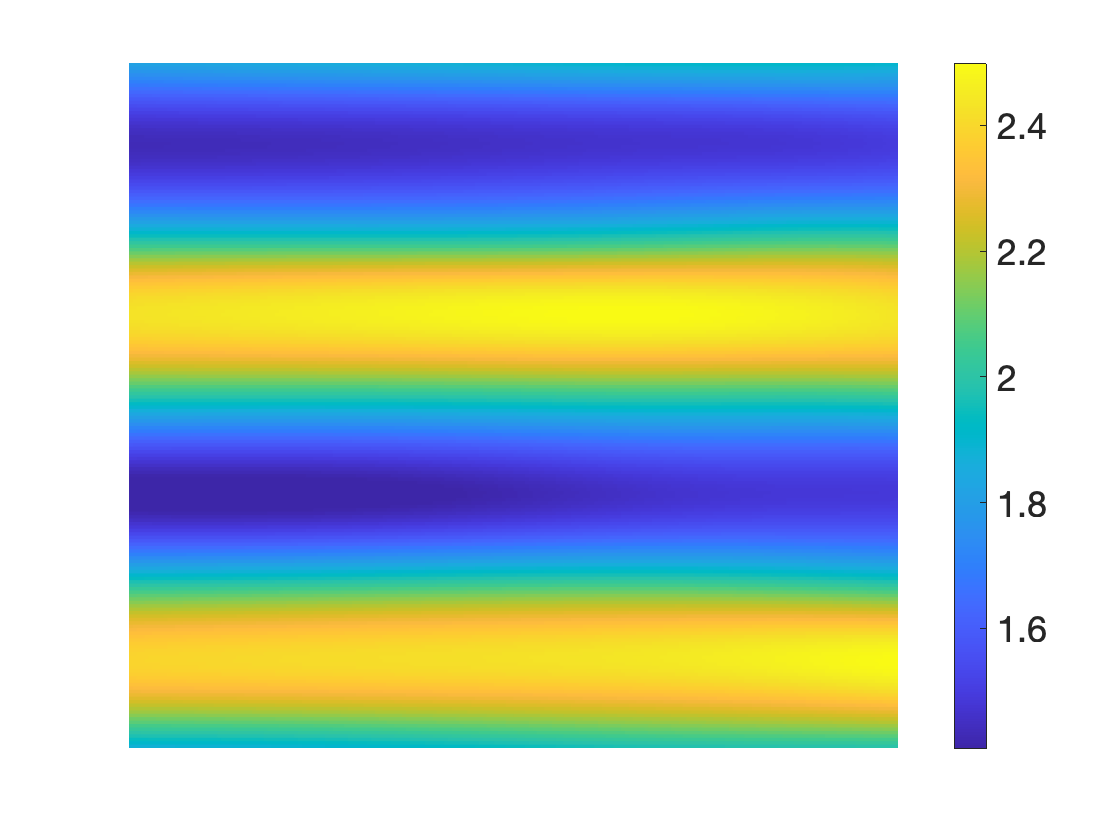} &
\includegraphics[width=0.199\textwidth]{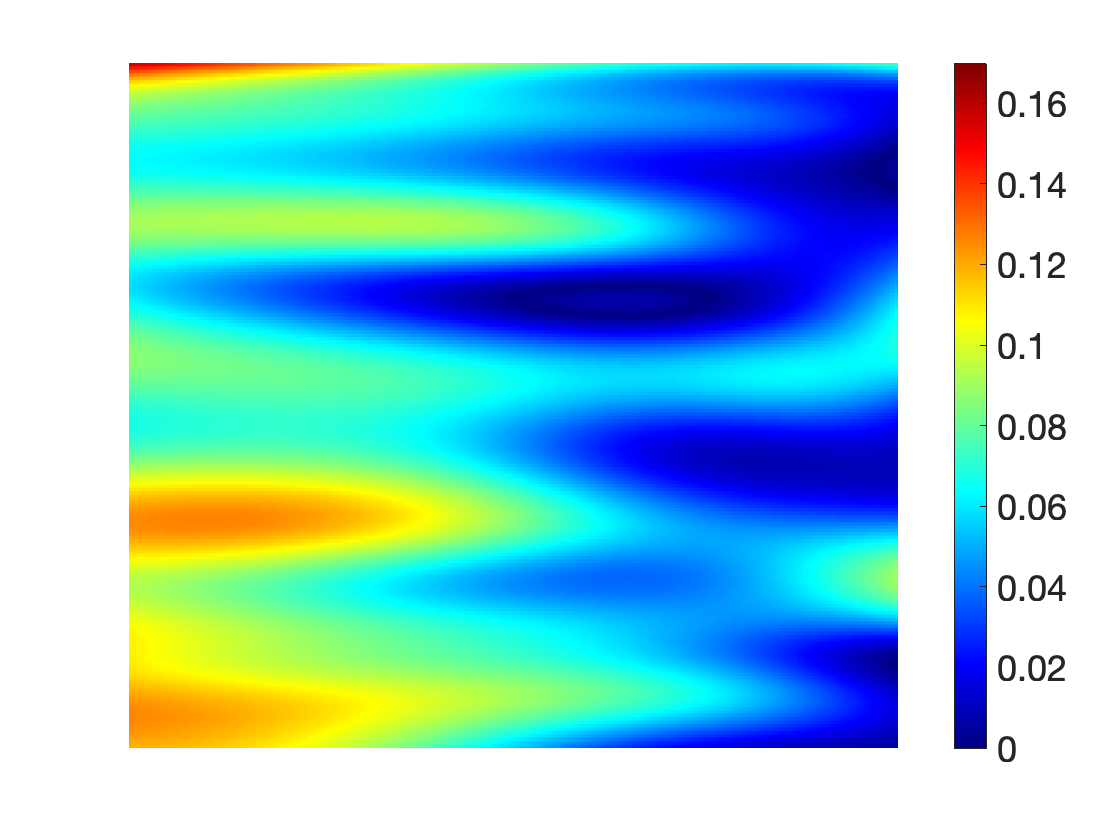} \\
\includegraphics[width=0.199\textwidth]{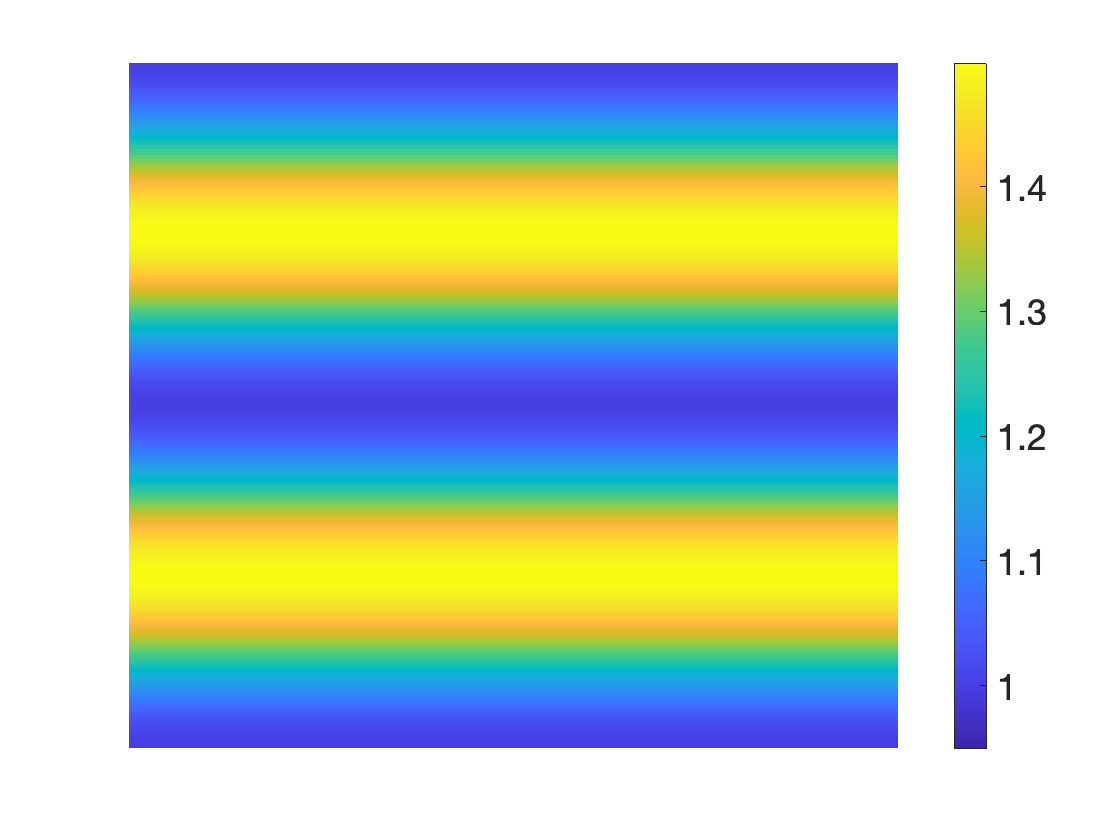} &
\includegraphics[width=0.199\textwidth]{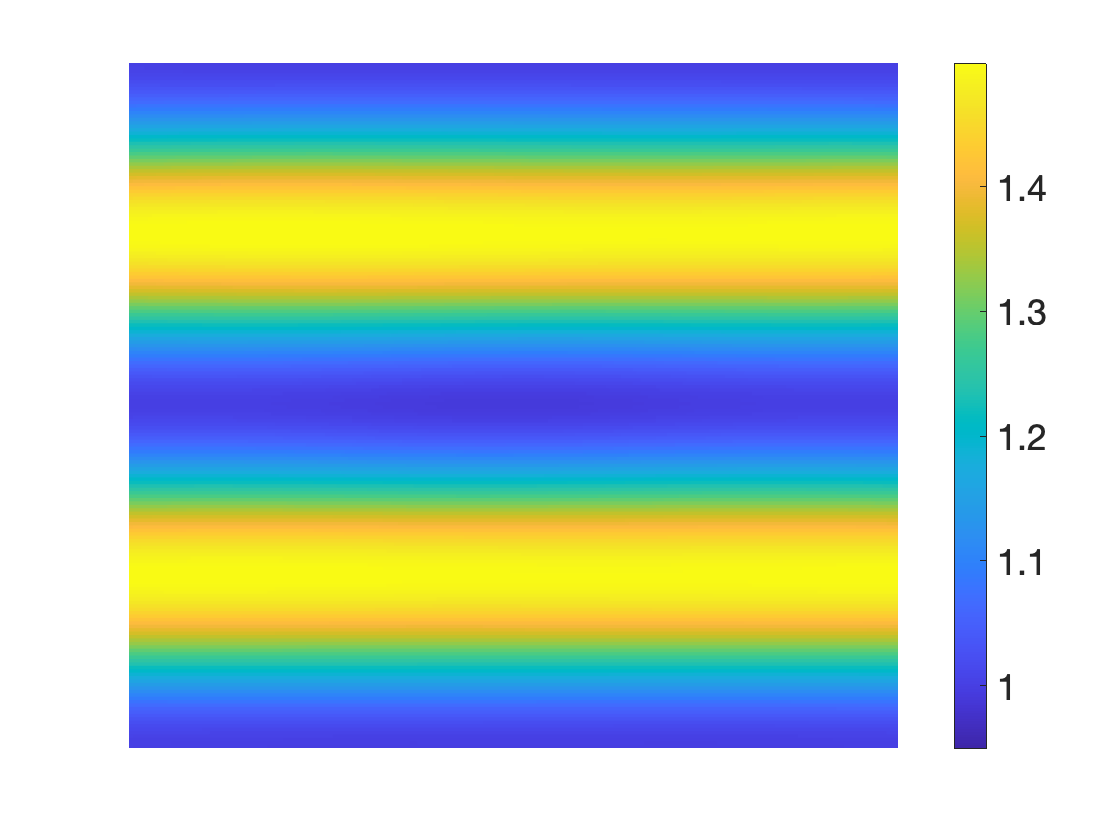} &
\includegraphics[width=0.199\textwidth]{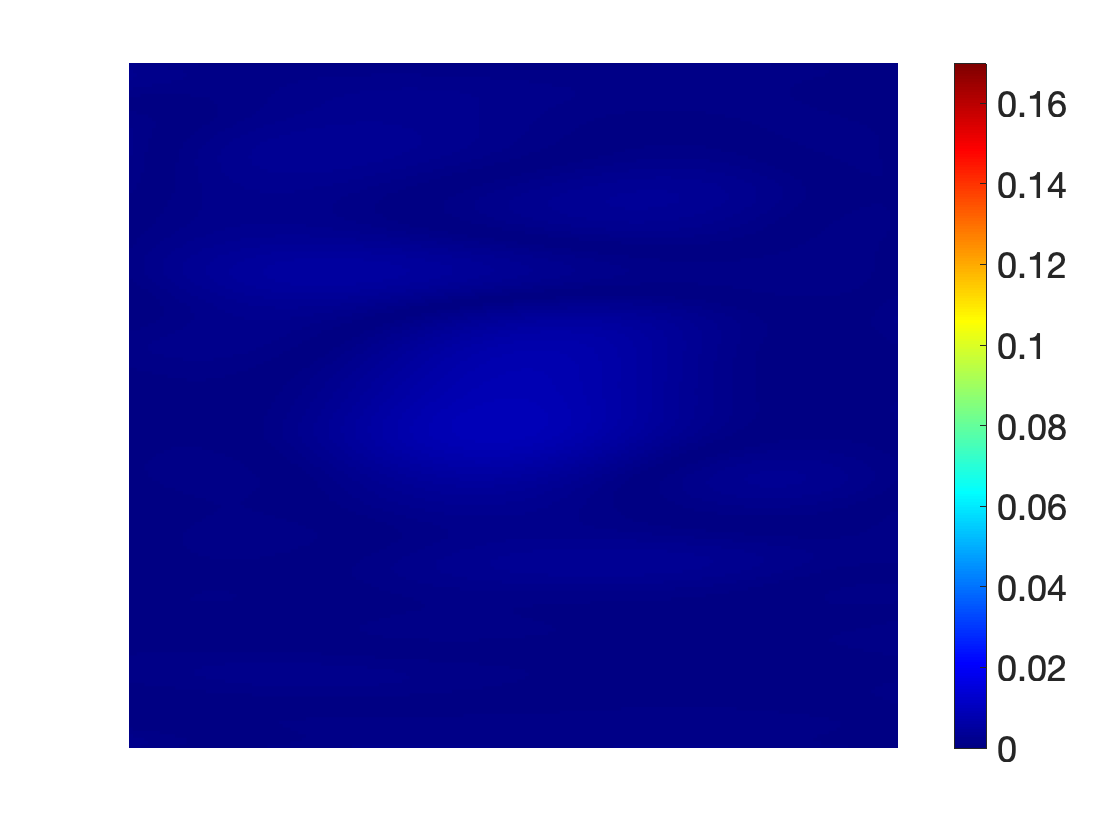} &
\includegraphics[width=0.199\textwidth]{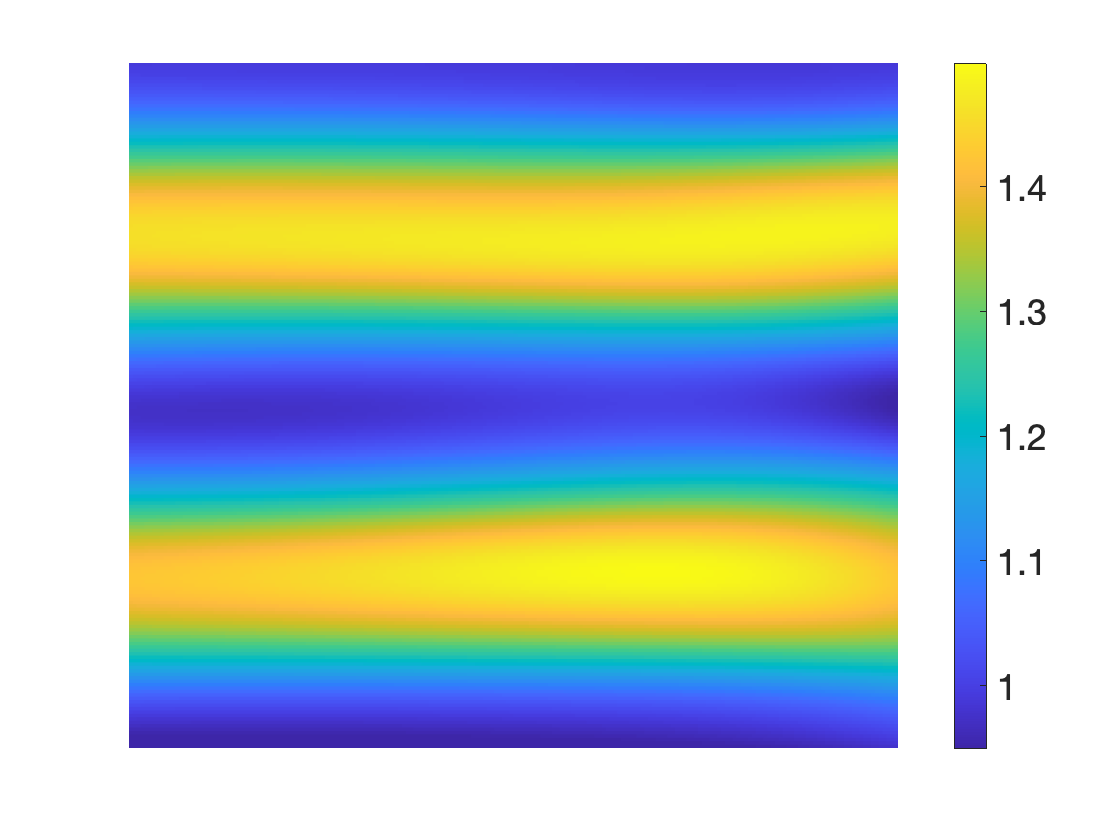} &
\includegraphics[width=0.199\textwidth]{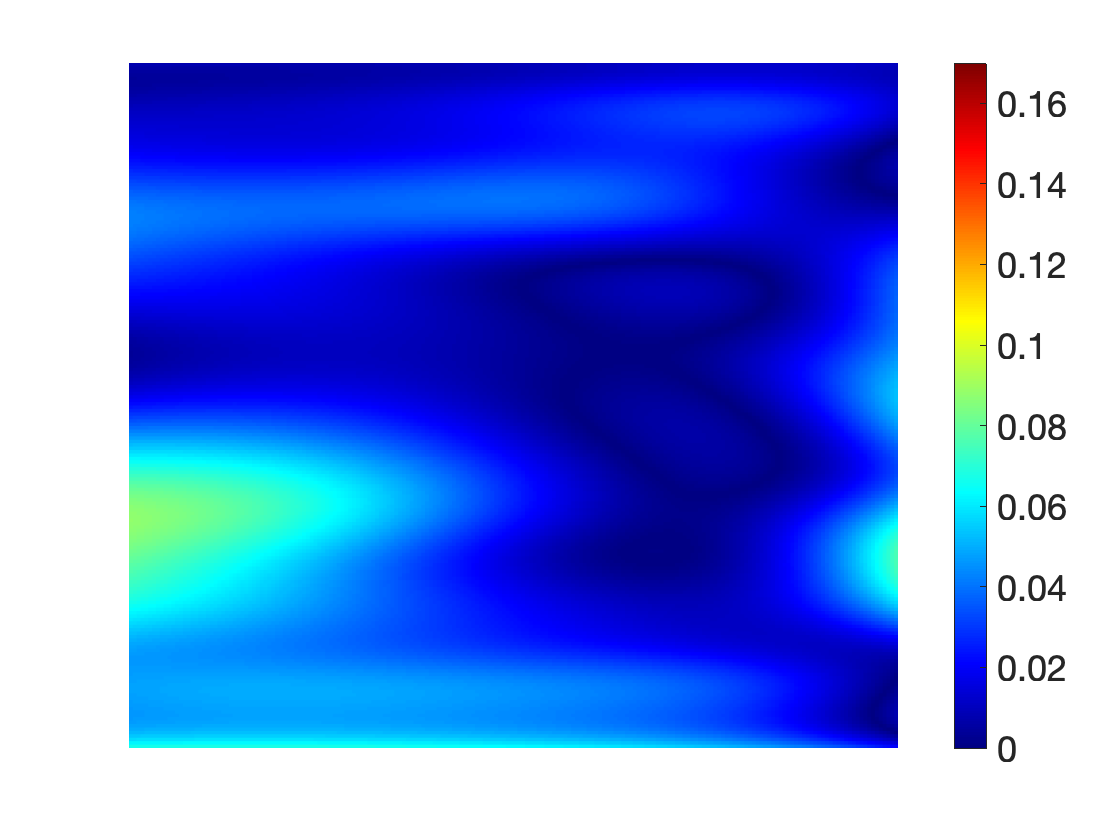} \\
\includegraphics[width=0.199\textwidth]{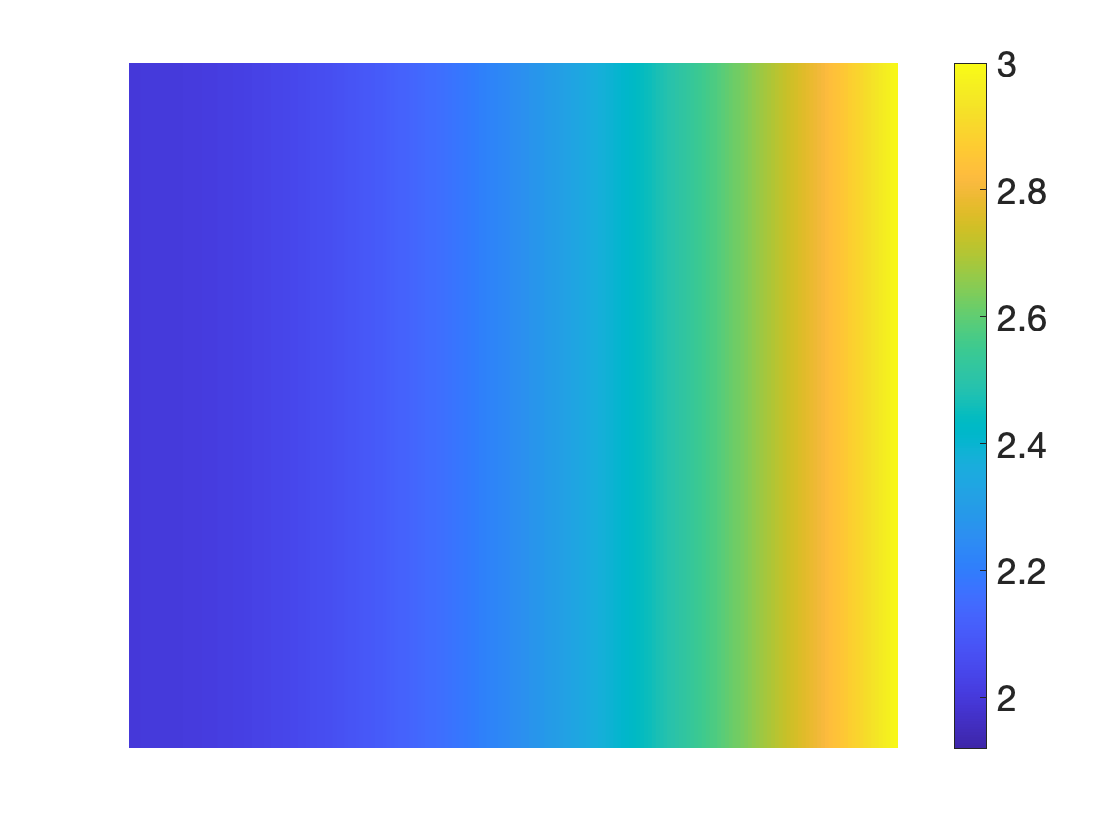} &
\includegraphics[width=0.199\textwidth]{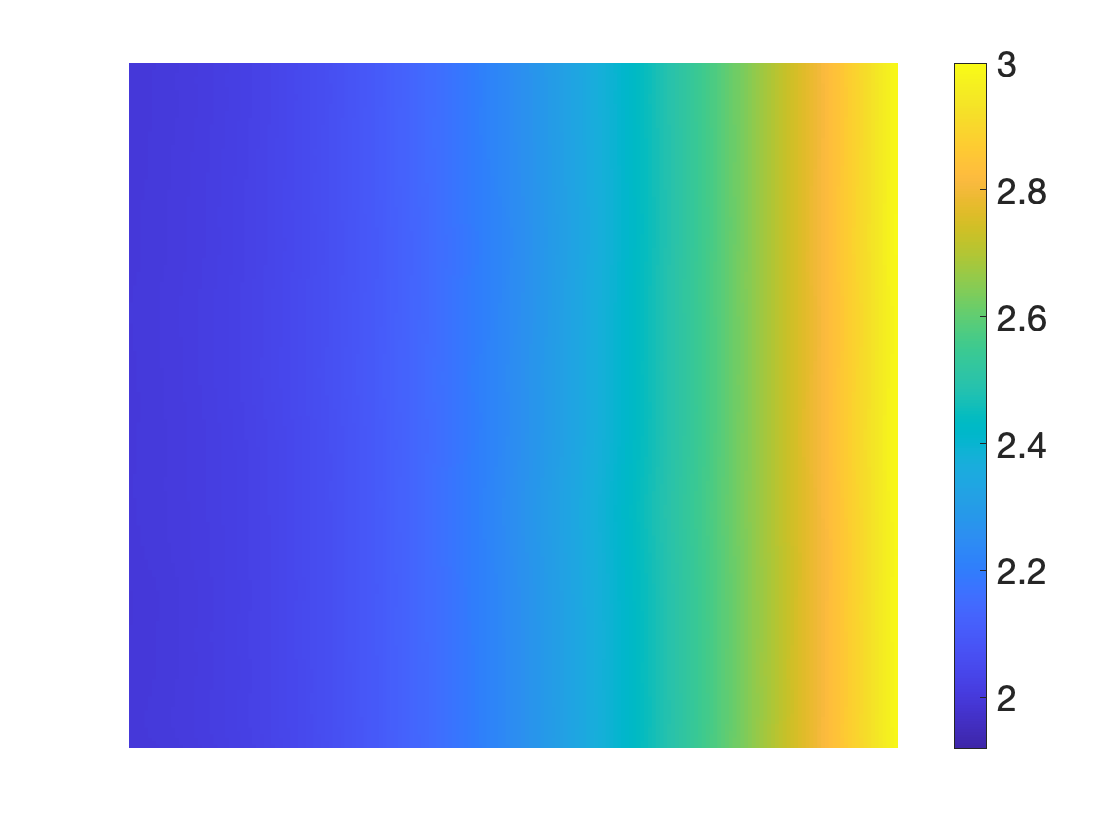} &
\includegraphics[width=0.199\textwidth]{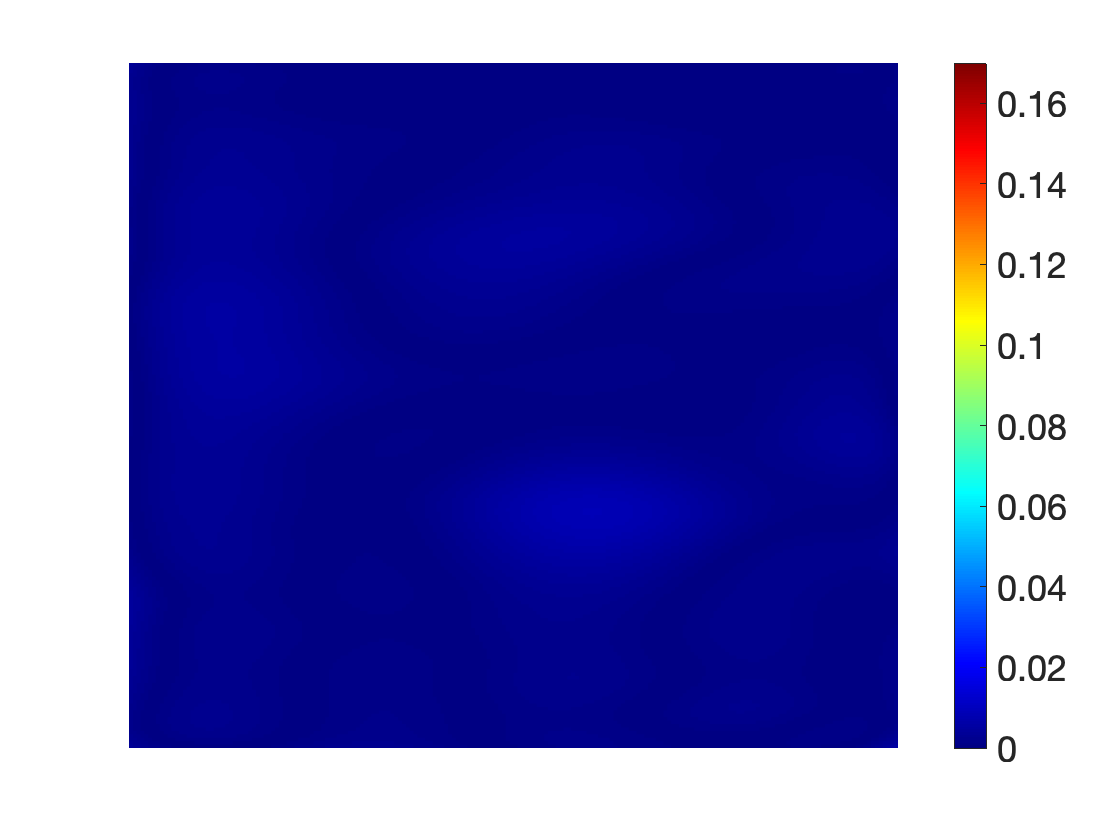} &
\includegraphics[width=0.199\textwidth]{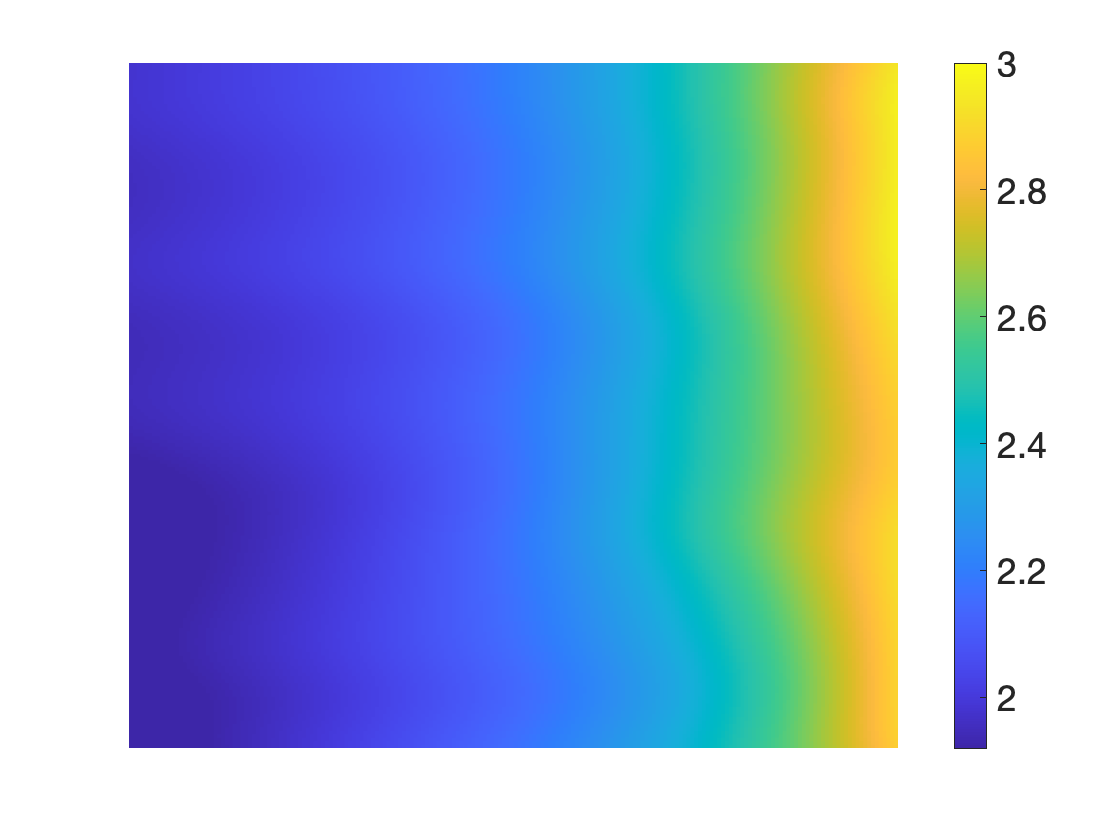} &
\includegraphics[width=0.199\textwidth]{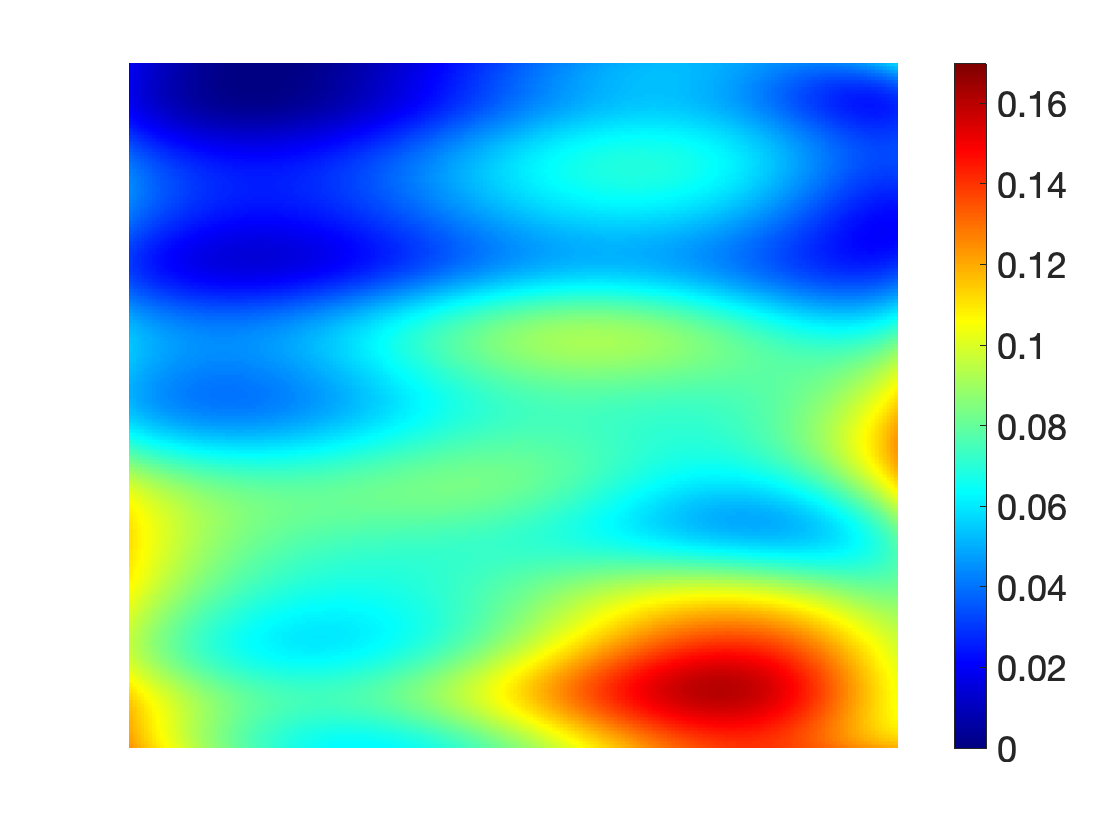} \\
(a) $A^\dag$  & (b) $\hat A$ & (c) $|\hat A-A^\dag|$ & (d) $\hat A$ & (e) $|\hat A-A^\dag|$
\end{tabular}
\caption{The reconstructions for Example \ref{exam:diri2d3} with exact data in (b) and noisy data $(\delta=10\%)$ in (d). From the top to bottom, the results are for $A_{11}$, $A_{12}$ and $A_{22}$, respectively.}
\label{fig:diri2d3}
\end{figure}

Next we show the reconstruction of the conductivity matrix from partial internal data.
\begin{example}
The domain $\Omega = (0,1)^2$, the measurement $\nabla  z_i^\delta$ on the region $\omega=\Omega\setminus (0.2,0.8)^2$, $A^\dag= \begin{pmatrix}
    2+\frac{\sin(4\pi x_2)}{2}&1+\frac{\sin^2(2\pi x_2)}{2}\\
    1+\frac{\sin^2(2\pi x_2)}{2}&2+x_1^2\\
\end{pmatrix},$ $u_1^\dag=x_1+x_2+\frac{1}{3}(x_1^3+x_2^3)$, $u_2^\dag=x_1-x_2+\frac{1}{3}(x_1^3-x_2^3)$, $u_3^\dag=-u_1^\dagger$, $u_4^\dag=-u_2^\dagger$.
\label{exam:diri2d3}
\end{example}
The approach minimizes an empirical version of the following loss:
\begin{equation}\label{obj:loss_diripartial}
    \begin{aligned}
            J_{\gamma}(\theta,\kappa)=\sum_{i=1}^{N}&\left(\|\sigma_{i,\kappa}-P_{\mathcal{K}}(A_{\theta})\nabla z_i^{\delta}\|_{L^2(\omega)^d}^2+\gamma_{\sigma}\|\nabla\cdot\sigma_{i,\kappa}+f_i\|_{L^2(\Omega)}^2\right.\\
            &\left.+\gamma_b\| \sigma_{i,\kappa}-A^{\dag}\nabla z_i^{\delta} \|_{L^2(\partial\Omega)^d}^2+\gamma_A\|P_{\mathcal{K}}(A_{\theta})\|_{L^2(\Omega)^{d,d}}^2\right).
    \end{aligned}
\end{equation}
Fig. \ref{fig:diri2d3} shows that the reconstruction results are accurate across the entire domain $\Omega$ for both exact and noisy data, including the central region where no observational data $ \nabla z_i^\delta$ is provided. This example again shows the stability of the MLS-DNN approach for partial internal data.

\begin{figure}[htb!]
\centering
\setlength{\tabcolsep}{0em}
\begin{tabular}{ccccc}
\includegraphics[width=0.21\textwidth]{a11diri2dpex.png} &
\includegraphics[width=0.199\textwidth]{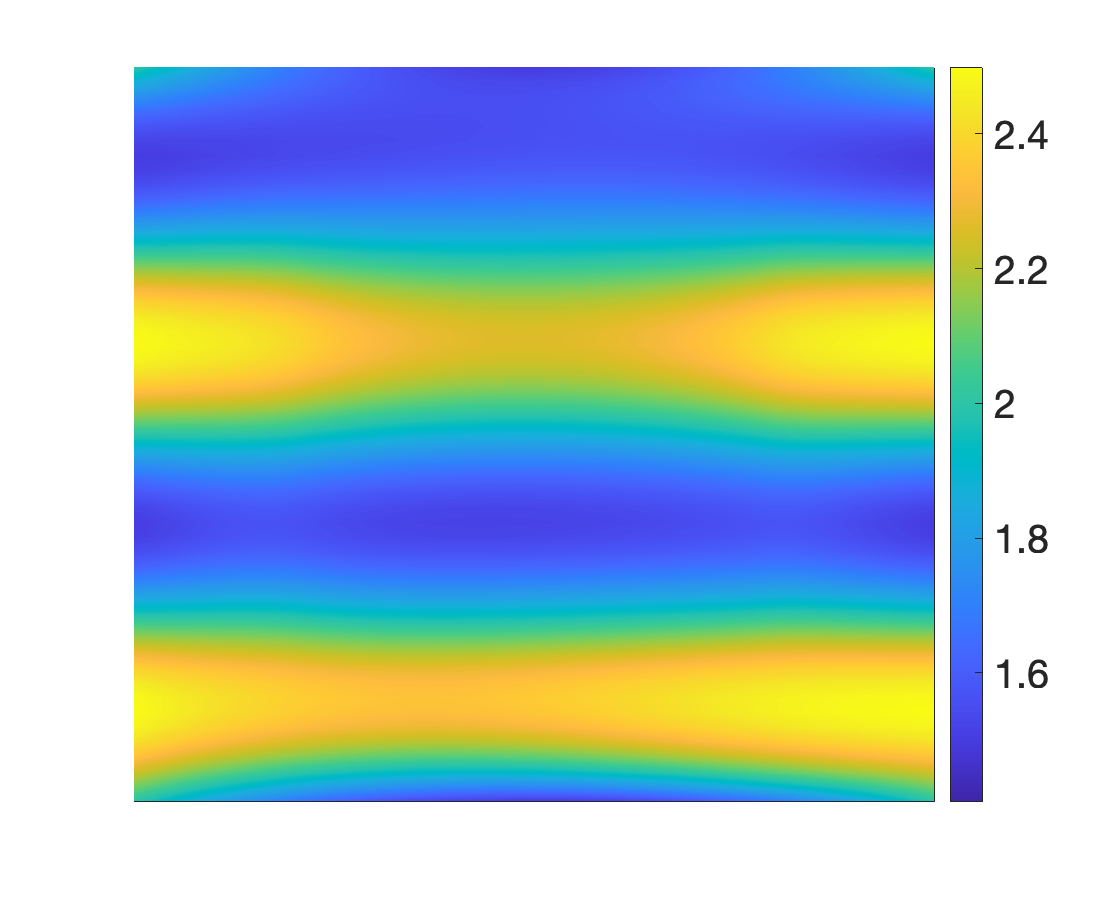} &
\includegraphics[width=0.199\textwidth]{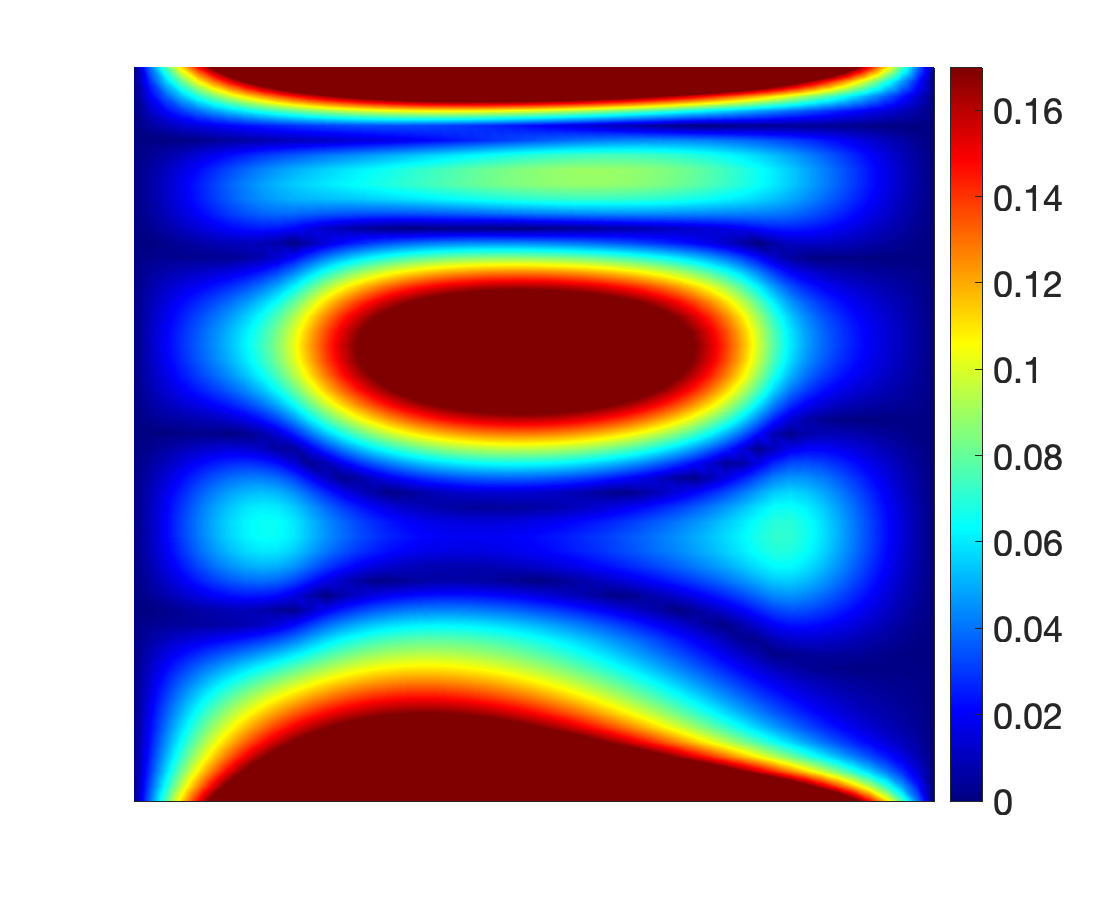} &
\includegraphics[width=0.199\textwidth]{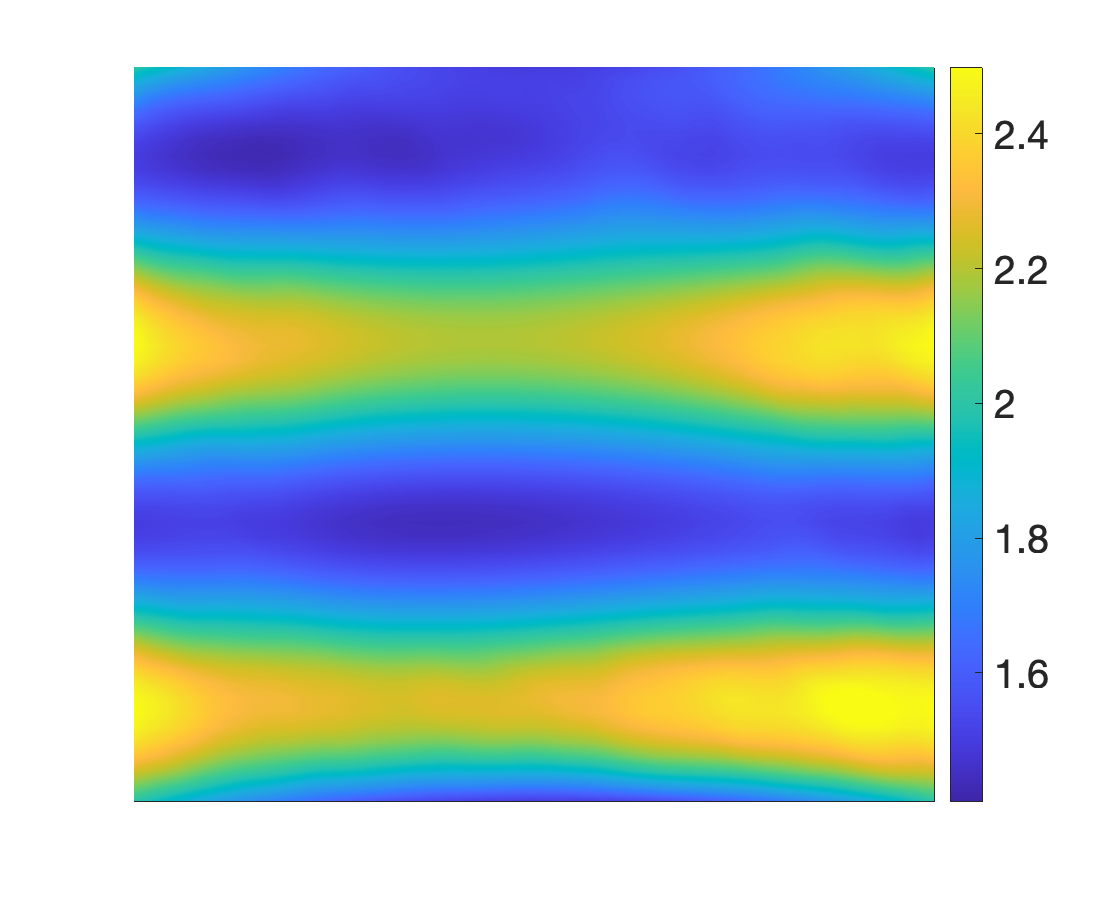} &
\includegraphics[width=0.199\textwidth]{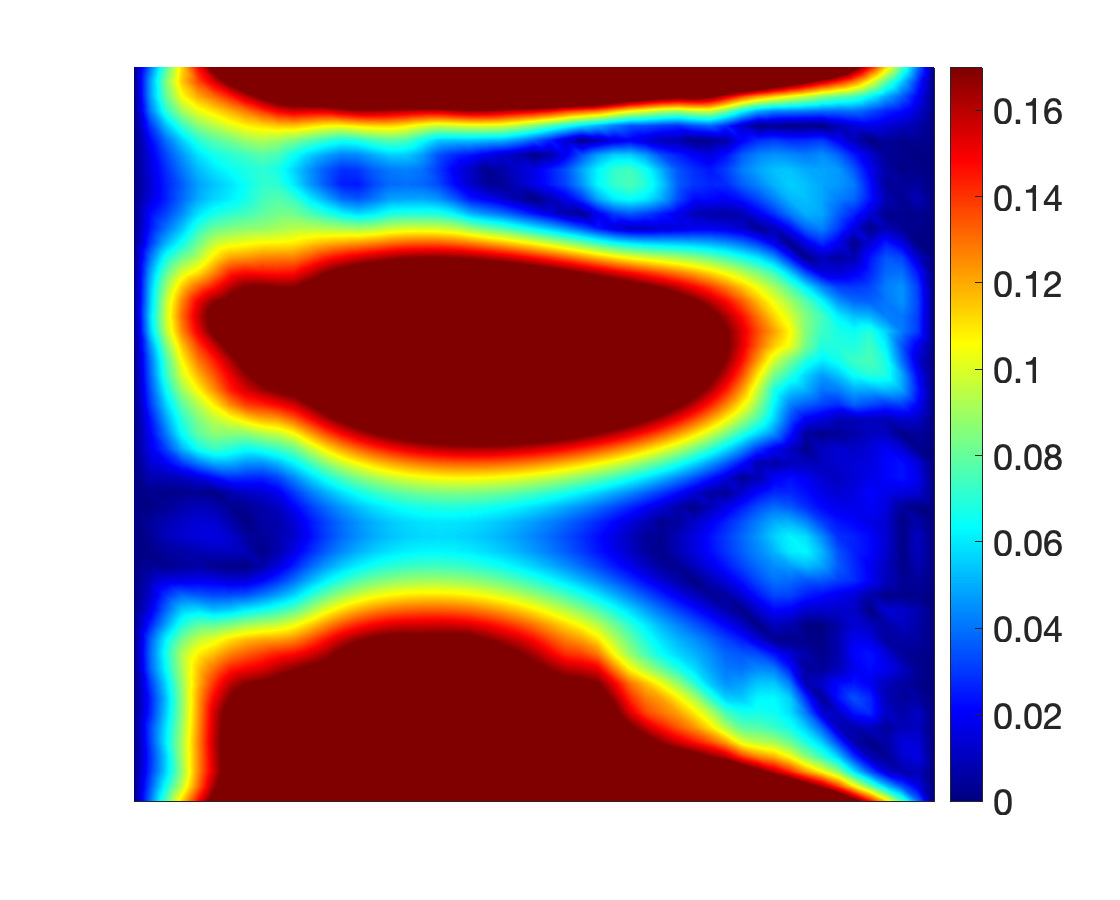} \\
\includegraphics[width=0.21\textwidth]{a12diri2dpex.png} &
\includegraphics[width=0.199\textwidth]{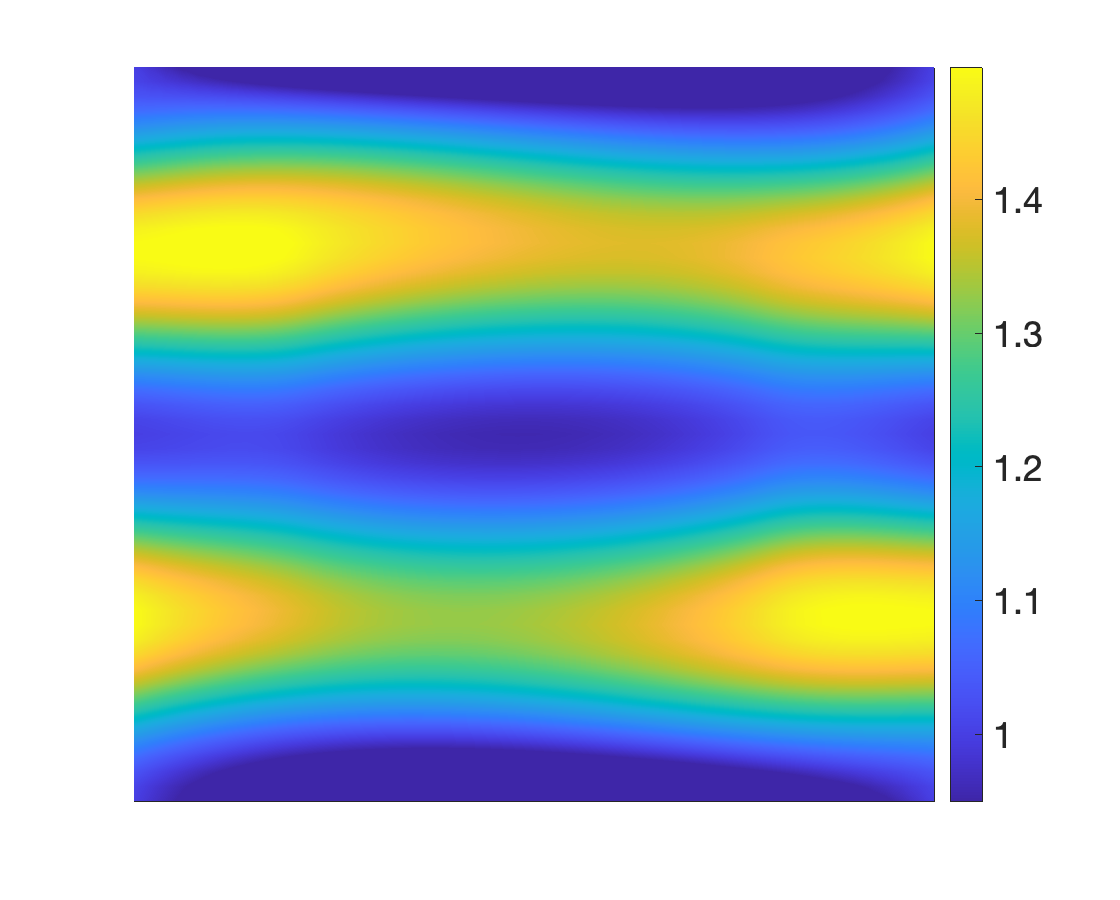} &
\includegraphics[width=0.199\textwidth]{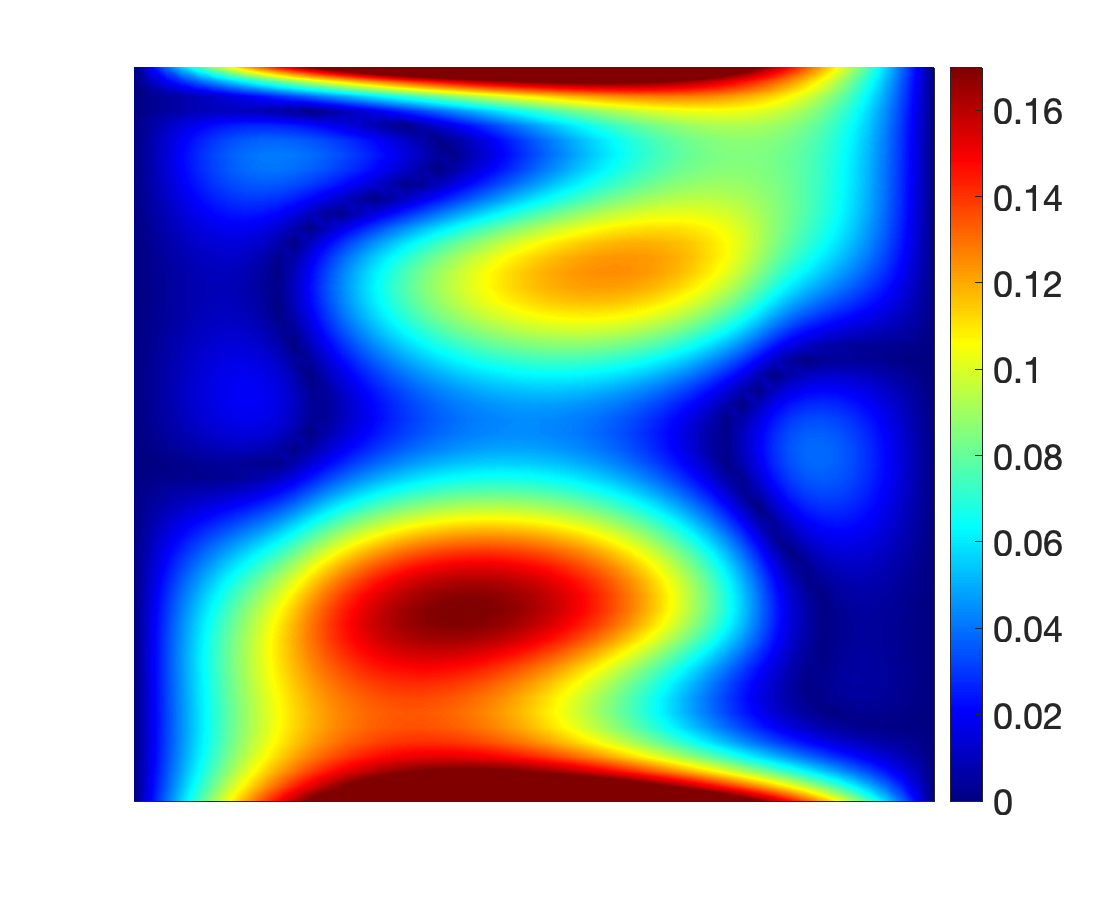} &
\includegraphics[width=0.199\textwidth]{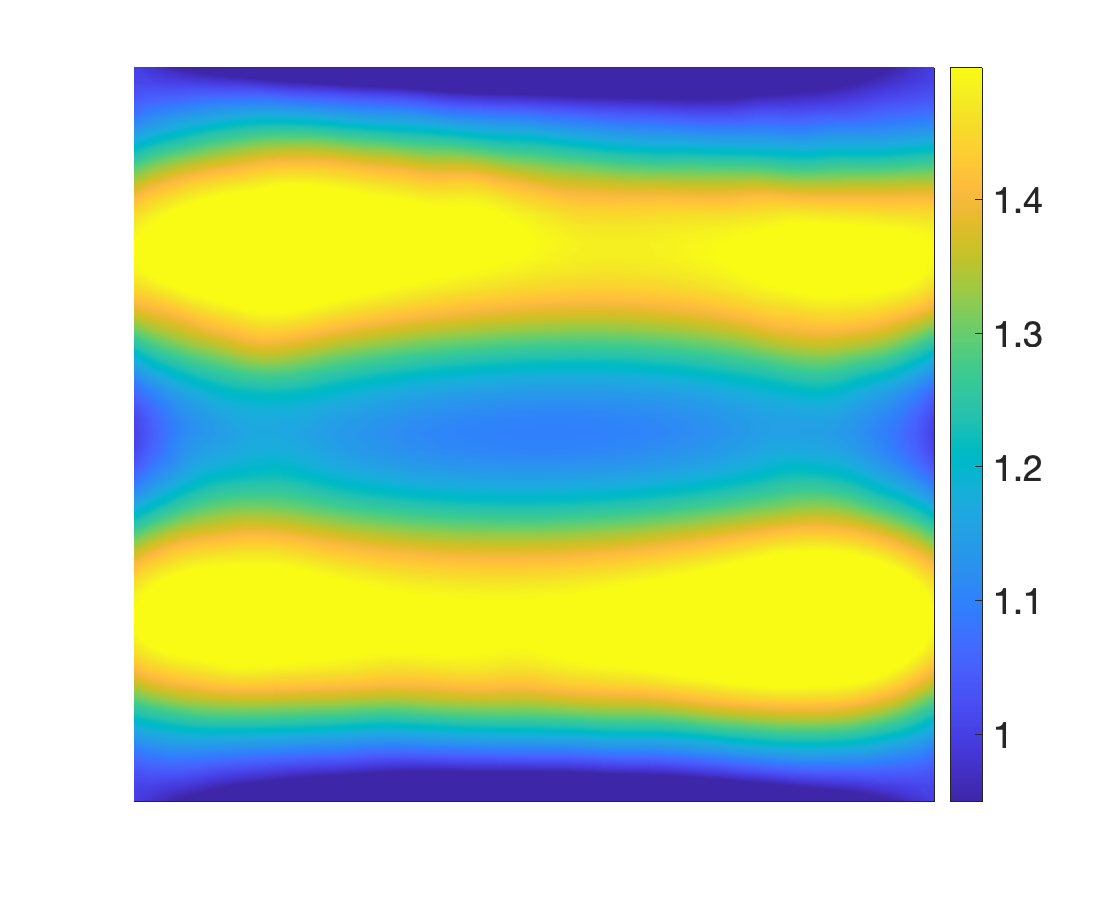} &
\includegraphics[width=0.199\textwidth]{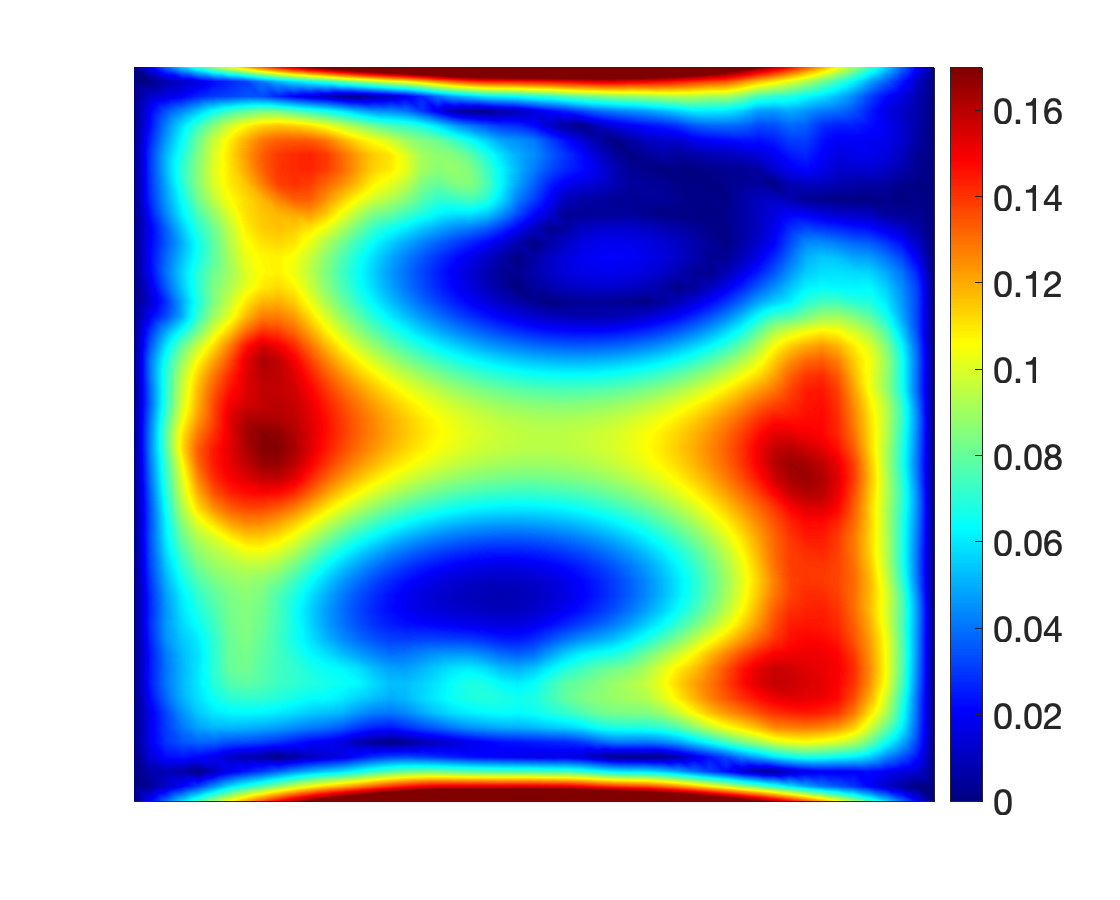} \\
\includegraphics[width=0.21\textwidth]{a22diri2dpex.png} &
\includegraphics[width=0.199\textwidth]{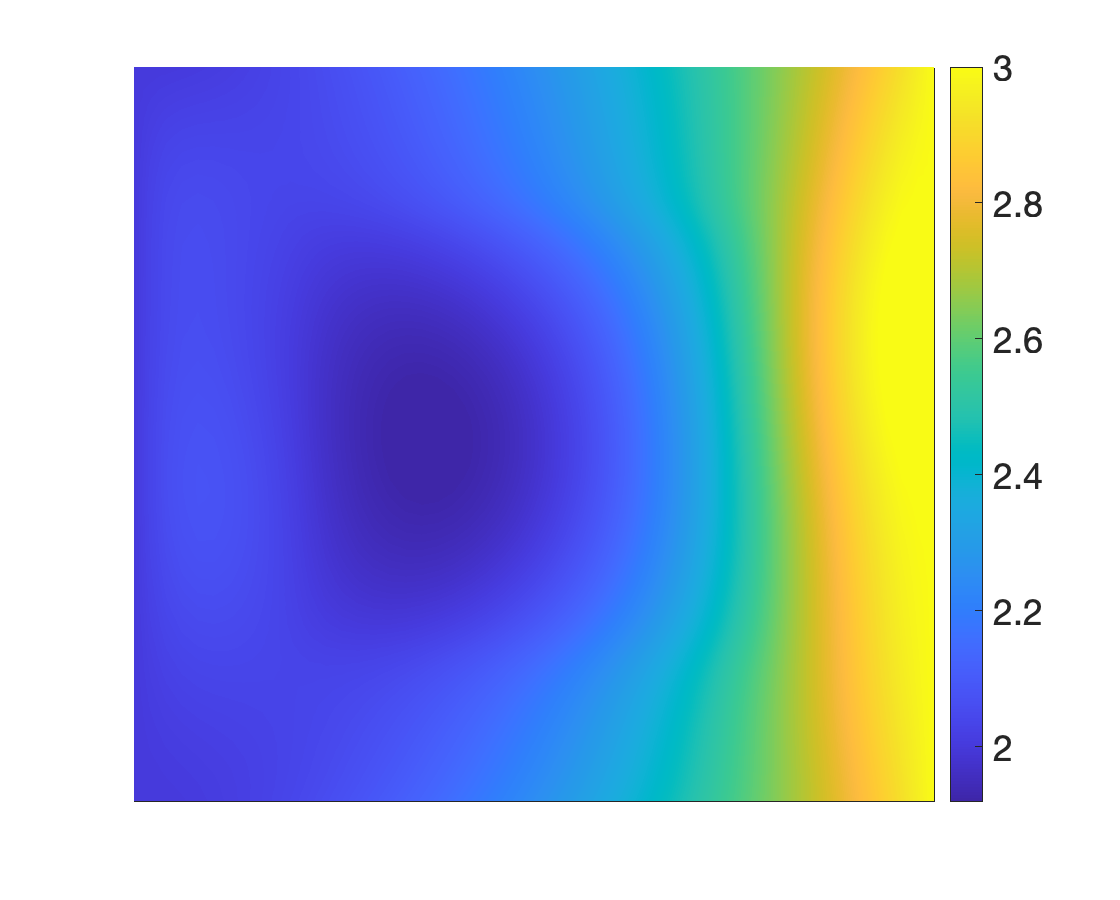} &
\includegraphics[width=0.199\textwidth]{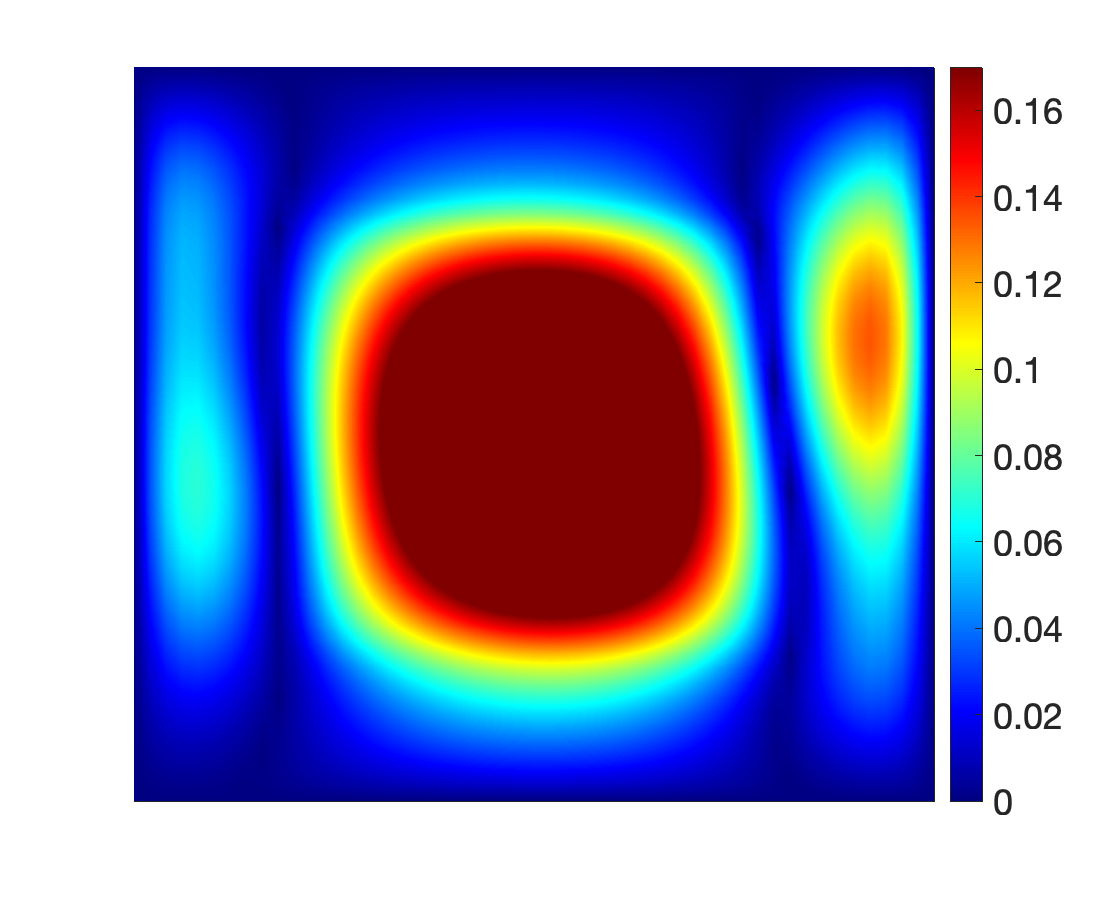} &
\includegraphics[width=0.199\textwidth]{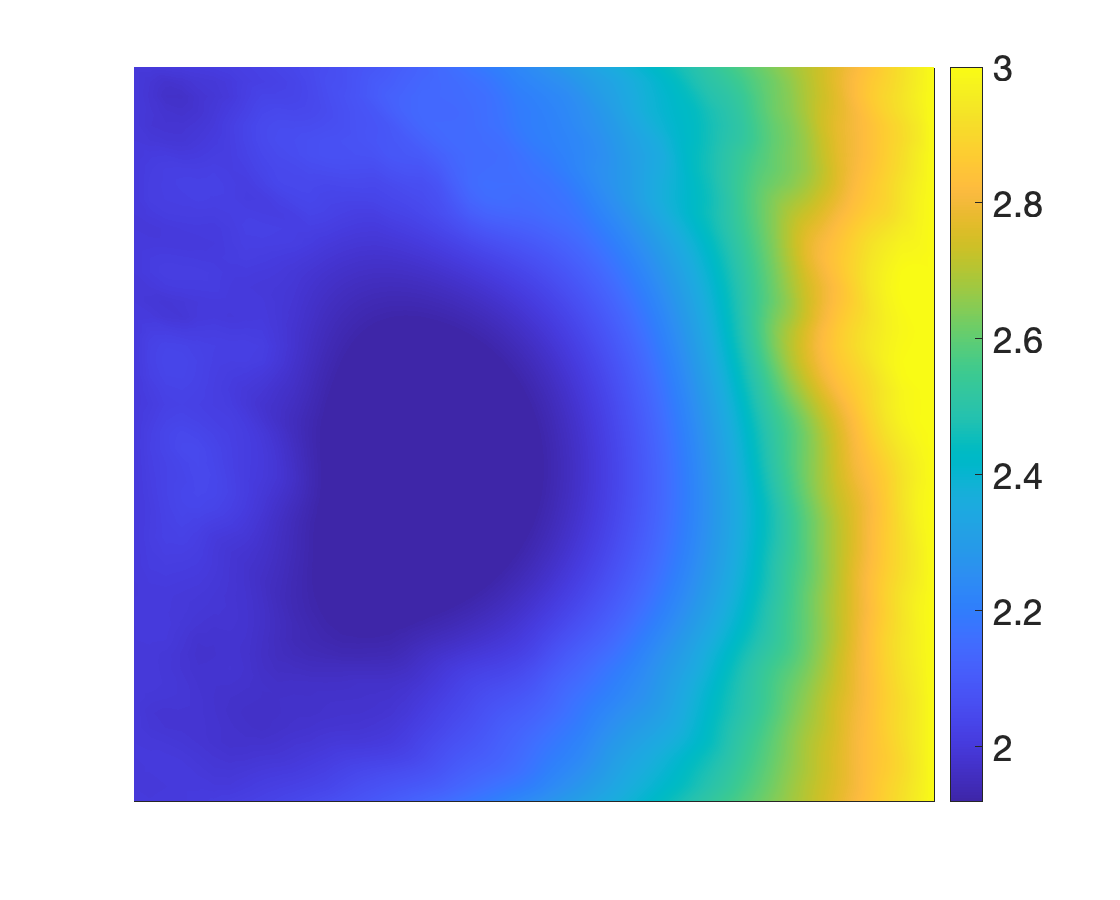} &
\includegraphics[width=0.199\textwidth]{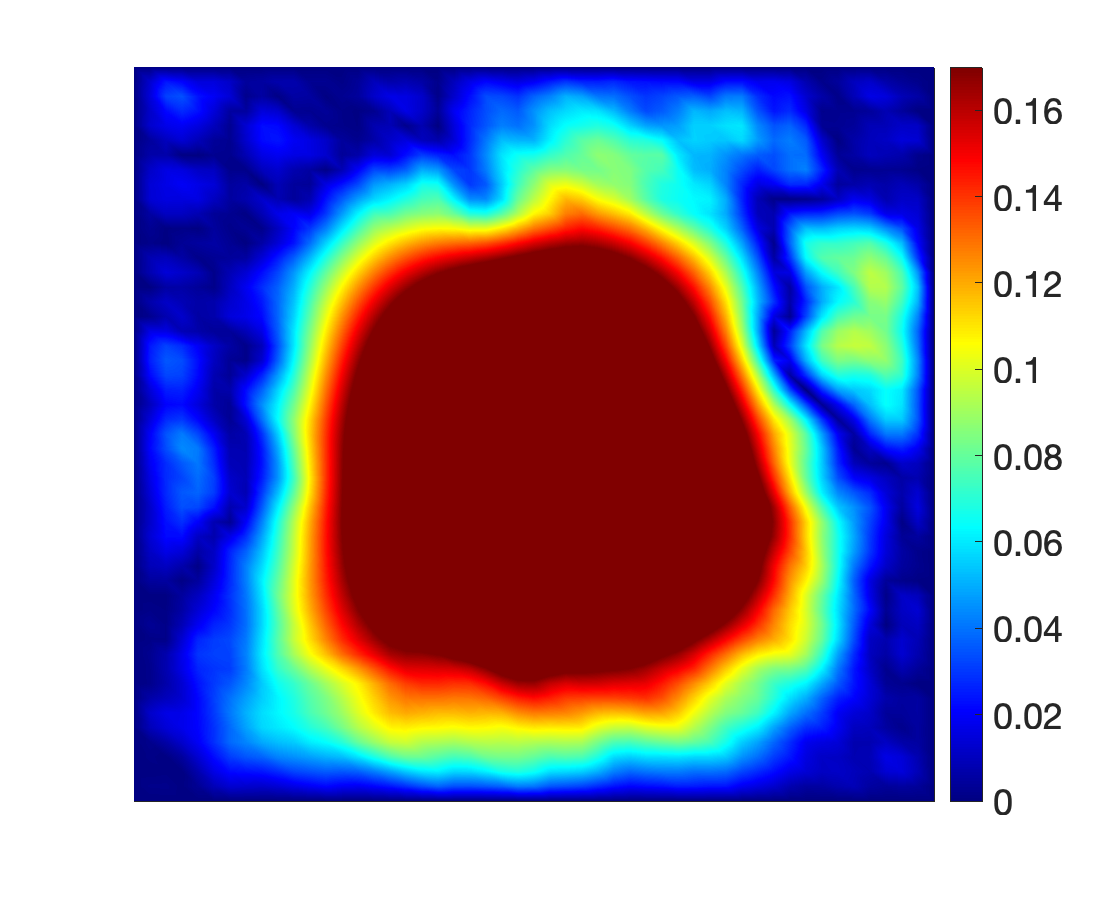} \\
(a) $A^\dag$  & (b) $\hat A$ & (c) $|\hat A-A^\dag|$ & (d) $\hat A$ & (e) $|\hat A-A^\dag|$
\end{tabular}
\caption{The reconstructions for Example \ref{exam:diri2d3} with exact data in (b) and noisy data $(\delta=10\%)$ in (d) using the FEM. From the top to bottom, the results are for $A_{11}$, $A_{12}$ and $A_{22}$, respectively.}
\label{fig:dirifem}
\end{figure}

Next we present the results obtained using the Galerkin FEM for Example \ref{exam:diri2d3}, in which the loss \eqref{obj:loss_diripartial} is discretized using the standard P1 FEM. The initial guess is set as follows: $A_{11}^0= (2+\sin(4\pi x_2)/2)(1+x_1(x_1-1)) $, $A_{12}^0= (1+\sin(2\pi x_2)^2/2)(1+x_1(x_1-1)) $, $A_{22}^0= (2+x_1^2)(1+4x_1x_2(x_1-1)(1-x_2)) $ and $\sigma_i^0=A^0\nabla u_i^\dagger$. The penalty parameters are chosen to be $\gamma_\sigma=\gamma_b=1$ and $\gamma_A=$1e-6 for reconstructions with both exact and noisy data. Fig. \ref{fig:dirifem} shows the reconstruction results using the FEM, which further demonstrates the superiority of the MLS-DNN approach over the traditional FEM for partial internal data.

\begin{figure}[htb!]
\centering
\setlength{\tabcolsep}{0em}
\begin{tabular}{ccccc}
\includegraphics[width=0.199\textwidth]{a11diri2dpex.png} &
\includegraphics[width=0.199\textwidth]{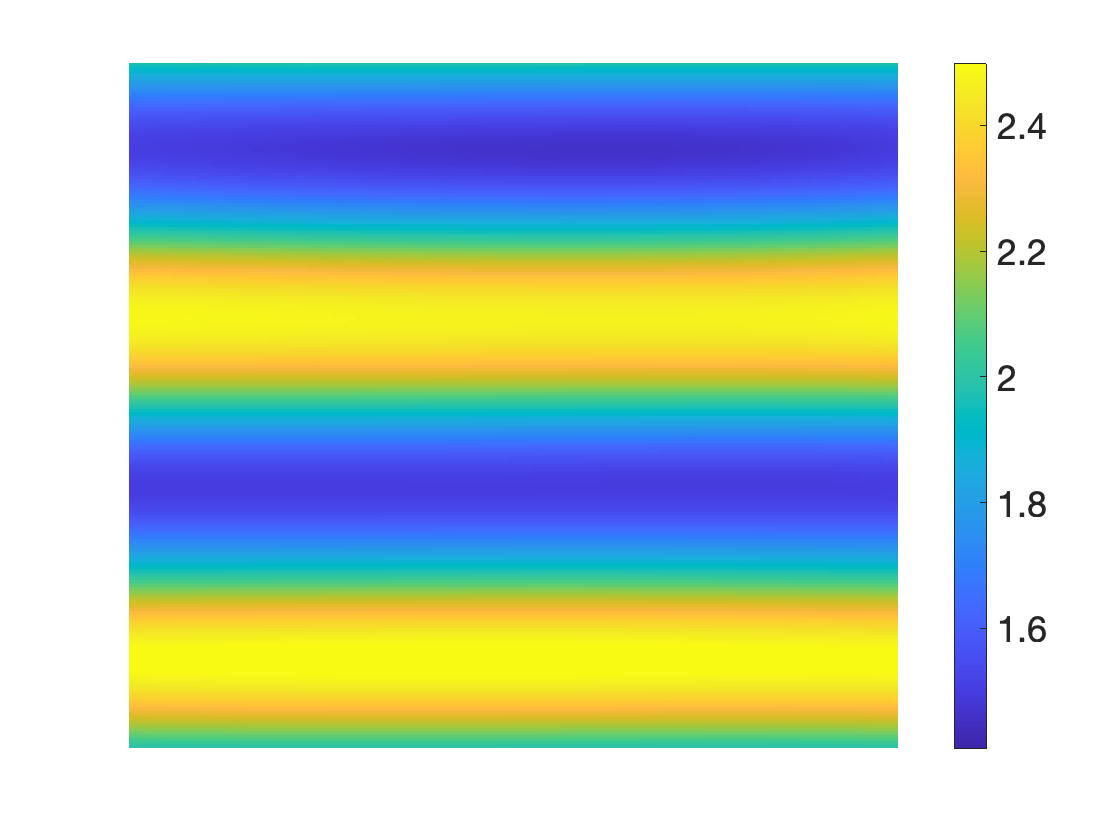} &
\includegraphics[width=0.199\textwidth]{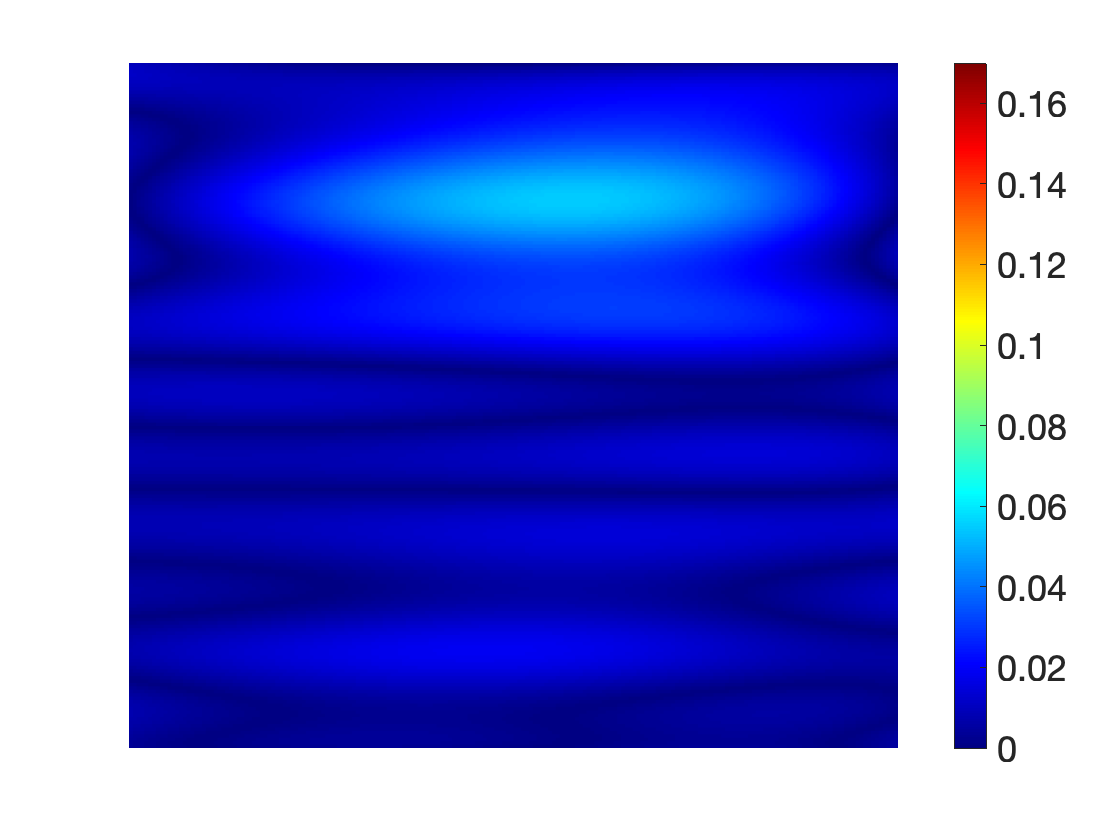} &
\includegraphics[width=0.199\textwidth]{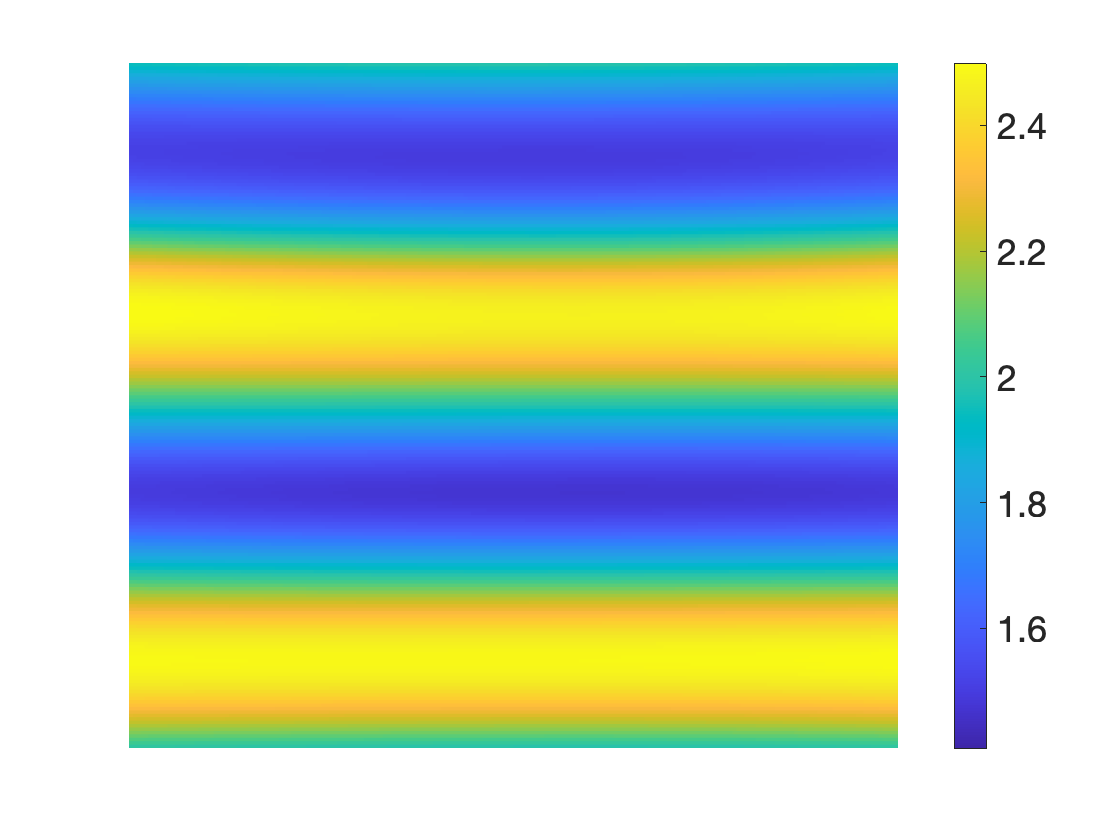} &
\includegraphics[width=0.199\textwidth]{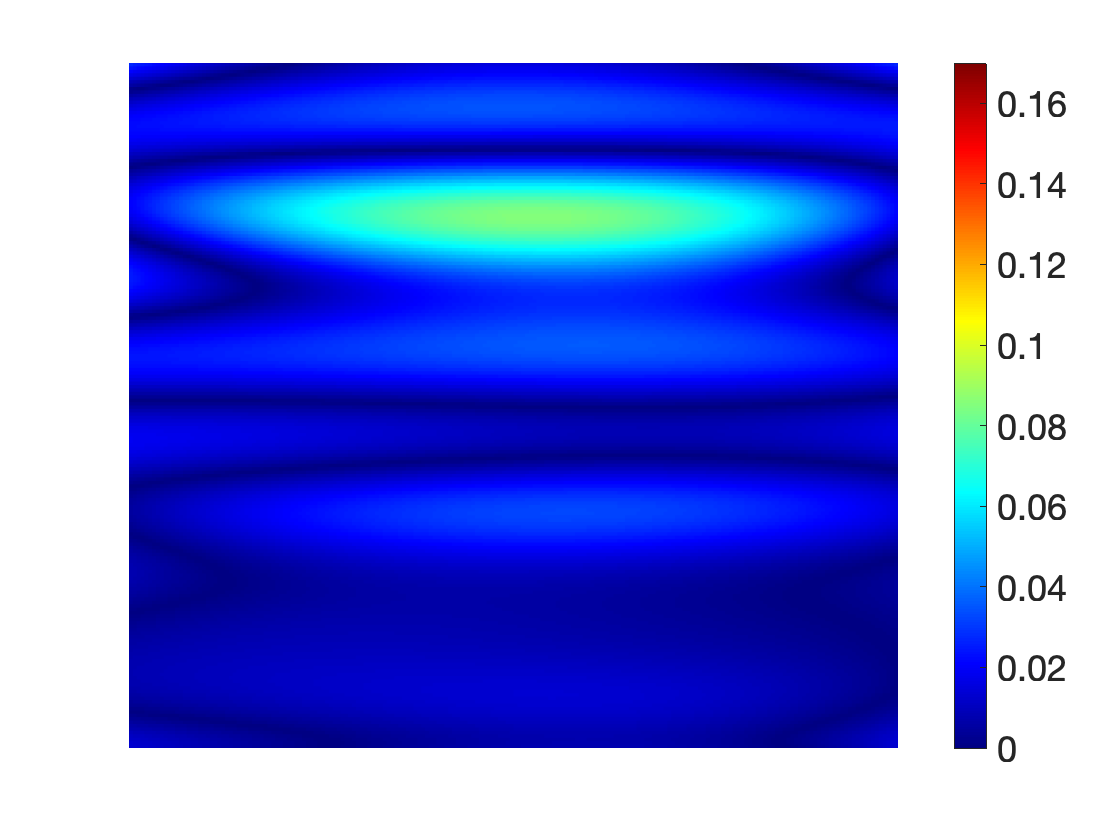} \\
\includegraphics[width=0.199\textwidth]{a12diri2dpex.png} &
\includegraphics[width=0.199\textwidth]{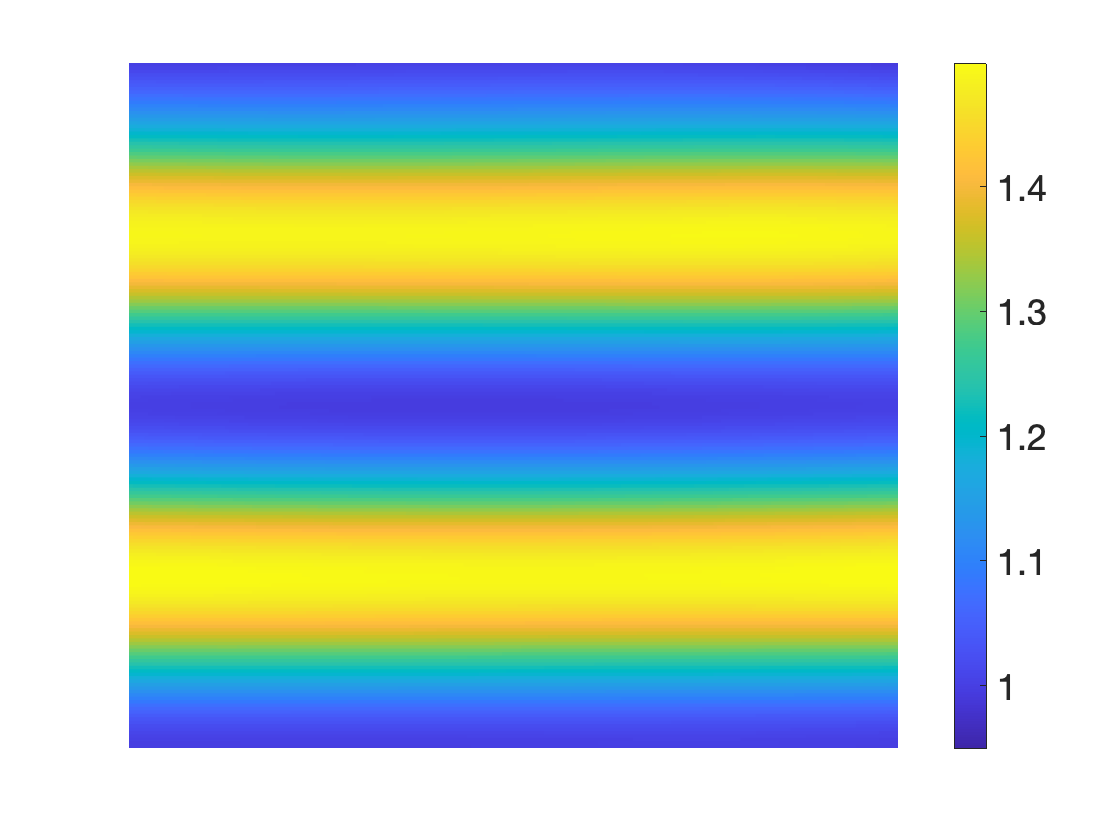} &
\includegraphics[width=0.199\textwidth]{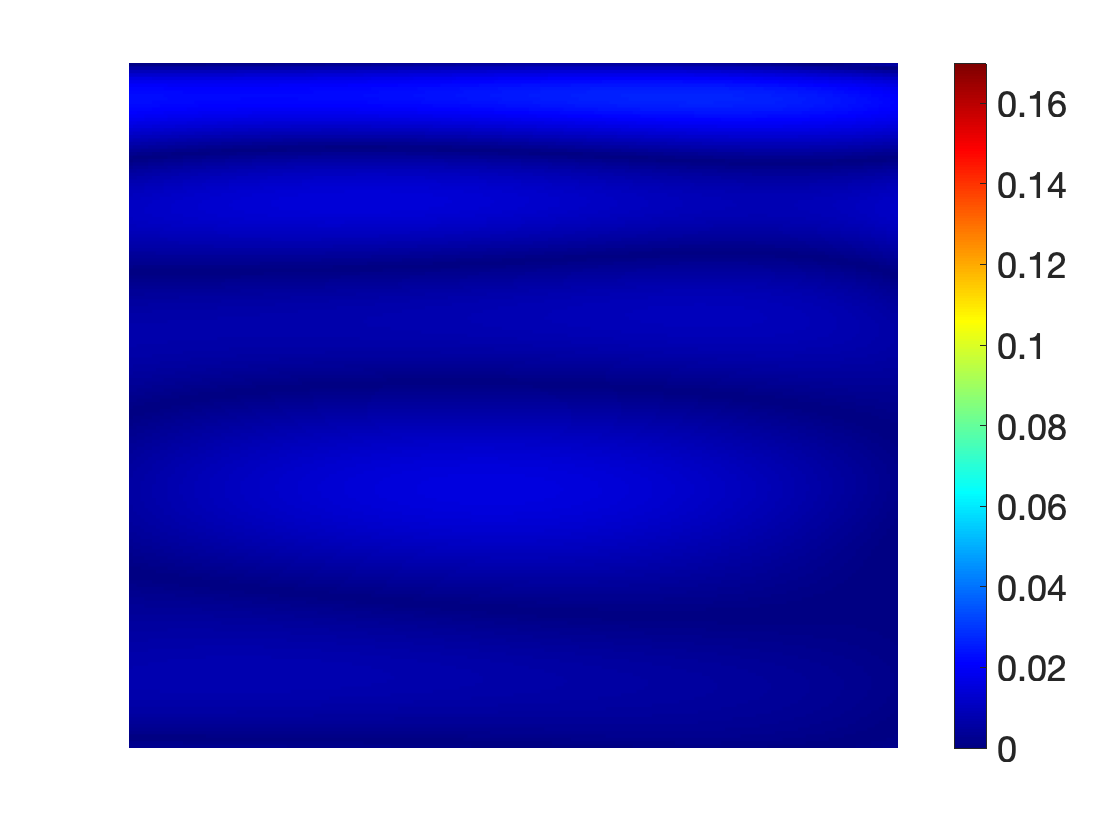} &
\includegraphics[width=0.199\textwidth]{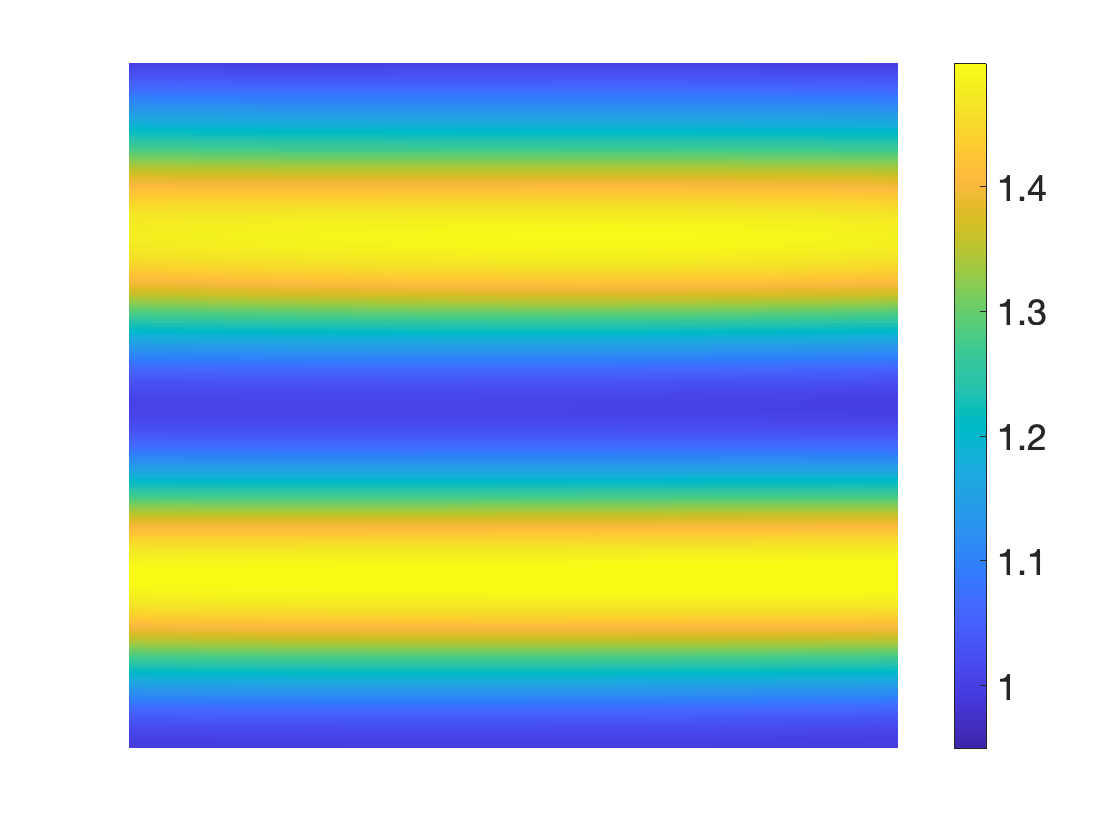} &
\includegraphics[width=0.199\textwidth]{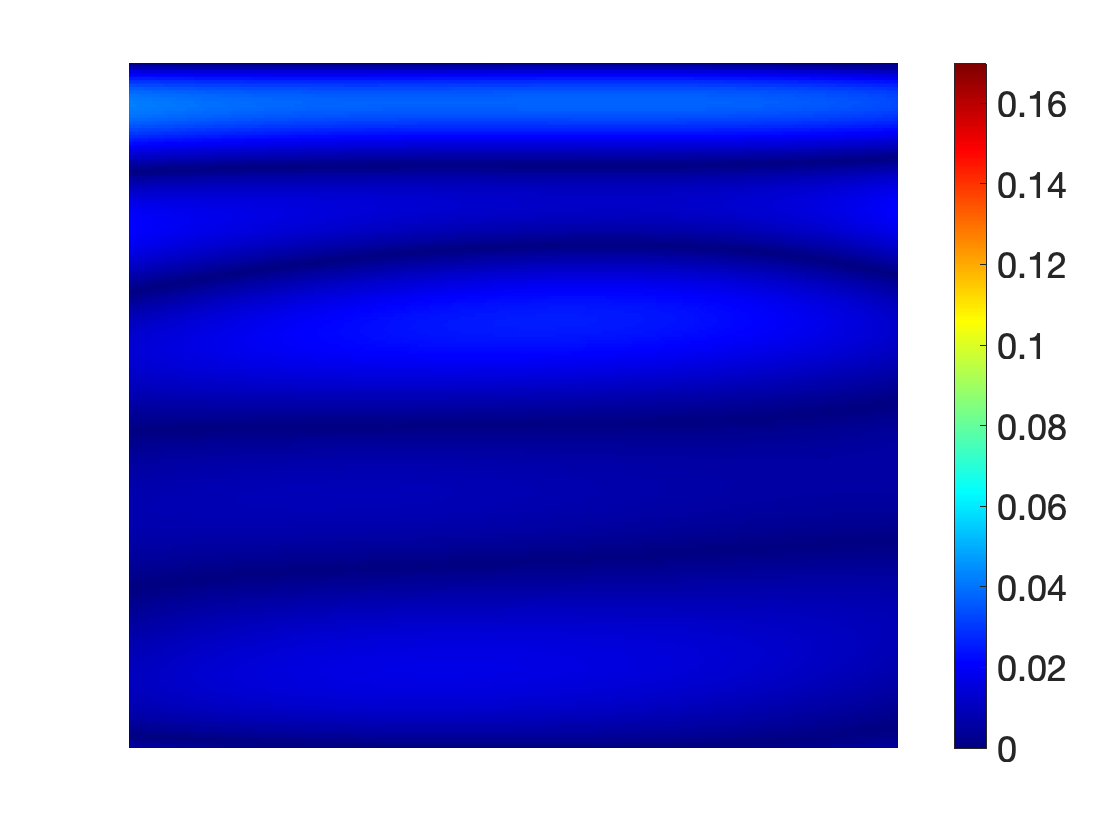} \\
\includegraphics[width=0.199\textwidth]{a22diri2dpex.png} &
\includegraphics[width=0.199\textwidth]{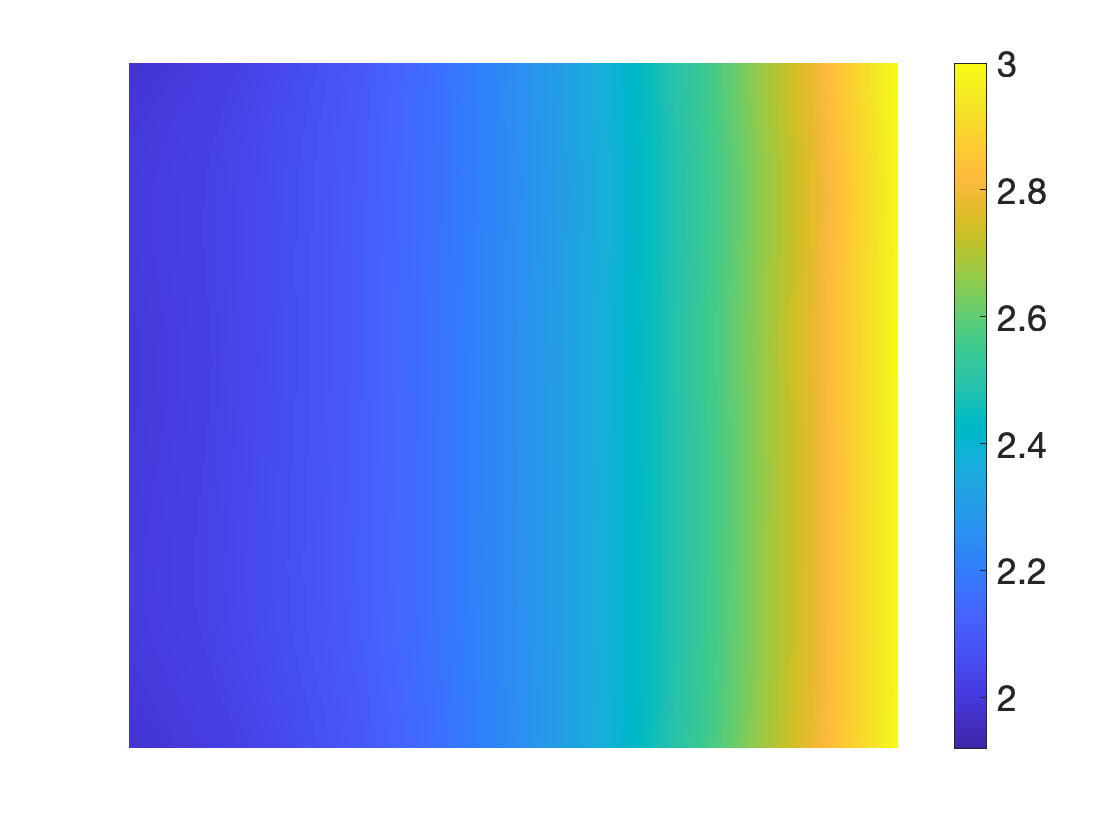} &
\includegraphics[width=0.199\textwidth]{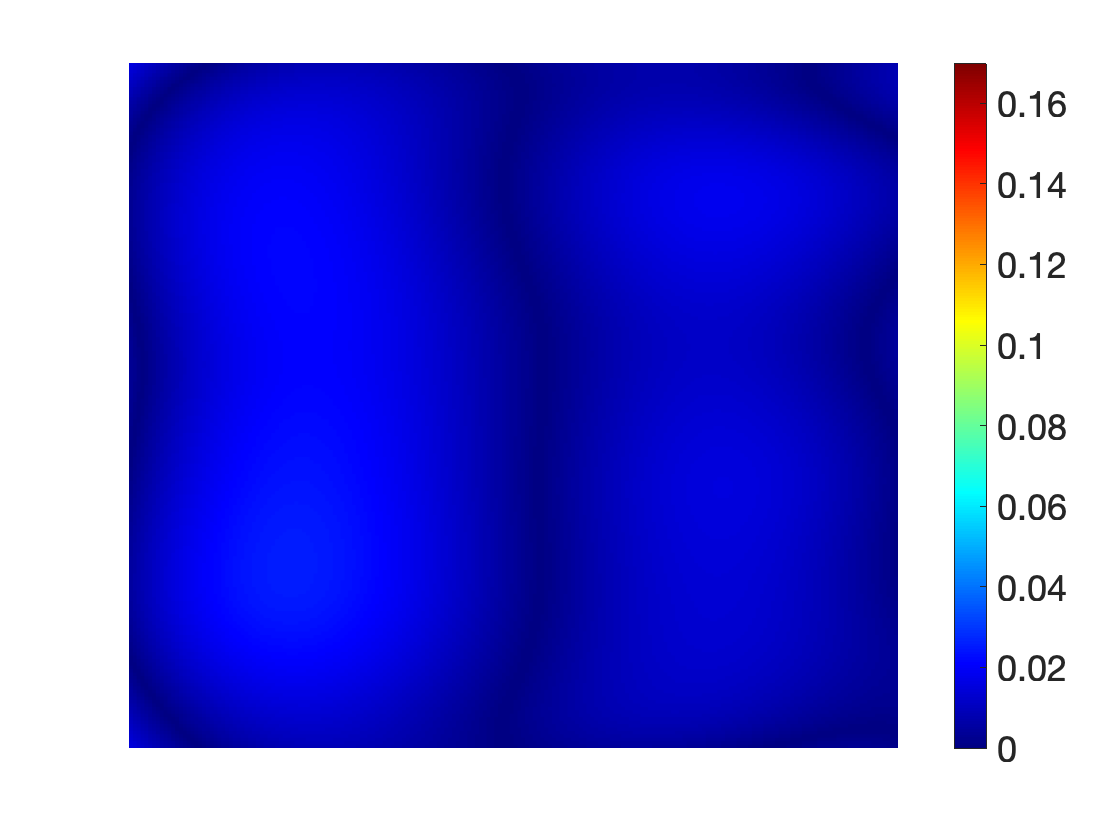} &
\includegraphics[width=0.199\textwidth]{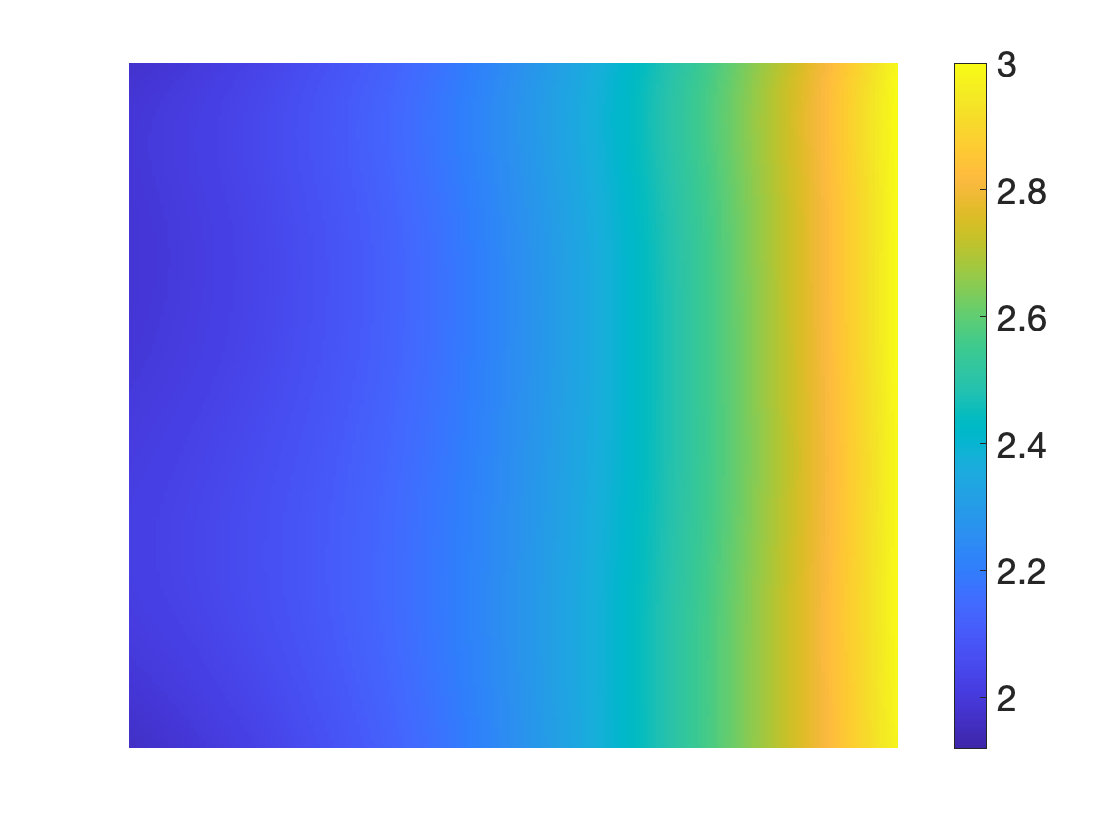} &
\includegraphics[width=0.199\textwidth]{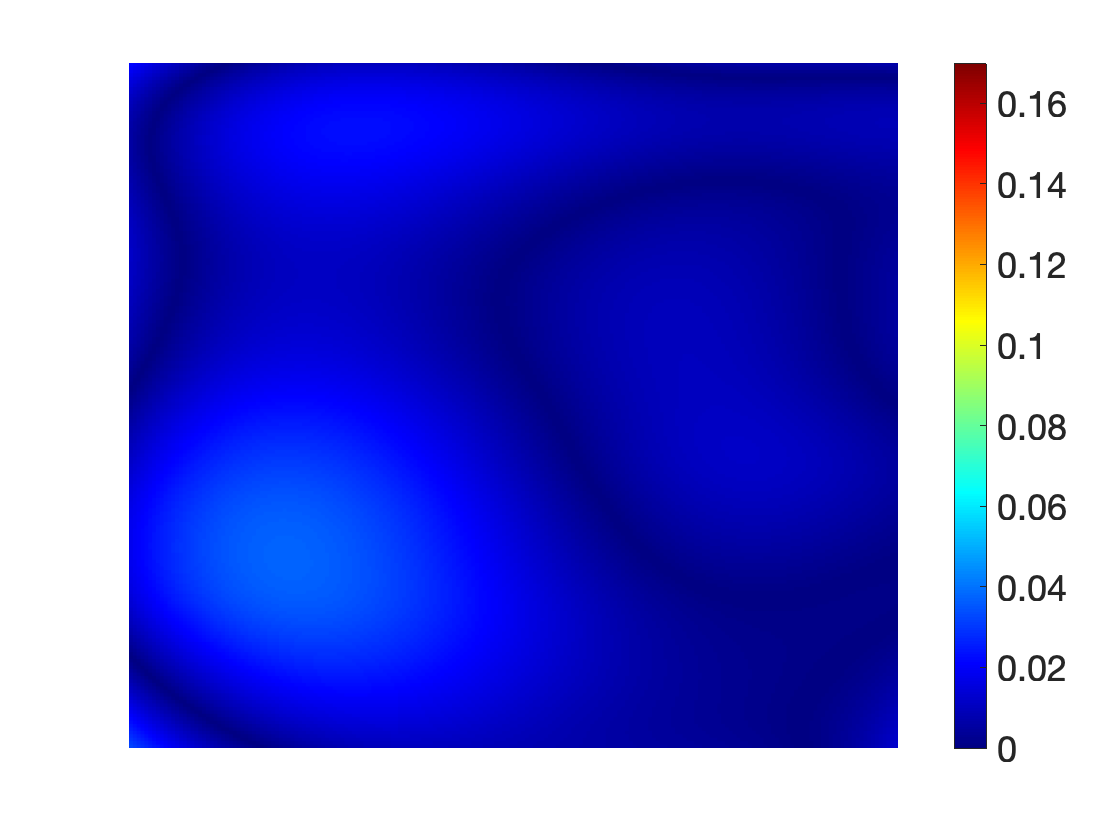} \\
(a) $A^\dag$  & (b) $\hat A$ & (c) $|\hat A-A^\dag|$ & (d) $\hat A$ & (e) $|\hat A-A^\dag|$
\end{tabular}
\caption{The reconstructions for Example \ref{exam:diri2d3} with exact data in (b) and noisy data $(\delta=10\%)$ in (d) using PINN. From the top to bottom, the results are for $A_{11}$, $A_{12}$ and $A_{22}$, respectively.}
\label{fig:diripinn}
\end{figure}

Also we present the results by PINN for Example \ref{exam:diri2d3}, which involves minimizing an empirical version of the following population loss
\begin{equation}\label{obj:diripinnloss}
    \begin{aligned}
            J_{\gamma}(\theta,\kappa)=\sum_{i=1}^{N}&\left(\| \nabla u_{i,\kappa} - \nabla z_i^{\delta} \|_{L^2(\omega)^d}^2+\gamma_{1}\| \nabla \cdot (A_{\theta} \nabla u_{i,\kappa}) + f_i \|_{L^2(\Omega)}^2\right.\\
            &\left.+\gamma_2 \|u_{i,\kappa}-g_i\|_{L^2(\partial \Omega)}^2+\gamma_3 \| A_{\theta}\|_{L^2(\Omega)^{d,d}}^2+\gamma_4 \|A_\theta-A^\dagger\|_{L^2(\partial \Omega)^{d,d}}^2 \right).
    \end{aligned}
\end{equation}
Similar to the Neumann case (i.e., Example \ref{exam:neu2d3}), the knowledge of the conductivity tensor $A^\dagger|_{\partial\Omega}$ is essential for achieving a reliable reconstruction of $A$. Fig. \ref{fig:diripinn} shows that the PINN reconstructions are comparable to that by MLS-DNN in terms of the reconstruction accuracy for both exact and noisy data. Similar to the Neumann case, PINN takes longer per epoch to compute $J_{\gamma}(\theta,\kappa)$, primarily due to the storage and computation of second-order spatial derivatives in the term $\| \nabla \cdot (A_{\theta} \nabla u_{i,\kappa}) + f_i \|_{L^2(\Omega)}^2$.

\begin{figure}[htb!]
\centering
\setlength{\tabcolsep}{0em}
\begin{tabular}{ccccc}
\includegraphics[width=0.199\textwidth]{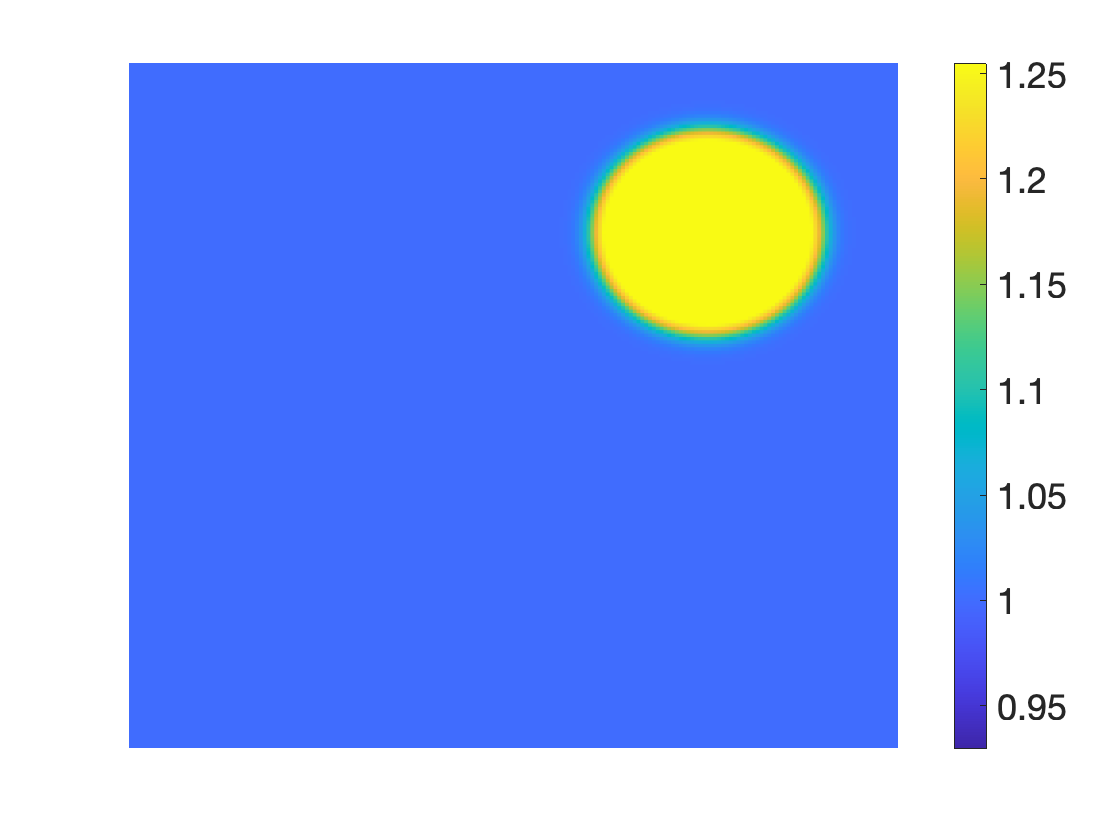} &
\includegraphics[width=0.199\textwidth]{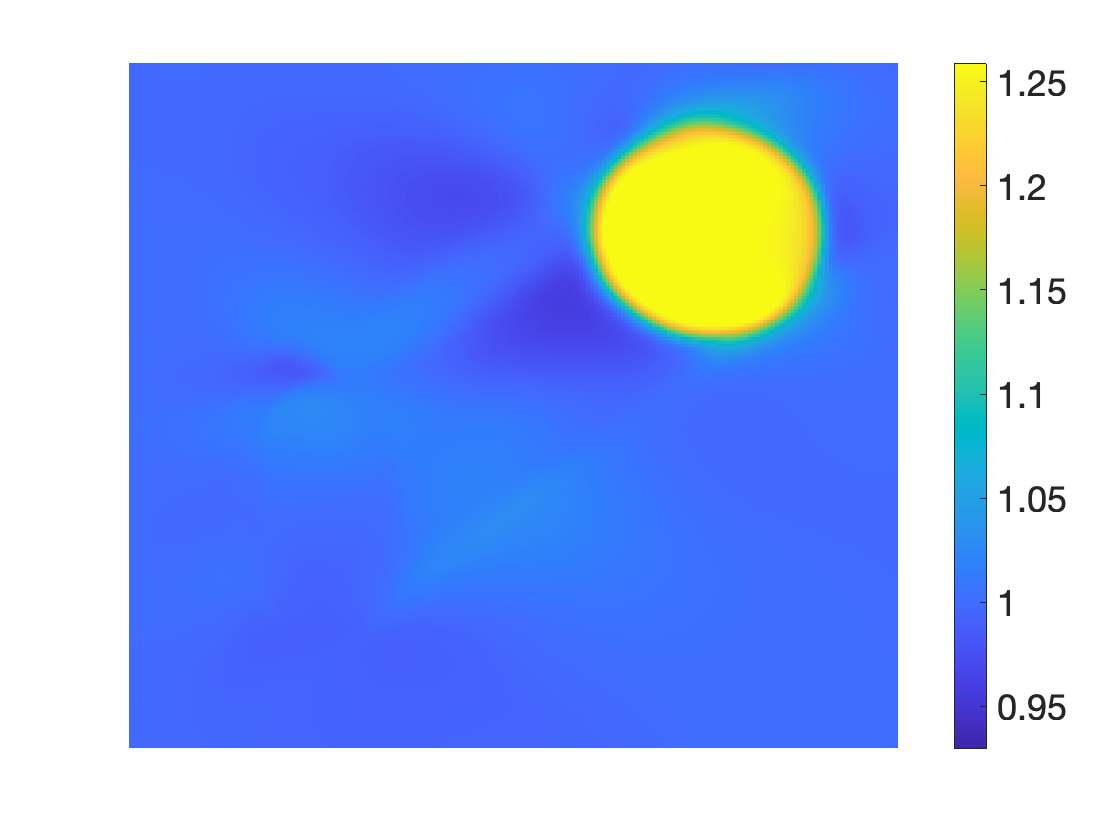} &
\includegraphics[width=0.199\textwidth]{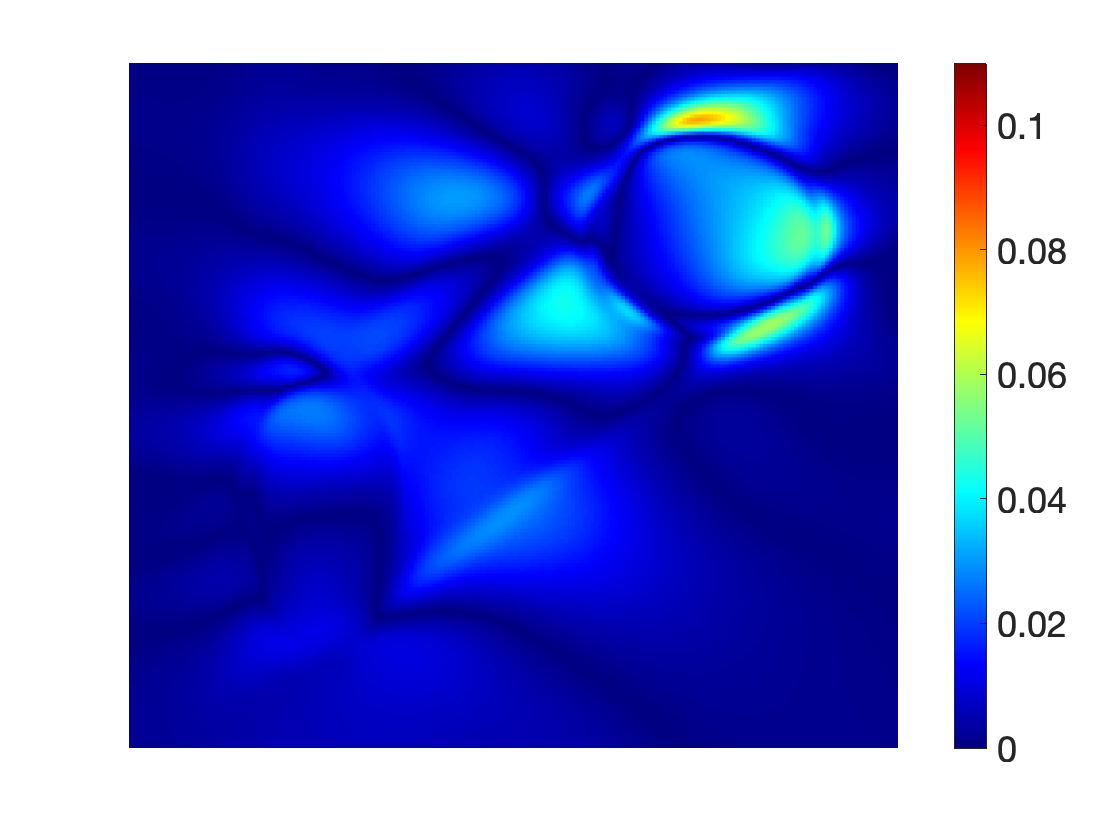} &
\includegraphics[width=0.199\textwidth]{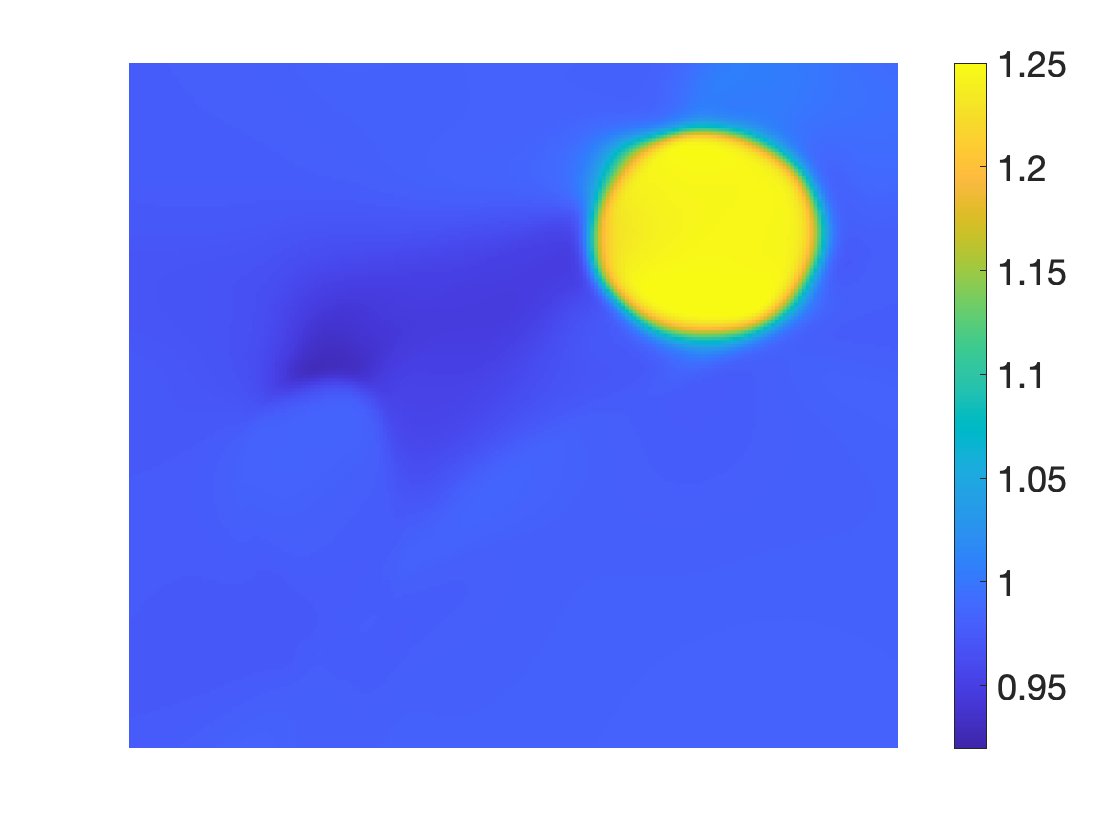} &
\includegraphics[width=0.199\textwidth]{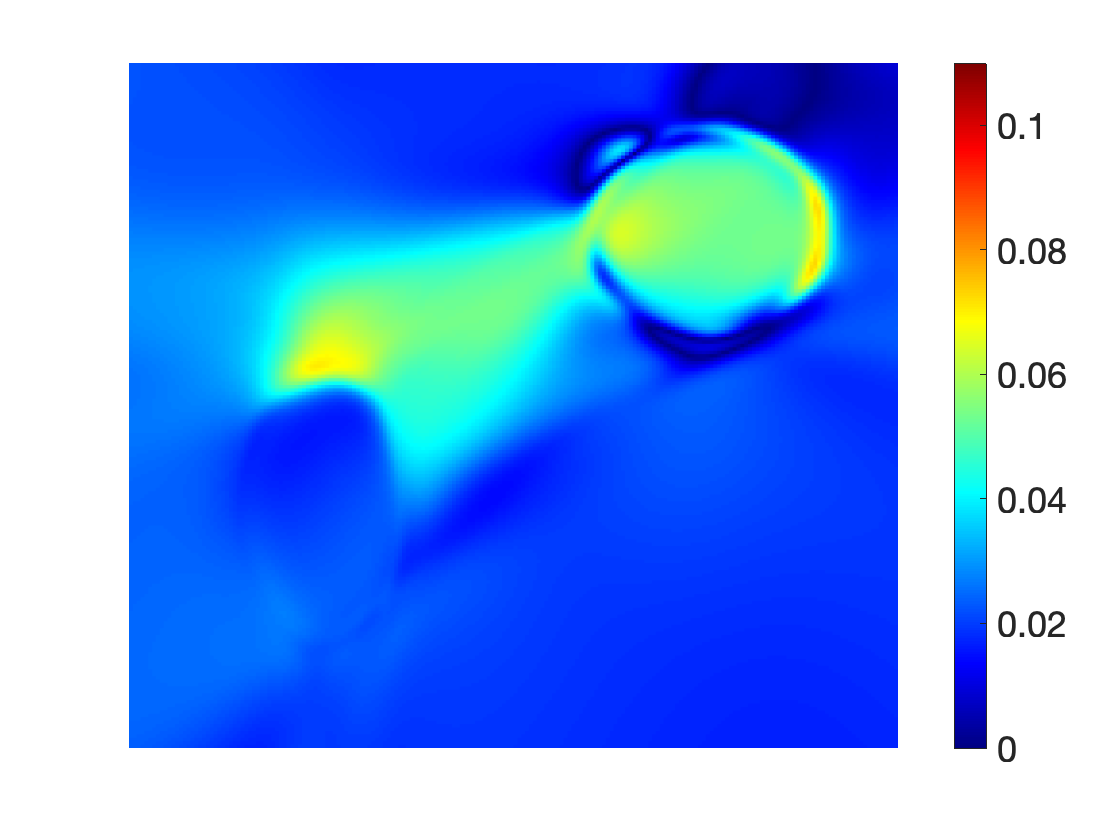} \\
\includegraphics[width=0.199\textwidth]{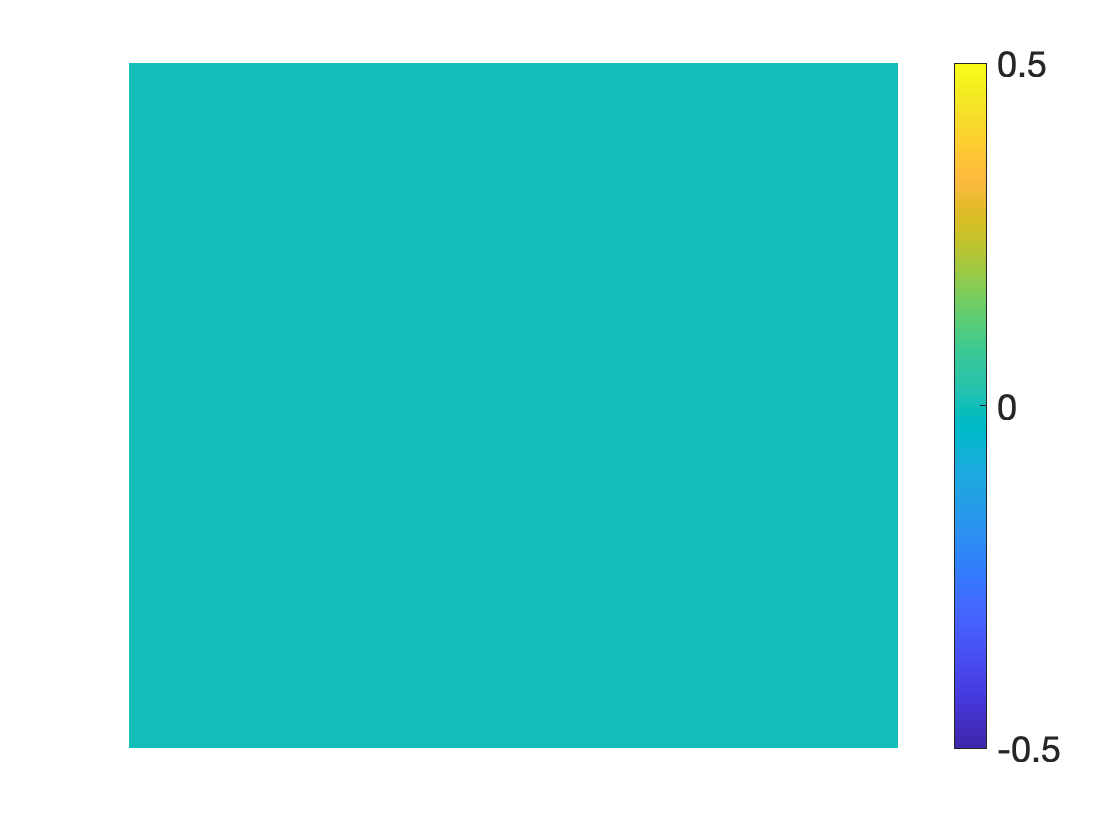} &
\includegraphics[width=0.199\textwidth]{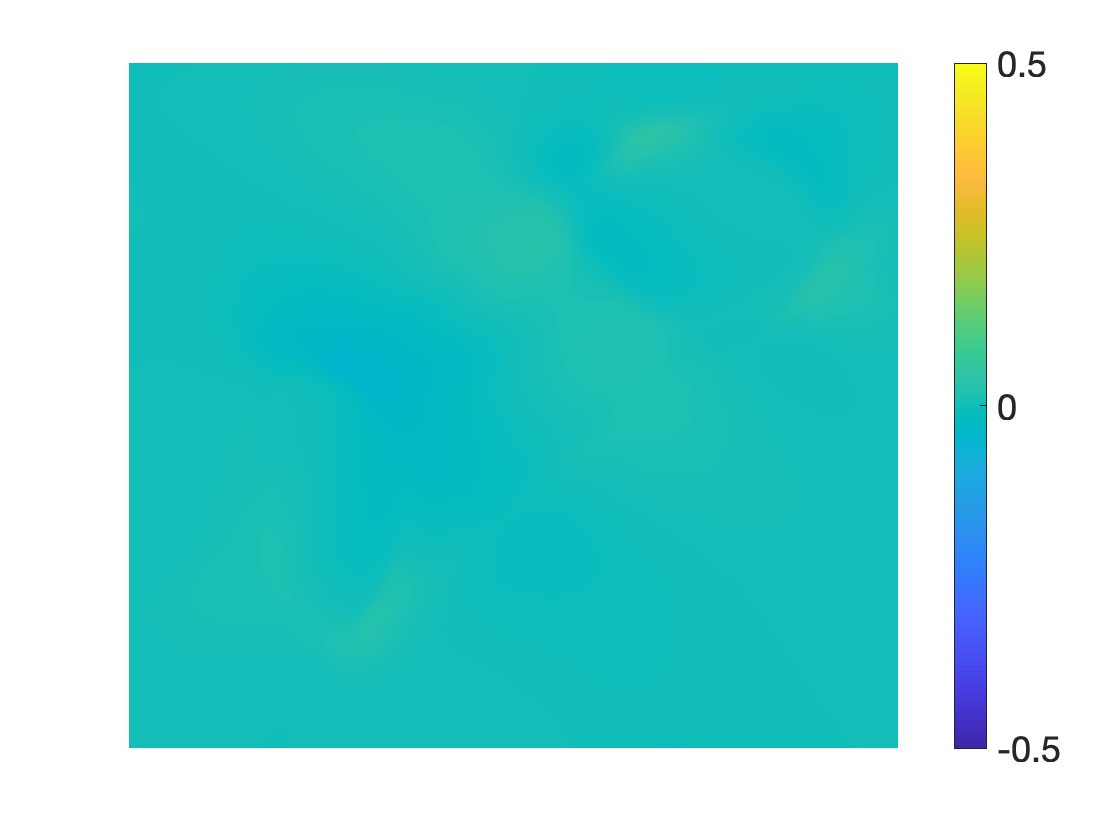} &
\includegraphics[width=0.199\textwidth]{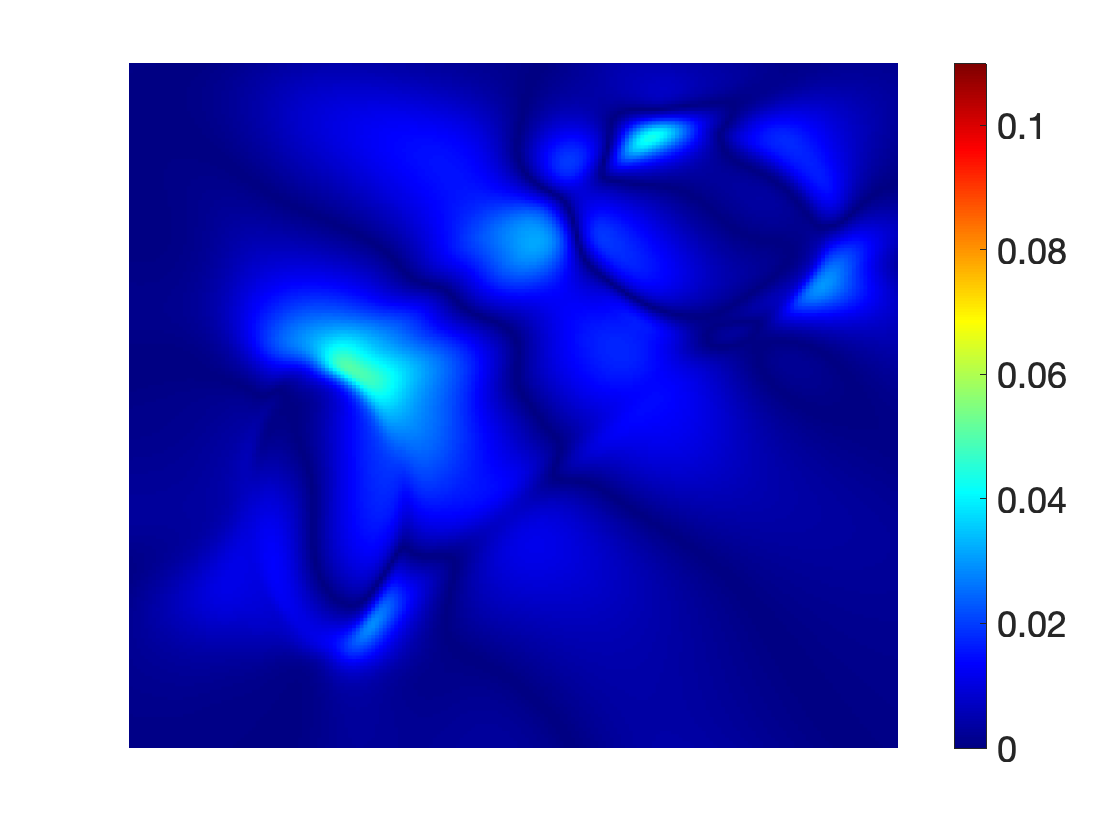} &
\includegraphics[width=0.199\textwidth]{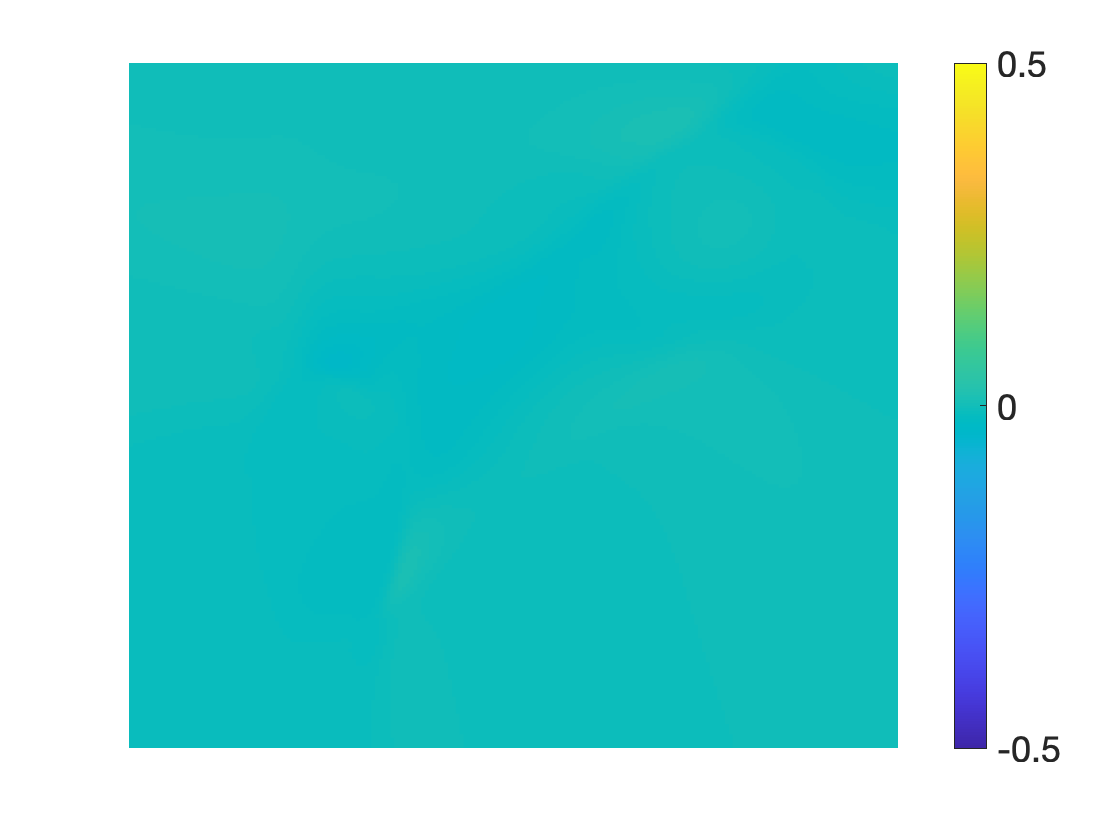} &
\includegraphics[width=0.199\textwidth]{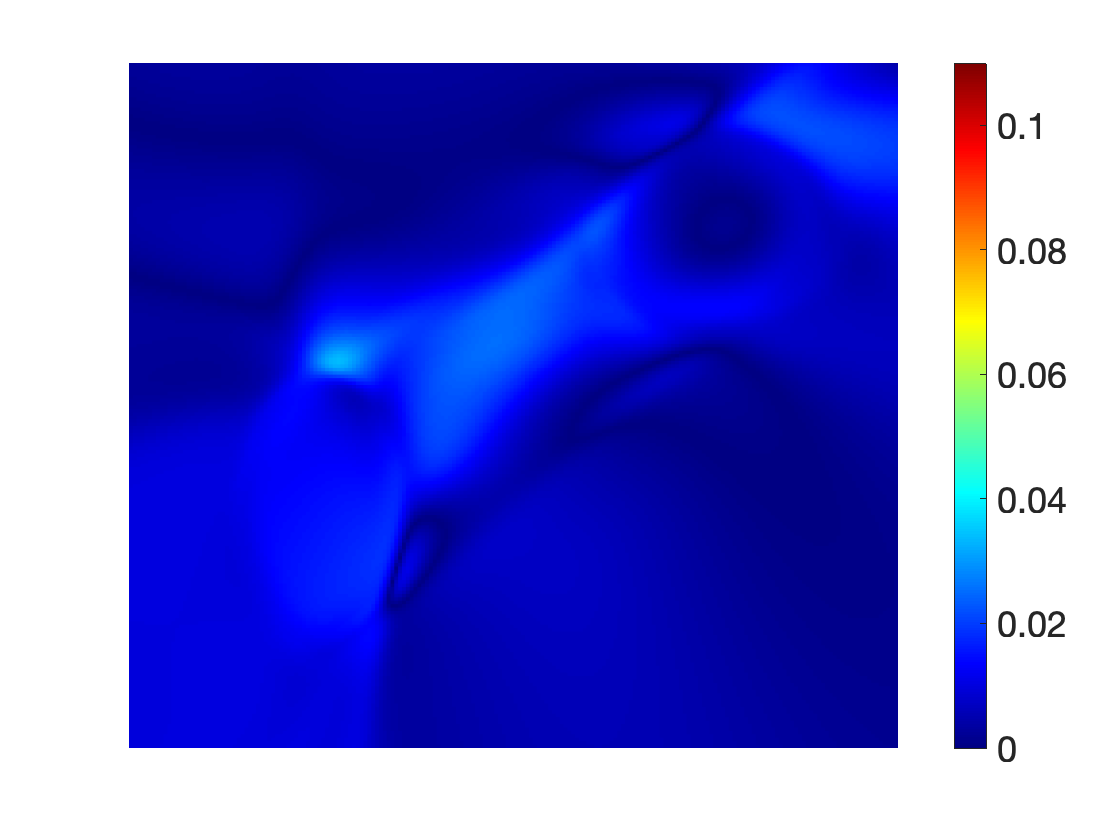} \\
\includegraphics[width=0.199\textwidth]{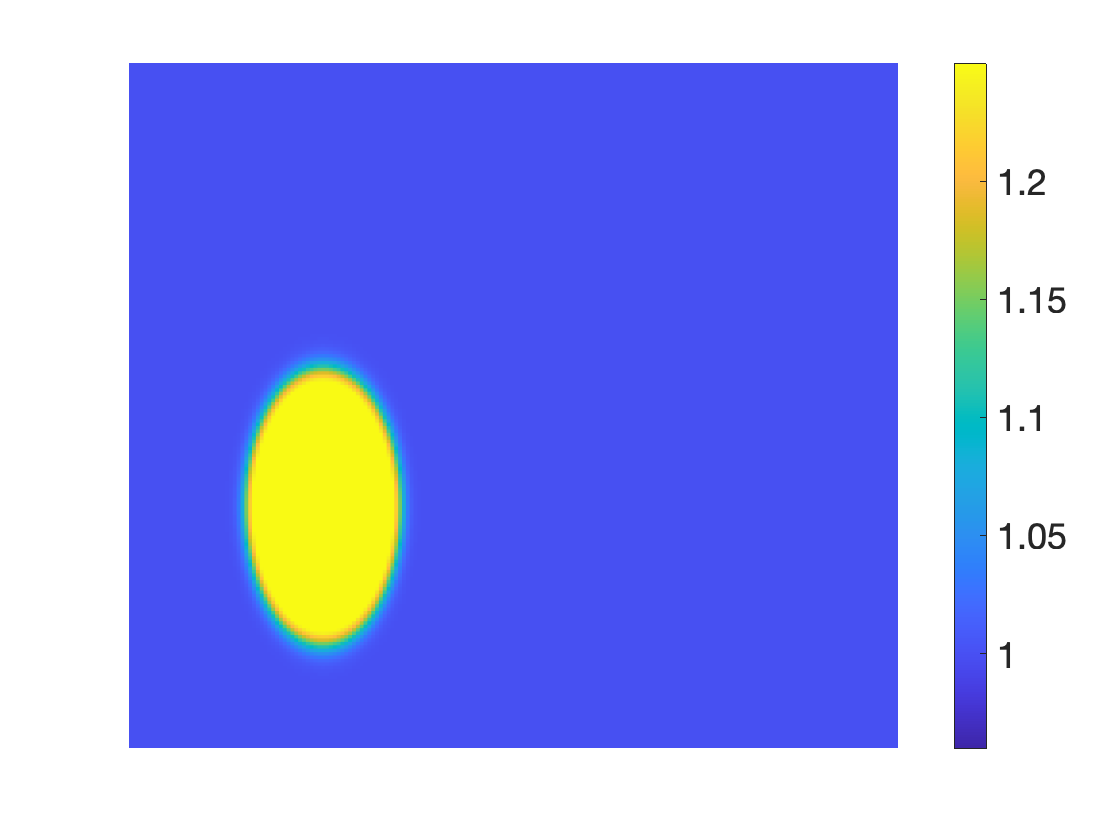} &
\includegraphics[width=0.199\textwidth]{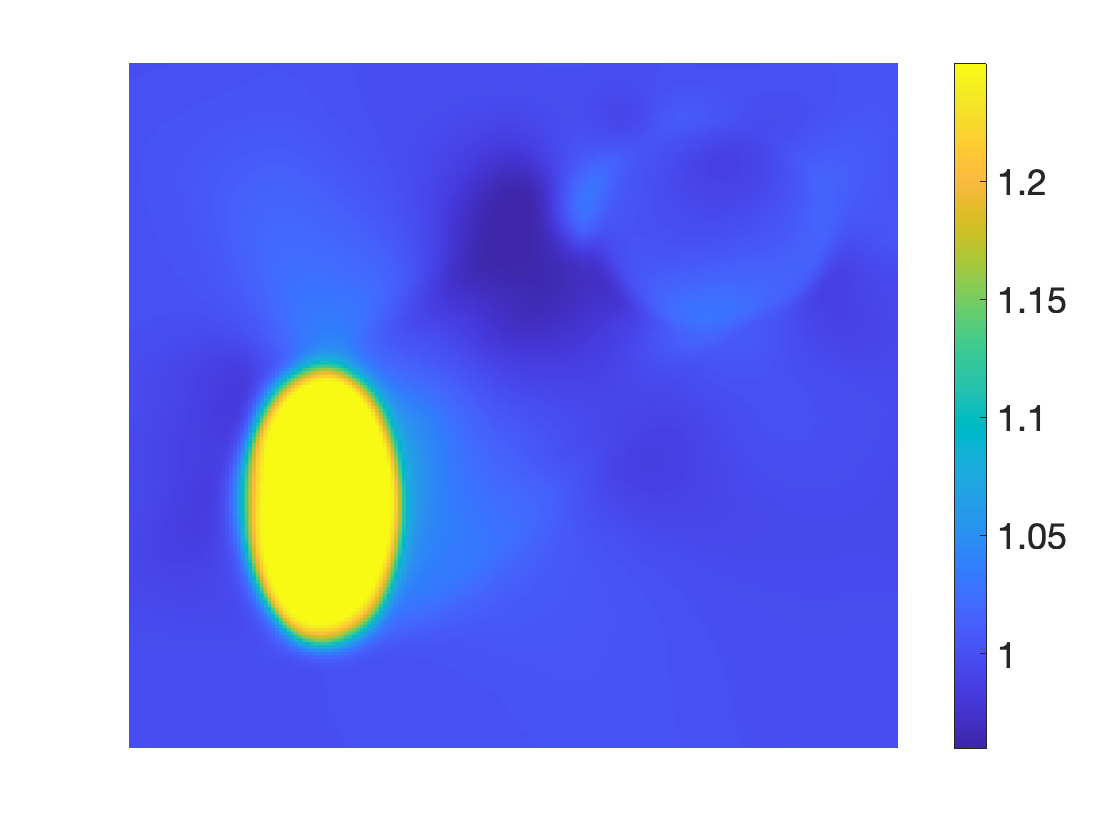} &
\includegraphics[width=0.199\textwidth]{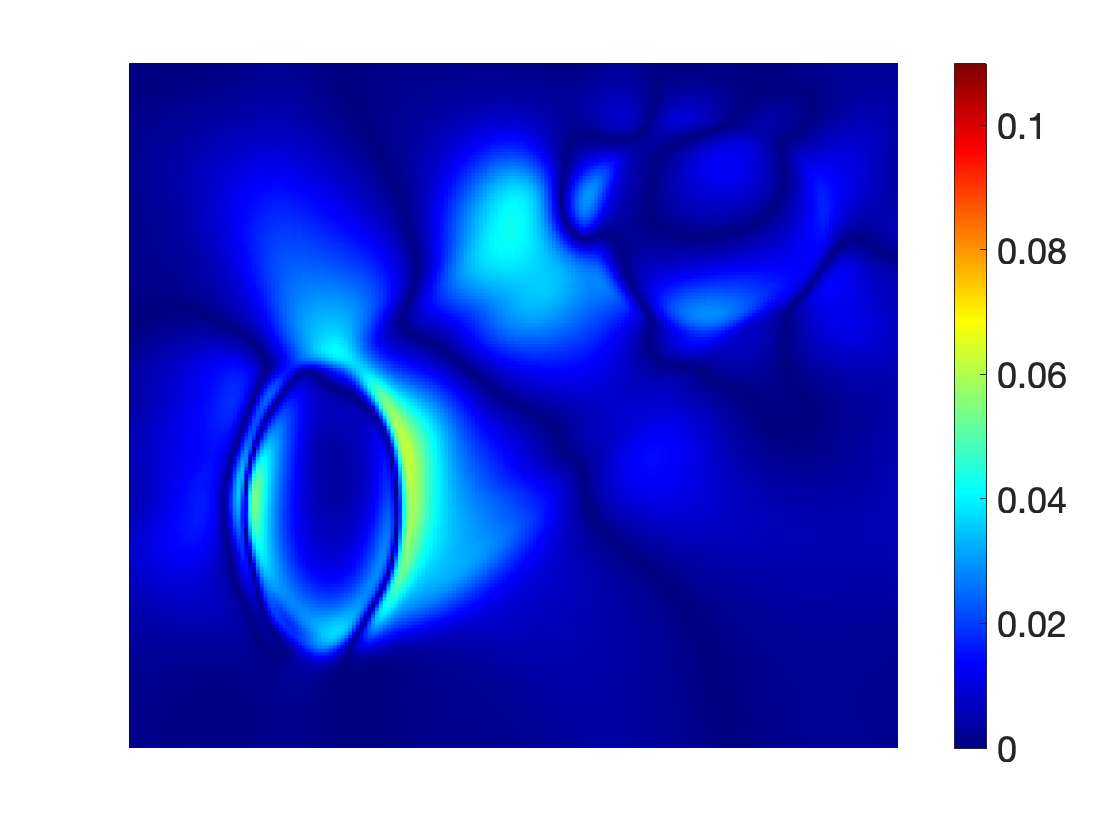} &
\includegraphics[width=0.199\textwidth]{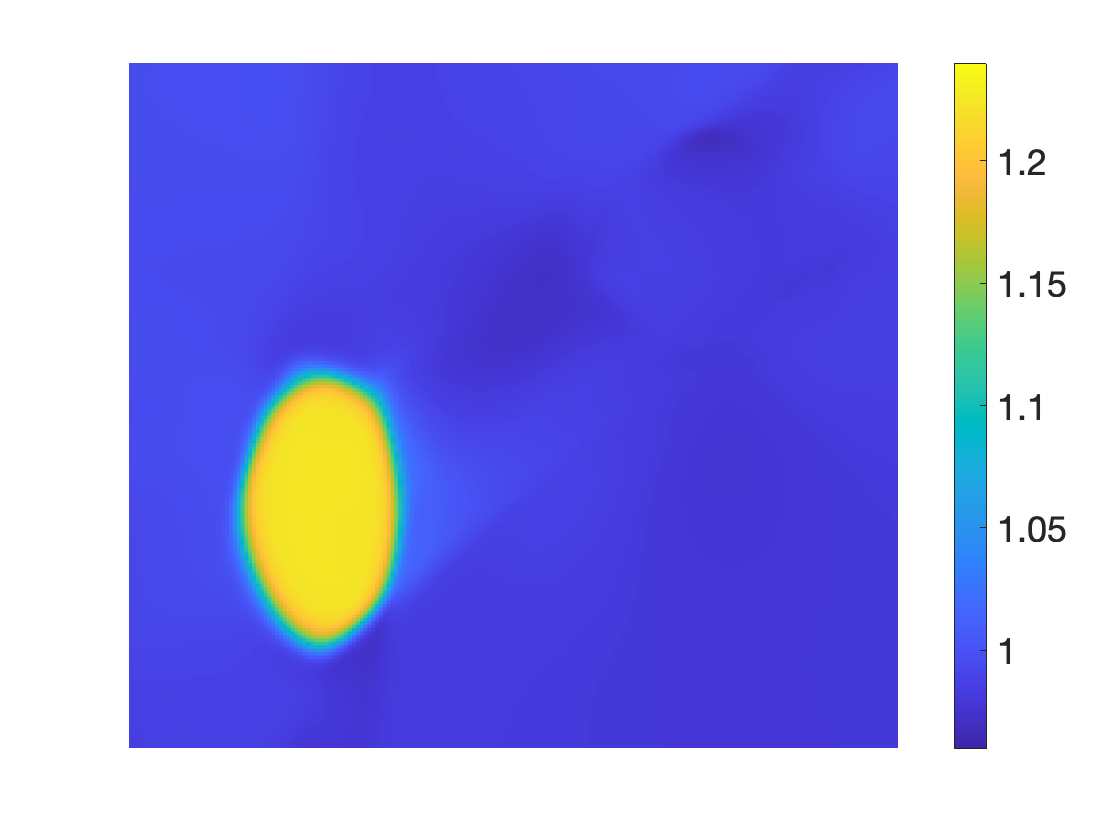} &
\includegraphics[width=0.199\textwidth]{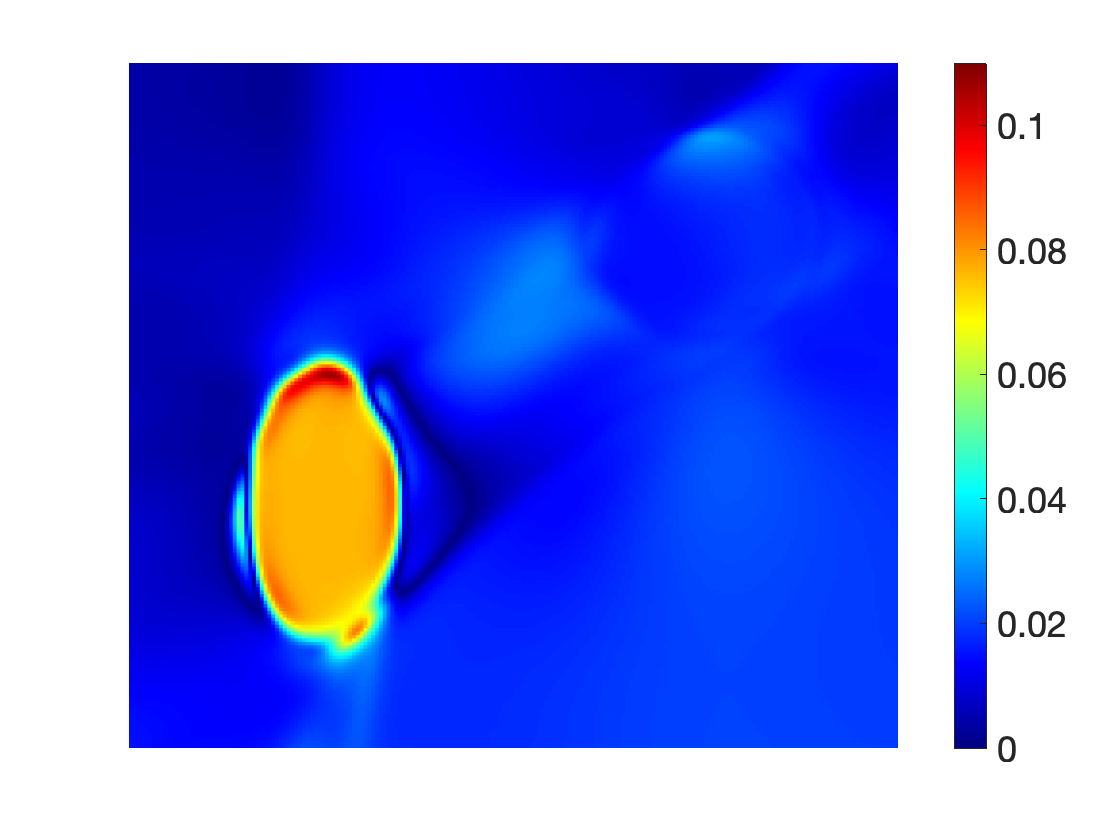} \\
(a) $A^\dag$  & (b) $\hat A$ & (c) $|\hat A-A^\dag|$ & (d) $\hat A$ & (e) $|\hat A-A^\dag|$
\end{tabular}
\caption{The reconstructions for Example \ref{exam:diri2d4} with exact data in (b) and noisy data $(\delta=5\%)$ in (d). From the top to bottom, the results are for $A_{11}$, $A_{12}$ and $A_{22}$, respectively.}
\label{fig:diri2d4}
\end{figure}

The fourth example is about recovering a diagonal conductivity matrix $A^\dag$ with nearly piecewise constant component functions.
\begin{example}
The domain $\Omega = (0,1)^2$, $A^\dag= \begin{pmatrix}
    \zeta_1(x_1,x_2) & 0\\
   0 & \zeta_2(x_1,x_2)\\
\end{pmatrix},$ where $\zeta_1(x_1,x_2)=1+0.3/(1+\exp(400((x_1-0.75)^2+(x_2-0.75)^2-0.0225)))$ and $\zeta_2(x_1,x_2)=1+0.3/(1+\exp(250(4(x_1-0.25)^2+(x_2-0.35)^2-0.04)))$. $u_1^\dag=x_1+x_2+\frac{1}{3}(x_1^3+x_2^3)$, $u_2^\dag=x_1-x_2+\frac{1}{3}(x_1^3-x_2^3)$, $u_3^\dag=-u_1^\dagger$, $u_4^\dag=-u_2^\dagger$.
\label{exam:diri2d4}
\end{example}
Similar to Example \ref{exam:neu2d4}, an additional total variation penalty,
i.e., $\gamma_{tv}\sum_{i,j=1}^2| A_{ij}|_{\rm TV}$, is added to the loss $J_{\bsgamma}(\theta,\kappa)$ in order to promote piecewise constancy in the reconstruction. Fig. \ref{fig:diri2d4} presents the reconstruction results for both exact and noisy data ($\delta = 5\%$), where $\gamma_{tv}=\ $8e-4 and 5e-3, respectively. For the exact data, the reconstruction can accurately capture the location and shape of the circular and elliptic bumps. When 5\% noise is present, the reconstruction remains accurate, showing only a slight underestimation of the bump peaks and minor artifacts in the background.

The fifth example is about recovering a 3D anisotropic conductivity matrix with mixed oscillatory and polynomial entries.
\begin{example}\label{exam:diri3d1}
    The domain $\Omega=(0,1)^3$, $A^\dagger = \begin{pmatrix}
    3& 1+\frac{\sin(4\pi x_3)}{2}& 1+\frac{\sin(4\pi x_2)}{2}\\
    1+\frac{\sin(4\pi x_3)}{2}& 2& 1+\frac{\sin(4\pi x_1)}{2}\\
    1+\frac{\sin(4\pi x_2)}{2}& 1+\frac{\sin(4\pi x_1)}{2}& 2+\frac{\sin(2\pi x_1)\sin(2\pi x_2)}{2}\\
    \end{pmatrix}$, $u_1^\dag=x_1+x_2+x_3+\frac{1}{3}(x_1^3+x_2^3+x_3^3)$, $u_2^\dag=x_1-x_2+x_3+\frac{1}{3}(x_1^3-x_2^3+x_3^3)$, $u_3^\dag=x_1+x_2-x_3+\frac{1}{3}(x_1^3+x_2^3-x_3^3)$, $u_4^\dag=-x_1+x_2+x_3+\frac{1}{3}(-x_1^3+x_2^3+x_3^3)$, $u_5^\dagger=-u_3^\dagger$ and $u_6^\dagger=-u_2^\dagger$.
\end{example}

\begin{figure}[htb!]
\centering
\setlength{\tabcolsep}{0em}
\begin{tabular}{ccccc}
\includegraphics[width=0.199\textwidth]{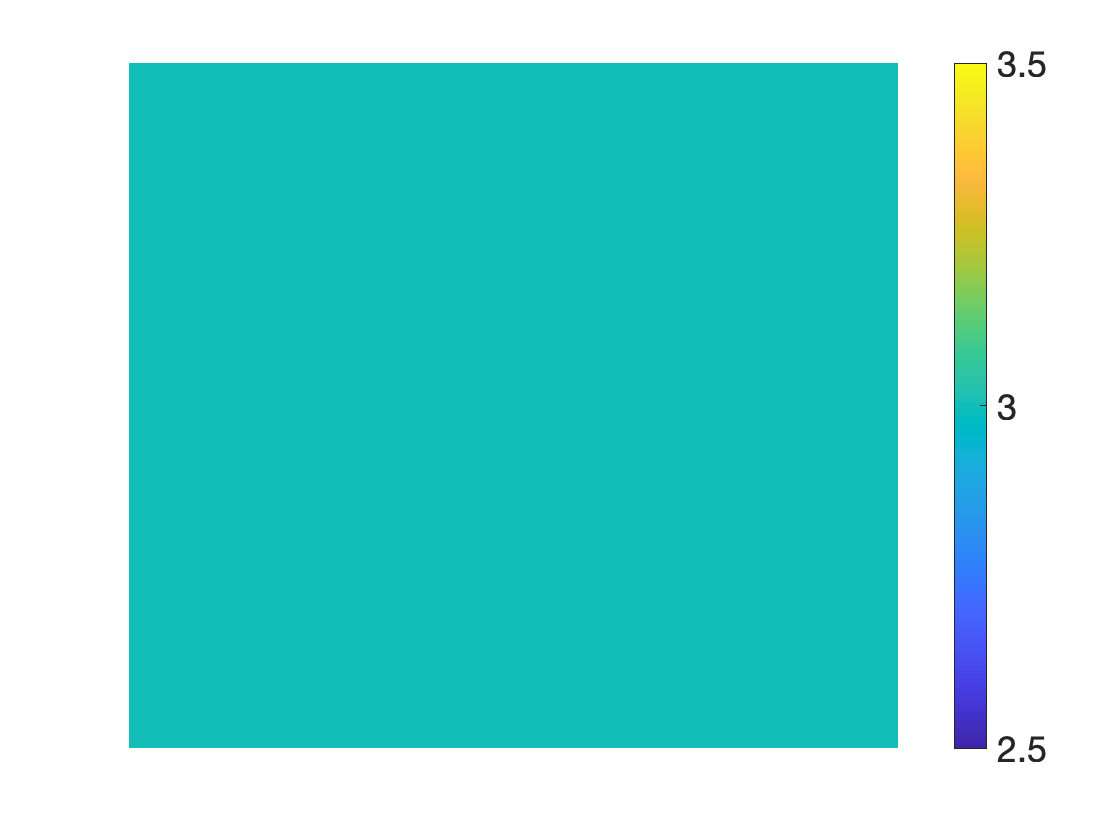} &
\includegraphics[width=0.199\textwidth]{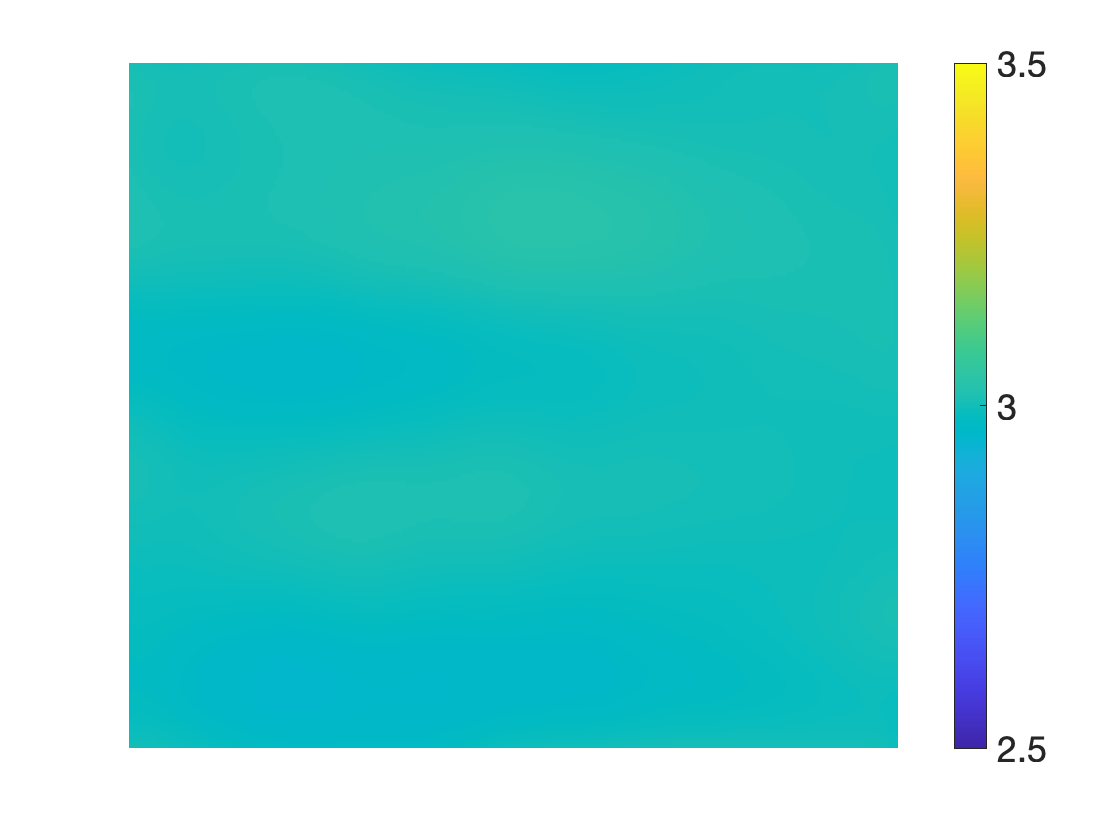} &
\includegraphics[width=0.199\textwidth]{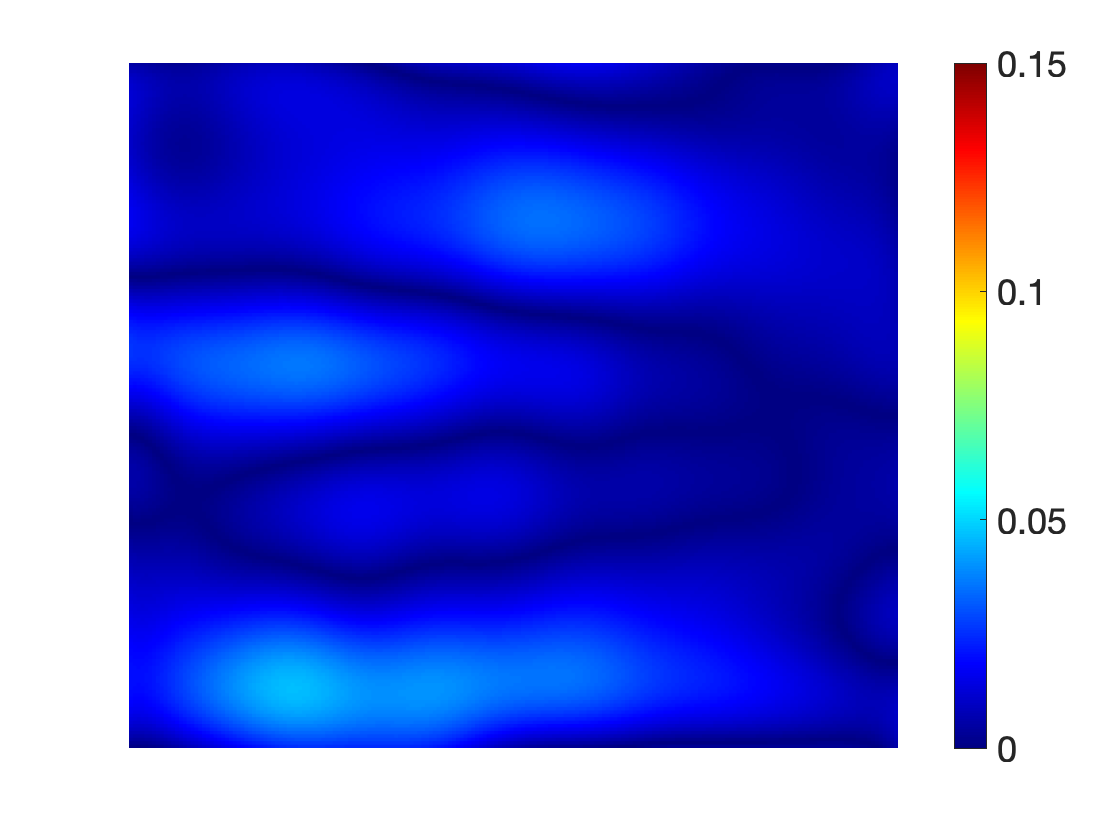} &
\includegraphics[width=0.199\textwidth]{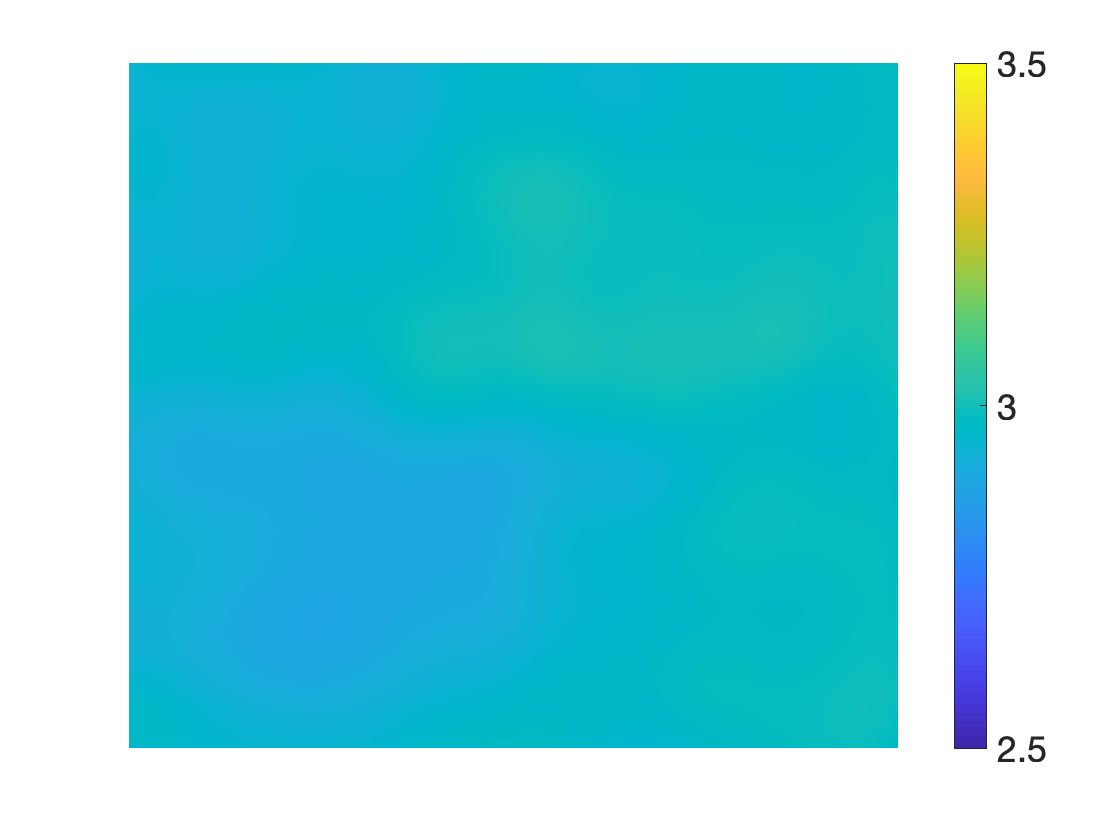} &
\includegraphics[width=0.199\textwidth]{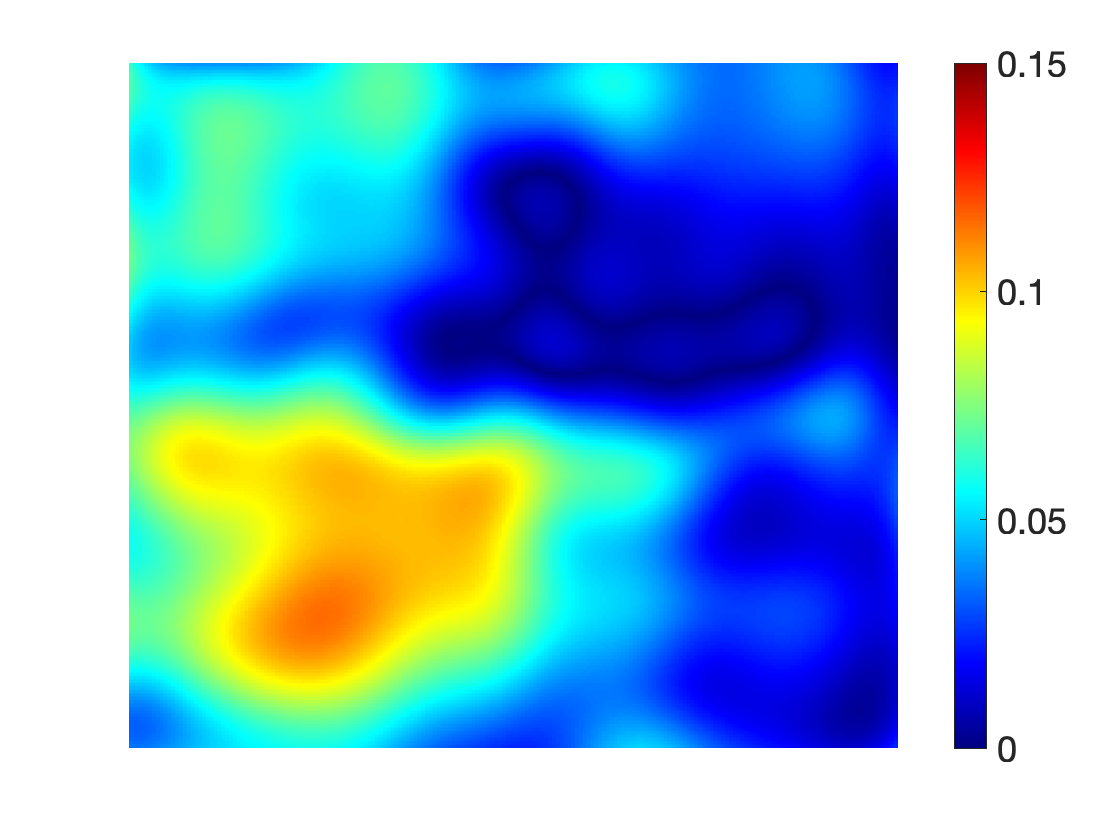}\\
\includegraphics[width=0.199\textwidth]{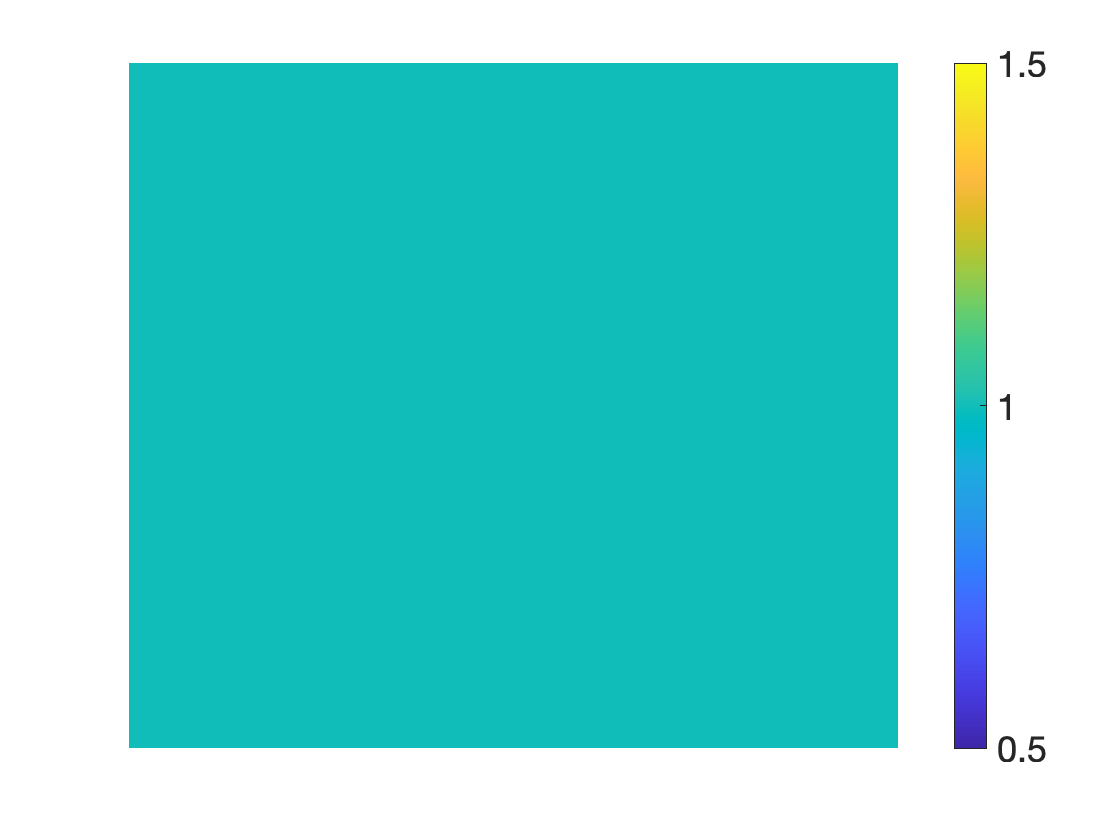} &
\includegraphics[width=0.199\textwidth]{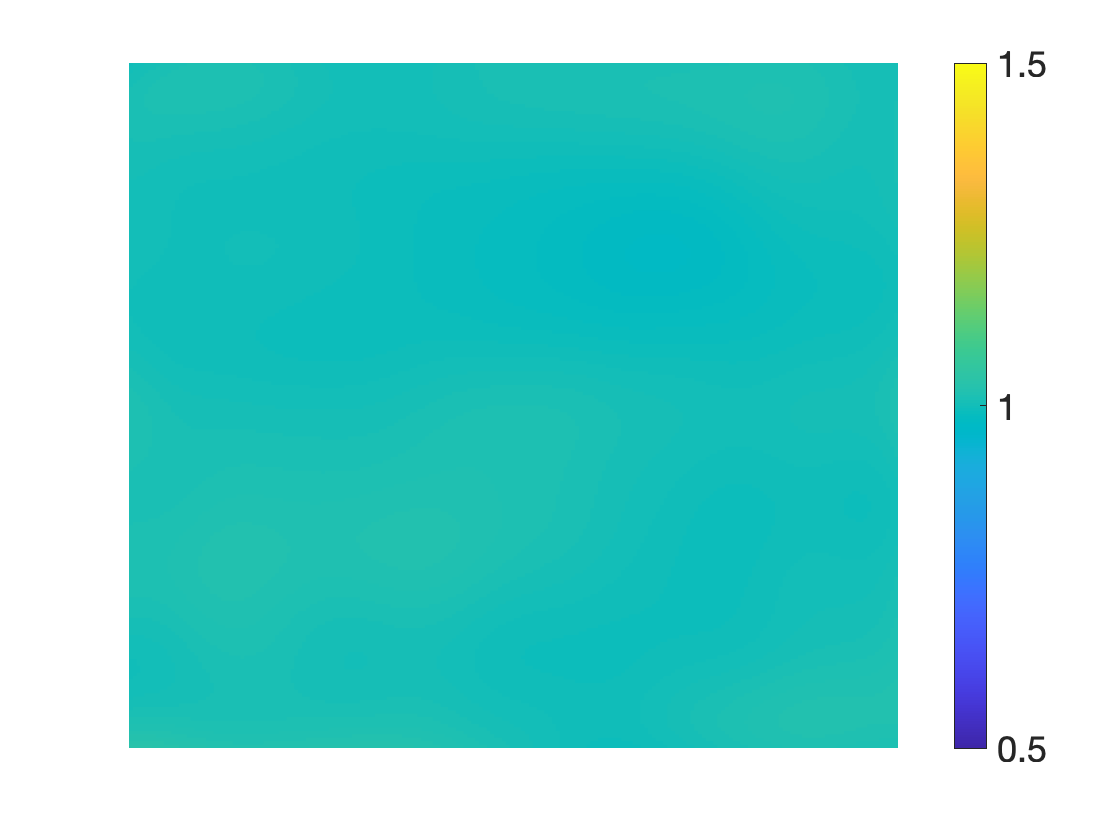} &
\includegraphics[width=0.199\textwidth]{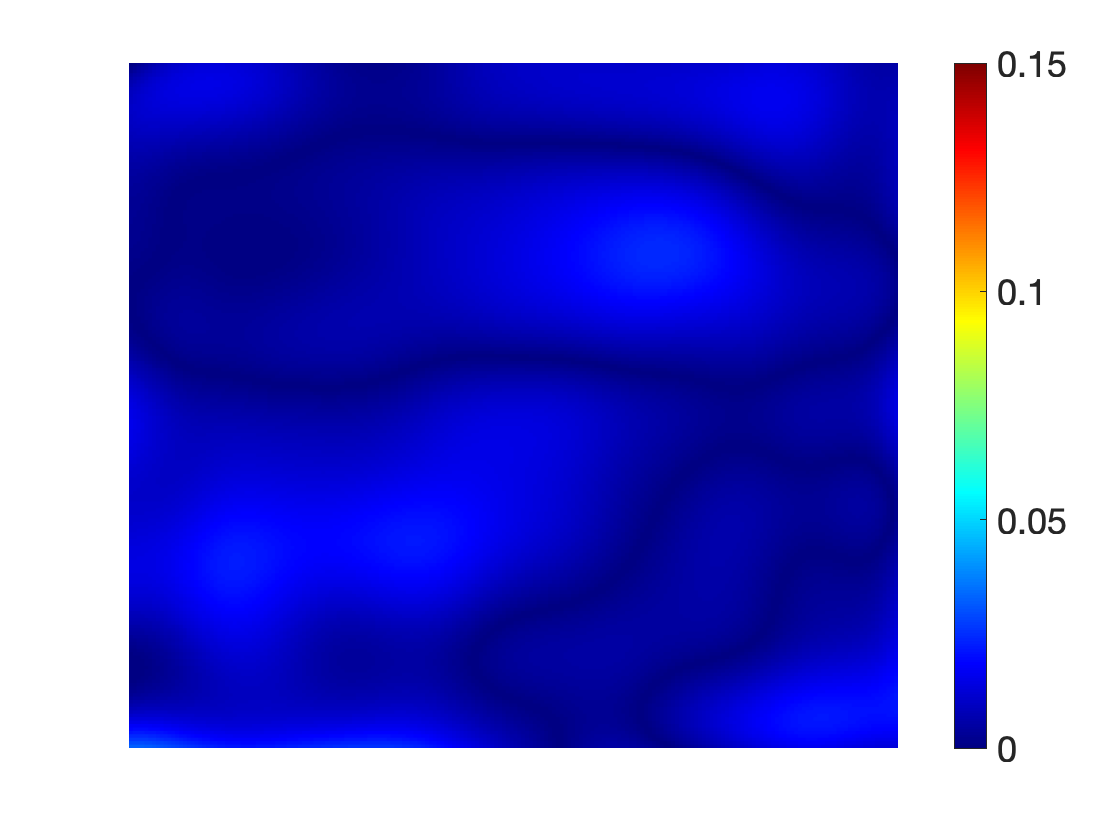}&
\includegraphics[width=0.199\textwidth]{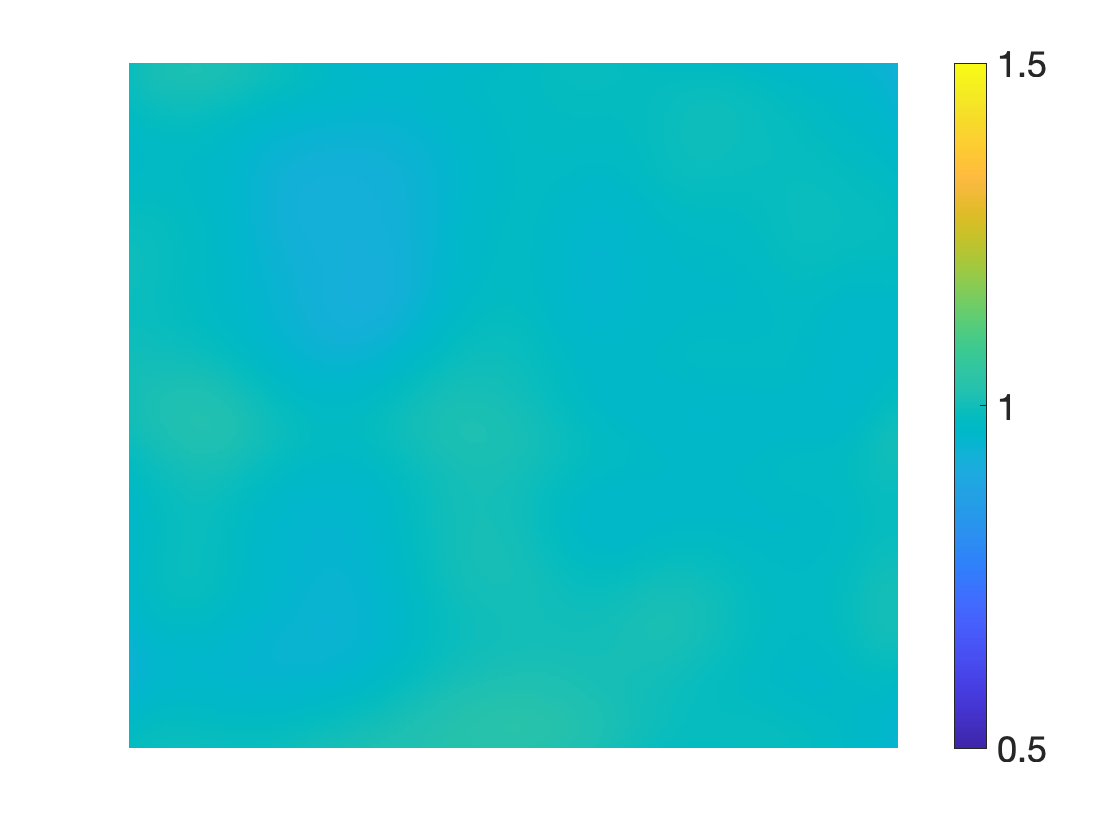} &
\includegraphics[width=0.199\textwidth]{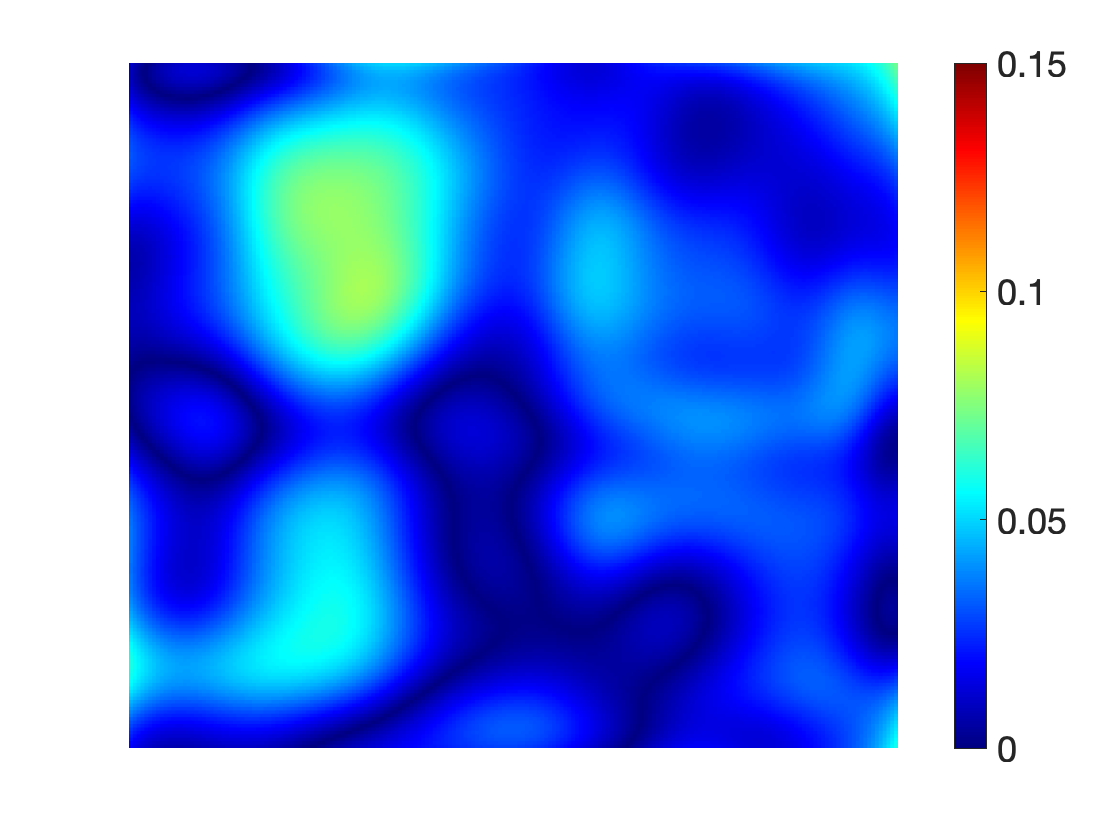}\\
\includegraphics[width=0.199\textwidth]{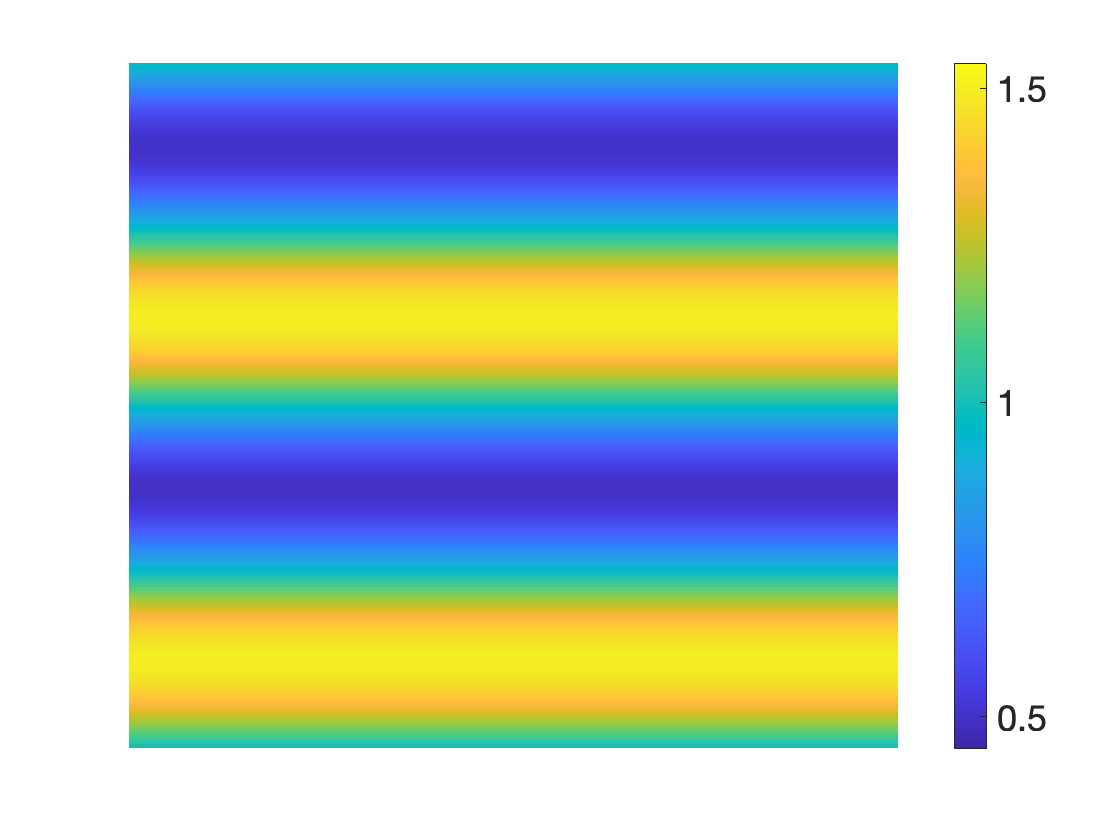} &
\includegraphics[width=0.199\textwidth]{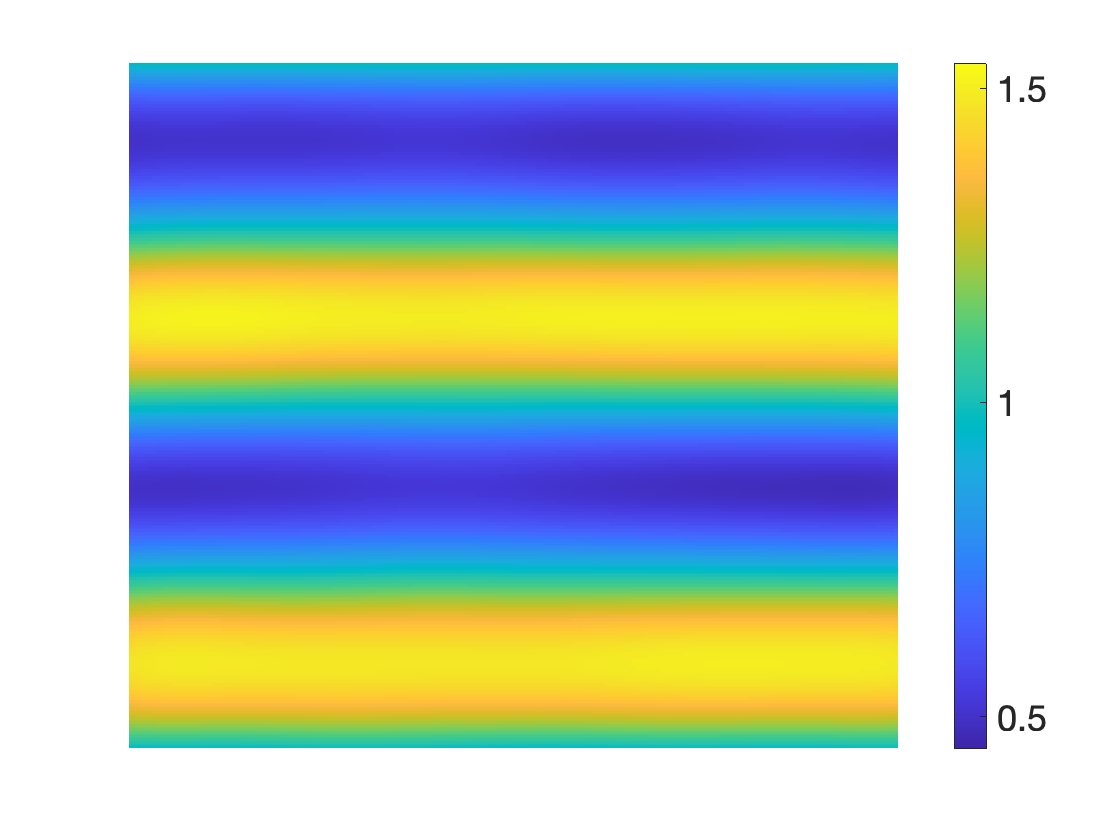} &
\includegraphics[width=0.199\textwidth]{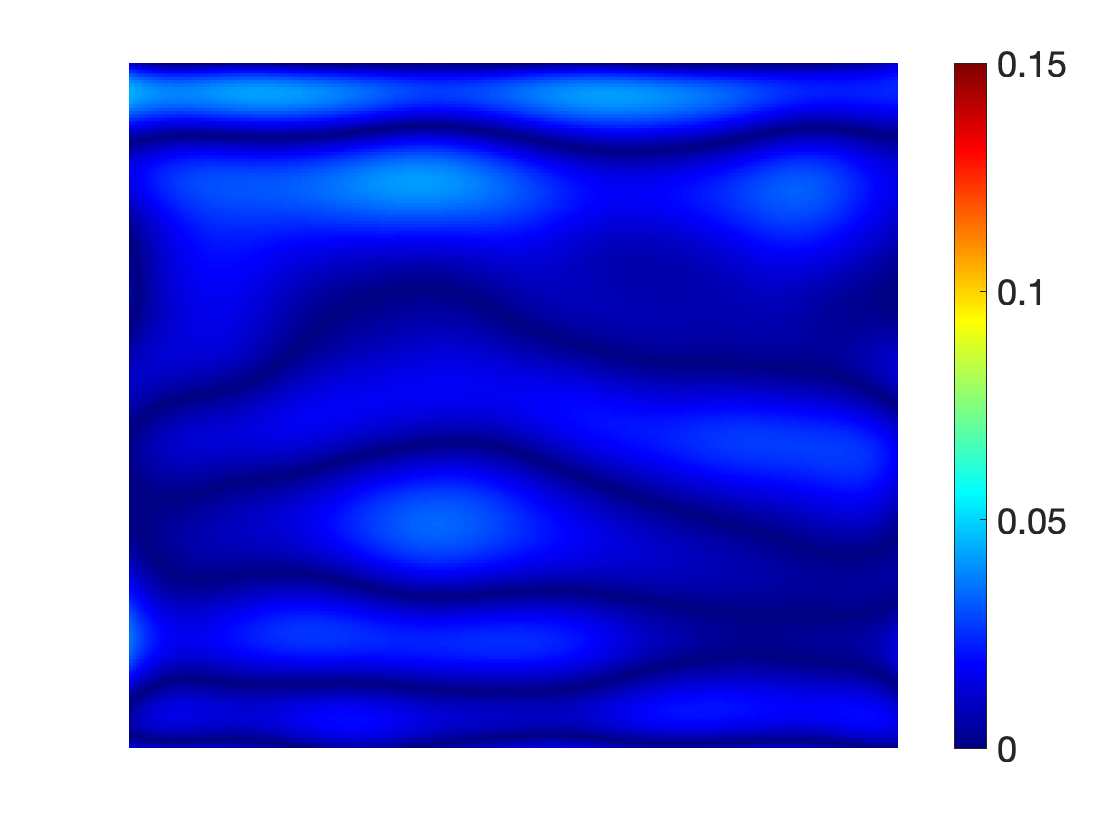} &
\includegraphics[width=0.199\textwidth]{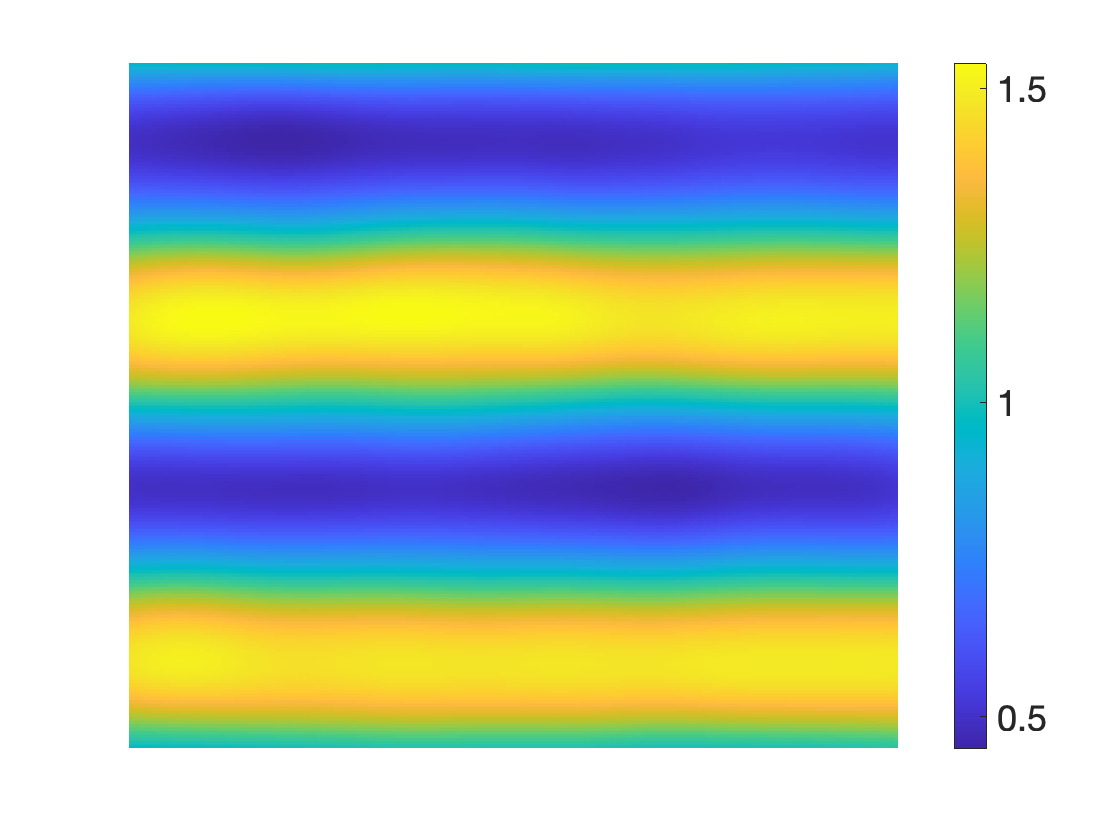} &
\includegraphics[width=0.199\textwidth]{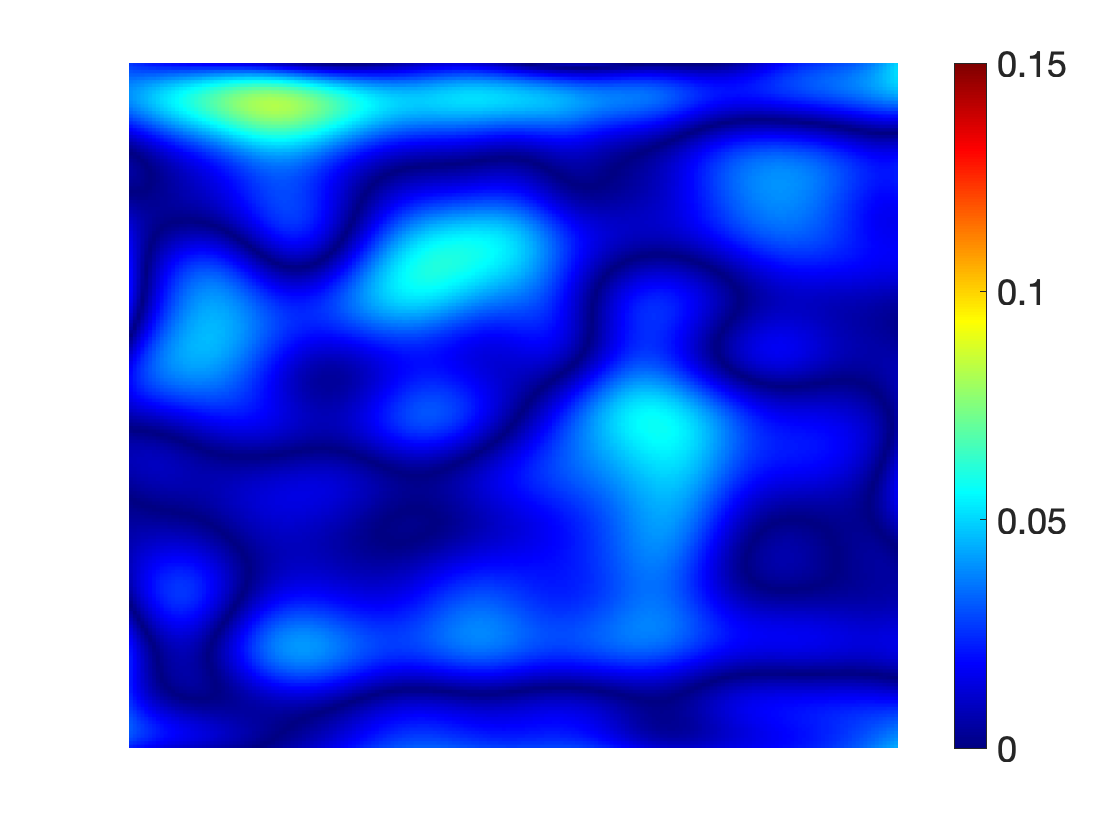}\\
\includegraphics[width=0.199\textwidth]{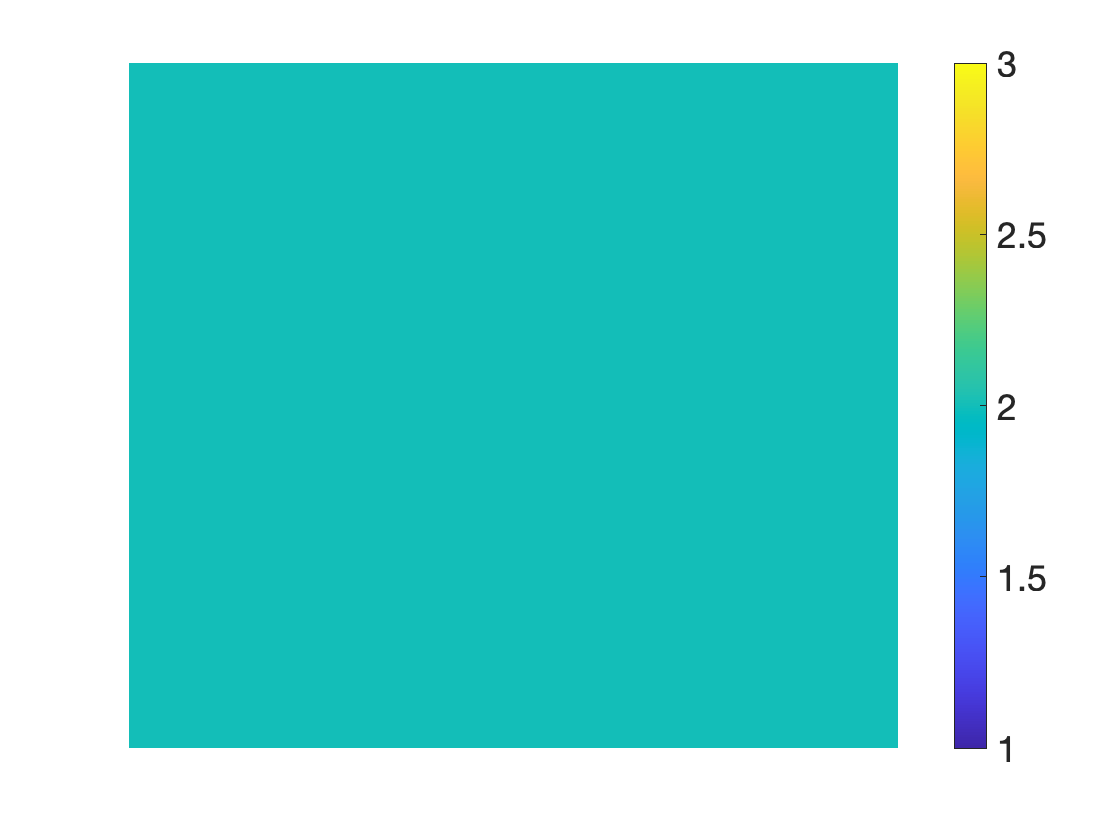} &
\includegraphics[width=0.199\textwidth]{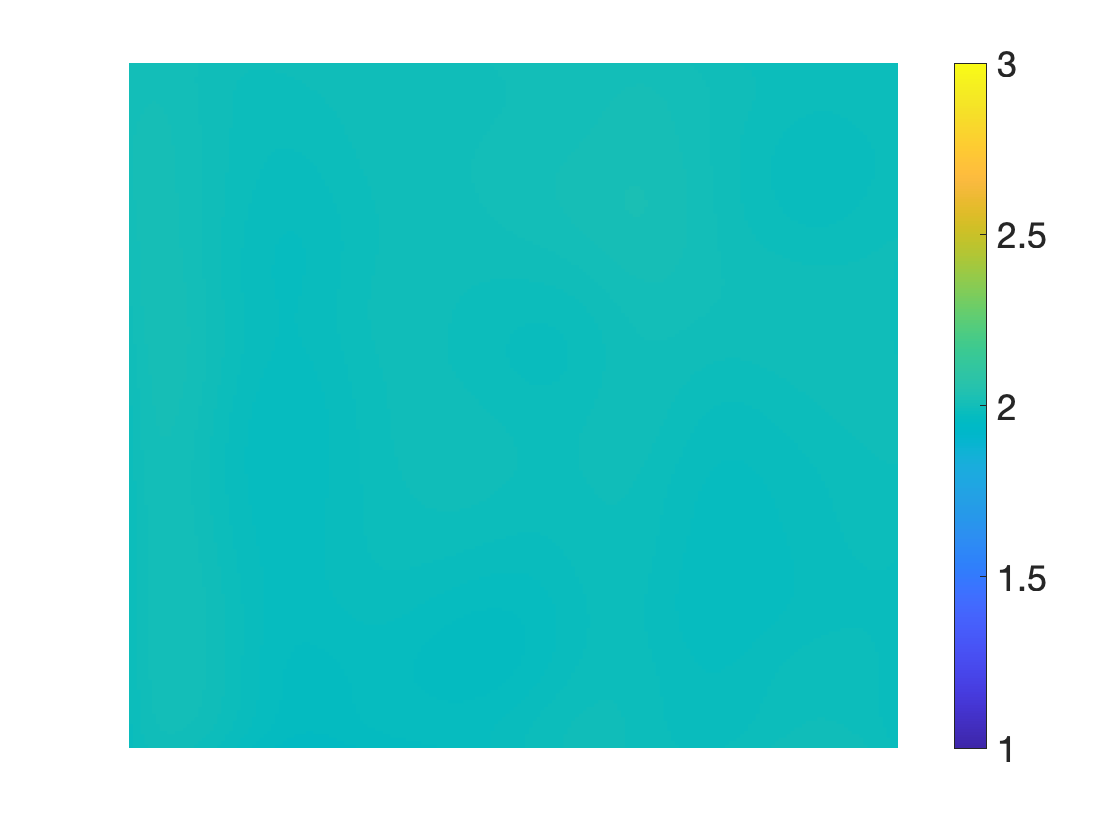} &
\includegraphics[width=0.199\textwidth]{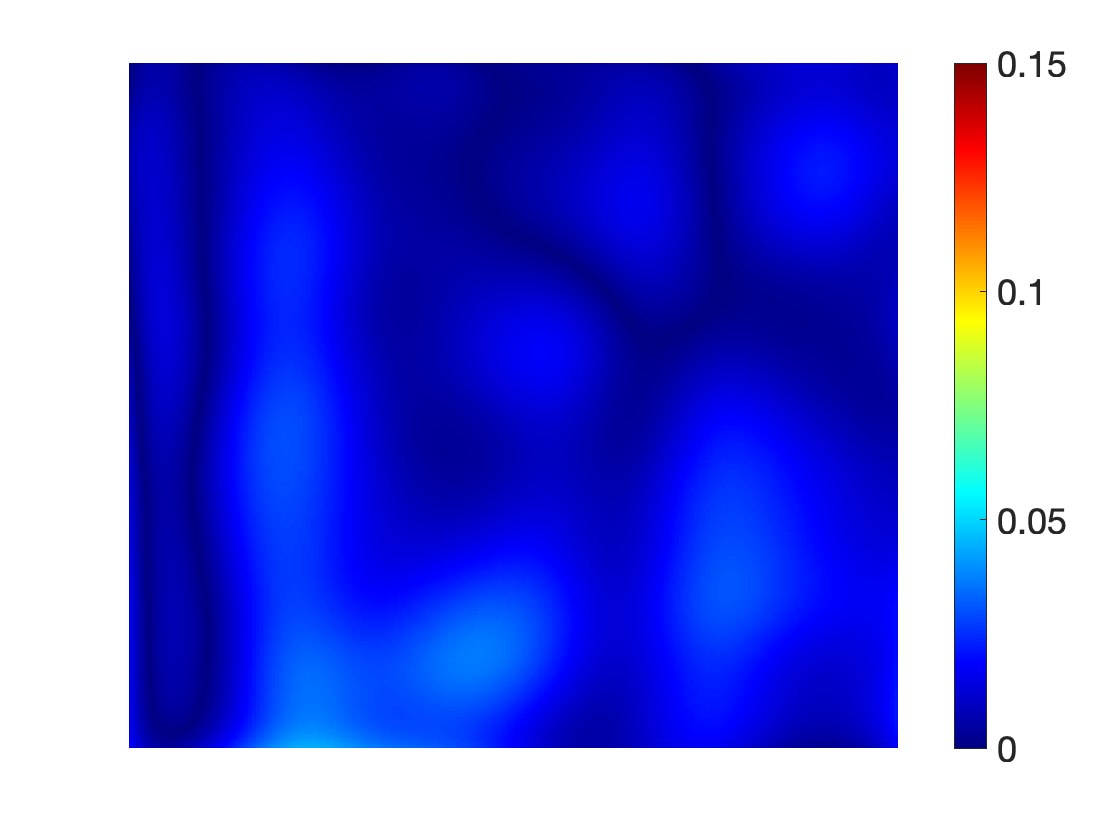} &
\includegraphics[width=0.199\textwidth]{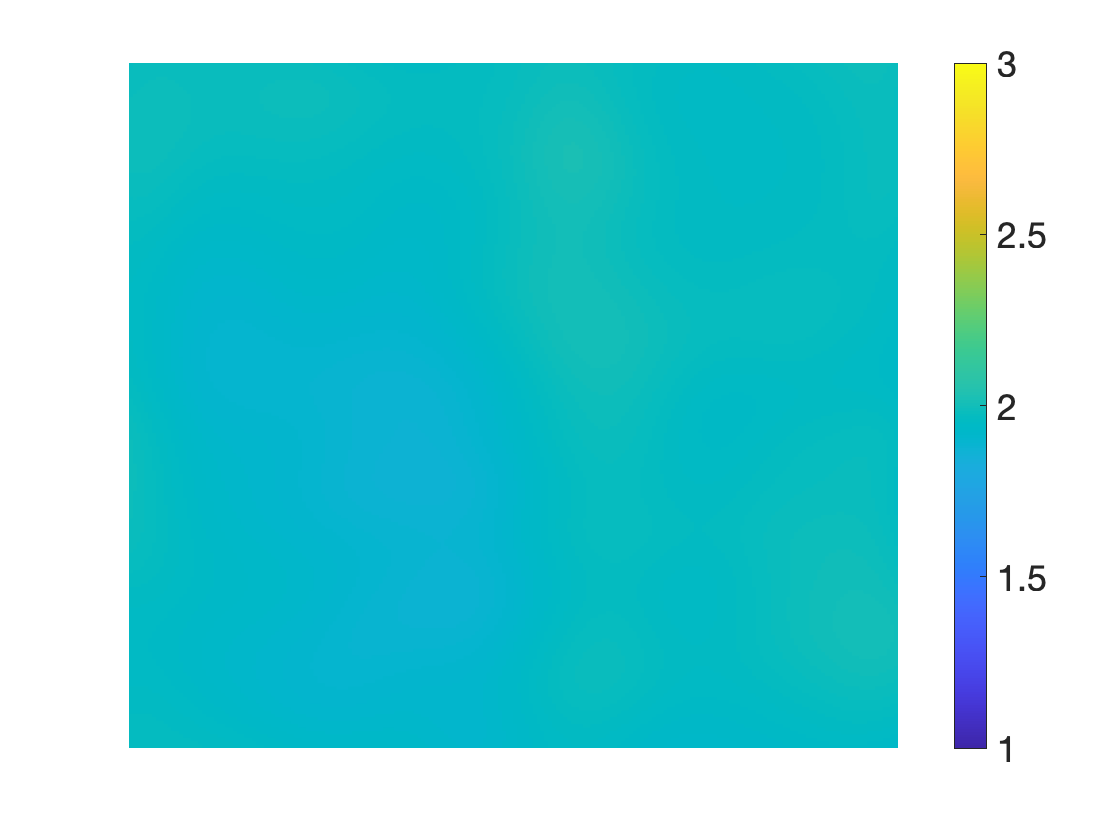} &
\includegraphics[width=0.199\textwidth]{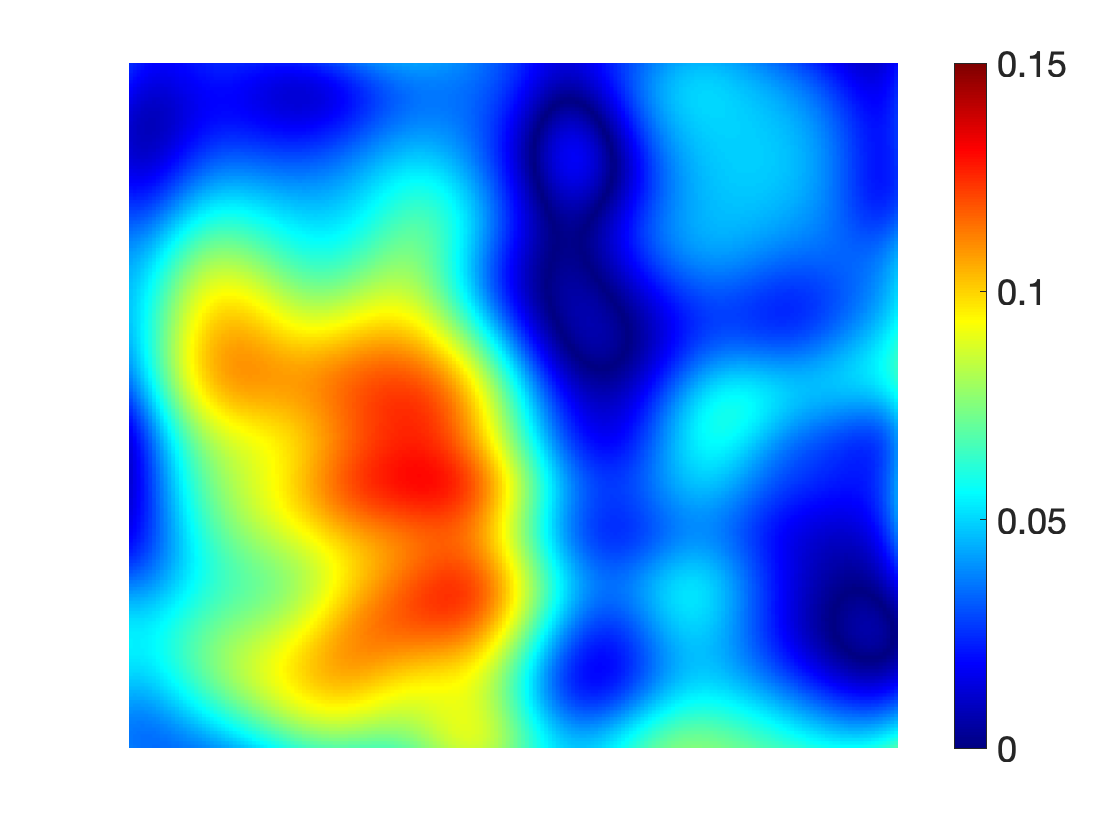}\\
\includegraphics[width=0.199\textwidth]{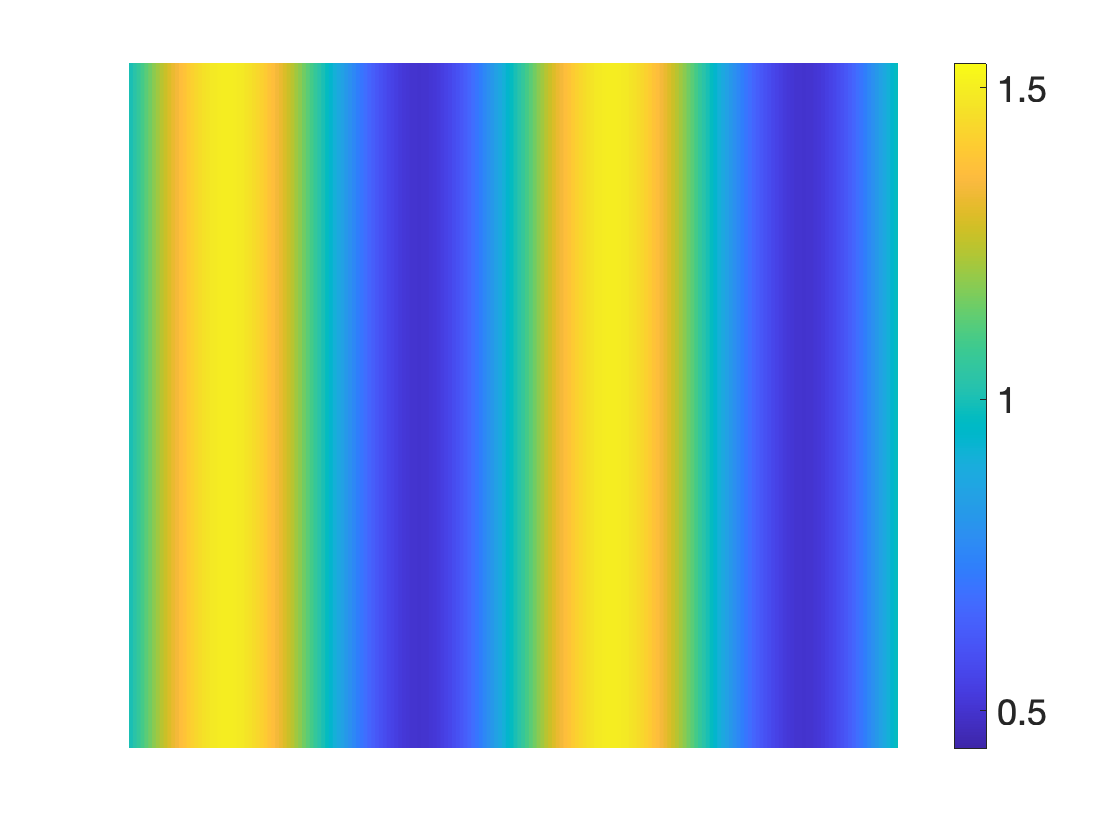} &
\includegraphics[width=0.199\textwidth]{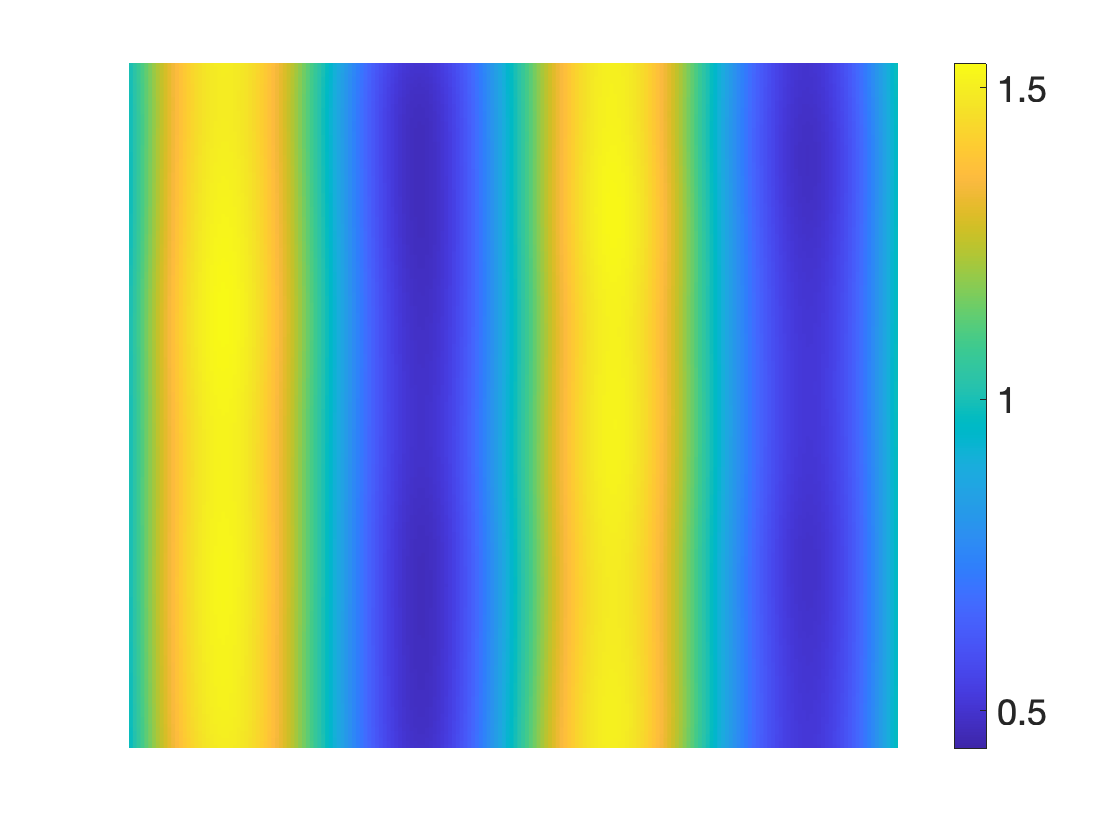} &
\includegraphics[width=0.199\textwidth]{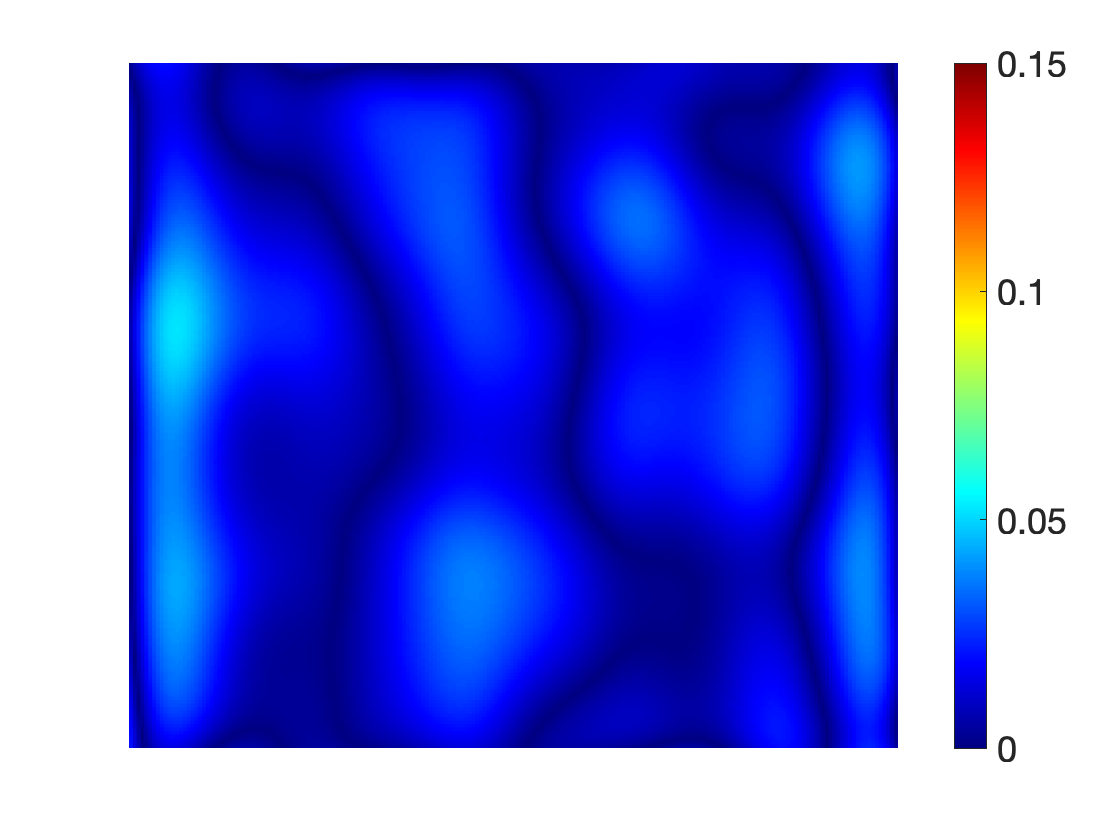} &
\includegraphics[width=0.199\textwidth]{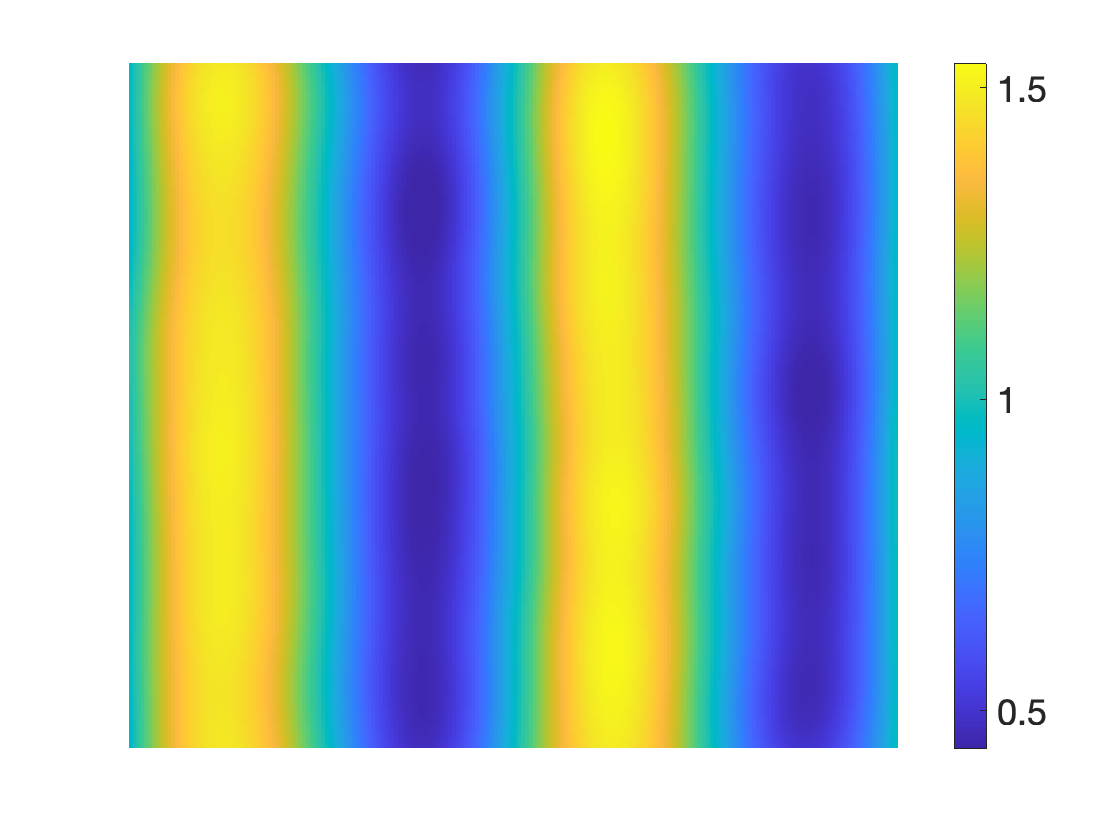} &
\includegraphics[width=0.199\textwidth]{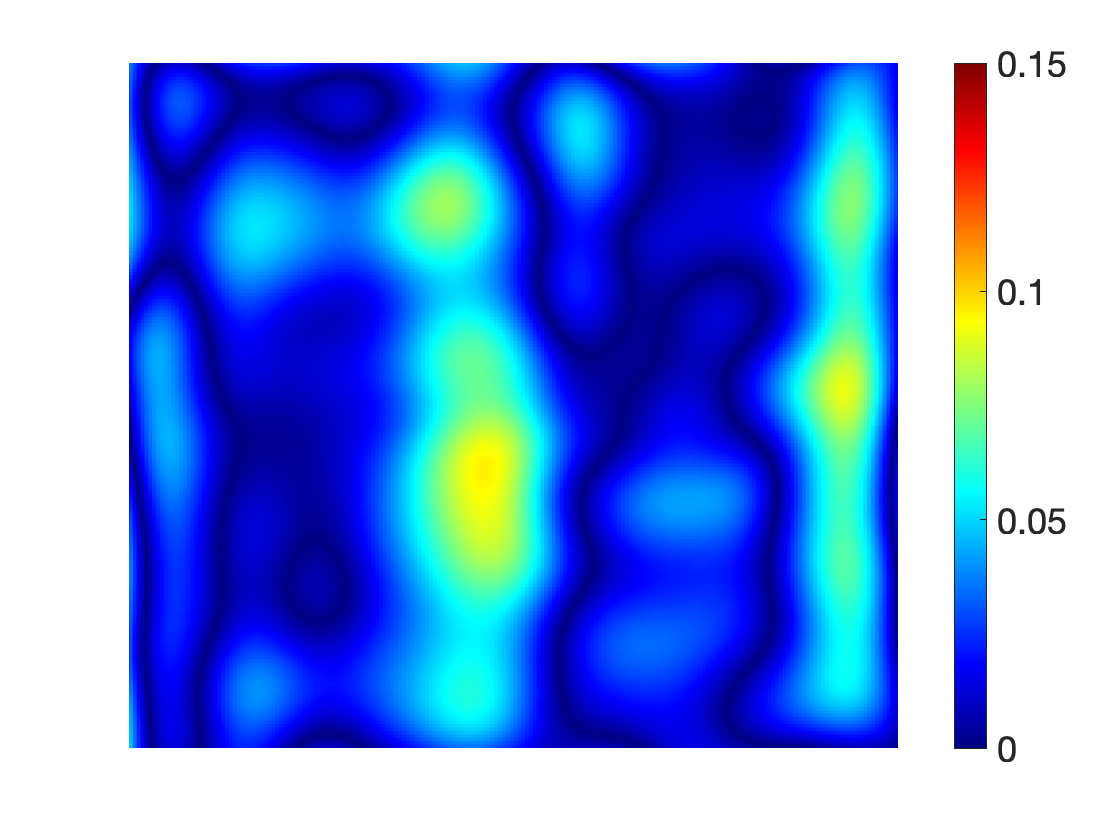}\\
\includegraphics[width=0.199\textwidth]{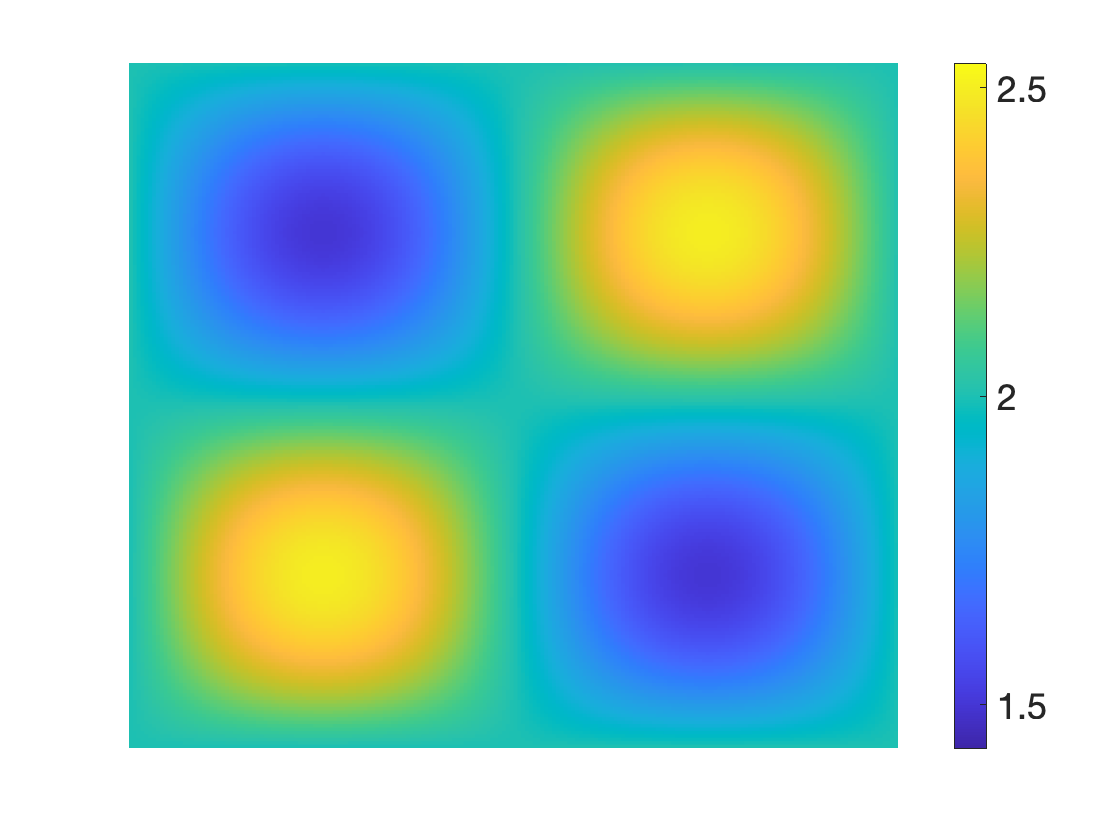} &
\includegraphics[width=0.199\textwidth]{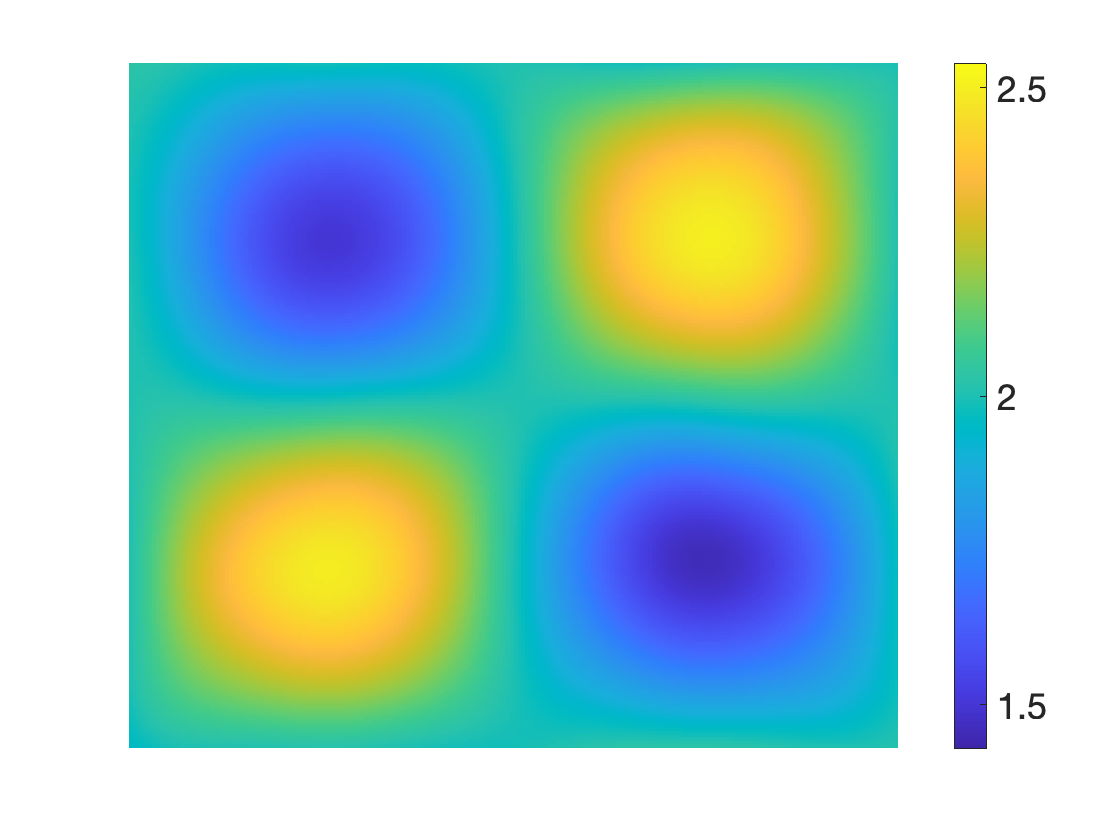} &
\includegraphics[width=0.199\textwidth]{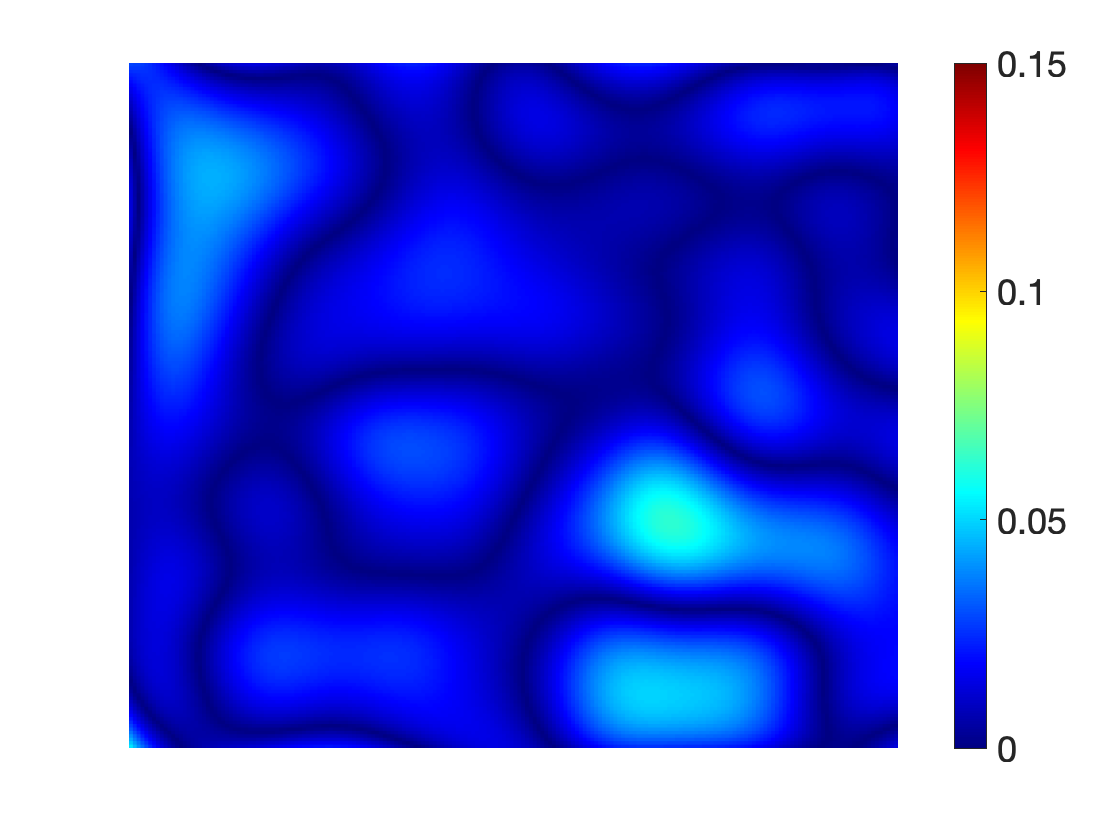} &
\includegraphics[width=0.199\textwidth]{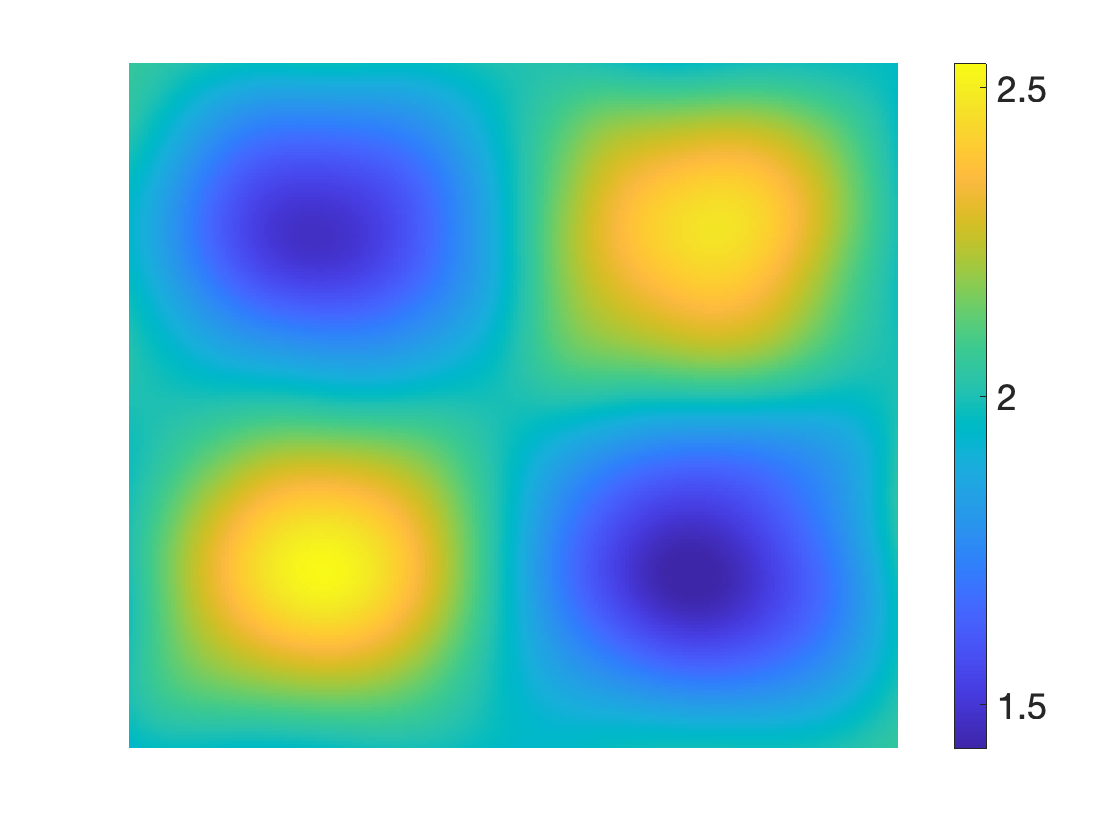} &
\includegraphics[width=0.199\textwidth]{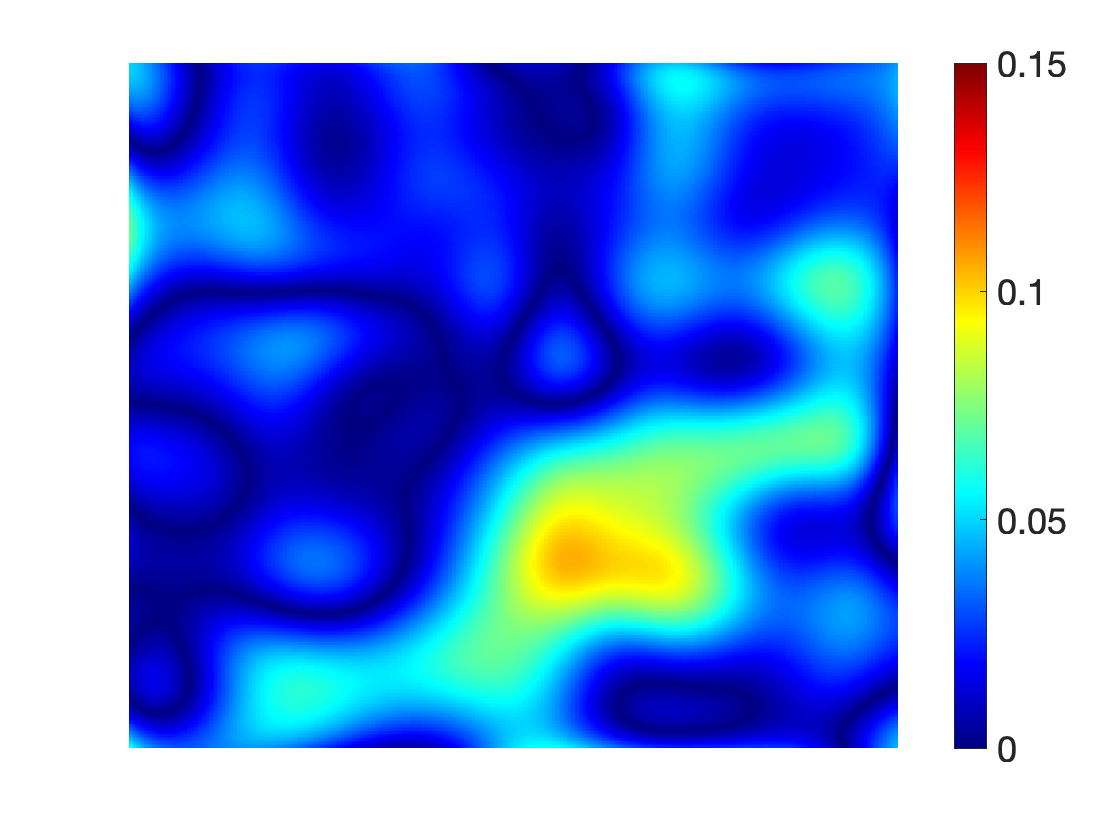}\\
(a) $A^\dag$  & (b) $\hat A$ & (c) $|\hat A-A^\dag|$ & (d) $\hat A$ & (e) $|\hat A-A^\dag|$
\end{tabular}
\caption{The reconstructions for Example \ref{exam:diri3d1} with exact data in (b) and noisy data $(\delta=5\%)$ in (d). From the top to bottom, the results are for $A_{11}$, $A_{12}$, $A_{13}$, $A_{22}$, $A_{23}$ and $A_{33}$, respectively.}
\label{fig:diri3d1}
\end{figure}

The reconstruction results on a 2D cross section at $x_3=0.5$ for both exact and noisy data ($\delta=5\%$) are shown in Fig. \ref{fig:diri3d1}. The overall features of the anisotropic conductivity tensor are well recovered at both noise levels, demonstrating again the robustness of the DNN approach to data noise and its effectiveness in the three-dimensional case.

The last example is about recovering a 3D conductivity matrix from partial internal data.
\begin{example}\label{exam:diri3d2}
The domain $\Omega=(0,1)^3$, the measurement $\nabla z_i^\delta$ on the region $\omega=\Omega\setminus (0.2,0.8)^3$, $A^\dagger = \begin{pmatrix}
    2+x_2^2& 1& 1+\frac{\sin(4\pi x_2)}{2}\\
    1& 2+x_1^2& 1+\frac{\sin(4\pi x_1)}{2}\\
    1+\frac{\sin(4\pi x_2)}{2}& 1+\frac{\sin(4\pi x_1)}{2}& 2+\frac{\sin(\pi x_1)\sin(\pi x_2)}{2}\\
    \end{pmatrix}$, $u_1^\dag=x_1+x_2+x_3+\frac{1}{3}(x_1^3+x_2^3+x_3^3)$, $u_2^\dag=x_1-x_2+x_3+\frac{1}{3}(x_1^3-x_2^3+x_3^3)$, $u_3^\dag=x_1+x_2-x_3+\frac{1}{3}(x_1^3+x_2^3-x_3^3)$, $u_4^\dag=-x_1+x_2+x_3+\frac{1}{3}(-x_1^3+x_2^3+x_3^3)$, $u_5^\dagger=-u_3^\dagger$ and $u_6^\dagger=-u_2^\dagger$.
\end{example}

\begin{figure}[htb!]
\centering
\setlength{\tabcolsep}{0em}
\begin{tabular}{ccccc}
\includegraphics[width=0.199\textwidth]{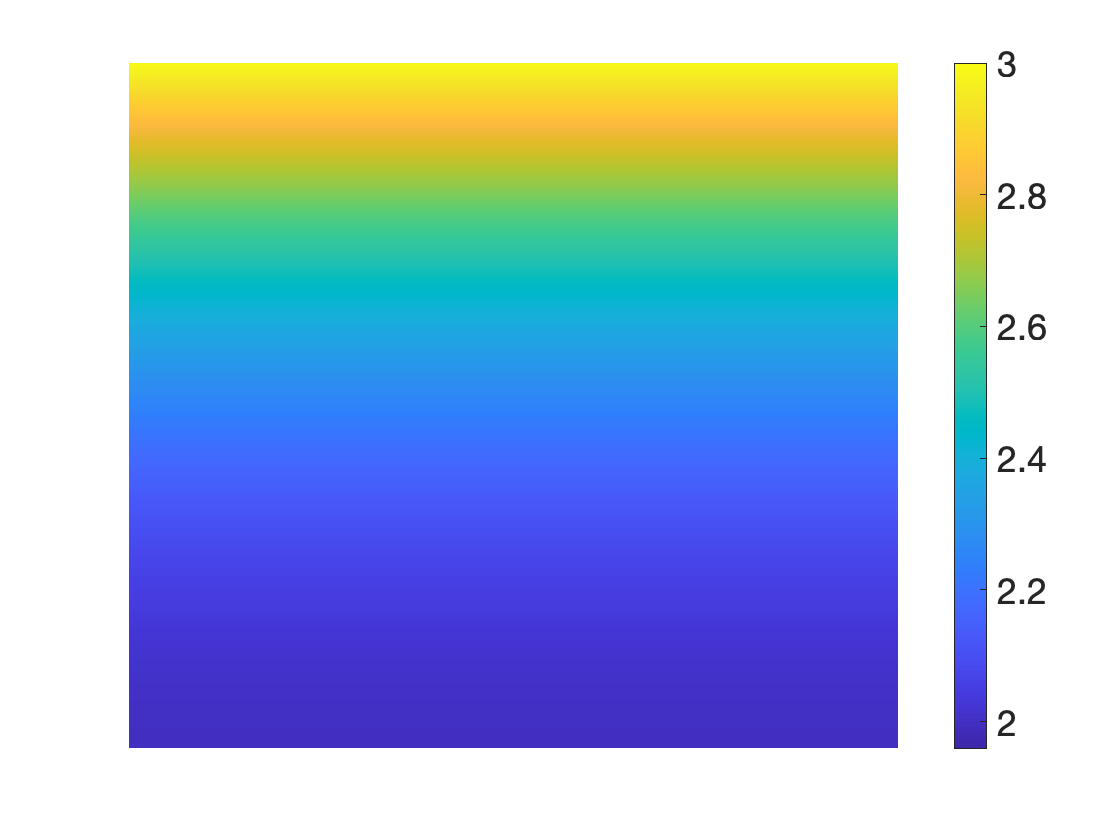} &
\includegraphics[width=0.199\textwidth]{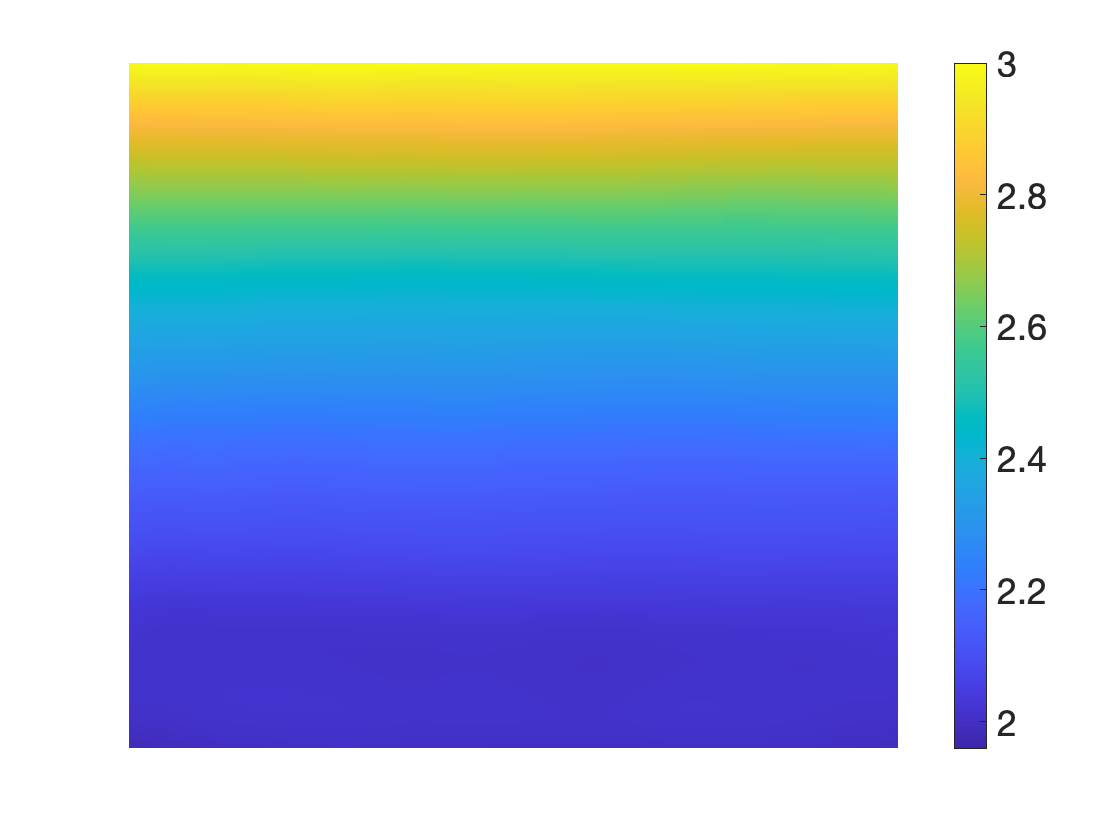} &
\includegraphics[width=0.199\textwidth]{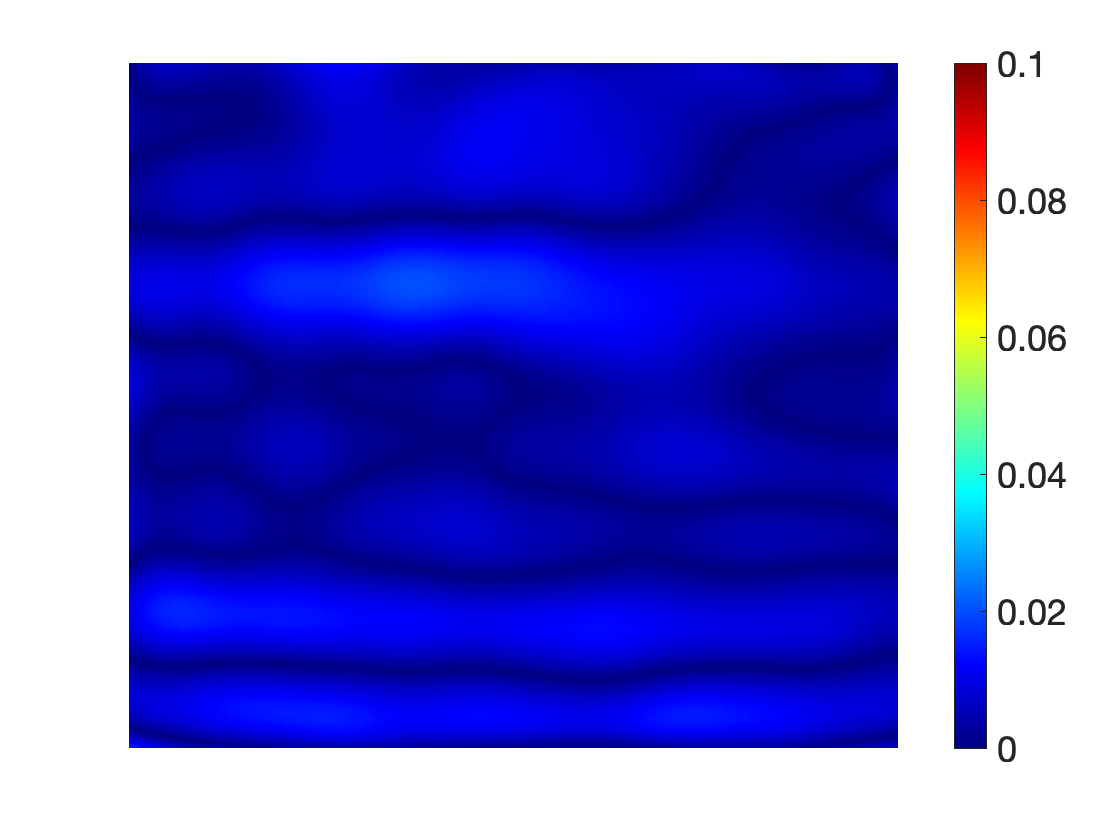} &
\includegraphics[width=0.199\textwidth]{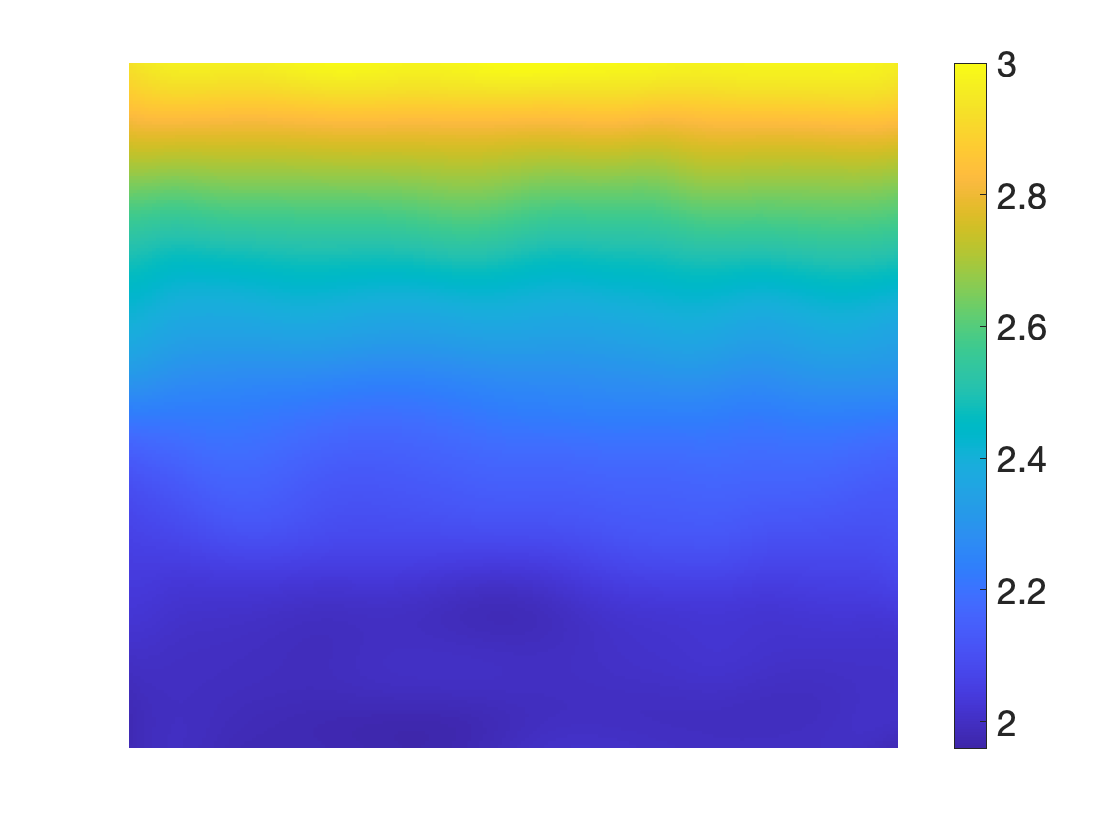} &
\includegraphics[width=0.199\textwidth]{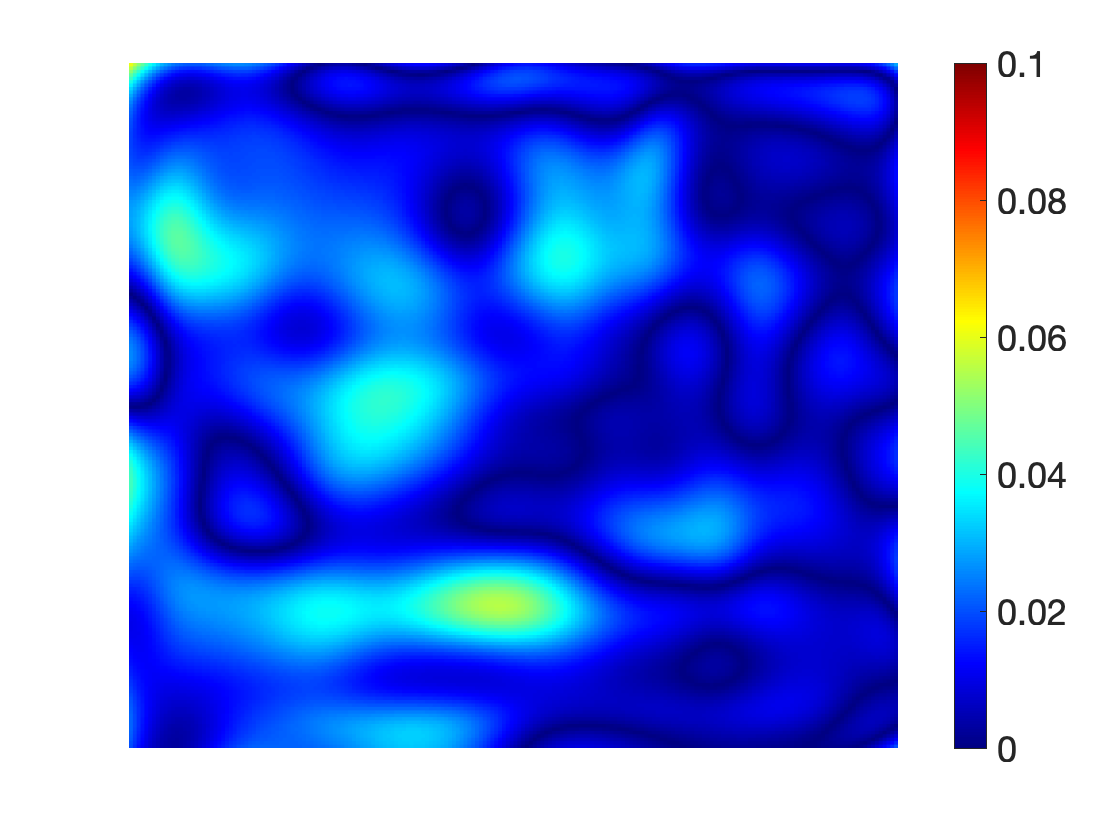}\\
\includegraphics[width=0.199\textwidth]{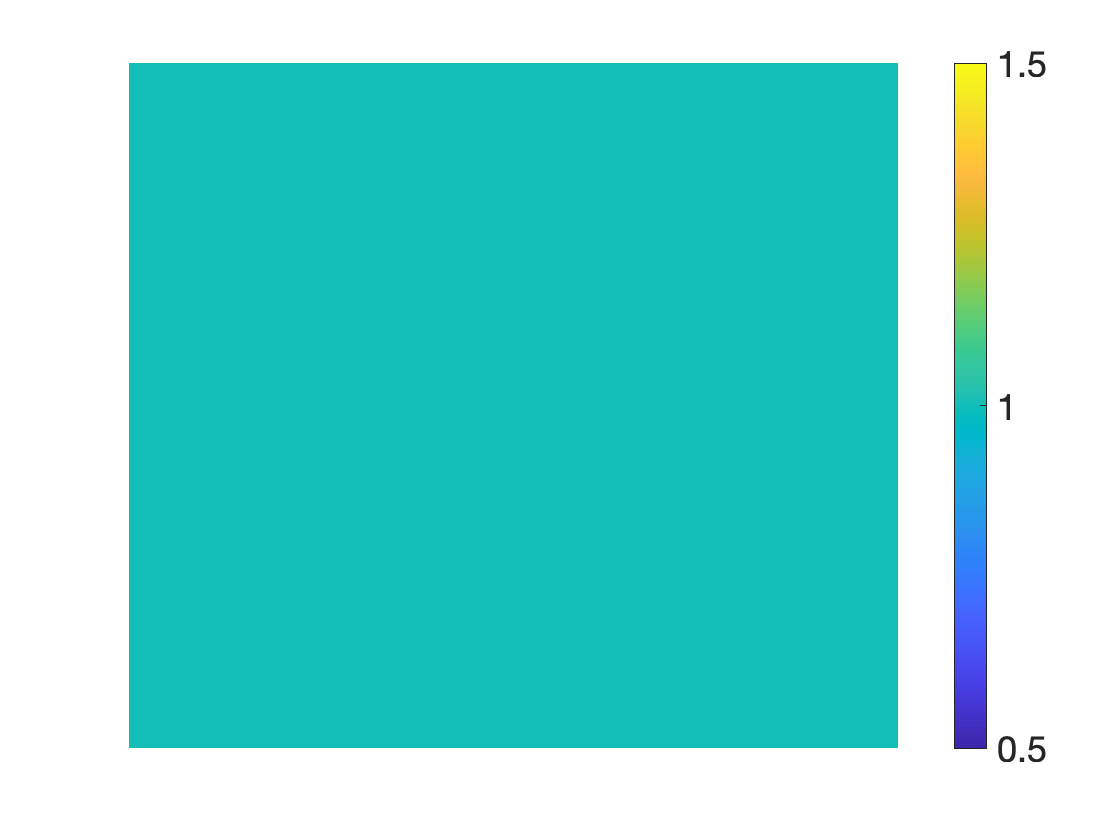} &
\includegraphics[width=0.199\textwidth]{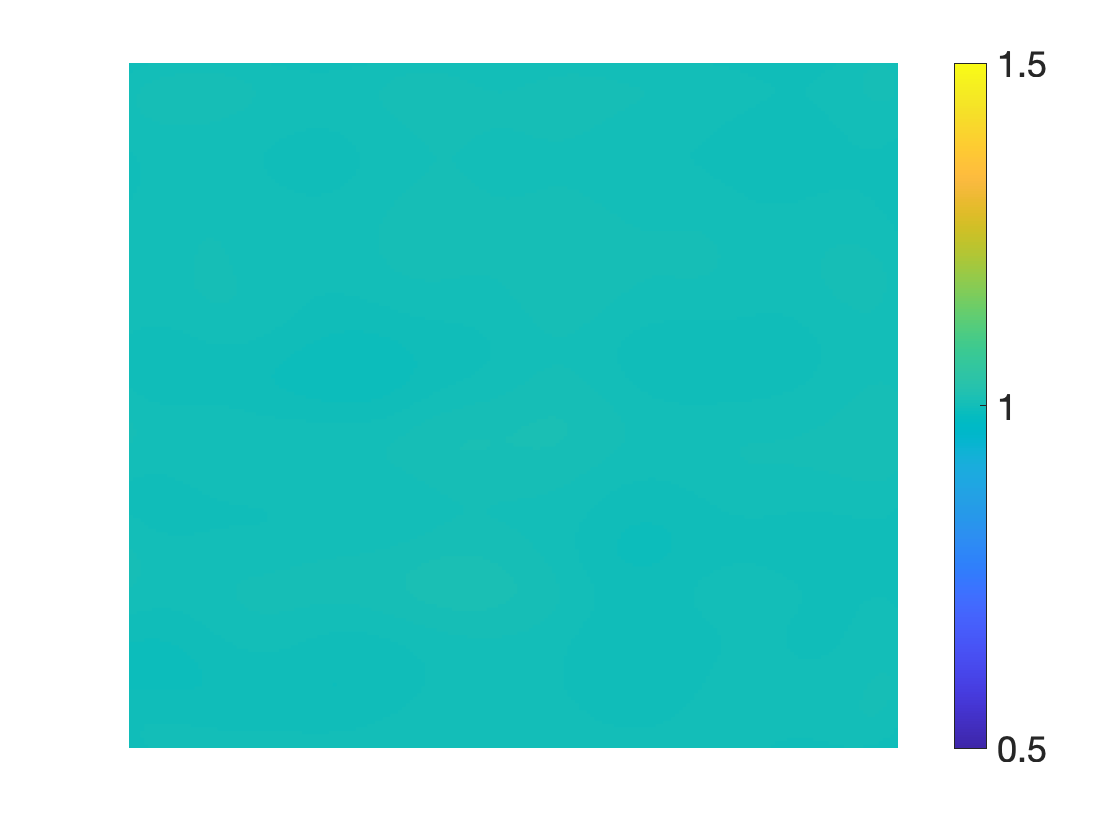} &
\includegraphics[width=0.199\textwidth]{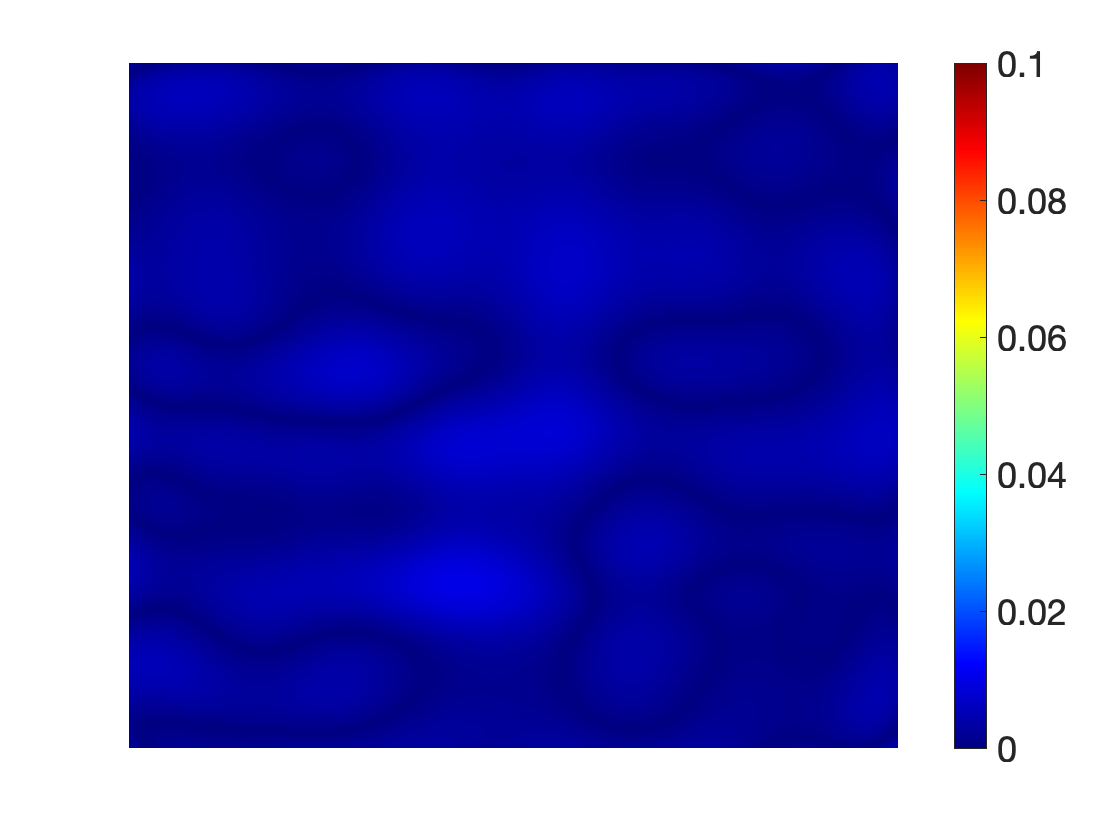}&
\includegraphics[width=0.199\textwidth]{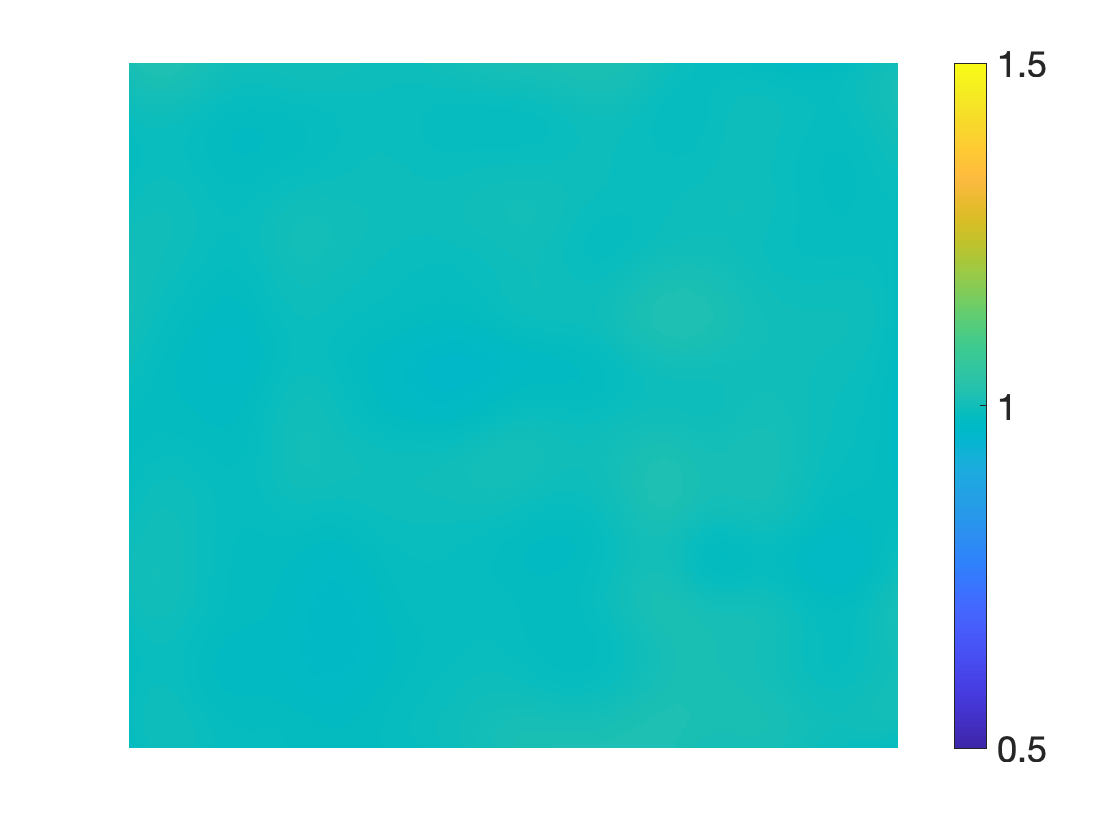} &
\includegraphics[width=0.199\textwidth]{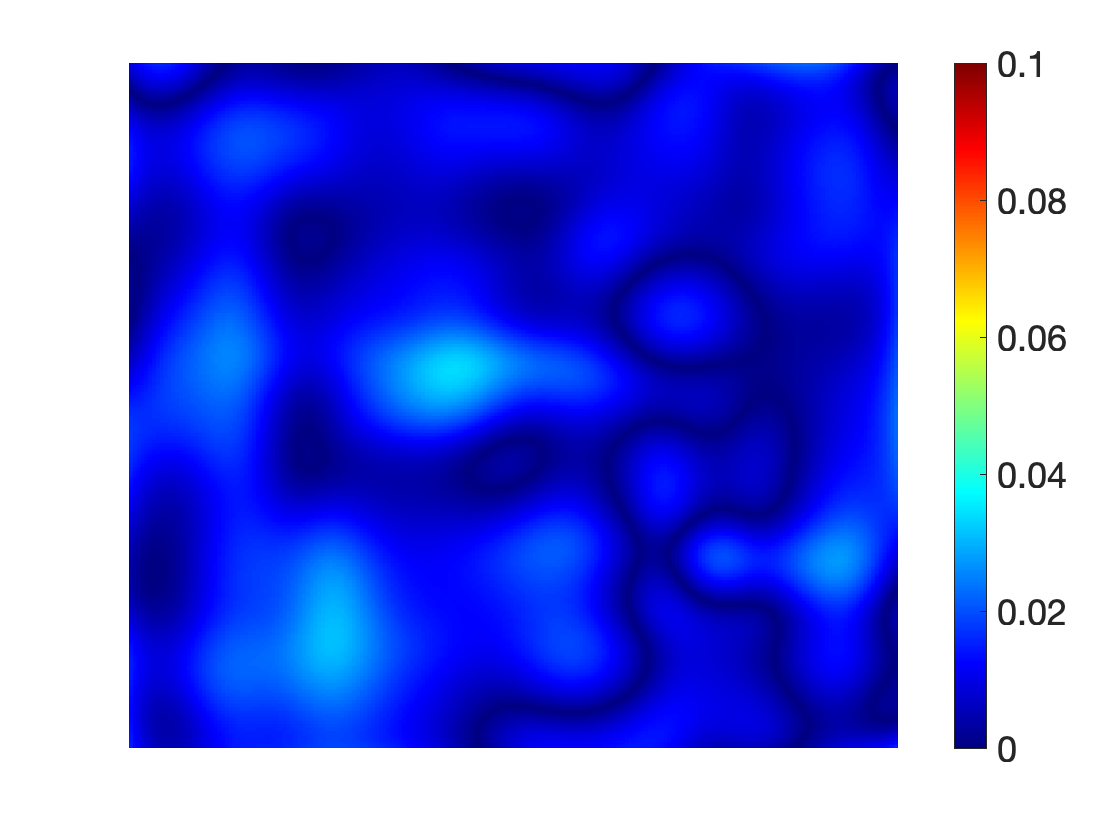}\\
\includegraphics[width=0.199\textwidth]{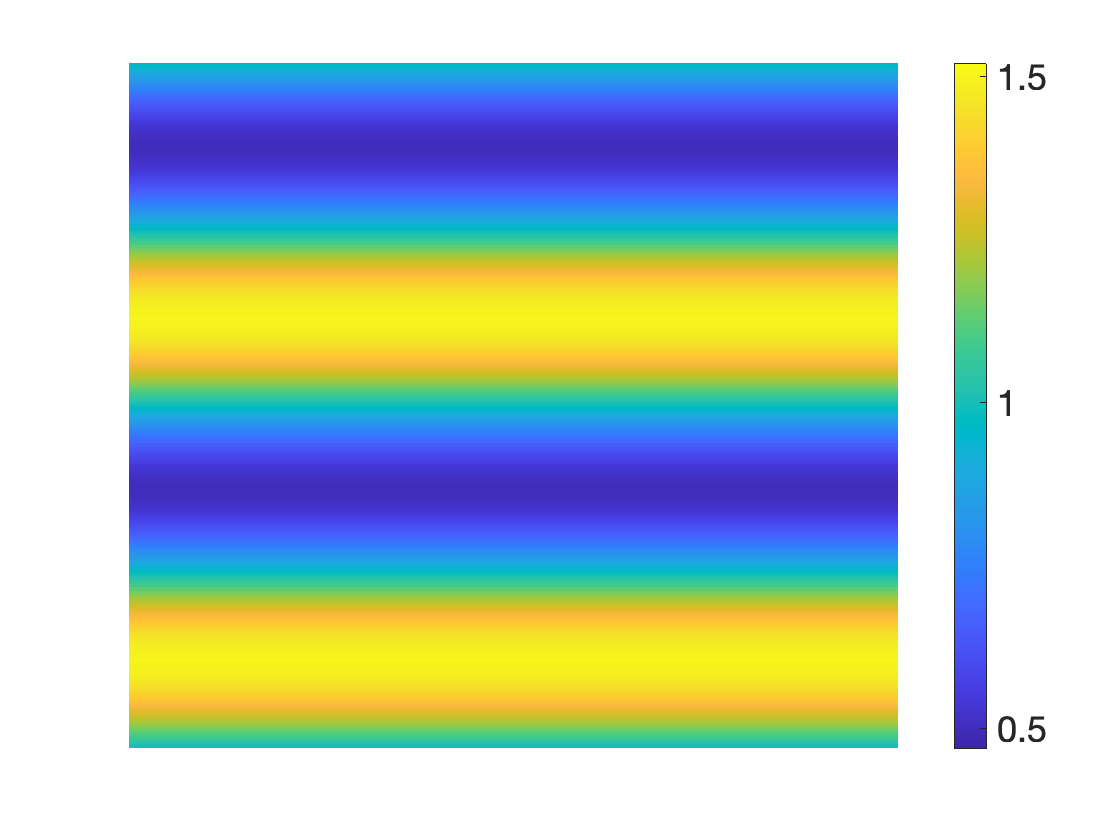} &
\includegraphics[width=0.199\textwidth]{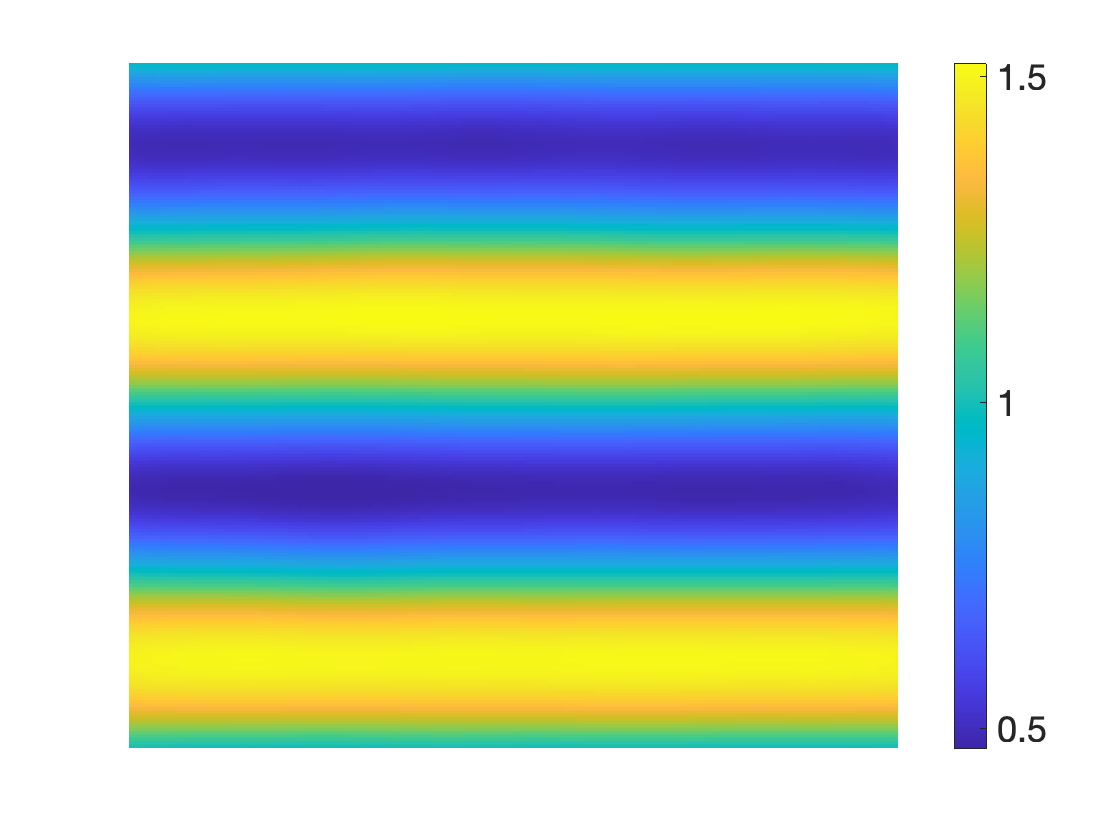} &
\includegraphics[width=0.199\textwidth]{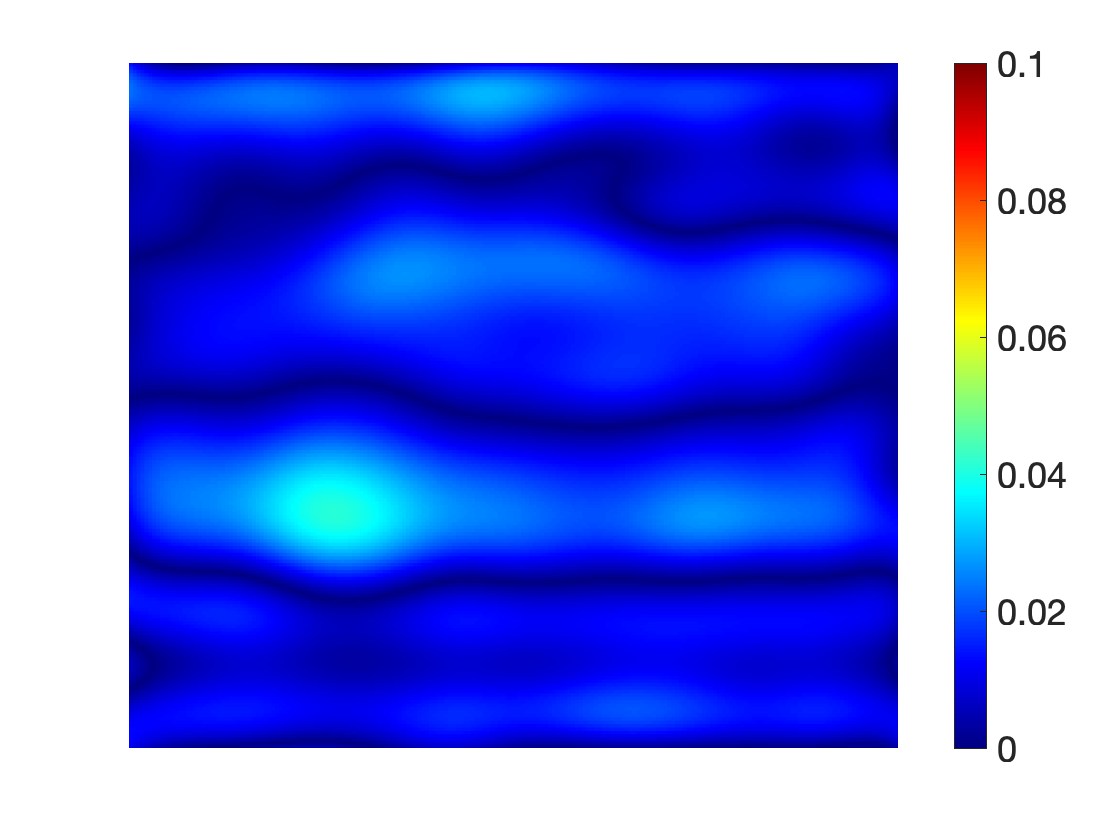} &
\includegraphics[width=0.199\textwidth]{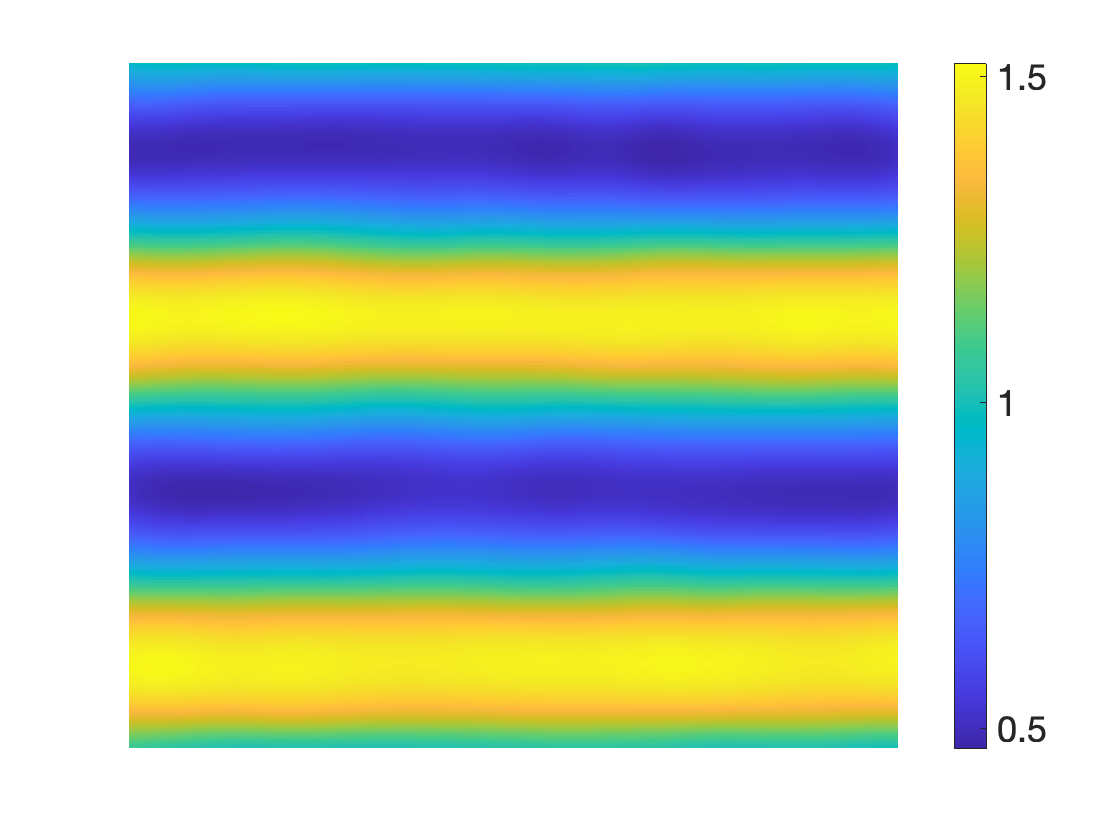} &
\includegraphics[width=0.199\textwidth]{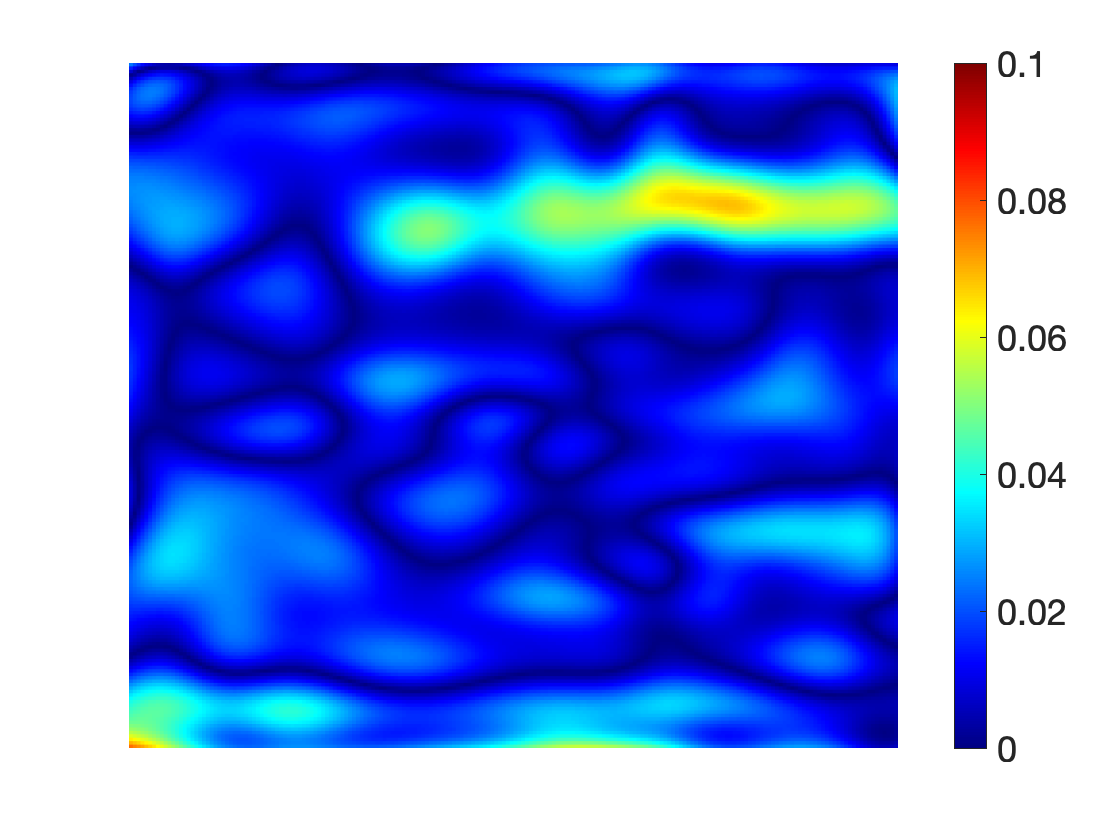}\\
\includegraphics[width=0.199\textwidth]{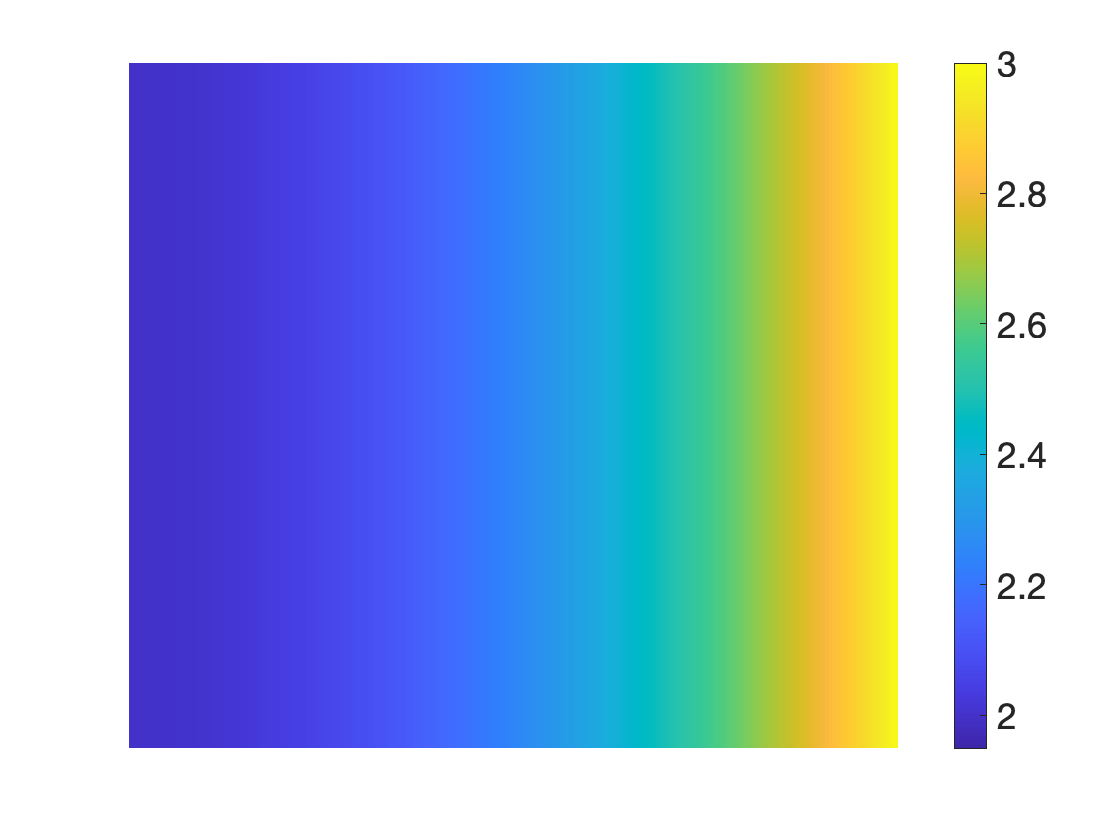} &
\includegraphics[width=0.199\textwidth]{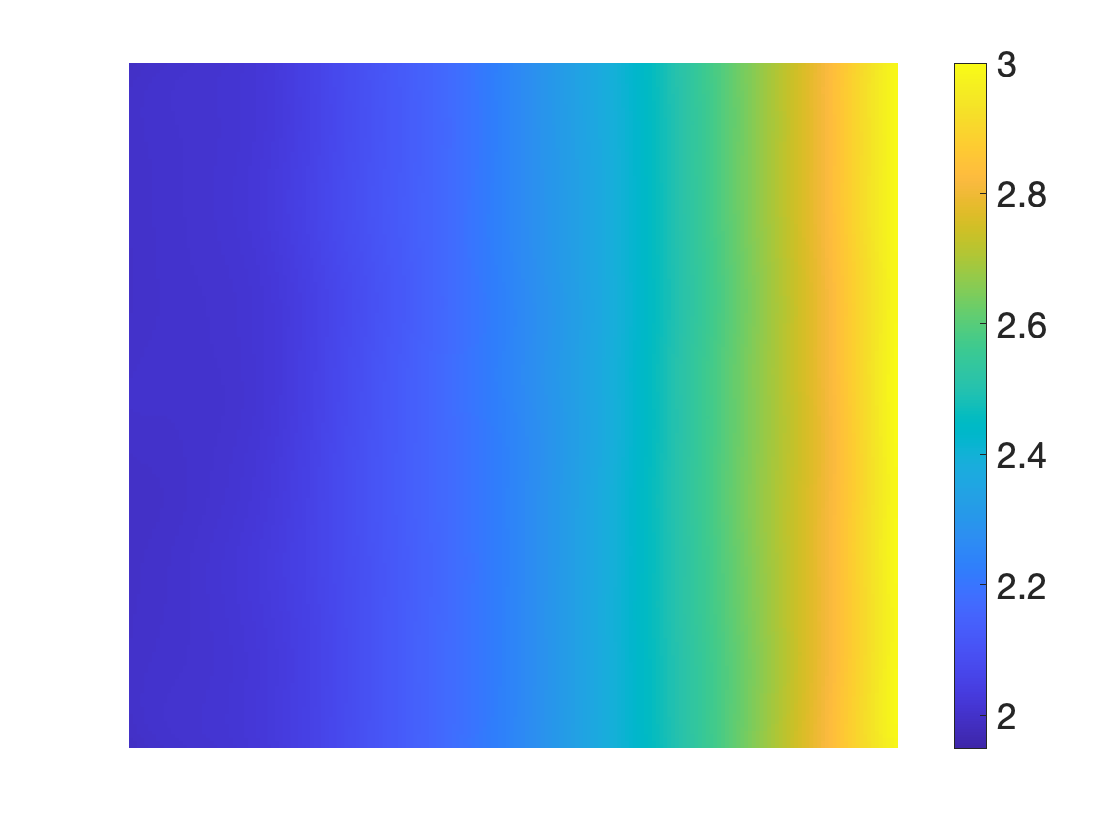} &
\includegraphics[width=0.199\textwidth]{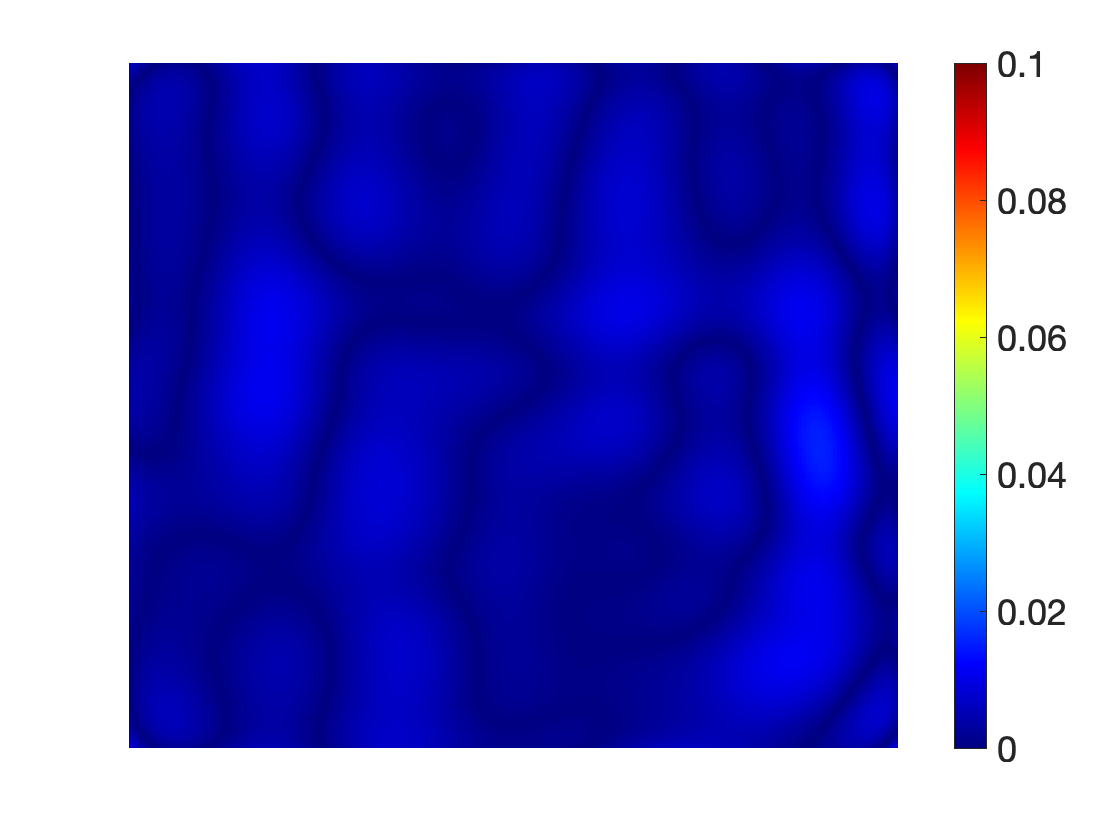} &
\includegraphics[width=0.199\textwidth]{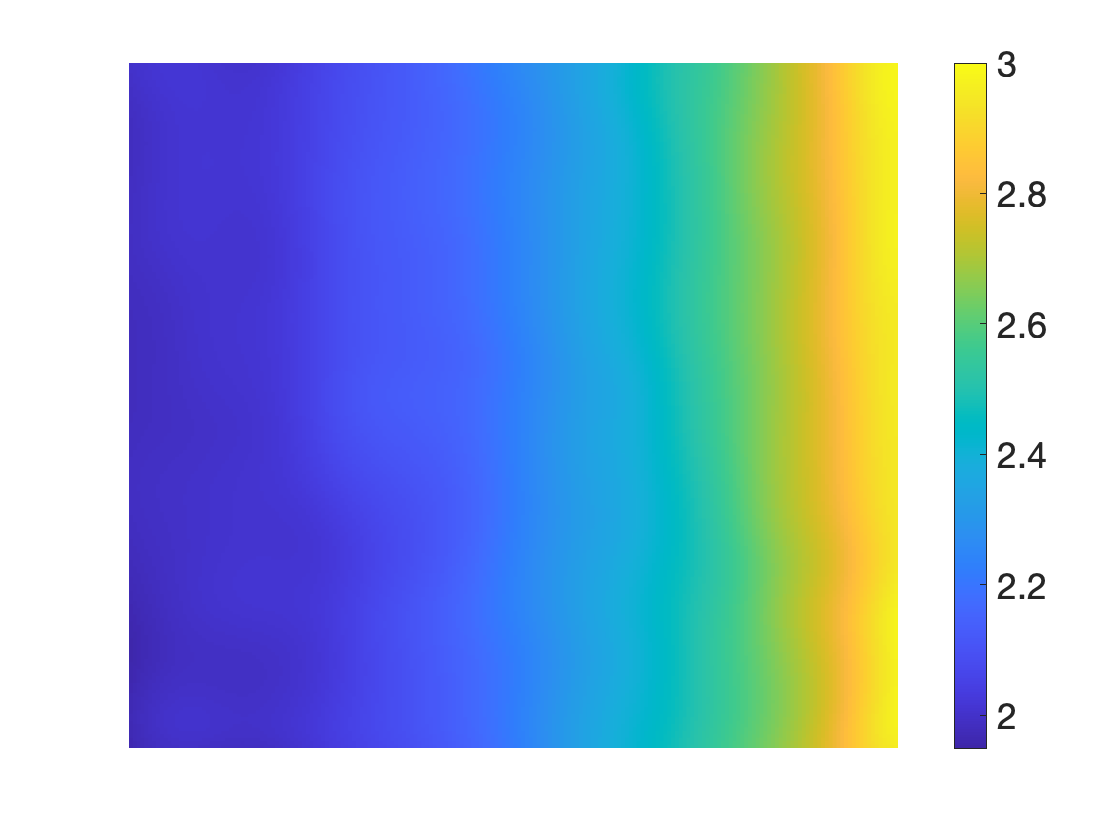} &
\includegraphics[width=0.199\textwidth]{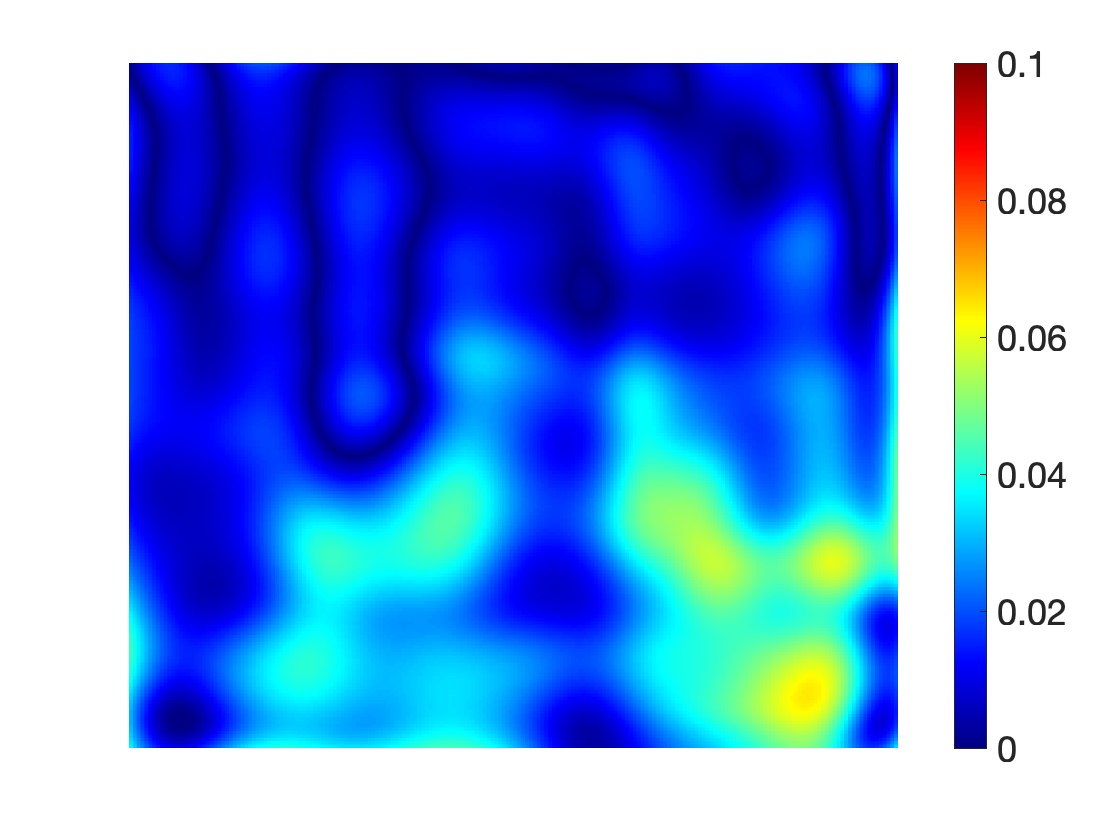}\\
\includegraphics[width=0.199\textwidth]{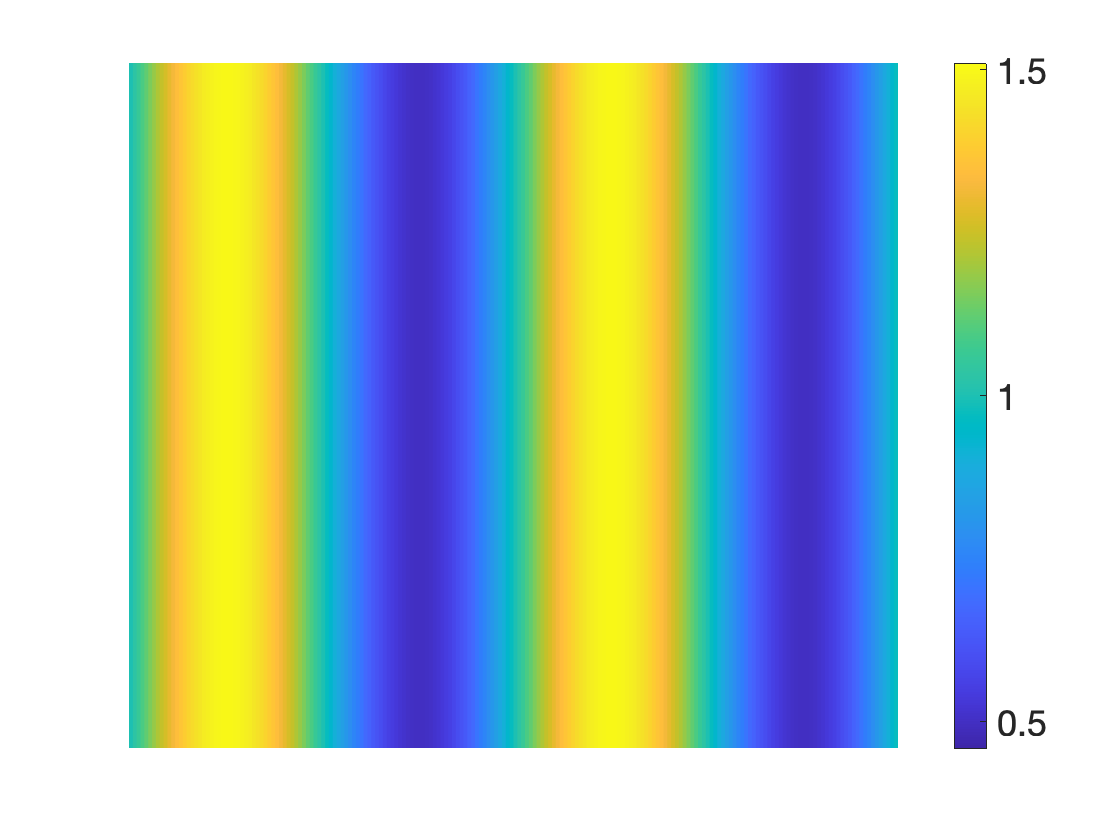} &
\includegraphics[width=0.199\textwidth]{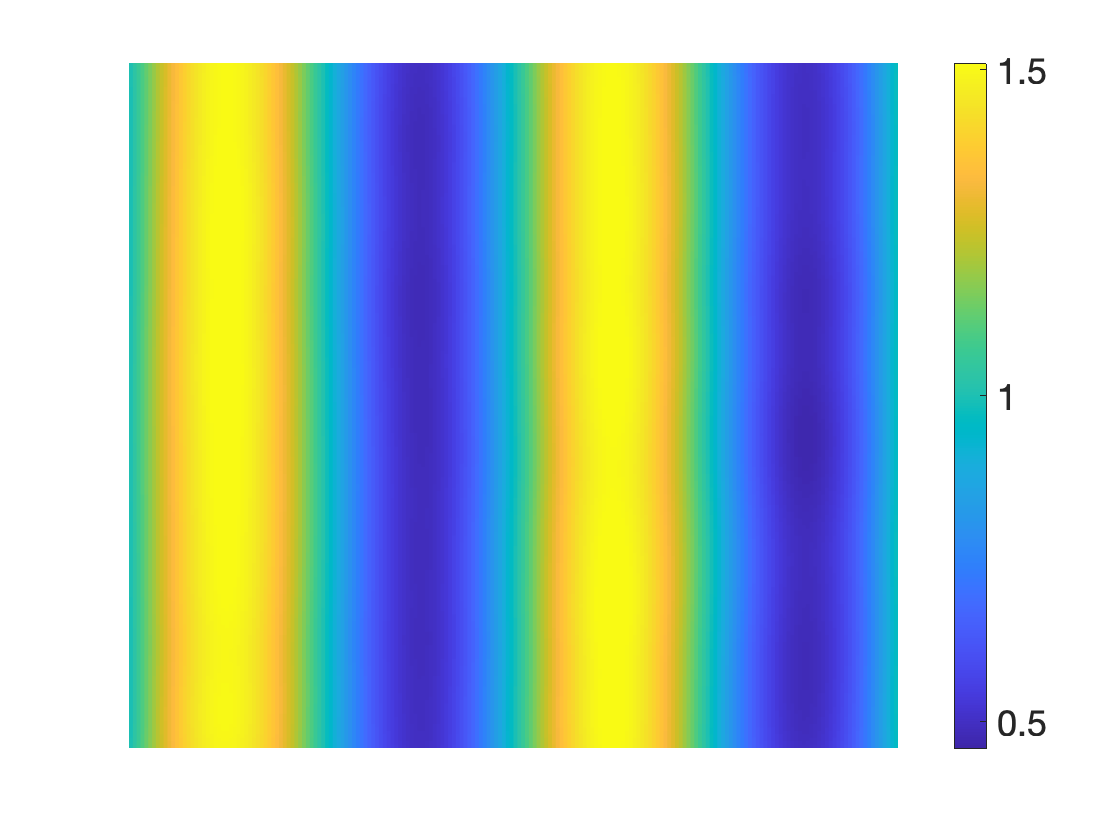} &
\includegraphics[width=0.199\textwidth]{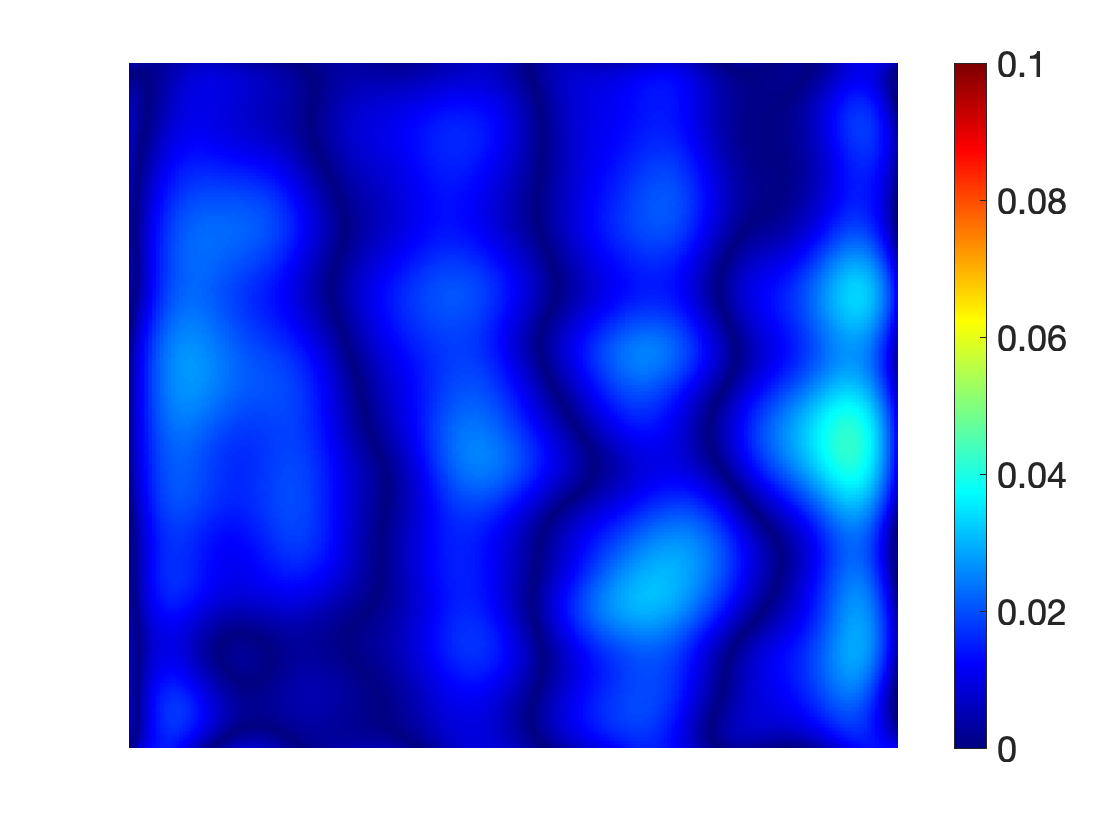} &
\includegraphics[width=0.199\textwidth]{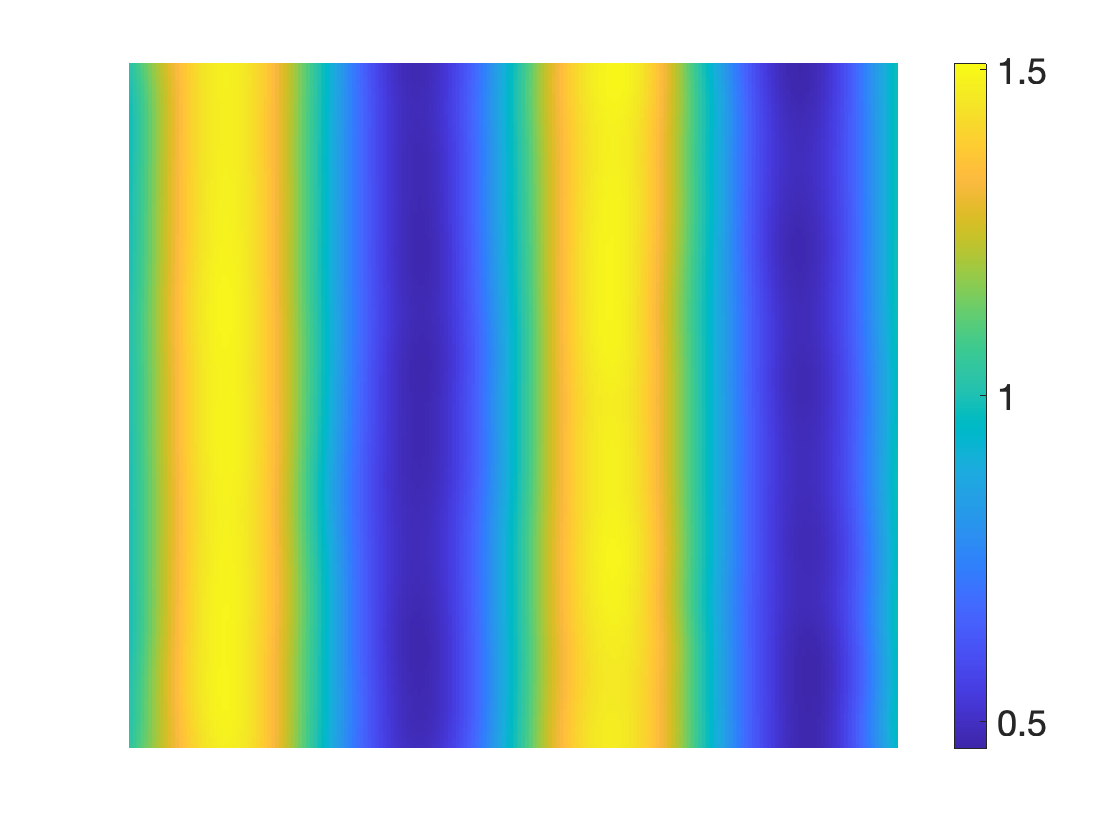} &
\includegraphics[width=0.199\textwidth]{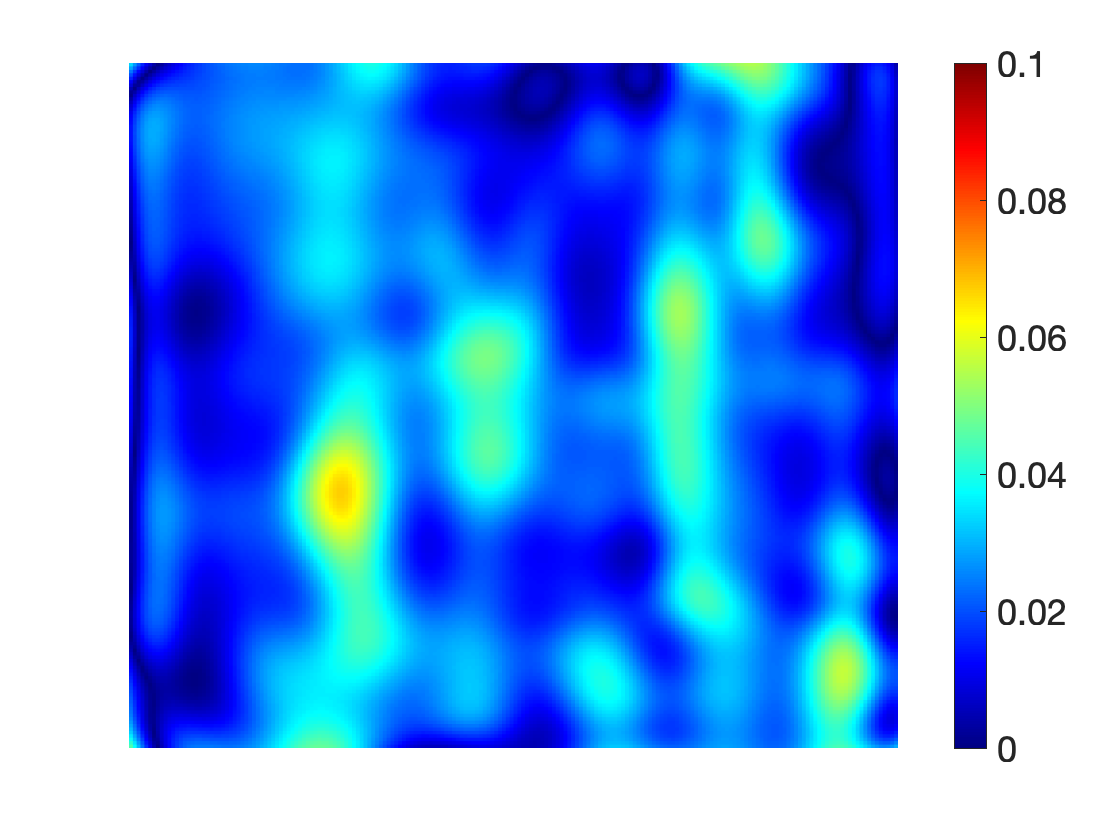}\\
\includegraphics[width=0.199\textwidth]{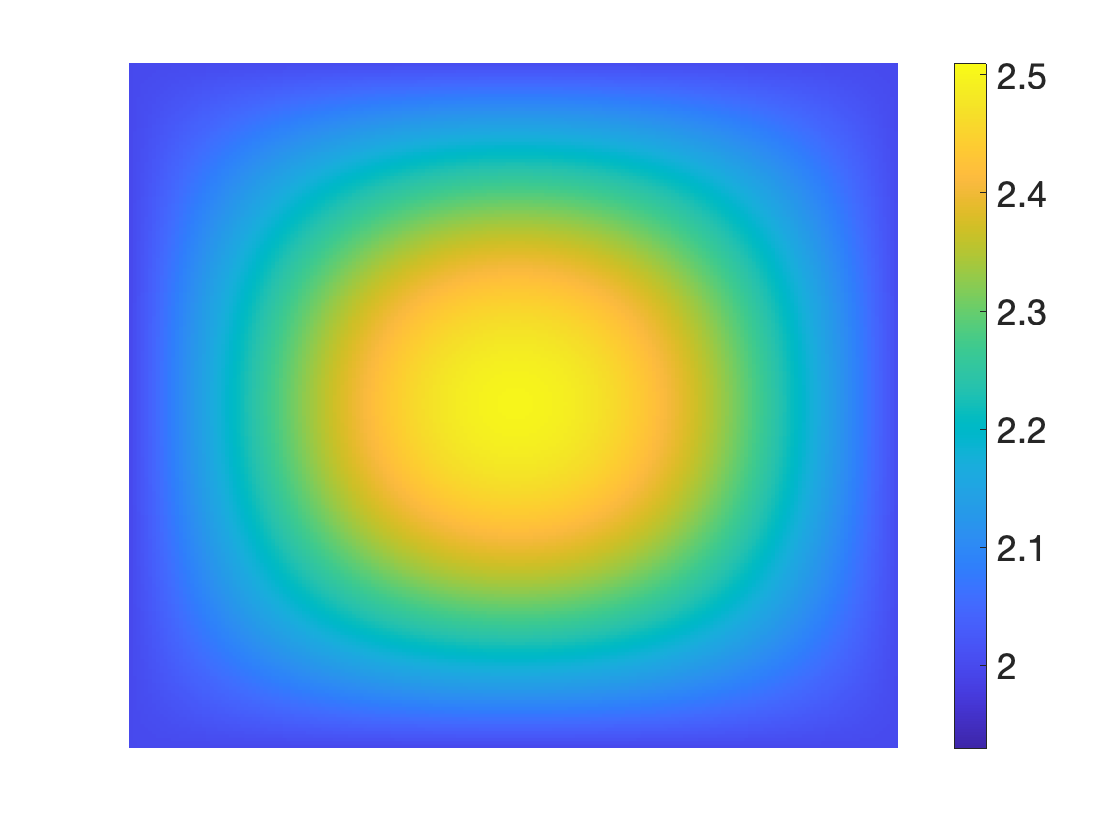} &
\includegraphics[width=0.199\textwidth]{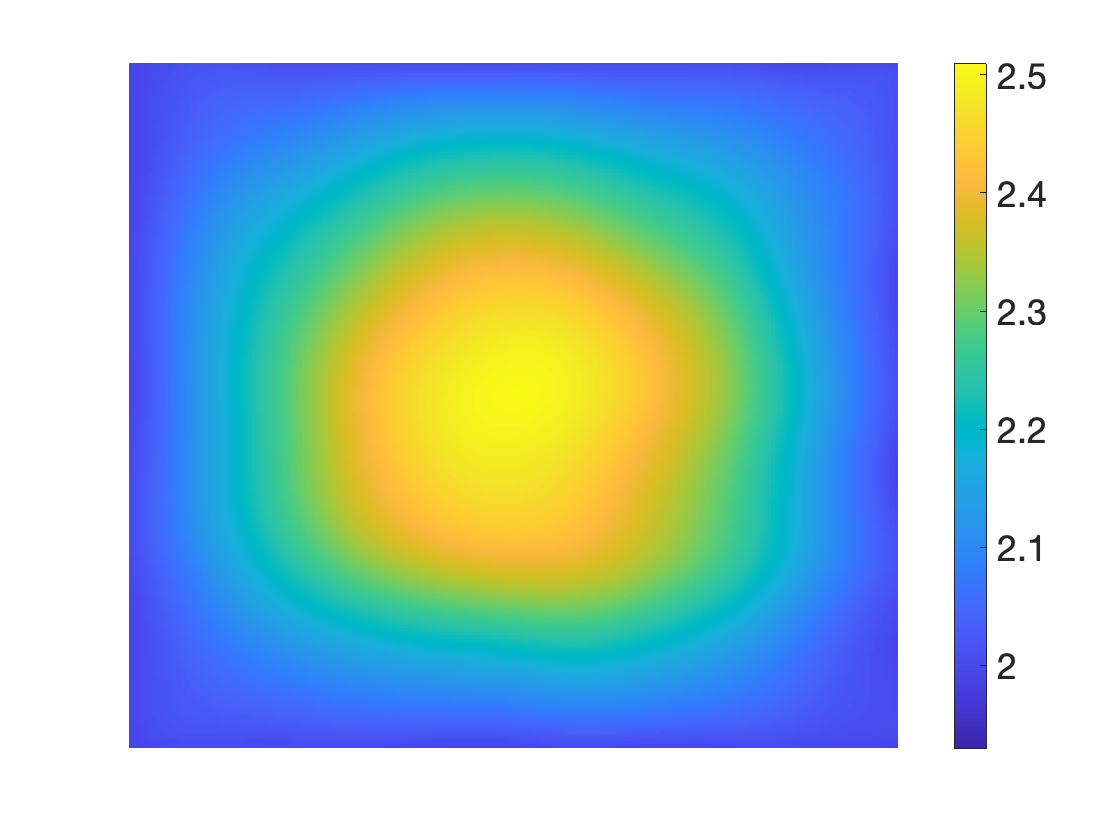} &
\includegraphics[width=0.199\textwidth]{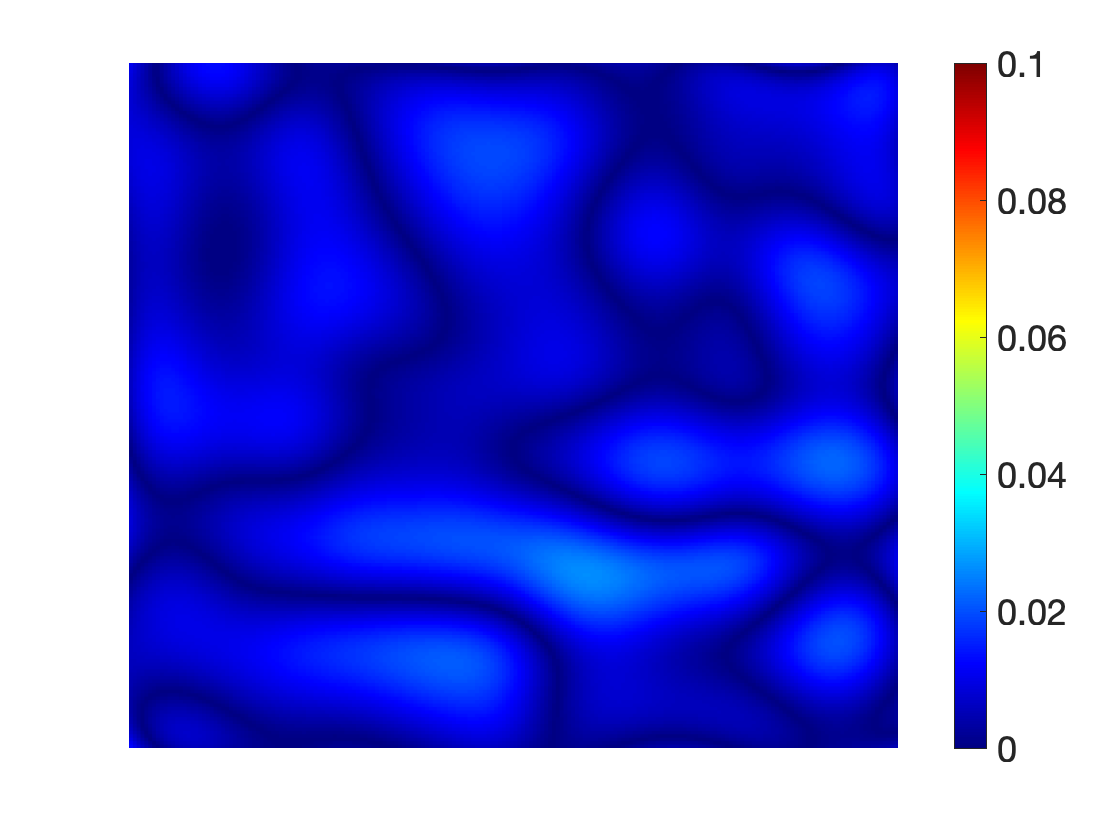} &
\includegraphics[width=0.199\textwidth]{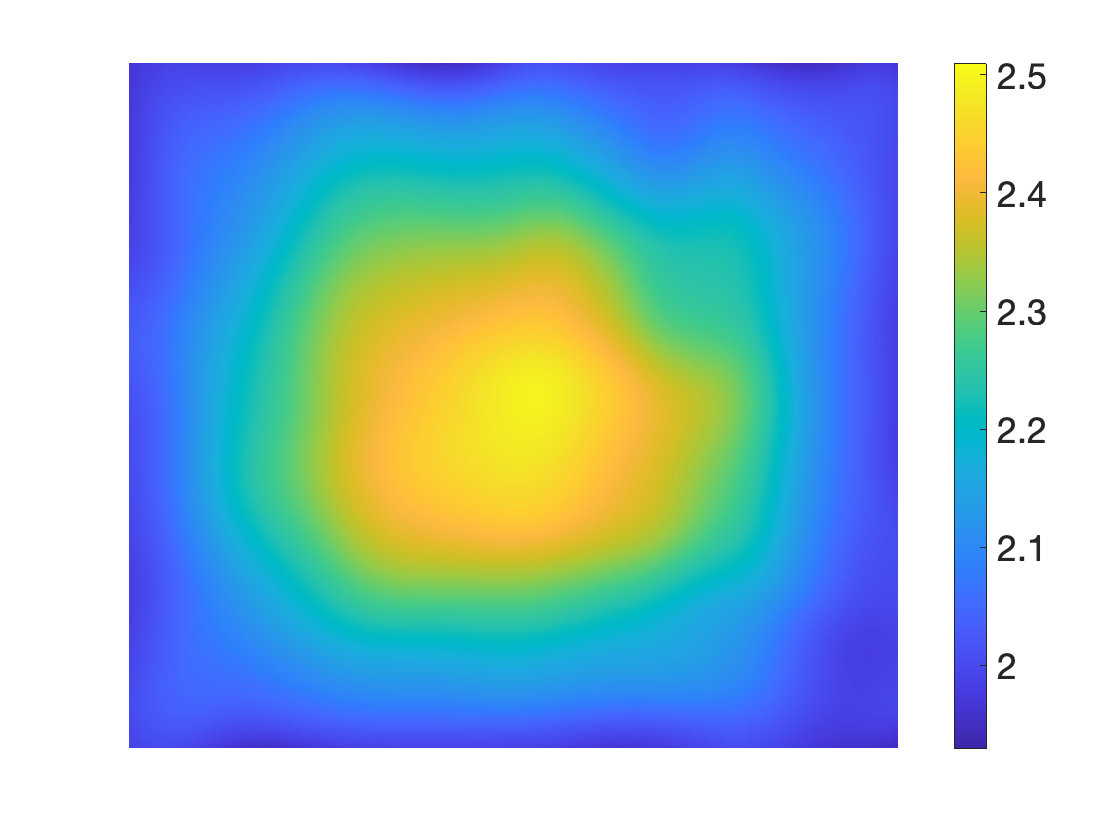} &
\includegraphics[width=0.199\textwidth]{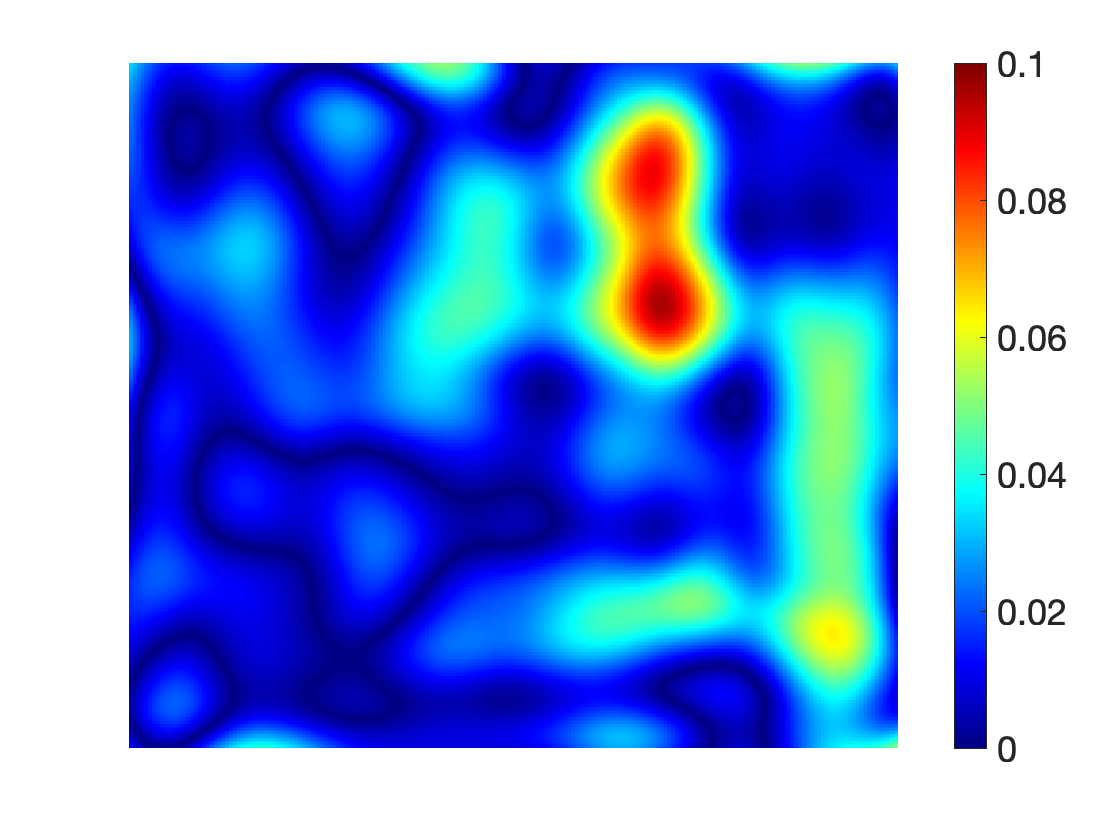}\\
(a) $A^\dag$  & (b) $\hat A$ & (c) $|\hat A-A^\dag|$ & (d) $\hat A$ & (e) $|\hat A-A^\dag|$
\end{tabular}
\caption{The reconstructions for Example \ref{exam:diri3d2} with exact data in (b) and noisy data $(\delta=5\%)$ in (d). From the top to bottom, the results are for $A_{11}$, $A_{12}$, $A_{13}$, $A_{22}$, $A_{23}$ and $A_{33}$, respectively.}
\label{fig:diri3d2}
\end{figure}

Fig. \ref{fig:diri3d2} shows the results on a 2D cross section at $x_3 = 0.5$. The overall features of $A^\dagger$ are accurately reconstructed across the entire domain $\Omega$, including the central region $\Omega \setminus \omega$ where there is no observational data. The accuracy does not deteriorate much for up to $5\%$ noise, showing only a slight deformation in the reconstruction of $A_{33}$. This example again illustrates the effectiveness of the DNN approach in the three-dimensional case and with partial internal data.

These experiments clearly demonstrate that the proposed DNN approach can accurately reconstruct the anisotropic conductivity in both two- and three-dimensional cases, for up to 10\% noise in the data, indicating its high robustness with respect to data noise. Also due to the mesh-free nature of the method, its implementation is direct and nearly identical for the 2D and 3D cases. Thus, the method holds potential for the concerned inverse problem.

\bibliographystyle{abbrv}
\bibliography{reference}
\end{document}